\documentclass[12pt]{article}

\usepackage[utf8]{inputenc}
\usepackage[russian, english]{babel}
 \usepackage[T2A]{fontenc}

\usepackage{amsmath, amssymb,latexsym,amsthm, epigraph}
\usepackage{authblk}

\usepackage{enumitem}
\usepackage{moreenum}

\usepackage{hyperref}

\usepackage{tikz}
\usetikzlibrary{decorations.markings}
\usetikzlibrary{decorations.pathreplacing}
\usetikzlibrary{shapes.geometric,positioning,calc}
\usetikzlibrary{intersections}
\usetikzlibrary{positioning}

\usepackage{caption}
\usepackage{array}

\usepackage{multicol}

\theoremstyle{definition}
\newtheorem{definition}{Definition}[section]
\newtheorem{remark}{Remark}[section]
\newtheorem{example}{Example}[section]
\newtheorem{condition}{Condition}
\newtheorem*{condition3a}{Condition 3a}
\newtheorem*{condition3b}{Condition 3b}

\theoremstyle{plain}
\newtheorem{theorem}{Theorem}
\newtheorem{lemma}{Lemma}[section]
\newtheorem{proposition}[lemma]{Proposition}
\newtheorem{corollary}[lemma]{Corollary}

\newcommand{\Fr}{\mathcal{F}}
\newcommand{\Ideal}{\mathcal{I}}
\newcommand{\Rel}{\mathcal{R}}
\newcommand{\Mon}{\mathcal{M}}
\newcommand{\SP}{\mathcal{S}}
\newcommand{\SPM}{\Lambda}
\newcommand{\fld}{\mathit{k}}

\newcommand{\Galg}{\fld\Fr}
\newcommand{\Qalg}{\fld\Fr / \Ideal}

\newcommand{\Dp}{\mathrm{Dp}}
\newcommand{\EDp}{\mathrm{EDp}}
\newcommand{\GDp}{\mathcal{T}}
\newcommand{\Ft}{\mathrm{F}}
\newcommand{\Low}{\mathrm{L}}

\newcommand{\Gr}{\mathrm{Gr}}

\newcommand{\MLow}{\text{$\mathrm{Lower}$-$\mathrm{f}$}}
\newcommand{\MUp}{\text{$\mathrm{Equal}$-$\mathrm{f}$}}

\newcommand{\Add}{\mathrm{Add}}
\newcommand{\GreedyAlg}{\mathrm{GreedyAlg}}

\newcommand{\DSpace}[1]{\langle{#1}\rangle_d}
\newcommand{\DMUp}[1]{\MUp(#1)_d}
\newcommand{\DMLow}[1]{\MLow(#1)_d}

\newcommand{\mo}[1]{\mathsf{Max}({#1})}
\newcommand{\nfc}[1]{\mathsf{Max}^{\mathsf{nfc}}({#1})}
\newcommand{\fc}[1]{\mathsf{Max}^{\mathsf{fc}}({#1})}
\newcommand{\virt}[1]{\mathcal{V}({#1})}

\newcommand{\shortmo}[1]{\mathsf{Max}^{\leqslant 2}({#1})}
\newcommand{\longmo}[1]{\mathsf{Max}^{\geqslant 3}({#1})}

\newcommand{\mincov}[1]{\mathsf{MinCov}({#1})}
\newcommand{\mincovk}[1]{\mathsf{MinCov}_k({#1})}
\newcommand{\nvirt}[1]{\mathsf{NVirt}({#1})}

\newcommand{\dbtilde}[1]{\widetilde{\raisebox{0pt}[0.85\height]{$\widetilde{#1}$}}}
\newcommand{\dbhat}[1]{\widehat{\raisebox{0pt}[0.85\height]{$\widehat{#1}$}}}

\tikzset{arrow/.style={
        decoration={markings,
            mark= at position #1 with {\arrow{stealth}},
        },
        postaction={decorate}
    }
}

\tikzset{bigarrow/.style={
        decoration={markings,
            mark= at position #1 with {\arrow[scale=1.5]{stealth}},
        },
        postaction={decorate}
    }
}

\tikzset{reversearrow/.style={
        decoration={markings,
            mark= at position #1 with {\arrow{stealth reversed}},
        },
        postaction={decorate}
    }
}

\tikzset{reversebigarrow/.style={
        decoration={markings,
            mark= at position #1 with {\arrow[scale=1.5]{stealth reversed}},
        },
        postaction={decorate}
    }
}

\title{Group-like Small Cancellation Theory for Rings}

\author{A.\,Atkarskaya}
\affil{Department of Mathematics, Bar-Ilan University, 5290002 Ramat Gan, Israel}
\affil{Department of Mathematics, The Hebrew University of Jerusalem, Givat Ram, 9190401 Jerusalem, Israel \\ atkarskaya.agatha@gmail.com}

\author{A.\,Kanel-Belov}
\affil{Department of Mathematics, Bar-Ilan University, 5290002 Ramat Gan, Israel}
\affil{Department of Discrete Mathematics, Moscow Institute of Physics and Technology, Dolgoprudnyi, Institutskiy Pereulok, 141700 Moscow Oblast, Russia}
\affil{College of Mathematics and Statistics, Shenzhen University, Shenzhen 518061, China \\ beloval@macs.biu.ac.il}

\author{E.\,Plotkin}
\affil{Department of Mathematics, Bar-Ilan University, 5290002 Ramat Gan, Israel \\ plotkin.evgeny@gmail.com}

\author{E.\,Rips}
\affil{Department of Mathematics, The Hebrew University of Jerusalem, Givat Ram, 9190401 Jerusalem, Israel \\ eliyahu.rips@mail.huji.ac.il}

\date{}

\begin{document}
\maketitle

\setlength{\epigraphwidth}{0.83\textwidth}
\epigraph{
\smallskip
\begin{multicols}{2}
{\it In memoriam of \\
 Sergei Ivanovich Adian \\
 whose groundbreaking work \\
 remains a permanent source \\
 of mathematical inspiration}\\
 \newpage
{\it Памяти Сергея Ивановича Адяна чья первопроходческая деятельность является постоянным источником математического вдохновения}\\
 \end{multicols}
 {\it
\ }}


\begin{abstract}
In the present paper we develop a small cancellation theory for associative algebras with a basis of invertible elements. Namely, we study quotients of a group algebra of a free group and introduce three axioms for the corresponding defining relations. We show that the obtained ring is non-trivial. Moreover, we show  that this ring enjoys a global filtration that agrees with relations, find a basis of the ring as a vector space  and establish the corresponding structure theorems. We also provide a revision of a concept of  Gr\"{o}bner basis for our rings and establish a greedy algorithm for the Ideal Membership Problem.
\end{abstract}

\medskip


\medskip

\tableofcontents

\section{Introduction}\label{intro}

For any variety of algebraic systems, one can define a specific algebraic system by generators and defining relations~\cite{Co}. When the free algebra of this variety is well enough understood, we may pose the following (vague) question:

\medskip

{\it If the interactions among the defining relations are weak in a certain sense, 
will the resulting algebra bear some resemblance to the free algebra?}

\medskip

For instance,  may it be non-trivial,  may it have a reasonable structure theory,  may it have  a solvable equality problem or even may its equality problem enjoy a greedy algorithm.

In the case of groups (or semigroups and monoids) the Small Cancellation Theory provides an answer to this question. In the present paper we develop a similar theory for associative algebras with a basis of invertible elements.  In fact, in course of studying the question:
\medskip

\centerline{\it ``what is a small cancellation  associative ring?"}

\medskip
\noindent
we axiomatically define  a ring, which can reasonably be called a ring with  small cancellation. We also determine the structure and properties of this ring.

General theory presented in this paper is modeled after a particular case we have treated in our previous paper~\cite{AKPR}.

\subsection{Motivation, objectives, results}
\label{intro_motivations_section}

The motivation for developing a ring-theoretical analog of small cancellation comes from the fact that small cancellation for groups and, especially, its more far-reaching versions, provides a very powerful technique for constructing groups with unusual, and even exotic, properties, like for example,  infinite Burnside groups \cite{NA1}--\cite{NA3}, \cite{A}, \cite{Ol3},\cite{Ivanov}, \cite{Lys},\cite{DelGr},\cite{Co},  Tarski monster \cite{Ol2}, finitely generated infinite divisible groups~\cite{Gu},  and many others~\cite{Ol1}.

On the other hand, there is a conceptual desire to understand what negative curvature could mean for ring theory.

For any group with the fixed system of generators, its Cayley graph can be considered as a metric space.  This leads to Gromov's program ``Groups as geometric objects" \cite{Gr1}, see also \cite{Gr2}. In particular, a finitely generated group is word-hyperbolic when its Cayley graph is $\delta$-hyperbolic for $\delta > 0$ (see \cite{Bo}, \cite{DK} for modern exposition and references).

So far, we do not know a way to associate a geometric object to a ring. Thus, having in mind the negative curvature as a heuristic and indirect hint for our considerations, we, nevertheless,  follow a more accessible  combinatorial line of studying rings. Therefore, small cancellation groups appear naturally on the stage.   

Finitely generated small cancellation groups turned out to be word hyperbolic (when every relation needs at least 7 pieces). So, if we could generalize small cancellation to the ring theoretic situation, it would provide examples to the yet undefined concept of a ring with a negative curvature. Another source  of potential examples should be group algebras of hyperbolic groups.

Following this reasoning,  we introduce in the paper the small cancellation axioms for rings. We develop the theory of rings that satisfy these axioms. It is a major problem to establish that the resulting quotient is non-zero. Only in the end of the paper we are in a position to verify this claim. We also construct an explicit linear basis for such rings. In parallel we show that the equality problem in our ring  is solvable. In fact it possesses algorithmic properties similar to ones for groups with small cancellation.

Group algebras of groups can be defined by ring presentations in which all polynomials are in fact binomials. If we take small cancellation groups with appropriate constants then their group algebras  satisfy  Small Cancellation Axioms, see Subsection~\ref{group_algebras_of_sc_groups_section}. Surface groups are among the well-known examples of  small cancellation groups and, therefore, their group algebras provide examples of rings with small cancellation.

We do hope that being of interest as a ring of new type by itself, this ring will inherit useful practical properties known for small cancellation groups and, thus, it can be used for obtaining complicated algebras with the very specific properties.

We shall note  another direction in the construction of algebras with unusual properties. This refers to the breakthrough made by A.~Smoktunowicz. Her innovatory approach to controlling relations in rings led to the construction of a simple nil-ring and other important examples of nil algebras (see \cite{Sm1}, \cite{Sm2}, \cite{LSm},\cite{LSmY}).


\subsection{Overview of the work}
\label{intro_overview_section}

In the present paper we develop a small cancellation theory for associative algebras with a basis of invertible elements. In this introductory overview we explain why such a generalization exists at all, why it turns out to be so technical and difficult, and what is the strategy to overcome difficulties.

\subsubsection{Small cancellation for groups}
\label{intro_smallgroups_section}

The small cancellation theory for groups, in its simplest version, is quite an elementary theory (see~\cite{LS}), and let us shortly outline it here. Consider a group presentation $G = \langle \mathcal{X} \mid  \Rel\rangle$ where we assume that the set of relations $\Rel$ is closed under cyclic permutations and inverses.

The interaction between the defining relations is described in terms of small pieces (see the definition below).

Throughout the paper we use the following notations. Let $\Fr$ be the free group with the set of free generators~$\mathcal{X}$. Assume $A$ and $B$ are two elements of $\Fr$. We write the product of the monomials $A$ and $B$ as $A\cdot B$. There may occur cancellations between $A$ and $B$ in $A\cdot B$. We write the product of $A$ and $B$ in the form $AB$ when there are no cancellations in~$A\cdot B$.

A word $s \in \Fr$ is called a \emph{small piece with respect to $\Rel$} if there are relations of the form $l_1 s r_1$ and  $l_2 s r_2$ in $\Rel$ such that  $r_1 l_1 \neq r_2 l_2$ as words in the free group, and any conjugate of $(r_1 l_1) \cdot (l_2^{-1} r_2^{-1})$ is not contained in $\Rel$, after possible cancellations (cf.~\cite{R}).

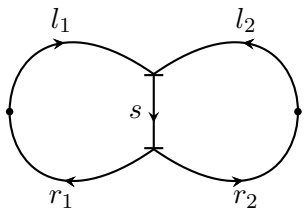
\begin{center}
\begin{tikzpicture}
\tikzstyle{every node}=[font=\small];
\draw[|-|, black, thick, reversearrow=0.5] (0,0)--(0,1) node [midway, left] {$s$};
\draw[black, thick, arrow=0.25, arrow=0.75] (0,0) to [bend left=125, looseness=8] coordinate[pos=0.5] (P1) (0,1);
\draw[black, thick, arrow=0.25, arrow=0.75] (0,0) to [bend right=125, looseness=8] coordinate[pos=0.5] (P2) (0,1);
\path (0,0) to [bend left=125, looseness=8] node[pos=0.2, below] {$r_1$} node[pos=0.8, above] {$l_1$} (0,1);
\path (0,0) to [bend right=125, looseness=8] node[pos=0.2, below] {$r_2$} node[pos=0.8, above] {$l_2$} (0,1);
\node[circle, fill, inner sep=1] at (P1) {};
\node[circle, fill, inner sep=1] at (P2) {};
\end{tikzpicture}
\end{center}

\begin{remark}
The geometric way to think about small pieces is seeing them as words that may appear on the common boundary between two cells in the van Kampen diagram~\cite{LS},~\cite{Ol1}. In particular, if $(r_1l_1)\cdot (l_2^{-1}r_2^{-1}) \in \Rel$, then we can substitute these cells by a simple cell, so we are entitled to assume from the beginning that $(r_1l_1)\cdot (l_2^{-1}r_2^{-1}) \notin \Rel$.
\end{remark}

The \emph{small cancellation condition} says that any relation in $\Rel$ cannot be written as a product of too few small pieces. For most of purposes seven small pieces suffice, since the discrete Euler characteristic per cell becomes negative~\cite{LS},~\cite{Ly}.

To ensure this,  we can assume that the length of any small piece is less than one sixth of the length of the whole relation in which it appears.

Note that in the simpler version of the Small Cancellation Theory the last condition $(r_1 l_1)\cdot (l_2^{-1} r_2^{-1}) \notin \Rel$ is omitted. So in the simpler version more words  are qualified  as small pieces and, correspondingly, this class is smaller.

Every element of $G$ can be presented by a word in the generators $\mathcal{X}$ that does not contain occurrences of more than a half of one of the relations (otherwise there is a shorter word presenting the same element). However this presentation need not be unique. The question is to which extent it is not unique. The Main Theorem of Small Cancellation Theory provides the following answer.

Let $w_1$, $w_2$ be two words that do not contain occurrences of more than a half of a relation. They represent the same element of $G$ if and only if they can be connected by a one-layer diagram (\cite{LS}, especially see Greendlinger's Lemma).

As an illustration, consider the following example
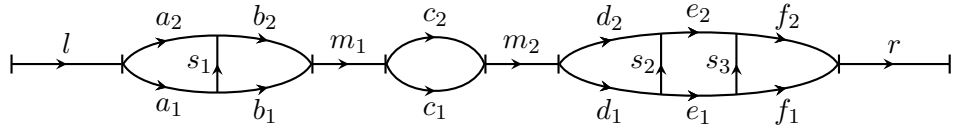
\begin{center}
\begin{tikzpicture}
\tikzstyle{every node}=[font=\small];
\draw[|-|, black, thick, arrow=0.5] (0.5,0)--(2,0) node [midway, above] {$l$};
\draw[black, thick, arrow=0.25, arrow=0.75] (2,0) to [bend left=60, looseness=0.6] coordinate[pos=0.5] (P1)  (4.5,0);
\draw[black, thick,  arrow=0.25, arrow=0.75] (2,0) to [bend right=60, looseness=0.6] coordinate[pos=0.5] (P2) (4.5,0);
\path (2,0) to [bend left=60, looseness=0.6] node[pos=0.3, above] {$a_2$} node[pos=0.7, above] {$b_2$} (4.5,0);
\path (2,0) to [bend right=60, looseness=0.6] node[pos=0.3, below] {$a_1$} node[pos=0.7, below] {$b_1$}  (4.5,0);
\path (2,0) to [bend right=60, looseness=0.6] (4.5,0);
\draw[black, thick, arrow=0.5] (P2)--(P1) node[midway, left, xshift=2.5] {$s_1$};
\draw[|-|, black, thick, arrow=0.5] (4.5,0)--(5.5,0) node [midway, above] {$m_1$};
\draw[black, thick, arrow=0.5] (5.5,0) to [bend left=70] node [above] {$c_2$} (6.8,0);
\draw[black, thick, arrow=0.5] (5.5,0) to [bend right=70] node [below] {$c_1$} (6.8,0);
\draw[|-|, black, thick, arrow=0.5] (6.8,0)--(7.8,0) node [midway, above] {$m_2$};
\draw[black, thick, arrow=0.2, arrow=0.5, arrow=0.8] (7.8,0) to [bend left=60, looseness=0.45] coordinate[pos=0.4] (P3) coordinate[pos=0.6] (P5) (11.5,0);
\draw[black, thick,  arrow=0.2, arrow=0.5, arrow=0.8] (7.8,0) to [bend right=60, looseness=0.45] coordinate[pos=0.4] (P4) coordinate[pos=0.6] (P6) (11.5,0);
\path (7.8,0) to [bend left=60, looseness=0.45] node[pos=0.25, above] {$d_2$} node[pos=0.5, above] {$e_2$} node[pos=0.75, above] {$f_2$}(11.5,0);
\path (7.8,0) to [bend right=60, looseness=0.45] node[pos=0.25, below] {$d_1$} node[pos=0.5, below] {$e_1$} node[pos=0.75, below] {$f_1$} (11.5,0);
\draw[black, thick, arrow=0.5] (P4)--(P3) node[midway, left, xshift=2.5] {$s_2$};
\draw[black, thick, arrow=0.5] (P6)--(P5) node[midway, left, xshift=2.5] {$s_3$};
\draw[|-|, black, thick, arrow=0.5] (11.5,0)--(13,0) node [midway, above] {$r$};
\end{tikzpicture}
\end{center}

\noindent
where
\begin{align*}
&w_1=la_1 b_1 m_1 c_1 m_2d_1e_1f_1r,\\
&w_2=la_2 b_2 m_1 c_2 m_2d_2e_2f_2r.
\end{align*}

\noindent
and $a_1s_1a_2^{-1}$, $b_1b_2^{-1}s_1^{-1}$,  $c_1c_2^{-1}$, $d_1s_2d_2^{-1}$, $e_1s_3 e_2^{-1}s_2^{-1}$, $f_1f_2^{-1}s_3^{-1}$ are relations, and $s_1$, $s_2$, $s_3$ are small pieces.

The transition from $w_1$ to $w_2$ can be divided into a sequence of  elementary steps called \emph{turns} (\cite{NA1}---\cite{NA3}). Each turn reverses just  one cell. For example,
\begin{equation*}
w_1 = v_0\mapsto v_1\mapsto v_2\mapsto v_3\mapsto v_4\mapsto v_5\mapsto v_6 = u,
\end{equation*}
where
\begin{align*}
&v_0 = w_1 = la_1 b_1 m_1 c_1 m_2d_1e_1f_1r,\\
&v_1 = la_2s_1^{-1}b_1m_1c_1m_2d_1e_1f_1r,\\
&v_2 = la_2b_2m_1c_1m_2d_1e_1f_1r,
\end{align*}
etc. The turns can be performed in any order.

In general, if $m = pq\in \Rel$, $w = lpr$, then $w = lpr \mapsto lq^{-1}r = w^{\prime}$ is a turn.

\subsubsection{Main definitions and examples for the ring case}
\label{intro_defex_section}

Since we start with small cancellation of groups (or, not, say, of semigroups and monoids, see for example~\cite{Hig},~\cite{GS},~\cite{S}), we deal with quotients of the group algebra of a free group and not with quotients of a free associative algebra.

Given a field $\fld$ and the  free group $\Fr$, denote by $\Galg$ the corresponding group algebra. Let a set of polynomials $\Rel$ from $\Galg$ be fixed. Define $\Ideal$ to be the ideal generated by the elements of $\Rel$ as an ideal.

As above, we fix the following notations. Assume $X$ and $Y$ are two elements of $\Galg$. We write their product as $X \cdot Y$. Assume $A$ and $B$ are two elements of $\Fr$. We write the product of the monomials $A$ and $B$ as $A\cdot B$. There may occur cancellations between $A$ and $B$ in $A\cdot B$. We write the product of $A$ and $B$ in the form $AB$ when there are no cancellations in~$A\cdot B$.

First of all we introduce three types of conditions on elements of $\Rel$ called Compatibility Axiom, Small Cancellation Axiom and Isolation Axiom. These restrictions are our analogue of Small Cancellation Condition for group presentations (see~\cite{LS}).

Now we indicate counterparts of the notions of a small piece and of a turn.

Our definition of \emph{a small piece for presentations of associative algebras with a basis of invertible elements} is based  on the definition of a small piece for group presentations mentioned above. However, this generalization is not at all straightforward; rather it constitutes a major conceptual novelty (see Definition~\ref{sp}).

The notion of \emph{turn} is replaced by the notion of \emph{multi-turn} (Definition~\ref{multiturn_def}). Throughout the paper we reserve small Greek letters for non-zero elements of the field~$\fld$. Let $\sum_{j = 1}^n \alpha_j m_j \in \Rel$, with all $\alpha_j \neq 0$ given. Let $v$ be a monomial of the form $v = l m_h r$ for some $h$, $ 1 \leqslant h \leqslant n$. The transition from $v = l m_h r$ to $\sum_{\substack{j = 1\\j \neq h}}^n ( -\alpha_h^{-1} \alpha_j) l m_j r$ is called \emph{a multi-turn}. It is extended by linearity to $\beta v$, and then extended to polynomials, such that one of their terms is of the form $\beta v$. The corresponding polynomial $\sum_{j = 1}^n\alpha_jl m_j r$ is called \emph{the layout of this multi-turn}.

\noindent
{\bf  Examples.}
In our examples we assume for simplicity that the ground field is the field with two elements.

\medskip

A. Let $v = la_1r$ and let our polynomial be $a_1 + a_2$ from $\Rel$. In this case we have a transition from $la_1r$ to $la_2r$ (see the picture)
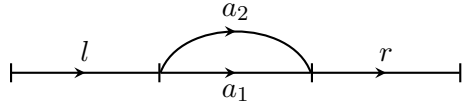
\begin{center}
\begin{tikzpicture}
\tikzstyle{every node}=[font=\small];
\draw[|-|, black, thick, arrow=0.5] (0,0)--(2,0) node [midway, above] {$l$};
\draw[black, thick, arrow=0.5] (2,0) to [bend left=70] node [above] {$a_2$} (4,0);
\draw[black, thick, arrow=0.5] (2,0) to node [below] {$a_1$} (4,0);
\draw[|-|, black, thick, arrow=0.5] (4,0)--(6,0) node [midway, above] {$r$};
\end{tikzpicture}
\end{center}
The transition from $l a_1 r$ to $l a_2 r$ shows that a turn is a particular case of a multi-turn.

\medskip

B. Let $v = la_1r$ and let our polynomial be $a_1 + a_2 + a_3$ from $\Rel$. Then we have a transition from $la_1r$ to $la_2r + la_3r$. We use the picture
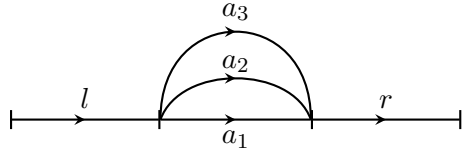
\begin{center}
\begin{tikzpicture}
\tikzstyle{every node}=[font=\small];
\draw[|-|, black, thick, arrow=0.5] (0,0)--(2,0) node [midway, above] {$l$};
\draw[black, thick, arrow=0.5] (2,0) to [bend left=70] node [above, yshift=-2] {$a_2$} (4,0);
\draw[black, thick, arrow=0.5] (2,0) to [bend left=90, looseness=2] node [above] {$a_3$} (4,0);
\draw[black, thick, arrow=0.5] (2,0) to node [below] {$a_1$} (4,0);
\draw[|-|, black, thick, arrow=0.5] (4,0)--(6,0) node [midway, above] {$r$};
\end{tikzpicture}
\end{center}

\medskip

C. Let the monomial $v = lqa_1q^{-1}r$ be given, and let our polynomial be $a_1 + a_2 + 1$ from $\Rel$. Then the multi-turn is a transition from $l q a_1 q^{-1} r$ to $l q a_2 q^{-1} r + l\cdot r$. Note that after insertion of $1$ instead of $a_1$ the factor $q$ cancels out with $q^{-1}$. The $1$ is not reflected in the picture.
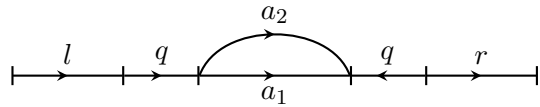
\begin{center}
\begin{tikzpicture}
\tikzstyle{every node}=[font=\small];
\draw[|-|, black, thick, arrow=0.5] (0,0)--(1.5,0) node [midway, above] {$l$};
\draw[-|, black, thick, arrow=0.5] (1.5,0)--(2.5,0) node [midway, above] {$q$};
\draw[black, thick, arrow=0.5] (2.5,0) to [bend left=70] node [above] {$a_2$} (4.5,0);
\draw[black, thick, arrow=0.5] (2.5,0) to node [below] {$a_1$} (4.5,0);
\draw[|-, black, thick, reversearrow=0.5] (4.5,0)--(5.5,0) node [midway, above] {$q$};
\draw[|-|, black, thick, arrow=0.5] (5.5,0)--(7,0) node [midway, above] {$r$};
\end{tikzpicture}
\end{center}

\medskip

D. Let the monomial $v = la_1b_1r$ be given. Consider the polynomials  $a_1 + a_2 + a_3 s_1^{-1}$ and $b_1 + b_2 + 1$ from $\Rel$. We have two adjacent multi-turns: $a_1$ is replaced by $a_2 + a_3 s_1^{-1}$ and $b_1$ is replaced by $s_1 b_2 + 1$.  We have a picture
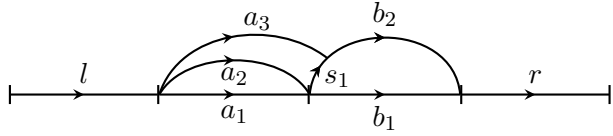
\begin{center}
\begin{tikzpicture}
\tikzstyle{every node}=[font=\small];
\draw[|-|, black, thick, arrow=0.5] (0,0)--(2,0) node [midway, above] {$l$};
\draw[black, thick, arrow=0.5] (2,0) to [bend left=60, looseness=0.9] node [below, yshift=1.5] {$a_2$} (4,0);
\draw[-|, black, thick, arrow=0.5] (2,0) to node [below] {$a_1$} (4,0);
\draw[black, thick, arrow=0.5] (4,0) to node [below] {$b_1$} (6,0);
\draw[black, thick, arrow=0.5, arrow=0.15] (4,0) to [bend left=85, looseness=1.3] node [above] {$b_2$} (6,0);
\path (4,0) to [bend left=85, looseness=1.3] coordinate[pos=0.2] (S) (6,0);
\path (4,0)--(S) node[midway, right, yshift=-1, xshift=-2] {$s_1$};
\draw[black, thick, arrow=0.5] (2,0) to [bend left=60, in=130] node [above, midway, xshift=10, yshift=-1.5] {$a_3$} (S);
\draw[|-|, black, thick, arrow=0.5] (6,0)--(8,0) node [midway, above] {$r$};
\end{tikzpicture}
\end{center}
Performing both multi-turns, we obtain that $l a_1 b_1 r$ is replaced by $l a_2 b_1 r + l a_3 s_1^{-1}b_1 r$ and then by $l a_2 s_1 b_2 r + l a_2 \cdot r + l a_3 b_2 r + l a_3 s_1 ^{-1}\cdot r$. There may be further cancellations in the monomials $l a_2 r$ and $l a_3 s_1 ^{-1} r$. Performing the multi-turns in the opposite order, we obtain that $l a_1 b_1 r$ is replaced by $l a_1 s_1 b_2 r + l a_1\cdot r$, and then by $ la_2 s_1 b_2 r + l a_3 b_2 r + l a_2 \cdot r + l a_3 s_1 ^{-1} \cdot r$, also with further cancellations. Notice that we needed to modify the second multi-turn before applying it, and notice that the result does not depend on the order of performing the multi-turns.

These examples show that, in comparison with the group case, a number of new phenomena occur. For instance, in Examples C and D unavoidable cancellations may happen. Occurrences of short words cause additional difficulties.

Now let us explain how we cope with this situation. Like in the group case, we introduce the notion of \emph{the chart} of a monomial. The chart of a monomial is a distinguished set of occurrences of its submonomials from a given set. In this paper, unless specified otherwise, all these occurrences are maximal occurrences of monomials from the set $\Mon$ of the monomials appearing in the set $\Rel$ of polynomial relations (see page~\pageref{mon_def}). Notice that according to the Compatibility Axiom, the set $\Mon$ is closed under taking submonomials. In our setting, the multi-turns we consider always start with an appropriate member of the chart of the monomial. Given a monomial, its chart, and a set of multi-turns applicable to this monomial, we define a new central notion of the \emph{space of linear dependencies} associated with this monomial and the set of multi-turns. If we apply one multi-turn generated by $\sum _{j = 1}^n \alpha_j m_j \in \Rel$ to the monomial $v = lm_hr$, then we consider the linear subspace of $\Galg$ spanned by the monomials $lm_jr,$  $j = 1,\ldots ,n$, after possible cancellations. In the image of this subspace in $\Qalg$, the linear dependence $\sum _{j = 1}^n \alpha_j l m_j r$ holds. If we apply several multi-turns to the monomial $v$, then we consider all monomials we obtain from $v$, and the linear dependencies between them in $\Qalg$ induced by all these multi-turns. In our last example D, we have the monomials
\begin{align*}
&la_1b_1r, \quad la_2b_1r, \quad la_3s_1^{-1}b_1r,\\
&la_1s_1b_2r, \quad la_2s_1b_2r, \quad la_3b_2r,\\
&la_1\cdot r, \quad la_2 \cdot r, \quad la_3 s_1^{-1} \cdot r
\end{align*}
after appropriate cancellations in the last three monomials and the linear dependencies
\begin{align*}
&la_1b_1r + la_2b_1r + la_3s_1^{-1}b_1r = 0,\\
&la_1s_1b_2r + la_2s_1b_2r + la_3b_2r = 0,\\
&la_1\cdot r + la_2 \cdot r + la_3s_1^{-1} \cdot r = 0,\\
&la_1b_1r + la_1s_1b_2r + la_1 \cdot r = 0,\\
&la_2b_1r + la_2s_1b_2r + la_2\cdot r = 0,\\
&la_3s_1^{-1}b_1r + la_3b_2r + la_3s_1^{-1}\cdot r = 0.
\end{align*}
Notice that here we have only $5$ linear constraints (and not $6$), so the resulting linear space is (at most) $4$-dimensional. It is easy to see that in the group case the corresponding linear space is $1$-dimensional, and so this phenomenon degenerates in the group case.

\subsubsection{Structure of small cancellation algebras}
\label{intro_structure_section}

Like in the group case, starting with small pieces, we introduce a $\SPM$-measure on the monomials from~$\Mon$. If a monomial $v$ can be written as a product of small pieces, then $\SPM(v)$ is the minimal number of small pieces needed for such presentation, otherwise $\SPM(v)=\infty$.

We fix a threshold constant $\tau$ and say that  $\Qalg$ is \emph{$C(\tau)$-small cancellation ring} if it satisfies Compatibility Axiom, Small Cancellation Axiom with respect to $\tau + 1$ small pieces and Isolation Axiom. In the further argument we require $\tau \geqslant 10$.

In our paper, we perform a very detailed study of the influence of multi-turns on members of the charts of monomials. Here we have to treat several caveats. When we define members of a chart in the terms of their $SPM$-measure, such definition is not stable enough under multi-turns. So, we have to define a very delicate notion of \emph{a virtual member of the chart} (Section~\ref{virtual_members_section}). In order to do this, one has to verify that the cumulative influence of a long sequence of multi-turns has no adverse effects. In the course of this verification, we introduce several invariants of the chart, and finally, we introduce the notion of $f$-characteristic (Definition~\ref{f_char_def}). Using the $f$-characteristic, we introduce a global increasing filtration $\Ft_n (\Galg)$, $n \geqslant 0$, on $\Galg$, which gives rise to the filtration on the quotient algebra $\Qalg$, and then consider the corresponding graded objects. All the members of our filtration are spanned by monomials. Hence,  we can define the space of dependencies of a given member of the filtration as a linear subspace spanned by the spaces of dependencies of all its monomials. It is important to stress that when multi-turns cause degeneracies, like it happened in our Example C, the corresponding monomials descend to the lower members of the filtration, so that they disappear in the graded objects. It turns out, that it is sufficient to study dependencies in the following particular local situation. Namely, our Main Lemma~\ref{fall_through_linear_dep_one} describes the interaction between the filtration and the spaces of linear dependencies:
\begin{equation*}
\Dp(\DSpace{U}) \cap \Low(\DSpace{U}) \subseteq \Ft_{n - 1}(\Galg),
\end{equation*}
where $U$ is a monomial from $\Ft_{n}(\Galg) \setminus \Ft_{n - 1}(\Galg)$, $\DSpace{U}$ stands for the current level of the filtration, $\Low(\DSpace{U})$ denotes the next level of the filtration (in the descending order), and $\Dp$ denotes the spaces of dependencies (see Section~\ref{structure_calc} for precise definitions).

Here is the place to make  some comments. In essence, the main result above claims that in the quotient algebra $\Qalg$ \emph{there are no unexpected linear dependencies}. But, first, one has to explain what are \emph{the expected linear dependencies}.

Consider the filtration $\Ft_n(\Galg)$, $n \geqslant 0$, on $\Galg$ introduced above. Let $U \in \Ft_n (\Galg)$ be a monomial such that its chart has $m$ virtual members, $U = L^{(i)} u^{(i)} R^{(i)}$, $i=1, 2, \ldots , m$. For any $p \in \Rel$ of the form $p = \alpha u^{(i)} + \sum_{j = 1}^k \alpha_j a_j$, $\alpha \neq 0$, we consider the polynomial $L^{(i)}\cdot p\cdot R^{(i)} \in \Galg$. All such polynomials obviously belong to $\Ft_n (\Galg) \cap \Ideal$ and are regarded as expected dependencies.

In fact, our main result claims that the opposite is also true.

\begin{theorem}
$\Ft_n (\Galg) \cap \Ideal$
is linearly spanned by all the polynomials of the form $L^{(i)} \cdot p\cdot R^{(i)}$, $i = 1,\ldots , m$ for all monomials $U \in \Ft_n (\Galg)$ and polynomials $p \in \Rel$ as above, $n \geqslant 0$.
\end{theorem}

Given the theorem, the graded object $\Gr(\Qalg)$ can be explicitly determined by local calculations (see, in particular, Proposition~\ref{component_subspaces_structure}). This allows us to show that the quotient algebra $\Qalg$ is non-trivial, construct a linear basis of $\Qalg$, and prove that the equality problem in $\Qalg$ is solvable by a counterpart of Dehn's algorithm.

\subsubsection{{\bf {\it A posteriori}} insight. Gr\"{o}bner basis }
\label{intro_Grobner_section}

With the benefit of hindsight, we observe that our theory bears an intimate relation to the theory of Gr\"{o}bner basis and Diamond Lemma (see, e.g.,~\cite{BKW}, \cite{BCKK}, \cite{U}, \cite{Be} and references therein). In particular, the polynomial relations~$\Rel$ satisfying our Small Cancellation axioms form a Gr\"{o}bner basis of the ideal $\Ideal$. Moreover, Main Lemma~\ref{fall_through_linear_dep_one} can be considered as a sophisticated analog of Diamond Lemma.

One should emphasize that  Gr\"{o}bner basis philosophy was not a ruling engine of the calculations. Nevertheless, it seems  that it sheds a lot of light on what is going on behind the technicalities of proofs. We recall here some known facts in order to make the exposition self-contained.

Let us start with the  classical setting. Let $K\langle a_1,\ldots,a_n\rangle$ be the free associative algebra with generators $a_1,\ldots,a_n$. Elements of the algebra $K\langle a_1,\ldots,a_n \rangle$ are polynomials with non-commutative monomials. The elements of the corresponding monoid are called monomials. We fix a linear ordering of the monomials. Usually the monomials are ordered by length, and then lexicographically. This ordering is called \emph{$DegLex$}. Take the ideal $\Ideal$ generated as an ideal by elements $\lbrace f_i\rbrace$. Let  $\overline{g}$ denote the highest monomial of a polynomial $g$ in this order.

We consider the following  natural \emph{greedy algorithm}. Its steps are as follows: given a polynomial $g$ in $K \langle a_1,\ldots,a_n \rangle$, we take its highest monomial $\overline{g}$. If $\overline{g}$ contains a submonomial of the form $\overline{f}_i$ for some $i$, then $\overline{g}$ has the form $\overline{g} = l\overline{f}_i r$. Then subtract from $g$ the product $l\cdot f_i\cdot r$ with an appropriate coefficient, thus cancelling the highest term of $g$. If we come to a polynomial such that its highest monomial cannot be cancelled using this procedure, then the algorithm terminates; otherwise we continue reduction. By definition of the greedy algorithm, termination at zero means that the element $g$ belongs to ideal $\Ideal$.

The family $\lbrace f_i\rbrace$ is called \emph{a Gr\"{o}bner Basis of the ideal $\Ideal$} if the following property holds: if the greedy algorithm  which starts at some $g$ in $K \langle a_1,\ldots,a_n \rangle $ terminates at some non-zero element $g_0$, then $g$ does not belong to ideal $\Ideal$.

Not every system of generators of an ideal is a Gr\"{o}bner basis. Moreover, in general the ideal membership problem is not algorithmically solvable and, in particular, not every ideal has a finite Gr\"{o}bner basis. However, there is a criterion for a set of generators of the ideal to be a Gr\"{o}bner basis. This criterion is provided by Bergman's Diamond Lemma~\cite{Be}. It works as follows. Suppose we are given a finite set of polynomials $\lbrace f_i\rbrace$ of the ideal $\Ideal$ and a greedy reduction algorithm. We want to check if $\lbrace f_i\rbrace$ is a Gr\"{o}bner basis. Fix two monomials $f_j$ and $f_k$ from our set. Look at all monomials  $M$ such that $\overline{f}_j$ is a prefix of $M$ and $\overline{f}_k$ is a suffix of $M$, that is $M = h_l\overline{f}_k = \overline{f}_j h_r$. Then $M$ can be reduced to a polynomial with smaller highest monomials in two ways, namely using the reduction via the submonomial $\overline{f}_j$ or via the submonomial $\overline{f}_k$. Take the difference $f_{jk}$ of the results. Clearly, $f_{jk} \in \Ideal$. The set of polynomials $\lbrace f_i\rbrace$ is a Gr\"{o}bner basis if and only if every difference can be reduced to zero by the greedy algorithm.

If the overlap of $\overline{f}_j$ and  $\overline{f}_k$ is empty, then $f_{jk}$ is reduced to zero by the greedy algorithm in a trivial way. Suppose that $\overline{f}_j$ and  $\overline{f}_k$ have a non-empty overlap. Clearly, there is a finite number of such $M$. So, in order to show that a finite family $\lbrace f_i\rbrace$ is a Gr\"{o}bner basis we need to check only the possibility of reduction to $0$ of a finite number of polynomials. This criterion is called  Bergman's Diamond Lemma.

One can look at the condition in Diamond Lemma from the following point of view. Assume we have two polynomials $T_j = l_j\cdot f_j \cdot r_j$, $T_k = l_k\cdot f_k\cdot r_k$, where $l_j, r_j, l_k, r_k$ are monomials, such that the highest monomials of $T_j$ and $T_k$ are equal and have equal coefficients. Clearly, this highest monomial cancels in $T_{jk} = T_j - T_k$. The condition of Diamond Lemma (aka Diamond Condition) states that every $T_{jk}$ can be reduced to~$0$. This condition can be replaced by a stronger one. Namely, instead of two polynomials we take a number of polynomials $T_1, \ldots, T_n$, $T_i = l_i\cdot f_i\cdot r_i$. Let $m$ be the biggest monomial among all $\overline{T}_i$. We consider linear combinations $\sum\gamma_iT_i$ such that the monomial $m$ cancels out. The new condition states that $\sum\gamma_iT_i$ can be reduced to~$0$. Then it is easy to see that $\lbrace f_k\rbrace$ is a Gr\"{o}bner basis if and only if the new condition holds. In particular, this means that the new stronger condition is, in fact, equivalent to the initial condition of Diamond Lemma.

The above modification of Diamond Condition allows us to consider our Main Lemma as a verification of this modified condition for a certain set of generators of the space~$\Dp\DSpace{U}$. Indeed, consider a linear combination of layouts of multi-turns $\sum\gamma_iT_i$ that belongs to $\Ft_{n - 1}(\Galg)$, where $U \in \Ft_n(\Galg)\setminus \Ft_{n - 1}(\Galg)$. This means that all monomials of $T_i$ that do not belong to $\Ft_{n - 1}(\Galg)$ cancel out in this sum. This is exactly a particular form of cancelling of the highest monomial. Main Lemma~\ref{fall_through_linear_dep_one} states that such $\sum\gamma_iT_i$ can be represented as a linear combination of layouts of multi-turns of the monomials from $\Ft_{n - 1}(\Galg)$. This is an analogue of the possibility of reduction.

Let us notice that inside the proof of our Main Lemma~\ref{fall_through_linear_dep_one} we use a machinery which is very similar to the one from Diamond Lemma in the classical sense. In particular, using the small cancellation conditions, we introduce a linear ordering on monomials of $\DSpace{U}$ that do not belong to $\Ft_{n - 1}(\Galg)$. Moreover, using the same ideas that we use in the argument of Main Lemma~\ref{fall_through_linear_dep_one}, we can introduce a linear ordering for all monomials and prove a classic version of Diamond Lemma for layouts of multi-turns and this ordering. Thereby, one can show that $\Rel$ is a Gr\"{o}bner basis of the ideal $\Ideal$ with respect to this special linear ordering.

\subsection{Route map of the paper}
\label{intro_Route_section}

For the general idea of the paper we refer to the preceding  Subsection~\ref{intro_overview_section} ``Overview of the work''. We study the quotient ring $\Qalg$, where $\Galg$ is the group algebra over the free finitely generated group $\Fr$, $k$ is a field, and  $\Ideal$ is an ideal of $\Galg$ generated as an ideal by set of polynomials $\Rel$ that satisfies certain conditions (Axioms~\ref{comp_ax}, \ref{sc_ax}, 3a, 3b in Section~\ref{ring_sc_ax_section}).

Since the paper is long and technical we aim to provide the reader of ``Route map'' with some familiarity with the essential notions and results.

Numbering of statements is as follows. Theorems have through numbers throughout the paper. Lemmas, propositions, corollaries, definitions, examples and remarks are enumerated within every section. Lemmas, propositions and corollaries have through numbers within every section. Definitions, examples and remarks are enumerated independently.

In {\bf Section~\ref{ring_sc_ax_section}} we fix our setting. We develop a ring-theoretic counterpart of the group-theoretic Small Cancellation theory, so a part of our axioms is modeled after the corresponding group-theoretic axioms. Condition~\ref{comp_ax}, the  Compatibility Axiom, corresponds to group-theoretic condition that a set of relations is closed under cyclic shifts. Condition~\ref{sc_ax} called Small Cancellation Axiom is a counterpart of the group-theoretic Small Cancellation Axiom. To state it we define for rings the notion of \emph{a small piece} (Definition~\ref{sp}). The connection to the corresponding group-theoretic notion of a small piece is not obvious, so we recommend to look at the example from Subsection~\ref{group_algebras_of_sc_groups_section} and Remark~\ref{small_pieces_comment} from Section~\ref{ring_sc_ax_section}, that may clarify the relations between these notions.

Elements of $\Fr$ and $\Galg$ are called \emph{monomials (or words)} and \emph{polynomials}, respectively. We denote by $\Mon$ the set of all monomials that are non-zero summands of polynomials from the set of relations $\Rel$. Polynomials from $\Rel$ generate $\Ideal$ as an ideal and, as it was said, are subject to three axioms mentioned above.

We introduce  \emph{a measure} on the monomials from $\Mon$ defined as the minimal number of small pieces needed to represent the monomial as their product. We call it $\SPM$-measure.

Next we define \emph{incident monomials}.  Monomials $a_i$ and $a_j$ are incident if they appear with non-zero coefficients in some polynomial $p\in \Rel$.

{\bf Section~\ref{basic_def}} deals with other basic definitions.  We start with the notion of \emph{occurrence of a submonomial in a given monomial}. We  are interested  in occurrences of monomials from the set $\Mon$ and define maximal occurrences of such monomials. In what follows we call them for short \emph{maximal occurrences}. Then we consider \emph{overlaps of occurrences} and show the important property that maximal occurrences from $\Mon$ can overlap only via a small piece.  Throughout the paper we graphically represent monomials as segments and their subwords as subsegments.

Definition~\ref{chart_def} is devoted to the notion of \emph{the chart of a monomial} and \emph{members of the chart}. We postulate a threshold value $\tau$ of the $\SPM$-measure of the monomials. The \emph{chart of a monomial $U$} is defined to be the set of all maximal occurrences of monomials of $\Mon$ in $U$.  The maximal occurrences $m \in \Mon$ in $U$ such that $\SPM(m) \geqslant \tau$ are called \emph{members of the chart}.

It follows from Small Cancellation Axiom that the ideal $\Ideal$ is linearly generated by the polynomials such that some of their monomials have maximal occurrences with $\SPM$-measure  $\geqslant\tau$ (see Proposition~\ref{the_ideal_characterisation}).

Definition~\ref{multiturn_def} introduces an operation on monomials and polynomials which we call  \emph{multi-turn}. This operation is identical modulo ideal $\Ideal$. On the other hand, Compatibility Axiom guarantees that any two polynomials equal modulo ideal $\Ideal$ can be transformed one to the other by a finite sequence of multi-turns.

Let $p = \sum\limits_{j = 1}^{n} \alpha_{j} a_{j}$ be a polynomial from $\Rel$. For $h$, $1\leqslant h \leqslant n$, the transformation
\begin{equation*}
a_h \longmapsto \sum_{\substack{j = 1\\ j\neq h}}^{n}(-\alpha_h^{-1}\alpha_j) a_j
\end{equation*}
is called \emph{an elementary multi-turn}. Let $U = La_hR$ be a monomial with the maximal occurrence $a_h$. The transformation\begin{equation*}
U = La_hR \longmapsto \sum\limits_{\substack{j = 0\\ j\neq h}}^{n} (-\alpha_{h}^{-1} \alpha_{j})La_{j}R,
\end{equation*}
with the further cancellations if there are any, is called \emph{a multi-turn} of the occurrence $a_h$ in $U$. Given a polynomial we apply a multi-turn to one of its components. The transformations of individual monomials $U_h = La_hR \mapsto U_j = La_jR$ are called \emph{replacements}. In the end of this section we introduce the notion of \emph{a layout}. Given a monomial $U = La_hR$ and a multi-turn as above, we call the polynomial $L\cdot p\cdot R$ the layout of this multi-turn. Obviously, a layout of a multi-turn always belongs to the ideal $\Ideal$.

We have set the stage for the central  topic of this paper: the interplay between the charts of monomials and sequences of multi-turns.

{\bf Section~\ref{ideal_linear_descr_section}} is devoted to description of the ideal $\Ideal$ as a linear subspace in $\Galg$. For every monomial of $\Fr$ we do all multi-turns of all members of the chart and take the set $\GDp$ of all layouts of these multi-turns. Denote by $\langle \GDp\rangle$ the linear span of $\GDp$. The main result of this chapter states that our ideal $\Ideal$ coincides with $\langle \GDp\rangle$. From the point of view of the ring $\Qalg$ the elements of $\GDp$ can be considered as linear dependencies and $\langle \GDp\rangle$ as \emph{the space of linear dependencies}.

In  {\bf Section~\ref{mt_configurations}} we figure out how a multi-turn influences on the chart of a monomial. So, consider replacement $U_h = La_hR \mapsto U_j = La_jR$ under condition $\SPM(a_h)\geqslant 3$, where $a_h$ and $a_j$ are incident monomials. To every maximal occurrence in $U_h$ there corresponds the set of \emph{images} in $U_j$ under the given replacement. In most cases of interest there is precisely a unique image or no images. In this section we describe how the images look like.

We consider three variants for the resulting monomial $U_j = La_jR$: $a_j$ is not a small piece; $a_j$ is  a small piece; $a_j$ is $1$. We show that in the first case the structure of the chart remains almost stable after a multi-turn. In particular, only maximal occurrences that are \emph{not separated} from $a_h$ may change at most by one small piece. In the second case the replacement $a_h$ by $a_j$ can cause merging and  restructuring of the chart, and in the third case massive cancellations resulting in complete modification of the chart are possible.

We produce the full list of all possible changes of maximal occurrences. The calculations are based on thorough analysis of all combinatorial possibilities.

Let us compare our general situation with the case of the group algebra, where the group relations may be expressed as binomials (see Subsection~\ref{group_algebras_of_sc_groups_section}). Hence, we have a unique resulting monomial $U_j$ of every multi-turn, where the merging either happens or not. In the ring case a multi-turn is associated with a polynomial relation and the structure of an emerging chart depends on the particular resulting monomial $U_j$. Thus,  one and the same multi-turn can produce significantly different charts of resulting monomials.

In {\bf Section~\ref{virtual_members_section}} we apply the results of the previous sections to introduce \emph{numerical invariants of the  charts}. We define the notions of \emph{a minimal covering}, of \emph{an admissible replacement} and of \emph{virtual members of the chart of a monomial}. The number of elements of a minimal covering and the number of virtual members of a chart turn out to be useful numerical invariants of a monomial that interact in a satisfactory way with replacements $U_h = L a_h R \mapsto U_j = L a_j R$.

Recall that a member of a chart is a maximal occurrence with $\SPM$-measure $\geqslant \tau$. Example~\ref{example_members_not_stable} shows that the total number of members of the chart does not behave in a satisfactory way with replacements $U_h = L a_h R \mapsto U_j = L a_j R$. So, we have to look on more subtle invariants. One such invariant is the number of elements in \emph{a minimal covering of $U_h$} (see Subsection~\ref{coverings_subsection}). Let $\mo{U_h}$ be the set of all maximal occurrences in $U_h$. Consider subsets of $\mo{U_h}$ that cover the same letters in $U_h$ as the whole $\mo{U_h}$. A covering of such type consisting of the smallest number of elements is called \emph{a minimal covering}. Of course, such covering is not, necessarily, unique. We denote the size of the minimal covering by $\mincov{U_h}$. The precise behavior of the minimal covering number is stated in Lemma~\ref{minimal_coverings_property}, which is the main result of Subsection~\ref{coverings_subsection}.

Now we replace the notion of a member of the chart by a quite delicate notion of \emph{a virtual member} of the chart. First, we introduce the notion of \emph{an admissible replacement}. We say the replacement $U_h = La_h R \mapsto U_j = L a_j R$ is admissible if $\SPM(a_h) \geqslant \tau - 2$ and $a_j \neq 1$ and $a_j$ is not fully covered by the images of the elements of $\mo{U_h} \setminus \lbrace a_h \rbrace$ with $\SPM$-measure $\geqslant 3$. So, this is a rather technical notion. Roughly speaking, virtual members of the chart are those maximal occurrences which originally are not necessarily members of the chart, but under a series of admissible replacements their images are members of the chart (see Definition~\ref{virtual_members_def}). Denote the number of virtual members of $U$ by $\nvirt{U}$.

It is worth mentioning that while member of a chart is a \emph{local} notion, a virtual member of a chart is a \emph{global} notion, that is, by changing one single letter at the end of a monomial some maximal occurrences inside the monomial may become (or seize to be) a virtual member of the chart.

Corollary~\ref{virtual_members_number} claims that the behavior of virtual members of charts is good enough, namely $\nvirt{U_j} \leqslant \nvirt{U_h}$ for an admissible replacement $U_h = L a_h R \mapsto U_j = L a_j R$. If $a_j$ is not a virtual member of the chart of $U_j$, then $\nvirt{U_j} < \nvirt{U_h}$. Moreover, given a replacement $a_h\mapsto a_j$, where $a_h$ is a virtual member of the chart of $U_h$, $a_j$ is a virtual member of the chart of $U_j$, taking images gives a bijective correspondence between all virtual members of the chart of $U_h$ and all  virtual members of the chart of $U_j$. In particular, $\nvirt{U_h} = \nvirt{U_j}$, see Corollary~\ref{full_virtual_members_corresponding}. That is, the number of virtual members of the chart can not increase after an admissible replacement, unlike the number of members of the chart.

Now, exactly in the same way as in Section~\ref{ideal_linear_descr_section}, we define the linear subspace $\langle \GDp^{\prime} \rangle$ of $\Galg$ of linear dependencies induced by multi-turns of virtual members of the chart of monomials. Proposition~\ref{the_ideal_characterisation2} states that $\langle\GDp\rangle = \langle\GDp^{\prime}\rangle = \Ideal$, i.e., all layouts of multi-turns of virtual members generate linearly the same ideal as all layouts of multi-turns of members of the chart (cf. Proposition~\ref{the_ideal_characterisation}).

Next Proposition~\ref{non_increasing_parameter} aggregates all properties of the numerical characteristics $\mincov{U_h}$ and $\nvirt{U_h}$  obtained before. Assume $U_h$ is a monomial, $a_h$ is a virtual member of the chart of $U_h$. Let $a_h$ and $a_j$ be incident monomials. Consider the replacement $U_h = La_hR \mapsto U_j = La_jR$ in $U_h$. If $a_j$ is a virtual member of the chart of $U_j$, then $\mincov{U_h} = \mincov{U_j}$ and $\nvirt{U_h} = \nvirt{U_j}$. If $a_j$ is not a virtual member of the chart of $U_j$, then either $\mincov{U_j} < \mincov{U_h}$, or $\mincov{U_j} = \mincov{U_h}$ and $\nvirt{U_j} < \nvirt{U_h}$.

In fact, Section~\ref{mt_configurations} and Section~\ref{virtual_members_section} lay ground for constructing a special partial ordering on monomials. This order relation is introduced in Subsection~\ref{derived_monom_def_section} via the invariant called $f$-characteristic, but, in fact, all necessary results are formulated in Subsection~\ref{virtual_members_subsection},  Proposition~\ref{non_increasing_parameter}. In what follows it serves as a background of a special total order, which is introduced in Section~\ref{algorithm_section}.
 Once  Proposition~\ref{non_increasing_parameter} is proven, we will start the second part of the paper devoted to combinatorial structure of our ring.

In {\bf Section~\ref{derived_monom_section}} we look for further details of transformation of monomials. It consists of three parts described below.

{\bf Subsection~\ref{derived_monom_def_section}} is devoted solely to definition of two notions. First, we define \emph{derived monomials} of $U$ as a result of application of a sequence of replacements of virtual members of the chart by incident monomials, starting from a given monomial $U$. Second, in Definition~\ref{f_char_def} we introduce the important \emph{$f$-characteristic} of a monomial. Given a monomial $U$ we define $f(U)$ to be the pair $(\mincov{U}, \nvirt{U})$. Then we have a partial order on monomials defined by $f(U_1) < f(U_2)$ if and only if either $\mincov{U_1} < \mincov{U_2}$, or $\mincov{U_1} = \mincov{U_2}$ and $\nvirt{U_1} < \nvirt{U_2}$. The next Lemma~\ref{estimation_value_property} explains that $f$-characteristic works properly with respect to transformation of the monomials, i.e., the $f$-characteristic of a derived monomial is less than or equal to $f$-characteristic of the original one. Moreover, we know  effectively when we have inequality.

{\bf Subsection~\ref{replacements_by_incident_section}} and {\bf Subsection~\ref{replacements_by_U_incident_section}} are devoted to further thorough study of the replacements of virtual members of the chart by  incident monomials. This study is a quite long series of technical statements used in the next  section in order to describe linear spaces generated by derived monomials in terms of tensor products. In particular, we introduce the notion of $U$-incident monomials widely used in Section~\ref{structure_calc}.

In {\bf Section~\ref{structure_calc}} we define a filtration and a corresponding grading on $\Galg$ that agrees with layouts of multi-turns. We describe the corresponding homogeneous components in terms of tensor products of some smaller spaces. This is the place where we start constructing the combinatorial structure of~$\Qalg$.

In particular, in {\bf Subsection~\ref{filtration_def_subsection}} we define an increasing filtration on $\Galg$, using $f$-characteristic of monomials introduced in Subsection~\ref{derived_monom_def_section}. Namely, the $f$-characteristic gives rise to a certain function $t$ on natural numbers defined as follows. We put $t(0) = (0, 0)$. Assume $t(n) = (r, s)$, then we put
\begin{equation*}
t(n + 1) = \begin{cases}
(r, s + 1) &\textit{ if } r > s,\\
(r + 1, 0) &\textit{ if } r = s.
\end{cases}
\end{equation*}
We define an increasing filtration on $\Galg$ by the rule:
\begin{equation*}
\Ft_n(\Galg) = \left\langle Z \mid Z \in \Fr, f(Z) \leqslant t(n)\right\rangle.
\end{equation*}
That is, the space $\Ft_n(\Galg)$ is generated by all monomials with $f$-characteristics not greater than $t(n)$. We show that each
$\Ft_n(\Galg)$ is closed under taking derived monomials.

We need a set of new notions. Let $U$ be a monomial. By $\DSpace{U}$ we denote a linear subspace of $\Ft_n(\Galg)$ generated by all derived monomials of $U$. By $\Low\DSpace{U}$ we denote the subspace generated by all derived monomials of $U$ with $f$-characteristic smaller than~$f(U)$. The next principal object is \emph{the space of dependencies}, defined in Definition~\ref{subspace_of_dependencies} as follows. Suppose $Y$ is a subspace of $\Galg$ linearly generated by a set of monomials and closed under taking derived monomials. We take the set of all the layouts of multi-turns of virtual members of the chart of monomials of $Y$ and look at its linear envelope $\Dp(Y)$, which is our set of dependencies related to $Y$. Note that in this terms Proposition~\ref{the_ideal_characterisation2} claims that $\Dp(\Galg) = \Ideal$. The  key statement of Subsection~\ref{filtration_def_subsection} is Proposition~\ref{fall_through_linear_dep}. It describes the nice interaction between dependencies and filtration:
\begin{equation*}
\Dp(\Ft_n(\Galg)) \cap \Ft_{n - 1}(\Galg) = \Dp(\Ft_{n - 1}(\Galg)).
\end{equation*}
Its proof is based on Main Lemma~\ref{fall_through_linear_dep_one} which deals with a linear space generated by a single monomial and its derived monomials. Namely, let $U$ be an arbitrary monomial, $U \in \Ft_n(\Galg) \setminus \Ft_{n - 1}(\Galg)$. Then Main Lemma~\ref{fall_through_linear_dep_one} claims that
\begin{equation*}
\Dp\DSpace{U} \cap \Low\DSpace{U} \subseteq \Dp(\Ft_{n - 1}(\Galg)).
\end{equation*}
The proof of Main Lemma is postponed till the next Subsection~\ref{tensor_product_section}, since we need first to introduce some calculations related to tensor products.

In the end of this subsection we obtain the following important extension of Proposition~\ref{fall_through_linear_dep_one}. Suppose $X, Y$ are subspaces of $\Galg$ generated by monomials and closed under taking derived monomials, $Y \subseteq X$. Then $\Dp(X)\cap Y = \Dp(Y)$ (see Proposition~\ref{fall_to_smaller_subspace}).

In {\bf Subsection~\ref{tensor_product_section}} we study tensor products of special linear spaces which are tightly related to the charts of monomials. The idea behind that is as follows. Assume $U$ is a monomial. Derived monomials of $U$ are defined with the use of certain sequences of replacements of virtual members of the chart (see Definition~\ref{derived_monomials}). When we perform replacements that preserve $f$-characteristics of monomials, they preserve, roughly speaking, the structure of the chart. Moreover, there is no interaction between the replaced occurrence and the separated virtual members of the chart and there is a very small interaction between the replaced occurrence and its neighbours. This kind of behaviour gave rise the idea of considering a tensor product of linear spaces that correspond to each place of the chart of $U$.

Assume a monomial $U$ has $m$ virtual members of the chart, that is, $\nvirt{U} = m$. We enumerate all the virtual members of the chart of $U$ from left to right. Let $u^{(i)}$ be the $i$-th virtual member of the chart of $U$. We define a linear space $A_i[U]$ by the following formula
\begin{equation*}
A_i[U] = \left\langle a^{(i)} \mid u^{(i)} \textit{ and } a^{(i)} \textit{ are } U\textit{-incident monomials}\right\rangle.
\end{equation*}
Then, given $U \in \Ft_{n}(\Galg) \setminus \Ft_{n - 1}(\Galg)$, we construct in Definition~\ref{mu_def} a linear mapping
\begin{equation*}
\mu[U] : A_1[U] \otimes \ldots \otimes A_m[U] \to \DSpace{U} +  \Ft_{n - 1}(\Galg).
\end{equation*}
Lemma~\ref{mu_properties} states important properties of the mapping $\mu[U]$. Here is the unique place where we use Isolation Axiom. Now we are in a position to prove Main Lemma~\ref{fall_through_linear_dep_one}. The proof of Main Lemma~\ref{fall_through_linear_dep_one} completes Subsection~\ref{tensor_product_section}.

The next {\bf Subsection~\ref{grading_subsection}} is devoted to grading of the space~$\Qalg$. First we define the corresponding filtration on $ \Qalg$ in the following way:
\begin{equation*}
\Ft_n (\Qalg) = (\Ft_n(\Galg) + \Dp(\Galg)) / \Dp(\Galg) = (\Ft_n(\Galg) + \Ideal) / \Ideal.
\end{equation*}
We define a grading on $\Qalg$ by the rule:
\begin{equation*}
\Gr(\Qalg) = \bigoplus\limits_{n = 0}^{\infty}\Gr_n(\Qalg) = \bigoplus\limits_{n = 0}^{\infty}\Ft_n (\Qalg) / \Ft_{n - 1} (\Qalg).
\end{equation*}
Theorem~\ref{structure_of_quotient_space} establishes the compatibility of the filtration and the corresponding grading on~$\Qalg$ with the space of dependencies~$\Dp(\Galg)$. Namely, it states that
\begin{equation*}
\Gr_n(\Qalg) \cong \Ft_n(\Galg) / (\Dp(\Ft_n(\Galg)) + \Ft_{n - 1}(\Galg)).
\end{equation*}
Proposition~\ref{component_subspaces_structure} in the end of the section provides a kind of semi-simplicity property for $\Gr_n$.

{\bf Section~\ref{basis_section}} contains the main results that summarize all previous work. In {\bf Subsection~\ref{quotient_space_non_triviality_subsection}} we show that the ring $\Qalg$ is non-trivial. First, we notice that the space $\DSpace{X} / (\Dp\DSpace{X} + \Low\DSpace{X})$, where $X$ is a monomial with no virtual members of the chart, is always non-trivial, and of dimension~$1$ (Lemma~\ref{non_trivial_spaces_existence}). By definition the empty monomial~$1$ is always a small piece and thus it has no virtual members of the chart. Then Corollary~\ref{non_trivial_quotient}, stating that $\Qalg$ is non-trivial, is just a simple combination of Proposition~\ref{fall_to_smaller_subspace} and Lemma~\ref{non_trivial_spaces_existence}.

In {\bf Subsection~\ref{construction_of_basis_subsection}} we are able, at last, to describe a basis of $\Qalg$ as a vector space. This is done in two steps. First, we construct in Proposition~\ref{graded_component_basis} a basis for non-trivial graded components of our filtration on $\Qalg$:
\begin{equation*}
\Gr_n(\Qalg) = \Ft_n(\Qalg) / \Ft_{n - 1}(\Qalg).
\end{equation*}
Given $n$ we consider the set of spaces $\lbrace\DSpace{Z} \mid Z \in \Fr, Z \in \Ft_n(\Galg)  \setminus \Ft_{n - 1}(\Galg)\rbrace$, such that $\DSpace{Z} / (\Dp\DSpace{Z} + \Low\DSpace{Z}) \neq 0$. Let $\lbrace V_{i}^{(n)}\rbrace_{i \in I^{(n)}}$ be all different spaces from this set. Then, the semi-simplicity property from Subsection~\ref{grading_subsection} implies that
\begin{equation*}
\Gr_n(\Qalg) \cong \bigoplus\limits_{i \in I^{(n)}} V_{i}^{(n)} /(\Dp(V_{i}^{(n)}) + \Low(V_{i}^{(n)})).
\end{equation*}
Assume $\lbrace \overline{W}^{(i, n)}_j\rbrace_j$ is a basis of $V_{i}^{(n)}/(\Dp(V_{i}^{(n)}) + \Low(V_{i}^{(n)}))$, $i \in I^{(n)}$. Let $W^{(i, n)}_j \in V_{i}^{(n)}$ be an arbitrary representative of the coset $\overline{W}^{(i, n)}_j$. Then
\begin{equation*}
\bigcup\limits_{i \in I^{(n)}}\left\lbrace W^{(i, n)}_j + \Ideal + \Ft_{n - 1}(\Qalg)\right\rbrace_j
\end{equation*}
is a basis of $\Gr_n(\Qalg)$. Finally, Theorem~\ref{whole_quotient_ring_structure} describes a basis of $\Qalg$.  Let $\lbrace V_i\rbrace_{i \in I}$ be all different spaces $\lbrace\DSpace{Z} \mid Z \in \Fr\rbrace$.  We have
\begin{equation*}
\Qalg \cong \bigoplus\limits_{i \in I} V_i/(\Dp(V_i) + \Low(V_i)),
\end{equation*}
as vector spaces, and the right-hand side is explicitly described in Proposition~\ref{correspondence_to_tensor_product}.
Assume $\lbrace \overline{W}^{(i)}_j \rbrace_j$ is a basis of $V_i/(\Dp(V_i) + \Low(V_i))$, $i \in I$. Let $W^{(i)}_j \in V_i$ be an arbitrary representative of the coset $\overline{W}^{(i)}_j$. Then
\begin{equation*}
\bigcup\limits_{i \in I}\left\lbrace W^{(i)}_j  + \Ideal\right\rbrace_j
\end{equation*}
is a basis of $\Qalg$.

In {\bf Section~\ref{algorithm_section}} we study algorithmic properties of our ring. We show that they are similar in a sense to the ones valid for small cancellation groups. However, in the ring case  the essential differences arise everywhere.

Recall that small cancellation groups enjoy Dehn's algorithm~\cite{LS}. In this section we define and study a greedy algorithm for rings which plays a similar role as Dehn's algorithm does for groups.

Let a ring $\Qalg$ with small cancellation conditions be given. First of all we need to extend a bit our set of relations $\Rel$ to a certain  additive closure $\Add(\Rel)$. It is important that for the natural examples of {\bf Section~\ref{examples_section}} we have $\Rel = \Add(\Rel)$. Then we define a linear order on all monomials, based on $f$-characteristic and the description in terms of tensor products (see Subsections~\ref{derived_monom_def_section} and~\ref{tensor_product_section}), and denote it by $<_f$ (see Definition~\ref{f_char_order_def}). Then given the  order $<_f$ and the set $\Add(\Rel)$, we define a non-deterministic greedy algorithm with external source of knowledge and denote it by $\GreedyAlg(<_f, \Add(\Rel))$ (see Definition~\ref{greedy_alg_rings}).

Recall that given a small cancellation group $G = \langle \mathcal{X} \mid \Rel_G\rangle$, a word $W$ from a free group is equal to $1$ in $G$ if and only if Dehn's algorithm, starting from $W$, terminates at $1$,~\cite{LS}. Our Theorem~\ref{possibility_of_reduction} establishes the similar properties of $\GreedyAlg(<_f, \Add(\Rel))$ in much more complicated situation of rings.

Namely, assume $W_1, \ldots, W_k$ are different monomials. We take an element $\sum_{i = 1}^k \gamma_i W_i \in \Galg$, $\gamma_i \neq 0$. Then the following statements are equivalent:
\begin{itemize}
\item
some branch of the algorithm $\GreedyAlg(<_f, \Add(\Rel))$, starting from $\sum_{i = 1}^k \gamma_i W_i$, terminates at~$0$;
\item
$\sum_{i = 1}^k \gamma_i W_i \in \Ideal$;
\item
every branch of the algorithm $\GreedyAlg(<_f, \Add(\Rel))$, starting from $\sum_{i = 1}^k \gamma_i W_i$, terminates at~$0$.
\end{itemize}
Hence,
\begin{itemize}
\item
$\GreedyAlg(<_f, \Add(\Rel))$ solves the Ideal Membership Problem for~$\Ideal$,
\item
$\Add(\Rel)$ is a Gr\"{o}bner basis of the ideal $\Ideal$ with respect to monomial ordering $<_f$.
\end{itemize}

In {\bf Section~\ref{examples_section}} we give two examples of small cancellation rings. First, in {\bf Subsection~\ref{group_algebras_of_sc_groups_section}} we check that the group algebra of a small cancellation group that satisfies condition  $C(m)$ for $m\geqslant 22$ is a small cancellation ring. Here the small cancellation property appears twice and, correspondingly, has two faces. The first one is group theoretic while the second one is, of course, ring theoretic. In this case we show that there is no  need to distinguish between small pieces in a group sense and in a ring sense.

Another example is a ring constructed in~\cite{AKPR}, see {\bf Subsection ~\ref{equating_to_monomial_subsection}}. This is  a quotient ring $\Qalg$, where $\Galg$ is the group algebra of the free group $\Fr$ over the field $k$, and the ideal $\Ideal$ is generated by a single trinomial $v^{-1}-(1+w)$, where $v$ is a complicated word depending on $w$. The ring $\Qalg$ is of special interest because in this ring we have $(1 + w)^{-1} = v$. Thus, $1 + w$ becomes invertible. Here we show how to construct a set of defining relations equivalent to the relation $v^{-1} = 1 + w$ in $\Qalg$, that satisfies our three axioms.

\section{Group-like small cancellation axioms}
\label{ring_sc_ax_section}

Let $\Galg$\label{galg_def} be the group algebra of a free group $\Fr$ over some field $\fld$. Assume $\Fr$ has a fixed system of generators. Then its elements are reduced words in these generators and their inverses. We call these words monomials. Then the elements of the group algebra are linear combinations of monomials. We call them polynomials. Let $\Ideal$ be an ideal of $\Galg$ generated by a set of polynomials and let $\Qalg$ be the corresponding quotient algebra.

We state conditions on these polynomials that will enable a combinatorial description of the quotient algebra similar to small cancellation quotients of a free group. These axioms emerged when we studied the particular case described in~\cite{AKPR}.

Let us move on to formal definitions. Let the free group $\Fr$\label{fr_def} be freely generated by an alphabet $S$\label{alphabet_def}. Assume
\begin{equation*}
\label{rel_def}
\Rel = \left\lbrace p_i =  \sum\limits_{j = 1}^{n(i)} \alpha_{ij} m_{ij} \mid \alpha_{ij}\in \fld, m_{ij} \in \Fr, i \in I \right\rbrace
\end{equation*}
is a (finite or infinite) set of polynomials that generates the ideal $\Ideal$ (as an ideal). We denote this way of generating by $\langle \rangle_i$. So,
\begin{equation*}
\label{ideal_def}
\Ideal = \left\langle \Rel \right\rangle_i = \left\langle p_i = \sum\limits_{j = 1}^{n(i)} \alpha_{ij} m_{ij} \mid \alpha_{ij}\in \fld, m_{ij} \in \Fr, i \in I \right\rangle_i.
\end{equation*}
We assume that the monomials $m_{ij}$ are reduced, the polynomials $p_i$ are additively reduced, $I \subseteq \mathbb{N}$ is some index set. In particular, we assume that all the coefficients $\alpha_{ij}$ are non-zero. Denote the set of all monomials $m_{ij}$ of $\Rel$ by $\Mon$\label{mon_def}.

Throughout the paper we reserve small Greek letters for non-zero elements of the field~$\fld$.

\label{product_def}Throughout the paper we use the following notations. Assume $X$ and $Y$ are two elements of $\Galg$. We write their product as $X \cdot Y$. Assume $A$ and $B$ are two elements of $\Fr$. We write the product of the monomials $A$ and $B$ as $A\cdot B$. There may occur cancellations between $A$ and $B$ in $A\cdot B$. We write the product of $A$ and $B$ in the form $AB$ when there are no cancellations in~$A\cdot B$.

\begin{condition}[Compatibility Axiom]
\label{comp_ax}
The axiom consists of the following two conditions.
\begin{enumerate}
\item
If $p = \sum\limits_{j = 1}^{n} \alpha_{j} m_{j} \in \Rel$, then $\beta p = \sum\limits_{j = 1}^{n} \beta\alpha_{j} m_{j} \in \Rel$ for every $\beta \in \fld, \beta \neq 0$.

\item
Let $x \in S \cup S^{-1}$, $p = \sum\limits_{j = 1}^{n} \alpha_{j} m_{j} \in \Rel$. Suppose there exists $j_0 \in \lbrace 1, \ldots, n \rbrace$ such that $x^{-1}$ is the initial symbol of $m_{j_0}$. Then
\begin{equation*}
x\cdot p = \sum\limits_{j = 1}^{n} \alpha_{j} x\cdot m_{j} \in \Rel
\end{equation*}
(after the cancellations in the monomials $x\cdot m_{j}$).

We require the same condition from the right side as well. Suppose there exists $j_0 \in \lbrace 1, \ldots, n\rbrace$ such that $x^{-1}$ is the final symbol of $m_{j_0}$, then
\begin{equation*}
p\cdot x = \sum\limits_{j = 1}^{n} \alpha_{j} m_{j}\cdot x \in \Rel
\end{equation*}
(after the cancellations in the monomials $m_{j}\cdot x$).
\end{enumerate}
\end{condition}

Notice that taking any set of polynomials $\Rel_0$, $\langle\Rel_0\rangle_i = \Ideal$, one can construct a set of polynomials $\Rel \supseteq \Rel_0$ that satisfies Compatibility Axiom and generates the same ideal $\Ideal$.

From the second condition of Compatibility Axiom it immediately follows that the set $\Mon$ is closed under taking subwords. In particular, $1$ always belongs to $\Mon$.

Let $p = \sum_{j = 1}^{n} \alpha_{j} m_{j} \in \Rel$. Assume $c$ is a monomial and there exists $j_0 \in \lbrace 1, \ldots, n\rbrace$ such that $c^{-1}$ is a prefix of $m_{j_0}$, that is, $m_{j_0} = c^{-1}m_{j_0}^{\prime}$ ($c^{-1}$ is a suffix of $m_{j_0}$, that is, $m_{j_0} = m_{j_0}^{\prime}c^{-1}$). Then it easily follows from the second condition of Compatibility Axiom that $c\cdot p \in \Rel$ ($p\cdot c \in \Rel$).

Now we state a definition of \emph{a small piece}. It plays a central role in the further argument.

\begin{definition}
\label{sp}
Let $c \in \Mon$. Assume there exist two polynomials
\begin{align*}
&p = \sum\limits_{j = 1}^{n_1} \alpha_j a_j + \alpha a \in \Rel,\\
&q = \sum\limits_{j = 1}^{n_2} \beta_j b_j + \beta b \in \Rel,
\end{align*}
such that $c$ is a subword of $a$ and a subword of $b$. Namely,
\begin{align*}
&a = \widehat{a}_1c\widehat{a}_2,\\
&b = \widehat{b}_1c\widehat{b}_2,
\end{align*}
where $\widehat{a}_1$, $\widehat{a}_2$, $\widehat{b}_1$, $\widehat{b}_2$ are allowed to be empty. Assume that
\begin{equation*}
\widehat{b}_1\cdot \widehat{a}_1^{-1}\cdot p = \widehat{b}_1\cdot  \widehat{a}_1^{-1}\cdot \left(\sum\limits_{j = 1}^{n_1} \alpha_j a_j + \alpha \widehat{a}_1c\widehat{a}_2\right) = \sum\limits_{j = 1}^{n_1} \alpha_j \widehat{b}_1\cdot \widehat{a}_1^{-1}\cdot a_j + \alpha \widehat{b}_1c\widehat{a}_2 \notin \Rel
\end{equation*}
(even after the cancellations), or
\begin{equation*}
p\cdot \widehat{a}_2^{-1}\cdot \widehat{b}_2 = \left(\sum\limits_{j = 1}^{n_1} \alpha_j a_j  + \alpha \widehat{a}_1c\widehat{a}_2\right)\cdot \widehat{a}_2^{-1}\cdot \widehat{b}_2 = \sum\limits_{j = 1}^{n_1} \alpha_j a_j\cdot \widehat{a}_2^{-1}\cdot \widehat{b}_2 + \alpha \widehat{a}_1c\widehat{b}_2 \notin \Rel
\end{equation*}
(even after the cancellations). Then the monomial $c$ is called \emph{a small piece with respect to $\Rel$}.


We denote the set of all small pieces with respect to $\Rel$ by $\SP$\label{set_sp_def}. Clearly, $\SP \subseteq \Mon$. From the definition it follows that the set $\SP$ is closed under taking subwords. In particular, if the set $\SP$ is non-empty, the monomial $1$ is always a small piece. If the set $\SP$ is turned out to be empty, then \emph{we still assign $1$ to be a small piece}.
\end{definition}

The key property of monomials of $\Mon \setminus \SP$ is the following:
\begin{lemma}
\label{not_sp_prolongation}
Let $c \in \Mon \setminus \SP$, $p = \sum_{j = 1}^{n} \gamma_j c_j  + \gamma c \in \Rel$. Assume $a, b$ are monomials such that the monomial $ac$ has no cancellations, $ac \in \Mon$, and the monomial $cb$ has no cancellations, $cb \in \Mon$. Then $a\cdot p\cdot b \in \Rel$ (possibly after the cancellations). In particular, $acb \in \Mon$.
\end{lemma}
\begin{proof}
Since $ac \in \Mon$, there exists a polynomial $q_1 \in \Rel$ such that $ac$ is s monomial in $q_1$. That is, $q_1 = \sum_{j = 1}^{k_1}\eta_j d_j + \eta ac$. Since $ac$ has no cancellations, $c$ is a subword of monomials in the polynomials $p \in \Rel$ and $q_1 \in \Rel$. Assume $a\cdot p = \sum_{j = 1}^{n} \gamma_j a\cdot c_j  + \gamma ac \notin \Rel$. Then, by definition, $c$ is a small piece. This contradicts  the conditions of Lemma~\ref{not_sp_prolongation}. Hence, $a\cdot p \in \Rel$.

Since $cb \in \Mon$, there exists a polynomial $q_2 \in \Rel$ such that $cb$ is s monomial in $q_2$. That is, $q_2 = \sum_{j = 1}^{k_2}\delta_j t_j + \delta cb$. Since $cb$ and $ac$ have no cancellations, $c$ is a subword of monomials in the polynomials $a\cdot p \in \Rel$ and $q_2 \in \Rel$. Assume $a\cdot p\cdot b = \sum_{j = 1}^{n} \gamma_j a\cdot c_j\cdot b  + \gamma acb \notin \Rel$. Then again, by definition, $c$ is a small piece. This contradicts  the conditions of Lemma~\ref{not_sp_prolongation}. Thus, we finally obtain $a\cdot p\cdot b \in \Rel$.
\end{proof}

\begin{corollary}
\label{not_sp_prolongation2}
Let $c \in \Mon \setminus \SP$, $p = \sum_{j = 1}^{n} \gamma_j c_j  + \gamma c \in \Rel$. Assume $a, b$ are monomials such that the monomial $acb$ has no cancellations and $acb \in \Mon$. Then $a\cdot p\cdot b \in \Rel$ (possibly after the cancellations).
\end{corollary}
\begin{proof}
Assume $acb$ has no cancellations and $acb \in \Mon$. Since $\Mon$ is closed under taking subwords, we obtain $ac \in \Mon$ and $cb \in \Mon$. Hence, we can apply Lemma~\ref{not_sp_prolongation} and obtain $a\cdot p\cdot b \in \Rel$.
\end{proof}

In Subsection~\ref{tensor_product_section} we widely use the following obvious corollary of Compatibility Axiom and Corollary~\ref{not_sp_prolongation2}.
\begin{corollary}
\label{not_sp_cancellation_prolongation}
Let $p = \sum_{j = 1}^{n} \gamma_j c_j  + \gamma c \in \Rel$. Assume that $c = amb$, where $a$ is a prefix of $c$, $b$ is a suffix of $c$, and $m$ is not a small piece. Assume that $\widetilde{a}$ and $\widetilde{b}$ are monomials such that the monomial $\widetilde{a}m\widetilde{b}$ has no cancellations and $\widetilde{a}m\widetilde{b} \in \Mon$. Then $\widetilde{a}\cdot a^{-1}\cdot p\cdot b^{-1}\cdot \widetilde{b} \in \Rel$.
\end{corollary}
\begin{proof}
Since $c = amb$, it follows from Compatibility Axiom that $a^{-1}\cdot p\cdot b^{-1} \in \Rel$. So, we have $a^{-1}\cdot p\cdot b^{-1} = \sum_{j = 1}^{n} \gamma_j a^{-1}\cdot c_j\cdot b^{-1} + \gamma m \in \Rel$. We assumed that $\widetilde{a}m\widetilde{b}$ has no cancellations and $\widetilde{a}m\widetilde{b} \in \Mon$. Since $m$ is not a small piece, Corollary~\ref{not_sp_prolongation2} implies that
\begin{equation*}
\widetilde{a}\cdot a^{-1}\cdot p\cdot b^{-1}\cdot \widetilde{b} = \sum\limits_{j = 1}^{n} \gamma_j \widetilde{a}\cdot a^{-1}\cdot c_j\cdot b^{-1}\cdot \widetilde{b} + \gamma \widetilde{a}m\widetilde{b} \in \Rel.
\end{equation*}
\end{proof}

Let $u \in \Mon$. Then either $u = p_1\cdots p_n$, where $p_1, \ldots, p_n$ are small pieces, or $u$ can not be represented as a product of small pieces. We introduce a measure on monomials of $\Mon$ (aka \emph{$\SPM$-measure}). We say that\label{measure_def}
\begin{align*}
\SPM(u) = n \textit{ if } u &\textit{ can be represented as a product of small pieces}\\
&\textit{and minimal possible number of small pieces}\\
&\textit{in such representation is equal to } n.
\end{align*}
We say that
\begin{equation*}
\SPM(u) = \infty \textit{ if } u \textit{ can not be represented as a product of small pieces}.
\end{equation*}

\begin{remark}
We assume standard arithmetic rules for $\infty$. Namely,
\begin{align*}
&n < \infty \textit{ for every } n \in \mathbb{Z},\\
&\infty + n = n + \infty = \infty \textit{ for every } n \in \mathbb{Z},\\
&\infty + \infty = \infty.
\end{align*}
\end{remark}

\begin{lemma}
\label{measure_properties}
Let $a \in \Mon$. Assume $c$ is a prefix of $a$, that is, $a = cd$. Then we always have either $\SPM(a) = \SPM(c) + \SPM(d)$, or $\SPM(a) = \SPM(c) + \SPM(d) - 1$. In particular, if $c \in \SP$, then either $\SPM(a) = \SPM(d) + 1$, or $\SPM(a) = \SPM(d)$; if $d \in \SP$, then either $\SPM(a) = \SPM(c) + 1$, or $\SPM(a) = \SPM(c)$.

If $m$ is some subword of $a$, then $\SPM(m) \leqslant \SPM(a)$. In particular, if $\SPM(a)$ is finite, then $\SPM(m)$ is finite as well; if $\SPM(a) = \infty$, then $\SPM(m)$ may be both finite or infinite.
\end{lemma}
\begin{proof}
Assume $\SPM(a) = n < \infty$, then $a = p_1\cdots p_n$, where $p_1, \ldots, p_n$ are small pieces. Let $m$ be a subword of $a$. Then, clearly,
\begin{equation*}
m = p_{l_1}^{\prime}\cdot p_{l_1 + 1}\cdots p_{l_2}^{\prime},
\end{equation*}
where $p_{l_1}^{\prime}$ is a suffix of $p_{l_1}$ and $p_{l_2}^{\prime}$ is a prefix of $p_{l_2}$. Hence, $\SPM(m) \leqslant l_2 - l_1 + 1 \leqslant \SPM(a)$.

Assume $a = cd$ and $\SPM(a) = n < \infty$. Then $a = p_1\cdots p_n$, where $p_1, \ldots, p_n$ are small pieces. Hence,
\begin{align*}
&c = p_1\cdots p_l^{\prime},\\
&d = p_l^{\prime\prime}\cdots p_n,
\end{align*}
where $p_l = p_l^{\prime}p_l^{\prime\prime}$, one of $p_l^{\prime}$ or $p_l^{\prime\prime}$ can be empty. Then $\SPM(c) \leqslant l$, $\SPM(d) \leqslant n - l + 1$. Hence, $\SPM(c) + \SPM(d) \leqslant l + n - l + 1 = \SPM(a) + 1$.

On the other hand, if $\SPM(c) = n_1$ and $\SPM(d) = n_2$, then
\begin{align*}
&c = s_1\cdots s_{n_1},\\
&d = t_1\cdots t_{n_2},
\end{align*}
where $s_1,\ldots, s_{n_1}, t_1,\ldots, t_{n_2}$ are small pieces. So,
\begin{equation*}
a = cd = s_1\cdots s_{n_1}\cdot t_1\cdots t_{n_2}.
\end{equation*}
Hence, $\SPM(a) \leqslant n_1 + n_2 =  \SPM(c) + \SPM(d)$.

So, finally we obtain
\begin{equation*}
\SPM(c) + \SPM(d) - 1\leqslant \SPM(a) \leqslant \SPM(c) + \SPM(d).
\end{equation*}
Since values of $\SPM$-measure are natural numbers, we obtain
\begin{equation}
\label{measure_value}
\textit{either } \SPM(a) = \SPM(c) + \SPM(d), \textit{ or } \SPM(a) = \SPM(c) + \SPM(d) - 1.
\end{equation}

From the above argument it is clear that if $\SPM(c) < \infty$ and $\SPM(d) < \infty$, then also $\SPM(a) < \infty$. So, if $\SPM(a) = \infty$, then at least one of $\SPM(c)$ and $\SPM(d)$ is infinite. So, formula~\eqref{measure_value} is applicable for this case as well.
\end{proof}

We fix a constant $\tau \in \mathbb{N}$. In the further argument we require $\tau \geqslant 10$\label{tau_def}.
\begin{condition}[Small Cancellation Axiom]
\label{sc_ax}
Assume $q_1, \ldots, q_n \in \Rel$ and a linear combination $\sum_{l = 1}^{n} \gamma_l q_l$ is non-zero after additive cancellations. Then there exists a monomial $a$ in $\sum_{l= 1}^{n} \gamma_l q_l$ with a non-zero coefficient after additive cancellations such that either $a$ can not be represented as a product of small pieces or every representation of $a$ as a product of small pieces contains at least $\tau + 1$ small pieces. That is, $\SPM(a) \geqslant \tau + 1$, including $\SPM(a) = \infty$.
\end{condition}

\begin{remark}
\label{small_pieces_comment}
We want to explain informally a source of a concept of a small piece in the ring case (see Definition~\ref{sp}).

Let $G = \langle \mathcal{X}\mid \Rel_G\rangle$, $\Rel_G = \lbrace R_j\rbrace_{j\in J}$ be a group given by generators and defining relations. Assume $\Rel_G$ is closed under taking inverses and cyclic shifts of relators, and every relator $R_j$ is a cyclically reduced word. Basically, we take the idea of Definition~\ref{sp} of a small piece from the following concept. Let $s$ be a prefix of $R_{j_1} \in \Rel_G$ and of $R_{j_2} \in \Rel_G$, $R_{j_1} = sR_{j_1}^{\prime}$, $R_{j_2} = sR_{j_2}^{\prime}$. Then $s$ is called \emph{a small piece (in a generalized group sense)} if $R_{j_1}^{\prime} \cdot {R_{j_2}^{\prime}}^{-1} \neq 1$ and $R_{j_1}^{\prime} \cdot {R_{j_2}^{\prime}}^{-1}$ is not a conjugate of a relator from $\Rel_G$ in the corresponding free group (even after the cancellations) (see~\cite{LS}, page~240 and page~271, condition~(1), and~\cite{R}, page~5). That is, if $s$ is not a small piece, we obtain that either $R_{j_1} = R_{j_2}$, or $R_{j_1}^{\prime} \cdot {R_{j_2}^{\prime}}^{-1}$ is equal to a conjugate of a relator from $\Rel_G$ in the corresponding free group (possibly after the cancellations).

We consider the Cayley graph of the group $G$ with respect to the set of generators $\mathcal{X}$. Then every $R_j \in \Rel_G$ corresponds to a closed path in the Cayley graph. Assume $s$ is a maximal common prefix of $R_{j_1} = sR_{j_1}^{\prime}$ and $R_{j_2} = sR_{j_2}^{\prime}$, and $s$ is a maximal common suffix of $R_{j_1}^{\prime}s$ and $R_{j_2}^{\prime}s$. Then $R_{j_1}^{\prime} \cdot {R_{j_2}^{\prime}}^{-1}$ has no cancellations. If we glue together the graphs that correspond to $R_{j_1}$ and $R_{j_2}$ by $s$, then $R_{j_1}^{\prime}{R_{j_2}^{\prime}}^{-1}$ corresponds to the path around the outer bypass of the obtained graph.
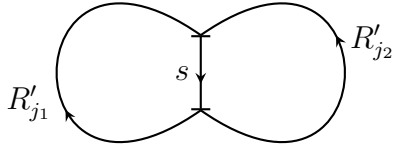
\begin{center}
\begin{tikzpicture}
\draw[|-|, black, thick, reversearrow=0.5] (0,0)--(0,1) node [midway, left] {$s$};
\draw[black, thick, arrow=0.4] (0,0) to [bend left=125, looseness=8] (0,1);
\draw[black, thick, arrow=0.6] (0,0) to [bend right=125, looseness=8] (0,1);
\path (0,0) to [bend left=125, looseness=8] node[pos=0.4, left] {$R_{j_1}^{\prime}$} (0,1);
\path (0,0) to [bend right=125, looseness=8] node[pos=0.6, right] {$R_{j_2}^{\prime}$} (0,1);
\end{tikzpicture}
\end{center}

However, in Definition~\ref{sp} there are certain modifications comparatively with small pieces in a generalized group sense. Now we want to compare informally Definition~\ref{sp} of a small piece in the ring case with the given above definition of a small piece in the group case.
\begin{enumerate}
\item
Let us analyze Definition~\ref{sp}. Assume $c \in \Mon$ is not a small piece in a sense of Definition~\ref{sp}. Let
\begin{align*}
&\widehat{a}_1c\widehat{a}_2 + \sum_{j = 1}^{n_1}\alpha_ja_j \in \Rel,\\
&c + \sum_{j = 1}^{n_2}\beta_jb_j \in \Rel.
\end{align*}
Clearly, $c = -\sum_{j = 1}^{n_2}\beta_jb_j \mod \Ideal$. We want to replace $c$ by the linear combination $\left(-\sum_{j = 1}^{n_2}\beta_jb_j\right)$ in the monomial $\widehat{a}_1c\widehat{a}_2$. If we follow the group analogue completely, we should require that
\begin{align*}
\widehat{a}_1\cdot \left(-\sum_{j = 1}^{n_2}\beta_jb_j\right)\cdot \widehat{a}_2 &+ \sum_{j = 1}^{n_1}\alpha_ja_j =\\
&= \widehat{a}_1c\widehat{a}_2 + \sum_{j = 1}^{n_1}\alpha_ja_j - \widehat{a}_1\cdot \left(c + \sum_{j = 1}^{n_2}\beta_jb_j\right)\cdot \widehat{a}_2 \in \Rel.
\end{align*}
But this condition is not convenient for the further work for some reasons. So, in Definition~\ref{sp} we require its modification
\begin{equation*}
\widehat{a}_1\cdot \left(c + \sum_{j = 1}^{n_2}\beta_jb_j\right)\cdot \widehat{a}_2 = \widehat{a}_1c\widehat{a}_2 + \sum_{j = 1}^{n_2}\beta_j\widehat{a}_1\cdot b_j\cdot \widehat{a}_2 \in \Rel
\end{equation*}
instead.

Now let us show that in the above condition we still follow the same ideology as in the case of groups. Let us explain this in more detail. Assume $AsBD \in \Rel_G$, $sE \in \Rel_G$, and $s$ is not a small piece in generalized group sense. Consider the prefix $AsB$ of $AsBD$. Let us replace $s$ by $E^{-1}$ in $AsB$, then we obtain $AE^{-1}B$ as a result. Assume the word $AE^{-1}B$ has no cancellations. Obviously, $AE^{-1}B$ is a subword of $AE^{-1}BD$. Since $s$ is not a small piece, we get that $AE^{-1}BD \in \Rel_G$. Therefore, the resulting word $AE^{-1}B$ is a subword of a relator from~$\Rel_G$ as the initial word $AsB$.
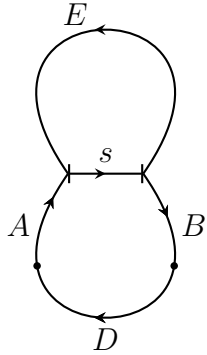
\begin{equation}
\begin{tikzpicture}
\draw[|-|, black, thick, arrow=0.5] (0,0)--(1,0) node [midway, above] {$s$};
\draw[black, thick, reversearrow=0.5] (0,0) to [bend left=125, looseness=8] (1,0);
\draw[black, thick, reversearrow=0.1, reversearrow=0.5, reversearrow=0.9] (0,0) to [bend right=125, looseness=8] coordinate[pos=0.2] (P1) coordinate[pos=0.8] (P2) (1,0);
\node[circle, fill, inner sep=1] at (P1) {};
\node[circle, fill, inner sep=1] at (P2) {};
\path (0,0) to [bend left=125, looseness=8] node[pos=0.4, above] {$E$} (1,0);
\path (0,0) to [bend right=125, looseness=8] node[pos=0.1, left] {$A$} (1,0);
\path (0,0) to [bend right=125, looseness=8] node[pos=0.9, right] {$B$} (1,0);
\path (0,0) to [bend right=125, looseness=8] node[pos=0.5, below] {$D$} (1,0);
\end{tikzpicture}
\label{replacement_of_non_sp_in_group}
\end{equation}

Similar situation happens for the case of rings. Assume $c \in \Mon$ is not a small piece in the sense of Definition~\ref{sp}. We replace $c$ by the linear combination $\left(-\sum_{j = 1}^{n_2}\beta_jb_j\right)$ in $\widehat{a}_1c\widehat{a}_2$. By Definition~\ref{sp},
\begin{equation*}
\widehat{a}_1\cdot \left(c + \sum_{j = 1}^{n_2}\beta_jb_j\right)\cdot \widehat{a}_2 = \widehat{a}_1c\widehat{a}_2 + \sum_{j = 1}^{n_2}\beta_j\widehat{a}_1\cdot b_j\cdot \widehat{a}_2 \in \Rel.
\end{equation*}
Therefore, we do not obtain any completely new resulting monomials in $\sum_{j = 1}^{n_2}\beta_j\widehat{a}_1\cdot b_j\cdot \widehat{a}_2$,  which means that all monomials $\widehat{a}_1b_j\widehat{a}_2 \in \Mon$ (possibly after the cancellations).

\item
Assume $c\in \Mon$ is not a small piece in the sense of Definition~\ref{sp}. Let
\begin{align*}
&\widehat{a}_1c + \sum_{j = 1}^{n_1}\alpha_ja_j \in \Rel,\\
&c\widehat{b}_2 + \sum_{j = 1}^{n_2}\beta_jb_j \in \Rel.
\end{align*}
Then in Definition~\ref{sp} we require that
\begin{align*}
&\widehat{a}_1c\widehat{b}_2 + \sum_{j = 1}^{n_1}\alpha_ja_j\cdot \widehat{b}_2 \in \Rel,\\
&\widehat{a}_1c\widehat{b}_2 + \sum_{j = 1}^{n_2}\beta_j\widehat{a}_1\cdot b_j \in \Rel.
\end{align*}
This allows us to glue via $c$ the monomials $\widehat{a}_1c \in \Mon$ and $c\widehat{b}_2 \in \Mon$ to one monomial $\widehat{a}_1c\widehat{b}_2$ from $\Mon$. Notice that this property does not necessarily hold in the case of groups for subwords of relators from~$\Rel_G$. For example, looking at picture~\eqref{replacement_of_non_sp_in_group}, $AsE$ is not necessarily a subword of a relator from~$\Rel_G$.

The reader will notice wide applications of this property already in Section~\ref{basic_def} when we start to discuss overlaps of occurrences.
\end{enumerate}

The group $G = \langle \mathcal{X}\mid \Rel_G\rangle$ satisfies \emph{small cancellation condition $C(m)$ in a generalized sense}\label{cond_cm_1} if every $R_j \in \Rel_G$ can not be written as a product of less than $m$ small pieces (in a generalized group sense). One can see that Small Cancellation Axiom stated above is an analogue of condition~$C(m)$.
\end{remark}

The next lemma follows from Compatibility Axiom and Small Cancellation Axiom.
\begin{lemma}
\label{shifts_of_relations}
Let $p = \sum_{j = 1}^{n} \gamma_j c_j \in \Rel$, $a$, $b$ be some monomials. Assume $a\cdot c_j\cdot b \in \Mon$ for all $j = 1, \ldots, n$ (possibly after the cancellations). Then $a\cdot p\cdot b \in \Rel$.
\end{lemma}
\begin{proof}
First assume $a$ and $b$ have no cancellation with all monomials $c_j$, $j = 1, \ldots, n$. From Small Cancellation Axiom it follows that there exists a monomial $c_{j_0}$ that is not a small piece. Therefore, since $ac_{j_0}b \in \Mon$ and $ac_{j_0}b$ have no cancellations, Lemma~\ref{not_sp_prolongation} implies $a\cdot p\cdot b \in \Rel$.

Further we argue by induction on $\vert a\vert + \vert b\vert$, where $\vert \cdot \vert$ is number of symbols of $S \cup S^{-1}$ in a reduced word. Assume $a = a_1x$, $x\in S\cup S^{-1}$, $x$ cancels with at least one monomial in $p$. Let $p_1 = x\cdot p$, then $p_1 \in \Rel$, by Compatibility Axiom. Since $a_1\cdot p_1\cdot b = a\cdot p\cdot b$ (after the cancellations from both sides), all monomials of $a_1\cdot p_1\cdot b$ belong to $\Mon$. Hence, $a_1\cdot p_1\cdot b \in \Rel$, by the induction hypothesis. Clearly, if $a$ is empty or does not cancel with all monomials $c_j$, $j = 1, \ldots, n$, we can argue in the same way with $b$.
\end{proof}

\begin{definition}
\label{incident_momom}
Let $p = \sum_{j = 1}^n \alpha_j a_j \in \Rel$. Then we call the monomials $a_{j_1}, a_{j_2}$, $1 \leqslant j_1, j_2 \leqslant n$, \emph{incident monomials} (including the case $a_{j_1} = a_{j_2}$). Recall that $\alpha_j \neq 0$, $j = 1, \ldots, n$.
\end{definition}

The next straightforward lemma follows directly from Compatibility Axiom and Lemma~\ref{not_sp_prolongation}.
\begin{lemma}
\label{compatibility_for_incedent}
Let $m, m^{\prime} \in \Mon$ be incident monomials, $c$ be some monomial.
\begin{enumerate}
\item
If $c^{-1}$ is a suffix of $m$ or $m^{\prime}$ (a prefix of $m$ or $m^{\prime}$), then $m\cdot c, m^{\prime}\cdot c \in \Mon$ ($c\cdot m, c\cdot m^{\prime} \in \Mon$) and $m\cdot c, m^{\prime}\cdot c$ ($c\cdot m, c\cdot m^{\prime}$) are incident monomials after the cancellations.
\item
If $m$ is not a small piece and $c$ has no cancellations with $m$ from the right side (from the left side) and $mc \in \Mon$ ($cm \in \Mon$), then $m^{\prime}\cdot c \in \Mon$ ($c\cdot m^{\prime} \in \Mon$) and $mc, m^{\prime}\cdot c$ ($cm, c\cdot m^{\prime}$) are incident monomials possibly after the cancellations.
\end{enumerate}
\end{lemma}
\begin{proof}
Since $m$ and $m^{\prime}$ are incident monomials, there exists $p = \sum_{j = 1}^n \alpha_ja_j \in \Rel$ such that $m = a_{j_1}$ and $m^{\prime} = a_{j_2}$, $j_1, j_2 \in \lbrace 1, \ldots, n\rbrace$. If $c^{-1}$ is a suffix of $m$ or $m^{\prime}$ (a prefix of $m$ or $m^{\prime}$), then, by Compatibility Axiom, we obtain $p\cdot c \in \Rel$ ($cp \in \Rel$) after the cancellations. In particular, $m\cdot c, m^{\prime}\cdot c \in \Mon$ ($c\cdot m, c\cdot m^{\prime} \in \Mon$).

If $m$ is not a small piece, $c$ has no cancellations with $m$ from the right side (from the left side) and $mc \in \Mon$ ($cm \in \Mon$), then Lemma~\ref{not_sp_prolongation} implies that $p\cdot c \in \Rel$ ($c\cdot p \in \Rel$). In particular, $m^{\prime}\cdot c \in \Mon$ ($c\cdot m^{\prime} \in \Mon$).

So, in both cases $p\cdot c \in \Rel$ ($cp \in \Rel$). Clearly, $a_{j_1}\cdot c = m\cdot c$ and  $a_{j_2}\cdot c = m^{\prime}\cdot c$ ($c\cdot a_{j_1} = c\cdot m$ and $c\cdot a_{j_2} = c\cdot m^{\prime}$) are monomials in $p\cdot c \in \Rel$ ($c\cdot p \in \Rel$). Therefore, by definition, $m\cdot c$ and $m^{\prime}\cdot c$ ($c\cdot m$ and $c\cdot m^{\prime}$) are incident monomials after the cancellations.
\end{proof}

Now we introduce the last condition, we call it Isolation Axiom. Unlike previous two axioms, this is entirely ring-theoretic condition. Here we use notions of a maximal occurrence and an overlap that we introduce in Section~\ref{basic_def} (see Definition~\ref{max_occurrence_def}, Definition~\ref{overlap_def} and the list of possibilities~\ref{separated}---\ref{overlap} on page~\pageref{separated}). We wish to draw the reader's attention to the  fact  that Isolation Axiom is used only in order to prove statement~\ref{mu_properties_up_correspondence} of Lemma~\ref{mu_properties} in Subsection~\ref{tensor_product_section}. Before this place we develop the theory without this condition. We urge the reader to skip the statement of Isolation Axiom for the first reading and to return to it when it becomes necessary.

\begin{condition3a}[Isolation Axiom, left-sided]
Let $m_1, m_2, \ldots, m_{k}$ be a sequence of monomials of $\Mon$ such that $m_1 \neq m_k$ and $m_i, m_{i + 1}$ are incident monomials for all $i = 1, \ldots, k - 1$, and $\SPM(m_i) \geqslant \tau - 2$ for all $i = 1, \ldots, k$. Let us take a monomial $a \in \Mon$ with the following properties.
\begin{itemize}
\item
$\SPM(a) \geqslant \tau - 2$;
\item
$am_1, am_k \notin \Mon$, $am_1$ has no cancellations, $am_k$ has no cancellations;
\item
$m_1$ is a maximal occurrence in $am_1$, $m_k$ is a maximal occurrence in $am_k$.
\item
Let $ap_1(a)$ be a maximal occurrence in $am_1$ that contains $a$, let $ap_k(a)$ be a maximal occurrence in $am_k$ that contains $a$ (that is, $p_1(a)$ is the overlap of $ap_1(a)$ and $m_1$, $p_1(a)$ may be empty, and $p_k(a)$ is the overlap of $ap_k(a)$ and $m_k$, $p_k(a)$ may be empty). Assume that there exist monomials $l$, $l^{\prime} \in \Mon$ such that
\begin{itemize}
\item
$l$, $l^{\prime}$ are small pieces;
\item
$la, l^{\prime}a \in \Mon$, $la$ has no cancellations, $l^{\prime}a$ has no cancellations;
\item
there exists a sequence of monomials $b_1, \ldots, b_n$ from $\Mon$ such that $b_1 = lap_1(a)$, $b_n = l^{\prime}ap_k(a)$, $b_i, b_{i + 1}$ are incident monomials for all $i = 1, \ldots, n - 1$, and $\SPM(b_i) \geqslant \tau - 2$ for all $i = 1, \ldots, n$.
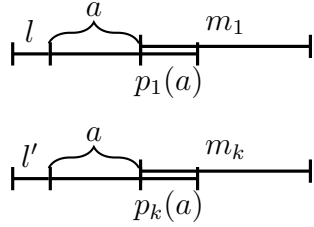
\begin{center}
\begin{tikzpicture}
\draw[|-|, black, very thick] (2.2,0.1)--(4.5,0.1) node[midway, above] {$m_1$};
\draw[|-|, black, very thick] (2.2, 0)--(3, 0) node[midway, below] {$p_1(a)$};
\draw[|-, black, very thick] (1,0)--(2.2,0);
\draw[|-, black, very thick] (0.5, 0)--(1,0) node[midway, above] {$l$};
\draw [thick, decorate, decoration={brace, amplitude=8pt, raise=4pt}] (1, 0) to node[midway, above, yshift=10pt] {$a$} (2.2, 0);
\end{tikzpicture}

\begin{tikzpicture}
\draw[|-|, black, very thick] (2.2,0.1)--(4.5,0.1) node[midway, above] {$m_k$};
\draw[|-|, black, very thick] (2.2, 0)--(3, 0) node[midway, below] {$p_k(a)$};
\draw[|-, black, very thick] (1,0)--(2.2,0);
\draw[|-, black, very thick] (0.5, 0)--(1,0) node[midway, above] {$l^{\prime}$};
\draw [thick, decorate, decoration={brace, amplitude=8pt, raise=4pt}] (1, 0) to node[midway, above, yshift=10pt] {$a$} (2.2, 0);
\end{tikzpicture}
\end{center}
Notice that since $a$ is not a small piece, by Lemma~\ref{max_overlap}, we get that $lap_1(a), l^{\prime}ap_k(a) \in \Mon$, and $lap_1(a)$ is a maximal occurrence in $lap_1(a)m_1$, $l^{\prime}ap_k(a)$ is a maximal occurrence in $l^{\prime}ap_k(a)m_k$.
\end{itemize}
\end{itemize}
Then we require that ${p_1(a)}^{-1}\cdot m_1 \neq {p_k(a)}^{-1}\cdot m_k$ for every such $a \in \Mon$.
\end{condition3a}

\begin{condition3b}[Isolation Axiom, right-sided]
Let $m_1, m_2, \ldots, m_{k}$ be a sequence of monomials of $\Mon$ such that $m_1 \neq m_k$ and $m_i, m_{i + 1}$ are incident monomials for all $i = 1, \ldots, k - 1$, and $\SPM(m_i) \geqslant \tau - 2$ for all $i = 1, \ldots, k$. Let us take a monomial $a \in \Mon$ such that
\begin{itemize}
\item
$\SPM(a) \geqslant \tau - 2$;
\item
$m_1a, m_ka \notin \Mon$, $m_1a$ has no cancellations, $m_ka$ has no cancellations;
\item
$m_1$ is a maximal occurrence in $m_1a$, $m_k$ is a maximal occurrence in $m_ka$.
\item
Let $s_1(a)a$ be a maximal occurrence in $m_1a$ that contains $a$, let $s_k(a)a$ be a maximal occurrence in $m_ka$ that contains $a$ (that is, $s_1(a)$ is the overlap of $m_1$ and $s_1(a)a$, $s_1(a)$ may be empty, and $s_k(a)$ is the overlap of $m_k$ and $s_k(a)a$, $s_k(a)$ may be empty). Assume that there exist monomials $r$, $r^{\prime} \in \Mon$ such that
\begin{itemize}
\item
$r$, $r^{\prime}$ are small pieces;
\item
$ar, ar^{\prime} \in \Mon$, $ar$ has no cancellations, $ar^{\prime}$ has no cancellations;
\item
there exists a sequence of monomials $b_1, \ldots, b_n$ from $\Mon$ such that $b_1 = s_1(a)ar$, $b_n = s_k(a)ar^{\prime}$, $b_i, b_{i + 1}$ are incident monomials for all $i = 1, \ldots, n - 1$, and $\SPM(b_i) \geqslant \tau - 2$ for all $i = 1, \ldots, n$.
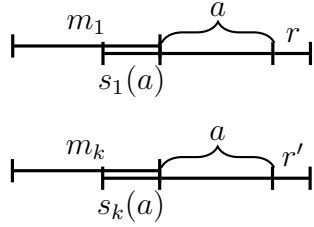
\begin{center}
\begin{tikzpicture}
\draw[-|, black, very thick] (3,0)--(4.5,0);
\draw[|-|, black, very thick] (2.2, 0)--(3, 0) node[midway, below] {$s_1(a)$};
\draw[|-|, black, very thick] (1,0.1)--(3,0.1) node[midway, above] {$m_1$};
\draw[-|, black, very thick] (4.5,0)--(5,0) node[midway, above] {$r$};
\draw [thick, decorate, decoration={brace, amplitude=8pt, raise=4pt}] (3, 0) to node[midway, above, yshift=10pt] {$a$} (4.5, 0);
\end{tikzpicture}

\begin{tikzpicture}
\draw[-|, black, very thick] (3, 0)--(4.5, 0);
\draw[|-|, black, very thick] (2.2, 0)--(3, 0) node[midway, below] {$s_k(a)$};
\draw[|-|, black, very thick] (1,0.1)--(3,0.1) node[midway, above] {$m_k$};
\draw[-|, black, very thick] (4.5,0)--(5,0) node[midway, above] {$r^{\prime}$};
\draw [thick, decorate, decoration={brace, amplitude=8pt, raise=4pt}] (3, 0) to node[midway, above, yshift=10pt] {$a$} (4.5, 0);
\end{tikzpicture}
\end{center}
Notice that since $a$ is not a small piece, by Lemma~\ref{max_overlap}, we get that $s_1(a)ar, s_k(a)ar^{\prime} \in \Mon$, and $s_1(a)ar$ is a maximal occurrence in $m_1s_1(a)ar$, $s_k(a)ar^{\prime}$ is a maximal occurrence in $m_ks_k(a)ar^{\prime}$.
\end{itemize}
\end{itemize}
Then we require that $m_1\cdot {s_1(a)}^{-1} \neq m_k\cdot {s_k(a)}^{-1}$ for every such $a \in \Mon$.
\end{condition3b}

\begin{remark}
\label{isolation_ax_meaning}
We shall informally explain the essence of Isolation Axioms. Given a monomial $U$, consider the set of its non-degenerate derived monomials (see Subsection~\ref{derived_monom_def_section} for the definition of derived monomials). Every derived monomial can be imagined as a result of a sequence of replacements of virtual members of the chart by $U$-incident monomials (see Subsection~\ref{replacements_by_U_incident_section} for the definition of $U$-incident monomials). If two essentially different sequences of replacements result in one and the same derived monomial, the exotic dependencies appear in the ideal $\Ideal$. Isolation Axiom guarantees that essentially different sequences of replacements result in different monomials. Hence, exotic dependencies are not present in~$\Ideal$.
\end{remark}

\section{Basic definitions}
\label{basic_def}
In this section we introduce notions of \emph{the chart of a monomial} and \emph{a multi-turn}. Both notions play an important role in the further argument.

Let $U$ be a word and $\widehat{U}$ be its subword. We call the triple that consists of $U$, $\widehat{U}$ and the position of $\widehat{U}$ in $U$ \emph{an occurrence of $\widehat{U}$ in $U$}.

Let $U$ be a monomial, $a \in \Mon$ be an occurrence in $U$, that is, $U = LaR$, where $L$, $R$ can be empty. Since $a \in \Mon$, there exists a polynomial $p \in \Rel$ such that $a$ is a monomial of $p$. Assume $L$ is not empty, $X$ is a suffix of $L$, $L = L_1X$, $L_1$ is possibly empty. If $X\cdot p \in \Rel$ (possibly after the cancellations), then we say that \emph{$X$ prolongs $a$ in $U$ from the left with respect to $p$}. In particular, $Xa \in \Mon$ in this case. Assume $R$ is not empty, $Y$ is a prefix of $R$, $R = YR_1$, $R_1$ is possibly empty. If $p\cdot Y \in \Rel$ (possibly after the cancellations), then we say that \emph{$Y$ prolongs $a$ in $U$ from the right with respect to $p$}. In particular, $aY \in \Mon$ in this case.

\begin{remark}
Let $U$ be a monomial, $a \in \Mon$ be an occurrence in $U$, $U = LaR$. Let $a$ be a monomial in a polynomial $p \in \Rel$ and in some other polynomial $q \in \Rel$.

Assume $a$ is not a small piece. Assume $X$ is a suffix of $L$ such that $X$ prolongs $a$ from the left with respect to $p$ (that is, $X\cdot p \in \Rel$). Since $a$ is not a small piece and $Xa \in \Mon$ and $Xa$ has no cancellations, it follows from Lemma~\ref{not_sp_prolongation} that $X\cdot q \in \Rel$. That is, $X$ also prolongs $a$ from the left with respect to $q$. So, if $a$ is not a small piece, its possible prolongations in $U$ from the left do not depend on a particular relation in which $a$ is a monomial. Similarly, if $a$ is not a small piece, its possible prolongations in $U$ from the right do not depend on a particular relation in which $a$ is a monomial.

On the contrary, if $a$ is a small piece, then it may happen that $a$ is prolonged in $U$ with respect to $p$ and is not prolonged in $U$ with respect to~$q$.
\end{remark}

Now we give a definition of \emph{a maximal occurrence of a monomial of $\Mon$ in $U$}.
\begin{definition}
\label{max_occurrence_one_rel_def}
Let $U$ be a monomial, $a \in \Mon$ be an occurrence in $U$. Let $p \in \Rel$, $a$ be a monomial in $p$. The occurrence $a$ is called \emph{maximal in $U$ with respect to $p$} if $a$ can not be prolonged neither to the left nor to the right (even by a single letter) in $U$ with respect to $p$.

Let $p = \sum_{j = 1}^n \alpha_ja_j + \alpha a \in \Rel$, $U = LaR$. In more detail, this definition means the following.
\begin{enumerate}
\item
If both $L$ and $R$ are not empty, $L = L_1x$, $R = yR_1$, $x, y \in S \cup S^{-1}$, then both
\begin{align*}
&x\cdot p = \sum\limits_{j = 1}^n \alpha_j x\cdot a_j + \alpha xa \notin \Rel,\\
&p\cdot y = \sum\limits_{j = 1}^n \alpha_j a_j\cdot y + \alpha ay \notin \Rel.
\end{align*}
\item
If $L$ is not empty and $R$ is empty, that is, $U = La$, $L = L_1x$, $x \in S \cup S^{-1}$, then
\begin{equation*}
x\cdot p = \sum\limits_{j = 1}^n \alpha_j x\cdot a_j + \alpha xa \notin \Rel.
\end{equation*}
\item
If $L$ is empty and $R$ is not empty, that is, $U = aR$, $R = yR_1$, $y \in S \cup S^{-1}$, then
\begin{equation*}
p\cdot y = \sum\limits_{j = 1}^n \alpha_j a_j\cdot y + \alpha ay \notin \Rel.
\end{equation*}
\item
If both $L$ and $R$ are empty, that is, $U = a$, then there are no additional conditions, $a$ is always a maximal occurrence in $U$ with respect to $p$.
\end{enumerate}
\end{definition}

\begin{definition}
\label{max_occurrence_def}
Let $U$ be a monomial, $a \in \Mon$ be an occurrence in $U$. The occurrence $a$ is called \emph{maximal in $U$} if for every $p \in \Rel$ such that $a$ is a monomial in $p$ the occurrence $a$ is maximal in $U$ with respect to $p$.
\end{definition}

The following lemma gives a very natural characterisation of maximal occurrences of monomials of $\Mon$.
\begin{lemma}
\label{max_can_not_be_inside}
Let $U$ be a monomial, $a \in \Mon$ be an occurrence in $U$. Then $a$ is a maximal occurrence of a monomial of $\Mon$ in $U$ if and only if $a$ is not properly contained in any other occurrence of a monomial of $\Mon$ in $U$.
\end{lemma}
\begin{proof}
Assume $b \in \Mon$ is an occurrence in $U$, $a$ is properly contained inside $b$. That is $b = cad$, where $c$ is a prefix of $b$, $d$ is a suffix of $d$, one of $c$ or $d$ can be empty. Since $b \in \Mon$, there exists $p = \sum_{j = 1}^n \beta_j b_j + \beta b \in \Rel$. Since $b = cad$, it follows from Compatibility Axiom that $c^{-1}\cdot p \in \Rel$, $p\cdot d^{-1} \in \Rel$ and $c^{-1}\cdot p\cdot d^{-1} \in \Rel$. Denote $c^{-1}\cdot p\cdot d^{-1}$ by $q$, $q\in \Rel$. Then $q = \sum_{j = 1}^n \beta_j c^{-1}\cdot b_j\cdot d^{-1} + \beta a$, that is, $a$ is a monomial in $q \in \Rel$. We have $c\cdot q = p\cdot d^{-1} \in \Rel$ and $q\cdot d = c^{-1}\cdot p \in \Rel$, hence, by definition, $c$ prolongs $a$ from the left in $U$ with respect to $q$ and $d$ prolongs $a$ from the right in $U$ with respect to $q$. Thus, $a$ is not a maximal occurrence of a monomial of $\Mon$ in $U$. This contradiction completes the proof.
\end{proof}

\begin{remark}
Let $U$ be a monomial, $a \in \Mon$ be an occurrence in $U$. We proved that $a$ is a maximal occurrence of a monomial of $\Mon$ in $U$ if and only if $a$ is not properly contained in any other occurrence of a monomial of $\Mon$ in $U$. It can also be stated in the following way.
\begin{enumerate}
\item
If both $L$ and $R$ are not empty, $L = L_1x$, $R = yR_1$, $x, y \in S \cup S^{-1}$, then $a$ is a maximal occurrence in $U$ if and only if $xa \notin \Mon$ and $ay \notin \Mon$.
\item
If $L$ is not empty and $R$ is empty, that is, $U = La$, $L = L_1x$, $x \in S \cup S^{-1}$, then $a$ is a maximal occurrence in $U$ if and only if $xa \notin \Mon$.
\item
If $L$ is empty and $R$ is not empty, that is, $U = aR$, $R = yR_1$, $y \in S \cup S^{-1}$, then $a$ is a maximal occurrence in $U$ if and only if $ay \notin \Mon$.
\item
If both $L$ and $R$ are empty, that is, $U = a$, then there are no additional conditions, $a$ is always a maximal occurrence in $U$.
\end{enumerate}
\end{remark}

\begin{lemma}
Let $U$ be a monomial, $a \in \Mon$ be an occurrence in $U.$ Assume $a$ is not a small piece, $a$ is a monomial in $p \in \Rel$. Let $a$ be a maximal occurrence of a monomial of $\Mon$ in $U$ with respect to $p$. Then $a$ is a maximal occurrence of a monomial of $\Mon$ in $U$ (that is, $a$ is a maximal occurrence in $U$ with respect to every $q \in \Rel$ such that $a$ is a monomial in $q$).
\end{lemma}
\begin{proof}
Let $U = LaR$. Let $a$ be a monomial of a polynomial $q \in \Rel$ different from $p$. Assume $a$ is not a maximal occurrence of a monomial of $\Mon$ in $U$ with respect to $q$. Assume $a$ can be prolonged in $U$ from the left with respect to $q$. That is, $L = L_1x$, $x \in S \cup S^{-1}$, $x\cdot q \in \Rel$. In particular, we obtain $xa \in \Mon$. Since $a$ is not a small piece, from Lemma~\ref{not_sp_prolongation} it follows that $x\cdot p \in \Rel$. This contradicts the assumption that $a$ is a maximal occurrence in $U$ with respect to $p$. The case when $a$ can be prolonged in $U$ from the right with respect to $q$ is considered in the same way. Thus, $a$ is a maximal occurrence in $U$ with respect to $q$. Since $q$ is an arbitrary polynomial of $\Rel$ such that $a$ is a monomial of $q$, finally we obtain that $a$ is a maximal occurrence in $U$.
\end{proof}

\begin{remark}
Let $U$ be a monomial, $a \in \Mon$ be an occurrence in $U$. If $a$ is a small piece, the situation is different. Let $a$ be a monomial in $p \in \Rel$ and in $q \in \Rel$. Then it is possible that $a$ can not be prolonged in $U$ with respect to $p$ and can be prolonged in $U$ with respect to $q$. Then $a$ is a maximal occurrence in $U$ with respect to $p$ and is not a maximal occurrence in $U$ with respect to $q$. Hence, if $a$ is a maximal occurrence in $U$ with respect to $p$ and $a$ is a small piece, $a$ can be contained in some different occurrence of a monomial of $\Mon$ in $U$.
\end{remark}

Further we speak only about maximal occurrences of monomials of~$\Mon$ in~$U$. We call them for short ``maximal occurrences in $U$''.

\begin{definition}
\label{overlap_def}
Let $U$ be a monomial. Let $a$ and $b$ be two different occurrences of monomials of $\Mon$ in $U$. Assume $a$ is not contained inside $b$, $b$ is not contained inside $a$, and $a$ and $b$ have a non-empty common subword in $U$. Then we call this common subword \emph{an overlap of $a$ and $b$}.
\end{definition}

\begin{remark}
\label{max_overlap}
Let $U$ be a monomial. Let $a$ and $b$ be two different maximal occurrences of monomials of $\Mon$ in $U$. Then, by Lemma~\ref{max_can_not_be_inside}, $a$ is not contained inside $b$ and $b$ is not contained inside $a$. Hence, if $a$ and $b$ have a non-empty common subword in $U$, then $a$ and $b$ have an overlap.

We denote by $\SP^{\prime}$\label{max_overlaps_set} the set of all overlaps of maximal occurrences in all monomials, including the empty word.
\end{remark}

\begin{lemma}
\label{one_max_occurrence_interactions}
Let $U$ be a monomial. Let $a$ be a maximal occurrence of a monomial of $\Mon$ in $U$, $b$ be some occurrence of a monomial of $\Mon$ in $U$ (not necessarily maximal). If $a$ and $b$ have an overlap, then this overlap is a small piece.
\end{lemma}
\begin{proof}
Assume $a$ starts from the left of the beginning of $b$. The case when $a$ starts from the right of the beginning of $b$ is considered in the same way. Assume $a$ and $b$ have an overlap $c$ in $U$. That is, $a = c_1c$, $b = cc_2$, where $c_1$ and $c_2$ are non-empty. Assume $c$ is not a small piece. Assume $a$ is a monomial in a polynomial $p \in \Rel$, $b$ is a monomial in a polynomial $q \in \Rel$. Let
\begin{align*}
&p = \sum_{j = 1}^{n_1}\alpha_j a_j + \alpha a = \sum_{j = 1}^{n_1}\alpha_j a_j + \alpha c_1c,\\
&q = \sum_{j = 1}^{n_2}\beta_j b_j + \beta b = \sum_{j = 1}^{n_2}\beta_j b_j + \beta cc_2.
\end{align*}
Then, since $c$ is not a small piece, by definition, we obtain $p\cdot c_2 \in \Rel$. Since this holds for every polynomial $p \in \Rel$ such that $a$ is a monomial of $p$, we obtain that $a$ is not a maximal occurrence in $U$. A contradiction.
\end{proof}

\begin{corollary}
\label{max_occurrences_interactions}
Let $U$ be a monomial. Let $a$ and $b$ be two different maximal occurrences of monomials of $\Mon$ in $U$. If $a$ and $b$ have an overlap, then this overlap is a small piece.
\end{corollary}
\begin{proof}
The proof is trivial.
\end{proof}

\begin{corollary}
\label{max_occurrences_contains}
Let $U$ be a monomial. Let $a$ be a maximal occurrence of a monomial of $\Mon$ in $U$, $b$ be some occurrence of a monomial of $\Mon$ in $U$ (not necessarily maximal). Assume $a$ and $b$ have a non-empty common subword and this subword is not a small piece. Then $b$ is contained inside $a$ (not necessarily properly, $b = a$ is also possible).
\end{corollary}
\begin{proof}
Since $a$ and $b$ have a non-empty common subword and this subword is not a small piece, it follows from Lemma~\ref{one_max_occurrence_interactions} that this subword is not an overlap of $a$ and $b$. Therefore, either $a$ is properly contained inside $b$, or $b$ is properly contained inside $a$, or $a = b$. Since $a$ is a maximal occurrence, Lemma~\ref{max_can_not_be_inside} implies that the first variant is not possible. Thus, $b$ is contained inside $a$.
\end{proof}

\begin{corollary}
\label{max_occurrences_coinside}
Let $U$ be a monomial. Let $a$ and $b$ be two maximal occurrences of monomials of $\Mon$ in $U$. If $a$ and $b$ have a non-empty common subword in $U$ and this subword is not a small piece, then the occurrences $a$ and $b$ coincide.
\end{corollary}
\begin{proof}
It follows directly from Lemma~\ref{max_can_not_be_inside} and Corollary~\ref{max_occurrences_contains}.
\end{proof}

Throughout the paper we graphically represent monomials as segments and their subwords as subsegments. In particular, we represent maximal occurrences in a monomial as its subsegments.

Assume $a$ is a maximal occurrence of a monomial of $\Mon$ in $U$, $b$ is some occurrence of a monomial of $\Mon$ in $U$ (not necessarily maximal), $b$ is not contained in $a$ (in particular, $b \neq a$). To be definite, assume that $a$ starts from the left of the beginning of $b$. Lemma~\ref{max_can_not_be_inside} and Lemma~\ref{one_max_occurrence_interactions} imply that possible configurations of $a$ and $b$ are the following:
\begin{enumerate}
\item
\label{separated}
There exists some non empty subword between $a$ and $b$ in $U$.
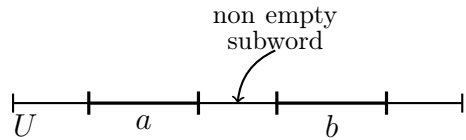
\begin{center}
\begin{tikzpicture}
\draw[|-|, black, thick] (0,0)--(6,0);
\node[below, xshift=5] at (0, 0) {$U$};
\draw[|-|, black, very thick] (1,0)--(2.5,0) node[midway, below] {$a$};
\draw[|-|, black, very thick] (3.5,0)--(5,0) node[midway, below] {$b$};
\node[text width=2cm, align=center] at (3.5,1) {\footnotesize{\baselineskip=10ptnon empty subword \par}};
\draw[->, black, thick] (3.5, 0.7) to [bend right] (3, 0);
\end{tikzpicture}
\end{center}
In this case we say that $a$ and $b$ are \emph{separated}.
\item
\label{touch}
$a$ and $b$ are adjacent and have no common non empty subword.
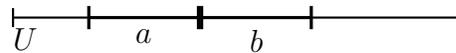
\begin{center}
\begin{tikzpicture}
\draw[|-|, black, thick] (0,0)--(6,0);
\node[below, xshift=5] at (0, 0) {$U$};
\draw[|-|, black, very thick] (1,0)--(2.5,0) node[midway, below] {$a$};
\draw[|-|, black, very thick] (2.5,0)--(4,0) node[midway, below] {$b$};
\end{tikzpicture}
\end{center}
In this case we say that $a$ and $b$ \emph{touch at a point}.
\item
\label{overlap}
$a$ and $b$ have a non empty common subword, wherein $b$ is not contained inside $a$.
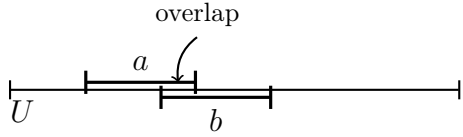
\begin{center}
\begin{tikzpicture}
\draw[|-|, black, thick] (0,0)--(6,0);
\node[below, xshift=5] at (0, 0) {$U$};
\draw[|-|, black, very thick] (1,0.1)--(2.5,0.1) node[midway, above] {$a$};
\draw[|-|, black, very thick] (2,-0.1)--(3.5,-0.1) node[midway, below] {$b$};
\node[text width=2cm, align=center] at (2.5,1) {\footnotesize{\baselineskip=10ptoverlap \par}};
\draw[->, black, thick] (2.5, 0.7) to [bend right] (2.25, 0.1);
\end{tikzpicture}
\end{center}
In this case we say that $a$ and $b$ \emph{have an overlap}. Then, by Lemma~\ref{one_max_occurrence_interactions}, the overlap is a small piece.
\end{enumerate}
So, one can see that the positions of maximal occurrences in $U$ are linearly ordered.

\begin{definition}
\label{chart_def}
Let $U$ be a monomial. We define \emph{the chart of $U$} as the set of all maximal occurrences of monomials of $\Mon$ in $U$. The maximal occurrences $m_i \in \Mon$ in $U$ such that $\SPM(m_i) \geqslant \tau$ are called \emph{members of the chart}.
\end{definition}

\begin{remark}
Assume $U$ is a monomial. Note that the notion of members of the chart of $U$ in Definition~\ref{chart_def} in fact depends on two parameters: the measure $\SPM$ and the constant $\tau$. So, formally they should be called the $(\SPM, \tau)$-members of the chart of~$U$. However, since $\SPM$ and $\tau$ are fixed throughout the paper, we omit these parameters and call them members of the chart of $U$.
\end{remark}

\medskip

Let $G$ be a small cancellation group, $R_i = M_1M_2^{-1}$ be a relator of its small cancellation presentation. Assume $LM_1R$ and $LM_2R$ are two words, then the transition from $LM_1R$ to $LM_2R$
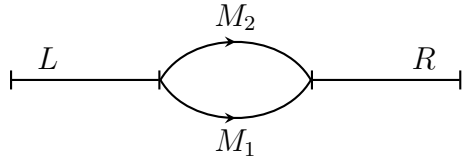
\begin{center}
\begin{tikzpicture}
\draw[|-|, black, thick] (0,0)--(2,0) node [near start, above] {$L$};
\draw[black, thick, arrow=0.5] (2,0) to [bend left=60] node [above] {$M_2$} (4,0);
\draw[black, thick, arrow=0.5] (2,0) to [bend right=60] node [below] {$M_1$} (4,0);
\draw[|-|, black, thick] (4,0)--(6,0) node [near end, above] {$R$};
\end{tikzpicture}
\end{center}
is called \emph{a turn} of an occurrence of the subrelation $M_1$ (to its complement $M_2$), see~\cite{NA1}. Analogously, in our case we define a multi-turn.
\begin{definition}
\label{multiturn_def}
Let $p = \sum_{j = 1}^{n} \alpha_{j} a_{j} \in \Rel$. For every $h = 1, \ldots, n$ we call the transition
\begin{equation*}
a_h\longmapsto \sum\limits_{\substack{j = 1 \\ j\neq h}}^{n} (-\alpha_{h}^{-1} \alpha_{j}a_{j})
\end{equation*}
an \emph{elementary multi-turn of $a_h$ with respect to $p$}.

Let $p = \sum_{j = 1}^{n} \alpha_{j} a_{j} \in \Rel$. Let $a_h$ be a maximal occurrence in $U$, $U = La_hR$. The transformation
\begin{equation}
\label{multi_turn_def}
U = La_hR \longmapsto \sum\limits_{\substack{j = 0 \\ j\neq h}}^{n} (-\alpha_{h}^{-1} \alpha_{j}L\cdot a_{j}\cdot R),
\end{equation}
with the further cancellations if there are any, is called \emph{a multi-turn of the maximal occurrence $a_h$ in $U$ that comes from an elementary multi-turn $a_h \mapsto \sum_{\substack{j = 1 \\ j\neq h}}^{n}(-\alpha_h^{-1}\alpha_j a_j)$}. Obviously,
\begin{equation*}
U - \sum\limits_{\substack{j = 0 \\ j\neq h}}^{n} (-\alpha_{h}^{-1} \alpha_{j}L\cdot a_{j}\cdot R) =  L\cdot (\alpha_{h}^{-1}p)\cdot R \in \Ideal.
\end{equation*}
In this case the polynomial $L\cdot p\cdot R = \sum_{j = 1}^n\alpha_j L\cdot a_j\cdot R$ (after the cancellations) is called \emph{a layout of the multi-turn~\eqref{multi_turn_def}}.
\medskip

We will show in Section~\ref{mt_configurations} that if $a_j \neq 1$, then the monomial $L\cdot a_{j}\cdot R$ has no cancellations. So, cancellations are possible only in the monomial $L\cdot R$ if $a_j = 1$ is a monomial of $p$. That is, in fact we can write the multi-turn above in the way
\begin{equation*}
U = La_hR \longmapsto \sum\limits_{\substack{j = 0 \\ j\neq h}}^{n} (-\alpha_{h}^{-1} \alpha_{j}La_{j}R),
\end{equation*}
and its layout $L\cdot p\cdot R$ in the form
\begin{equation*}
L\cdot p\cdot R = \sum_{j = 1}^n\alpha_j La_j R.
\end{equation*}
\end{definition}

\section{The description of the ideal $\Ideal$ as a linear subspace of $\Galg$}
\label{ideal_linear_descr_section}
Let us define a subspace of $\Galg$ of linear dependencies induced by multi-turns of members of the chart of monomials. For every monomial of $\Fr$ we do all multi-turns of all members of the chart and we consider all layouts of these multi-turns. As a result, we obtain the set of expressions
\begin{multline}
\label{linear_dep_members}
\GDp = \Bigg\lbrace \sum\limits_{j = 1}^n \alpha_jU_j \mid U_j \in \Fr, \textit{there exists an index }1\leqslant h \leqslant n \textit{ such that}\\
U_h \mapsto \sum_{\substack{j = 1 \\ j\neq h}}^{n} (-\alpha_h^{-1} \alpha_j U_j) \textit{ is a multi-turn of a member of the chart of } U_h\Bigg\rbrace.
\end{multline}
Assume $U$ runs through all monomials of $\Fr$ with non-empty charts, $a$ runs through all members of the chart of $U$. Clearly, for every fixed $U$ and $a$ we can write $U = L_{U, a}aR_{U, a}$, where $L_{U, a}$ is a prefix of $U$ and $R_{U, a}$ is a suffix of $U$. Then $\GDp$ consists of all polynomials of the form $\sum_{j = 1}^n \alpha_j L_{U, a} a_j R_{U, a}$ for different $U$ and $a$, where $\sum_{j = 1}^{n} \alpha_j a_j$ runs through all polynomials of $\Rel$ such that $a$ is a monomial in these polynomials. We denote the linear span of $\GDp$ by $\langle\GDp\rangle$.

\begin{proposition}
\label{the_ideal_characterisation}
The linear subspace $\langle\GDp \rangle \subseteq \Galg$ is equal to the ideal $\Ideal$.
\end{proposition}
\begin{proof}
First we will show that $\langle \GDp \rangle$ is an ideal of $\Galg$. Let $T \in \GDp$. We have to check that if $Z$ is a monomial, then $ZT \in \langle \GDp \rangle$ and $TZ \in \langle \GDp \rangle$. Clearly, it is sufficient to check this property only for a monomial $Z$ that consists of only one symbol of $S \cup S^{-1}$. In this particular case denote it by $z$. We will show an even stronger property, namely that $z\cdot T \in \GDp$ and $T\cdot z \in \GDp$.

Assume $U_h = La_hR$, $a_h$ is a member of the chart of $U_h$. Assume $T$ is a layout of a multi-turn of $a_h$ in $U_h$ that comes from an elementary multi-turn $a_h \mapsto \sum_{\substack{j = 1 \\ j\neq h}}^{n}\alpha_h^{-1}\alpha_ja_j$. That is,
\begin{equation*}
T = \sum\limits_{j = 1}^n \alpha_jU_j = \sum\limits_{j = 1}^n \alpha_j La_jR.
\end{equation*}
First consider the case when $L$ is not empty. Since $L$ is not empty, $a_h$ remains unchanged in $z\cdot U_h = (z\cdot L)a_hR$ both when $z$ does not cancel out with $L$, or when it does. Since $a_h$ is a member of the chart of $U_h$, $a_h$ is a maximal occurrence in $U_h$ and $\SPM(a_h) \geqslant \tau$. Then, clearly, $a_h$ is also a maximal occurrence in $z\cdot U_h$. Therefore, since $\SPM(a_h) \geqslant \tau$, $a_h$ is a member of the chart of $z\cdot U_h$. So, $z\cdot U_h = (z\cdot L)a_hR \mapsto \sum_{\substack{j = 1 \\ j\neq h}}^n \alpha_h^{-1}\alpha_j(z\cdot L)a_jR$ is a multi-turn of a member of the chart of $z\cdot U_h$. Clearly,
\begin{equation*}
z\cdot T = \sum\limits_{j = 1}^n \alpha_j z\cdot U_j = \sum\limits_{j = 1}^n \alpha_j(z\cdot L)a_jR
\end{equation*}
is a layout of this multi-turn. Hence, $z\cdot T\in \GDp$.

Now consider the case when $L$ is empty, that is, $U_j = a_jR$ and
\begin{equation*}
z\cdot T = \sum\limits_{j = 1}^n \alpha_j z\cdot a_jR.
\end{equation*}
First assume $z$ does not cancel with $a_h$ and $za_h \notin \Mon$. Clearly, in this case $a_h$ is a maximal occurrence in $zU_h = za_hR$. Since $\SPM(a_h) \geqslant \tau$, $a_h$ is a member of the chart of $za_hR$. Hence, $za_hR \mapsto \sum_{\substack{j = 1 \\ j\neq h}}^{n} \alpha_h^{-1}\alpha_j za_jR$ is a multi-turn of a member of the chart of $zU_h$. Since $z\cdot T$ is its layout, we have $z\cdot T \in \GDp$.

Assume $z$ does not cancel with $a_h$ and $za_h \in \Mon$. Then $za_h$ is a maximal occurrence of a monomial of $\Mon$ in $zU_h = za_hR$. It follows from Lemma~\ref{measure_properties} that $\SPM(za_h) \geqslant \SPM(a_h) \geqslant \tau$. Therefore, $za_h$ is a member of the chart of $zU_h$. Since $a_h$ is not a small piece, by Lemma~\ref{not_sp_prolongation}, we obtain $\sum_{j = 1}^n \alpha_j za_j \in \Rel$ (after the cancellations if there are any). Hence, $za_hR \mapsto \sum_{\substack{j\neq h \\ j = 1}}^n \alpha_h^{-1}\alpha_jza_jR$ is a multi-turn of a member of the chart of $zU_h$ and $z\cdot T$ is its layout. So, $z\cdot T \in \GDp$.

Assume $z$ cancels with $a_h$. Then, by Compatibility Axiom, we obtain $\sum_{j = 1}^n \alpha_j z\cdot a_j \in \Rel$ (after the cancellations). Since $\SPM(z\cdot a_h) \leqslant \SPM(a_h)$, we distinguish two possibilities: $\SPM(z\cdot a_h) \geqslant \tau$ and $\SPM(z\cdot a_h) < \tau$.

Consider the case $\SPM(z\cdot a_h) \geqslant \tau$. Then $z\cdot a_h$ is a member of the chart of $z\cdot U_h$ after the cancellations. So, $z\cdot a_hR \mapsto \sum_{\substack{j = 1 \\ j\neq h}}^n \alpha_h^{-1}\alpha_j z\cdot a_jR$ is a multi-turn of a member of the chart of $z\cdot U_h$ and $z\cdot T$ is its layout. Hence, $z\cdot T \in \GDp$.

Now consider the case $\SPM(z\cdot a_h) < \tau$. In this case $z\cdot a_h$ is not a member of the chart of $z\cdot U_h$. However, from Small Cancellation Axiom it follows that in the polynomial $\sum_{j = 1}^n \alpha_j z\cdot a_j$ there exists a monomial $z\cdot a_{h_0}$ such that $\SPM(z\cdot a_{h_0}) \geqslant \tau + 1$ (possibly after the cancellations). Then $z\cdot a_{h_0}$ is a member of the chart of $z\cdot a_{h_0}R$ (possibly after the cancellations). Hence, $z\cdot a_{h_0}R \mapsto \sum_{\substack{j = 1 \\ j\neq h_0}}^n \alpha_{h_0}^{-1}\alpha_j z\cdot a_jR$ is a multi-turn of a member of the chart of $z\cdot a_{h_0}R$. Clearly, a layout of this multi-turn is also $z\cdot T$. Therefore, $z\cdot T \in \GDp$.

Summarising all of the above, we obtain $z\cdot T \in \GDp$ in all cases. Clearly, for the same reason we obtain $T\cdot z \in \GDp$. Hence, $\langle \GDp\rangle$ is an ideal of $\Galg$.

Consider a polynomial $\sum_{j = 1}^n \alpha_j U_j = \sum_{j = 1}^n \alpha_j La_jR \in \GDp$ that comes from an elementary multi-turn $a_h \mapsto \sum_{\substack{j = 1 \\ j\neq h}}^n\alpha_{h}^{-1}\alpha_j a_j$. Since $\Rel \subseteq \Ideal$, we have $\sum_{j = 1}^n\alpha_j a_j \in \Ideal$. Hence,
\begin{equation*}
\sum_{j = 1}^n \alpha_j U_j = L\cdot (\sum_{j = 1}^n \alpha_j a_j)\cdot R \in \Ideal.
\end{equation*}
So, $\GDp \subseteq \Ideal$ and $\langle\GDp\rangle \subseteq \Ideal$.

Assume $p = \sum_{j = 1}^n \alpha_ja_j \in \Rel$. By Small Cancellation Axiom, in this polynomial there exists a monomial $a_{h_0}$ such that $\SPM(a_{h_0}) \geqslant \tau + 1$. Hence, $a_{h_0}$ is the single member of the chart of the monomial $a_{h_0}$. Then $a_{h_0} \mapsto \sum_{\substack{j = 1 \\ j\neq h_0}}^n\alpha_{h_0}^{-1}\alpha_j a_j$ is a multi-turn of the member of the chart of $a_{h_0}$. Clearly, the polynomial $p$ is a layout of this multi-turn, so, $p \in \GDp$. Therefore, $\Rel \subseteq \GDp$. So, the ideal of $\Galg$ generated by $\Rel$ is contained in the ideal of $\Galg$ generated by $\GDp$. Since $\langle\GDp\rangle$ is an ideal of $\Galg$, we get that the ideal of $\Galg$ generated by $\GDp$ is equal to $\langle\GDp\rangle$. By definition, $\Ideal = \langle \Rel\rangle_i$. Combining these facts, we obtain $\Ideal \subseteq \langle\GDp\rangle$. Thus, $\Ideal =\langle\GDp\rangle$.
\end{proof}

\section{How multi-turns influence maximal occurrences}
\label{mt_configurations}
Let $U_h$ be a monomial, $a_h \in \Mon$ be a maximal occurrence in $U_h$, $U_h = La_hR$. We suppose that $a_h$ is not too short, namely, that $\SPM(a_h) \geqslant 3$. Assume $U_h \mapsto \sum_{\substack{j = 1 \\ j\neq h}}^n (-\alpha_h^{-1}\alpha_j U_j)$ is a multi-turn that comes from an elementary multi-turn $a_h \mapsto \sum_{\substack{j = 1 \\ j\neq h}}^n (-\alpha_h^{-1}\alpha_j a_j)$, so, $U_j = La_jR$. It is intuitively clear that for every maximal occurrence in $U_h$ there exists a corresponding maximal occurrence in every $U_j$. In this section we describe precisely how that corresponding maximal occurrences look like.

\begin{definition}
\label{intersection_def}
Let $U_h$ be a monomial, $b \in \Mon$ be a maximal occurrence in $U_h$, $A$ be an occurrence in $U_h$ ($A$ does not necessarily belong to $\Mon$). We define \emph{the intersection of $b$ and $A$ in $U_h$} as the occurrence that consists of letters that are contained both in $b$ and in $A$.
\end{definition}
Evidently, there are the following five possibilities:
\begin{enumerate}
\item
$b$ and $A$ have the empty intersection;
\begin{center}
\begin{tikzpicture}
\draw[|-|, black, thick] (0,0)--(6,0);
\node[below, xshift=10] at (0, 0) {$U_h$};
\draw[|-|, black, very thick] (1,0)--(2.5,0) node[midway, above] {$A$};
\draw[|-|, black, very thick] (3.5,0)--(5,0) node[midway, below] {$b$};
\end{tikzpicture}

\begin{tikzpicture}
\draw[|-|, black, thick] (0,0)--(6,0);
\node[below, xshift=10] at (0, 0) {$U_h$};
\draw[|-|, black, very thick] (1,0)--(2.5,0) node[midway, above] {$b$};
\draw[|-|, black, very thick] (3.5,0)--(5,0) node[midway, below] {$A$};
\end{tikzpicture}

\begin{tikzpicture}
\draw[|-|, black, thick] (0,0)--(6,0);
\node[below, xshift=10] at (0, 0) {$U_h$};
\draw[|-|, black, very thick] (1.5,0)--(3,0) node[midway, above] {$A$};
\draw[|-|, black, very thick] (3,0)--(4.5,0) node[midway, below] {$b$};
\end{tikzpicture}

\begin{tikzpicture}
\draw[|-|, black, thick] (0,0)--(6,0);
\node[below, xshift=10] at (0, 0) {$U_h$};
\draw[|-|, black, very thick] (1.5,0)--(3,0) node[midway, below] {$b$};
\draw[|-|, black, very thick] (3,0)--(4.5,0) node[midway, above] {$A$};
\end{tikzpicture}
\end{center}
\item
the intersection of $b$ and $A$ is equal to $b$, that is, $b$ is contained in $A$;
\begin{center}
\begin{tikzpicture}
\draw[|-|, black, thick] (0,0)--(6,0);
\node[below, xshift=10] at (0, 0) {$U_h$};
\draw[|-|, black, very thick] (2,0)--(5,0) node[near start, above] {$A$};
\draw[|-|, black, very thick] (3,-0.1)--(4.5,-0.1) node[midway, below] {$b$};
\end{tikzpicture}
\end{center}
\item
the intersection of $b$ and $A$ is equal to $A$, that is, $A$ is contained in $b$;
\begin{center}
\begin{tikzpicture}
\draw[|-|, black, thick] (0,0)--(6,0);
\node[below, xshift=10] at (0, 0) {$U_h$};
\draw[|-|, black, very thick] (2,0)--(5,0) node[near start, below] {$b$};
\draw[|-|, black, very thick] (3,0.1)--(4.5,0.1) node[midway, above] {$A$};
\end{tikzpicture}
\end{center}
\item
the intersection of $b$ and $A$ is equal to some proper prefix of $b$;
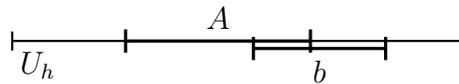
\begin{center}
\begin{tikzpicture}
\draw[|-|, black, thick] (0,0)--(6,0);
\node[below, xshift=10] at (0, 0) {$U_h$};
\draw[|-|, black, very thick] (1.5,0)--(4,0) node[midway, above] {$A$};
\draw[|-|, black, very thick] (3.2,-0.1)--(5,-0.1) node[midway, below] {$b$};
\end{tikzpicture}
\end{center}
\item
the intersection of $b$ and $A$ is equal to some proper suffix of $b$.
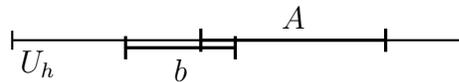
\begin{center}
\begin{tikzpicture}
\draw[|-|, black, thick] (0,0)--(6,0);
\node[below, xshift=10] at (0, 0) {$U_h$};
\draw[|-|, black, very thick] (2.5,0)--(5,0) node[midway, above] {$A$};
\draw[|-|, black, very thick] (1.5,-0.1)--(3,-0.1) node[midway, below] {$b$};
\end{tikzpicture}
\end{center}
\end{enumerate}

There are the following four possibilities for $a_h$:
\begin{itemize}
\item
there exist maximal occurrences in $U_h$ that begin from the left of the beginning of $a_h$ and that begin from the right of the beginning of $a_h$;
\item
there are no maximal occurrences in $U_h$ that begin from the left of the beginning of $a_h$;
\item
there are no maximal occurrences in $U_h$ that begin from the right of the beginning of $a_h$;
\item
$a_h$ is a single maximal occurrence in $U_h$.
\end{itemize}
In this section we consider only the first possibility as the most interesting. Other cases are treated in the similar (but simpler) way.

Let us fix the notations. Throughout this section we assume that $b$ and $c$ are two maximal occurrences in $U_h$ such that $b$ starts from the left of the beginning of $a_h$ and $c$ starts from the right of the beginning of $a_h$. Since $\SPM(a_h) \geqslant 3$, $b$ and $c$ are separated in $U_h$.
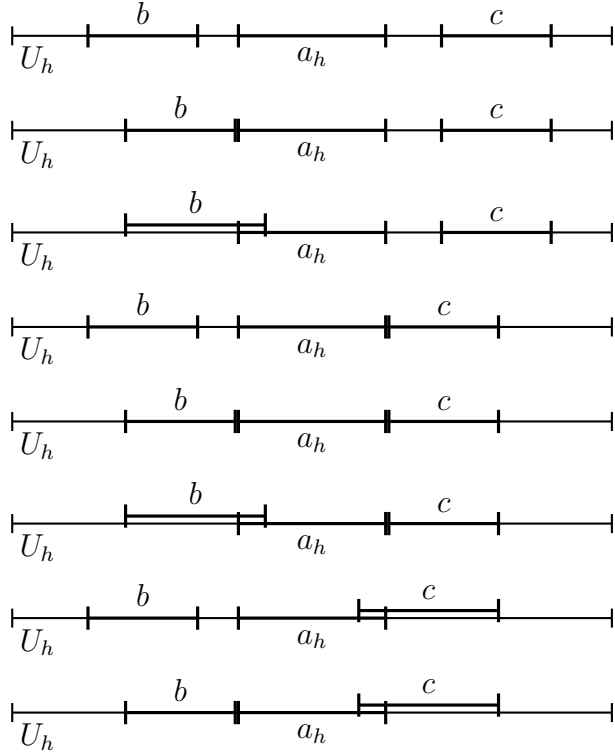
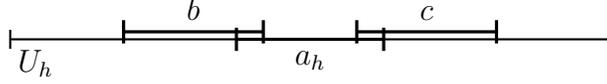
\begin{center}
\begin{tikzpicture}
\draw[|-|, black, thick] (0,0)--(8,0);
\node[below, xshift=10] at (0, 0) {$U_h$};
\draw[|-|, black, very thick] (3,0)--(5,0) node[midway, below] {$a_h$};
\draw[|-|, black, very thick] (1,0)--(2.5,0) node[midway, above] {$b$};
\draw[|-|, black, very thick] (5.7,0)--(7.2,0) node[midway, above] {$c$};
\end{tikzpicture}

\begin{tikzpicture}
\draw[|-|, black, thick] (0,0)--(8,0);
\node[below, xshift=10] at (0, 0) {$U_h$};
\draw[|-|, black, very thick] (3,0)--(5,0) node[midway, below] {$a_h$};
\draw[|-|, black, very thick] (1.5,0)--(3,0) node[midway, above] {$b$};
\draw[|-|, black, very thick] (5.7,0)--(7.2,0) node[midway, above] {$c$};
\end{tikzpicture}

\begin{tikzpicture}
\draw[|-|, black, thick] (0,0)--(8,0);
\node[below, xshift=10] at (0, 0) {$U_h$};
\draw[|-|, black, very thick] (3,0)--(5,0) node[midway, below] {$a_h$};
\draw[|-|, black, very thick] (1.5,0.1)--(3.4,0.1) node[midway, above] {$b$};
\draw[|-|, black, very thick] (5.7,0)--(7.2,0) node[midway, above] {$c$};
\end{tikzpicture}

\begin{tikzpicture}
\draw[|-|, black, thick] (0,0)--(8,0);
\node[below, xshift=10] at (0, 0) {$U_h$};
\draw[|-|, black, very thick] (3,0)--(5,0) node[midway, below] {$a_h$};
\draw[|-|, black, very thick] (1,0)--(2.5,0) node[midway, above] {$b$};
\draw[|-|, black, very thick] (5,0)--(6.5,0) node[midway, above] {$c$};
\end{tikzpicture}

\begin{tikzpicture}
\draw[|-|, black, thick] (0,0)--(8,0);
\node[below, xshift=10] at (0, 0) {$U_h$};
\draw[|-|, black, very thick] (3,0)--(5,0) node[midway, below] {$a_h$};
\draw[|-|, black, very thick] (1.5,0)--(3,0) node[midway, above] {$b$};
\draw[|-|, black, very thick] (5,0)--(6.5,0) node[midway, above] {$c$};
\end{tikzpicture}

\begin{tikzpicture}
\draw[|-|, black, thick] (0,0)--(8,0);
\node[below, xshift=10] at (0, 0) {$U_h$};
\draw[|-|, black, very thick] (3,0)--(5,0) node[midway, below] {$a_h$};
\draw[|-|, black, very thick] (1.5,0.1)--(3.4,0.1) node[midway, above] {$b$};
\draw[|-|, black, very thick] (5,0)--(6.5,0) node[midway, above] {$c$};
\end{tikzpicture}

\begin{tikzpicture}
\draw[|-|, black, thick] (0,0)--(8,0);
\node[below, xshift=10] at (0, 0) {$U_h$};
\draw[|-|, black, very thick] (3,0)--(5,0) node[midway, below] {$a_h$};
\draw[|-|, black, very thick] (1,0)--(2.5,0) node[midway, above] {$b$};
\draw[|-|, black, very thick] (4.6,0.1)--(6.5,0.1) node[midway, above] {$c$};
\end{tikzpicture}

\begin{tikzpicture}
\draw[|-|, black, thick] (0,0)--(8,0);
\node[below, xshift=10] at (0, 0) {$U_h$};
\draw[|-|, black, very thick] (3,0)--(5,0) node[midway, below] {$a_h$};
\draw[|-|, black, very thick] (1.5,0)--(3,0) node[midway, above] {$b$};
\draw[|-|, black, very thick] (4.6,0.1)--(6.5,0.1) node[midway, above] {$c$};
\end{tikzpicture}

\begin{tikzpicture}
\draw[|-|, black, thick] (0,0)--(8,0);
\node[below, xshift=10] at (0, 0) {$U_h$};
\draw[|-|, black, very thick] (3,0)--(5,0) node[midway, below] {$a_h$};
\draw[|-|, black, very thick] (1.5,0.1)--(3.4,0.1) node[midway, above] {$b$};
\draw[|-|, black, very thick] (4.6,0.1)--(6.5,0.1) node[midway, above] {$c$};
\end{tikzpicture}
\end{center}

We distinguish three types of monomials in the sum
\begin{equation*}
\sum_{\substack{j = 1\\ j\neq h}}^n(-\alpha_h^{-1}\alpha_j U_j) = \sum_{\substack{j = 1\\ j\neq h}}^n(-\alpha_h^{-1}\alpha_j La_jR):
\end{equation*}
\begin{enumerate}
\item
\label{keep_structure}
$La_jR$, $j \neq h$, such that $a_j$ is not a small piece;
\item
\label{donot_keep_structure1}
$La_jR$, $j \neq h$, such that $a_j$ is a small piece but $a_j \neq 1$;
\item
\label{donot_keep_structure2}
$La_jR$, $j \neq h$, such that $a_j = 1$.
\end{enumerate}

Let $U_j = La_jR$ be a monomial of type~\ref{keep_structure} or~\ref{donot_keep_structure1}. First of all notice that the monomial $U_j = La_jR$ is reduced. Indeed, assume $U_j = L\cdot a_j\cdot R$ is not reduced. Since $L$ and $R$ are reduced monomials and $a_j \neq 1$, we obtain that at least one of $L\cdot a_j$ and $a_j\cdot R$ is not reduced. Assume $L\cdot a_j$ is not reduced, so, $L = L^{\prime}X$, $a_j = X^{-1}a_j^{\prime}$, $L\cdot a_j = \left(L^{\prime}X\right)\cdot \left(X^{-1}a_j^{\prime}\right)$, where $L^{\prime}$, $a_j^{\prime}$ can be empty.
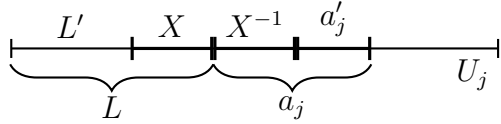
\begin{center}
\begin{tikzpicture}
\draw[|-|, black, thick] (0,0)--(6.5,0);
\node[below, xshift=-10] at (6.5, 0) {$U_j$};
\draw [thick, decorate, decoration={brace, amplitude=10pt, raise=6pt, mirror}] (0,0) to node[midway, below, yshift=-14pt] {$L$} (2.7, 0);
\draw [thick, decorate, decoration={brace, amplitude=10pt, raise=6pt, mirror}] (2.7,0) to node[midway, below, yshift=-14pt] {$a_j$} (4.8, 0);
\draw[|-|, black, very thick] (3.8,0)--(4.8,0) node[midway, above] {$a_j^{\prime}$};
\draw[|-|, black, very thick] (2.7,0)--(3.8,0) node[midway, above] {$X^{-1}$};
\draw[|-|, black, very thick] (1.6,0)--(2.7,0) node[midway, above] {$X$};
\path (0,0)--(1.6,0) node[midway, above] {$L^{\prime}$};
\end{tikzpicture}
\end{center}
Assume $p$ is a layout of the elementary multi-turn $a_h \mapsto \sum_{\substack{j = 1 \\ j\neq h}}^n (-\alpha_h^{-1}\alpha_j a_j)$. Then $p \in \Rel$ and $a_h$ and $a_j$ are monomials in $p$. Since $X^{-1}$ is a prefix of $a_j$, it follows from Compatibility Axiom that $X\cdot p \in \Rel$. Therefore, $Xa_h \in \Mon$. Since $U_h$ is of the form $U_h = L^{\prime}Xa_hR$, $a_h$ is not a maximal occurrence in $U_h$, a contradiction.
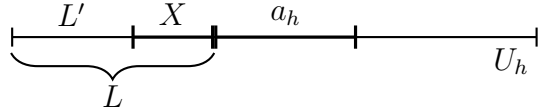
\begin{center}
\begin{tikzpicture}
\draw[|-|, black, thick] (0,0)--(7,0);
\node[below, xshift=-10] at (7, 0) {$U_h$};
\draw [thick, decorate, decoration={brace, amplitude=10pt, raise=6pt, mirror}] (0,0) to node[midway, below, yshift=-14pt] {$L$} (2.7, 0);
\draw[|-|, black, very thick] (2.7,0)--(4.6,0) node[midway, above] {$a_h$};
\draw[|-|, black, very thick] (1.6,0)--(2.7,0) node[midway, above] {$X$};
\path (0,0)--(1.6,0) node[midway, above] {$L^{\prime}$};
\end{tikzpicture}
\end{center}
The case when $a_j\cdot R$ is not reduced is considered similarly.

Let $\widehat{b}$ be the intersection of $L$ and $b$ in $U_h$, $\widehat{c}$ be the intersection of $R$ and $c$ in $U_h$. Notice that since $b$ begins from the left of the beginning of $a_h$, $\widehat{b}$ is not empty. Similarly, since $c$ begins from the right of the beginning of $a_h$, $\widehat{c}$ is not empty. Then $\widehat{b}$ and $\widehat{c}$ can be considered as occurrences in $U_j$ in a natural way.
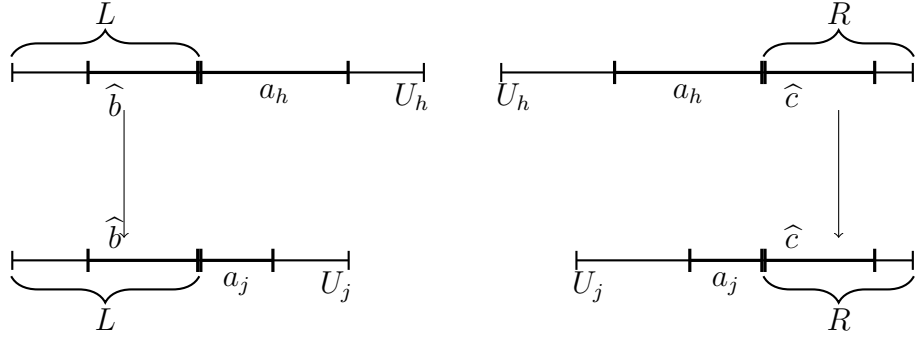
\begin{center}
\begin{tikzpicture}
\draw[|-|, black, thick] (0,0)--(5.5,0);
\node[below, xshift=-5] at (5.5, 0) {$U_h$};
\draw [thick, decorate, decoration={brace, amplitude=10pt, raise=6pt}] (0,0) to node[midway, above, yshift=14pt] {$L$} (2.5, 0);
\draw[|-|, black, very thick] (2.5,0)--(4.5,0) node[midway, below] {$a_h$};
\draw[|-|, black, very thick] (1,0)--(2.5,0) node[near start, below] {$\widehat{b}$};

\draw[|-|, black, thick] (0,-2.5)--(4.5,-2.5);
\node[below, xshift=-5] at (4.5, -2.5) {$U_j$};
\draw [thick, decorate, decoration={brace, amplitude=10pt, raise=6pt, mirror}] (0,-2.5) to node[midway, below, yshift=-14pt] {$L$} (2.5, -2.5);
\draw[|-|, black, very thick] (2.5,-2.5)--(3.5,-2.5) node[midway, below] {$a_j$};
\draw[|-|, black, very thick] (1,-2.5)--(2.5,-2.5) node[near start, above] {$\widehat{b}$};

\draw[->, black] (1.5, -0.5)--(1.5,-2.2);

\draw[|-|, black, thick] (6.5,0)--(12,0);
\node[below, xshift=5] at (0 + 6.5, 0) {$U_h$};
\draw [thick, decorate, decoration={brace, amplitude=10pt, raise=6pt}] (10, 0) to node[midway, above, yshift=14pt] {$R$} (12, 0);
\draw[|-|, black, very thick] (8,0)--(10,0) node[midway, below] {$a_h$};
\draw[|-|, black, very thick] (10,0)--(11.5,0) node[near start, below] {$\widehat{c}$};

\draw[|-|, black, thick] (7.5,-2.5)--(12,-2.5);
\node[below, xshift=5] at (0 + 7.5, -2.5) {$U_j$};
\draw [thick, decorate, decoration={brace, amplitude=10pt, raise=6pt, mirror}] (10, -2.5) to node[midway, below, yshift=-14pt] {$R$} (12, -2.5);
\draw[|-|, black, very thick] (9,-2.5)--(10,-2.5) node[midway, below] {$a_j$};
\draw[|-|, black, very thick] (10,-2.5)--(11.5,-2.5) node[near start, above] {$\widehat{c}$};

\draw[->, black] (11, -0.5)--(11,-2.2);
\end{tikzpicture}
\end{center}
The occurrences $\widehat{b}$ and $\widehat{c}$ do not have to be maximal in $U_j$. We study only $b$ and $c$ such that $\widehat{b}$ and $\widehat{c}$ are not small pieces. Then, by Corollary~\ref{max_occurrences_coinside}, there exists a unique maximal occurrence in $U_j$ that contains $\widehat{b}$ and there exists a unique maximal occurrence in $U_j$ that contains $\widehat{c}$. We denote them by $b^{\prime}$ and $c^{\prime}$ respectively. So, it is natural to consider $b$ in $U_h$ and $b^{\prime}$ in $U_j$ as maximal occurrences that correspond to each other and $c$ in $U_h$ and $c^{\prime}$ in $U_j$ as maximal occurrences that correspond to each other.

We study how $b^{\prime}$ and $c^{\prime}$ look like and how they are mutually arranged. Assume $b$ and $a_h$ are separated in $U_h$. Then $b$ is contained inside $L$, so, $\widehat{b} = b$ in $U_h$ and $\widehat{b}$ is separated from $a_j$ in $U_j$. Clearly, in this case $\widehat{b}$ itself is a maximal occurrence in $U_j$, that is, $b^{\prime} = \widehat{b}$ in $U_j$. Informally speaking, we can say that $b$ stays unchanged in $U_j$ in this case.
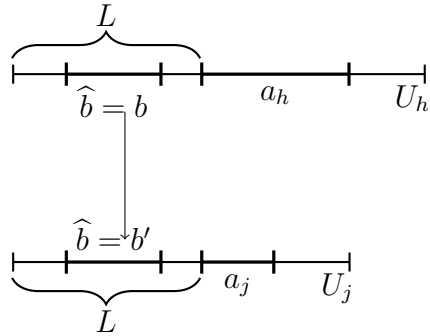
\begin{center}
\begin{tikzpicture}
\draw[|-|, black, thick] (0,0)--(5.5,0);
\node[below, xshift=-5] at (5.5, 0) {$U_h$};
\draw [thick, decorate, decoration={brace, amplitude=10pt, raise=6pt}] (0,0) to node[midway, above, yshift=14pt] {$L$} (2.5, 0);
\draw[|-|, black, very thick] (2.5,0)--(4.5,0) node[midway, below] {$a_h$};
\draw[|-|, black, very thick] (0.7,0)--(2,0) node[midway, below] {$\widehat{b} = b$};

\draw[|-|, black, thick] (0,-2.5)--(4.5,-2.5);
\node[below, xshift=-5] at (4.5, -2.5) {$U_j$};
\draw [thick, decorate, decoration={brace, amplitude=10pt, raise=6pt, mirror}] (0,-2.5) to node[midway, below, yshift=-14pt] {$L$} (2.5, -2.5);
\draw[|-|, black, very thick] (2.5,-2.5)--(3.5,-2.5) node[midway, below] {$a_j$};
\draw[|-|, black, very thick] (0.7,-2.5)--(2,-2.5) node[midway, above] {$\widehat{b} = b^{\prime}$};

\draw[->, black] (1.5, -0.5)--(1.5,-2.2);
\end{tikzpicture}
\end{center}
The same for $c$. Namely, if $c$ and $a_h$ are separated in $U_h$, then $c$ is contained inside $R$, $\widehat{c} = c$ in $U_h$ and $\widehat{c}$ is separated from $a_j$ in $U_j$. Then $\widehat{c}$ itself is a maximal occurrence in $U_j$, that is, $c^{\prime} = \widehat{c}$ in $U_j$.
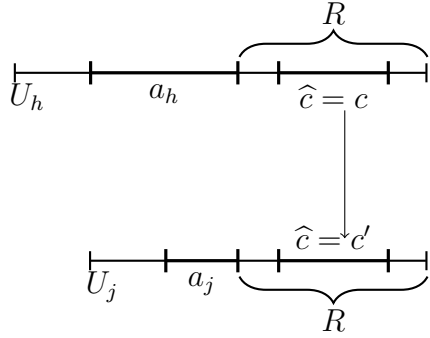
\begin{center}
\begin{tikzpicture}
\draw[|-|, black, thick] (0,0)--(5.5,0);
\node[below, xshift=5] at (0, 0) {$U_h$};
\draw [thick, decorate, decoration={brace, amplitude=10pt, raise=6pt}] (3, 0) to node[midway, above, yshift=14pt] {$R$} (5.5, 0);
\draw[|-|, black, very thick] (1,0)--(3,0) node[midway, below] {$a_h$};
\draw[|-|, black, very thick] (3.5,0)--(5,0) node[midway, below] {$\widehat{c} = c$};

\draw[|-|, black, thick] (1,-2.5)--(5.5,-2.5);
\node[below, xshift=5] at (1, -2.5) {$U_j$};
\draw [thick, decorate, decoration={brace, amplitude=10pt, raise=6pt, mirror}] (3, -2.5) to node[midway, below, yshift=-14pt] {$R$} (5.5, -2.5);
\draw[|-|, black, very thick] (2,-2.5)--(3,-2.5) node[midway, below] {$a_j$};
\draw[|-|, black, very thick] (3.5,-2.5)--(5,-2.5) node[midway, above] {$\widehat{c} = c^{\prime}$};

\draw[->, black] (4.4, -0.5)--(4.4,-2.2);
\end{tikzpicture}
\end{center}

Assume $b$ is not separated from $a_h$ in $U_h$. Then $\widehat{b}$ is a suffix of $L$, $L = L_1\widehat{b}$. Let us show that $b^{\prime}$ is of the form $b^{\prime} = \widehat{b}Y$, where $Y$ is a suffix of $b^{\prime}$. Indeed, assume $b^{\prime} = X\widehat{b}Y$, where $X$ is a non-empty prefix of $b^{\prime}$. First suppose $b$ and $a_h$ touch at a point. Then $\widehat{b} = b$ in $U_h$ and $b^{\prime} = XbY$.
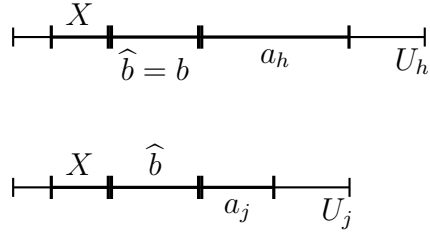
\begin{center}
\begin{tikzpicture}
\draw[|-|, black, thick] (0,0)--(5.5,0);
\node[below, xshift=-5] at (5.5, 0) {$U_h$};
\draw[|-|, black, very thick] (0.5,0)--(1.3,0) node[midway, above] {$X$};
\draw[|-|, black, very thick] (2.5,0)--(4.5,0) node[midway, below] {$a_h$};
\draw[|-|, black, very thick] (1.3,0)--(2.5,0) node[midway, below] {$\widehat{b} = b$};

\draw[|-|, black, thick] (0,-2)--(4.5,-2);
\node[below, xshift=-5] at (4.5, -2) {$U_j$};
\draw[|-|, black, very thick] (2.5,-2)--(3.5,-2) node[midway, below] {$a_j$};
\draw[|-|, black, very thick] (0.5,-2)--(1.3,-2) node[midway, above] {$X$};
\draw[|-|, black, very thick] (1.3,-2)--(2.5,-2) node[midway, above] {$\widehat{b}$};
\end{tikzpicture}
\end{center}
Hence, $Xb \in \Mon$, so, $b$ is not a maximal occurrence in $U_h$. Contradiction. Now let $b$ and $a_h$ have an overlap $d$, then $b = \widehat{b}d$ and $b^{\prime} = X\widehat{b}Y$.
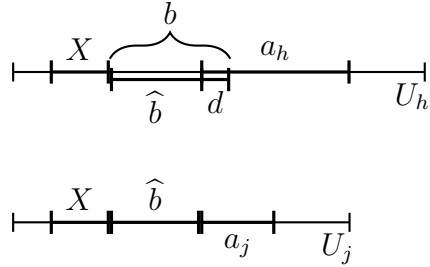
\begin{center}
\begin{tikzpicture}
\draw[|-|, black, thick] (0,0)--(5.5,0);
\node[below, xshift=-5] at (5.5, 0) {$U_h$};
\draw[|-|, black, very thick] (0.5,0)--(1.3,0) node[midway, above] {$X$};
\draw[|-|, black, very thick] (2.5,0)--(4.5,0) node[midway, above] {$a_h$};
\draw[|-, black, very thick] (1.3,-0.1)--(2.5,-0.1) node[midway, below] {$\widehat{b}$};
\draw[|-|, black, very thick] (2.5,-0.1)--(2.9,-0.1) node[midway, below] {$d$};
\draw [thick, decorate, decoration={brace, amplitude=10pt, raise=5pt}] (1.3, 0) to node[midway, above, yshift=14pt] {$b$} (2.9, 0);

\draw[|-|, black, thick] (0,-2)--(4.5,-2);
\node[below, xshift=-5] at (4.5, -2) {$U_j$};
\draw[|-|, black, very thick] (2.5,-2)--(3.5,-2) node[midway, below] {$a_j$};
\draw[|-|, black, very thick] (0.5,-2)--(1.3,-2) node[midway, above] {$X$};
\draw[|-|, black, very thick] (1.3,-2)--(2.5,-2) node[midway, above] {$\widehat{b}$};
\end{tikzpicture}
\end{center}
So we have $X\widehat{b} \in \Mon$ and $\widehat{b}d \in \Mon$. Therefore, since $\widehat{b}$ is not a small piece, Lemma~\ref{not_sp_prolongation} implies that $Xb = X\widehat{b}d \in \Mon$. This again contradicts with the assumption that $b$ is a maximal occurrence in $U_h$. Assume $c$ is not separated from $a_h$ in $U_h$. Then $\widehat{c}$ is a prefix of $R$, $R = \widehat{c}R_1$. By the same argument, we obtain $c^{\prime} = Z\widehat{c}$, where $Z$ is a prefix of $c^{\prime}$.

\begin{remark}
On the other hand, if $\widehat{b}$ ($\widehat{c}$) is a small piece, then in general it can be contained in several maximal occurrences in $U_j$. Since we do not need this for the further  argument, we do not consider the case when $\widehat{b}$ ($\widehat{c}$) is a small piece in detail in this section and state here this remark only for general information.

The occurrence $b$ ($c$) can be considered as a maximal occurrence in $U_h$ with respect to some fixed polynomial $p = \sum_{j = 1}^n \beta_j b_j \in \Rel$, $b = b_h$, $1 \leqslant h \leqslant n$ ($q = \sum_{j = 1}^k \gamma_j c_j \in \Rel$, $c = c_h$, $1 \leqslant h \leqslant k$). Then, in fact, the situation is the same as above. Namely, if we consider prolongations of $\widehat{b}$ ($\widehat{c}$) only with respect to $p$ ($q$), then we obtain that $\widehat{b}$ ($\widehat{c}$) is contained in a unique maximal occurrence in $U_j$ with respect to $p$ ($q$) and this occurrence is of the form $bY$ ($Zc$). However, in order to show this precisely, we need to use some additional considerations. Since we do not use this for the further  argument, we do not prove this statement here.
\end{remark}

Let $m_1$ and $m_2$ be maximal occurrences in $U_h$ such that the intersection of $m_1$ and $L$ is not a small piece and the the intersection of $m_2$ and $L$ is not a small piece. Let $\widehat{m}_1$ be the intersection of $L$ and $m_1$, $\widehat{m}_2$ be the intersection of $L$ and $m_2$. Assume both $m_1$ and $m_2$ are not separated from $a_h$ in $U_h$. Then, clearly, either $\widehat{m}_1$ is contained in $\widehat{m}_2$, or $\widehat{m}_2$ is contained in $\widehat{m}_1$.
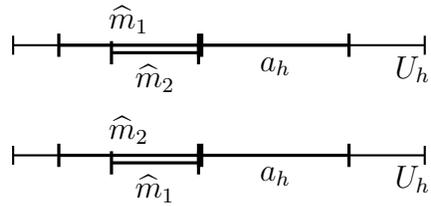
\begin{center}
\begin{tikzpicture}
\draw[|-|, black, thick] (0,0)--(5.5,0);
\node[below, xshift=-5] at (5.5, 0) {$U_h$};
\draw[|-|, black, very thick] (2.5,0)--(4.5,0) node[midway, below] {$a_h$};
\draw[|-|, black, very thick] (0.6,0)--(2.5,0) node[midway, above] {$\widehat{m}_1$};
\draw[|-|, black, very thick] (1.3,-0.1)--(2.5,-0.1) node[midway, below] {$\widehat{m}_2$};
\end{tikzpicture}

\begin{tikzpicture}
\draw[|-|, black, thick] (0,0)--(5.5,0);
\node[below, xshift=-5] at (5.5, 0) {$U_h$};
\draw[|-|, black, very thick] (2.5,0)--(4.5,0) node[midway, below] {$a_h$};
\draw[|-|, black, very thick] (0.6,0)--(2.5,0) node[midway, above] {$\widehat{m}_2$};
\draw[|-|, black, very thick] (1.3,-0.1)--(2.5,-0.1) node[midway, below] {$\widehat{m}_1$};
\end{tikzpicture}
\end{center}
Since $\widehat{m}_1$ and $\widehat{m}_2$ are not small pieces, by Corollary~\ref{max_occurrences_coinside}, we obtain $m_1 = m_2$. Hence, there exists not more than one maximal occurrence in $U_h$ such that it is not separated from $a_h$ from the left side and its intersection with $L$ is not a small piece. By the same argument, there exists not more than one maximal occurrence in $U_h$ such that it is not separated from $a_h$ from the right side and its intersection with $R$ is not a small piece.

We already showed what happens with maximal occurrences that are separated from $a_h$ in $U_h$. So we have to study only what happens with maximal occurrences that are not separated from $a_h$ in $U_h$. There are the following possibilities.
\begin{enumerate}[label=(\alph*)]
\item
\label{no_not_separated}
There does not exist a maximal occurrence in $U_h$ such that it is not separated from $a_h$ from the left side and such that its intersection with $L$ is not a small piece. And similarly, there does not exist a maximal occurrence in $U_h$ such that it is not separated from $a_h$ from the right side and such that its intersection with $R$ is not a small piece.
\item
\label{not_separated_left}
There exists a maximal occurrence in $U_h$ such that it is not separated from $a_h$ from the left side and such that its intersection with $L$ is not a small piece. Let this be our $b$. As before, $\widehat{b}$ is its intersection with $L$. But there does not exist a maximal occurrence in $U_h$ such that it is not separated from $a_h$ from the right side and such that its intersection with $R$ is not a small piece.
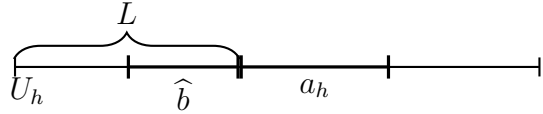
\begin{center}
\begin{tikzpicture}
\draw[|-|, black, thick] (0,0)--(7,0);
\node[below, xshift=5] at (0, 0) {$U_h$};
\draw[|-|, black, very thick] (3,0)--(5,0) node[midway, below] {$a_h$};
\draw[|-|, black, very thick] (1.5,0)--(3,0) node[midway, below] {$\widehat{b}$};
\draw [thick, decorate, decoration={brace, amplitude=10pt, raise=3pt}] (0, 0) to node[midway, above, yshift=12pt] {$L$} (3, 0);
\end{tikzpicture}
\end{center}
\item
\label{not_separated_right}
There does not exist a maximal occurrence in $U_h$ such that it is not separated from $a_h$ from the left side and such that its intersection with $L$ is not a small piece. But there exists a maximal occurrence in $U_h$ such that it is not separated from $a_h$ from the right side and such that its intersection with $R$ is not a small piece. Let this be our $c$. As before, $\widehat{c}$ is its intersection with $R$.
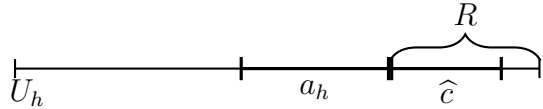
\begin{center}
\begin{tikzpicture}
\draw[|-|, black, thick] (0,0)--(7,0);
\node[below, xshift=5] at (0, 0) {$U_h$};
\draw[|-|, black, very thick] (3,0)--(5,0) node[midway, below] {$a_h$};
\draw[|-|, black, very thick] (5,0)--(6.5,0) node[midway, below] {$\widehat{c}$};
\draw [thick, decorate, decoration={brace, amplitude=10pt, raise=3pt}] (5, 0) to node[midway, above, yshift=12pt] {$R$} (7, 0);
\end{tikzpicture}
\end{center}
\item
\label{not_separated_left_right}
There exists a maximal occurrence in $U_h$ such that it is not separated from $a_h$ from the left side and such that its intersection with $L$ is not a small piece. Let this be our $b$. As before, $\widehat{b}$ is its intersection with $L$. And there exists a maximal occurrence in $U_h$ such that it is not separated from $a_h$ from the right side and such that its intersection with $R$ is not a small piece. Let this be our $c$. As before, $\widehat{c}$ is its intersection with $R$.
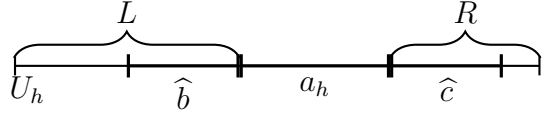
\begin{center}
\begin{tikzpicture}
\draw[|-|, black, thick] (0,0)--(7,0);
\node[below, xshift=5] at (0, 0) {$U_h$};
\draw[|-|, black, very thick] (3,0)--(5,0) node[midway, below] {$a_h$};
\draw[|-|, black, very thick] (5,0)--(6.5,0) node[midway, below] {$\widehat{c}$};
\draw[|-|, black, very thick] (1.5,0)--(3,0) node[midway, below] {$\widehat{b}$};
\draw [thick, decorate, decoration={brace, amplitude=10pt, raise=3pt}] (0, 0) to node[midway, above, yshift=12pt] {$L$} (3, 0);
\draw [thick, decorate, decoration={brace, amplitude=10pt, raise=3pt}] (5, 0) to node[midway, above, yshift=12pt] {$R$} (7, 0);
\end{tikzpicture}
\end{center}
\end{enumerate}
We have already considered case~\ref{no_not_separated}, so, we have to study only cases~\ref{not_separated_left}---\ref{not_separated_left_right}.

Now we treat monomials of types~\ref{keep_structure} and~\ref{donot_keep_structure1} separately (see page~\pageref{keep_structure}). First consider a monomial $U_j = La_jR$ of type~\ref{keep_structure}. Let us study what happen with $a_j$. There is the only possibility:
\begin{enumerate}[label=A\,\ref{keep_structure}.\arabic*]
\item
\label{a_j_keep_structure}
$a_j$ is a maximal occurrence in $U_j$.

Assume the contrary, that is, $La_jR = L^{\prime}X a_j YR^{\prime}$, $Xa_jY \in \Mon$, where $X$ is a suffix of $L$, $Y$ is a prefix of $R$, at least one of $X$ and $Y$ is non-empty. Since $a_j$ is not a small piece, we apply Corollary~\ref{not_sp_prolongation2} and obtain $Xa_hY \in \Mon$. Thus, $a_h$ is not a maximal occurrence in $U_h$. Contradiction.
\end{enumerate}

Now let us study possible forms of $b^{\prime}$ and $c^{\prime}$ in a monomial of type~\ref{keep_structure} in detail. First we study what happens with $b$ in cases~\ref{not_separated_left} and~\ref{not_separated_left_right}. That is, $b$ is not separated from $a_h$, $L = L_1\widehat{b}$. Then there are the following possibilities for $b^{\prime}$. Notice that they do not depend on a configuration of $c^{\prime}$ in the same monomial.
\begin{enumerate}[label=L\,\ref{keep_structure}.\arabic*]
\item
\label{keep_structure_case_l1}
$\widehat{b}$ is a maximal occurrence in $U_j$, that is, $b^{\prime} = \widehat{b}$. If $b$ and $a_h$ touch at a point in $U_h$, we obtain $b = \widehat{b}$ and $b^{\prime} = \widehat{b}$.
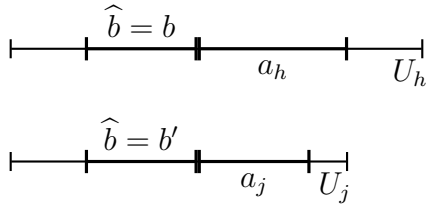
\begin{center}
\begin{tikzpicture}
\draw[|-|, black, thick] (0,0)--(5.5,0);
\node[below, xshift=-5] at (5.5, 0) {$U_h$};
\draw[|-|, black, very thick] (2.5,0)--(4.5,0) node[midway, below] {$a_h$};
\draw[|-|, black, very thick] (1,0)--(2.5,0) node[midway, above] {$\widehat{b} = b$};

\draw[|-|, black, thick] (0,-1.5)--(4.5,-1.5);
\node[below, xshift=-5] at (4.5, -1.5) {$U_j$};
\draw[|-|, black, very thick] (2.5,-1.5)--(4,-1.5) node[midway, below] {$a_j$};
\draw[|-|, black, very thick] (1,-1.5)--(2.5,-1.5) node[midway, above] {$\widehat{b} = b^{\prime}$};
\end{tikzpicture}
\end{center}
If $b$ and $a_h$ have an overlap $d_1$ in $U_h$, then $b= \widehat{b}d_1$, $b^{\prime} = \widehat{b}$.
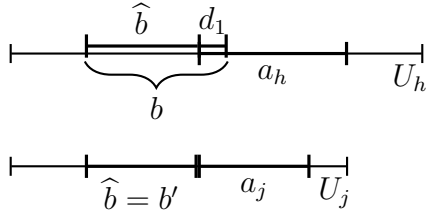
\begin{center}
\begin{tikzpicture}
\draw[|-|, black, thick] (0,0)--(5.5,0);
\node[below, xshift=-5] at (5.5, 0) {$U_h$};
\draw[|-|, black, very thick] (2.5,0)--(4.5,0) node[midway, below] {$a_h$};
\draw[|-, black, very thick] (1,0.1)--(2.5,0.1) node[midway, above] {$\widehat{b}$};
\draw[|-|, black, very thick] (2.5,0.1)--(2.9,0.1) node[midway, above] {$d_1$};
\draw [thick, decorate, decoration={brace, amplitude=10pt, raise=3pt, mirror}] (1, 0) to node[midway, below, yshift=-12pt] {$b$} (2.9, 0);

\draw[|-|, black, thick] (0,-1.5)--(4.5,-1.5);
\node[below, xshift=-5] at (4.5, -1.5) {$U_j$};
\draw[|-|, black, very thick] (2.5,-1.5)--(4,-1.5) node[midway, below] {$a_j$};
\draw[|-|, black, very thick] (1,-1.5)--(2.5,-1.5) node[midway, below] {$\widehat{b} = b^{\prime}$};
\end{tikzpicture}
\end{center}
\item
\label{keep_structure_case_l2}
$\widehat{b}$ is prolonged from the right in $U_j$ by a small piece, that is, $b^{\prime} = \widehat{b}e_1$, where $e_1$ is a non-empty small piece. If $e_1$ is not a small piece, then, by Lemma~\ref{not_sp_prolongation}, we obtain that $\widehat{b}a_j \in \Mon$. But this is not possible, since $a_j$ is a maximal occurrence in $U_j$, a contradiction.

If $b$ and $a_h$ touch at a point in $U_h$, we obtain $b = \widehat{b}$ and $b^{\prime} = \widehat{b}e_1$.
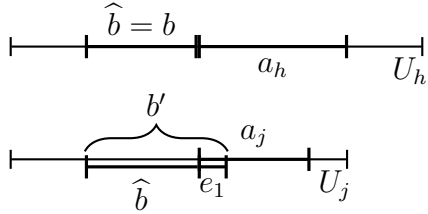
\begin{center}
\begin{tikzpicture}
\draw[|-|, black, thick] (0,0)--(5.5,0);
\node[below, xshift=-5] at (5.5, 0) {$U_h$};
\draw[|-|, black, very thick] (2.5,0)--(4.5,0) node[midway, below] {$a_h$};
\draw[|-|, black, very thick] (1,0)--(2.5,0) node[midway, above] {$\widehat{b} = b$};

\draw[|-|, black, thick] (0,-1.5)--(4.5,-1.5);
\node[below, xshift=-5] at (4.5, -1.5) {$U_j$};
\draw[|-|, black, very thick] (2.5,-1.5)--(4,-1.5) node[midway, above] {$a_j$};
\draw[|-, black, very thick] (1,-1.6)--(2.5,-1.6) node[midway, below] {$\widehat{b}$};
\draw[|-|, black, very thick] (2.5,-1.6)--(2.9,-1.6) node[midway, below] {$e_1$};
\draw [thick, decorate, decoration={brace, amplitude=10pt, raise=3pt}] (1, -1.5) to node[midway, above, yshift=12pt] {$b^{\prime}$} (2.9, -1.5);
\end{tikzpicture}
\end{center}

If $b$ and $a_h$ have an overlap $d_1$ in $U_h$, then $b = \widehat{b}d_1$, $b^{\prime} = \widehat{b}e_1$.
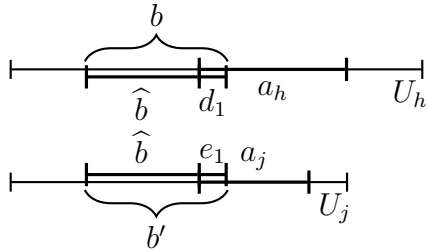
\begin{center}
\begin{tikzpicture}
\draw[|-|, black, thick] (0,0)--(5.5,0);
\node[below, xshift=-5] at (5.5, 0) {$U_h$};
\draw[|-|, black, very thick] (2.5,0)--(4.5,0) node[midway, below] {$a_h$};
\draw[|-, black, very thick] (1,-0.1)--(2.5,-0.1) node[midway, below] {$\widehat{b}$};
\draw[|-|, black, very thick] (2.5,-0.1)--(2.9,-0.1) node[midway, below] {$d_1$};
\draw [thick, decorate, decoration={brace, amplitude=10pt, raise=3pt}] (1, 0) to node[midway, above, yshift=12pt] {$b$} (2.9, 0);

\draw[|-|, black, thick] (0,-1.5)--(4.5,-1.5);
\node[below, xshift=-5] at (4.5, -1.5) {$U_j$};
\draw[|-|, black, very thick] (2.5,-1.5)--(4,-1.5) node[midway, above] {$a_j$};
\draw[|-, black, very thick] (1,-1.4)--(2.5,-1.4) node[midway, above] {$\widehat{b}$};
\draw[|-|, black, very thick] (2.5,-1.4)--(2.9,-1.4) node[midway, above] {$e_1$};
\draw [thick, decorate, decoration={brace, amplitude=10pt, raise=3pt, mirror}] (1, -1.5) to node[midway, below, yshift=-12pt] {$b^{\prime}$} (2.9, -1.5);
\end{tikzpicture}
\end{center}
\end{enumerate}

Now consider what happens with $c$ in cases~\ref{not_separated_right} and~\ref{not_separated_left_right}. That is, $c$ is not separated from $a_h$, $R = \widehat{c}R_1$. Then we obtain the following possibilities for $c^{\prime}$. They do not depend on a configuration of $b^{\prime}$ in the same monomial.
\begin{enumerate}[label=R\,\ref{keep_structure}.\arabic*]
\item
\label{keep_structure_case_r1}
$\widehat{c}$ is a maximal occurrence in $U_j$, that is, $c^{\prime} = \widehat{c}$. If $c$ and $a_h$ touch at a point in $U_h$, then $c = \widehat{c}$ and $c^{\prime} = \widehat{c}$.
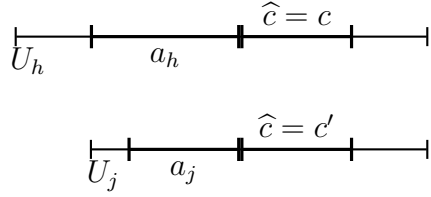
\begin{center}
\begin{tikzpicture}
\draw[|-|, black, thick] (0,0)--(5.5,0);
\node[below, xshift=5] at (0, 0) {$U_h$};
\draw[|-|, black, very thick] (1,0)--(3,0) node[midway, below] {$a_h$};
\draw[|-|, black, very thick] (3,0)--(4.5,0) node[midway, above] {$\widehat{c} = c$};

\draw[|-|, black, thick] (1,-1.5)--(5.5,-1.5);
\node[below, xshift=5] at (1, -1.5) {$U_j$};
\draw[|-|, black, very thick] (1.5,-1.5)--(3,-1.5) node[midway, below] {$a_j$};
\draw[|-|, black, very thick] (3,-1.5)--(4.5,-1.5) node[midway, above] {$\widehat{c} = c^{\prime}$};
\end{tikzpicture}
\end{center}
If $c$ and $a_h$ have an overlap $d_2$ in $U_h$, then $c = d_2\widehat{c}$ and $c^{\prime} = \widehat{c}$.
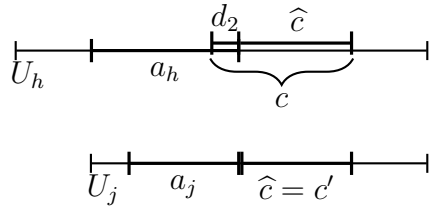
\begin{center}
\begin{tikzpicture}
\draw[|-|, black, thick] (0,0)--(5.5,0);
\node[below, xshift=5] at (0, 0) {$U_h$};
\draw[|-|, black, very thick] (1,0)--(3,0) node[midway, below] {$a_h$};
\draw[-|, black, very thick] (3,0.1)--(4.5,0.1) node[midway, above] {$\widehat{c}$};
\draw[|-|, black, very thick] (2.6,0.1)--(3,0.1) node[midway, above] {$d_2$};
\draw [thick, decorate, decoration={brace, amplitude=10pt, raise=3pt, mirror}] (2.6, 0) to node[midway, below, yshift=-12pt] {$c$} (4.5, 0);

\draw[|-|, black, thick] (1,-1.5)--(5.5,-1.5);
\node[below, xshift=5] at (1, -1.5) {$U_j$};
\draw[|-|, black, very thick] (1.5,-1.5)--(3,-1.5) node[midway, below] {$a_j$};
\draw[|-|, black, very thick] (3,-1.5)--(4.5,-1.5) node[midway, below] {$\widehat{c} = c^{\prime}$};
\end{tikzpicture}
\end{center}
\item
\label{keep_structure_case_r2}
$\widehat{c}$ is prolonged from the left in $U_j$ by a small piece, that is, $c^{\prime} = e_2\widehat{c}$, where $e_2$ is a non-empty small piece. If $e_2$ is not a small piece, then, by Lemma~\ref{not_sp_prolongation}, we obtain $a_j\widehat{c} \in \Mon$. But this is not possible, since $a_j$ is a maximal occurrence in $U_j$, a contradiction. Since $a_j$ is a maximal occurrence in $U_j$, $a_j$ is not contained in $e_2$.

If $c$ and $a_h$ touch at a point in $U_h$, then $c = \widehat{c}$ and $c^{\prime} = e_2\widehat{c}$.
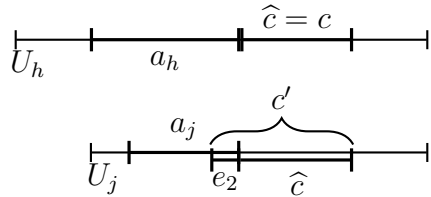
\begin{center}
\begin{tikzpicture}
\draw[|-|, black, thick] (0,0)--(5.5,0);
\node[below, xshift=5] at (0, 0) {$U_h$};
\draw[|-|, black, very thick] (1,0)--(3,0) node[midway, below] {$a_h$};
\draw[|-|, black, very thick] (3,0)--(4.5,0) node[midway, above] {$\widehat{c} = c$};

\draw[|-|, black, thick] (1,-1.5)--(5.5,-1.5);
\node[below, xshift=5] at (1, -1.5) {$U_j$};
\draw[|-|, black, very thick] (1.5,-1.5)--(3,-1.5) node[midway, above] {$a_j$};
\draw[-|, black, very thick] (3,-1.6)--(4.5,-1.6) node[midway, below] {$\widehat{c}$};
\draw[|-|, black, very thick] (2.6,-1.6)--(3,-1.6) node[midway, below] {$e_2$};
\draw [thick, decorate, decoration={brace, amplitude=10pt, raise=3pt}] (2.6, -1.5) to node[midway, above, yshift=12pt] {$c^{\prime}$} (4.5, -1.5);
\end{tikzpicture}
\end{center}
If $c$ and $a_h$ have an overlap $d_2$ in $U_h$, then $c = d_2\widehat{c}$ and $c^{\prime} = e_2\widehat{c}$.
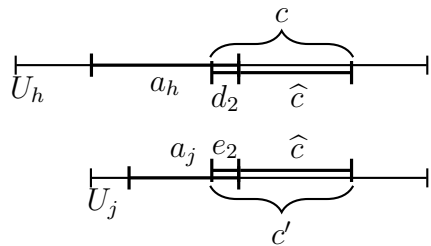
\begin{center}
\begin{tikzpicture}
\draw[|-|, black, thick] (0,0)--(5.5,0);
\node[below, xshift=5] at (0, 0) {$U_h$};
\draw[|-|, black, very thick] (1,0)--(3,0) node[midway, below] {$a_h$};
\draw[-|, black, very thick] (3,-0.1)--(4.5,-0.1) node[midway, below] {$\widehat{c}$};
\draw[|-|, black, very thick] (2.6,-0.1)--(3,-0.1) node[midway, below] {$d_2$};
\draw [thick, decorate, decoration={brace, amplitude=10pt, raise=3pt}] (2.6, 0) to node[midway, above, yshift=12pt] {$c$} (4.5, 0);

\draw[|-|, black, thick] (1,-1.5)--(5.5,-1.5);
\node[below, xshift=5] at (1, -1.5) {$U_j$};
\draw[|-|, black, very thick] (1.5,-1.5)--(3,-1.5) node[midway, above] {$a_j$};
\draw[-|, black, very thick] (3,-1.4)--(4.5,-1.4) node[midway, above] {$\widehat{c}$};
\draw[|-|, black, very thick] (2.6,-1.4)--(3,-1.4) node[midway, above] {$e_2$};
\draw [thick, decorate, decoration={brace, amplitude=10pt, raise=3pt, mirror}] (2.6, -1.5) to node[midway, below, yshift=-12pt] {$c^{\prime}$} (4.5, -1.5);
\end{tikzpicture}
\end{center}
\end{enumerate}

If both $b$ and $c$ are not separated from $a_h$ in $U_h$, $L = L_1\widehat{b}$, $R = \widehat{c}R_1$, then any combination of configurations~\ref{keep_structure_case_l1}, \ref{keep_structure_case_l2} from the left and \ref{keep_structure_case_r1}, \ref{keep_structure_case_r2} from the right is possible. Therefore, we obtain the following mutual positions of $b^{\prime}$ and $c^{\prime}$ in case~\ref{not_separated_left_right}.
\begin{enumerate}[label=\ref{keep_structure}.\arabic*]
\item
\label{keep_structure_case_m1}
$\widehat{b}$ is a maximal occurrence in $U_j$, that is, $b^{\prime} = \widehat{b}$, $\widehat{c}$ is a maximal occurrence in $U_j$, that is, $c^{\prime} = \widehat{c}$.
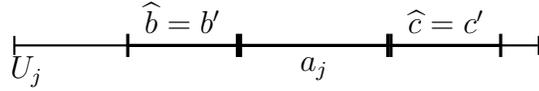
\begin{center}
\begin{tikzpicture}
\draw[|-|, black, thick] (0,0)--(7,0);
\node[below, xshift=5] at (0, 0) {$U_j$};
\draw[|-|, black, very thick] (3,0)--(5,0) node[midway, below] {$a_j$};
\draw[|-|, black, very thick] (5,0)--(6.5,0) node[midway, above] {$\widehat{c} = c^{\prime}$};
\draw[|-|, black, very thick] (1.5,0)--(3,0) node[midway, above] {$\widehat{b} = b^{\prime}$};
\end{tikzpicture}
\end{center}
\item
\label{keep_structure_case_m2}
$\widehat{b}$ is a maximal occurrence in $U_j$, that is, $b^{\prime} = \widehat{b}$, $\widehat{c}$ is prolonged in $U_j$ by a small piece, that is, $c^{\prime} = e_2\widehat{c}$, where $e_2$ is a non-empty small piece. Since $a_j$ is not a small piece, $b^{\prime}$ and $c^{\prime}$ are separated in $U_j$.
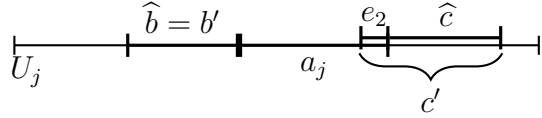
\begin{center}
\begin{tikzpicture}
\draw[|-|, black, thick] (0,0)--(7,0);
\node[below, xshift=5] at (0, 0) {$U_j$};
\draw[|-|, black, very thick] (3,0)--(5,0) node[midway, below] {$a_j$};
\draw[-|, black, very thick] (5,0.1)--(6.5,0.1) node[midway, above] {$\widehat{c}$};
\draw[|-|, black, very thick] (4.6,0.1)--(5,0.1) node[midway, above] {$e_2$};
\draw [thick, decorate, decoration={brace, amplitude=10pt, raise=3pt, mirror}] (4.6, 0) to node[midway, below, yshift=-12pt] {$c^{\prime}$} (6.5, 0);
\draw[|-|, black, very thick] (1.5,0)--(3,0) node[midway, above] {$\widehat{b} = b^{\prime}$};
\end{tikzpicture}
\end{center}
\item
\label{keep_structure_case_m3}
$\widehat{b}$ is prolonged by a small piece in $U_j$, that is, $b^{\prime} = \widehat{b}e_1$, where $e_1$ is a non-empty small piece, $\widehat{c}$ is a maximal occurrence in $U_j$, that is, $c^{\prime} = \widehat{c}$. Since $a_j$ is not a small piece, $b^{\prime}$ and $c^{\prime}$ are separated in $U_j$.
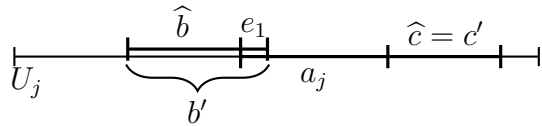
\begin{center}
\begin{tikzpicture}
\draw[|-|, black, thick] (0,0)--(7,0);
\node[below, xshift=5] at (0, 0) {$U_j$};
\draw[|-|, black, very thick] (3,0)--(5,0) node[midway, below] {$a_j$};
\draw[-|, black, very thick] (5,0)--(6.5,0) node[midway, above] {$\widehat{c} = c^{\prime}$};
\draw [thick, decorate, decoration={brace, amplitude=10pt, raise=3pt, mirror}] (1.5, 0) to node[midway, below, yshift=-12pt] {$b^{\prime}$} (3.4, 0);
\draw[|-, black, very thick] (1.5,0.1)--(3,0.1) node[midway, above] {$\widehat{b}$};
\draw[|-|, black, very thick] (3,0.1)--(3.4,0.1) node[midway, above] {$e_1$};
\end{tikzpicture}
\end{center}
\item
\label{keep_structure_case_m4}
both $\widehat{b}$ and $\widehat{c}$ are prolonged by small pieces in $U_j$, that is, $b^{\prime} = \widehat{b}e_1$, $c^{\prime} = e_2\widehat{c}$, where $e_1$ and $e_2$ are non-empty small pieces, and $b^{\prime}$ and $c^{\prime}$ are separated in $U_j$.
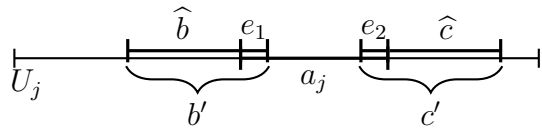
\begin{center}
\begin{tikzpicture}
\draw[|-|, black, thick] (0,0)--(7,0);
\node[below, xshift=5] at (0, 0) {$U_j$};
\draw[|-|, black, very thick] (3,0)--(5,0) node[midway, below] {$a_j$};
\draw[-|, black, very thick] (5,0.1)--(6.5,0.1) node[midway, above] {$\widehat{c}$};
\draw[|-|, black, very thick] (4.6,0.1)--(5,0.1) node[midway, above] {$e_2$};
\draw [thick, decorate, decoration={brace, amplitude=10pt, raise=3pt, mirror}] (1.5, 0) to node[midway, below, yshift=-12pt] {$b^{\prime}$} (3.4, 0);
\draw[|-, black, very thick] (1.5,0.1)--(3,0.1) node[midway, above] {$\widehat{b}$};
\draw[|-|, black, very thick] (3,0.1)--(3.4,0.1) node[midway, above] {$e_1$};
\draw [thick, decorate, decoration={brace, amplitude=10pt, raise=3pt, mirror}] (4.6, 0) to node[midway, below, yshift=-12pt] {$c^{\prime}$} (6.5, 0);
\end{tikzpicture}
\end{center}
\item
\label{keep_structure_case_m5}
both $\widehat{b}$ and $\widehat{c}$ are prolonged by small pieces in $U_j$, that is, $b^{\prime} = \widehat{b}e_1$, $c^{\prime} = e_2\widehat{c}$, where $e_1$ and $e_2$ are non-empty small pieces, and $b^{\prime}$ and $c^{\prime}$ touch at a point in $U_j$.
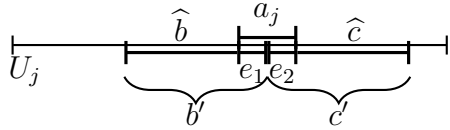
\begin{center}
\begin{tikzpicture}
\draw[|-|, black, thick] (0,0)--(7- 1.2,0);
\node[below, xshift=5] at (0, 0) {$U_j$};
\draw[|-|, black, very thick] (3,0.1)--(5 - 1.2,0.1) node[midway, above] {$a_j$};
\draw[-|, black, very thick] (5- 1.2,-0.1)--(6.5- 1.2,-0.1) node[midway, above] {$\widehat{c}$};
\draw[|-|, black, very thick] (4.6- 1.2,-0.1)--(5- 1.2,-0.1) node[midway, below] {$e_2$};
\draw [thick, decorate, decoration={brace, amplitude=10pt, raise=12pt, mirror}] (1.5, 0) to node[midway, below, yshift=-17pt] {$b^{\prime}$} (3.4, 0);
\draw[|-, black, very thick] (1.5,-0.1)--(3,-0.1) node[midway, above] {$\widehat{b}$};
\draw[|-|, black, very thick] (3,-0.1)--(3.4,-0.1) node[midway, below] {$e_1$};
\draw [thick, decorate, decoration={brace, amplitude=10pt, raise=12pt, mirror}] (4.6- 1.2, 0) to node[midway, below, yshift=-17pt] {$c^{\prime}$} (6.5- 1.2, 0);
\end{tikzpicture}
\end{center}
This may happen only if $\SPM(a_j) \leqslant 2$.
\item
\label{keep_structure_case_m6}
both $\widehat{b}$ and $\widehat{c}$ are prolonged by small pieces in $U_j$, that is, $b^{\prime} = \widehat{b}e_1$, $c^{\prime} = e_2\widehat{c}$, where $e_1$ and $e_2$ are non-empty small pieces, and $b^{\prime}$ and $c^{\prime}$ have an overlap in $U_j$.
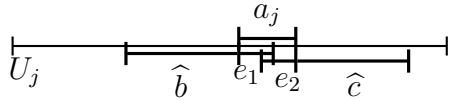
\begin{center}
\begin{tikzpicture}
\draw[|-|, black, thick] (0,0)--(7- 1.2,0);
\node[below, xshift=5] at (0, 0) {$U_j$};
\draw[|-|, black, very thick] (3,0.1)--(5 - 1.2,0.1) node[midway, above] {$a_j$};
\draw[-|, black, very thick] (5- 1.2,-0.2)--(6.5- 1.2,-0.2) node[midway, below] {$\widehat{c}$};
\draw[|-|, black, very thick] (4.5- 1.2,-0.2)--(5- 1.2,-0.2) node[near end, below] {$e_2$};
\draw[|-, black, very thick] (1.5,-0.1)--(3,-0.1) node[midway, below] {$\widehat{b}$};
\draw[|-|, black, very thick] (3,-0.1)--(3.5,-0.1) node[near start, below] {$e_1$};
\end{tikzpicture}
\end{center}
This may happen only if $\SPM(a_j) \leqslant 2$.
\end{enumerate}

Consider a monomial $U_j = La_jR$ of type~\ref{donot_keep_structure1} (see page~\pageref{donot_keep_structure1}). Let us study what happens with $b$ in cases~\ref{not_separated_left} and~\ref{not_separated_left_right}. That is, $b$ is not separated from $a_h$, $L = L_1\widehat{b}$. First we enumerate possibilities that do not depend on $c^{\prime}$. They can be obtained together in cases~\ref{not_separated_left} and~\ref{not_separated_left_right}. Then there are the following configurations for $b^{\prime}$:
\begin{enumerate}[label=L\,\ref{donot_keep_structure1}.\arabic*]
\item
\label{donot_keep_structure1_l1}
$\widehat{b}$ is a maximal occurrence in $U_j$, that is, $b^{\prime} = \widehat{b}$.
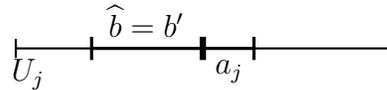
\begin{center}
\begin{tikzpicture}
\draw[|-|, black, thick] (0,0)--(5,0);
\node[below, xshift=5] at (0, 0) {$U_j$};
\draw[|-|, black, very thick] (2.5,0)--(3.2,0) node[midway, below] {$a_j$};
\draw[|-|, black, very thick] (1,0)--(2.5,0) node[midway, above] {$\widehat{b} = b^{\prime}$};
\end{tikzpicture}
\end{center}
\item
\label{donot_keep_structure1_l2}
$\widehat{b}$ is prolonged in $U_j$ by a small piece, that is, $b^{\prime} = \widehat{b}e_1$, where $e_1$ is a non-empty small piece, wherein $b^{\prime}$ does not cover $a_j$.
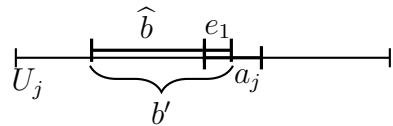
\begin{center}
\begin{tikzpicture}
\draw[|-|, black, thick] (0,0)--(5,0);
\node[below, xshift=5] at (0, 0) {$U_j$};
\draw[|-|, black, very thick] (2.5,0)--(3.3,0) node[near end, below] {$a_j$};
\draw[|-, black, very thick] (1,0.1)--(2.5,0.1) node[midway, above] {$\widehat{b}$};
\draw[|-|, black, very thick] (2.5,0.1)--(2.9,0.1) node[midway, above] {$e_1$};
\draw [thick, decorate, decoration={brace, amplitude=10pt, raise=3pt, mirror}] (1, 0) to node[midway, below, yshift=-12pt] {$b^{\prime}$} (2.9, 0);
\end{tikzpicture}
\end{center}
\item
\label{donot_keep_structure1_l3}
$\widehat{b}$ is prolonged in $U_j$ by a small piece, $b^{\prime} = \widehat{b}a_j$.
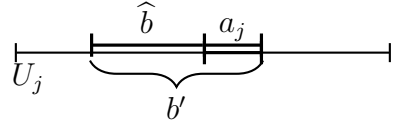
\begin{center}
\begin{tikzpicture}
\draw[|-|, black, thick] (0,0)--(5,0);
\node[below, xshift=5] at (0, 0) {$U_j$};
\draw[|-|, black, very thick] (2.5,0)--(3.3,0) node[midway, above] {$a_j$};
\draw[|-, black, very thick] (1,0.1)--(2.5,0.1) node[midway, above] {$\widehat{b}$};
\draw[|-|, black, very thick] (2.5,0.1)--(3.3,0.1);
\draw [thick, decorate, decoration={brace, amplitude=10pt, raise=3pt, mirror}] (1, 0) to node[midway, below, yshift=-12pt] {$b^{\prime}$} (3.3, 0);
\end{tikzpicture}
\end{center}
\item
\label{donot_keep_structure1_l4}
$\widehat{b}$ is prolonged in $U_j$, $b^{\prime}$ covers $a_j$, $b^{\prime} = \widehat{b}a_je_1$, where $e_1$ is a non-empty small piece.
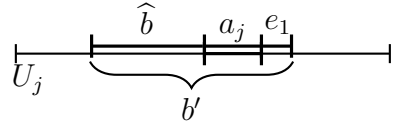
\begin{center}
\begin{tikzpicture}
\draw[|-|, black, thick] (0,0)--(5,0);
\node[below, xshift=5] at (0, 0) {$U_j$};
\draw[|-|, black, very thick] (2.5,0)--(3.3,0) node[midway, above] {$a_j$};
\draw[|-, black, very thick] (1,0.1)--(2.5,0.1) node[midway, above] {$\widehat{b}$};
\draw[|-|, black, very thick] (2.5,0.1)--(3.3,0.1);
\draw[-|, black, very thick] (3.3,0.1)--(3.7,0.1) node[midway, above] {$e_1$};;
\draw [thick, decorate, decoration={brace, amplitude=10pt, raise=3pt, mirror}] (1, 0) to node[midway, below, yshift=-12pt] {$b^{\prime}$} (3.7, 0);
\end{tikzpicture}
\end{center}
\end{enumerate}
\begin{remark}
\label{remark_possible_prolongations}
Assume $b^{\prime} = \widehat{b}a_je_1$ and $e_1$ is not a small piece.
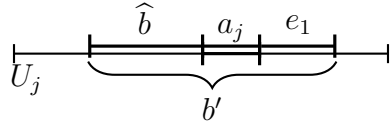
\begin{center}
\begin{tikzpicture}
\draw[|-|, black, thick] (0,0)--(5,0);
\node[below, xshift=5] at (0, 0) {$U_j$};
\draw[|-|, black, very thick] (2.5,0)--(3.3,0) node[midway, above] {$a_j$};
\draw[|-, black, very thick] (1,0.1)--(2.5,0.1) node[midway, above] {$\widehat{b}$};
\draw[|-|, black, very thick] (2.5,0.1)--(3.3,0.1);
\draw[-|, black, very thick] (3.3,0.1)--(4.3,0.1) node[midway, above] {$e_1$};
\draw [thick, decorate, decoration={brace, amplitude=10pt, raise=3pt, mirror}] (1, 0) to node[midway, below, yshift=-12pt] {$b^{\prime}$} (4.3, 0);
\end{tikzpicture}
\end{center}
Since $e_1$ is contained inside $R$, $e_1$ can be considered as an occurrence in $U_h$.
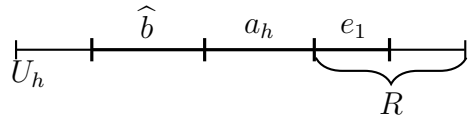
\begin{center}
\begin{tikzpicture}
\draw[|-|, black, thick] (0,0)--(6,0);
\node[below, xshift=5] at (0, 0) {$U_h$};
\draw[|-|, black, very thick] (2.5,0)--(4,0) node[midway, above] {$a_h$};
\draw[|-, black, very thick] (1,0)--(2.5,0) node[midway, above] {$\widehat{b}$};
\draw[-|, black, very thick] (4,0)--(5,0) node[midway, above] {$e_1$};
\draw [thick, decorate, decoration={brace, amplitude=10pt, raise=3pt, mirror}] (4, 0) to node[midway, below, yshift=-12pt] {$R$} (6, 0);
\end{tikzpicture}
\end{center}
Since $e_1$ is not a small piece, there exists a unique maximal occurrence $e$ in $U_h$ that contains $e_1$. Clearly, $e$ is not separated from $a_h$ and its intersection with $R$ is equal to $e_1$. There are two possibilities:
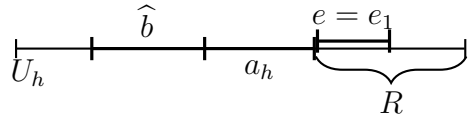
\begin{center}
\begin{tikzpicture}
\draw[|-|, black, thick] (0,0)--(6,0);
\node[below, xshift=5] at (0, 0) {$U_h$};
\draw[|-|, black, very thick] (2.5,0)--(4,0) node[midway, below] {$a_h$};
\draw[|-, black, very thick] (1,0)--(2.5,0) node[midway, above] {$\widehat{b}$};
\draw[|-|, black, very thick] (4,0.1)--(5,0.1) node[midway, above] {$e = e_1$};
\draw [thick, decorate, decoration={brace, amplitude=10pt, raise=3pt, mirror}] (4, 0) to node[midway, below, yshift=-12pt] {$R$} (6, 0);
\end{tikzpicture}

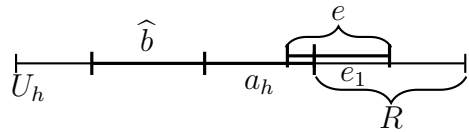
\begin{tikzpicture}
\draw[|-|, black, thick] (0,0)--(6,0);
\node[below, xshift=5] at (0, 0) {$U_h$};
\draw[|-|, black, very thick] (2.5,0)--(4,0) node[midway, below] {$a_h$};
\draw[|-, black, very thick] (1,0)--(2.5,0) node[midway, above] {$\widehat{b}$};
\draw[-|, black, very thick] (4,0.1)--(5,0.1) node[midway, below] {$e_1$};
\draw[|-|, black, very thick] (3.6,0.1)--(4,0.1);
\draw [thick, decorate, decoration={brace, amplitude=10pt, raise=5pt, mirror}] (4, 0) to node[midway, below, yshift=-12pt] {$R$} (6, 0);
\draw [thick, decorate, decoration={brace, amplitude=8pt, raise=6pt}] (3.6, 0) to node[midway, above, yshift=12pt] {$e$} (5, 0);
\end{tikzpicture}
\end{center}
Then, since $e_1$ is not a small piece, we obtain $e = c$, $e_1 = \widehat{c}$ and $b^{\prime} = \widehat{b}a_j\widehat{c}$.
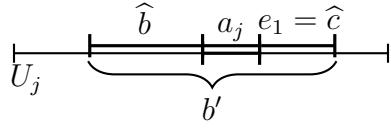
\begin{center}
\begin{tikzpicture}
\draw[|-|, black, thick] (0,0)--(5,0);
\node[below, xshift=5] at (0, 0) {$U_j$};
\draw[|-|, black, very thick] (2.5,0)--(3.3,0) node[midway, above] {$a_j$};
\draw[|-, black, very thick] (1,0.1)--(2.5,0.1) node[midway, above] {$\widehat{b}$};
\draw[|-|, black, very thick] (2.5,0.1)--(3.3,0.1);
\draw[-|, black, very thick] (3.3,0.1)--(4.3,0.1) node[midway, above] {$e_1 = \widehat{c}$};
\draw [thick, decorate, decoration={brace, amplitude=10pt, raise=3pt, mirror}] (1, 0) to node[midway, below, yshift=-12pt] {$b^{\prime}$} (4.3, 0);
\end{tikzpicture}
\end{center}
Therefore, $b^{\prime} = \widehat{b}a_je_1$, where $e_1$ is not a small piece, can be obtained only in case~\ref{not_separated_left_right}, and moreover $b^{\prime} = \widehat{b}a_j\widehat{c}$ in this case. So, \ref{donot_keep_structure1_l1}---\ref{donot_keep_structure1_l4} exhaust all possibilities that can be obtained both in cases~\ref{not_separated_left} and~\ref{not_separated_left_right}.

Clearly, the same argument is applicable for configurations~\ref{donot_keep_structure1_r1}---\ref{donot_keep_structure1_r4} and~\ref{donot_keep_structure2_l1}---\ref{donot_keep_structure2_l2}, \ref{donot_keep_structure2_r1}---\ref{donot_keep_structure2_r2}.
\end{remark}

Now we consider what happens with $c$ in cases~\ref{not_separated_right} and~\ref{not_separated_left_right}. So $c$ is not separated from $a_h$, $R = \widehat{c}R_1$. Let us enumerate possibilities that can be obtained together in cases~\ref{not_separated_right} and~\ref{not_separated_left_right}. They do not depend on a configuration of $b^{\prime}$ in the same monomial $U_j$. Then there are the following configurations for $c^{\prime}$ (see Remark~\ref{remark_possible_prolongations}):
\begin{enumerate}[label=R\,\ref{donot_keep_structure1}.\arabic*]
\item
\label{donot_keep_structure1_r1}
$\widehat{c}$ is a maximal occurrence in $U_j$, that is, $c^{\prime} = \widehat{c}$.
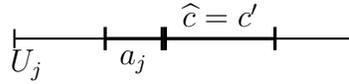
\begin{center}
\begin{tikzpicture}
\draw[|-|, black, thick] (1,0)--(5.5,0);
\node[below, xshift=5] at (1, 0) {$U_j$};
\draw[|-|, black, very thick] (2.2,0)--(3,0) node[midway, below] {$a_j$};
\draw[|-|, black, very thick] (3,0)--(4.5,0) node[midway, above] {$\widehat{c} = c^{\prime}$};
\end{tikzpicture}
\end{center}
\item
\label{donot_keep_structure1_r2}
$\widehat{c}$ is prolonged in $U_j$, that is, $c^{\prime} = e_2\widehat{c}$, where $e_2$ is a non-empty small piece, wherein $c^{\prime}$ does not cover $a_j$.
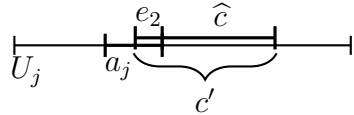
\begin{center}
\begin{tikzpicture}
\draw[|-|, black, thick] (1,0)--(5.5,0);
\node[below, xshift=5] at (1, 0) {$U_j$};
\draw[|-|, black, very thick] (2.2,0)--(3,0) node[near start, below] {$a_j$};
\draw[-|, black, very thick] (3,0.1)--(4.5,0.1) node[midway, above] {$\widehat{c}$};
\draw[|-|, black, very thick] (2.6,0.1)--(3,0.1) node[midway, above] {$e_2$};
\draw [thick, decorate, decoration={brace, amplitude=10pt, raise=3pt, mirror}] (2.6, 0) to node[midway, below, yshift=-12pt] {$c^{\prime}$} (4.5, 0);
\end{tikzpicture}
\end{center}
\item
\label{donot_keep_structure1_r3}
$\widehat{c}$ is prolonged in $U_j$, $c^{\prime} = a_j\widehat{c}$.
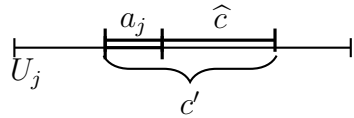
\begin{center}
\begin{tikzpicture}
\draw[|-|, black, thick] (1,0)--(5.5,0);
\node[below, xshift=5] at (1, 0) {$U_j$};
\draw[|-|, black, very thick] (2.2,0)--(3,0) node[midway, above] {$a_j$};
\draw[-|, black, very thick] (3,0.1)--(4.5,0.1) node[midway, above] {$\widehat{c}$};
\draw[|-|, black, very thick] (2.2,0.1)--(3,0.1);
\draw [thick, decorate, decoration={brace, amplitude=10pt, raise=3pt, mirror}] (2.2, 0) to node[midway, below, yshift=-12pt] {$c^{\prime}$} (4.5, 0);
\end{tikzpicture}
\end{center}
\item
\label{donot_keep_structure1_r4}
$\widehat{c}$ is prolonged in $U_j$, $c^{\prime}$ covers $a_j$, $c^{\prime} = e_2a_j\widehat{c}$, where $e_2$ is a non-empty small piece.
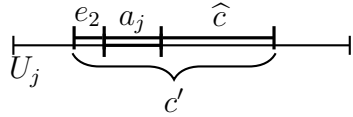
\begin{center}
\begin{tikzpicture}
\draw[|-|, black, thick] (1,0)--(5.5,0);
\node[below, xshift=5] at (1, 0) {$U_j$};
\draw[|-|, black, very thick] (2.2,0)--(3,0) node[midway, above] {$a_j$};
\draw[-|, black, very thick] (3,0.1)--(4.5,0.1) node[midway, above] {$\widehat{c}$};
\draw[|-|, black, very thick] (2.2,0.1)--(3,0.1);
\draw[|-, black, very thick] (1.8,0.1)--(2.2,0.1) node[midway, above] {$e_2$};
\draw [thick, decorate, decoration={brace, amplitude=10pt, raise=3pt, mirror}] (1.8, 0) to node[midway, below, yshift=-12pt] {$c^{\prime}$} (4.5, 0);
\end{tikzpicture}
\end{center}
\end{enumerate}

Assume both $b$ and $c$ are not separated from $a_h$ in $U_h$, $L = L_1\widehat{b}$, $R = \widehat{c}R_1$ (case~\ref{not_separated_left_right}). Then we can have any combination of configurations~\ref{donot_keep_structure1_l1}---\ref{donot_keep_structure1_l4} for $b^{\prime}$ and configurations~\ref{donot_keep_structure1_r1}---\ref{donot_keep_structure1_r4} for $c^{\prime}$. In these cases $b^{\prime}$ and $c^{\prime}$ are different maximal occurrences in $U_j$ and they may have any mutual position: be separated, touch at a point or have an overlap. If they have an overlap, then their overlap may contain or may not contain $a_j$. But here we have an extra possibility. Unlike monomials of type~\ref{keep_structure}, $b^{\prime}$ and $c^{\prime}$ may also be the same maximal occurrence in $U_j$:
\begin{enumerate}[label=\ref{donot_keep_structure1}.\arabic*]
\item
\label{donot_keep_structure1_merge}
$\widehat{b}$ and $\widehat{c}$ merge to the maximal occurrence $\widehat{b}a_j\widehat{c}$ in $U_j$, that is, $b^{\prime} = c^{\prime} = \widehat{b}a_j\widehat{c}$.
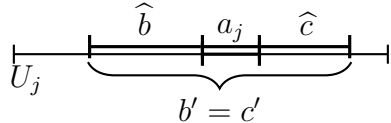
\begin{center}
\begin{tikzpicture}
\draw[|-|, black, thick] (0,0)--(5,0);
\node[below, xshift=5] at (0, 0) {$U_j$};
\draw[|-|, black, very thick] (2.5,0)--(3.3,0) node[midway, above] {$a_j$};
\draw[|-, black, very thick] (1,0.1)--(2.5,0.1) node[midway, above] {$\widehat{b}$};
\draw[|-|, black, very thick] (2.5,0.1)--(3.3,0.1);
\draw[-|, black, very thick] (3.3,0.1)--(4.5,0.1) node[midway, above] {$\widehat{c}$};
\draw [thick, decorate, decoration={brace, amplitude=10pt, raise=3pt, mirror}] (1, 0) to node[midway, below, yshift=-12pt] {$b^{\prime} = c^{\prime}$} (4.5, 0);
\end{tikzpicture}
\end{center}
\end{enumerate}

We continue to study a monomial $U_j = La_jR$ of type~\ref{donot_keep_structure1}. Let us consider what happens with $a_j$. Assume $a_j^{\prime}$ is a maximal occurrence in $U_j$ such that $a_j^{\prime}$ contains $a_j$. Clearly, $a_j^{\prime}$ is not necessarily unique. The first option is
\begin{enumerate}[label=A\,\ref{donot_keep_structure1}.\arabic*]
\item
\label{a_j_donot_keep_structure1_1}
$a_j = a_j^{\prime}$.
\end{enumerate}
Let $a_j \neq a_j^{\prime}$. Almost all configurations are already described above, namely, cases~\ref{donot_keep_structure1_l3}, \ref{donot_keep_structure1_l4}, \ref{donot_keep_structure1_r3}, \ref{donot_keep_structure1_r4}, \ref{donot_keep_structure1_merge}. The rest of the possibilities are the following:
\begin{enumerate}[label=A\,\ref{donot_keep_structure1}.\arabic*]
\setcounter{enumi}{1}
\item
\label{a_j_donot_keep_structure1_2}
$a_j^{\prime} = a_je_2$, where $e_2$ is a non-empty small piece.
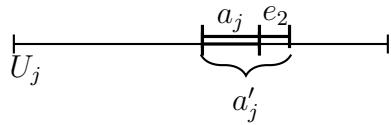
\begin{center}
\begin{tikzpicture}
\draw[|-|, black, thick] (0,0)--(5,0);
\node[below, xshift=5] at (0, 0) {$U_j$};
\draw[|-|, black, very thick] (2.5,0)--(3.3,0) node[midway, above] {$a_j$};
\draw[|-|, black, very thick] (2.5,0.1)--(3.3,0.1);
\draw[-|, black, very thick] (3.3,0.1)--(3.7,0.1) node[midway, above] {$e_2$};
\draw [thick, decorate, decoration={brace, amplitude=8pt, raise=4pt, mirror}] (2.5, 0) to node[midway, below, yshift=-12pt] {$a_j^{\prime}$} (3.7, 0);
\end{tikzpicture}
\end{center}
\item
\label{a_j_donot_keep_structure1_3}
$a_j^{\prime} = e_1a_j$, where $e_1$ is a non-empty small piece.
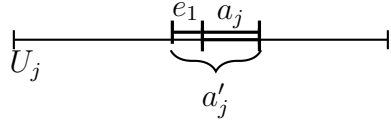
\begin{center}
\begin{tikzpicture}
\draw[|-|, black, thick] (0,0)--(5,0);
\node[below, xshift=5] at (0, 0) {$U_j$};
\draw[|-|, black, very thick] (2.5,0)--(3.3,0) node[midway, above] {$a_j$};
\draw[|-, black, very thick] (2.1,0.1)--(2.5,0.1) node[midway, above] {$e_1$};
\draw[|-|, black, very thick] (2.5,0.1)--(3.3,0.1);
\draw [thick, decorate, decoration={brace, amplitude=8pt, raise=4pt, mirror}] (2.1, 0) to node[midway, below, yshift=-12pt] {$a_j^{\prime}$} (3.3, 0);
\end{tikzpicture}
\end{center}
\item
\label{a_j_donot_keep_structure1_4}
$a_j^{\prime} = e_1a_je_2$, where $e_1$ and $e_2$ are non-empty small pieces.
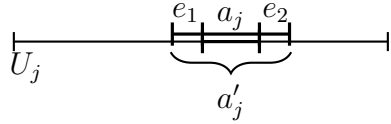
\begin{center}
\begin{tikzpicture}
\draw[|-|, black, thick] (0,0)--(5,0);
\node[below, xshift=5] at (0, 0) {$U_j$};
\draw[|-|, black, very thick] (2.5,0)--(3.3,0) node[midway, above] {$a_j$};
\draw[|-, black, very thick] (2.1,0.1)--(2.5,0.1) node[midway, above] {$e_1$};
\draw[|-|, black, very thick] (2.5,0.1)--(3.3,0.1);
\draw[-|, black, very thick] (3.3,0.1)--(3.7,0.1) node[midway, above] {$e_2$};
\draw [thick, decorate, decoration={brace, amplitude=8pt, raise=4pt, mirror}] (2.1, 0) to node[midway, below, yshift=-12pt] {$a_j^{\prime}$} (3.7, 0);
\end{tikzpicture}
\end{center}
\end{enumerate}

Now we study a monomial of type~\ref{donot_keep_structure2}, that is, $a_j = 1$, $U_j = L\cdot R$. In this case cancellations between $L$ and $R$ may occur. Suppose $L = L^{\prime}C$, $R = C^{-1}R^{\prime}$ and $L^{\prime}R^{\prime}$ does not have further cancellations.
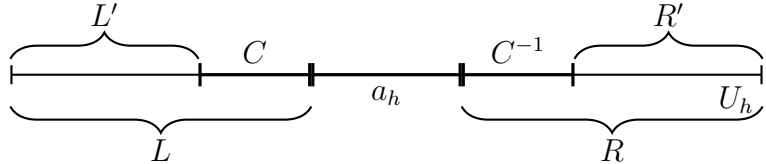
\begin{center}
\begin{tikzpicture}
\draw[|-|, black, thick] (0,0)--(10,0) node[at end, below, xshift=-10] {$U_h$};
\draw [thick, decorate, decoration={brace, amplitude=10pt, raise=6pt}] (0,0) to node[midway, above, yshift=14pt] {$L^{\prime}$} (2.5, 0);
\draw [thick, decorate, decoration={brace, amplitude=10pt, raise=6pt}] (7.5,0) to node[midway, above, yshift=14pt] {$R^{\prime}$} (10, 0);
\draw[|-|, black, very thick] (4,0)--(6,0) node[midway, below] {$a_h$};
\draw[|-|, black, very thick] (6,0)--(7.5,0) node[midway, above] {$C^{-1}$};
\draw[|-|, black, very thick] (2.5,0)--(4,0) node[midway, above] {$C$};
\draw [thick, decorate, decoration={brace, amplitude=10pt, raise=12pt, mirror}] (0,0) to node[midway, below, yshift=-20pt] {$L$} (4, 0);
\draw [thick, decorate, decoration={brace, amplitude=10pt, raise=12pt, mirror}] (6,0) to node[midway, below, yshift=-20pt] {$R$} (10, 0);
\end{tikzpicture}
\end{center}
If both $L^{\prime}$ and $R^{\prime}$ are empty, then $U_j = 1$ and there are no maximal occurrences to study. From now on we assume that at least one of $L^{\prime}$ and $R^{\prime}$ is non-empty. If some maximal occurrence in $U_h$ is fully contained in $C$, then, clearly, it is cancelled in $U_j$. If some maximal occurrence in $U_h$ is fully contained in $C^{-1}$, then it is cancelled in $U_j$. We have the following possibilities:
\begin{itemize}
\item
there exist maximal occurrences in $U_h$ that have a non-empty intersection with $L^{\prime}$ and that have a non-empty intersection with $R^{\prime}$;
\item
there are no maximal occurrences in $U_h$ that have a non-empty intersection with $L^{\prime}$;
\item
there are no maximal occurrences in $U_h$ that have a non-empty intersection with $R^{\prime}$;
\item
there are neither maximal occurrences in $U_h$ that have a non-empty intersection with $L^{\prime}$, nor maximal occurrences in $U_h$ that have a non-empty intersection with $R^{\prime}$.
\end{itemize}
We consider only the first possibility as the most interesting. Other cases are considered in the similar (but simpler) way.

We consider $b$ such that it has a non trivial intersection with $L^{\prime}$ and $c$ such that it has a non trivial intersection with $R^{\prime}$. Denote by $\widehat{b}$ the intersection of $b$ and $L^{\prime}$ and by $\widehat{c}$ the intersection of $c$ and $R^{\prime}$. Clearly, as above, $\widehat{b}$ and $\widehat{c}$ can be considered as occurrences in $U_j = L^{\prime}R^{\prime}$.
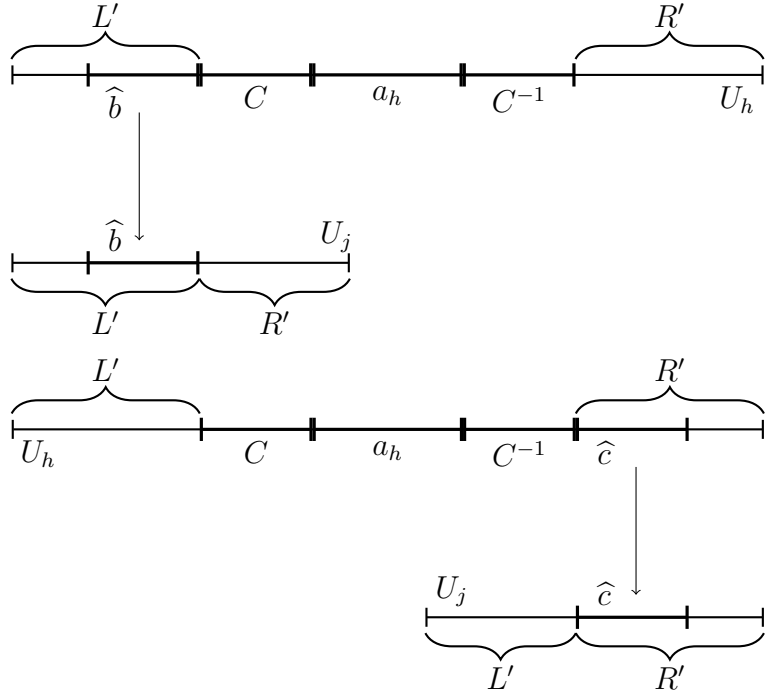
\begin{center}
\begin{tikzpicture}
\draw[|-|, black, thick] (0,0)--(10,0) node[at end, below, xshift=-10] {$U_h$};
\draw [thick, decorate, decoration={brace, amplitude=10pt, raise=6pt}] (0,0) to node[midway, above, yshift=14pt] {$L^{\prime}$} (2.5, 0);
\draw [thick, decorate, decoration={brace, amplitude=10pt, raise=6pt}] (7.5,0) to node[midway, above, yshift=14pt] {$R^{\prime}$} (10, 0);
\draw[|-|, black, very thick] (4,0)--(6,0) node[midway, below] {$a_h$};
\draw[|-|, black, very thick] (6,0)--(7.5,0) node[midway, below] {$C^{-1}$};
\draw[|-|, black, very thick] (2.5,0)--(4,0) node[midway, below] {$C$};
\draw[|-|, black, very thick] (1,0)--(2.5,0) node[near start, below] {$\widehat{b}$};

\draw[|-|, black, thick] (0,-2.5)--(4.5,-2.5);
\node[above, xshift=-5] at (4.5, -2.5) {$U_j$};
\draw [thick, decorate, decoration={brace, amplitude=10pt, raise=6pt, mirror}] (0,-2.5) to node[midway, below, yshift=-14pt] {$L^{\prime}$} (2.5, -2.5);
\draw[|-|, black, very thick] (1,-2.5)--(2.5,-2.5) node[near start, above] {$\widehat{b}$};
\draw [thick, decorate, decoration={brace, amplitude=10pt, raise=6pt, mirror}] (2.5,-2.5) to node[midway, below, yshift=-14pt] {$R^{\prime}$} (4.5, -2.5);

\draw[->, black] (1.7, -0.5)--(1.7,-2.2);
\end{tikzpicture}

\begin{tikzpicture}
\draw[|-|, black, thick] (0,0)--(10,0) node[at start, below, xshift=10] {$U_h$};
\draw [thick, decorate, decoration={brace, amplitude=10pt, raise=6pt}] (0,0) to node[midway, above, yshift=14pt] {$L^{\prime}$} (2.5, 0);
\draw [thick, decorate, decoration={brace, amplitude=10pt, raise=6pt}] (7.5,0) to node[midway, above, yshift=14pt] {$R^{\prime}$} (10, 0);
\draw[|-|, black, very thick] (4,0)--(6,0) node[midway, below] {$a_h$};
\draw[|-|, black, very thick] (6,0)--(7.5,0) node[midway, below] {$C^{-1}$};
\draw[|-|, black, very thick] (2.5,0)--(4,0) node[midway, below] {$C$};
\draw[|-|, black, very thick] (7.5,0)--(9,0) node[near start, below] {$\widehat{c}$};

\draw[|-|, black, thick] (5.5,-2.5)--(10,-2.5);
\node[above, xshift=10] at (5.5, -2.5) {$U_j$};
\draw [thick, decorate, decoration={brace, amplitude=10pt, raise=6pt, mirror}] (5.5,-2.5) to node[midway, below, yshift=-14pt] {$L^{\prime}$} (7.5, -2.5);
\draw[|-|, black, very thick] (7.5,-2.5)--(9,-2.5) node[near start, above] {$\widehat{c}$};
\draw [thick, decorate, decoration={brace, amplitude=10pt, raise=6pt, mirror}] (7.5,-2.5) to node[midway, below, yshift=-14pt] {$R^{\prime}$} (10, -2.5);

\draw[->, black] (8.3, -0.5)--(8.3,-2.2);
\end{tikzpicture}
\end{center}
The occurrences $\widehat{b}$ and $\widehat{c}$ do not have to be maximal occurrences in $U_j$. As above, we consider only $b$ and $c$ such that both $\widehat{b}$ and $\widehat{c}$ are not small pieces. Then, by Corollary~\ref{max_occurrences_coinside}, there exists a unique maximal occurrence in $U_j$ that contains $\widehat{b}$, we denote it by $b^{\prime}$, and there exists a unique maximal occurrence in $U_j$ that contains $\widehat{c}$, we denote it by $c^{\prime}$.

Similar to the above, if $\widehat{b}$ is not a terminal subword of $L^{\prime}$, then $\widehat{b} = b$ in $U_h$ and $\widehat{b}$ is separated from the end of $L^{\prime}$ in $U_j$. So, $\widehat{b}$ is a maximal occurrence in $U_j$. Hence, $b$ stays unchanged in $U_j$.
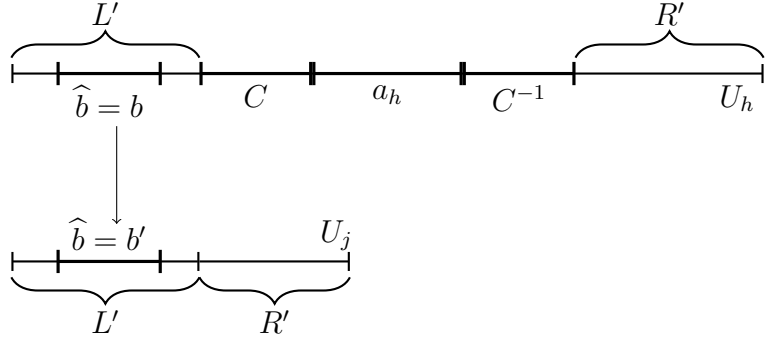
\begin{center}
\begin{tikzpicture}
\draw[|-|, black, thick] (0,0)--(10,0) node[at end, below, xshift=-10] {$U_h$};
\draw [thick, decorate, decoration={brace, amplitude=10pt, raise=6pt}] (0,0) to node[midway, above, yshift=14pt] {$L^{\prime}$} (2.5, 0);
\draw[|-|, black, very thick] (4,0)--(6,0) node[midway, below] {$a_h$};
\draw[|-|, black, very thick] (6,0)--(7.5,0) node[midway, below] {$C^{-1}$};
\draw[|-|, black, very thick] (2.5,0)--(4,0) node[midway, below] {$C$};
\draw[|-|, black, very thick] (0.6,0)--(2,0) node[midway, below] {$\widehat{b} = b$};
\draw [thick, decorate, decoration={brace, amplitude=10pt, raise=6pt}] (7.5,0) to node[midway, above, yshift=14pt] {$R^{\prime}$} (10, 0);

\draw[-|, black, thick] (2.5,-2.5)--(4.5,-2.5);
\draw[|-|, black, thick] (0,-2.5)--(2.5,-2.5);
\node[above, xshift=-5] at (4.5, -2.5) {$U_j$};
\draw [thick, decorate, decoration={brace, amplitude=10pt, raise=6pt, mirror}] (0,-2.5) to node[midway, below, yshift=-14pt] {$L^{\prime}$} (2.5, -2.5);
\draw[|-|, black, very thick] (0.6,-2.5)--(2,-2.5) node[midway, above] {$\widehat{b} = b^{\prime}$};
\draw [thick, decorate, decoration={brace, amplitude=10pt, raise=6pt, mirror}] (2.5,-2.5) to node[midway, below, yshift=-14pt] {$R^{\prime}$} (4.5, -2.5);

\draw[->, black] (1.4, -0.7)--(1.4,-2);
\end{tikzpicture}
\end{center}
If $\widehat{c}$ is not an initial subword of $R^{\prime}$, then $\widehat{c} = c$ in $U_h$ and $\widehat{c}$ is separated from the end of $L^{\prime}$ in $U_j$. So, $\widehat{c}$ is a maximal occurrence in $U_j$. Therefore, $c$ stays unchanged in $U_j$.
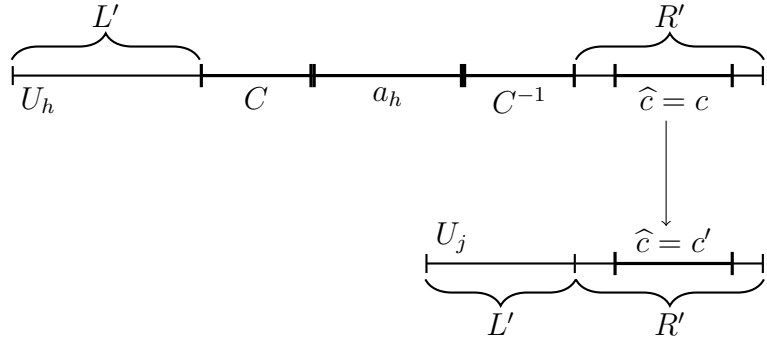
\begin{center}
\begin{tikzpicture}
\draw[|-|, black, thick] (0,0)--(10,0) node[at start, below, xshift=10] {$U_h$};
\draw [thick, decorate, decoration={brace, amplitude=10pt, raise=6pt}] (0,0) to node[midway, above, yshift=14pt] {$L^{\prime}$} (2.5, 0);
\draw [thick, decorate, decoration={brace, amplitude=10pt, raise=6pt}] (7.5,0) to node[midway, above, yshift=14pt] {$R^{\prime}$} (10, 0);
\draw[|-|, black, very thick] (4,0)--(6,0) node[midway, below] {$a_h$};
\draw[|-|, black, very thick] (6,0)--(7.5,0) node[midway, below] {$C^{-1}$};
\draw[|-|, black, very thick] (2.5,0)--(4,0) node[midway, below] {$C$};
\draw[|-|, black, very thick] (8,0)--(9.6,0) node[midway, below] {$\widehat{c} = c$};

\draw[|-|, black, thick] (5.5,-2.5)--(7.5,-2.5);
\draw[-|, black, thick] (7.5,-2.5)--(10,-2.5);
\node[above, xshift=10] at (5.5, -2.5) {$U_j$};
\draw [thick, decorate, decoration={brace, amplitude=10pt, raise=6pt, mirror}] (5.5,-2.5) to node[midway, below, yshift=-14pt] {$L^{\prime}$} (7.5, -2.5);
\draw[|-|, black, very thick] (8,-2.5)--(9.6,-2.5) node[midway, above] {$\widehat{c} = c^{\prime}$};
\draw [thick, decorate, decoration={brace, amplitude=10pt, raise=6pt, mirror}] (7.5,-2.5) to node[midway, below, yshift=-14pt] {$R^{\prime}$} (10, -2.5);

\draw[->, black] (8.7, -0.6)--(8.7,-2);
\end{tikzpicture}
\end{center}

If $\widehat{b}$ is a terminal subword of $L^{\prime}$, then, using the same argument as we used for monomials of types~\ref{keep_structure} and~\ref{donot_keep_structure1}, we can show that $b^{\prime} = \widehat{b}Y$, where $Y$ is a suffix of $b^{\prime}$. The same for $c^{\prime}$, if $\widehat{c}$ is an initial subword of $R^{\prime}$, then we obtain $c^{\prime} = Z\widehat{c}$, where $Z$ is a prefix of $c^{\prime}$.

Using the same argument as above, we can show that there exists no more than one maximal occurrence in $U_h$ such that it is not separated from the end of $L^{\prime}$ from the left side and its intersection with $L^{\prime}$ is not a small piece. The same from the right side, namely, there exists no more than one maximal occurrence in $U_h$ such that it is not separated from the beginning of $R^{\prime}$ from the right side and its intersection with $R^{\prime}$ is not a small piece.

We already described what happened with maximal occurrences that are separated from the end of $L^{\prime}$. So we have to study what happens with maximal occurrences that are not separated from the end of $L^{\prime}$. There are the following possibilities.
\begin{enumerate}[label=($\text{\alph*}^{\prime}$)]
\item
\label{no_not_separated2}
There does not exist a maximal occurrence in $U_h$ such that it is not separated from the end of $L^{\prime}$ from the left side and such that its intersection with $L^{\prime}$ is not a small piece. And similarly, there does not exist a maximal occurrence in $U_h$ such that it is not separated from the beginning of $R^{\prime}$ from the right side and such that its intersection with $R^{\prime}$ is not a small piece.
\item
\label{not_separated_left2}
There exists a maximal occurrence in $U_h$ such that it is not separated from the end of $L^{\prime}$ from left side and such that its intersection with $L^{\prime}$ is not a small piece. Let this be our $b$. As before, $\widehat{b}$ is its intersection with $L^{\prime}$. At the same time there does not exist a maximal occurrence in $U_h$ such that it is not separated from the beginning of $R^{\prime}$ from the right side and such that its intersection with $R^{\prime}$ is not a small piece.
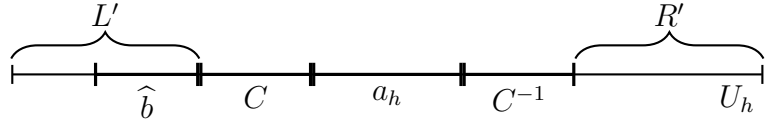
\begin{center}
\begin{tikzpicture}
\draw[|-|, black, thick] (0,0)--(10,0) node[at end, below, xshift=-10] {$U_h$};
\draw [thick, decorate, decoration={brace, amplitude=10pt, raise=6pt}] (0,0) to node[midway, above, yshift=14pt] {$L^{\prime}$} (2.5, 0);
\draw[|-|, black, very thick] (4,0)--(6,0) node[midway, below] {$a_h$};
\draw[|-|, black, very thick] (6,0)--(7.5,0) node[midway, below] {$C^{-1}$};
\draw[|-|, black, very thick] (2.5,0)--(4,0) node[midway, below] {$C$};
\draw[|-|, black, very thick] (1.1,0)--(2.5,0) node[midway, below] {$\widehat{b}$};
\draw [thick, decorate, decoration={brace, amplitude=10pt, raise=6pt}] (7.5,0) to node[midway, above, yshift=14pt] {$R^{\prime}$} (10, 0);
\end{tikzpicture}
\end{center}
\item
\label{not_separated_right2}
There does not exist a maximal occurrence in $U_h$ such that it is not separated from the end of $L^{\prime}$ from the left side and such that its intersection with $L^{\prime}$ is not a small piece. But there exists a maximal occurrence in $U_h$ such that it is not separated from the beginning of $R^{\prime}$ from the right side and such that its intersection with $R^{\prime}$ is not a small piece. Let this be our $c$. As before, $\widehat{c}$ is its intersection with $R^{\prime}$.
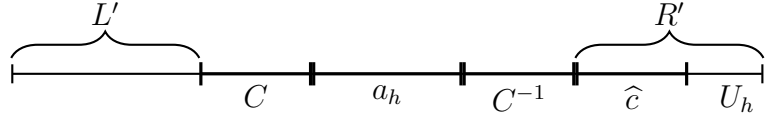
\begin{center}
\begin{tikzpicture}
\draw[|-|, black, thick] (0,0)--(10,0) node[at end, below, xshift=-10] {$U_h$};
\draw [thick, decorate, decoration={brace, amplitude=10pt, raise=6pt}] (0,0) to node[midway, above, yshift=14pt] {$L^{\prime}$} (2.5, 0);
\draw[|-|, black, very thick] (4,0)--(6,0) node[midway, below] {$a_h$};
\draw[|-|, black, very thick] (6,0)--(7.5,0) node[midway, below] {$C^{-1}$};
\draw[|-|, black, very thick] (2.5,0)--(4,0) node[midway, below] {$C$};
\draw [thick, decorate, decoration={brace, amplitude=10pt, raise=6pt}] (7.5,0) to node[midway, above, yshift=14pt] {$R^{\prime}$} (10, 0);
\draw[|-|, black, very thick] (7.5,0)--(9,0) node[midway, below] {$\widehat{c}$};
\end{tikzpicture}
\end{center}
\item
\label{not_separated_left_right2}
There exists a maximal occurrence in $U_h$ such that it is not separated from the end of $L^{\prime}$ from the left side and such that its intersection with $L^{\prime}$ is not a small piece. Let this be $b$. As before, $\widehat{b}$ is its intersection with $L^{\prime}$. But there exists a maximal occurrence in $U_h$ such that it is not separated from the end of $L^{\prime}$ from the right side and such that its intersection with $R^{\prime}$ is not a small piece. Let this be $c$. As before, $\widehat{c}$ is its intersection with $R^{\prime}$.
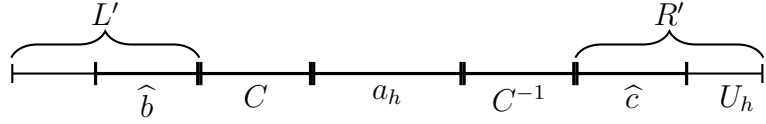
\begin{center}
\begin{tikzpicture}
\draw[|-|, black, thick] (0,0)--(10,0) node[at end, below, xshift=-10] {$U_h$};
\draw [thick, decorate, decoration={brace, amplitude=10pt, raise=6pt}] (0,0) to node[midway, above, yshift=14pt] {$L^{\prime}$} (2.5, 0);
\draw[|-|, black, very thick] (4,0)--(6,0) node[midway, below] {$a_h$};
\draw[|-|, black, very thick] (6,0)--(7.5,0) node[midway, below] {$C^{-1}$};
\draw[|-|, black, very thick] (2.5,0)--(4,0) node[midway, below] {$C$};
\draw [thick, decorate, decoration={brace, amplitude=10pt, raise=6pt}] (7.5,0) to node[midway, above, yshift=14pt] {$R^{\prime}$} (10, 0);
\draw[|-|, black, very thick] (7.5,0)--(9,0) node[midway, below] {$\widehat{c}$};
\draw[|-|, black, very thick] (1.1,0)--(2.5,0) node[midway, below] {$\widehat{b}$};
\end{tikzpicture}
\end{center}
\end{enumerate}
We have already considered case~\ref{no_not_separated2}, so, we have to consider only cases~\ref{not_separated_left2}---\ref{not_separated_left_right2}.

Now let us study $b^{\prime}$ and $c^{\prime}$ in a monomial of type~\ref{donot_keep_structure2} in detail. Consider cases~\ref{not_separated_left2} and~\ref{not_separated_left_right2}. That is, $\widehat{b}$ is a terminal subword of $L^{\prime}$, $L^{\prime} = L_1^{\prime}\widehat{b}$. First we describe possibilities that can be obtained together in cases~\ref{not_separated_left2} and~\ref{not_separated_left_right2}. They do not depend on a configuration of $c^{\prime}$ in the same monomial $U_j$. Then there are the following configurations for $b^{\prime}$ (see Remark~\ref{remark_possible_prolongations}):
\begin{enumerate}[label=L\,\ref{donot_keep_structure2}.\arabic*]
\item
\label{donot_keep_structure2_l1}
$\widehat{b}$ is a maximal occurrence in $U_j$, that is $b^{\prime} = \widehat{b}$.
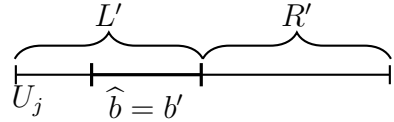
\begin{center}
\begin{tikzpicture}
\draw[|-|, black, thick] (0,0)--(5,0);
\node[below, xshift=5] at (0, 0) {$U_j$};
\draw[|-|, black, very thick] (1,0)--(2.5,0) node[midway, below] {$\widehat{b} = b^{\prime}$};
\draw [thick, decorate, decoration={brace, amplitude=10pt, raise=6pt}] (0, 0) to node[midway, above, yshift=14pt] {$L^{\prime}$} (2.5, 0);
\draw [thick, decorate, decoration={brace, amplitude=10pt, raise=6pt}] (2.5, 0) to node[midway, above, yshift=14pt] {$R^{\prime}$} (5, 0);
\end{tikzpicture}
\end{center}
\item
\label{donot_keep_structure2_l2}
$\widehat{b}$ can be prolonged in $U_j$ by a small piece, that is, $b^{\prime} = \widehat{b}e_1$, where $e_1$ is a non-empty small piece.
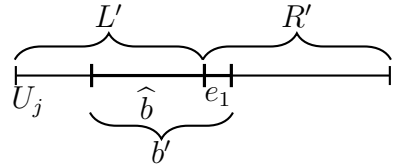
\begin{center}
\begin{tikzpicture}
\draw[|-|, black, thick] (0,0)--(5,0);
\node[below, xshift=5] at (0, 0) {$U_j$};
\draw[|-, black, very thick] (1,0)--(2.5,0) node[midway, below] {$\widehat{b}$};
\draw[|-|, black, very thick] (2.5,0)--(2.9,0) node[midway, below] {$e_1$};
\draw [thick, decorate, decoration={brace, amplitude=10pt, raise=14pt, mirror}] (1, 0) to node[midway, below, yshift=-20pt] {$b^{\prime}$} (2.9, 0);
\draw [thick, decorate, decoration={brace, amplitude=10pt, raise=6pt}] (0, 0) to node[midway, above, yshift=14pt] {$L^{\prime}$} (2.5, 0);
\draw [thick, decorate, decoration={brace, amplitude=10pt, raise=6pt}] (2.5, 0) to node[midway, above, yshift=14pt] {$R^{\prime}$} (5, 0);
\end{tikzpicture}
\end{center}
\end{enumerate}

Let us treat cases~\ref{not_separated_left_right2} and~\ref{not_separated_left_right2}. So $\widehat{c}$ is an initial subword of $R^{\prime}$, $R^{\prime} = \widehat{c}R_1^{\prime}$. We enumerate configurations that can be obtained together in cases~\ref{not_separated_left_right2} and~\ref{not_separated_left_right2}. They do not depend on a configuration of $b^{\prime}$ in the same monomial $U_j$. Then there are the following possibilities for $c^{\prime}$ (see Remark~\ref{remark_possible_prolongations}):
\begin{enumerate}[label=R\,\ref{donot_keep_structure2}.\arabic*]
\item
\label{donot_keep_structure2_r1}
$\widehat{c}$ is a maximal occurrence in $U_j$, that is, $c^{\prime} = \widehat{c}$.
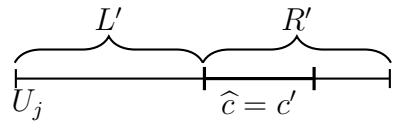
\begin{center}
\begin{tikzpicture}
\draw[|-|, black, thick] (0,0)--(5,0);
\node[below, xshift=5] at (0, 0) {$U_j$};
\draw[|-|, black, very thick] (2.5,0)--(4,0) node[midway, below] {$\widehat{c} = c^{\prime}$};
\draw [thick, decorate, decoration={brace, amplitude=10pt, raise=6pt}] (0, 0) to node[midway, above, yshift=14pt] {$L^{\prime}$} (2.5, 0);
\draw [thick, decorate, decoration={brace, amplitude=10pt, raise=6pt}] (2.5, 0) to node[midway, above, yshift=14pt] {$R^{\prime}$} (5, 0);
\end{tikzpicture}
\end{center}
\item
\label{donot_keep_structure2_r2}
$\widehat{c}$ can be prolonged in $U_j$, that is, $c^{\prime} = e_2\widehat{c}$, where $e_2$ is a non-empty small piece.
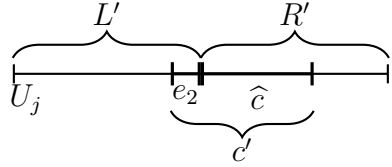
\begin{center}
\begin{tikzpicture}
\draw[|-|, black, thick] (0,0)--(5,0);
\node[below, xshift=5] at (0, 0) {$U_j$};
\draw[|-|, black, very thick] (2.5,0)--(4,0) node[midway, below] {$\widehat{c}$};
\draw[|-|, black, very thick] (2.1,0)--(2.5,0) node[midway, below] {$e_2$};
\draw [thick, decorate, decoration={brace, amplitude=10pt, raise=14pt, mirror}] (2.1, 0) to node[midway, below, yshift=-20pt] {$c^{\prime}$} (4, 0);
\draw [thick, decorate, decoration={brace, amplitude=10pt, raise=6pt}] (0, 0) to node[midway, above, yshift=14pt] {$L^{\prime}$} (2.5, 0);
\draw [thick, decorate, decoration={brace, amplitude=10pt, raise=6pt}] (2.5, 0) to node[midway, above, yshift=14pt] {$R^{\prime}$} (5, 0);
\end{tikzpicture}
\end{center}
\end{enumerate}

Now we will finish with case~\ref{not_separated_left_right2}. So, we assume that $\widehat{b}$ is a terminal subword of $L^{\prime}$, $L^{\prime} = L_1^{\prime}\widehat{b}$, $\widehat{c}$ is an initial subword of $R^{\prime}$, $R^{\prime} = \widehat{c}R_1^{\prime}$. Then we can have any combination of configurations~\ref{donot_keep_structure2_l1}, \ref{donot_keep_structure2_l2} for $b^{\prime}$ and ~\ref{donot_keep_structure2_r1}, \ref{donot_keep_structure2_r2} for $c^{\prime}$ in the monomial $U_j$. In these cases $b^{\prime}$ and $c^{\prime}$ are different maximal occurrences in $U_j$. The last possibility is that $b^{\prime}$ and $c^{\prime}$ is the same maximal occurrence in $U_j$:
\begin{enumerate}[label=\ref{donot_keep_structure2}.\arabic*]
\item
\label{donot_keep_structure2_merge}
$\widehat{b}$ and $\widehat{c}$ merge to one maximal occurrence $\widehat{b}\widehat{c}$, that is, $b^{\prime} = c^{\prime} = \widehat{b}\widehat{c}$.
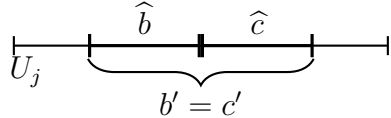
\begin{center}
\begin{tikzpicture}
\draw[|-|, black, thick] (0,0)--(5,0);
\node[below, xshift=5] at (0, 0) {$U_j$};
\draw[|-|, black, very thick] (1,0)--(2.5,0) node[midway, above] {$\widehat{b}$};
\draw[|-|, black, very thick] (2.5,0)--(4,0) node[midway, above] {$\widehat{c}$};
\draw [thick, decorate, decoration={brace, amplitude=10pt, raise=5pt, mirror}] (1, 0) to node[midway, below, yshift=-13pt] {$b^{\prime} = c^{\prime}$} (4, 0);
\end{tikzpicture}
\end{center}
\end{enumerate}


\section{Virtual members of the chart}
\label{virtual_members_section}
Let $U_h$ be a monomial. We denote by $\mo{U_h}$\label{mo_def} the set of all maximal occurrences of monomials of $\Mon$ in $U_h$. Assume $a_h$ is a member of the chart of $U_h$ (see Definition~\ref{chart_def}), $U_h = La_hR$. Consider a multi-turn  $U_h \mapsto \sum_{\substack{j = 1 \\ j\neq h}}^{n}U_j$ that comes from an elementary multi-turn $a_h \mapsto \sum_{\substack{j = 1 \\ j\neq h}}^{n}a_j$. That is,
\begin{equation*}
U_h = La_hR \mapsto \sum_{\substack{j = 1 \\ j\neq h}}^{n}U_j = \sum_{\substack{j = 1 \\ j\neq h}}^{n}La_jR.
\end{equation*}
In Section~\ref{structure_calc} we employ an inductive argument for $U_h$ and $U_j$. For the sake of this argument, we need a numerical parameter of a monomial that does not increase when we go from $U_h$ to $U_j$. The following examples show that neither the number of maximal occurrences in a monomial, nor the number of members of the chart of a monomial do not suite this purpose.

\begin{example}
\label{example_members_not_stable}
Let $a_h$ be a member of the chart of $U_h$. Assume $b, c \in \mo{U_h}$, $b$ starts from the left of the beginning of $a_h$, $c$ starts from the right of the beginning of $a_h$. Assume $b$ and $a_h$ touch at a point, and $a_h$ and $c$ also touch at a point. Let $\SPM(b) < \tau$ and $\SPM(c) < \tau$. Then $b$ and $c$ are not members of the chart of $U_h$. According to the previous section, the element $b$ may be prolonged in $U_j$ and the $\SPM$-measure of the corresponding prolonged element $b^{\prime}$ in $U_j$ may increase and become $\geqslant \tau$. The same may happen with $c$, that is, the element $c$ may be prolonged in $U_j$ and the $\SPM$-measure of the corresponding prolonged element $c^{\prime}$ in $U_j$ may increase and become $\geqslant \tau$. If this happens simultaneously, then both $b^{\prime}$ and $c^{\prime}$ are members of the chart of $U_j$. In this case, the number of members of the chart of $U_j$ becomes greater than the number of members of the chart of $U_h$ even if $\SPM(a_j) < \tau$.
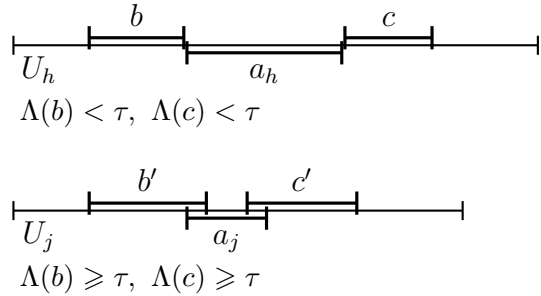
\begin{center}
\begin{tikzpicture}
\node at (0, 0) [below, xshift=10] {$U_h$};
\draw[|-|, black, thick] (0,0)--(7,0);
\draw[|-|, black, very thick] (1,0.1)--(2.3,0.1) node[midway, above] {$b$};
\draw[|-|, black, very thick] (2.3,-0.1)--(4.4,-0.1) node[midway, below] {$a_h$};
\draw[|-|, black, very thick] (4.4,0.1)--(5.6,0.1) node[midway, above] {$c$};
\node[text width=4cm, align=left] at (2.1, -0.9) {\small\baselineskip=10pt \textit{$\SPM(b) < \tau,\ \SPM(c) < \tau$}};

\node at (0, -2.2) [below, xshift=10] {$U_j$};
\draw[|-|, black, thick] (0,-2.2)--(6,-2.2);
\draw[|-|, black, very thick] (1,0.1 - 2.2)--(2.6,0.1-2.2) node[midway, above] {$b^{\prime}$};
\draw[|-|, black, very thick] (2.3,-0.1-2.2)--(3.4,-0.1-2.2) node[midway, below] {$a_j$};
\draw[|-|, black, very thick] (3.1,0.1-2.2)--(4.6,0.1-2.2) node[midway, above] {$c^{\prime}$};
\node[text width=4cm, align=left] at (2.1, -0.9-2.2) {\small\baselineskip=10pt \textit{$\SPM(b) \geqslant \tau,\ \SPM(c) \geqslant \tau$}};
\end{tikzpicture}
\end{center}
\end{example}

\begin{example}
\label{example_mo_not_stable}
Let $a_h$ be a member of the chart of $U_h$. Assume $b\in \mo{U_h}$, $b$ starts from the left of the beginning of $a_h$. Assume $b$ and $a_h$ touch at a point. Assume $b$ is the only element of $\mo{U_h}$ that is not separated from $a_h$ from the left side. Let $b = \widetilde{b}d^{\prime}$, where $d^{\prime}$ is a proper suffix of $b$, $d^{\prime}$ is a small piece. Let $a_j = d^{\prime\prime}\widetilde{a}_j$, where $d^{\prime\prime}$ is a proper prefix of $a_j$, $d^{\prime\prime}$ is a small piece. It may happen that $d = d^{\prime}d^{\prime\prime}$ is a maximal occurrence in $U_j$. If this happens, we obtain a new element in $\fc{U_j}$, which grows from two occurrences $d^{\prime} \notin \mo{U_h}$ and $d^{\prime\prime} \notin \mo{U_j}$. So, assume $d = d^{\prime}d^{\prime\prime}$ is a maximal occurrence in $U_j$. Assume additionally that $b$ stays unchanged in $U_j$, $a_j \in \mo{U_j}$, and there are no elements of $\mo{U_h}$ that are not separated from $a_h$ from the right side. Then $\vert \mo{U_j}\vert > \vert \mo{U_h}\vert$.
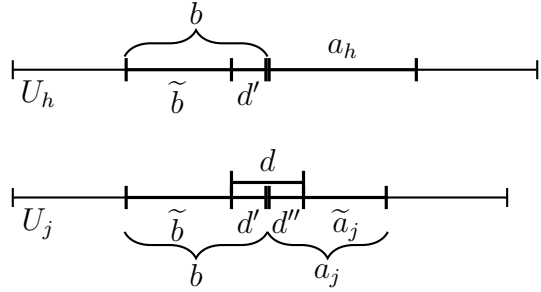
\begin{center}
\begin{tikzpicture}
\node at (0, 0) [below, xshift=10] {$U_h$};
\draw[|-|, black, thick] (0, 0)--(7, 0);
\draw[|-, black, very thick] (1.5, 0)--(2.9, 0) node[midway, below] {$\widetilde{b}$};
\draw[|-|, black, very thick] (2.9, 0)--(3.4, 0) node[midway, below] {$d^{\prime}$};
\draw[|-|, black, very thick] (3.4, 0)--(5.4, 0) node[midway, above] {$a_h$};
\draw [thick, decorate, decoration={brace, amplitude=10pt, raise=5pt}] (1.5, 0) to node[midway, above, yshift=13pt] {$b$} (3.4, 0);

\node at (0, -1.7) [below, xshift=10] {$U_j$};
\draw[|-|, black, thick] (0, -1.7)--(6.6, -1.7);
\draw[|-, black, very thick] (1.5, -1.7)--(2.9, -1.7) node[midway, below] {$\widetilde{b}$};
\draw[|-|, black, very thick] (2.9, -1.7)--(3.4, -1.7) node[midway, below] {$d^{\prime}$};
\draw[|-|, black, very thick] (3.4, -1.7)--(3.9, -1.7) node[midway, below] {$d^{\prime\prime}$};
\draw[-|, black, very thick] (3.9, -1.7)--(5, -1.7) node[midway, below] {$\widetilde{a}_j$};
\draw[|-|, black, very thick] (2.9, 0.2-1.7)--(3.9, 0.2-1.7) node[midway, above] {$d$};
\draw [thick, decorate, decoration={brace, amplitude=10pt, raise=13pt, mirror}] (1.5, -1.7) to node[midway, below, yshift=-21pt] {$b$} (3.4, -1.7);
\draw [thick, decorate, decoration={brace, amplitude=10pt, raise=13pt, mirror}] (3.4, -1.7) to node[midway, below, yshift=-21pt] {$a_j$} (5, -1.7);
\end{tikzpicture}
\end{center}
\end{example}

So, the number of members of the chart of $U_h$ and $\vert \mo{U_h}\vert$ are not appropriate parameters for the induction. If we want to prove that the induction is nevertheless finite, we need to introduce a parameter with more stable properties with respect to multi-turns. In order to do this, we refine the notion of a member of the chart in this section. After that we introduce a function that guarantees finiteness of the inductive process (see Proposition~\ref{non_increasing_parameter}).

\subsection{Images of maximal occurrences}
\begin{definition}
\label{images_def}
Let $U_h$ be a monomial, $a_h \in \mo{U_h}$, $\SPM(a_h) \geqslant 3$, $U_h = La_hR$. Let $a_h$ and $a_j$ be incident monomials (that is, there exists $p \in \Rel$ such that $a_h$ and $a_j$ are monomials in $p$). Consider the replacement $a_h \mapsto a_j$\label{incidents_replacement_def} in $U_h$. Then $La_jR$ is the resulting monomial, we denote it by $U_j$.
\begin{enumerate}
\item
First we define \emph{images of $a_h$ in $U_j$}. If $a_j = 1$  we say that \emph{the set of images of $a_j$ is empty}. If $a_j \neq 1$,
an element of $\mo{U_j}$ that contains $a_j$ is called \emph{an image of $a_h$ in $U_j$}. The set of all such elements is called \emph{the set of images of $a_j$ in $U_j$}.

\item
Assume $a_j \neq 1$. Let $b \in \mo{U_h}$, $b \neq a_h$, $b$ starts from the left of the beginning of $a_h$. Let $\widehat{b}$ be the intersection of $b$ and $L$. Then an element of $\mo{U_j}$ that contains $\widehat{b}$ is called \emph{an image of $b$ in $U_j$}. The set of all such elements is called \emph{the set of images of $b$ in $U_j$}.

Let $c \in \mo{U_h}$, $c \neq a_h$, $c$ starts from the right of the beginning of $a_h$. Let $\widehat{c}$ be the intersection of $c$ and $R$. Then an element of $\mo{U_j}$ that contains $\widehat{c}$ is called \emph{an image of $c$ in $U_j$}. The set of all such elements is called \emph{the set of images of $c$ in $U_j$}.

\item
Let $a_j = 1$. Assume $L = L^{\prime}C$, $R = C^{-1}R^{\prime}$, $L^{\prime}R^{\prime}$ has no further cancellations. Let $b \in \mo{U_h}$ such that $b \neq a_h$, $b$ starts from the left of the beginning of $a_h$. If $b$ is contained in $C$, then $b$ cancels in $U_j$ and we say that \emph{the set of images of $b$ in $U_j$ is empty}. Assume $b$ is not contained in $C$. Let $\widehat{b}$ be the intersection of $b$ and $L^{\prime}$. Then an element of $\mo{U_j}$ that contains $\widehat{b}$ is called \emph{an image of $b$ in $U_j$}. The set of all such elements is called \emph{the set of images of $b$ in $U_j$}.

Let $c \in \mo{U_h}$ such that $c \neq a_h$, $c$ starts from the right of the beginning of $a_h$. If $c$ is contained in $C^{-1}$, then $c$ cancels in $U_j$ and we say that \emph{the set of images of $c$ in $U_j$ is empty}. Assume $c$ is not contained in $C^{-1}$. Let $\widehat{c}$ be the intersection of $c$ and $R^{\prime}$. Then an element of $\mo{U_j}$ that contains $\widehat{c}$ is called \emph{an image of $c$ in $U_j$}. The set of all such elements is called \emph{the set of images of $c$ in $U_j$}.
\end{enumerate}
\end{definition}

\begin{remark}
Let us state several observations on images.
\begin{enumerate}
\item
Consider images of $a_h$ in $U_j$. All possible images of $a_h$ in $U_j$ are described in Section~\ref{mt_configurations}, see cases~\ref{a_j_keep_structure}, \ref{donot_keep_structure1_l3}, \ref{donot_keep_structure1_l4}, \ref{donot_keep_structure1_r3}, \ref{donot_keep_structure1_r4}, \ref{donot_keep_structure1_merge}, \ref{a_j_donot_keep_structure1_1}---\ref{a_j_donot_keep_structure1_4}.

\item
Consider the case $a_j \neq 1$. Let $b \in \mo{U_h}$ such that $b \neq a_h$, $b$ starts from the left of the beginning of $a_h$. Let $\widehat{b}$ be the intersection of $b$ and $L$. Let $c \in \mo{U_h}$ such that $c \neq a_h$, $c$ starts from the right of the beginning of $a_h$. Let $\widehat{c}$ be the intersection of $c$ and $R$.

By Corollary~\ref{max_occurrences_coinside}, if $\widehat{b}$ is not a small piece, then $b$ has a single image in $U_j$. If $\widehat{b}$ is a small piece, then it may have more than one image in $U_j$. Similarly, if $\widehat{c}$ is not a small piece, then $c$ has a single image in $U_j$. If $\widehat{c}$ is a small piece, then it may have more than one image in $U_j$.

Recall that in Section~\ref{mt_configurations} we studied $b$ such that $\widehat{b}$ is not a small piece and $c$ such that $\widehat{c}$ is not a small piece. In fact, we described all possible forms of their images in $U_j$. We proved that if $b$ is separated from $a_h$ in $U_h$, then $\widehat{b} = b$ in $U_h$ and the image of $b$ is the corresponding occurrence of $\widehat{b}$ in $U_j$. If $b$ is not separated from $a_h$ in $U_h$, then all  possible forms of its image are described in Section~\ref{mt_configurations} in~\ref{keep_structure_case_l1}, \ref{keep_structure_case_l2}, \ref{donot_keep_structure1_l1}---\ref{donot_keep_structure1_l4}, \ref{donot_keep_structure1_merge}. Similarly for $c$, we proved that if $c$ is separated from $a_h$ in $U_h$, then $\widehat{c} = c$ in $U_h$ and the image of $c$ is the corresponding occurrence of $\widehat{c}$ in $U_j$. If $c$ is not separated from $a_h$ in $U_h$, then all  possible forms of its image are described in Section~\ref{mt_configurations} in~\ref{keep_structure_case_r1}, \ref{keep_structure_case_r2}, \ref{donot_keep_structure1_r1}---\ref{donot_keep_structure1_r4}, \ref{donot_keep_structure1_merge}.

\item
Consider the case $a_j = 1$, $L = L^{\prime}C$, $R = C^{-1}R^{\prime}$, $L^{\prime}R^{\prime}$ has no further cancellations. Let $b \in \mo{U_h}$ such that $b \neq a_h$, $b$ starts from the left of the beginning of $a_h$, $b$ is not contained in $C$. Assume $\widehat{b}$ is the intersection of $b$ and $L^{\prime}$. Let $c \in \mo{U_h}$ such that $c \neq a_h$, $c$ starts from the right of the beginning of $a_h$, $c$ is not contained in $C^{-1}$. Assume $\widehat{c}$ is the intersection of $c$ and $R^{\prime}$.

By Corollary~\ref{max_occurrences_coinside}, if $\widehat{b}$ is not a small piece, then $b$ has a single image in $U_j$. If $\widehat{b}$ is a small piece, then it may have more than one image in $U_j$. Similarly, if $\widehat{c}$ is not a small piece, then $c$ has a single image in $U_j$. If $\widehat{c}$ is a small piece, then it may have more than one image in $U_j$.

In Section~\ref{mt_configurations} we studied $b$ such that $\widehat{b}$ is not a small piece and $c$ such that $\widehat{c}$ is not a small piece. We described all possible forms of their images in $U_j$. We proved that if $b$ is separated from the end of $L^{\prime}$ in $U_h$, then $\widehat{b} = b$ in $U_h$ and the image of $b$ is the corresponding occurrence of $\widehat{b}$ in $U_j$. If $b$ is not separated from the end of $L^{\prime}$ in $U_h$, then all  possible forms of its image are described in Section~\ref{mt_configurations} in~\ref{donot_keep_structure2_l1}, \ref{donot_keep_structure2_l2}, \ref{donot_keep_structure2_merge}. Similarly for $c$, we proved that if $c$ is separated from the beginning of $R^{\prime}$ in $U_h$, then $\widehat{c} = c$ in $U_h$ and the image of $c$ is the corresponding occurrence of $\widehat{c}$ in $U_j$. If $c$ is not separated from the beginning of $R^{\prime}$ in $U_h$, then all  possible forms of its image are described in Section~\ref{mt_configurations} in~\ref{donot_keep_structure2_r1}, \ref{donot_keep_structure2_r2}, \ref{donot_keep_structure2_merge}.

\item
Notice that taking  images is not an injective mapping. That is, different elements of $\mo{U_h}$ may have the same set of images in $U_j$. For example, let us mention case~\ref{donot_keep_structure1_merge} in Section~\ref{mt_configurations}. It follows directly from our definition of images, that in this case the maximal occurrences $b$, $c$ and $a_h$ in $U_h$ have the same set of images in $U_j$ (and this set consists of the single element $\widehat{b}a_j\widehat{c}$).
\end{enumerate}
\end{remark}

\subsection{Minimal coverings of a monomial}
\label{coverings_subsection}
Let $U$ be a monomial. Assume $A, B_1, \ldots, B_n$ are occurrences in $U$. If every occurrence of a letter in $A$ is contained in some $B_i$, $1 \leqslant i \leqslant n$, then we say that \emph{$B_1, \ldots, B_n$ fully cover $A$} (or for short \emph{cover $A$}). Let us consider the monomial $U$ as a segment and its subwords as its subsegments. Assume $B_1, \ldots, B_n$ fully cover $A$ in the sense of the above definition. Then, clearly, the segments $B_1, \ldots, B_n$ cover the segment $A$.

Recall that $\mo{U}$ is the set of all maximal occurrences in $U$. There are two types of such occurrences: occurrences that are not fully covered by other occurrences of $\mo{U}$ and occurrences that are fully covered by other occurrences of $\mo{U}$. We denote the first set by $\nfc{U}$\label{nfc_def} and the second set by $\fc{U}$\label{fc_def}.

\begin{lemma}
\label{fc_measue}
Let $U$ be a monomial. Assume $a \in \fc{U}$. Then $\SPM(a) \leqslant 2$.
\end{lemma}
\begin{proof}
Let $b, c \in \mo{U_h}$. Assume $b$ starts from the left of the beginning of $a$ and $b$ has a non-empty intersection with $a$. Then $b$ covers some proper initial subword of $a$, since $a$ is not contained in $b$.
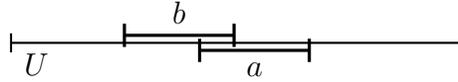
\begin{center}
\begin{tikzpicture}
\draw[|-|, black, thick] (0,0)--(6,0);
\node[below, xshift=10] at (0, 0) {$U$};
\draw[|-|, black, very thick] (1.5,0.1)--(3,0.1) node[midway, above] {$b$};
\draw[|-|, black, very thick] (2.5,-0.1)--(4,-0.1)node[midway, below] {$a$};
\end{tikzpicture}
\end{center}
Assume $c$ starts from the right of the beginning of $a$ and $c$ has a non-empty intersection with $a$. Then $c$ covers some proper final subword of $a$, since $c$ is not contained in $a$.
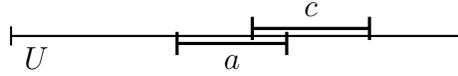
\begin{center}
\begin{tikzpicture}
\draw[|-|, black, thick] (0,0)--(6,0);
\node[below, xshift=10] at (0, 0) {$U$};
\draw[|-|, black, very thick] (2.2,-0.1)--(3.7,-0.1)node[midway, below] {$a$};
\draw[|-|, black, very thick] (3.2,0.1)--(4.8,0.1) node[midway, above] {$c$};
\end{tikzpicture}
\end{center}

Assume $d_1, \ldots, d_n \in \mo{U}$ cover $a$. Clearly, every $d_i$, $1 \leqslant i \leqslant n$, starts either from the left or from the right of the beginning of $a$. Therefore, using the above observations, once can easily see that there exist $d_{i_1}, d_{i_2}$, $1 \leqslant i_1, i_2 \leqslant n$, such that $d_{i_1}$ covers a proper initial subword $a_1$ of $a$, $d_{i_2}$ covers a proper final subword $a_2$ of $a$, $a_1$ and $a_2$ cover the whole $a$.
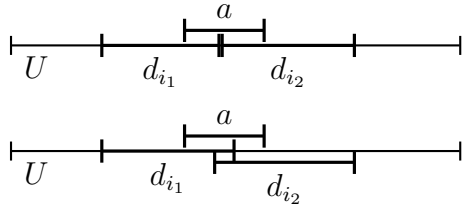
\begin{center}
\begin{tikzpicture}
\node at (-0.5, 0) [below, xshift=10] {$U$};
\draw[|-|, black, thick] (-0.5,0)--(5.5,0);
\draw[|-|, black, very thick] (0.7, 0)--(2.3, 0) node[midway, below] {$d_{i_1}$};
\draw[|-|, black, very thick] (2.3, 0)--(4.1, 0) node[midway, below] {$d_{i_2}$};
\draw[|-|, black, very thick] (1.8, 0.2)--(2.9, 0.2) node[midway, above] {$a$};
\end{tikzpicture}

\begin{tikzpicture}
\node at (-0.5, 0) [below, xshift=10] {$U$};
\draw[|-|, black, thick] (-0.5,0)--(5.5,0);
\draw[|-|, black, very thick] (0.7, 0)--(2.5, 0) node[midway, below] {$d_{i_1}$};
\draw[|-|, black, very thick] (2.2, -0.15)--(4.1, -0.15) node[midway, below] {$d_{i_2}$};
\draw[|-|, black, very thick] (1.8, 0.2)--(2.9, 0.2) node[midway, above] {$a$};
\end{tikzpicture}
\end{center}
Since $a_1$ and $a_2$ are overlaps of maximal occurrences in $U$, $a_1$ and $a_2$ are small pieces. Since $a_1$ and $a_2$ cover $a$, we obtain $\SPM(a) \leqslant 2$.
\end{proof}

Let us enumerate the beginnings of all the elements of $\mo{U}$ according to their positions in $U$ in ascending order. We consider this order on the beginnings as an order on $\mo{U}$. Notice that since elements of $\mo{U}$ are not contained inside each other, the order on their ends is the same as the order on their beginnings. Namely, assume $A_1B_1$ and $A_2B_2$ are two segments of our set, where $A_1$ and $A_2$ are starting points, $B_1$ and $B_2$ are ending points. If $A_1 < A_2$, then $B_1 < B_2$.
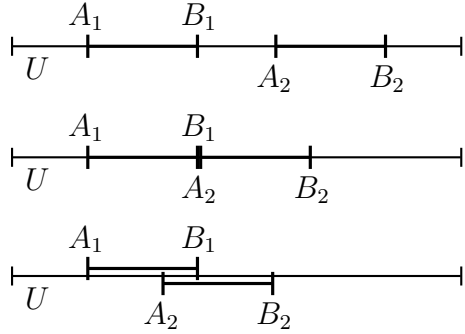
\begin{center}
\begin{tikzpicture}
\draw[|-|, black, thick] (0,0)--(6,0);
\node[below, xshift=10] at (0, 0) {$U$};
\draw[|-|, black, very thick] (1,0)--(2.5,0);
\node[above, yshift=2] at (1, 0) {$A_1$};
\node[above, yshift=2] at (2.5, 0) {$B_1$};
\draw[|-|, black, very thick] (3.5,0)--(5,0);
\node[below, yshift=-3] at (3.5, 0) {$A_2$};
\node[below, yshift=-3] at (5, 0) {$B_2$};
\end{tikzpicture}

\begin{tikzpicture}
\draw[|-|, black, thick] (0,0)--(6,0);
\node[below, xshift=10] at (0, 0) {$U$};
\draw[|-|, black, very thick] (1,0)--(2.5,0);
\node[above, yshift=2] at (1, 0) {$A_1$};
\node[above, yshift=2] at (2.5, 0) {$B_1$};
\draw[|-|, black, very thick] (2.5,0)--(4,0);
\node[below, yshift=-3] at (2.5, 0) {$A_2$};
\node[below, yshift=-3] at (4, 0) {$B_2$};
\end{tikzpicture}

\begin{tikzpicture}
\draw[|-|, black, thick] (0,0)--(6,0);
\node[below, xshift=10] at (0, 0) {$U$};
\draw[|-|, black, very thick] (1,0.1)--(2.5,0.1);
\node[above, yshift=2] at (1, 0.1) {$A_1$};
\node[above, yshift=2] at (2.5, 0.1) {$B_1$};
\draw[|-|, black, very thick] (2,-0.1)--(3.5,-0.1);
\node[below, yshift=-3] at (2, -0.1) {$A_2$};
\node[below, yshift=-3] at (3.5, -0.1) {$B_2$};
\end{tikzpicture}
\end{center}

Let us consider a graph with a set of vertices $\mo{U}$. Two vertices are connected with an edge in this graph if and only if the corresponding maximal occurrences in $U$ are not separated. Let $\mo{U}_1, \ldots, \mo{U}_{m}$ be all the maximal connected components of this graph. Then, clearly,
\begin{equation*}
\mo{U} = \mo{U}_1 \sqcup \ldots \sqcup \mo{U}_{m},
\end{equation*}
where $\sqcup$ is a disjoint union. Assume $a \in \mo{U}_{i_1}$ and $b \in \mo{U}_{i_2}$, $i_1 \neq i_2$. Then $a$ and $b$ are separated in $U$. If $a < b$, then without loss of generality we can assume that $i_1 < i_2$.
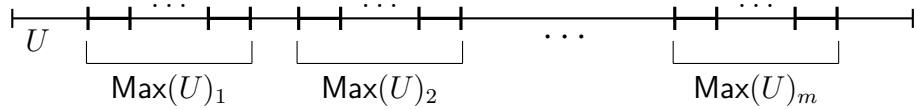
\begin{center}
\begin{tikzpicture}
\draw[|-|, black, thick] (0,0) to node[at start, below, xshift=10] {$U$} (12,0);
\draw[|-|, black, very thick] (1,0)--(1.6,0);
\path (1.6, 0)--(2.6, 0) node[midway, above] {$\ldots$};
\draw[|-|, black, very thick] (2.6,0)--(3.2,0);
\draw (1, -0.3)--(1, -0.6);
\draw (3.2, -0.3)--(3.2, -0.6);
\draw (1, -0.6)--(3.2, -0.6) node[midway, below] {$\mo{U}_1$};
\draw[|-|, black, very thick] (3.8,0)--(4.4,0);
\path (4.4, 0)--(5.4, 0) node[midway, above] {$\ldots$};
\draw[|-|, black, very thick] (5.4,0)--(6,0);
\draw (3.8, -0.3)--(3.8, -0.6);
\draw (6, -0.3)--(6, -0.6);
\draw (3.8, -0.6)--(6, -0.6) node[midway, below] {$\mo{U}_2$};
\path (6, 0)--(8.8, 0) node[midway, below] {\large $\cdots$};
\draw[|-|, black, very thick] (8.8,0)--(9.4,0);
\path (9.4, 0)--(10.4, 0) node[midway, above] {$\ldots$};
\draw[|-|, black, very thick] (10.4,0)--(11,0);
\draw (8.8, -0.3)--(8.8, -0.6);
\draw (11, -0.3)--(11, -0.6);
\draw (8.8, -0.6)--(11, -0.6) node[midway, below] {$\mo{U}_{m}$};
\end{tikzpicture}
\end{center}

Let us consider a set $\mo{U}_k$, $1 \leqslant k \leqslant m$, in detail. Every element of $\mo{U}_k$ belongs either to $\nfc{U}$, or to $\fc{U}$. Clearly, the first element and the last element of $\mo{U}_k$ belong to $\nfc{U}$. So, if we enumerate elements of $\mo{U}_k$ in ascending order, we obtain
\begin{align*}
&d_1, c_{11}, \ldots, c_{1t_1}, d_2, \ldots, d_{l - 1}, c_{l-1,1}, \ldots, c_{l-1, t_{l - 1}}, d_l,\\
&\textit{where } d_i \in \nfc{U}, c_{ij} \in \fc{U}, j = 1, \ldots, t_i, i = 1, \ldots, l,\\
&\textit{a sequence } c_{i1}, \ldots, c_{it_i} \textit{ can be empty}.
\end{align*}
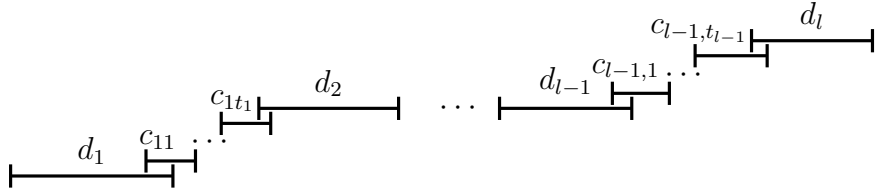
\begin{center}
\begin{tikzpicture}
\draw[|-|, black, very thick] (0, 0)--(2.2, 0) node[midway, above] {$d_1$};
\draw[|-|, black, very thick] (1.8, 0.2)--(2.5, 0.2) node[near start, above] {$c_{11}$};
\path[|-|, black, very thick] (2.4, 0.45)--(3, 0.45) node[midway] {$\cdots$};
\draw[|-|, black, very thick] (2.8, 0.7)--(3.5, 0.7) node[near start, above] {$c_{1t_1}$};
\draw[|-|, black, very thick] (3.3, 0.9)--(5.2, 0.9) node[midway, above] {$d_2$};

\path[|-|, black, very thick] (5.2, 0.9)--(6.8, 0.9) node[midway] {$\cdots$};

\draw[|-|, black, very thick] (6.5, 0.9)--(8.3, 0.9) node[midway, above] {$d_{l - 1}$};
\draw[|-|, black, very thick] (8, 1.1)--(8.8, 1.1) node[near start, above] {$c_{l-1,1}$};
\path[|-|, black, very thick] (8.8, 1.35)--(9.2, 1.35) node[midway] {$\cdots$};
\draw[|-|, black, very thick] (9.1, 1.6)--(10.1, 1.6) node[near start, above, xshift=-5] {$c_{l-1,t_{l-1}}$};
\draw[|-|, black, very thick] (9.85, 1.8)--(11.5, 1.8) node[midway, above] {$d_l$};
\end{tikzpicture}
\end{center}
The elements $d_i$ and $d_{i + 1}$ can either be separated, or touch at a point, or have an overlap. The elements $d_i$ and $c_{i1}$, $c_{1k_i}$ and $d_{i + 1}$, $c_{ij}$ and $c_{i,j+1}$ always have an overlap. Non-successive elements of $\mo{U}_k$ are not necessarily separated. They still can touch at a point or have an overlap. For example,
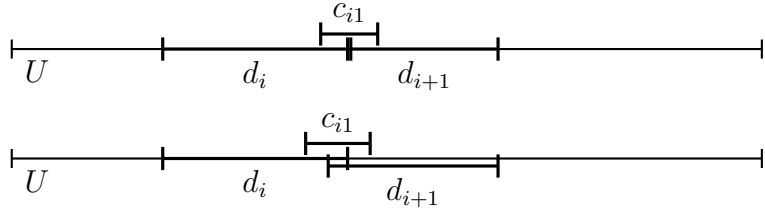
\begin{center}
\begin{tikzpicture}
\draw[|-|, black, thick] (0,0) to node[at start, below, xshift=10] {$U$} (10,0);
\draw[|-|, black, very thick] (2, 0)--(4.5, 0) node[midway, below] {$d_i$};
\draw[|-|, black, very thick] (4.5, 0)--(6.5, 0) node[midway, below] {$d_{i+1}$};
\draw[|-|, black, very thick] (4.1, 0.2)--(4.9, 0.2) node[midway, above] {$c_{i1}$};
\end{tikzpicture}

\begin{tikzpicture}
\draw[|-|, black, thick] (0,0) to node[at start, below, xshift=10] {$U$} (10,0);
\draw[|-|, black, very thick] (2, 0)--(4.5, 0) node[midway, below] {$d_i$};
\draw[|-|, black, very thick] (4.2, -0.1)--(6.5, -0.1) node[midway, below] {$d_{i+1}$};
\draw[|-|, black, very thick] (3.9, 0.2)--(4.8, 0.2) node[midway, above] {$c_{i1}$};
\end{tikzpicture}
\end{center}

Let $P_k$ be the starting point of the first element of $\mo{U}_k$ and $Q_k$ be the end point of the last element of $\mo{U}_k$, $k = 1, \ldots, m$. Clearly, segments $P_{k_1}Q_{k_1}$ and $P_{k_2}Q_{k_2}$ are separated if $k_1 \neq k_2$.
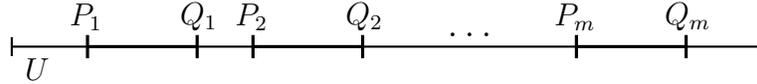
\begin{center}
\begin{tikzpicture}
\draw[|-|, black, thick] (0,0) to node[at start, below, xshift=10] {$U$} (10,0);
\draw[|-|, black, very thick] (1,0)--(2.5,0);
\node[above, yshift=2] at (1, 0) {$P_1$};
\node[above, yshift=2] at (2.5, 0) {$Q_1$};
\draw[|-|, black, very thick] (3.2,0)--(4.7,0);
\node[above, yshift=2] at (3.2, 0) {$P_2$};
\node[above, yshift=2] at (4.7, 0) {$Q_2$};
\path (4.7, 0) to node[midway, above] {\large $\ldots$} (7.5, 0);
\draw[|-|, black, very thick] (7.5,0)--(9,0);
\node[above, yshift=2] at (7.5, 0) {$P_m$};
\node[above, yshift=2] at (9, 0) {$Q_m$};
\end{tikzpicture}
\end{center}
In the other words, if we consider elements of $\mo{U}$ as segments and glue all non-separated segments, as a result we obtain the segments $P_1Q_1, \ldots, P_mQ_m$. So, $\mo{U}$ covers all the segments $P_1Q_1, \ldots, P_mQ_m$ (and only them) in $U$. We consider a subcovering of every segment $P_kQ_k$, $k = 1, \ldots, m$, by elements of $\mo{U}$. Clearly, a subcovering of $P_kQ_k$ is a subset of $\mo{U}_k$. We call the union of these subcoverings by $k = 1, \ldots, m$ \emph{a covering of $U$} (despite the fact that it covers only certain subsegments of $U$ and may not cover the whole $U$). There are finitely many different coverings of $U$, we denote them by $\mathcal{C}_i(U)$, $i = 1, \ldots , n(U)$\label{cov_def}.

If $Z$ is some occurrence in $U$ ($Z$ not necessarily belongs to $\Mon$) and $\mathcal{C}_i(U)$ is a covering of $U$, we can consider all the elements of $\mathcal{C}_i(U)$ that have a non-empty intersection with $Z$ (but not necessarily fully contained in $Z$). We denote this set by $\mathcal{C}_i(Z, U)$\label{cov_intersection_def}.

Let us consider a subcovering of every segment $P_kQ_k$, $k = 1, \ldots, m$, that consists of the smallest number of elements and denote this number by $\mincovk{U}$. We call the union of these subcoverings by $k = 1, \ldots, m$ \emph{a minimal covering of $U$}. Clearly, a minimal covering of $U$ is not necessarily unique. Let $\mincov{U}$ be the number of elements in a minimal covering of $U$, that is,
\begin{equation*}
\label{mincov_def}
\mincov{U} = \sum\limits_{k = 1}^{m}\mincovk{U}.
\end{equation*}
Notice that all elements of $\nfc{U}$ are contained in any covering of $U$, including a minimal covering. If some element of $\fc{U}$ is fully covered by elements of $\nfc{U}$, it is never contained in a minimal covering of $U$.

Assume $U_h$ is a monomial, $a_h \in \mo{U_h}$. Let $a_h$ and $a_j$ be incident monomials and consider the replacement $a_h \mapsto a_j$ in $U_h$. Let $U_j$ be the resulting monomial. We study how a minimal covering of $U_h$ is connected with a minimal covering of $U_j$. Let us prove the following auxiliary lemma.
\begin{lemma}
\label{single_image_lemma}
Let $U_h$ be a monomial, $a_h \in \mo{U_h}$, $U_h = La_hR$. Let $a_h$ and $a_j$ be incident monomials, $a_j \neq 1$. Consider the replacement $a_h \mapsto a_j$ in $U_h$. Let $U_j = La_jR$ be the resulting monomial. Assume $b \in \nfc{U_h}$, $b \neq a_h$, $b$ starts from the left of the beginning of $a_h$. Then $b$ has a single image in $U_j$. Moreover, if $\widehat{b}$ is the intersection of $b$ and $L$, then the single image of $b$ in $U_j$ is of the form $\widehat{b}Y$, where $Y$ is a suffix of the image.

Under the same conditions assume $c \in \nfc{U_h}$, $c \neq a_h$, $c$ starts from the right of the beginning of $a_h$. Then $c$ has a single image in $U_j$. Moreover, if $\widehat{c}$ is the intersection of $c$ and $R$, then the single image of $c$ in $U_j$ is of the form $X\widehat{b}$, where $X$ is a prefix of the image.
\end{lemma}
\begin{proof}
We consider the case of $b$ that starts from the left of the beginning of $a_h$. Since $a_j \neq 1$, the monomial $La_jR$ has no cancellations and $b$ has a non-empty set of images in $U_j$. If $\widehat{b}$ is not a small piece, then, by Corollary~\ref{max_occurrences_coinside}, $b$ has a single image without any additional conditions. So, in fact, the condition $b \in \nfc{U_h}$ is essential only if $\widehat{b}$ is a small piece.

If $b$ is separated from $a_h$ in $U_h$, then the corresponding occurrence of $\widehat{b}$ in $U_j$ is a maximal occurrence in $U_j$. Hence, it is the single image of $b$ in $U_j$.

Assume $b$ is not separated from $a_h$ in $U_h$. Then $L = L_1\widehat{b}$. Assume $d^{\prime}$ is an image of $b$ in $U_j$. That is, $d^{\prime} \in \mo{U_j}$, $d^{\prime}$ contains $\widehat{b}$ in $U_j$. Let us show that $\widehat{b}$ is a prefix of $d^{\prime}$, that is, $d^{\prime} = \widehat{b}Y$, where $Y$ is a suffix of $d^{\prime}$. Assume the contrary, namely that $d^{\prime} = X\widehat{b}Y$, where $X$ is a non-empty prefix of $d^{\prime}$. Assume $L = L_1^{\prime}X\widehat{b}$.
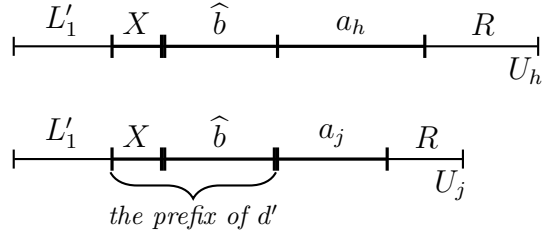
\begin{center}
\begin{tikzpicture}
\draw[|-|, black, thick] (-1,0)--(6,0);
\node[below, xshift=-5] at (6, 0) {$U_h$};
\draw[|-|, black, very thick] (2.5,0)--(4.5,0) node[midway, above] {$a_h$};
\draw[|-, black, very thick] (1,0)--(2.5,0) node[midway, above] {$\widehat{b}$};
\draw[|-|, black, very thick] (0.3,0)--(1,0) node[midway, above] {$X$};
\path (-1,0)--(0.3,0) node[midway, above] {$L_1^{\prime}$};
\path (4.5,0)--(6,0) node[midway, above] {$R$};

\draw[|-|, black, thick] (-1,-1.5)--(5,-1.5);
\node[below, xshift=-5] at (5, -1.5) {$U_j$};
\path (-1,-1.5)--(0.3,-1.5) node[midway, above] {$L_1^{\prime}$};
\draw[|-|, black, very thick] (2.5,-1.5)--(4,-1.5) node[midway, above] {$a_j$};
\draw[|-|, black, very thick] (1,-1.5)--(2.5,-1.5) node[midway, above] {$\widehat{b}$};
\draw[|-|, black, very thick] (0.3,-1.5)--(1,-1.5) node[midway, above] {$X$};
\path (4,-1.5)--(5,-1.5) node[midway, above] {$R$};
\draw [thick, decorate, decoration={brace, amplitude=10pt, raise=5pt, mirror}] (0.3, -1.5) to (2.5, -1.5);
\node[text width=3cm, align=center, below, yshift=-13] at (1.4, -1.5) {\footnotesize{\baselineskip=10pt\textit{the prefix of $d^{\prime}$}\par}};
\end{tikzpicture}
\end{center}

Consider the case when $a_h$ and $b$ touch at a point. Then $\widehat{b} = b$ in $U_h$ and $L = L_1b$.
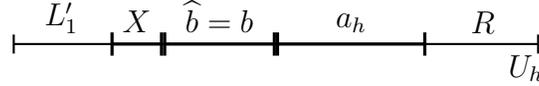
\begin{center}
\begin{tikzpicture}
\draw[|-|, black, thick] (-1,0)--(6,0);
\node[below, xshift=-5] at (6, 0) {$U_h$};
\path (-1,0)--(0.3,0) node[midway, above] {$L_1^{\prime}$};
\draw[|-|, black, very thick] (2.5,0)--(4.5,0) node[midway, above] {$a_h$};
\draw[|-|, black, very thick] (1,0)--(2.5,0) node[midway, above] {$\widehat{b} = b$};
\draw[|-|, black, very thick] (0.3,0)--(1,0) node[midway, above] {$X$};
\path (4.5,0)--(6,0) node[midway, above] {$R$};
\end{tikzpicture}
\end{center}
We assumed $d^{\prime} = X\widehat{b}Y \in \Mon$, therefore $X\widehat{b} \in \Mon$. Since $\widehat{b} = b$ in $U_h$, we obtain $Xb \in \Mon$. So, since $U_h = L_1^{\prime}Xba_hR$, $b$ is not a maximal occurrence in $U_h$, a contradiction.

Assume $a_h$ and $b$ have an overlap $e$ in $U_h$. Since $X\widehat{b} \in \Mon$, there exists a maximal occurrence in $U_h$ that contains $X\widehat{b}$. Denote this maximal occurrence in $U_h$ by $d$, so $d \in \mo{U_h}$.
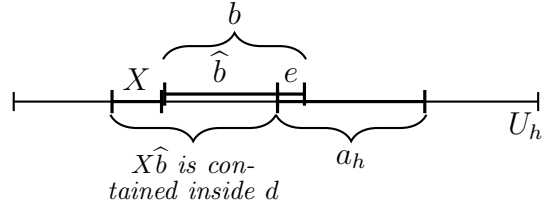
\begin{center}
\begin{tikzpicture}
\draw[|-|, black, thick] (-1,0)--(6,0);
\node[below, xshift=-5] at (6, 0) {$U_h$};
\draw[|-|, black, very thick] (2.5,0)--(4.5,0);
\draw[|-, black, very thick] (1,0.1)--(2.5,0.1) node[midway, above] {$\widehat{b}$};
\draw[|-|, black, very thick] (2.5,0.1)--(2.9,0.1) node[midway, above] {$e$};
\draw [thick, decorate, decoration={brace, amplitude=10pt, raise=16pt}] (1, 0) to node[midway, above, yshift=25pt] {$b$} (2.9, 0);
\draw [thick, decorate, decoration={brace, amplitude=10pt, raise=5pt, mirror}] (2.5, 0) to node[midway, below, yshift=-14pt] {$a_h$} (4.5, 0);
\draw[|-|, black, very thick] (0.3,0)--(1,0) node[midway, above] {$X$};
\draw [thick, decorate, decoration={brace, amplitude=10pt, raise=5pt, mirror}] (0.3, 0) to (2.5, 0);
\node[text width=3cm, align=center, below, yshift=-13] at (1.4, 0) {\footnotesize{\baselineskip=10pt\textit{$X\widehat{b}$ is contained inside $d$}\par}};
\end{tikzpicture}
\end{center}
We assumed that $X$ is non empty, hence, $d$ and $b$ have different beginnings in $U_h$. Therefore $d$ and $b$ are different maximal occurrences in $U_h$. Then $d$ together with $a_h$ cover $b$ in $U_h$. Hence, $b \in \fc{U_h}$, a contradiction.

So finally we obtain $d^{\prime} = \widehat{b}Y$, where $d^{\prime} \in \mo{U_j}$ is an image of $b$ in $U_j$. Since $d^{\prime}$ is an arbitrary image of $b$ in $U_j$, this means that all images of $b$ in $U_j$ have the same beginning. Moreover, every image of $b$ in $U_j$ has the prefix $\widehat{b}$. Clearly, there exists a single maximal occurrence in $U_j$ with the prefix $\widehat{b}$. Thus, $b$ has a single image in $U_j$. The case when $b$ starts from the left of the beginning of $a_h$ in $U_h$ is considered in the same way.
\end{proof}

\begin{lemma}
\label{minimal_coverings_property}
Let $U_h$ be a monomial, $a_h \in \nfc{U_h}$, $U_h = La_hR$. Assume $a_h$ and $a_j$ are incident monomials, $U_j = La_jR$. Then we have $\mincov{U_j)} \leqslant \mincov{U_h}$. Moreover, if $a_j = 1$ or $a_j$ is fully covered by images of  $\nfc{U_h} \setminus \lbrace a_h\rbrace$ in $U_h$, then $\mincov{U_j} < \mincov{U_h}$.
\end{lemma}
\begin{proof}
Let $\mathcal{C}_i(U_h)$ be a minimal covering of $U_h$. Since $a_h \in \nfc{U_h}$, it necessarily belongs to $\mathcal{C}_i(U_h)$. Hence, the covering $\mathcal{C}_i(U_h)$ can be written as
\begin{equation*}
\mathcal{C}_i(U_h) = \mathcal{C}_i(L, U_h) \sqcup \lbrace a_h\rbrace \sqcup \mathcal{C}_i(R, U_h),
\end{equation*}
where $\sqcup$ is a disjoint union. First consider the case $a_j \neq 1$. Let $a_j^{\prime}$ be some image of $a_h$. Consider images of elements of $\mathcal{C}_i(L, U_h)$. Let us take one arbitrary image of every element of $\mathcal{C}_i(L, U_h)$ and denote the obtained set by $\mathcal{C}^{\prime}$. Similarly for $\mathcal{C}_i(R, U_h)$, let us take one arbitrary image of every element of $\mathcal{C}_i(L, U_h)$ and denote the obtained set by $\mathcal{C}^{\prime\prime}$.

As above, let $P_1Q_1, \ldots, P_mQ_m$ be the subsegments of $U_h$ that are actually covered by $\mathcal{C}_i(U_h)$. Recall that if we consider elements of $\mo{U_h}$ as segments and glue all non-separated segments, we obtain the segments $P_1Q_1, \ldots, P_mQ_m$. Then $a_h$ is a subsegment of some $P_tQ_t$, $1 \leqslant t \leqslant m$.
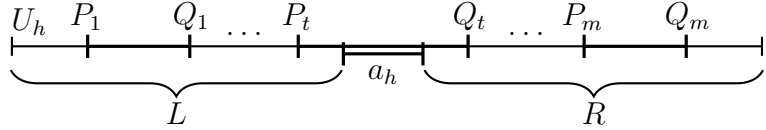
\begin{center}
\begin{tikzpicture}
\draw[|-|, black, thick] (0,0) to node[at start, above, xshift=7] {$U_h$} (10,0);
\draw[|-|, black, very thick] (1,0)--(2.4,0);
\node[above, yshift=2] at (1, 0) {$P_1$};
\node[above, yshift=2] at (2.4, 0) {$Q_1$};
\path (2.4, 0) to node[midway, above] {$\ldots$} (3.8, 0);
\draw[|-|, black, very thick] (3.8,0)--(6.1,0);
\node[above, yshift=2] at (3.8, 0) {$P_t$};
\node[above, yshift=2] at (6.1, 0) {$Q_t$};
\draw[|-|, black, very thick] (4.4,-0.1)--(5.5,-0.1) node[midway, below] {$a_h$};
\path (6.1, 0) to node[midway, above] {$\ldots$} (7.6, 0);
\draw[|-|, black, very thick] (7.6,0)--(9,0);
\node[above, yshift=2] at (7.6, 0) {$P_m$};
\node[above, yshift=2] at (9, 0) {$Q_m$};
\draw [thick, decorate, decoration={brace, amplitude=10pt, raise=9pt, mirror}] (0, 0) to node[midway, below, yshift=-17pt] {$L$} (4.4, 0);
\draw [thick, decorate, decoration={brace, amplitude=10pt, raise=9pt, mirror}] (5.5, 0) to node[midway, below, yshift=-17pt] {$R$} (10, 0);
\end{tikzpicture}
\end{center}
So, in $U_j$ we obtain
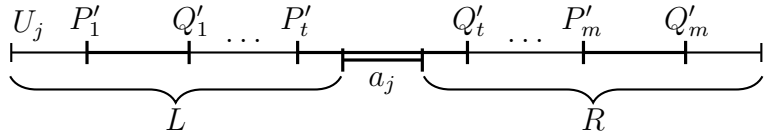
\begin{center}
\begin{tikzpicture}
\draw[|-|, black, thick] (0,0) to node[at start, above, xshift=7] {$U_j$} (10,0);
\draw[|-|, black, very thick] (1,0)--(2.4,0);
\node[above, yshift=2] at (1, 0) {$P_1^{\prime}$};
\node[above, yshift=2] at (2.4, 0) {$Q_1^{\prime}$};
\path (2.4, 0) to node[midway, above] {$\ldots$} (3.8, 0);
\draw[|-|, black, very thick] (3.8,0)--(6.1,0);
\node[above, yshift=2] at (3.8, 0) {$P_t^{\prime}$};
\node[above, yshift=2] at (6.1, 0) {$Q_t^{\prime}$};
\draw[|-|, black, very thick] (4.4,-0.1)--(5.5,-0.1) node[midway, below] {$a_j$};
\path (6.1, 0) to node[midway, above] {$\ldots$} (7.6, 0);
\draw[|-|, black, very thick] (7.6,0)--(9,0);
\node[above, yshift=2] at (7.6, 0) {$P_m^{\prime}$};
\node[above, yshift=2] at (9, 0) {$Q_m^{\prime}$};
\draw [thick, decorate, decoration={brace, amplitude=10pt, raise=9pt, mirror}] (0, 0) to node[midway, below, yshift=-17pt] {$L$} (4.4, 0);
\draw [thick, decorate, decoration={brace, amplitude=10pt, raise=9pt, mirror}] (5.5, 0) to node[midway, below, yshift=-17pt] {$R$} (10, 0);
\end{tikzpicture}
\end{center}
Here the points $P_1^{\prime}, Q_1^{\prime}, \ldots, P_t^{\prime}$ have the same positions in $L$ in $U_j$ as the points $P_1, Q_1, \ldots ,P_t$ have in $L$ in $U_h$, respectively. The points $Q_t^{\prime},P_{t+1}^{\prime}, \ldots, Q_m^{\prime}$ have the same positions in $R$ in $U_j$ as the points $Q_t,P_{t+1},\ldots, Q_m$ have in $R$ in $U_h$, respectively. Clearly, the replacement $a_h \mapsto a_j$ in $U_h$ can not produce any occurrences of monomials of $\Mon$ in $U_j$ outside of the segments $P_1^{\prime}Q_1^{\prime}, \ldots, P_m^{\prime}Q_m^{\prime}$. Hence, a covering of $U_j$ covers the segments $P_1^{\prime}Q_1^{\prime}, \ldots, P_m^{\prime}Q_m^{\prime}$ and only them. Let $\mathcal{C}_l(U_j)$ be a covering of $U_j$. Then we actually proved that an occurrence of a letter in $L$ or in $R$ in $U_h$ is covered by $\mathcal{C}_i(U_h)$ if and only if the corresponding occurrence of this letter in $U_j$ is covered by $\mathcal{C}_l(U_j)$.

Let $x \in S \cup S^{-1}$ be an occurrence of a letter in $L$. Assume $x$ is covered by some $b \in \mathcal{C}_i(L, U_h)$, then the corresponding occurrence of $x$ in $U_j$ is contained in the intersection of $b$ and $L$. Hence, by the definition of images, the corresponding occurrence of $x$ in $U_j$ is covered by any image of $b$ in $U_j$. Hence, $x$ is covered by $\mathcal{C}^{\prime}$ in $U_j$. The same for $R$. Namely, let $y \in S \cup S^{-1}$ be an occurrence of a letter in $R$. Assume $y$ is covered by some $c \in \mathcal{C}_i(R, U_h)$, then the corresponding occurrence of $y$ in $U_j$ is contained in the intersection of $c$ and $R$. Hence, by the definition of images, the corresponding occurrence of $y$ in $U_j$ is covered by any image of $c$ in $U_j$. Hence, $y$ is covered by $\mathcal{C}^{\prime\prime}$ in $U_j$. Therefore, every letter of $L$ that is covered by $\mathcal{C}_i(L, U_h)$ in $U_h$ is covered by $\mathcal{C}^{\prime}$ in $U_j$ and every letter of $R$ that is covered by $\mathcal{C}_i(R, U_h)$ in $U_h$ is covered by $\mathcal{C}^{\prime\prime}$ in $U_j$. The occurrence $a_j$ in $U_j$ is covered by $a_j^{\prime}$. Hence,
\begin{equation*}
\mathcal{C}^{\prime} \cup \lbrace a_j^{\prime}\rbrace \cup \mathcal{C}^{\prime\prime} \textit{ is a covering of } U_j.
\end{equation*}
Since we take one image of every element of $\mathcal{C}_i(L, U_h)$ and $\mathcal{C}_i(R, U_h)$,
\begin{equation*}
\vert \mathcal{C}^{\prime}\vert \leqslant \vert \mathcal{C}_i(L, U_h)\vert, \vert \mathcal{C}^{\prime\prime}\vert \leqslant \vert \mathcal{C}_i(R, U_h)\vert,
\end{equation*}
where $\vert \cdot \vert$ is the number of elements in a set. Hence,
\begin{multline*}
\vert \mathcal{C}^{\prime} \cup \lbrace a_j^{\prime}\rbrace \cup \mathcal{C}^{\prime\prime}\vert \leqslant \vert \mathcal{C}^{\prime}\vert + \vert \mathcal{C}^{\prime\prime}\vert + \vert \lbrace a_j^{\prime}\rbrace\vert \leqslant \vert \mathcal{C}_i(L, U_h)\vert + \vert \mathcal{C}_i(R, U_h)\vert + \vert \lbrace a_j\rbrace\vert  = \\
= \vert\mathcal{C}_i(L, U_h)\vert + \vert \mathcal{C}_i(R, U_h)\vert + 1 = \vert \mathcal{C}_i(U_h) \vert = \mincov{U_h}.
\end{multline*}
So, we constructed a covering of $U_j$ that consists of no more than $\mincov{U_h}$ elements. Thus, $\mincov{U_j} \leqslant \mincov{U_h}$.

Assume $a_j$ is fully covered by images of elements of $\nfc{U_h} \setminus \lbrace a_h\rbrace$ in $U_j$. Every element of $\nfc{U_h}$ is contained in any covering of $U_h$. Clearly, every element of $\nfc{U_h} \setminus \lbrace a_h\rbrace$ has a non empty intersection either with $L$, or with $R$. Therefore, every element of $\nfc{U_h} \setminus \lbrace a_h\rbrace$ is contained either in $\mathcal{C}_i(L, U_h)$, or in $\mathcal{C}_i(R, U_h)$. By Lemma~\ref{single_image_lemma}, every element of $\nfc{U_h} \setminus \lbrace a_h\rbrace$ has only one image. Hence, all images of $\nfc{U_h} \setminus \lbrace a_h\rbrace$ are contained in $\mathcal{C}^{\prime} \cup \mathcal{C}^{\prime\prime}$ (regardless of images that we take for elements that have more than one image). So, $\mathcal{C}^{\prime} \cup \mathcal{C}^{\prime\prime}$ covers $a_j$. Clearly, $a_j^{\prime} = a^{\prime}a_ja^{\prime\prime}$, where $a^{\prime}$ is a suffix of $L$, $a^{\prime\prime}$ is a prefix of $R$ ($a^{\prime}$, $a^{\prime\prime}$ are possibly empty). Since $\mathcal{C}^{\prime}$ covers $L$, we see that $\mathcal{C}^{\prime}$ covers $a^{\prime}$. Since $\mathcal{C}^{\prime\prime}$ covers $R$, we see that $\mathcal{C}^{\prime\prime}$ covers $a^{\prime\prime}$. Therefore, finally we obtain  $\mathcal{C}^{\prime} \cup \mathcal{C}^{\prime\prime}$ covers $a_j^{\prime}$. Thus,
\begin{equation*}
\mathcal{C}^{\prime} \cup \mathcal{C}^{\prime\prime} \textit{ is a covering of } U_j.
\end{equation*}
So, since $\vert \mathcal{C}_i(U_h)\vert = \vert\mathcal{C}_i(L, U_h)\vert + \vert \mathcal{C}_i(R, U_h)\vert + 1$, we have
\begin{multline*}
\mincov{U_j} \leqslant \vert \mathcal{C}^{\prime} \cup \mathcal{C}^{\prime\prime} \vert \leqslant \vert \mathcal{C}^{\prime}\vert + \vert \mathcal{C}^{\prime\prime}\vert \leqslant \\
\leqslant \vert \mathcal{C}_i(L, U_h)\vert + \vert \mathcal{C}_i(R, U_h)\vert < \vert \mathcal{C}_i(U_h)\vert = \mincov{U_h}.
\end{multline*}

Now consider the case $a_j = 1$. Assume $U_j = L\cdot R = \left(L^{\prime}C\right)\cdot \left(C^{-1}R^{\prime}\right) = L^{\prime}R^{\prime}$, where $L^{\prime}R^{\prime}$ has no further cancellations. Assume $\mathcal{C}_{i^{\prime}}(U_j)$ is some covering of $U_j$. Using the same argument as above, we obtain that an occurrence of a letter in $L^{\prime}$ or in $R^{\prime}$ in $U_h$ is covered by $\mathcal{C}_i(U_h)$ if and only if the corresponding occurrence of this letter in $U_j$ is covered by $\mathcal{C}_{i^{\prime}}(U_j)$.

Consider images of elements of $\mathcal{C}_i(L^{\prime}, U_h)$ in $U_j$. Let us take one arbitrary image of every element of $\mathcal{C}_i(L^{\prime}, U_h)$ and denote the obtained set by $\mathcal{C}^{\prime}$. Similarly for $\mathcal{C}_i(R^{\prime}, U_h)$, let us take one arbitrary image of every element of $\mathcal{C}_i(L^{\prime}, U_h)$ and denote the obtained set by $\mathcal{C}^{\prime\prime}$. As above, we see that every letter of $L^{\prime}$ that is covered by $\mathcal{C}_i(L^{\prime}, U_h)$ in $U_h$ is covered by $\mathcal{C}^{\prime}$ in $U_j$. Every letter of $R^{\prime}$ that is covered by $\mathcal{C}_i(R^{\prime}, U_h)$ in $U_h$ is covered by $\mathcal{C}^{\prime\prime}$ in $U_j$. So, since $U_j = L^{\prime}R^{\prime}$,
\begin{equation*}
\mathcal{C}^{\prime} \cup \mathcal{C}^{\prime\prime} \textit{ is a covering of } U_j.
\end{equation*}
Since we take one image of every element of $\mathcal{C}_i(L^{\prime}, U_h)$ and $\mathcal{C}_i(R^{\prime}, U_h)$, we have
\begin{equation*}
\vert \mathcal{C}^{\prime}\vert \leqslant \vert \mathcal{C}_{i}(L^{\prime}, U_h)\vert, \vert \mathcal{C}^{\prime\prime}\vert \leqslant \vert \mathcal{C}_{i}(R^{\prime}, U_h)\vert.
\end{equation*}
Since $L^{\prime}$ is contained in $L$ in $U_h$ and $R^{\prime}$ is contained in $R$ in $U_h$, we have
\begin{equation*}
\mathcal{C}_i(L^{\prime}, U_h) \subseteq \mathcal{C}_i(L, U_h),\ \mathcal{C}_i(R^{\prime}, U_h) \subseteq \mathcal{C}_i(R, U_h).
\end{equation*}
So,
\begin{equation*}
\vert \mathcal{C}_i(L^{\prime}, U_h)\vert \leqslant \vert\mathcal{C}_i(L, U_h)\vert,\ \vert\mathcal{C}_i(R^{\prime}, U_h)\vert \leqslant \mathcal{C}_i(R, U_h)\vert.
\end{equation*}
Recall that $\vert \mathcal{C}_i(U_h)\vert = \vert\mathcal{C}_i(L, U_h)\vert + \vert \mathcal{C}_i(R, U_h)\vert + 1$. Therefore, finally we obtain
\begin{multline*}
\mincov{U_j} \leqslant \vert \mathcal{C}^{\prime} \cup \mathcal{C}^{\prime\prime} \vert \leqslant \vert \mathcal{C}^{\prime}\vert + \vert \mathcal{C}^{\prime\prime}\vert \leqslant \vert \mathcal{C}_i(L^{\prime}, U_h)\vert + \vert \mathcal{C}_i(R^{\prime}, U_h)\vert \leqslant\\
\leqslant \vert \mathcal{C}_i(L, U_h)\vert + \vert \mathcal{C}_i(R, U_h)\vert < \vert \mathcal{C}_i(U_h)\vert = \mincov{U_h}.
\end{multline*}
\end{proof}

\subsection{Admissible replacements of incident monomials}
\label{admissible_replacements_section}
Let $U$ be a monomial. We separate the set $\mo{U}$ in two distinct parts according to $\SPM$-measure of maximal occurrences in $U$. Namely, we consider all the maximal occurrences in $U$ of $\SPM$-measure $\geqslant 3$ and denote this set by $\longmo{U}$\label{longmo_def}, and we consider all the maximal occurrences in $U$ of $\SPM$-measure $\leqslant 2$ and denote this set by $\shortmo{U}$\label{shortmo_def}. Since values of $\SPM$-measure are natural numbers, we see that $\mo{U} = \longmo{U} \sqcup \shortmo{U}$, where $\sqcup$ is a disjoint union. It follows directly from Lemma~\ref{fc_measue} that $\fc{U} \subseteq \shortmo{U}$ and $\longmo{U} \subseteq \nfc{U}$.

\begin{definition}
\label{admissible_replacement}
Let $U_h$ be a monomial, $a_h \in \mo{U_h}$. Let $a_h$ and $a_j$ be incident monomials. Consider a replacement $a_h \mapsto a_j$ in $U_h$. Let $U_j$ be the resulting monomial. We say that the replacement $a_h \mapsto a_j$ in $U_h$ is \emph{admissible} if $\SPM(a_h) \geqslant \tau - 2$ and $a_j \neq 1$ and $a_j$ is not fully covered by images of elements of $\longmo{U_h} \setminus \lbrace a_h\rbrace$ in $U_j$.
\end{definition}

\begin{remark}
Assume $\SPM(a_h) \geqslant \tau - 2$. Then, by Lemma~\ref{minimal_coverings_property}, we have $\mincov{U_j} \leqslant \mincov{U_h}$. Notice that if $a_h \mapsto a_j$ is an admissible replacement in $U_h$, we can obtain both $\mincov{U_j} = \mincov{U_h}$ or $\mincov{U_j} < \mincov{U_h}$. Let $\SPM(a_h) \geqslant \tau - 2$ and $a_h \mapsto a_j$ be not an admissible replacement in $U_h$. Then, by definition, $a_j$ is covered by images of element of $\longmo{U_h} \setminus \lbrace a_h\rbrace$ in $U_j$. This yields that $a_j$ is covered by images of elements $\nfc{U_h} \setminus \lbrace a_h\rbrace$ in $U_j$, because $\longmo{U_h} \subseteq \nfc{U_h}$. Hence, by Lemma~\ref{minimal_coverings_property}, we obtain $\mincov{U_j} < \mincov{U_h}$ in this case.
\end{remark}

\begin{definition}
\label{neighbours_def}
Assume $U$ is a monomial, $a \in \mo{U}$. Assume $b\in \mo{U}$, $b$ starts from the left of the beginning of $a$, and $b$ is not separated from $a$. Then $b$ is called \emph{a left neighbour of $a$ in $U$}. Assume $c \in \mo{U}$, $c$ starts from the right of the beginning of $a$, and $c$ is not separated from $a$. Then $c$ is called \emph{a right neighbour of $a$ in $U$}.

Notice that the maximal occurrence $a$ may have several left neighbours and several right neighbours.
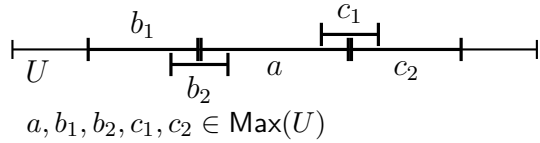
\begin{center}
\begin{tikzpicture}
\node at (-0.5, 0) [below, xshift=10] {$U$};
\draw[|-|, black, thick] (-0.5,0)--(6.5,0);
\draw[|-|, black, very thick] (2, 0)--(4, 0) node[midway, below] {$a$};
\draw[|-|, black, very thick] (0.5, 0)--(2, 0) node[midway, above] {$b_1$};
\draw[|-|, black, very thick] (1.6, -0.2)--(2.4, -0.2) node[midway, below] {$b_2$};
\draw[|-|, black, very thick] (4, 0)--(5.5, 0) node[midway, below] {$c_2$};
\draw[|-|, black, very thick] (3.6, 0.2)--(4.4, 0.2) node[midway, above] {$c_1$};
\node[text width=4cm, align=left] at (1.7, -1) {\small\baselineskip=10pt \textit{$a, b_1, b_2, c_1, c_2 \in \mo{U}$}};
\end{tikzpicture}
\end{center}
\end{definition}

\begin{remark}
\label{admissible_replacement_propoerties}
Let $U_h$ be a monomial, $a_h \in \mo{U_h}$. Let $a_h \mapsto a_j$ be an admissible replacement in $U_h$, $U_j$ be a resulting monomial. Let $U_h = La_hR$. We need a notion of an admissible replacement to prevent the following two effects.
\begin{enumerate}
\item
\label{a_j_not_inside}
Assume $b \in \mo{U_h}$ is a left neighbour of $a_h$. Let $\widehat{b}$ be the intersection of $b$ and $L$. Assume $\SPM(\widehat{b}) \geqslant 3$. Then $b \in \longmo{U_h}$. Let $b^{\prime}$ be the image of $b$ in $U_j$. Since $b \in \longmo{U_h}$, it follows directly from the definition of admissible replacements that $a_j$ is not contained in $b^{\prime}$. Namely, the following configurations are not possible:
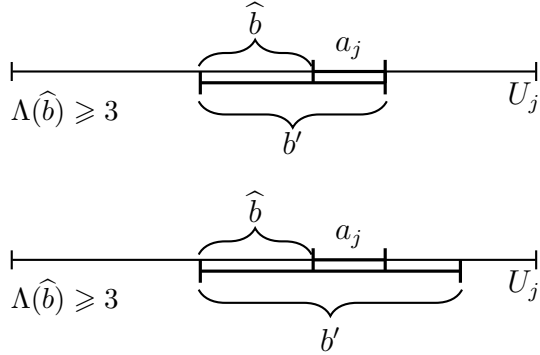
\begin{center}
\begin{tikzpicture}
\draw[|-|, black, thick] (0,0)--(7,0);
\node[below, xshift=-5] at (7, 0) {$U_j$};
\draw[|-|, black, very thick] (4,0)--(5,0) node[midway, above] {$a_j$};
\draw[|-|, black, very thick] (2.5,-0.15)--(5,-0.15);
\draw [thick, decorate, decoration={brace, amplitude=10pt, raise=7pt, mirror}] (2.5,-0.15) to node[midway, below, yshift=-14pt] {$b^{\prime}$} (5, -0.15);
\draw [thick, decorate, decoration={brace, amplitude=10pt, raise=2pt}] (2.5,0) to node[midway, above, yshift=10pt] {$\widehat{b}$} (4, 0);
\node[text width=4cm, align=left] at (2, -0.5) {\small\baselineskip=10pt \textit{$\SPM(\widehat{b}) \geqslant 3$}};
\end{tikzpicture}

\begin{tikzpicture}
\draw[|-|, black, thick] (0,0)--(7,0);
\node[below, xshift=-5] at (7, 0) {$U_j$};
\draw[|-|, black, very thick] (4,0)--(5,0) node[midway, above] {$a_j$};
\draw[|-|, black, very thick] (2.5,-0.15)--(6,-0.15);
\draw [thick, decorate, decoration={brace, amplitude=10pt, raise=7pt, mirror}] (2.5,-0.15) to node[midway, below, yshift=-17pt] {$b^{\prime}$} (6, -0.15);
\draw [thick, decorate, decoration={brace, amplitude=10pt, raise=2pt}] (2.5,0) to node[midway, above, yshift=10pt] {$\widehat{b}$} (4, 0);
\node[text width=4cm, align=left] at (2, -0.5) {\small\baselineskip=10pt \textit{$\SPM(\widehat{b}) \geqslant 3$}};
\end{tikzpicture}
\end{center}

Similarly, assume $c \in \mo{U_h}$ is a right neighbour of $a_h$. Let $\widehat{c}$ be the intersection of $c$ and $R$. Assume $\SPM(\widehat{c}) \geqslant 3$. Then $c \in \longmo{U_h}$. Let $c^{\prime}$ be the image of $b$ in $U_j$. Since $c \in \longmo{U_h}$, it follows directly from the definition that $a_j$ is not contained in $c^{\prime}$. Namely, the following configurations are not possible:
\begin{center}
\begin{tikzpicture}
\draw[|-|, black, thick] (0,0)--(7,0);
\node[below, xshift=10] at (0, 0) {$U_j$};
\draw[|-|, black, very thick] (2,0)--(3,0) node[midway, above] {$a_j$};
\draw[|-|, black, very thick] (2,-0.15)--(4.5,-0.15);
\draw [thick, decorate, decoration={brace, amplitude=10pt, raise=7pt, mirror}] (2,-0.15) to node[midway, below, yshift=-14pt] {$c^{\prime}$} (4.5, -0.15);
\draw [thick, decorate, decoration={brace, amplitude=10pt, raise=2pt}] (3,0) to node[midway, above, yshift=10pt] {$\widehat{c}$} (4.5, 0);
\node[text width=4cm, align=left] at (7, -0.5) {\small\baselineskip=10pt \textit{$\SPM(\widehat{c}) \geqslant 3$}};
\end{tikzpicture}

\begin{tikzpicture}
\draw[|-|, black, thick] (0,0)--(7,0);
\node[below, xshift=10] at (0, 0) {$U_j$};
\draw[|-|, black, very thick] (2,0)--(3,0) node[midway, above] {$a_j$};
\draw[|-|, black, very thick] (1,-0.15)--(4.5,-0.15);
\draw [thick, decorate, decoration={brace, amplitude=10pt, raise=7pt, mirror}] (1,-0.15) to node[midway, below, yshift=-17pt] {$c^{\prime}$} (4.5, -0.15);
\draw [thick, decorate, decoration={brace, amplitude=10pt, raise=2pt}] (3,0) to node[midway, above, yshift=10pt] {$\widehat{c}$} (4.5, 0);
\node[text width=4cm, align=left] at (7, -0.5) {\small\baselineskip=10pt \textit{$\SPM(\widehat{c}) \geqslant 3$}};
\end{tikzpicture}
\end{center}

\item
\label{a_j_not_covered_by_two}
Assume $b \in \mo{U_h}$ is a left neighbour if $a_h$, $c \in \mo{U_h}$ is a right neighbour of $a_h$. Let $\widehat{b}$ be the intersection of $b$ and $L$, $\widehat{c}$ be the intersection of $c$ and $R$. Assume $\SPM(\widehat{b}) \geqslant 3$ and $\SPM(\widehat{c}) \geqslant 3$. Then $b \in \longmo{U_h}$ and $c \in \longmo{U_h}$. Let $b^{\prime}$ be the image of $b$ in $U_j$, and $c^{\prime}$ be the image of $c$ in $U_j$. Then it follows directly from the definition of admissible replacements that $b^{\prime}$ and $c^{\prime}$ are separated. Namely, the following configurations are not possible:
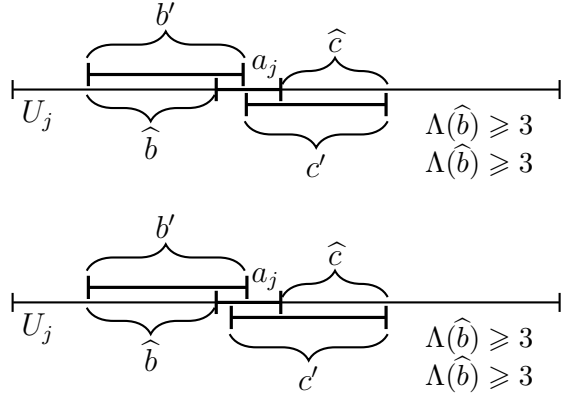
\begin{center}
\begin{tikzpicture}
\node at (0, 0) [below, xshift=10] {$U_j$};
\draw[|-|, black, thick] (0,0)--(7.3,0);
\draw[|-|, black, very thick] (1, 0.2)--(3.1, 0.2);
\draw [thick, decorate, decoration={brace, amplitude=10pt, raise=12pt}] (1, 0) to node[midway, above, yshift=20pt] {$b^{\prime}$} (3.1, 0);
\draw[|-|, black, very thick] (2.7, 0)--(3.6, 0) node[near end, above] {$a_j$};
\draw[|-|, black, very thick] (3.1, -0.2)--(5, -0.2);
\draw [thick, decorate, decoration={brace, amplitude=10pt, raise=12pt, mirror}] (3.1, 0) to node[midway, below, yshift=-20pt] {$c^{\prime}$} (5, 0);
\draw [thick, decorate, decoration={brace, amplitude=10pt, raise=2pt, mirror}] (1,0) to node[midway, below, yshift=-10pt] {$\widehat{b}$} (2.7, 0);
\draw [thick, decorate, decoration={brace, amplitude=10pt, raise=2pt}] (3.6,0) to node[midway, above, yshift=10pt] {$\widehat{c}$} (5, 0);
\node[text width=4cm, align=left] at (7.5, -0.7) {\small\baselineskip=10pt \textit{$\SPM(\widehat{b}) \geqslant 3$\\$\SPM(\widehat{b}) \geqslant 3$}};
\end{tikzpicture}

\begin{tikzpicture}
\node at (0, 0) [below, xshift=10] {$U_j$};
\draw[|-|, black, thick] (0,0)--(7.3,0);
\draw[|-|, black, very thick] (1, 0.2)--(3.15, 0.2);
\draw [thick, decorate, decoration={brace, amplitude=10pt, raise=12pt}] (1, 0) to node[midway, above, yshift=20pt] {$b^{\prime}$} (3.1, 0);
\draw[|-|, black, very thick] (2.7, 0)--(3.6, 0) node[near end, above] {$a_j$};
\draw[|-|, black, very thick] (2.9, -0.2)--(5, -0.2);
\draw [thick, decorate, decoration={brace, amplitude=10pt, raise=12pt, mirror}] (2.9, 0) to node[midway, below, yshift=-20pt] {$c^{\prime}$} (5, 0);
\draw [thick, decorate, decoration={brace, amplitude=10pt, raise=2pt, mirror}] (1,0) to node[midway, below, yshift=-10pt] {$\widehat{b}$} (2.7, 0);
\draw [thick, decorate, decoration={brace, amplitude=10pt, raise=2pt}] (3.6,0) to node[midway, above, yshift=10pt] {$\widehat{c}$} (5, 0);
\node[text width=4cm, align=left] at (7.5, -0.7) {\small\baselineskip=10pt \textit{$\SPM(\widehat{b}) \geqslant 3$\\$\SPM(\widehat{b}) \geqslant 3$}};
\end{tikzpicture}
\end{center}
\end{enumerate}
\end{remark}

\begin{remark}
\label{admissible_short_merging}
Let us consider Property~\ref{a_j_not_inside} from Remark~\ref{admissible_replacement_propoerties} from a different point of view and state it in other words.

Let $U_h$ be a monomial, $a_h \in \mo{U_h}$, $U_h = La_hR$, and $\SPM(a_h) \geqslant \tau - 2$. Assume $a_h$ and $a_j$ are incident monomials. Consider the replacement $a_h \mapsto a_j$ in $U_h$, let $U_j$ be the resulting monomial. Assume $b \in \mo{U_h}$, $b$ is a left neighbour of $a_h$. Let $\widehat{b}$ be the intersection of $b$ and $L$. Assume $c \in \mo{U_h}$, $c$ is a left neighbour of $a_h$. Let $\widehat{c}$ be the intersection of $c$ and $R$. Assume $\widehat{b}$ and $\widehat{c}$ merge to one maximal occurrence $\widehat{b}a_j\widehat{c}$ in $U_j$. Then it follows from the results of Section~\ref{mt_configurations} that $a_j$ is a small piece (see~\ref{a_j_keep_structure}).
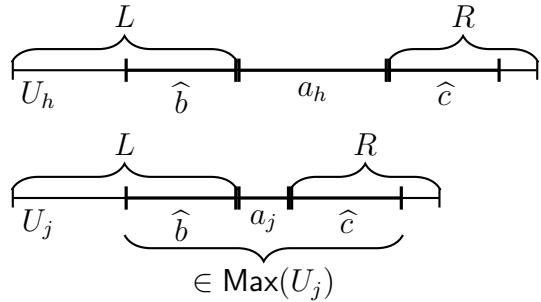
\begin{center}
\begin{tikzpicture}
\draw[|-|, black, thick] (0,0)--(7,0);
\node[below, xshift=10] at (0, 0) {$U_h$};
\draw[|-|, black, very thick] (3,0)--(5,0) node[midway, below] {$a_h$};
\draw[|-|, black, very thick] (5,0)--(6.5,0) node[midway, below] {$\widehat{c}$};
\draw[|-|, black, very thick] (1.5,0)--(3,0) node[midway, below] {$\widehat{b}$};
\draw [thick, decorate, decoration={brace, amplitude=10pt, raise=3pt}] (0, 0) to node[midway, above, yshift=12pt] {$L$} (3, 0);
\draw [thick, decorate, decoration={brace, amplitude=10pt, raise=3pt}] (5, 0) to node[midway, above, yshift=12pt] {$R$} (7, 0);

\draw[|-|, black, thick] (0,-1.7)--(5.7,-1.7);
\node[below, xshift=10] at (0, -1.7) {$U_j$};
\draw[|-|, black, very thick] (3,-1.7)--(3.7,-1.7) node[midway, below] {$a_j$};
\draw[|-|, black, very thick] (3.7,-1.7)--(5.2,-1.7) node[midway, below] {$\widehat{c}$};
\draw[|-|, black, very thick] (1.5,-1.7)--(3,-1.7) node[midway, below] {$\widehat{b}$};
\draw [thick, decorate, decoration={brace, amplitude=10pt, raise=3pt}] (0, -1.7) to node[midway, above, yshift=12pt] {$L$} (3, -1.7);
\draw [thick, decorate, decoration={brace, amplitude=10pt, raise=3pt}] (3.7, -1.7) to node[midway, above, yshift=12pt] {$R$} (5.7, -1.7);
\draw [thick, decorate, decoration={brace, amplitude=10pt, raise=14pt, mirror}] (1.5, -1.7) to node[midway, below, yshift=-22pt] {$\in \mo{U_j}$} (5.2, -1.7);
\end{tikzpicture}
\end{center}

Assume $\SPM(\widehat{b}a_j\widehat{c}) \geqslant 6$. Notice that $\SPM(a_j) = 1$ because $a_j$ is a small piece. Therefore,
\begin{equation*}
\SPM(\widehat{b}a_j\widehat{c}) \leqslant \SPM(\widehat{b}) + \SPM(\widehat{c}) + \SPM(a_j) = \SPM(\widehat{b}) + \SPM(\widehat{c}) + 1.
\end{equation*}
Since $\SPM(\widehat{b}a_j\widehat{c}) \geqslant 6$, we see that $\SPM(\widehat{b}) \geqslant 3$ or $\SPM(\widehat{c}) \geqslant 3$. To be definite, assume $\SPM(\widehat{b}) \geqslant 3$. Then $\SPM(b) \geqslant \SPM(\widehat{b}) \geqslant 3$. Hence, $b \in \longmo{U_h}$. By definition, $\widehat{b}a_j\widehat{c}$ is an image of $b$ in $U_j$. Therefore, $a_j$ is covered by an image of an element of $\longmo{U_h} \setminus \lbrace a_h\rbrace$. Clearly, we obtain the same result if $\SPM(\widehat{c}) \geqslant 3$. So, $a_h \mapsto a_j$ is not an admissible replacement in $U_h$ if $\widehat{b}a_j\widehat{c} \in \mo{U_j}$ and $\SPM(\widehat{b}a_j\widehat{c}) \geqslant 6$.

Thus, if $a_h \mapsto a_j$ is an admissible replacement in $U_h$, then $\SPM(\widehat{b}a_j\widehat{c}) \leqslant 5$. This means that we can not obtain ``very long'' maximal occurrences in $U_j$ as a result of merging if $a_h \mapsto a_j$ is an admissible replacement in $U_h$.
\end{remark}

\begin{lemma}
\label{single_nfc_neighbour}
Let $U$ be a monomial, $a \in \mo{U}$. Then either $a$ does not have left neighbours that belong to $\nfc{U}$, or it has just one left neighbour that belongs to $\nfc{U}$ (but it may have other left neighbours provided they do not belong to $\nfc{U}$). Moreover, if $a$ has a left neighbour that belongs to $\nfc{U}$, then this is a left neighbour of $a$ with the leftmost beginning point. The same holds for right neighbours, that is, either $a$ does not have right neighbours that belong to $\nfc{U}$, or it has just one right neighbour that belongs to $\nfc{U}$ (but it may have other right neighbours provided they do not belong to $\nfc{U}$). Moreover, if $a$ has a right neighbour that belongs to $\nfc{U}$, then this is a right neighbour of $a$ with the rightmost end point.
\end{lemma}
\begin{proof}
Assume $m_1, m_2 \in \mo{U}$. Let $m_1$ and $m_2$ be left neighbours of $a$. By definition, both $m_1$ and $m_2$ start from the left of the beginning of $a$. Without loss of generality we may assume that $m_1$ starts from the left of the beginning of $m_2$. Then, since both $m_1$ and $m_2$ are not separated from $a$, $m_2$ is covered by $m_1$ and $a$.
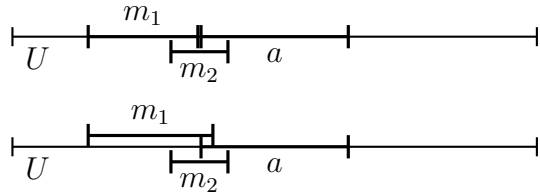
\begin{center}
\begin{tikzpicture}
\node at (-0.5, 0) [below, xshift=10] {$U$};
\draw[|-|, black, thick] (-0.5,0)--(6.5,0);
\draw[|-|, black, very thick] (2, 0)--(4, 0) node[midway, below] {$a$};
\draw[|-|, black, very thick] (0.5, 0)--(2, 0) node[midway, above] {$m_1$};
\draw[|-|, black, very thick] (1.6, -0.2)--(2.4, -0.2) node[midway, below] {$m_2$};
\end{tikzpicture}

\begin{tikzpicture}
\node at (-0.5, 0) [below, xshift=10] {$U$};
\draw[|-|, black, thick] (-0.5,0)--(6.5,0);
\draw[|-|, black, very thick] (2, 0)--(4, 0) node[midway, below] {$a$};
\draw[|-|, black, very thick] (0.5, 0.15)--(2.2, 0.15) node[midway, above] {$m_1$};
\draw[|-|, black, very thick] (1.6, -0.2)--(2.4, -0.2) node[midway, below] {$m_2$};
\end{tikzpicture}
\end{center}
So, $m_2 \in \fc{U}$. This means that $a$ either has a single left neighbour that belongs to $\nfc{U}$, or does not have left neighbours that belong to $\nfc{U}$ at all. Furthermore, only a left neighbour of $a$ with the leftmost beginning may belong to $\nfc{U}$. The statement about right neighbours of $a$ is proved in the same way.
\end{proof}

\begin{corollary}
\label{single_long_neighbour}
Let $U$ be a monomial, $a \in \mo{U}$. Then either $a$ does not have left neighbours that belong to $\longmo{U}$, or it has just one left neighbour that belongs to $\longmo{U}$ (but it may have other left neighbours provided they do not belong to $\longmo{U}$). Moreover, if $a$ has a left neighbour that belongs to $\longmo{U}$, then this is a left neighbour of $a$ with the leftmost beginning point. The same holds for right neighbours, that is, either $a$ does not have right neighbours that belong to $\longmo{U}$, or it has just one right neighbour that belongs to $\longmo{U}$ (but it may have other right neighbours provided they do not belong to $\longmo{U}$). Moreover, if $a$ has a right neighbour that belongs to $\longmo{U}$, then this is a right neighbour of $a$ with the rightmost end point.
\end{corollary}
\begin{proof}
Since $\longmo{U} \subseteq \nfc{U}$, the statement of Corollary~\ref{single_long_neighbour} follows directly from Lemma~\ref{single_nfc_neighbour}.
\end{proof}

\begin{lemma}
\label{non_covered_inverse_image}
Let $U_h$ be a monomial, $a_h \in \mo{U_h}$, $U_h = La_hR$. Let $a_h$ and $a_j$ be incident monomials, $a_j \neq 1$. Consider the replacement $a_h \mapsto a_j$ in $U_h$. Let $U_j$ be the resulting monomial, $U_j = La_jR$. Assume $b^{\prime} \in \mo{U_j}$. Let $\widehat{b}_1$ be the intersection of $b^{\prime}$ and $L$, $\widehat{b}_2$ be the intersection of $b^{\prime}$ and $R$.  Recall that we can consider $\widehat{b}_1$ and $\widehat{b}_2$ as occurrences in $U_h$.
\begin{enumerate}
\item
\label{non_covered_inverse_image1}
Assume $\SPM(\widehat{b}_1) \geqslant 3$. Then there exists $b \in \mo{U_h} \setminus \lbrace a_h \rbrace$ such that $b^{\prime}$ is a single image of $b$ in $U_j$. Furthermore, $b = \widehat{b}_1Y$, where $Y$ is a suffix of $b$, $Y$ is an overlap of $b$ and $a_h$ in $U_h$ ($Y$ is empty if $b$ and $a_h$ are separated or touch at a point), and $b \in \longmo{U_h}$.
\item
\label{non_covered_inverse_image2}
Assume $\SPM(\widehat{b}_2) \geqslant 3$. Then there exists $b \in \mo{U_h} \setminus \lbrace a_h \rbrace$ such that $b^{\prime}$ is a single image of $b$ in $U_j$. Furthermore, $b = X\widehat{b}_2$, where $X$ is a prefix of $b$, $X$ is an overlap of $b$ and $a_h$ in $U_h$ ($X$ is empty if $b$ and $a_h$ are separated or touch at a point), and $b \in \longmo{U_h}$.
\end{enumerate}

Under the same conditions suppose that $a_h \mapsto a_j$ is an admissible replacement in $U_h$. Then the element $b$ obtained above is a unique element of $\nfc{U_h}$ such that $b^{\prime}$ is its image in $U_j$.
\end{lemma}
\begin{proof}
Let us prove statement~\ref{non_covered_inverse_image1}. That is, $\SPM(\widehat{b}_1) \geqslant 3$. Statement~\ref{non_covered_inverse_image2} is proved similarly. Notice, that we can have $\SPM(\widehat{b}) \geqslant 3$ and $\SPM(\widehat{c}) \geqslant 3$ simultaneously.

By the condition of Lemma~\ref{non_covered_inverse_image}, the occurrence $\widehat{b}_1$ is the intersection of $b^{\prime}$ and $L$. If $b^{\prime}$ does not start from the left of the beginning of $a_j$, then $\widehat{b}_1$ is empty. So, the assumption $\SPM(\widehat{b}_1) \geqslant 3$ implies that $b^{\prime}$ starts from the left of the beginning of $a_j$.

Assume $b^{\prime}$ and $a_j$ are separated. Then $\widehat{b}_1 = b^{\prime}$. Denote by $b$ the occurrence $\widehat{b}_1$ in $U_h$. Then we can argue as earlier and obtain that $b$ is a maximal occurrence in $U_h$. Since $\SPM(\widehat{b}_1) \geqslant 3$, we obtain that $\SPM(b) \geqslant 3$. That is, $b \in \longmo{U_h}$. Clearly, $b$ and $a_h$ are separated. Hence, by the definition of images, $b^{\prime}$ is a single image of $b$ in $U_j$.
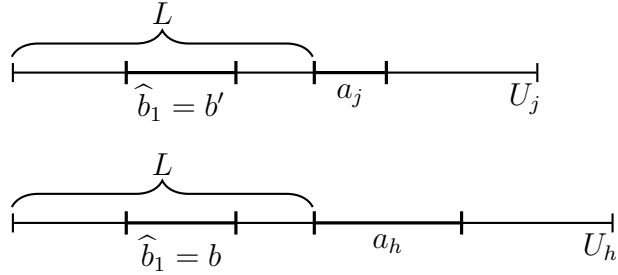
\begin{center}
\begin{tikzpicture}
\draw[|-|, black, thick] (0,0)--(7,0);
\node[below, xshift=-5] at (7, 0) {$U_j$};
\draw [thick, decorate, decoration={brace, amplitude=10pt, raise=6pt}] (0,0) to node[midway, above, yshift=14pt] {$L$} (4, 0);
\draw[|-|, black, very thick] (4,0)--(5,0) node[midway, below] {$a_j$};
\draw[|-|, black, very thick] (1.5,0)--(3,0) node[midway, below] {$\widehat{b}_1 = b^{\prime}$};

\draw[|-|, black, thick] (0,-2)--(8,-2);
\node[below, xshift=-5] at (8, -2) {$U_h$};
\draw [thick, decorate, decoration={brace, amplitude=10pt, raise=6pt}] (0,-2) to node[midway, above, yshift=14pt] {$L$} (4, -2);
\draw[|-|, black, very thick] (4,-2)--(6,-2) node[midway, below] {$a_h$};
\draw[|-|, black, very thick] (1.5,-2)--(3,-2) node[midway, below] {$\widehat{b}_1 = b$};
\end{tikzpicture}
\end{center}

Assume $b^{\prime}$ and $a_j$ are not separated. Then $b^{\prime}$ and $a_j$ either touch at a point, or have an overlap, or $a_j$ is contained inside $b^{\prime}$. Notice that if $a_j$ is contained in $b^{\prime}$, then $a_j$ is a small piece (see Section~\ref{mt_configurations},~\ref{a_j_keep_structure}).
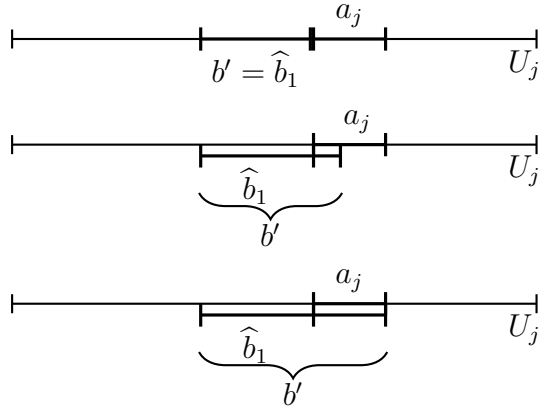
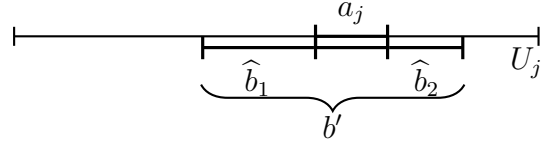
\begin{center}
\begin{tikzpicture}
\draw[|-|, black, thick] (0,0)--(7,0);
\node[below, xshift=-5] at (7, 0) {$U_j$};
\draw[|-|, black, very thick] (4,0)--(5,0) node[midway, above] {$a_j$};
\draw[|-|, black, very thick] (2.5,0)--(4,0) node[midway, below] {$b^{\prime} = \widehat{b}_1$};
\end{tikzpicture}

\begin{tikzpicture}
\draw[|-|, black, thick] (0,0)--(7,0);
\node[below, xshift=-5] at (7, 0) {$U_j$};
\draw[|-|, black, very thick] (4,0)--(5,0) node[midway, above, xshift=3] {$a_j$};
\draw[|-, black, very thick] (2.5,-0.15)--(4,-0.15) node[midway, below] {$\widehat{b}_1$};
\draw[|-|, black, very thick] (4,-0.15)--(4.4,-0.15);
\draw [thick, decorate, decoration={brace, amplitude=10pt, raise=14pt, mirror}] (2.5,-0.15) to node[midway, below, yshift=-21pt] {$b^{\prime}$} (4.4, -0.15);
\end{tikzpicture}

\begin{tikzpicture}
\draw[|-|, black, thick] (0,0)--(7,0);
\node[below, xshift=-5] at (7, 0) {$U_j$};
\draw[|-|, black, very thick] (4,0)--(5,0) node[midway, above] {$a_j$};
\draw[|-, black, very thick] (2.5,-0.15)--(4,-0.15) node[midway, below] {$\widehat{b}_1$};
\draw[|-|, black, very thick] (4,-0.15)--(5,-0.15);
\draw [thick, decorate, decoration={brace, amplitude=10pt, raise=14pt, mirror}] (2.5,-0.15) to node[midway, below, yshift=-21pt] {$b^{\prime}$} (5, -0.15);
\end{tikzpicture}

\begin{tikzpicture}
\draw[|-|, black, thick] (0,0)--(7,0);
\node[below, xshift=-5] at (7, 0) {$U_j$};
\draw[|-|, black, very thick] (4,0)--(5,0) node[midway, above] {$a_j$};
\draw[|-, black, very thick] (2.5,-0.15)--(4,-0.15) node[midway, below] {$\widehat{b}_1$};
\draw[|-|, black, very thick] (4,-0.15)--(5,-0.15);
\draw[-|, black, very thick] (5,-0.15)--(6,-0.15) node[midway, below] {$\widehat{b}_2$};
\draw [thick, decorate, decoration={brace, amplitude=10pt, raise=14pt, mirror}] (2.5,-0.15) to node[midway, below, yshift=-21pt] {$b^{\prime}$} (6, -0.15);
\end{tikzpicture}
\end{center}
We again consider $\widehat{b}_1$ as an occurrence in $U_h$.
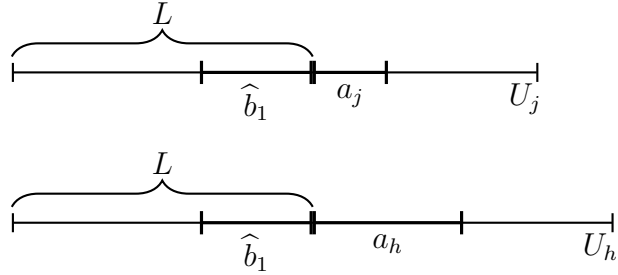
\begin{center}
\begin{tikzpicture}
\draw[|-|, black, thick] (0,0)--(7,0);
\node[below, xshift=-5] at (7, 0) {$U_j$};
\draw[|-|, black, very thick] (4,0)--(5,0) node[midway, below] {$a_j$};
\draw[|-|, black, very thick] (2.5,0)--(4,0) node[midway, below] {$\widehat{b}_1$};
\draw [thick, decorate, decoration={brace, amplitude=10pt, raise=6pt}] (0,0) to node[midway, above, yshift=14pt] {$L$} (4, 0);

\draw[|-|, black, thick] (0,-2)--(8,-2);
\node[below, xshift=-5] at (8, -2) {$U_h$};
\draw[|-|, black, very thick] (4,-2)--(6,-2) node[midway, below] {$a_h$};
\draw[|-|, black, very thick] (2.5,-2)--(4,-2) node[midway, below] {$\widehat{b}_1$};
\draw [thick, decorate, decoration={brace, amplitude=10pt, raise=6pt}] (0,-2) to node[midway, above, yshift=14pt] {$L$} (4, -2);
\end{tikzpicture}
\end{center}
Since we assume that $\SPM(\widehat{b}_1) \geqslant 3$, $\widehat{b}_1$ is not a small piece. Therefore, by Corollary~\ref{max_occurrences_coinside}, there exists a unique maximal occurrence in $U_h$ that contains $\widehat{b}_1$. Denote it by  $b$. Then $\SPM(b) \geqslant \SPM(\widehat{b}) \geqslant 3$, that is, $b \in \longmo{U_h}$.

Let us show that $b = \widehat{b}_1Y$, where $Y$ is a suffix of $b$. Indeed, assume $b = X\widehat{b}_1Y$, where $X$ is a non-empty prefix of $b$. Then $L = L_1X\widehat{b}_1$.
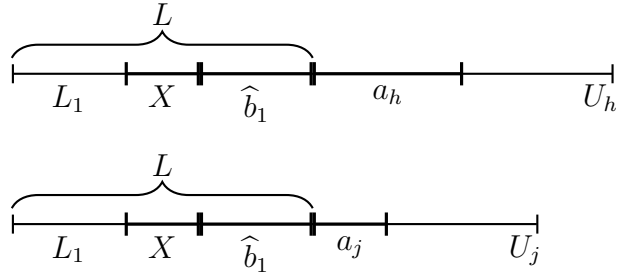
\begin{center}
\begin{tikzpicture}
\draw[|-|, black, thick] (0,0)--(8,0);
\node[below, xshift=-5] at (8, 0) {$U_h$};
\draw[|-|, black, very thick] (4,0)--(6,0) node[midway, below] {$a_h$};
\draw[|-|, black, very thick] (2.5,0)--(4,0) node[midway, below] {$\widehat{b}_1$};
\draw [thick, decorate, decoration={brace, amplitude=10pt, raise=6pt}] (0,0) to node[midway, above, yshift=14pt] {$L$} (4, 0);
\draw[|-|, black, very thick] (1.5,0)--(2.5,0) node[midway, below] {$X$};
\path[|-|, black, very thick] (0,0)--(1.5,0) node[midway, below] {$L_1$};

\draw[|-|, black, thick] (0,-2)--(7,-2);
\node[below, xshift=-5] at (7, -2) {$U_j$};
\draw[|-|, black, very thick] (4,-2)--(5,-2) node[midway, below] {$a_j$};
\draw[|-|, black, very thick] (2.5,-2)--(4,-2) node[midway, below] {$\widehat{b}_1$};
\draw [thick, decorate, decoration={brace, amplitude=10pt, raise=6pt}] (0,-2) to node[midway, above, yshift=14pt] {$L$} (4, -2);
\draw[|-|, black, very thick] (1.5,-2)--(2.5,-2) node[midway, below] {$X$};
\path[|-|, black, very thick] (0,-2)--(1.5,-2) node[midway, below] {$L_1$};
\end{tikzpicture}
\end{center}
Since $b = X\widehat{b}_1Y \in \Mon$, we obtain $X\widehat{b}_1 \in \Mon$. Since $\widehat{b}_1$ in $U_j$ is an intersection of $b^{\prime}$ and $L$, $\widehat{b}_1$ is an initial subword of $b^{\prime}$. Since $\widehat{b}_1$ is not a small piece, we apply Lemma~\ref{not_sp_prolongation} to $X\widehat{b}_1$ and $b^{\prime}$ and obtain that $Xb^{\prime} \in \Mon$.
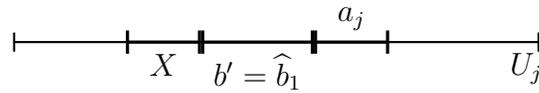
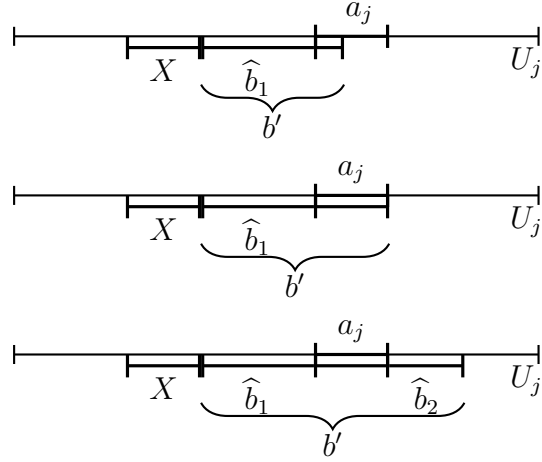
\begin{center}
\begin{tikzpicture}
\draw[|-|, black, thick] (0,0)--(7,0);
\node[below, xshift=-5] at (7, 0) {$U_j$};
\draw[|-|, black, very thick] (4,0)--(5,0) node[midway, above] {$a_j$};
\draw[|-|, black, very thick] (2.5,0)--(4,0) node[midway, below] {$b^{\prime} = \widehat{b}_1$};
\draw[|-|, black, very thick] (1.5,0)--(2.5,0) node[midway, below] {$X$};
\end{tikzpicture}

\begin{tikzpicture}
\draw[|-|, black, thick] (0,0)--(7,0);
\node[below, xshift=-5] at (7, 0) {$U_j$};
\draw[|-|, black, very thick] (4,0)--(5,0) node[midway, above, xshift=3] {$a_j$};
\draw[|-, black, very thick] (2.5,-0.15)--(4,-0.15) node[midway, below] {$\widehat{b}_1$};
\draw[|-|, black, very thick] (4,-0.15)--(4.4,-0.15);
\draw [thick, decorate, decoration={brace, amplitude=10pt, raise=14pt, mirror}] (2.5,-0.15) to node[midway, below, yshift=-21pt] {$b^{\prime}$} (4.4, -0.15);
\draw[|-|, black, very thick] (1.5,-0.15)--(2.5,-0.15) node[midway, below] {$X$};
\end{tikzpicture}

\begin{tikzpicture}
\draw[|-|, black, thick] (0,0)--(7,0);
\node[below, xshift=-5] at (7, 0) {$U_j$};
\draw[|-|, black, very thick] (4,0)--(5,0) node[midway, above] {$a_j$};
\draw[|-, black, very thick] (2.5,-0.15)--(4,-0.15) node[midway, below] {$\widehat{b}_1$};
\draw[|-|, black, very thick] (4,-0.15)--(5,-0.15);
\draw [thick, decorate, decoration={brace, amplitude=10pt, raise=14pt, mirror}] (2.5,-0.15) to node[midway, below, yshift=-21pt] {$b^{\prime}$} (5, -0.15);
\draw[|-|, black, very thick] (1.5,-0.15)--(2.5,-0.15) node[midway, below] {$X$};
\end{tikzpicture}

\begin{tikzpicture}
\draw[|-|, black, thick] (0,0)--(7,0);
\node[below, xshift=-5] at (7, 0) {$U_j$};
\draw[|-|, black, very thick] (4,0)--(5,0) node[midway, above] {$a_j$};
\draw[|-, black, very thick] (2.5,-0.15)--(4,-0.15) node[midway, below] {$\widehat{b}_1$};
\draw[|-|, black, very thick] (4,-0.15)--(5,-0.15);
\draw[-|, black, very thick] (5,-0.15)--(6,-0.15) node[midway, below] {$\widehat{b}_2$};
\draw [thick, decorate, decoration={brace, amplitude=10pt, raise=14pt, mirror}] (2.5,-0.15) to node[midway, below, yshift=-21pt] {$b^{\prime}$} (6, -0.15);
\draw[|-|, black, very thick] (1.5,-0.15)--(2.5,-0.15) node[midway, below] {$X$};
\end{tikzpicture}
\end{center}
Then $b^{\prime}$ is not a maximal occurrence in $U_j$, a contradiction. Therefore, $X$ is empty and $b = \widehat{b}_1Y$.

Since $\widehat{b}_1$ in $U_h$ is non-empty and is contained in the subword $L$ of $U_h$, $\widehat{b}_1$ in $U_h$ starts from the left of the beginning of $a_h$. Since $\widehat{b}_1$ in $U_h$ is contained in $b$, $b$ starts from the left of the beginning of $a_h$ as well. Therefore, $a_h \neq b$. Since $a_j$ and $\widehat{b}_1$ in $U_j$ touch at a point, $a_h$ and $\widehat{b}_1$ in $U_h$ touch at a point as well. So, since $b = \widehat{b}_1Y$, $Y$ is a common subword of $b$ and $a_h$. In other words, $Y$ is an overlap of $a_h$ and $b$ ($Y$ is empty if $a_h$ and $b$ touch at a point).
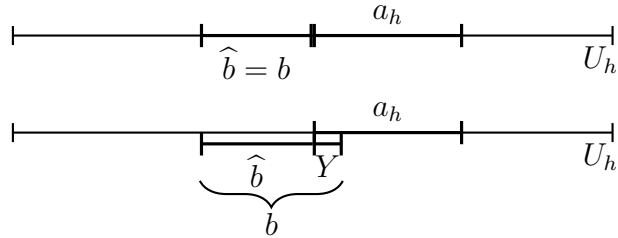
\begin{center}
\begin{tikzpicture}
\draw[|-|, black, thick] (0,0)--(8,0);
\node[below, xshift=-5] at (8, 0) {$U_h$};
\draw[|-|, black, very thick] (4,0)--(6,0) node[midway, above] {$a_h$};
\draw[|-|, black, very thick] (2.5,0)--(4,0) node[midway, below] {$\widehat{b} = b$};
\end{tikzpicture}

\begin{tikzpicture}
\draw[|-|, black, thick] (0,0)--(8,0);
\node[below, xshift=-5] at (8, 0) {$U_h$};
\draw[|-|, black, very thick] (4,0)--(6,0) node[midway, above] {$a_h$};
\draw[|-, black, very thick] (2.5,-0.15)--(4,-0.15) node[midway, below] {$\widehat{b}$};
\draw[|-|, black, very thick] (4,-0.15)--(4.4,-0.15) node[midway, below] {$Y$};
\draw [thick, decorate, decoration={brace, amplitude=10pt, raise=14pt, mirror}] (2.5,-0.15) to node[midway, below, yshift=-22pt] {$b$} (4.4, -0.15);
\end{tikzpicture}
\end{center}

Let us show that $b^{\prime}$ is an image of $b$. Indeed, $\widehat{b}_1$ is a suffix of $L$ and $b = \widehat{b}_1Y$. Hence, the intersection of $b$ and $L$ in $U_h$ is equal to $\widehat{b}_1$. Since $\widehat{b}_1$ in $U_j$ is the intersection of $b^{\prime}$ and $L$, $b^{\prime}$ contains $\widehat{b}$ in $U_j$. Thus, by definition, $b^{\prime}$ is an image of $b$ in $U_j$. Since $b \in \longmo{U_h} \subseteq \nfc{U_h}$, by Lemma~\ref{single_image_lemma}, $b^{\prime}$ is a single image of $b$ in $U_j$.

Now suppose that $a_h \mapsto a_j$ is an admissible replacement in $U_h$. Then we have to show that the element $b$ constructed above is a unique element of $\nfc{U_h}$ such that $b^{\prime}$ is an image of $b$. First of all notice that $a_j$ is not contained in $b^{\prime}$. Indeed, assume the contrary. We proved that the element $b$ constructed above belongs to $\longmo{U_h} \setminus \lbrace a_h\rbrace$. So, since $b^{\prime}$ is an image of $b$, we obtain that $a_j$ is covered by an image of $b \in \longmo{U_h} \setminus \lbrace a_h\rbrace$. This contradicts with the assumption that $a_h \mapsto a_j$ is an admissible replacement in $U_h$. Therefore, $a_j$ is not contained in $b^{\prime}$.

Assume $d \in \nfc{U_h}$, $d \neq b$, and $d^{\prime}$ is an image of $d$ in $U_j$. Let us show that $d^{\prime} \neq b^{\prime}$. There are three possibilities for $d$.
\begin{enumerate}
\item
$d = a_h$;
\item
$d$  starts from the right of the beginning of $a_h$;
\item
$d$ starts from the left of the beginning of $a_h$.
\end{enumerate}

Consider the first possibility, namely, $d = a_h$. Then, by definition, $a_j$ is contained in $d^{\prime}$. We proved above that $a_j$ is not contained in $b^{\prime}$. Hence, $d^{\prime} \neq b^{\prime}$.

Consider the second possibility, namely, that $d$ starts from the right of the beginning of $a_h$. Then the intersection of $d$ and $R$ is non-empty, and we denote it by $\widehat{d}_2$. Recall that $\widehat{d}_2$ can be considered as an occurrence in $U_j$. By the definition of images, $\widehat{d}_2$ in $U_j$ is contained in $d^{\prime}$. Since $\widehat{d}_2$ is non-empty and is contained in the subword $R$ in $U_j$, $\widehat{d}_2$ in $U_j$ ends from the right of the end of $a_j$. So, $d^{\prime}$ ends from the right of the end of $a_j$ as well. We proved above that $a_j$ is not contained in $b^{\prime}$. Hence, $b^{\prime}$ ends from the left of the end of $a_j$. That is, $b^{\prime}$ and $d^{\prime}$ have different end points. Thus, $b^{\prime} \neq d^{\prime}$.

Consider the third possibility, namely, that $d$ starts from the left of the beginning of $a_h$. Then the intersection of $d$ and $L$ is non-empty, and we denote it by $\widehat{d}_1$. We proved that $\widehat{b}_1$ in $U_h$ is an initial subword of $b$. By definition, $\widehat{d}_1$ in $U_h$ is an initial subword of $d$. Since $b$ and $d$ are different maximal occurrences in $U_h$, they have different beginnings. So, $\widehat{b}_1$ and $\widehat{d}_1$ have different beginnings in $U_h$. Recall that $\widehat{d}_1$ can be considered as an occurrence in $U_j$. Then, clearly, $\widehat{b}_1$ and $\widehat{d}_1$ have different beginnings in $U_j$.  Since  $d \in \nfc{U_h}$, it follows from Lemma~\ref{single_image_lemma} that $\widehat{d}_1$ in $U_j$ is an initial subword of $d^{\prime}$. By definition, $\widehat{b}_1$ is an initial subword of $b^{\prime}$. Thus, $b^{\prime}$ and $d^{\prime}$ have different beginnings. So, $b^{\prime} \neq d^{\prime}$. This completes the proof.
\end{proof}

\begin{lemma}
\label{neighbours_a_h_not_separated_right}
Let $U_h$ be a monomial, $b \in \longmo{U_h}$. Suppose $a_h \in \mo{U_h}$ is a right neighbour of $b$. Suppose $a_h \mapsto a_j$ is an admissible replacement in $U_h$. Let $U_j$ be the resulting monomial. Denote by $b^{\prime}$ the image of $b$ in $U_j$. Let $a_j \notin \longmo{U_j}$, that is, either $a_j \in \shortmo{U_j}$, or $a_j \notin \mo{U_j}$. Assume that there exists $c^{\prime}$ a right neighbour of $b^{\prime}$ such that $c^{\prime} \in \longmo{U_j}$. Then $\SPM(c^{\prime}) \leqslant 4$.
\end{lemma}
\begin{proof}
Assume $U_h = La_hR$, then $U_j = La_jR$. Let $\widehat{b}$ be the intersection of $b$ and $L$ in $U_h$. Recall that we also consider $\widehat{b}$ as an occurrence in $U_j$ and, by definition, $b^{\prime}$ contains $\widehat{b}$ in $U_j$. Moreover, since $b\in \longmo{U_h} \subseteq \nfc{U_h}$, it follows from Lemma~\ref{single_image_lemma} that $b^{\prime} = \widehat{b}Y$, where $Y$ is a suffix of $b^{\prime}$.

Since $a_h$ and $b$ are not separated in $U_h$, $a_h$ and $b$ either touch at a point or have an overlap.
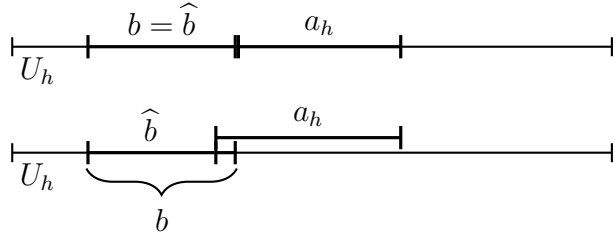
\begin{center}
\begin{tikzpicture}
\node at (0, 0) [below, xshift=10] {$U_h$};
\draw[|-|, black, thick] (0,0)--(8,0);
\draw[|-|, black, very thick] (1, 0)--(3, 0) node[midway, above] {$b = \widehat{b}$};
\draw[|-|, black, very thick] (3, 0)--(5.2, 0) node[midway, above] {$a_h$};
\end{tikzpicture}

\begin{tikzpicture}
\node at (0, 0) [below, xshift=10] {$U_h$};
\draw[|-|, black, thick] (0,0)--(8,0);
\draw[|-, black, very thick] (1, 0)--(2.7, 0) node[midway, above] {$\widehat{b}$};
\draw[|-|, black, very thick] (2.7, 0)--(3, 0);
\draw[|-|, black, very thick] (2.7, 0.2)--(5.2, 0.2) node[midway, above] {$a_h$};
\draw [thick, decorate, decoration={brace, amplitude=10pt, raise=6pt, mirror}] (1, 0) to node[midway, below, yshift=-17pt] {$b$} (3, 0);
\end{tikzpicture}
\end{center}
Clearly, $b^{\prime}$ and $a_j$ are not separated as well.
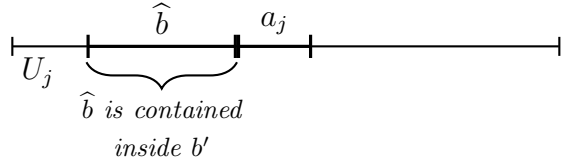
\begin{center}
\begin{tikzpicture}
\node at (0, 0) [below, xshift=10] {$U_j$};
\draw[|-|, black, thick] (0,0)--(7.3,0);
\draw[|-|, black, very thick] (1, 0)--(3, 0) node[midway, above] {$\widehat{b}$};
\draw [thick, decorate, decoration={brace, amplitude=10pt, raise=6pt, mirror}] (1, 0) to node[midway, below, yshift=-13pt, text width=2.5cm, align=center] {\footnotesize{\baselineskip=10pt\textit{$\widehat{b}$ is contained inside $b^{\prime}$}}} (3, 0);
\draw[|-|, black, very thick] (3, 0)--(4, 0) node[midway, above] {$a_j$};
\end{tikzpicture}
\end{center}
Moreover, since $b \in \longmo{U_h} \setminus \lbrace a_h\rbrace$ and $a_h \mapsto a_j$ is an admissible replacement in $U_h$, $a_j$ is not contained in $b^{\prime}$. So, either $b^{\prime}$ and $a_j$ touch at a point, or $b^{\prime}$ and $a_j$ have an overlap.
\begin{center}
\begin{tikzpicture}
\node at (0, 0) [below, xshift=10] {$U_j$};
\draw[|-|, black, thick] (0,0)--(7.3,0);
\draw[|-|, black, very thick] (1, 0)--(2.7, 0) node[midway, above] {$\widehat{b} = b^{\prime}$};
\draw[|-|, black, very thick] (2.7, 0)--(3.7, 0) node[midway, above] {$a_j$};
\end{tikzpicture}

\begin{tikzpicture}
\node at (0, 0) [below, xshift=10] {$U_j$};
\draw[|-|, black, thick] (0,0)--(7.3,0);
\draw[|-, black, very thick] (1, 0)--(2.7, 0) node[midway, above] {$\widehat{b}$};
\draw[|-|, black, very thick] (2.7, 0)--(3, 0);
\draw[|-|, black, very thick] (2.7, 0.2)--(3.7, 0.2) node[midway, above] {$a_j$};
\draw [thick, decorate, decoration={brace, amplitude=10pt, raise=6pt, mirror}] (1, 0) to node[midway, below, yshift=-17pt] {$b^{\prime}$} (3, 0);
\end{tikzpicture}
\end{center}

Since $b \in \longmo{U_h} \subseteq \nfc{U_h}$, by Lemma~\ref{single_image_lemma}, $b^{\prime} = \widehat{b}Y$, where $Y$ is a suffix of $b^{\prime}$. Since $c^{\prime}$ is a right neighbour of $b^{\prime}$, $c^{\prime}$ starts from the right of the beginning of $b^{\prime}$. Therefore, $c^{\prime}$ starts from the right of the beginning of $\widehat{b}$, and $\widehat{b}$ is not contained in $c^{\prime}$.

Let us show that either $c^{\prime}$ ends from the right of the end of $a_j$, or $c^{\prime}$ and $a_j$ have the same end point. Assume the contrary, namely, that $c^{\prime}$ ends from the left of the end of $a_j$. Then, obviously, $c^{\prime}$ starts from the left of the end of $a_j$ as well. First consider the case when $c^{\prime}$ starts from the left of the beginning of $a_j$.
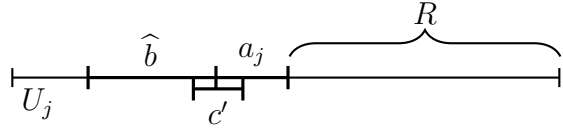
\begin{center}
\begin{tikzpicture}
\node at (0, 0) [below, xshift=10] {$U_j$};
\draw[|-|, black, thick] (0,0)--(7.3,0);
\draw[|-, black, very thick] (1, 0)--(2.7, 0) node[midway, above] {$\widehat{b}$};
\draw[|-|, black, very thick] (2.7, 0)--(3.7, 0) node[midway, above] {$a_j$};
\draw[|-|, black, very thick] (2.4, -0.15)--(3.1, -0.15) node[midway, below] {$c^{\prime}$};
\draw [thick, decorate, decoration={brace, amplitude=10pt, raise=8pt,}] (3.7, 0) to node[midway, above, yshift=16pt] {$R$} (7.3, 0);
\end{tikzpicture}
\end{center}
We see that $\widehat{b}$ and $a_j$ cover $c^{\prime}$. By definition, $b^{\prime}$ contains $\widehat{b}$ and every image of $a_h$ in $U_j$ contains $a_j$. Hence, $b^{\prime}$ and an arbitrary image of $a_h$ in $U_j$ cover $c^{\prime}$. On the other hand, $b^{\prime} \neq c^{\prime}$ and $c^{\prime}$ is not equal to any image of $a_h$ in $U_j$. Therefore, $c^{\prime} \in \fc{U_j}$. In particular, $c^{\prime} \notin \longmo{U_j}$, a contradiction.

Now assume $c^{\prime}$ ends from the left of the end of $a_j$, and $c^{\prime}$ either starts from the right of the beginning of $a_j$, or $c^{\prime}$ and $a_j$ have the same starting point.
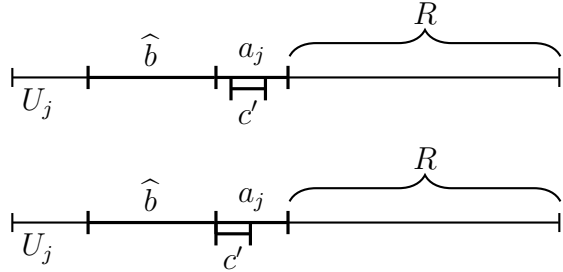
\begin{center}
\begin{tikzpicture}
\node at (0, 0) [below, xshift=10] {$U_j$};
\draw[|-|, black, thick] (0,0)--(7.3,0);
\draw[|-, black, very thick] (1, 0)--(2.7, 0) node[midway, above] {$\widehat{b}$};
\draw[|-|, black, very thick] (2.7, 0)--(3.7, 0) node[midway, above] {$a_j$};
\draw[|-|, black, very thick] (2.9, -0.15)--(3.4, -0.15) node[midway, below] {$c^{\prime}$};
\draw [thick, decorate, decoration={brace, amplitude=10pt, raise=8pt,}] (3.7, 0) to node[midway, above, yshift=16pt] {$R$} (7.3, 0);
\end{tikzpicture}

\begin{tikzpicture}
\node at (0, 0) [below, xshift=10] {$U_j$};
\draw[|-|, black, thick] (0,0)--(7.3,0);
\draw[|-, black, very thick] (1, 0)--(2.7, 0) node[midway, above] {$\widehat{b}$};
\draw[|-|, black, very thick] (2.7, 0)--(3.7, 0) node[midway, above] {$a_j$};
\draw[|-|, black, very thick] (2.7, -0.15)--(3.2, -0.15) node[midway, below] {$c^{\prime}$};
\draw [thick, decorate, decoration={brace, amplitude=10pt, raise=8pt,}] (3.7, 0) to node[midway, above, yshift=16pt] {$R$} (7.3, 0);
\end{tikzpicture}
\end{center}
Then $c^{\prime}$ is contained in $a_j$. So, $c^{\prime} \notin \mo{U_j}$, a contradiction. Therefore, $c^{\prime}$ ends either from the right of the end of $a_j$, or $c^{\prime}$ and $a_j$ have the same end point.

First consider the case when $c^{\prime}$ and $a_j$ have the same end point. Since $c^{\prime} \in \mo{U_j}$ and $c^{\prime} \neq a_j$, $c^{\prime}$ starts strictly from the left of the beginning of $a_j$ in this case. Since $\widehat{b}$ is not contained in $c^{\prime}$, we see that $\widehat{b}$ and $c^{\prime}$ have an overlap $d$, and $d$ is a small piece. Clearly, in this case $a_j$ is contained in $c^{\prime}$. Then $a_j$ is a small piece (see Section~\ref{mt_configurations},~\ref{a_j_keep_structure}).
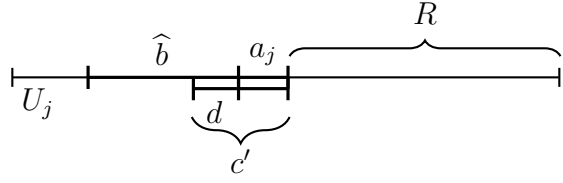
\begin{center}
\begin{tikzpicture}
\node at (0, 0) [below, xshift=10] {$U_j$};
\draw[|-|, black, thick] (0,0)--(7.3,0);
\draw[|-, black, very thick] (1, 0)--(3, 0) node[midway, above] {$\widehat{b}$};
\draw[|-|, black, very thick] (3, 0)--(3.7, 0) node[midway, above] {$a_j$};
\draw[|-, black, very thick] (2.4, -0.15)--(3, -0.15) node[midway, below] {$d$};
\draw[|-|, black, very thick] (3, -0.15)--(3.7, -0.15);
\draw [thick, decorate, decoration={brace, amplitude=10pt, raise=15pt, mirror}] (2.4, 0) to node[midway, below, yshift=-23pt] {$c^{\prime}$} (3.7, 0);
\draw [thick, decorate, decoration={brace, amplitude=6pt, raise=8pt}] (3.7, 0) to node[midway, above, yshift=16pt] {$R$} (7.3, 0);
\end{tikzpicture}
\end{center}
Therefore, $\SPM(c^{\prime}) \leqslant \SPM(d) + \SPM(a_j) \leqslant 2 \leqslant 4$.

Now assume that $c^{\prime}$ ends from the right of the end of $a_j$. Then $c^{\prime}$ has a non-empty intersection with $R$. Denote by $\widehat{c}$ the intersection of $c^{\prime}$ and $R$. First assume that either $c^{\prime}$ starts from the left of the beginning of $a_j$, or $c^{\prime}$ and $a_j$ have the same starting point. Then $a_j$ is contained in $c^{\prime}$. So, as above, $a_j$ is a small piece. Since $\widehat{b}$ is not contained in $c^{\prime}$, there are the following possibilities:
\begin{itemize}
\item
$\widehat{b}$ and $c^{\prime}$ touch at a point;
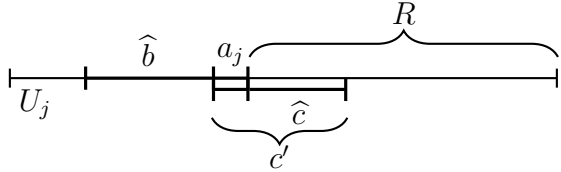
\begin{center}
\begin{tikzpicture}
\node at (0, 0) [below, xshift=10] {$U_j$};
\draw[|-|, black, thick] (0,0)--(7.3,0);
\draw[|-, black, very thick] (1, 0)--(2.7, 0) node[midway, above] {$\widehat{b}$};
\draw[|-|, black, very thick] (2.7, 0)--(3.2, 0) node[midway, above] {$a_j$};
\draw[|-|, black, very thick] (2.7, -0.15)--(3.2, -0.15);
\draw[-|, black, very thick] (3.2, -0.15)--(4.5, -0.15) node[midway, below] {$\widehat{c}$};
\draw [thick, decorate, decoration={brace, amplitude=10pt, raise=15pt, mirror}] (2.7, 0) to node[midway, below, yshift=-21pt] {$c^{\prime}$} (4.5, 0);
\draw [thick, decorate, decoration={brace, amplitude=10pt, raise=8pt,}] (3.2, 0) to node[midway, above, yshift=16pt] {$R$} (7.3, 0);
\end{tikzpicture}
\end{center}
\item
$\widehat{b}$ and $c^{\prime}$ have an overlap $d$, $d$ is a small piece;
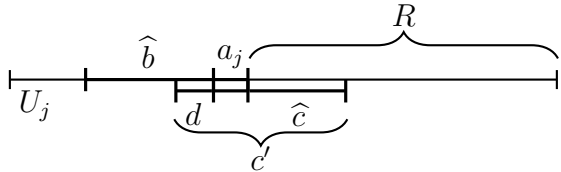
\begin{center}
\begin{tikzpicture}
\node at (0, 0) [below, xshift=10] {$U_j$};
\draw[|-|, black, thick] (0,0)--(7.3,0);
\draw[|-, black, very thick] (1, 0)--(2.7, 0) node[midway, above] {$\widehat{b}$};
\draw[|-|, black, very thick] (2.7, 0)--(3.2, 0) node[midway, above] {$a_j$};
\draw[|-|, black, very thick] (2.7, -0.15)--(3.2, -0.15);
\draw[|-, black, very thick] (2.2, -0.15)--(2.7, -0.15) node[midway, below] {$d$};
\draw[-|, black, very thick] (3.2, -0.15)--(4.5, -0.15) node[midway, below] {$\widehat{c}$};
\draw [thick, decorate, decoration={brace, amplitude=10pt, raise=15pt, mirror}] (2.2, 0) to node[midway, below, yshift=-21pt] {$c^{\prime}$} (4.5, 0);
\draw [thick, decorate, decoration={brace, amplitude=10pt, raise=8pt,}] (3.2, 0) to node[midway, above, yshift=16pt] {$R$} (7.3, 0);
\end{tikzpicture}
\end{center}
\end{itemize}
Assume $\SPM(\widehat{c}) \geqslant 3$. Then it follows from Lemma~\ref{non_covered_inverse_image} that there exists $c \in \longmo{U_h}\setminus \lbrace a_h \rbrace$ such that $c^{\prime}$ is an image of $c$. Since $a_j$ is contained in $c^{\prime}$, $a_j$ is fully covered by an image of the element $c \in \longmo{U_h} \setminus \lbrace a_h\rbrace$. However, by the conditions of Lemma~\ref{neighbours_a_h_not_separated_right}, $a_h \mapsto a_j$ is an admissible replacement in $U_h$, a contradiction. Therefore, $\SPM(\widehat{c}) \leqslant 2$. So, $\SPM(c^{\prime}) \leqslant \SPM(\widehat{c}) + \SPM(a_j) + 1 \leqslant 2 + 1 + 1 \leqslant 4$.

It remains to consider the case when $c^{\prime}$ starts from the right of the beginning of $a_j$ and ends from the right of the end of $a_j$. Since $b \in \longmo{U_h} \setminus \lbrace a_h\rbrace$ and $a_h \mapsto a_j$ is an admissible replacement in $U_h$, $a_j$ is not contained in $b^{\prime}$. So, $b^{\prime}$ ends from the left of the end of $a_j$. Since $c^{\prime}$ is a right neighbour of $b^{\prime}$ in $U_j$, $b^{\prime}$ and $c^{\prime}$ are not separated in $U_j$. Therefore, since $b^{\prime}$ ends from the left of the end of $a_j$ and $c^{\prime}$ starts from the right of the beginning of $a_j$, we obtain that the end point of $b^{\prime}$ is contained strictly inside $a_j$ and the beginning point of $c^{\prime}$ is contained strictly inside $a_j$. So, we obtain the following possibilities:
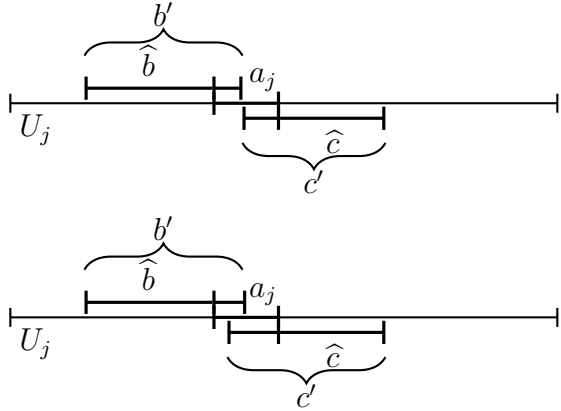
\begin{center}
\begin{tikzpicture}
\node at (0, 0) [below, xshift=10] {$U_j$};
\draw[|-|, black, thick] (0,0)--(7.3,0);
\draw[|-, black, very thick] (1, 0.2)--(2.7, 0.2) node[midway, above] {$\widehat{b}$};
\draw [thick, decorate, decoration={brace, amplitude=10pt, raise=18pt}] (1, 0) to node[midway, above, yshift=25pt] {$b^{\prime}$} (3.1, 0);
\draw[|-|, black, very thick] (2.7, 0.2)--(3.1, 0.2);
\draw[|-|, black, very thick] (2.7, 0)--(3.6, 0) node[near end, above] {$a_j$};
\draw[|-|, black, very thick] (3.1, -0.2)--(3.6, -0.2);
\draw[-|, black, very thick] (3.6, -0.2)--(5, -0.2) node[midway, below] {$\widehat{c}$};
\draw [thick, decorate, decoration={brace, amplitude=10pt, raise=15pt, mirror}] (3.1, 0) to node[midway, below, yshift=-21pt] {$c^{\prime}$} (5, 0);
\end{tikzpicture}

\begin{tikzpicture}
\node at (0, 0) [below, xshift=10] {$U_j$};
\draw[|-|, black, thick] (0,0)--(7.3,0);
\draw[|-, black, very thick] (1, 0.2)--(2.7, 0.2) node[midway, above] {$\widehat{b}$};
\draw [thick, decorate, decoration={brace, amplitude=10pt, raise=18pt}] (1, 0) to node[midway, above, yshift=25pt] {$b^{\prime}$} (3.1, 0);
\draw[|-|, black, very thick] (2.7, 0.2)--(3.15, 0.2);
\draw[|-|, black, very thick] (2.7, 0)--(3.6, 0) node[near end, above] {$a_j$};
\draw[|-|, black, very thick] (2.9, -0.2)--(3.6, -0.2);
\draw[-|, black, very thick] (3.6, -0.2)--(5, -0.2) node[midway, below] {$\widehat{c}$};
\draw [thick, decorate, decoration={brace, amplitude=10pt, raise=15pt, mirror}] (2.9, 0) to node[midway, below, yshift=-21pt] {$c^{\prime}$} (5, 0);
\end{tikzpicture}
\end{center}
We see that $b^{\prime}$ and $c^{\prime}$ cover $a_j$, wherein $b^{\prime} \neq a_j$ and $c^{\prime} \neq a_j$. Assume $\SPM(\widehat{c}) \geqslant 3$. Then, by Lemma~\ref{non_covered_inverse_image}, there exists $c \in \longmo{U_h} \setminus \lbrace a_h \rbrace$ such that $c^{\prime}$ is an image of $c$. Since $a_j$ is covered by $b^{\prime}$ and $c^{\prime}$ in $U_j$, we see that $a_j$ is fully covered by images of elements of $\longmo{U_h} \setminus \lbrace a_h\rbrace$. By the conditions of Lemma~\ref{neighbours_a_h_not_separated_right}, $a_h \mapsto a_j$ is an admissible replacement in $U_h$, a contradiction. Therefore, $\SPM(\widehat{c}) \leqslant 2$. Hence, $\SPM(c^{\prime}) \leqslant \SPM(\widehat{c}) + 1 \leqslant 3 \leqslant 4$.
\end{proof}

\begin{corollary}
\label{short_neighbours_right}
Let $U_h$ be a monomial, $b \in \longmo{U_h}$. Suppose $a_h \in \mo{U_h}$ is a right neighbour of $b$. Suppose $a_h \mapsto a_j$ is an admissible replacement in $U_h$. Let $U_j$ be the resulting monomial. Denote by $b^{\prime}$ the image of $b$ in $U_j$. Then every right neighbour of $b^{\prime}$ in $U_j$ is of $\SPM$-measure $\leqslant \max(4, \SPM(a_j))$.
\end{corollary}
\begin{proof}
Since $a_h$ and $b$ are not separated in $U_h$, clearly, $a_j$ and $b^{\prime}$ are not separated in $U_j$ as well.

Let $a_j \in \longmo{U_j}$. Then, in particular, $a_j \in \mo{U_h}$. So, by definition, $a_j$ is  a right neighbour of $b^{\prime}$ in $U_j$. Since $a_j \in \longmo{U_j}$, it follows from Corollary~\ref{single_long_neighbour} that all the rest of right neighbours of $b^{\prime}$ (if there are any) belong to $\shortmo{U_j}$. Therefore, in this case every right neighbour of $b^{\prime}$ is of $\SPM$-measure $\leqslant \SPM(a_j)$.

Let $a_j \notin \longmo{U_j}$. First assume that all the right neighbours of $b^{\prime}$ belong to $\shortmo{U_j}$. Then, by definition, every right neighbour of $b^{\prime}$ is of $\SPM$-measure $\leqslant 2$. Now assume there exists $c^{\prime}$ a right neighbour of $b^{\prime}$ such that $c^{\prime} \in \longmo{U_j}$. Then, by Lemma~\ref{neighbours_a_h_not_separated_right}, we obtain $\SPM(c^{\prime}) \leqslant 4$.

So, combining all the above, we see that every right neighbour of $b^{\prime}$ in $U_j$ is of $\SPM$-measure $\leqslant \max(4, \SPM(a_j))$.
\end{proof}

\begin{lemma}
\label{neighbours_a_h_not_separated_left}
Let $U_h$ be a monomial, $b \in \longmo{U_h}$. Suppose $a_h \in \mo{U_h}$ is a left neighbour of $b$. Suppose $a_h \mapsto a_j$ is an admissible replacement in $U_h$. Let $U_j$ be the resulting monomial. Denote by $b^{\prime}$ the image of $b$ in $U_j$. Let $a_j \notin \longmo{U_j}$, that is, either $a_j \in \shortmo{U_j}$, or $a_j \notin \mo{U_j}$. Assume that there exists $c^{\prime}$ a left neighbour of $b^{\prime}$ such that $c^{\prime} \in \longmo{U_j}$. Then $\SPM(c^{\prime}) \leqslant 4$.
\end{lemma}
\begin{proof}
This lemma is proved in the same way as Lemma~\ref{neighbours_a_h_not_separated_right}.
\end{proof}

\begin{corollary}
\label{short_neighbours_left}
Let $U_h$ be a monomial, $b \in \longmo{U_h}$. Suppose $a_h \in \mo{U_h}$ is a left neighbour of $b$. Suppose $a_h \mapsto a_j$ is an admissible replacement in $U_h$. Let $U_j$ be the resulting monomial. Denote by $b^{\prime}$ the image of $b$ in $U_j$. Then every left neighbour of $b^{\prime}$ in $U_j$ is of $\SPM$-measure $\leqslant \max(4, \SPM(a_j))$.
\end{corollary}
\begin{proof}
This corollary is proved in the same way as Corollary~\ref{short_neighbours_right}.
\end{proof}

\begin{lemma}
\label{neighbours_a_h_separated1_right}
Let $U_h$ be a monomial, $b \in \longmo{U_h}$. Suppose $a_h \in \mo{U_h}$ starts from the right of the beginning of $b$ and $a_h$ is separated from $b$. Suppose $a_h \mapsto a_j$ is an admissible replacement in $U_h$. Let $U_j$ be the resulting monomial. Denote by $b^{\prime}$ the image of $b$ in $U_j$. Assume there exists $c^{\prime} \in \mo{U_h}$  a right neighbour of $b^{\prime}$ such that $\SPM(c^{\prime}) \geqslant 6$. Then $a_j$ is not contained in $c^{\prime}$. That is, either $a_j$ and $c^{\prime}$ are separated, or they touch at a point, or they have an overlap.
\end{lemma}
\begin{proof}
Assume $U_h = La_hR$, then $U_j = La_jR$. Let $\widehat{b}$ be the intersection of $b$ and $L$ in $U_h$. As earlier, we consider $\widehat{b}$ as an occurrence in $U_j$. Since $a_h$ and $b$ are separated, $b$ stays unchanged in $U_j$, that is, $b = \widehat{b}$ in $U_h$ and $b^{\prime} = \widehat{b}$ in $U_j$. Hence, one can easily see that $a_j$ and $b^{\prime}$ are separated in $U_j$. So, $a_j$ is not a neighbour of $b^{\prime}$, and $b^{\prime}$ ends from the left of the beginning of $a_j$.
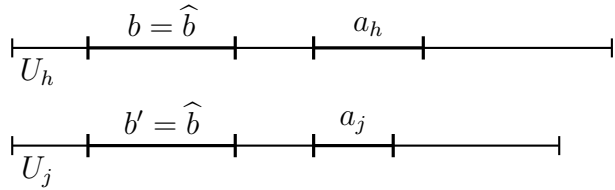
\begin{center}
\begin{tikzpicture}
\node at (0, 0) [below, xshift=10] {$U_h$};
\draw[|-|, black, thick] (0,0)--(8,0);
\draw[|-|, black, very thick] (1, 0)--(3, 0) node[midway, above] {$b = \widehat{b}$};
\draw[|-|, black, very thick] (4, 0)--(5.5, 0) node[midway, above] {$a_h$};

\node at (0, -1.3) [below, xshift=10] {$U_j$};
\draw[|-|, black, thick] (0,-1.3)--(7.3,-1.3);
\draw[|-|, black, very thick] (1, -1.3)--(3, -1.3) node[midway, above] {$b^{\prime} = \widehat{b}$};
\draw[|-|, black, very thick] (4, -1.3)--(5.1, -1.3) node[midway, above] {$a_j$};
\end{tikzpicture}
\end{center}

Let $\widehat{c}_1$ be the intersection of $c^{\prime}$ and $L$ in $U_j$ and $\widehat{c}_2$ be the intersection of $c^{\prime}$ and $R$ in $U_j$. Since $c^{\prime}$ and $b^{\prime}$ are not separated and $b^{\prime}$ ends from the left of the beginning of $a_j$, we see that $c^{\prime}$ starts from the left of the beginning of $a_j$. Hence, $\widehat{c}_1$ is always non-empty.

Assume $a_j$ is contained in $c^{\prime}$. Then $a_j$ is a small piece (see Section~\ref{mt_configurations},~\ref{a_j_keep_structure}). There are the following possibilities:
\begin{itemize}
\item
$a_j$ is a terminal subword of  $c^{\prime}$;
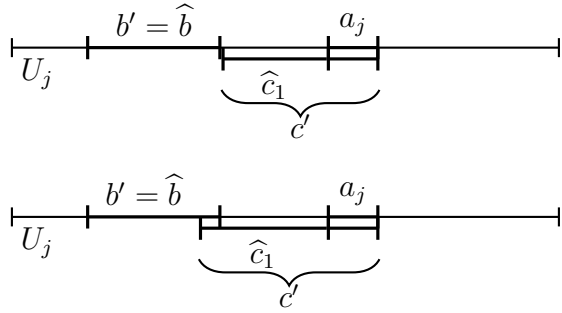
\begin{center}
\begin{tikzpicture}
\node at (-0.3, 0) [below, xshift=10] {$U_j$};
\draw[|-|, black, thick] (-0.3,0)--(7,0);
\draw[|-|, black, very thick] (0.7, 0)--(2.5, 0) node[midway, above] {$b^{\prime} = \widehat{b}$};
\draw[|-, black, very thick] (2.5, -0.15)--(3.9, -0.15) node[midway, below] {$\widehat{c}_1$};
\draw[|-|, black, very thick] (3.9, -0.15)--(4.6, -0.15);
\draw[|-|, black, very thick] (3.9, 0)--(4.6, 0) node[midway, above] {$a_j$};
\draw [thick, decorate, decoration={brace, amplitude=10pt, raise=16pt, mirror}] (2.5, 0) to node[midway, below, yshift=-21pt] {$c^{\prime}$} (4.6, 0);
\end{tikzpicture}

\begin{tikzpicture}
\node at (-0.3, 0) [below, xshift=10] {$U_j$};
\draw[|-|, black, thick] (-0.3,0)--(7,0);
\draw[|-|, black, very thick] (0.7, 0)--(2.5, 0) node[midway, above, xshift=-4] {$b^{\prime} = \widehat{b}$};
\draw[|-, black, very thick] (2.2, -0.15)--(3.9, -0.15) node[midway, below] {$\widehat{c}_1$};
\draw[|-|, black, very thick] (3.9, -0.15)--(4.6, -0.15);
\draw[|-|, black, very thick] (3.9, 0)--(4.6, 0) node[midway, above] {$a_j$};
\draw [thick, decorate, decoration={brace, amplitude=10pt, raise=16pt, mirror}] (2.2, 0) to node[midway, below, yshift=-21pt] {$c^{\prime}$} (4.6, 0);
\end{tikzpicture}
\end{center}
\item
$a_j$ is not a terminal subword of $c^{\prime}$;
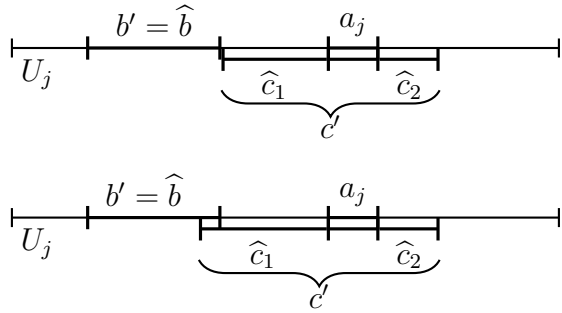
\begin{center}
\begin{tikzpicture}
\node at (-0.3, 0) [below, xshift=10] {$U_j$};
\draw[|-|, black, thick] (-0.3,0)--(7,0);
\draw[|-|, black, very thick] (0.7, 0)--(2.5, 0) node[midway, above] {$b^{\prime} = \widehat{b}$};
\draw[|-, black, very thick] (2.5, -0.15)--(3.9, -0.15) node[midway, below] {$\widehat{c}_1$};
\draw[|-|, black, very thick] (3.9, -0.15)--(4.6, -0.15);
\draw[-|, black, very thick] (4.6, -0.15)--(5.4, -0.15) node[midway, below] {$\widehat{c}_2$};
\draw[|-|, black, very thick] (3.9, 0)--(4.6, 0) node[midway, above] {$a_j$};
\draw [thick, decorate, decoration={brace, amplitude=10pt, raise=16pt, mirror}] (2.5, 0) to node[midway, below, yshift=-21pt] {$c^{\prime}$} (5.4, 0);
\end{tikzpicture}

\begin{tikzpicture}
\node at (-0.3, 0) [below, xshift=10] {$U_j$};
\draw[|-|, black, thick] (-0.3,0)--(7,0);
\draw[|-|, black, very thick] (0.7, 0)--(2.5, 0) node[midway, above, xshift=-4] {$b^{\prime} = \widehat{b}$};
\draw[|-, black, very thick] (2.2, -0.15)--(3.9, -0.15) node[midway, below] {$\widehat{c}_1$};
\draw[|-|, black, very thick] (3.9, -0.15)--(4.6, -0.15);
\draw[|-|, black, very thick] (3.9, 0)--(4.6, 0) node[midway, above] {$a_j$};
\draw[-|, black, very thick] (4.6, -0.15)--(5.4, -0.15) node[midway, below] {$\widehat{c}_2$};
\draw [thick, decorate, decoration={brace, amplitude=10pt, raise=16pt, mirror}] (2.2, 0) to node[midway, below, yshift=-21pt] {$c^{\prime}$} (5.4, 0);
\end{tikzpicture}
\end{center}
\end{itemize}
Since $a_j$ is a small piece, $\SPM(a_j) = 1$. Assume $\SPM(\widehat{c}_1) \leqslant 2$ and $\SPM(\widehat{c}_2) \leqslant 2$. Then $\SPM(c^{\prime}) \leqslant \SPM(\widehat{c}_1) + \SPM(\widehat{c}_2) + \SPM(a_j) \leqslant 2 + 2 + 1 = 5$. A contradiction with the condition $\SPM(c^{\prime}) \geqslant 6$. Therefore, $\SPM(\widehat{c}_1) \geqslant 3$ or $\SPM(\widehat{c}_2) \geqslant 3$. Assume that $\SPM(\widehat{c}_1) \geqslant 3$. Then it follows from Lemma~\ref{non_covered_inverse_image} that there exists $c \in \longmo{U_h}$ such that $c^{\prime}$ is an image of $c$.  Since $a_j$ is covered by $c^{\prime}$, we see that $a_j$ is fully covered by an image of an element $c \in \longmo{U_h} \setminus \lbrace a_h\rbrace$. However, by the conditions of Lemma~\ref{neighbours_a_h_separated1_right}, $a_h \mapsto a_j$ is an admissible replacement in $U_h$, a contradiction. Now assume that $\SPM(\widehat{c}_2) \geqslant 3$. Then, by the same argument, we obtain that $a_h \mapsto a_j$ is not an admissible replacement in $U_h$. This contradiction completes the proof.
\end{proof}

\begin{lemma}
\label{neighbours_a_h_separated1_left}
Let $U_h$ be a monomial, $b \in \longmo{U_h}$. Suppose $a_h \in \mo{U_h}$ starts from the left of the beginning of $b$ and $a_h$ is separated from $b$. Suppose $a_h \mapsto a_j$ is an admissible replacement in $U_h$. Let $U_j$ be the resulting monomial. Denote by $b^{\prime}$ the image of $b$ in $U_j$. Assume there exists $c^{\prime} \in \mo{U_h}$  a left neighbour of $b^{\prime}$ such that $\SPM(c^{\prime}) \geqslant 6$. Then $a_j$ is not contained in $c^{\prime}$. That is, either $a_j$ and $c^{\prime}$ are separated, or they touch at a point, or they have an overlap.
\end{lemma}

\begin{definition}
\label{admissible_sequence_def}
Let $U$ be a monomial, $b \in \longmo{U}$. Assume $U^{(1)},\ldots, U^{(K)}$, are monomials, $K \geqslant 2$, and
\begin{equation*}
U^{(1)} \mapsto \ldots \mapsto U^{(K)}
\end{equation*}
is a sequence of replacements of the following form:
\begin{itemize}
\item
$U = U^{(1)}$;
\item
$a_h^{(k)} \in \mo{U^{(k)}}$,  $a_h^{(k)}$ and $a_j^{(k)}$ are incident monomials and $a_h^{(k)} \mapsto a_j^{(k)}$ is an admissible replacement in $U^{(k)}$, $U^{(k + 1)}$ is the resulting monomial, $k = 1, \ldots, K - 1$;
\item
$b = b^{(1)}$;
\item
$b^{(k + 1)}$ is an image of $b^{(k)}$ in $U^{(k + 1)}$, $k = 1, \ldots, K - 1$;
\item
$a_h^{(k)} \neq b^{(k)}$, $k = 1, \ldots, K - 1$.
\end{itemize}
Then a triple that consists of
\begin{enumerate}
\item
the sequence $U = U^{(1)} \mapsto \ldots \mapsto U^{(K)}$,
\item
the tuple of pairs $((a_h^{(1)}, a_j^{(1)}), \ldots, (a_h^{(K - 1)}, a_j^{(K - 1)}))$,
\item
the tuple $(b = b^{(1)}, \ldots, b^{(K)})$
\end{enumerate}
is called \emph{$(b, U)$-admissible sequence}.
\end{definition}

\begin{remark}
\label{incident_same_indices_remark}
In Definition~\ref{admissible_sequence_def} we use same indices $h$ and $j$ for replacements $a_h^{(k_1)} \mapsto a_j^{(k_1)}$ and $a_h^{(k_2)} \mapsto a_j^{(k_2)}$ even if $k_1 \neq k_2$. Recall that, speaking formally, $h$ and $j$ are numbers of $a_h$ and $a_j$ in a polynomial of~$\Rel$. So, to be formal, we have to use $h(k)$ and $j(k)$ that depend on $k$. However, we never fixed any precise order of monomials in the polynomials of~$\Rel$. Throughout the paper we use indices $h$ and $j$ just to emphasise that $a_h$ and $a_j$ are incident monomials and we never focus on their precise numbers in a polynomial of~$\Rel$ in our proofs. That is, using same indices $h$ and $j$ for different pairs of incident monomials does not cause any additional restrictions to these monomials. Also, we do not want to make our denotations too bulky, so, we continue using same indices $h$ and $j$ for different pairs of incident monomials.
\end{remark}

Assume $U$ is a monomial, $b \in \longmo{U}$. Then there are the following configurations for left neighbours of $b$:
\begin{enumerate}[label=(L\arabic*)]
\item
\label{long_left_neighbour}
there exists a left neighbour of $b$ of $\SPM$-measure $\geqslant \tau - 3$;
\item
\label{short_left_neighbour}
all the left neighbours of $b$ are of $\SPM$-measure $< \tau - 3$;
\item
\label{no_left_neighbour}
$b$ does not have left neighbours at all.
\end{enumerate}

Let $t(b)$\label{terminal_subword_def} be a terminal subword of $b$ defined by the following rule. Assume we have configuration~\ref{long_left_neighbour} for $b$ and $c\in \mo{U}$ is a left neighbour of $b$ with $\SPM(c) \geqslant \tau - 3$. Then $t(b)$ is the suffix of $b$ such that its beginning point is equal to the end point of $c$.
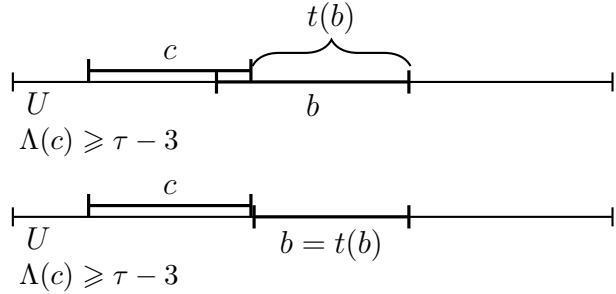
\begin{center}
\begin{tikzpicture}
\node at (0, 0) [below, xshift=10] {$U$};
\draw[|-|, black, thick] (0,0)--(8,0);
\draw[|-|, black, very thick] (1, 0.15)--(3.2, 0.15) node[midway, above] {$c$};
\draw[|-|, black, very thick] (2.7, 0)--(5.3, 0) node[midway, below] {$b$};
\draw [thick, decorate, decoration={brace, amplitude=10pt, raise=6pt}] (3.2, 0) to node[midway, above, yshift=14pt] {$t(b)$} (5.3, 0);
\node[text width=3cm, align=left] at (1.6, -0.8) {\small\baselineskip=10pt $\SPM(c) \geqslant \tau - 3$};
\end{tikzpicture}

\begin{tikzpicture}
\node at (0, 0) [below, xshift=10] {$U$};
\draw[|-|, black, thick] (0,0)--(8,0);
\draw[|-|, black, very thick] (1, 0.15)--(3.2, 0.15) node[midway, above] {$c$};
\draw[|-|, black, very thick] (3.2, 0)--(5.3, 0) node[midway, below] {$b = t(b)$};
\node[text width=3cm, align=left] at (1.6, -0.8) {\small\baselineskip=10pt $\SPM(c) \geqslant \tau - 3$};
\end{tikzpicture}
\end{center}
If we have configuration~\ref{short_left_neighbour} or~\ref{no_left_neighbour} for $b$, then we put $t(b)$ to be equal to $b$. So, in general $b = pt(b)$, where $p$ is a small piece (possibly empty).

\begin{lemma}
\label{admissible_replacements_left}
Let $U$ be a monomial. Assume $b \in \longmo{U}$ and, moreover, $\SPM(t(b)) \geqslant 3$. Assume
\begin{equation*}
(U = U^{(1)} \mapsto \ldots \mapsto U^{(K)}, ((a_h^{(1)}, a_j^{(1)}), \ldots, (a_h^{(K - 1)}, a_j^{(K - 1)})), (b = b^{(1)}, \ldots, b^{(K)}))
\end{equation*}
is a $(b, U)$-admissible sequence such that every $a_h^{(k)}$ starts from the left of the beginning of $b^{(k)}$, $k = 1, \ldots, K - 1$. Then $b^{(k)}$ is a single image of $b^{(k - 1)}$ in $U^{(k)}$, $k = 2, \ldots, K$, and
\begin{align*}
b^{(k)} = p^{(k)}t(b), &\textit{ where } p^{(k)} \textit{ is a prefix of } b^{(k)},\\
&p^{(k)} \textit{ is a small piece (possibly empty), } k = 1, \ldots, K.
\end{align*}
Moreover, if we start with~\ref{short_left_neighbour} or~\ref{no_left_neighbour} for $b$, then $a_h^{(k - 1)}$ is separated from $b^{(k - 1)}$ in $U^{(k - 1)}$, $k = 2, \ldots, K$, and $b^{(k)} = t(b)  = b$, $k = 1, \ldots, K$ (that is, $p^{(k)}$ is empty).

Under the same conditions, we have
\begin{equation*}
U^{(k)} = A^{(k)}b^{(k)}B = A^{(k)}p^{(k)}t(b)B, \ k = 1, \ldots, K,
\end{equation*}
where $A^{(k)}$ is a prefix of $U^{(k)}$ and $B$ is a suffix of $U^{(k)}$. That is, the suffix of $U^{(k)}$ that starts at the beginning point of $t(b)$ does not depend on $k$. Moreover, $a_h^{(k)}$ is a maximal occurrence in the subword $A^{(k)}p^{(k)}$, $k = 1, \ldots, K - 1$.
\end{lemma}
\begin{proof}
Consider the first replacement in the sequence, namely, $a_h^{(1)} \mapsto a_j^{(1)}$ in $U^{(1)}$, $U^{(2)}$ is a resulting monomial. Assume $U^{(1)} = La_h^{(1)}R$. Recall that, by the condition of Lemma~\ref{admissible_replacements_left}, $a_h^{(1)}$ starts from the left of the beginning of $b^{(1)}$. So, the intersection of $b^{(1)}$ and $R$ is non-empty and we denote it by $\widehat{b}^{(1)}$. As earlier, we can consider $\widehat{b}^{(1)}$ as an occurrence in the resulting monomial $U^{(2)} = La_j^{(1)}R$.
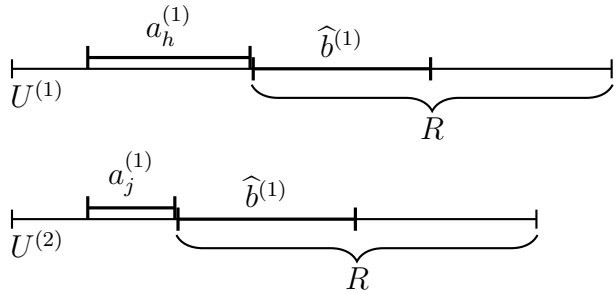
\begin{center}
\begin{tikzpicture}
\node at (0, 0) [below, xshift=10] {$U^{(1)}$};
\draw[|-|, black, thick] (0,0)--(8,0);
\draw[|-|, black, very thick] (1, 0.15)--(3.2, 0.15) node[midway, above] {$a_h^{(1)}$};
\draw[|-|, black, very thick] (3.2, 0)--(5.6, 0) node[midway, above] {$\widehat{b}^{(1)}$};
\draw [thick, decorate, decoration={brace, amplitude=10pt, raise=6pt, mirror}] (3.2, 0) to node[midway, below, yshift=-14pt] {$R$} (8, 0);

\node at (0, -2) [below, xshift=10] {$U^{(2)}$};
\draw[|-|, black, thick] (0, -2)--(7, -2);
\draw[|-|, black, very thick] (1, 0.15 - 2)--(2.2, 0.15 - 2) node[midway, above] {$a_j^{(1)}$};
\draw[|-|, black, very thick] (2.2, -2)--(4.6, -2) node[midway, above] {$\widehat{b}^{(1)}$};
\draw [thick, decorate, decoration={brace, amplitude=10pt, raise=6pt, mirror}] (2.2, -2) to node[midway, below, yshift=-14pt] {$R$} (7, -2);
\end{tikzpicture}

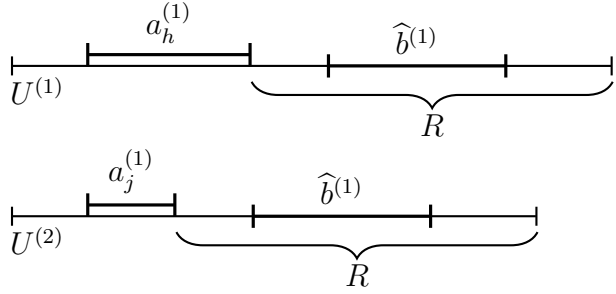
\begin{tikzpicture}
\node at (0, 0) [below, xshift=10] {$U^{(1)}$};
\draw[|-|, black, thick] (0,0)--(8,0);
\draw[|-|, black, very thick] (1, 0.15)--(3.2, 0.15) node[midway, above] {$a_h^{(1)}$};
\draw[|-|, black, very thick] (4.2, 0)--(6.6, 0) node[midway, above] {$\widehat{b}^{(1)}$};
\draw [thick, decorate, decoration={brace, amplitude=10pt, raise=6pt, mirror}] (3.2, 0) to node[midway, below, yshift=-14pt] {$R$} (8, 0);

\node at (0, -2) [below, xshift=10] {$U^{(2)}$};
\draw[|-|, black, thick] (0, -2)--(7, -2);
\draw[|-|, black, very thick] (1, 0.15 - 2)--(2.2, 0.15 - 2) node[midway, above] {$a_j^{(1)}$};
\draw[|-|, black, very thick] (3.2, -2)--(5.6, -2) node[midway, above] {$\widehat{b}^{(1)}$};
\draw [thick, decorate, decoration={brace, amplitude=10pt, raise=6pt, mirror}] (2.2, -2) to node[midway, below, yshift=-14pt] {$R$} (7, -2);
\end{tikzpicture}
\end{center}

First assume that $a_h^{(1)}$ is not separated from $b^{(1)}$. Since $\SPM(a_h^{(1)}) \geqslant \tau - 2$, we have configuration~\ref{long_left_neighbour} for $b^{(1)}$. Since $a_h^{(1)}$ is a left neighbour of $b^{(1)}$ of $\SPM$-measure $\geqslant \tau - 3$, it trivially follows from the definition that $t(b^{(1)}) = \widehat{b}^{(1)}$. Since $a_h^{(1)}$ is not separated from $b^{(1)}$, there are the following possibilities:
\begin{itemize}
\item
$a_h^{(1)}$ and $b^{(1)}$ touch at a point. Then $b^{(1)} = t(b^{(1)}) = \widehat{b}^{(1)}$.
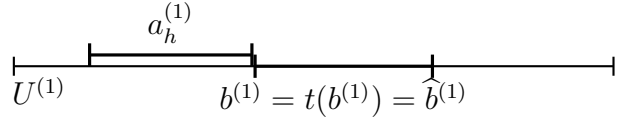
\begin{center}
\begin{tikzpicture}
\node at (0, 0) [below, xshift=10] {$U^{(1)}$};
\draw[|-|, black, thick] (0,0)--(8,0);
\draw[|-|, black, very thick] (1, 0.15)--(3.2, 0.15) node[midway, above] {$a_h^{(1)}$};
\draw[|-|, black, very thick] (3.2, 0)--(5.6, 0) node[midway, below] {$b^{(1)} = t(b^{(1)}) = \widehat{b}^{(1)}$};
\end{tikzpicture}
\end{center}
\item
$a_h^{(1)}$ and $b^{(1)}$ have an overlap $d$, $d$ is a small piece. Then $b^{(1)} = dt(b^{(1)}) = d\widehat{b}^{(1)}$.
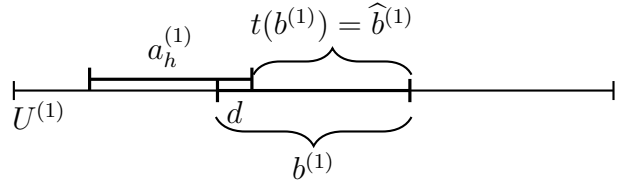
\begin{center}
\begin{tikzpicture}
\node at (0, 0) [below, xshift=10] {$U^{(1)}$};
\draw[|-|, black, thick] (0,0)--(8,0);
\draw[|-|, black, very thick] (1, 0.15)--(3.2, 0.15) node[midway, above] {$a_h^{(1)}$};
\draw[|-|, black, very thick] (2.7, 0)--(5.3, 0);
\path (2.7, 0)--(3.2, 0) node[midway, below] {$d$};
\draw [thick, decorate, decoration={brace, amplitude=10pt, raise=6pt}] (3.2, 0) to node[midway, above, yshift=14pt] {$t(b^{(1)}) = \widehat{b}^{(1)}$} (5.3, 0);
\draw [thick, decorate, decoration={brace, amplitude=11pt, raise=10pt, mirror}] (2.7, 0) to node[midway, below, yshift=-19pt] {$b^{(1)}$} (5.3, 0);
\end{tikzpicture}
\end{center}
\end{itemize}
So, we see that $b^{(1)} = p^{(1)}t(b)$, where $p^{(1)}$ is empty if $a_h^{(1)}$ and $b^{(1)}$ touch at a point, and $p^{(1)} = d$ if $a_h^{(1)}$ and $b^{(1)}$ have an overlap $d$.

Since $a_h^{(1)} \mapsto a_j^{(1)}$ is an admissible replacement in $U^{(1)}$ and $b^{(1)} \in \longmo{U^{(1)}} \setminus \lbrace a_h^{(1)}\rbrace$, $a_j^{(1)}$ is not contained in $b^{(2)}$. So, there are the following possibilities:
\begin{itemize}
\item
$a_j^{(1)}$ and $b^{(2)}$ touch at a point. Then $b^{(2)} = \widehat{b}^{(1)}$;
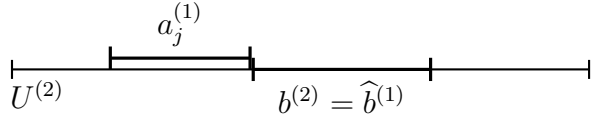
\begin{center}
\begin{tikzpicture}
\node at (0, 0) [below, xshift=10] {$U^{(2)}$};
\draw[|-|, black, thick] (0,0)--(7.7,0);
\draw[|-|, black, very thick] (1.3, 0.15)--(3.2, 0.15) node[midway, above] {$a_j^{(1)}$};
\draw[|-|, black, very thick] (3.2, 0)--(5.6, 0) node[midway, below] {$b^{(2)} = \widehat{b}^{(1)}$};
\end{tikzpicture}
\end{center}
\item
$a_j^{(1)}$ and $b^{(2)}$ have an overlap $d^{\prime}$, $d^{\prime}$ is a small piece. Then $b^{(2)} = d^{\prime}\widehat{b}^{(1)}$.
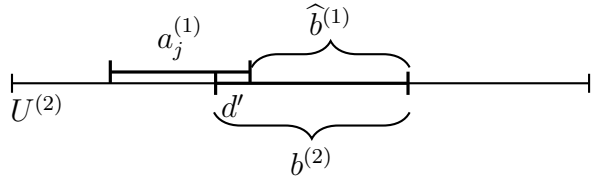
\begin{center}
\begin{tikzpicture}
\node at (0, 0) [below, xshift=10] {$U^{(2)}$};
\draw[|-|, black, thick] (0,0)--(7.7,0);
\draw[|-|, black, very thick] (1.3, 0.15)--(3.2, 0.15) node[midway, above] {$a_j^{(1)}$};
\draw[|-|, black, very thick] (2.7, 0)--(5.3, 0);
\path (2.7, 0)--(3.2, 0) node[midway, below] {$d^{\prime}$};
\draw [thick, decorate, decoration={brace, amplitude=10pt, raise=6pt}] (3.2, 0) to node[midway, above, yshift=14pt] {$\widehat{b}^{(1)}$} (5.3, 0);
\draw [thick, decorate, decoration={brace, amplitude=10pt, raise=11pt, mirror}] (2.7, 0) to node[midway, below, yshift=-20pt] {$b^{(2)}$} (5.3, 0);
\end{tikzpicture}
\end{center}
\end{itemize}
We proved above that $\widehat{b}^{(1)} = t(b^{(1)})$. Hence, $b^{(2)} = p^{(2)}t(b^{(1)})$, where $p^{(2)}$ is empty if $a_j^{(1)}$ and $b^{(2)}$ touch at a point, and $p^{(2)} = d^{\prime}$ if $a_j^{(1)}$ and $b^{(2)}$ have an overlap $d^{\prime}$.

Assume $a_h^{(1)}$ is separated from $b^{(1)}$. Then $b^{(1)} = \widehat{b}^{(1)}$ in $U^{(1)}$ and $b^{(2)} = \widehat{b}^{(1)}$ in $U^{(2)}$.
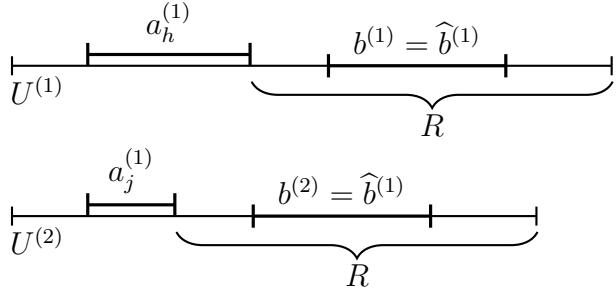
\begin{center}
\begin{tikzpicture}
\node at (0, 0) [below, xshift=10] {$U^{(1)}$};
\draw[|-|, black, thick] (0,0)--(8,0);
\draw[|-|, black, very thick] (1, 0.15)--(3.2, 0.15) node[midway, above] {$a_h^{(1)}$};
\draw[|-|, black, very thick] (4.2, 0)--(6.6, 0) node[midway, above] {$b^{(1)} = \widehat{b}^{(1)}$};
\draw [thick, decorate, decoration={brace, amplitude=10pt, raise=6pt, mirror}] (3.2, 0) to node[midway, below, yshift=-14pt] {$R$} (8, 0);

\node at (0, -2) [below, xshift=10] {$U^{(2)}$};
\draw[|-|, black, thick] (0, -2)--(7, -2);
\draw[|-|, black, very thick] (1, 0.15 - 2)--(2.2, 0.15 - 2) node[midway, above] {$a_j^{(1)}$};
\draw[|-|, black, very thick] (3.2, -2)--(5.6, -2) node[midway, above] {$b^{(2)} = \widehat{b}^{(1)}$};
\draw [thick, decorate, decoration={brace, amplitude=10pt, raise=6pt, mirror}] (2.2, -2) to node[midway, below, yshift=-14pt] {$R$} (7, -2);
\end{tikzpicture}
\end{center}
By the definition of $t(b^{(1)})$, $b^{(1)} = p^{(1)}t(b^{(1)})$, where $p^{(1)}$ is a small piece (possibly empty). Since $b^{(1)} = \widehat{b}^{(1)}$ in $U^{(1)}$ and $b^{(2)} = \widehat{b}^{(1)}$ in $U^{(2)}$, we obtain  $b^{(2)} = p^{(1)}t(b^{(1)})$. Moreover, by the definition of $t(b^{(1)})$, $p^{(1)}$ is empty if we have configuration~\ref{short_left_neighbour} or configuration~\ref{no_left_neighbour} for $b^{(1)}$.

Combining the above arguments, we see that we are done if $K = 2$. Further we argue by induction on $K$. The case $K = 2$ is a basis of induction.

Consider the sequence
\begin{equation}
\label{ind_hyp_sequence}
(U^{(2)} \mapsto \ldots \mapsto U^{(K)}, ((a_h^{(2)}, a_j^{(2)}), \ldots, (a_h^{(K - 1)}, a_j^{(K - 1)})), (b^{(2)}, \ldots, b^{(K)})).
\end{equation}
Clearly, this is a $(b^{(2)}, U^{(2)})$-admissible sequence. It has $K -1$ replacements. We want to apply the induction hypothesis to this sequence. In order to do this, first we calculate $t(b^{(2)})$.

First consider the case when $a_h^{(1)}$ is not separated from $b^{(1)}$ in $U^{(1)}$. Then $a_j^{(1)}$ is not separated from $b^{(2)}$ in $U^{(2)}$.
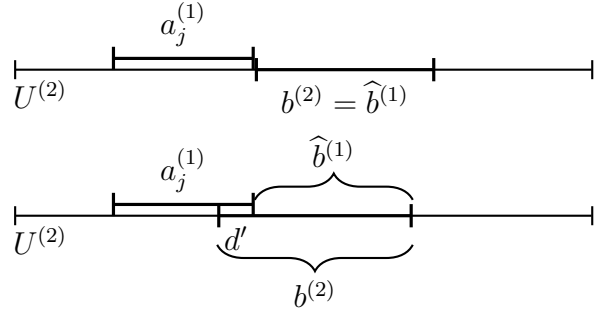
\begin{center}
\begin{tikzpicture}
\node at (0, 0) [below, xshift=10] {$U^{(2)}$};
\draw[|-|, black, thick] (0,0)--(7.7,0);
\draw[|-|, black, very thick] (1.3, 0.15)--(3.2, 0.15) node[midway, above] {$a_j^{(1)}$};
\draw[|-|, black, very thick] (3.2, 0)--(5.6, 0) node[midway, below] {$b^{(2)} = \widehat{b}^{(1)}$};
\end{tikzpicture}

\begin{tikzpicture}
\node at (0, 0) [below, xshift=10] {$U^{(2)}$};
\draw[|-|, black, thick] (0,0)--(7.7,0);
\draw[|-|, black, very thick] (1.3, 0.15)--(3.2, 0.15) node[midway, above] {$a_j^{(1)}$};
\draw[|-|, black, very thick] (2.7, 0)--(5.3, 0);
\path (2.7, 0)--(3.2, 0) node[midway, below] {$d^{\prime}$};
\draw [thick, decorate, decoration={brace, amplitude=10pt, raise=6pt}] (3.2, 0) to node[midway, above, yshift=14pt] {$\widehat{b}^{(1)}$} (5.3, 0);
\draw [thick, decorate, decoration={brace, amplitude=10pt, raise=11pt, mirror}] (2.7, 0) to node[midway, below, yshift=-20pt] {$b^{(2)}$} (5.3, 0);
\end{tikzpicture}
\end{center}
Assume $\SPM(a_j^{(1)}) \geqslant \tau - 3$. Since $a_j^{(1)}$ is not a small piece, $a_j^{(1)} \in \mo{U^{(2)}}$ (see Section~\ref{mt_configurations},~\ref{a_j_keep_structure}). Therefore, $a_j^{(1)} \in \longmo{U^{(2)}}$. Since $a_j^{(1)}$ is not separated from $b^{(2)}$ from the left, $a_j^{(1)}$ is a left neighbour of $b^{(2)}$. Therefore, we have configuration~\ref{long_left_neighbour} for $b^{(2)}$ in $U^{(2)}$. If $a_j^{(1)}$ and $b^{(2)}$ touch at a point, then $t(b^{(2)}) = b^{(2)} = \widehat{b}^{(1)}$.
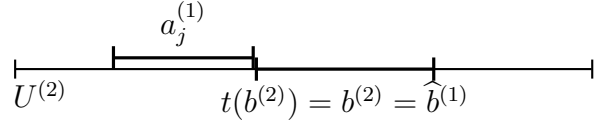
\begin{center}
\begin{tikzpicture}
\node at (0, 0) [below, xshift=10] {$U^{(2)}$};
\draw[|-|, black, thick] (0,0)--(7.7,0);
\draw[|-|, black, very thick] (1.3, 0.15)--(3.2, 0.15) node[midway, above] {$a_j^{(1)}$};
\draw[|-|, black, very thick] (3.2, 0)--(5.6, 0) node[midway, below] {$t(b^{(2)}) = b^{(2)} = \widehat{b}^{(1)}$};
\end{tikzpicture}
\end{center}
Assume $a_j^{(1)}$ and $b^{(2)}$ have an overlap $d^{\prime}$. Then $t(b^{(2)})$ starts at the end of $d^{\prime}$, because $\SPM(a_j^{(1)}) \geqslant \tau - 3$. Therefore, $t(b^{(2)}) = \widehat{b}^{(1)}$.
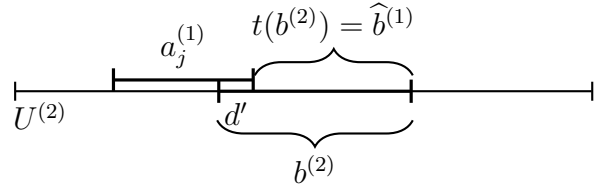
\begin{center}
\begin{tikzpicture}
\node at (0, 0) [below, xshift=10] {$U^{(2)}$};
\draw[|-|, black, thick] (0,0)--(7.7,0);
\draw[|-|, black, very thick] (1.3, 0.15)--(3.2, 0.15) node[midway, above] {$a_j^{(1)}$};
\draw[|-|, black, very thick] (2.7, 0)--(5.3, 0);
\path (2.7, 0)--(3.2, 0) node[midway, below] {$d^{\prime}$};
\draw [thick, decorate, decoration={brace, amplitude=10pt, raise=6pt}] (3.2, 0) to node[midway, above, yshift=14pt] {$t(b^{(2)}) = \widehat{b}^{(1)}$} (5.3, 0);
\draw [thick, decorate, decoration={brace, amplitude=10pt, raise=11pt, mirror}] (2.7, 0) to node[midway, below, yshift=-20pt] {$b^{(2)}$} (5.3, 0);
\end{tikzpicture}
\end{center}
We proved above that $\widehat{b}^{(1)} = t(b^{(1)})$ when $a_h$ is not separated from $b^{(1)}$. Therefore, in both cases we have $t(b^{(2)}) = \widehat{b}^{(1)} = t(b^{(1)})$. Since $\SPM(t(b)) \geqslant 3$, we have $\SPM(t(b^{(2)})) \geqslant 3$ as well. Hence, we can apply the induction hypothesis to the sequence~\eqref{ind_hyp_sequence}. Then we obtain $b^{(k)} = p^{(k)}t(b^{(2)}) = p^{(k)}t(b^{(1)}) = p^{(k)}t(b)$, where $p^{(k)}$ is a small piece (possibly empty), $k = 2, \ldots, K$. Since we have configuration~\ref{long_left_neighbour} for $b^{(1)}$ in $U^{(1)}$, this means that we are done for the case when $a_h^{(1)}$ is not separated from $b^{(1)}$ and $\SPM(a_j^{(1)}) \geqslant \tau - 3$.

Assume $a_h^{(1)}$ is not separated from $b^{(1)}$ in $U^{(1)}$ and $\SPM(a_j^{(1)}) < \tau - 3$. It follows from Corollary~\ref{short_neighbours_left} that every left neighbour of $b^{(2)}$ is of $\SPM$-measure $\leqslant \max(4, \SPM(a_j^{(1)}))$. Since $\tau \geqslant 10$, we have $\tau - 4 \geqslant 6 > 4$. Therefore, every left neighbour of $b^{(2)}$ is of $\SPM$-measure $\leqslant \tau - 4$. Hence, we obtain configuration~\ref{short_left_neighbour} for $b^{(2)}$ in $U^{(2)}$. Then, by definition, $t(b^{(2)}) = b^{(2)}$. We proved above that $b^{(2)} =  p^{(2)}t(b^{(1)})$, where $p^{(2)}$ is a small piece (possibly empty). Hence, $t(b^{(2)}) = p^{(2)}t(b^{(1)})$. So, $\SPM(t(b^{(2)})) \geqslant \SPM(t(b^{(1)})) \geqslant 3$. Therefore, we can apply the induction hypothesis to the sequence~\eqref{ind_hyp_sequence}. Then we obtain $b^{(k)} = t(b^{(2)}) = p^{(2)}t(b^{(1)}) = p^{(2)}t(b)$, where $p^{(2)}$ is a small piece (possibly empty), $k = 2, \ldots, K$. Since we have configuration~\ref{long_left_neighbour} for $b^{(1)}$ in $U^{(1)}$, this means that we are done for the case when $a_h^{(1)}$ is not separated from $b^{(1)}$ and $\SPM(a_j^{(1)}) < \tau - 3$.

Let us return to the case when $a_h^{(1)}$ is separated from $b^{(1)}$.
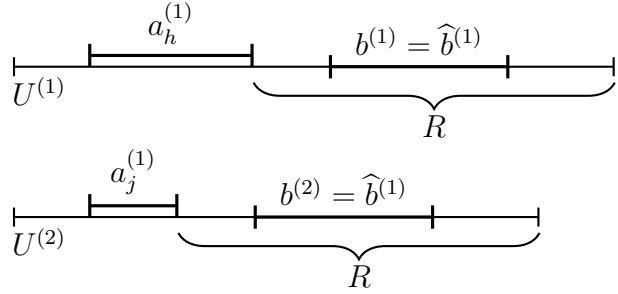
\begin{center}
\begin{tikzpicture}
\node at (0, 0) [below, xshift=10] {$U^{(1)}$};
\draw[|-|, black, thick] (0,0)--(8,0);
\draw[|-|, black, very thick] (1, 0.15)--(3.2, 0.15) node[midway, above] {$a_h^{(1)}$};
\draw[|-|, black, very thick] (4.2, 0)--(6.6, 0) node[midway, above] {$b^{(1)} = \widehat{b}^{(1)}$};
\draw [thick, decorate, decoration={brace, amplitude=10pt, raise=6pt, mirror}] (3.2, 0) to node[midway, below, yshift=-14pt] {$R$} (8, 0);

\node at (0, -2) [below, xshift=10] {$U^{(2)}$};
\draw[|-|, black, thick] (0, -2)--(7, -2);
\draw[|-|, black, very thick] (1, 0.15 - 2)--(2.2, 0.15 - 2) node[midway, above] {$a_j^{(1)}$};
\draw[|-|, black, very thick] (3.2, -2)--(5.6, -2) node[midway, above] {$b^{(2)} = \widehat{b}^{(1)}$};
\draw [thick, decorate, decoration={brace, amplitude=10pt, raise=6pt, mirror}] (2.2, -2) to node[midway, below, yshift=-14pt] {$R$} (7, -2);
\end{tikzpicture}
\end{center}
First consider configurations~\ref{short_left_neighbour} and~\ref{no_left_neighbour} for $b^{(2)}$ in $U^{(2)}$. That is, either all the left neighbours of $b^{(2)}$ are of $\SPM$-measure $< \tau - 3$, or $b^{(2)}$ does not have left neighbours at all. Then, by definition, $t(b^{(2)}) = b^{(2)}$. Since $a_h^{(1)}$ is separated from $b^{(1)}$, we obtain $b^{(2)} = \widehat{b}^{(1)} = b^{(1)}$. Hence, $t(b^{(2)}) = b^{(1)}$. So, $\SPM(t(b^{(2)})) = \SPM(b^{(1)}) \geqslant \SPM(t(b^{(1)})) \geqslant 3$. Therefore, we can apply the induction hypothesis to the sequence~\eqref{ind_hyp_sequence}. Since we have configuration~\ref{short_left_neighbour} or~\ref{no_left_neighbour} for $b^{(2)}$, we obtain $b^{(k)} = t(b^{(2)}) = b^{(1)}$, $k = 2, \ldots, K$. We proved above that $b^{(1)} = p^{(1)}t(b^{(1)})$, where $p^{(1)}$ is a small piece, and $p^{(1)}$ is empty if we have configuration~\ref{short_left_neighbour} or configuration~\ref{no_left_neighbour} for $b^{(1)}$. Therefore, $b^{(k)} = p^{(1)}t(b^{(1)})$, $k = 2, \ldots, K$. So, we are done for the case when $a_h^{(1)}$ is separated from $b^{(1)}$ and there are configurations~\ref{short_left_neighbour} or~\ref{no_left_neighbour} for $b^{(2)}$ in $U^{(2)}$.

It remains to consider the case when $a_h^{(1)}$ is separated from $b^{(1)}$ and we have configuration~\ref{long_left_neighbour} for $b^{(2)}$ in $U^{(2)}$. That is, there exists a left neighbour of $b^{(2)}$ of $\SPM$-measure $\geqslant \tau -3$. Recall that, by Corollary~\ref{single_long_neighbour}, $b^{(2)}$ can have only one left neighbour of $\SPM$-measure $\geqslant 3$. That is, in the considered case $b^{(2)}$ has a single left neighbour of $\SPM$-measure $\geqslant \tau - 3$. Let $c^{\prime}$ be a left neighbour of $b^{(2)}$ of $\SPM$-measure $\geqslant \tau - 3$. Since $a_h^{(1)} \mapsto a_j^{(1)}$ is an admissible replacement in $U^{(1)}$ and $\SPM(c^{\prime}) \geqslant \tau - 3 \geqslant 6$, by Lemma~\ref{neighbours_a_h_separated1_left}, we obtain that $a_j^{(1)}$ is not contained in $c^{\prime}$. Recall that $U^{(2)} = La_j^{(1)}R$. Let $\widehat{c}$ be the intersection of $R$ and $c^{\prime}$. First assume that $a_j^{(1)}$ is separated from $c^{\prime}$ in $U^{(2)}$. Then, clearly, $c^{\prime} = \widehat{c}$ in $U^{(2)}$, and $\widehat{c}$ is a maximal occurrence in $U^{(1)}$.
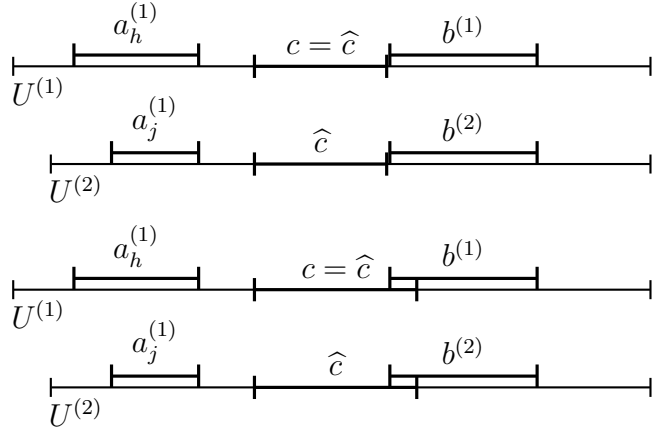
\begin{center}
\begin{tikzpicture}
\node at (0, 0) [below, xshift=10] {$U^{(1)}$};
\draw[|-|, black, thick] (0,0)--(8.5,0);
\draw[|-|, black, very thick] (0.8, 0.15)--(2.5, 0.15) node[midway, above] {$a_h^{(1)}$};
\draw[|-|, black, very thick] (3.2, 0)--(5, 0) node[midway, above] {$c = \widehat{c}$};
\draw[|-|, black, very thick] (5, 0.15)--(7, 0.15) node[midway, above] {$b^{(1)}$};

\node at (0 + 0.5, -1.3) [below, xshift=10] {$U^{(2)}$};
\draw[|-|, black, thick] (0.5,-1.3)--(8.5,-1.3);
\draw[|-|, black, very thick] (0.8 + 0.5, 0.15 - 1.3)--(2.5, 0.15 - 1.3) node[midway, above] {$a_j^{(1)}$};
\draw[|-|, black, very thick] (3.2, -1.3)--(5, -1.3) node[midway, above] {$\widehat{c}$};
\draw[|-|, black, very thick] (5, 0.15 - 1.3)--(7, 0.15 - 1.3) node[midway, above] {$b^{(2)}$};
\end{tikzpicture}

\begin{tikzpicture}
\node at (0, 0) [below, xshift=10] {$U^{(1)}$};
\draw[|-|, black, thick] (0,0)--(8.5,0);
\draw[|-|, black, very thick] (0.8, 0.15)--(2.5, 0.15) node[midway, above] {$a_h^{(1)}$};
\draw[|-|, black, very thick] (3.2, 0)--(5.4, 0) node[midway, above] {$c = \widehat{c}$};
\draw[|-|, black, very thick] (5, 0.15)--(7, 0.15) node[midway, above] {$b^{(1)}$};

\node at (0 + 0.5, -1.3) [below, xshift=10] {$U^{(2)}$};
\draw[|-|, black, thick] (0 + 0.5,-1.3)--(8.5,-1.3);
\draw[|-|, black, very thick] (0.8 + 0.5, 0.15 - 1.3)--(2.5, 0.15 - 1.3) node[midway, above] {$a_j^{(1)}$};
\draw[|-|, black, very thick] (3.2, -1.3)--(5.4, -1.3) node[midway, above] {$\widehat{c}$};
\draw[|-|, black, very thick] (5, 0.15 - 1.3)--(7, 0.15 - 1.3) node[midway, above] {$b^{(2)}$};
\end{tikzpicture}
\end{center}
It is clear that $\widehat{c}$ is a left neighbour of $b^{(1)}$ in $U^{(1)}$. Therefore, since $\SPM(\widehat{c}) = \SPM(c^{\prime}) \geqslant 3$, we have configuration~\ref{long_left_neighbour} for $b^{(1)}$ in $U^{(1)}$. Since $c^{\prime} = \widehat{c}$, we see that $t(b^{(2)}) = t(b^{(1)})$. Hence, we can apply the induction hypothesis to the sequence~\eqref{ind_hyp_sequence}. Then we obtain that $b^{(k)} = p^{(k)}t(b^{(1)})$, $k = 2, \ldots, K$. So, we are done for the case when $a_j^{(1)}$ is separated from $c^{\prime}$ in $U^{(2)}$.

Now assume that $a_j^{(1)}$ is not separated from $c^{\prime}$ in $U^{(2)}$. Recall that from Lemma~\ref{neighbours_a_h_separated1_left} it follows that either $a_j^{(1)}$ and $c^{\prime}$ touch at a point, or $a_j^{(1)}$ and $c^{\prime}$ have in overlap in this case.
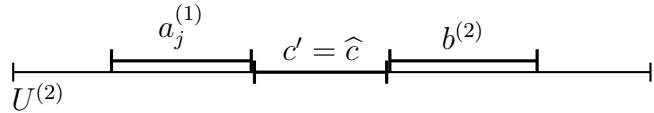
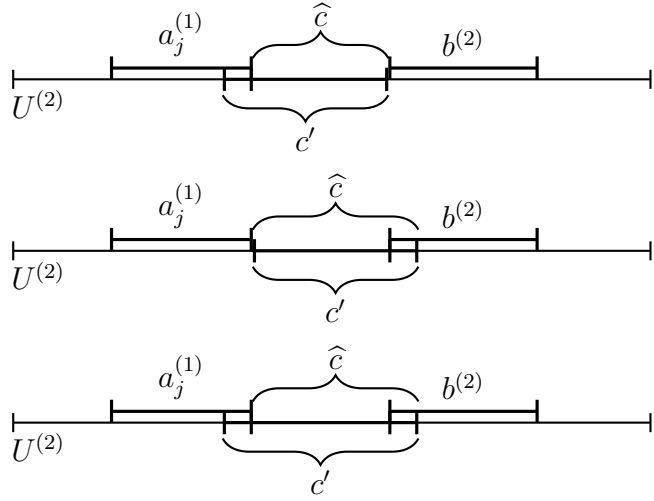
\begin{center}
\begin{tikzpicture}
\node at (0, 0) [below, xshift=10] {$U^{(2)}$};
\draw[|-|, black, thick] (0,0)--(8.5,0);
\draw[|-|, black, very thick] (1.3, 0.15)--(3.2, 0.15) node[midway, above] {$a_j^{(1)}$};
\draw[|-|, black, very thick] (3.2, 0)--(5, 0) node[midway, above] {$c^{\prime} = \widehat{c}$};
\draw[|-|, black, very thick] (5, 0.15)--(7, 0.15) node[midway, above] {$b^{(2)}$};
\end{tikzpicture}

\begin{tikzpicture}
\node at (0, 0) [below, xshift=10] {$U^{(2)}$};
\draw[|-|, black, thick] (0,0)--(8.5,0);
\draw[|-|, black, very thick] (1.3, 0.15)--(3.2, 0.15) node[midway, above] {$a_j^{(1)}$};
\draw[|-|, black, very thick] (2.8, 0)--(3.2, 0);
\draw[-|, black, very thick] (3.2, 0)--(5, 0);
\draw [thick, decorate, decoration={brace, amplitude=10pt, raise=8pt}] (3.2, 0) to node[midway, above, yshift=16pt] {$\widehat{c}$} (5, 0);
\draw[|-|, black, very thick] (5, 0.15)--(7, 0.15) node[midway, above] {$b^{(2)}$};
\draw [thick, decorate, decoration={brace, amplitude=10pt, raise=6pt, mirror}] (2.8, 0) to node[midway, below, yshift=-14pt] {$c^{\prime}$} (5, 0);
\end{tikzpicture}

\begin{tikzpicture}
\node at (0, 0) [below, xshift=10] {$U^{(2)}$};
\draw[|-|, black, thick] (0,0)--(8.5,0);
\draw[|-|, black, very thick] (1.3, 0.15)--(3.2, 0.15) node[midway, above] {$a_j^{(1)}$};
\draw[|-, black, very thick] (3.2, 0)--(5, 0);
\draw [thick, decorate, decoration={brace, amplitude=10pt, raise=8pt}] (3.2, 0) to node[midway, above, yshift=16pt] {$\widehat{c}$} (5.4, 0);
\draw[|-|, black, very thick] (5, 0)--(5.4, 0);
\draw[|-|, black, very thick] (5, 0.15)--(7, 0.15) node[midway, above] {$b^{(2)}$};
\draw [thick, decorate, decoration={brace, amplitude=10pt, raise=6pt, mirror}] (3.2, 0) to node[midway, below, yshift=-14pt] {$c^{\prime}$} (5.4, 0);
\end{tikzpicture}

\begin{tikzpicture}
\node at (0, 0) [below, xshift=10] {$U^{(2)}$};
\draw[|-|, black, thick] (0,0)--(8.5,0);
\draw[|-|, black, very thick] (1.3, 0.15)--(3.2, 0.15) node[midway, above] {$a_j^{(1)}$};
\draw[-, black, very thick] (3.2, 0)--(5, 0);
\draw[|-|, black, very thick] (2.8, 0)--(3.2, 0);
\draw[|-|, black, very thick] (5, 0)--(5.4, 0);
\draw [thick, decorate, decoration={brace, amplitude=10pt, raise=8pt}] (3.2, 0) to node[midway, above, yshift=16pt] {$\widehat{c}$} (5.4, 0);
\draw[|-|, black, very thick] (5, 0.15)--(7, 0.15) node[midway, above] {$b^{(2)}$};
\draw [thick, decorate, decoration={brace, amplitude=10pt, raise=6pt, mirror}] (2.8, 0) to node[midway, below, yshift=-14pt] {$c^{\prime}$} (5.4, 0);
\end{tikzpicture}
\end{center}
Since the overlap of $a_j^{(1)}$ and $c^{\prime}$ is a small piece (possibly empty), we see that
\begin{equation*}
\SPM(\widehat{c}) \geqslant \SPM(c^{\prime}) - 1 \geqslant \tau - 3 - 1 = \tau - 4.
\end{equation*}
Hence, it follows from Lemma~\ref{non_covered_inverse_image} that there exists $c \in \mo{U^{(1)}}$ such that $c^{\prime}$ is an image of $c$ in $U^{(2)}$. Moreover, $c = X\widehat{c}$, where $X$ is an overlap of $c$ and $a_h^{(1)}$. So, $c$ is a left neighbour of $b^{(1)}$.
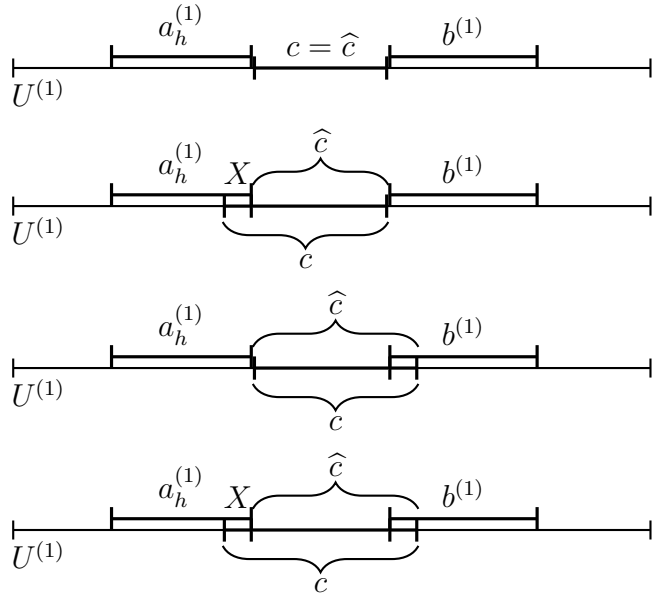
\begin{center}
\begin{tikzpicture}
\node at (0, 0) [below, xshift=10] {$U^{(1)}$};
\draw[|-|, black, thick] (0,0)--(8.5,0);
\draw[|-|, black, very thick] (1.3, 0.15)--(3.2, 0.15) node[midway, above] {$a_h^{(1)}$};
\draw[|-|, black, very thick] (3.2, 0)--(5, 0) node[midway, above] {$c = \widehat{c}$};
\draw[|-|, black, very thick] (5, 0.15)--(7, 0.15) node[midway, above] {$b^{(1)}$};
\end{tikzpicture}

\begin{tikzpicture}
\node at (0, 0) [below, xshift=10] {$U^{(1)}$};
\draw[|-|, black, thick] (0,0)--(8.5,0);
\draw[|-|, black, very thick] (1.3, 0.15)--(3.2, 0.15) node[midway, above] {$a_h^{(1)}$};
\draw[|-|, black, very thick] (2.8, 0)--(3.2, 0) node[midway, above, yshift=4pt] {$X$};
\draw[-|, black, very thick] (3.2, 0)--(5, 0);
\draw [thick, decorate, decoration={brace, amplitude=10pt, raise=8pt}] (3.2, 0) to node[midway, above, yshift=16pt] {$\widehat{c}$} (5, 0);
\draw[|-|, black, very thick] (5, 0.15)--(7, 0.15) node[midway, above] {$b^{(1)}$};
\draw [thick, decorate, decoration={brace, amplitude=10pt, raise=6pt, mirror}] (2.8, 0) to node[midway, below, yshift=-14pt] {$c$} (5, 0);
\end{tikzpicture}

\begin{tikzpicture}
\node at (0, 0) [below, xshift=10] {$U^{(1)}$};
\draw[|-|, black, thick] (0,0)--(8.5,0);
\draw[|-|, black, very thick] (1.3, 0.15)--(3.2, 0.15) node[midway, above] {$a_h^{(1)}$};
\draw[|-, black, very thick] (3.2, 0)--(5, 0);
\draw[|-|, black, very thick] (5, 0)--(5.4, 0);
\draw [thick, decorate, decoration={brace, amplitude=10pt, raise=8pt}] (3.2, 0) to node[midway, above, yshift=16pt] {$\widehat{c}$} (5.4, 0);
\draw[|-|, black, very thick] (5, 0.15)--(7, 0.15) node[midway, above] {$b^{(1)}$};
\draw [thick, decorate, decoration={brace, amplitude=10pt, raise=6pt, mirror}] (3.2, 0) to node[midway, below, yshift=-14pt] {$c$} (5.4, 0);
\end{tikzpicture}

\begin{tikzpicture}
\node at (0, 0) [below, xshift=10] {$U^{(1)}$};
\draw[|-|, black, thick] (0,0)--(8.5,0);
\draw[|-|, black, very thick] (1.3, 0.15)--(3.2, 0.15) node[midway, above] {$a_h^{(1)}$};
\draw[-, black, very thick] (3.2, 0)--(5, 0);
\draw [thick, decorate, decoration={brace, amplitude=10pt, raise=8pt}] (3.2, 0) to node[midway, above, yshift=16pt] {$\widehat{c}$} (5.4, 0);
\draw[|-|, black, very thick] (2.8, 0)--(3.2, 0) node[midway, above, yshift=4pt] {$X$};
\draw[|-|, black, very thick] (5, 0)--(5.4, 0);
\draw[|-|, black, very thick] (5, 0.15)--(7, 0.15) node[midway, above] {$b^{(1)}$};
\draw [thick, decorate, decoration={brace, amplitude=10pt, raise=6pt, mirror}] (2.8, 0) to node[midway, below, yshift=-14pt] {$c$} (5.4, 0);
\end{tikzpicture}
\end{center}
Clearly, $\SPM(c) \geqslant \SPM(\widehat{c})$. Therefore, since $\SPM(\widehat{c}) \geqslant \tau - 4$, we have $\SPM(c) \geqslant \tau - 4$ as well. We consider possibilities $\SPM(c) \geqslant \tau - 3$ and $\SPM(c) = \tau - 4$ separately.

Assume $\SPM(c) \geqslant \tau - 3$. Then we have configuration~\ref{long_left_neighbour} for $b^{(1)}$ in $U^{(1)}$. Since $c$ is a left neighbour of $b^{(1)}$ and $\SPM(c) \geqslant \tau - 3$, by definition, $t(b^{(1)})$ starts at the end of $\widehat{c}$ in $U^{(1)}$. Recall that we also have configuration~\ref{long_left_neighbour} for $b^{(2)}$ in $U^{(2)}$. Since $c^{\prime}$ is a left neighbour of $b^{(2)}$ and $\SPM(c^{\prime}) \geqslant \tau - 3$, by definition, $t(b^{(2)})$ starts at the end of $\widehat{c}$ in $U^{(2)}$. Therefore, since $b^{(1)} = b^{(2)}$, we obtain $t(b^{(1)}) = t(b^{(2)})$. Hence, $\SPM(t(b^{(2)})) = \SPM(t(b^{(1)})) \geqslant 3$. Therefore, we can apply the induction hypothesis to the sequence~\eqref{ind_hyp_sequence}. We obtain $b^{(k)} = p^{(k)}t(b^{(2)}) = p^{(k)}t(b^{(1)})$, where $p^{(k)}$ is a small piece (possibly empty), $k = 2, \ldots, K$. Since we have configuration~\ref{long_left_neighbour} for $b^{(1)}$ in $U^{(1)}$, this means that we are done for the case $\SPM(c) \geqslant \tau - 3$.

Assume $\SPM(c) = \tau - 4$. Clearly, $\vert \SPM(c) - \SPM(c^{\prime})\vert \leqslant 1$. By our assumption, $\SPM(c^{\prime}) \geqslant \tau - 3$. Hence, we obtain $\SPM(c^{\prime}) = \tau - 3$. So, we have configuration~\ref{short_left_neighbour} for $b^{(1)}$ and configuration~\ref{long_left_neighbour} for $b^{(2)}$. If we apply the induction hypothesis to the sequence~\eqref{ind_hyp_sequence} in this case, we do not obtain the desired result. Namely, since we have configuration~\ref{short_left_neighbour} for $b^{(1)}$, informally speaking, we have to show that there are no possibilities to change for $b^{(1)}$ in the initial sequence. However, since we have  configuration~\ref{long_left_neighbour} for $b^{(2)}$, by the induction hypothesis, we obtain that $b^{(2)}$ may change in the sequence~\eqref{ind_hyp_sequence}, and hence $b^{(1)}$ may change in the initial sequence. So, in what follows we do not apply the induction hypothesis directly to the sequence~\eqref{ind_hyp_sequence} and argue in a different way. We will directly show that $a_h^{(k)}$ is separated from $b^{(k)}$ for $k = 1, \ldots, K - 1$.

Let us prove an auxiliary lemma.
\begin{lemma}
\label{c_admissible_sequence}
We are under the conditions of Lemma~\ref{admissible_replacements_left} and are using notations introduced above. Assume $\SPM(c) = \tau - 4$ and $\SPM(c^{\prime}) = \tau - 3$. Consider the replacement $a_h^{(k - 1)} \mapsto a_j^{(k - 1)}$ in $U^{(k - 1)}$, $U^{(k)}$ is a resulting monomial, $k = 3, \ldots, K$. Let $c^{(2)} = c^{\prime}$, and let $c^{(k)}$ be an image of $c^{(k - 1)}$ in $U^{(k)}$, $k = 3, \ldots, K$. Then
\begin{itemize}
\item
$c^{(k)}$ is a left neighbour of $b^{(k)}$, $k = 2, \ldots, K$;
\item
$\tau - 4 \leqslant \SPM(c^{(k)}) \leqslant \tau - 3$, $k = 2, \ldots, K$, $c^{(k)} = q^{(k)}\widehat{c}$, where $q^{(k)}$ is a small piece (possibly empty);
\item
$a_h^{(k)}$ starts from the left of the beginning of $c^{(k)}$, $k = 2, \ldots, K - 1$. In particular, $a_h^{(k)} \neq c^{(k)}$.
\end{itemize}
\end{lemma}
\begin{proof}
We have the sequence of element $c^{(2)}, \ldots, c^{(K)}$. We prove Lemma~\ref{c_admissible_sequence} by induction on $k$ - the number of the element $c^{(k)}$ in this sequence.

Let us prove the basis of induction. That is, we consider $k = 2$. By our initial assumption, $\SPM(c^{(2)}) = \SPM(c^{\prime}) = \tau - 3$ and $c^{(2)} = c^{\prime}$ is a left neighbour of $b^{(2)}$. Since $a_h^{(2)} \mapsto a_j^{(2)}$ is an admissible replacement in $U^{(2)}$, we have $\SPM(a_h^{(2)}) \geqslant \tau - 2$. Therefore, $a_h^{(2)} \neq c^{(2)}$. Recall that, by the conditions of Lemma~\ref{admissible_replacements_left}, $a_h^{(2)}$ starts from the left of the beginning of $b^{(2)}$. Since $a_h^{(2)} \in \nfc{U^{(2)}}$ and $b^{(2)}$ and $c^{(2)}$ are not separated in $U^{(2)}$, $a_h^{(2)}$ can not start from the left of the beginning of $b^{(2)}$ and at the same time from the right of the beginning of $c^{(2)}$. Therefore, $a_h^{(2)}$ starts from the left of the beginning of $c^{(2)}$. So far, we are done with the basis of induction.

Let us prove the step of induction. In order to do this, let us first calculate $t(c^{(2)})$. Since $\SPM(c^{(2)}) = \tau - 3 \geqslant 6$, it follows from Lemma~\ref{neighbours_a_h_separated1_left} that $a_j^{(1)}$ is not contained in $c^{(2)}$. Therefore, since $a_j^{(1)}$ is not separated from $c^{(2)}$, $a_j^{(1)}$ and $c^{(2)}$ either touch at a point, or have an overlap. By our assumption, $\SPM(c) = \tau - 4$ and $\SPM(c^{(2)}) = \tau - 3$, that is, $\SPM(c^{(2)}) > \SPM(c)$. We proved above that $c = X\widehat{c}$, hence, $\SPM(c) \geqslant \SPM(\widehat{c})$. Therefore, $\SPM(c^{(2)}) > \SPM(c) \geqslant \SPM(\widehat{c})$. This implies that $a_j^{(1)}$ and $c^{(2)}$ have a non-empty overlap $X^{\prime}$ and $c^{(2)} = X^{\prime}\widehat{c}$.
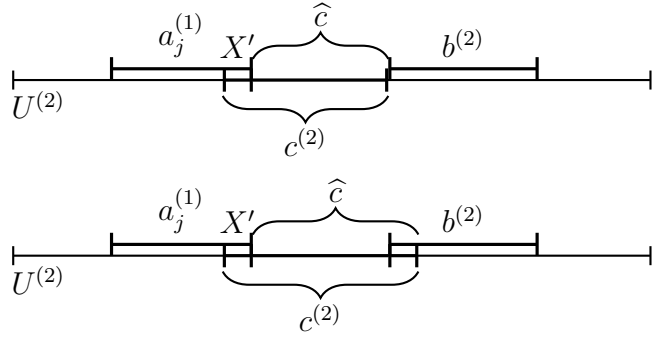
\begin{center}
\begin{tikzpicture}
\node at (0, 0) [below, xshift=10] {$U^{(2)}$};
\draw[|-|, black, thick] (0,0)--(8.5,0);
\draw[|-|, black, very thick] (1.3, 0.15)--(3.2, 0.15) node[midway, above] {$a_j^{(1)}$};
\draw[|-|, black, very thick] (2.8, 0)--(3.2, 0) node[midway, above, yshift=4] {$X^{\prime}$};
\draw[-|, black, very thick] (3.2, 0)--(5, 0);
\draw [thick, decorate, decoration={brace, amplitude=10pt, raise=8pt}] (3.2, 0) to node[midway, above, yshift=16pt] {$\widehat{c}$} (5, 0);
\draw[|-|, black, very thick] (5, 0.15)--(7, 0.15) node[midway, above] {$b^{(2)}$};
\draw [thick, decorate, decoration={brace, amplitude=10pt, raise=6pt, mirror}] (2.8, 0) to node[midway, below, yshift=-14pt] {$c^{(2)}$} (5, 0);
\end{tikzpicture}

\begin{tikzpicture}
\node at (0, 0) [below, xshift=10] {$U^{(2)}$};
\draw[|-|, black, thick] (0,0)--(8.5,0);
\draw[|-|, black, very thick] (1.3, 0.15)--(3.2, 0.15) node[midway, above] {$a_j^{(1)}$};
\draw[-, black, very thick] (3.2, 0)--(5, 0);
\draw[|-|, black, very thick] (2.8, 0)--(3.2, 0) node[midway, above, yshift=4] {$X^{\prime}$};
\draw[|-|, black, very thick] (5, 0)--(5.4, 0);
\draw [thick, decorate, decoration={brace, amplitude=10pt, raise=8pt}] (3.2, 0) to node[midway, above, yshift=16pt] {$\widehat{c}$} (5.4, 0);
\draw[|-|, black, very thick] (5, 0.15)--(7, 0.15) node[midway, above] {$b^{(2)}$};
\draw [thick, decorate, decoration={brace, amplitude=10pt, raise=6pt, mirror}] (2.8, 0) to node[midway, below, yshift=-14pt] {$c^{(2)}$} (5.4, 0);
\end{tikzpicture}
\end{center}
Assume $\SPM(a_j^{(1)}) \geqslant \tau - 3$. Recall that if $a_j^{(1)}$ is not a small piece, then $a_j^{(1)} \in \mo{U^{(2)}}$ (see Section~\ref{mt_configurations},~\ref{a_j_keep_structure}). Therefore, $a_j^{(1)}$ is a left neighbour of $c^{(2)}$. So, since $\SPM(a_j^{(1)}) \geqslant \tau - 3$, we have configuration~\ref{long_left_neighbour} for $c^{(2)}$. Then, by definition, $t(c^{(2)}) = \widehat{c}$.
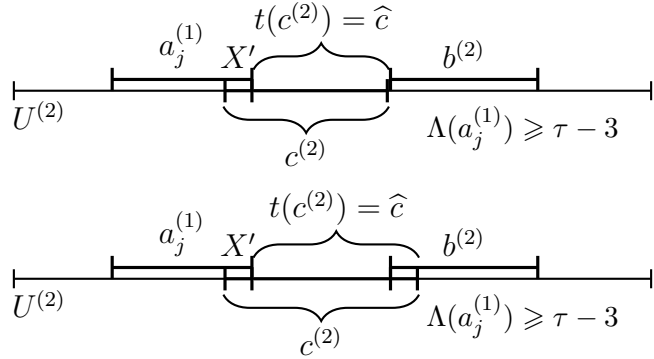
\begin{center}
\begin{tikzpicture}
\node at (0, 0) [below, xshift=10] {$U^{(2)}$};
\draw[|-|, black, thick] (0,0)--(8.5,0);
\draw[|-|, black, very thick] (1.3, 0.15)--(3.2, 0.15) node[midway, above] {$a_j^{(1)}$};
\draw[|-|, black, very thick] (2.8, 0)--(3.2, 0) node[midway, above, yshift=4] {$X^{\prime}$};
\draw[-|, black, very thick] (3.2, 0)--(5, 0);
\draw [thick, decorate, decoration={brace, amplitude=10pt, raise=8pt}] (3.2, 0) to node[midway, above, yshift=16pt] {$t(c^{(2)}) = \widehat{c}$} (5, 0);
\draw[|-|, black, very thick] (5, 0.15)--(7, 0.15) node[midway, above] {$b^{(2)}$};
\draw [thick, decorate, decoration={brace, amplitude=10pt, raise=6pt, mirror}] (2.8, 0) to node[midway, below, yshift=-14pt] {$c^{(2)}$} (5, 0);
\node[text width=4cm, align=left] at (7.5, -0.5) {\small\baselineskip=10pt \textit{$\SPM(a_j^{(1)}) \geqslant \tau - 3$}};
\end{tikzpicture}

\begin{tikzpicture}
\node at (0, 0) [below, xshift=10] {$U^{(2)}$};
\draw[|-|, black, thick] (0,0)--(8.5,0);
\draw[|-|, black, very thick] (1.3, 0.15)--(3.2, 0.15) node[midway, above] {$a_j^{(1)}$};
\draw[-, black, very thick] (3.2, 0)--(5, 0);
\draw[|-|, black, very thick] (2.8, 0)--(3.2, 0) node[midway, above, yshift=4] {$X^{\prime}$};
\draw[|-|, black, very thick] (5, 0)--(5.4, 0);
\draw [thick, decorate, decoration={brace, amplitude=10pt, raise=8pt}] (3.2, 0) to node[midway, above, yshift=16pt] {$t(c^{(2)}) = \widehat{c}$} (5.4, 0);
\draw[|-|, black, very thick] (5, 0.15)--(7, 0.15) node[midway, above] {$b^{(2)}$};
\draw [thick, decorate, decoration={brace, amplitude=10pt, raise=6pt, mirror}] (2.8, 0) to node[midway, below, yshift=-14pt] {$c^{(2)}$} (5.4, 0);
\node[text width=4cm, align=left] at (7.5, -0.5) {\small\baselineskip=10pt \textit{$\SPM(a_j^{(1)}) \geqslant \tau - 3$}};
\end{tikzpicture}
\end{center}
Assume $\SPM(a_j^{(1)}) < \tau - 3$. Then it follows from Corollary~\ref{short_neighbours_left} that $\SPM$-measure of all the left neighbours of $c^{(2)}$ is not greater than $\max(4, \SPM(a_j)) \leqslant \tau - 4$. That is, we have configuration~\ref{short_left_neighbour} for $c^{(2)}$. Then, by definition, $t(c^{(2)}) = c^{(2)}$.

Consider $c^{(N)}$, where $2 < N \leqslant K$. Consider the sequence
\begin{equation}
\label{c_ind_hyp_sequence}
(U^{(2)} \mapsto \ldots \mapsto U^{(N)}, ((a_h^{(2)}, a_j^{(2)}), \ldots, (a_h^{(N - 1)}, a_j^{(N - 1)})), (c^{(2)}, \ldots, c^{(N)})).
\end{equation}
By the induction hypothesis, $a_h^{(k)} \neq c^{(k)}$, $k = 2, \ldots, N - 1$. By the conditions of Lemma~\ref{admissible_replacements_left}, $a_h^{(k)} \mapsto a_j^{(k)}$ is an admissible replacement in $U^{(k)}$. Therefore, by definition, the sequence~\eqref{c_ind_hyp_sequence} is a $(c^{(2)}, U^{(2)})$-admissible sequence. By the induction hypothesis, $a_h^{(k)}$ starts from the left of the beginning of $c^{(k)}$ for $k = 2, \ldots, N - 1$. Therefore, since there are $N - 2 < K - 1$ replacements in the sequence~\eqref{c_ind_hyp_sequence}, we can apply the induction hypothesis of Lemma~\ref{admissible_replacements_left} to this sequence. Then we obtain $c^{(k)} = q^{(k)}t(c^{(2)})$, where $q^{(k)}$ is a small piece, $q^{(k)}$ is empty if we have configuration~\ref{short_left_neighbour} for $c^{(2)}$, $k = 2, \ldots, N$. Therefore, it follows from the above that
\begin{align*}
&c^{(k)} = \begin{cases}
c^{(2)} = X^{\prime}\widehat{c} &\textit{ if } \SPM(a_j^{(1)}) < \tau - 3,\\
q^{(k)}t(c^{(2)}) = q^{(k)}\widehat{c} &\textit{ if } \SPM(a_j^{(1)}) \geqslant \tau - 3,
\end{cases}\\
&k = 2, \ldots, N.
\end{align*}
This implies that $\widehat{c} \leqslant \SPM(c^{(N)}) \leqslant \widehat{c} + 1$.

Recall that $c = X\widehat{c}$ and $c^{(2)} = X^{\prime}\widehat{c}$, where $X$ and $X^{\prime}$ are small pieces ($X$ may be empty). Therefore, $\tau - 4 = \SPM(c) \geqslant \SPM(\widehat{c})$ and $\tau - 3 = \SPM(c^{(2)}) \leqslant \SPM(\widehat{c}) + 1$. So, we obtain $\SPM(\widehat{c}) = \tau - 4$. We proved above that $\widehat{c} \leqslant \SPM(c^{(N)}) \leqslant \widehat{c} + 1$. Hence,
\begin{equation*}
\tau - 4 \leqslant \SPM(c^{(N)}) \leqslant \tau - 3.
\end{equation*}
Since $c^{(2)}$ is a left neighbour of $b^{(2)}$, we have that $\widehat{c}$ is not separated from $b^{(2)}$ in $U^{(2)}$. We proved that $c^{(k)} = q^{(k)}\widehat{c}$ for $k = 2, \ldots, N$. This obviously implies that $c^{(k)}$ is a left neighbour of $b^{(k)}$ for $k = 2, \ldots, N$. Since $a_h^{(N)} \mapsto a_j^{(N)}$ is an admissible replacement in $U^{(N)}$, $\SPM(a_h^{(N)}) \geqslant \tau - 2$. Therefore, $a_h^{(N)} \neq c^{(N)}$. As above, since $a_h^{(N)} \in \nfc{U^{(N)}}$ and $c^{(N)}$ is not separated from $b^{(N)}$, $a_h^{(N)}$ can not start from the left of the beginning of $b^{(N)}$ and at the same time from the right of the beginning of $c^{(N)}$. Hence, $a_h^{(N)}$ starts from the left of the beginning of $c^{(N)}$. So, we are done with the step of induction.
\end{proof}

Lemma~\ref{c_admissible_sequence} implies that $a_h^{(k)}$ starts from the left of the beginning of $c^{(k)}$, $b^{(k)}$ starts from the right of the beginning of $c^{(k)}$, and $\SPM(c^{(k)}) \geqslant \tau - 4 \geqslant 6 > 3$ for $k = 2,\ldots, K - 1$. Therefore, $a_h^{(k)}$ is separated from $b^{(k)}$ in $U^{(k)}$, $k = 2, \ldots, K - 1$. The occurrence $a_h^{(1)}$ is separated from $b^{(1)}$ in $U^{(1)}$, by our initial assumption. Hence, $a_h^{(k)}$ is separated from $b^{(k)}$ for $k = 1, \ldots, K - 1$. Therefore, $b^{(k)} = b^{(k - 1)}$, $k = 2, \ldots, K$. This means that $b^{(k)} = b^{(1)}$, $k = 1, \ldots, K$. We assumed that $\SPM(c) = \tau - 4$, where $c$ is a left neighbour of $b^{(1)}$. That is, we have configuration~\ref{short_left_neighbour} for $b^{(1)}$, hence, $b^{(1)} = t(b^{(1)})$. Thus, $b^{(k)} = t(b^{(1)})$, $k = 1, \ldots, K$. So, we are done for the last case. This completes the proof of the first part of Lemma~\ref{admissible_replacements_left}.

Using the above argument, one can easily see that $a_h^{(k)}$ is always contained in the prefix of $U^{(k)}$ that ends at the beginning point of $t(b)$. Therefore, the suffix of $U^{(k)}$ that starts at the beginning point of $t(b)$ does not change after the replacement $a_h^{(k)} \mapsto a_j^{(k)}$, and, therefore, does not depend on $k$. This completes the proof of  Lemma~\ref{admissible_replacements_left}.
\end{proof}

Assume $U$ is a monomial, $b \in \longmo{U}$. Then there are the following possibilities for right neighbours of $b$:
\begin{enumerate}[label=(R\arabic*)]
\item
\label{long_right_neighbour}
there exists a right neighbour of $b$ of $\SPM$-measure $\geqslant \tau - 3$;
\item
\label{short_right_neighbour}
all the right neighbours of $b$ are of $\SPM$-measure $< \tau - 3$;
\item
\label{no_right_neighbour}
$b$ does not have a right neighbours at all.
\end{enumerate}

Let $i(b)$\label{initial_subword_def} be an initial subword of $b$ defined by the following rule. Assume we have configuration~\ref{long_right_neighbour} for $b$ and $d\in \longmo{U}$ is a right neighbour of $b$ with $\SPM(d) \geqslant \tau - 3$. Then $i(b)$ is the prefix of $b$ such that its end point is equal to the beginning point of $d$.
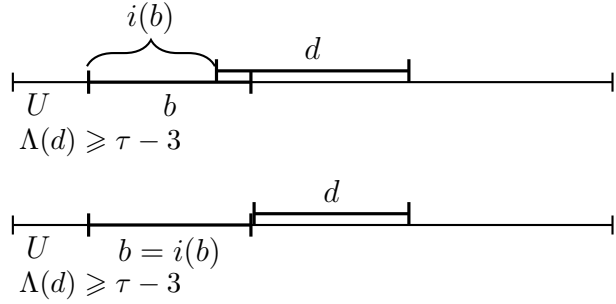
\begin{center}
\begin{tikzpicture}
\node at (0, 0) [below, xshift=10] {$U$};
\draw[|-|, black, thick] (0,0)--(8,0);
\draw[|-|, black, very thick] (1, 0)--(3.2, 0) node[midway, below] {$b$};
\draw[|-|, black, very thick] (2.7, 0.15)--(5.3, 0.15) node[midway, above] {$d$};
\draw [thick, decorate, decoration={brace, amplitude=10pt, raise=6pt}] (1, 0) to node[midway, above, yshift=14pt] {$i(b)$} (2.7, 0);
\node[text width=3cm, align=left] at (1.6, -0.8) {\small\baselineskip=10pt $\SPM(d) \geqslant \tau - 3$};
\end{tikzpicture}

\begin{tikzpicture}
\node at (0, 0) [below, xshift=10] {$U$};
\draw[|-|, black, thick] (0,0)--(8,0);
\draw[|-|, black, very thick] (1, 0)--(3.2, 0) node[midway, below] {$b = i(b)$};
\draw[|-|, black, very thick] (3.2, 0.15)--(5.3, 0.15) node[midway, above] {$d$};
\node[text width=3cm, align=left] at (1.6, -0.8) {\small\baselineskip=10pt $\SPM(d) \geqslant \tau - 3$};
\end{tikzpicture}
\end{center}
If we have configuration~\ref{short_right_neighbour} or~\ref{no_right_neighbour} for $b$, then we put $i(b)$ to be equal to $b$. So, in general, $b = i(b)p_2$, where $p_2$ is a small piece (possibly empty).

\begin{lemma}
\label{admissible_replacements_right}
Let $U$ be a monomial. Assume $b\in \longmo{U}$ and, moreover, $\SPM(i(b)) \geqslant 3$. Assume
\begin{equation*}
(U = U^{(1)} \mapsto \ldots \mapsto U^{(K)}, ((a_h^{(1)}, a_j^{(1)}), \ldots, (a_h^{(K)}, a_j^{(K)})), (b = b^{(1)}, \ldots, b^{(K)}))
\end{equation*}
is a $(b, U)$-admissible sequence such that every $a_h^{(k)}$ starts from the right of the beginning of $b^{(k)}$, $k = 1, \ldots, K - 1$ (that is, $a_h^{(k)}$ ends from the right of the end of $b^{(k)}$). Then $b^{(k)}$ is a single image of $b^{(k - 1)}$ in $U^{(k)}$, $k = 2, \ldots, K$, and
\begin{align*}
b^{(k)} = i(b)s^{(k)}, &\textit{ where } s^{(k)} \textit{ is a suffix of } b^{(k)},\\
&s^{(k)} \textit{ is a small piece (possibly empty), } k = 1, \ldots, K.
\end{align*}
Moreover, if we start with~\ref{short_right_neighbour} or~\ref{no_right_neighbour} for $b$, then $a_h^{(k - 1)}$ is separated from $b^{(k - 1)}$ in $U^{(k - 1)}$, $k = 2, \ldots, K$, and $b^{(k)} = b = i(b)$, $k = 1, \ldots, K$ (that is, $s^{(k)}$ is empty).

Under the same conditions, we have
\begin{equation*}
U^{(k)} = Ab^{(k)}B^{(k)} = Ai(b)s^{(k)}B^{(k)}, \ k = 1, \ldots, K,
\end{equation*}
where $A$ is a prefix of $U^{(k)}$ and $B^{(k)}$ is a suffix of $U^{(k)}$. That is, the prefix of $U^{(k)}$ that ends at the end point of $i(b)$ does not depend on $k$. Moreover, $a_h^{(k)}$ is a maximal occurrence in the subword $s^{(k)}B^{(k)}$, $k = 1, \ldots, K - 1$.
\end{lemma}
\begin{proof}
This lemma is proved in the same way as Lemma~\ref{admissible_replacements_left}.
\end{proof}

Let $m(b)$\label{middle_subword_def} be the intersection of $t(b)$ and $i(b)$. Namely, $m(b)$ is defined by the following rule.
\begin{itemize}
\item
\label{m_b_def_beginning}
Assume we have configuration~\ref{long_left_neighbour} for $b$ from the left side and $c \in \longmo{U}$ is a left neighbour of $b$ with $\SPM(c) \geqslant \tau - 3$. Assume we have configuration~\ref{long_right_neighbour} for $b$ from the right side and $d \in \longmo{U}$ is a right neighbour of $b$ with $\SPM(d) \geqslant \tau - 3$. Then the beginning point of $m(b)$ is equal to the end point of $c$ and the end point of $m(b)$ is equal to the beginning point of $d$.
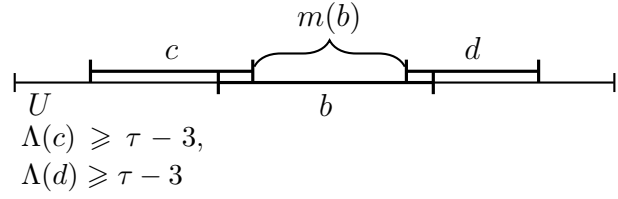
\begin{center}
\begin{tikzpicture}
\node at (0, 0) [below, xshift=10] {$U$};
\draw[|-|, black, thick] (0,0)--(8,0);
\draw[|-|, black, very thick] (1, 0.15)--(3.2, 0.15) node[midway, above] {$c$};
\draw[|-|, black, very thick] (2.7, 0)--(5.6, 0) node[midway, below] {$b$};
\draw[|-|, black, very thick] (5.2, 0.15)--(7, 0.15) node[midway, above] {$d$};
\draw [thick, decorate, decoration={brace, amplitude=10pt, raise=6pt}] (3.2, 0) to node[midway, above, yshift=14pt] {$m(b)$} (5.2, 0);
\node[text width=3cm, align=left] at (1.6, -1) {\small\baselineskip=10pt $\SPM(c) \geqslant \tau - 3$, $\SPM(d) \geqslant \tau - 3$};
\end{tikzpicture}

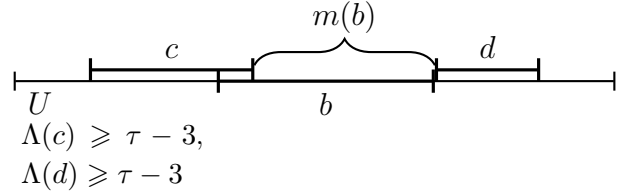
\begin{tikzpicture}
\node at (0, 0) [below, xshift=10] {$U$};
\draw[|-|, black, thick] (0,0)--(8,0);
\draw[|-|, black, very thick] (1, 0.15)--(3.2, 0.15) node[midway, above] {$c$};
\draw[|-|, black, very thick] (2.7, 0)--(5.6, 0) node[midway, below] {$b$};
\draw[|-|, black, very thick] (5.6, 0.15)--(7, 0.15) node[midway, above] {$d$};
\draw [thick, decorate, decoration={brace, amplitude=10pt, raise=6pt}] (3.2, 0) to node[midway, above, yshift=14pt] {$m(b)$} (5.6, 0);
\node[text width=3cm, align=left] at (1.6, -1) {\small\baselineskip=10pt $\SPM(c) \geqslant \tau - 3$, $\SPM(d) \geqslant \tau - 3$};
\end{tikzpicture}

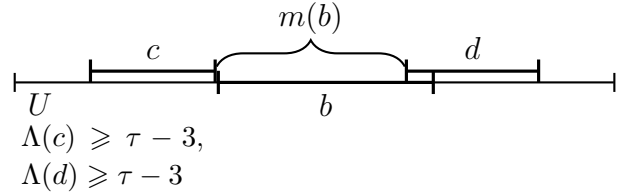
\begin{tikzpicture}
\node at (0, 0) [below, xshift=10] {$U$};
\draw[|-|, black, thick] (0,0)--(8,0);
\draw[|-|, black, very thick] (1, 0.15)--(2.7, 0.15) node[midway, above] {$c$};
\draw[|-|, black, very thick] (2.7, 0)--(5.6, 0) node[midway, below] {$b$};
\draw[|-|, black, very thick] (5.2, 0.15)--(7, 0.15) node[midway, above] {$d$};
\draw [thick, decorate, decoration={brace, amplitude=10pt, raise=6pt}] (2.7, 0) to node[midway, above, yshift=14pt] {$m(b)$} (5.2, 0);
\node[text width=3cm, align=left] at (1.6, -1) {\small\baselineskip=10pt $\SPM(c) \geqslant \tau - 3$, $\SPM(d) \geqslant \tau - 3$};
\end{tikzpicture}

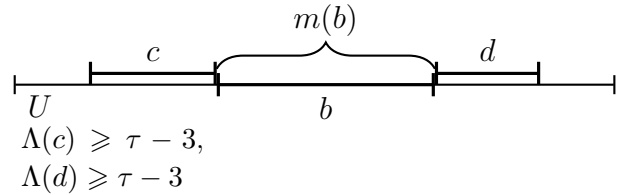
\begin{tikzpicture}
\node at (0, 0) [below, xshift=10] {$U$};
\draw[|-|, black, thick] (0,0)--(8,0);
\draw[|-|, black, very thick] (1, 0.15)--(2.7, 0.15) node[midway, above] {$c$};
\draw[|-|, black, very thick] (2.7, 0)--(5.6, 0) node[midway, below] {$b$};
\draw[|-|, black, very thick] (5.6, 0.15)--(7, 0.15) node[midway, above] {$d$};
\draw [thick, decorate, decoration={brace, amplitude=10pt, raise=6pt}] (2.7, 0) to node[midway, above, yshift=14pt] {$m(b)$} (5.6, 0);
\node[text width=3cm, align=left] at (1.6, -1) {\small\baselineskip=10pt $\SPM(c) \geqslant \tau - 3$, $\SPM(d) \geqslant \tau - 3$};
\end{tikzpicture}
\end{center}
\item
Assume we have configuration~\ref{long_left_neighbour} for $b$ from the left side and configuration~\ref{short_right_neighbour} or~\ref{no_right_neighbour} for $b$ from the right side. Then $m(b)$ is equal to $t(b)$.
\item
Assume we have configuration~\ref{short_left_neighbour} or~\ref{no_left_neighbour} for $b$ from the left side and configuration~\ref{long_right_neighbour} for $b$ from the right side. Then $m(b)$ is equal to $i(b)$.
\item
In the rest cases $m(b)$ is equal to $b$.
\end{itemize}
So, in general $b = pm(b)s$, where $p$ is a prefix of $b$, $s$ is a suffix of $b$, and $p$ and $s$ are small pieces (any of them may be empty). In particular, $\SPM(b) - 2 \leqslant \SPM(m(b)) \leqslant \SPM(b)$. Also notice that $t(b) = m(b)s$, $i(b) = pm(b)$.

\begin{lemma}
\label{admissible_replacement_in_subword}
Let $U_h$ be a monomial, $b \in \longmo{U_h}$. Let $m$ be an occurrence in $b$ such that $\SPM(m) \geqslant 3$. Assume $U_h = XmY$. Let $a_h \in \mo{U_h}$, $a_h \neq b$. Let $a_h$ and $a_j$ be incident monomials.

Assume $a_h$ is contained in $mY$. In this case the transformation $a_h \mapsto a_j$ can be considered as a replacement both in $U_h$ and in $mY$. Then $a_h \mapsto a_j$ is an admissible replacement in $U_h$ if and only if it is an admissible replacement in $mY$.

Assume $a_h$ is contained in $Xm$. In this case the transformation $a_h \mapsto a_j$ can be considered as a replacement both in $U_h$ and in $Xm$. Then $a_h \mapsto a_j$ is an admissible replacement in $U_h$ if and only if it is an admissible replacement in $Xm$.
\end{lemma}
\begin{proof}
Let us prove the first part of Lemma~\ref{admissible_replacement_in_subword}. That is, we study the case when $a_h$ is contained in $mY$. Let us denote $mY$ by $S_h$. Let us denote the beginning point of $m$ by $P$. Let $t$ be a terminal subword of $b$ that begins at the point $P$. Then $t \in \mo{S_h}$. Moreover, since $\SPM(m) \geqslant 3$, we have $t \in \longmo{S_h}$. Since $a_h \neq b$, we have the following configurations:
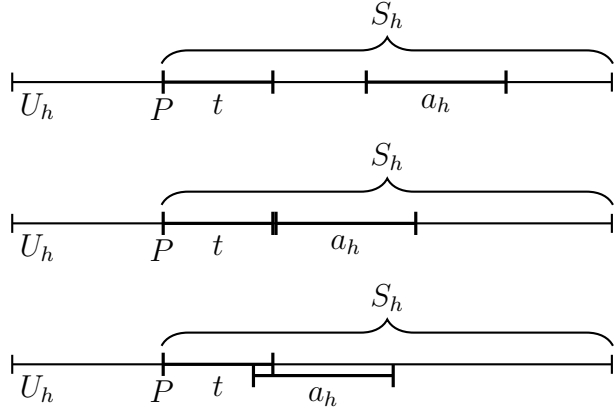
\begin{center}
\begin{tikzpicture}
\node at (0, 0) [below, xshift=10] {$U_h$};
\draw[|-|, black, thick] (0,0)--(8,0);
\draw[|-|, black, very thick] (2, 0)--(3.5, 0) node[midway, below] {$t$};
\node at (2, 0) [below, yshift=-2] {$P$};
\draw[|-|, black, very thick] (4.7, 0)--(6.6, 0) node[midway, below] {$a_h$};
\draw [thick, decorate, decoration={brace, amplitude=10pt, raise=7pt}] (2, 0) to node[midway, above, yshift=16pt] {$S_h$} (8, 0);
\end{tikzpicture}

\begin{tikzpicture}
\node at (0, 0) [below, xshift=10] {$U_h$};
\draw[|-|, black, thick] (0,0)--(8,0);
\draw[|-|, black, very thick] (2, 0)--(3.5, 0) node[midway, below] {$t$};
\node at (2, 0) [below, yshift=-2] {$P$};
\draw[|-|, black, very thick] (3.5, 0)--(5.4, 0) node[midway, below] {$a_h$};
\draw [thick, decorate, decoration={brace, amplitude=10pt, raise=7pt}] (2, 0) to node[midway, above, yshift=16pt] {$S_h$} (8, 0);
\end{tikzpicture}

\begin{tikzpicture}
\node at (0, 0) [below, xshift=10] {$U_h$};
\draw[|-|, black, thick] (0,0)--(8,0);
\draw[|-|, black, very thick] (2, 0)--(3.5, 0) node[midway, below] {$t$};
\node at (2, 0) [below, yshift=-2] {$P$};
\draw[|-|, black, very thick] (3.2, -0.15)--(5.1, -0.15) node[midway, below] {$a_h$};
\draw [thick, decorate, decoration={brace, amplitude=10pt, raise=7pt}] (2, 0) to node[midway, above, yshift=16pt] {$S_h$} (8, 0);
\end{tikzpicture}
\end{center}
So, $a_h \mapsto a_j$ can be considered as a replacement $S_h$. Let us denote the resulting monomial of the replacement $a_h \mapsto a_j$ in $S_h$ by $S_j$.

First assume $a_h \mapsto a_j$ is an admissible replacement in $U_h$. Let us show that $a_h \mapsto a_j$ is an admissible replacement in $S_h$. Assume the contrary. Namely, assume that $a_j$ is covered by images of elements of $\longmo{S_h} \setminus \lbrace a_h\rbrace$ in $S_j$. Let us show that every element of $\longmo{S_h}$ has the corresponding element in $\longmo{U_h}$. Indeed, if $a \in \longmo{S_h}$ and $a \neq t$, then $a$ can be considered as an element of $\longmo{U_h}$. And $t \in \longmo{S_h}$ naturally corresponds to $b \in \longmo{U_h}$. Hence, if $a_j$ is covered by images of elements of $\longmo{S_h} \setminus \lbrace a_h\rbrace$ in $S_j$, then $a_j$ is covered by images of the corresponding elements of $\longmo{U_h} \setminus \lbrace a_h\rbrace$ in $U_j$. So, $a_h \mapsto a_j$ is not an admissible replacement in $U_h$. A contradiction.

Assume $a_h \mapsto a_j$ is an admissible replacement in $S_h$. Let us show that $a_h \mapsto a_j$ is an admissible replacement in $U_h$. Assume the contrary. Namely, assume that $a_j$ is covered by images of elements of $\longmo{U_h} \setminus \lbrace a_h\rbrace$ in $U_j$. Then it follows from the results of Section~\ref{mt_configurations} that $a_j$ is covered by images of neighbours of $a_h$ that belong to $\longmo{U_h}$.

Let $c \in \longmo{U_h}$ be a right neighbour of $a_h$. Then, by definition, $c$ starts from the right of the beginning of $a_h$.
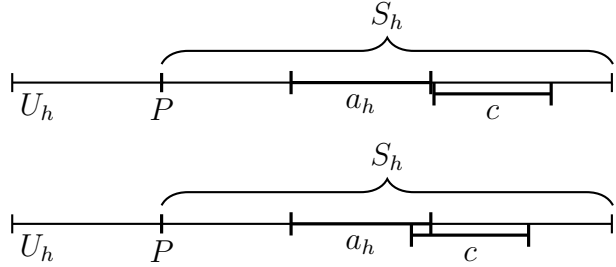
\begin{center}
\begin{tikzpicture}
\node at (0, 0) [below, xshift=10] {$U_h$};
\draw[|-|, black, thick] (0,0)--(8,0);
\draw[|-|, black, very thick] (3.7, 0)--(5.6, 0) node[midway, below] {$a_h$};
\draw[|-|, black, very thick] (5.6, -0.15)--(7.2, -0.15) node[midway, below] {$c$};
\draw [thick, decorate, decoration={brace, amplitude=10pt, raise=7pt}] (2, 0) to node[midway, above, yshift=16pt] {$S_h$} (8, 0);
\draw[black, very thick] (2, -0.128)--(2, 0.128);
\node[below, yshift=-2] at (2, 0) {$P$};
\end{tikzpicture}

\begin{tikzpicture}
\node at (0, 0) [below, xshift=10] {$U_h$};
\draw[|-|, black, thick] (0,0)--(8,0);
\draw[|-|, black, very thick] (3.7, 0)--(5.6, 0) node[midway, below] {$a_h$};
\draw[|-|, black, very thick] (5.3, -0.15)--(6.9, -0.15) node[midway, below] {$c$};
\draw [thick, decorate, decoration={brace, amplitude=10pt, raise=7pt}] (2, 0) to node[midway, above, yshift=16pt] {$S_h$} (8, 0);
\draw[black, very thick] (2, -0.128)--(2, 0.128);
\node[below, yshift=-2] at (2, 0) {$P$};
\end{tikzpicture}
\end{center}
Since $a_h$ is contained in $S_h$, we see that $c$ is contained in $S_h$. Hence, $c \in \longmo{S_h}$.

Let $d \in \longmo{U_h}$ be a left neighbour of $a_h$. First assume that $a_h$ is separated from $b$. Then $b$ is not a neighbour of $a_h$, so, $d \neq b$. Therefore, since $d$ is not separated from $a_h$, we see that $d$ starts from the right of the beginning of $b$. Since $m$ is not a small piece, we obtain that $d$ starts from the right of the beginning of $m$.
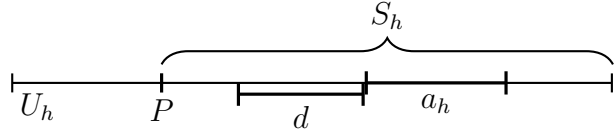
\begin{center}
\begin{tikzpicture}
\node at (0, 0) [below, xshift=10] {$U_h$};
\draw[|-|, black, thick] (0,0)--(8,0);
\draw[|-|, black, very thick] (3, -0.15)--(4.7, -0.15) node[midway, below] {$d$};
\draw[|-|, black, very thick] (4.7, 0)--(6.6, 0) node[midway, below] {$a_h$};
\draw [thick, decorate, decoration={brace, amplitude=10pt, raise=7pt}] (2, 0) to node[midway, above, yshift=16pt] {$S_h$} (8, 0);
\draw[black, very thick] (2, -0.128)--(2, 0.128);
\node[below, yshift=-2] at (2, 0) {$P$};
\end{tikzpicture}

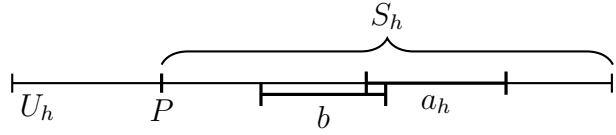
\begin{tikzpicture}
\node at (0, 0) [below, xshift=10] {$U_h$};
\draw[|-|, black, thick] (0,0)--(8,0);
\draw[|-|, black, very thick] (3.3, -0.15)--(5, -0.15) node[midway, below] {$b$};
\draw[|-|, black, very thick] (4.7, 0)--(6.6, 0) node[midway, below] {$a_h$};
\draw [thick, decorate, decoration={brace, amplitude=10pt, raise=7pt}] (2, 0) to node[midway, above, yshift=16pt] {$S_h$} (8, 0);
\draw[black, very thick] (2, -0.128)--(2, 0.128);
\node[below, yshift=-2] at (2, 0) {$P$};
\end{tikzpicture}
\end{center}
Therefore, $d$ is contained in $S_h$, so, $d \in \longmo{S_h}$. Now assume that $a_h$ is not separated from $b$.
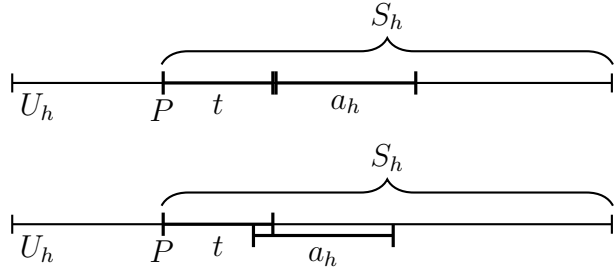
\begin{center}
\begin{tikzpicture}
\node at (0, 0) [below, xshift=10] {$U_h$};
\draw[|-|, black, thick] (0,0)--(8,0);
\draw[|-|, black, very thick] (2, 0)--(3.5, 0) node[midway, below] {$t$};
\node at (2, 0) [below, yshift=-2] {$P$};
\draw[|-|, black, very thick] (3.5, 0)--(5.4, 0) node[midway, below] {$a_h$};
\draw [thick, decorate, decoration={brace, amplitude=10pt, raise=7pt}] (2, 0) to node[midway, above, yshift=16pt] {$S_h$} (8, 0);
\end{tikzpicture}

\begin{tikzpicture}
\node at (0, 0) [below, xshift=10] {$U_h$};
\draw[|-|, black, thick] (0,0)--(8,0);
\draw[|-|, black, very thick] (2, 0)--(3.5, 0) node[midway, below] {$t$};
\node at (2, 0) [below, yshift=-2] {$P$};
\draw[|-|, black, very thick] (3.2, -0.15)--(5.1, -0.15) node[midway, below] {$a_h$};
\draw [thick, decorate, decoration={brace, amplitude=10pt, raise=7pt}] (2, 0) to node[midway, above, yshift=16pt] {$S_h$} (8, 0);
\end{tikzpicture}
\end{center}
Then $b$ is a left neighbour of $a_h$ in $U_h$. By Corollary~\ref{single_long_neighbour}, $a_h$ has a single left neighbour that belong to $\longmo{U_h}$. So, since $b \in \longmo{U_h}$, we obtain $b = d$. Therefore, $t \in \longmo{S_h}$ is a left neighbour of $a_h$ in $S_h$ that corresponds to $d$ in this case.

So, we proved that neighbours of $a_h$ in $U_h$ that belong to $\longmo{U_h}$ correspond to neighbours of $a_h$ in $S_h$ that belong to $\longmo{S_h}$. We assumed that $a_j$ in $U_j$ is covered by images of neighbours of $a_h$ in $U_h$ that belong to $\longmo{U_h}$. Therefore, $a_j$ is covered by images of the corresponding neighbours of $a_h$ in $S_h$ that belong to $\longmo{S_h}$. That is, $a_h \mapsto a_j$ is not an admissible replacement in $S_h$. A contradiction.

The second part of Lemma~\ref{admissible_replacement_in_subword} is proved in the same way.
\end{proof}

\begin{proposition}
\label{admissible_replacements_left_right}
Let $U$ be a monomial. Assume $b \in \longmo{U}$ and, moreover, $\SPM(m(b)) \geqslant 3$. Assume
\begin{equation}
\label{initial_sequence}
(U = U^{(1)} \mapsto \ldots \mapsto U^{(K)}, ((a_h^{(1)}, a_j^{(1)}), \ldots, (a_h^{(K - 1)}, a_j^{(K - 1)})), (b = b^{(1)}, \ldots, b^{(K)}))
\end{equation}
is a $(b, U)$-admissible sequence. Then $b^{(k)}$ is a single image of $b^{(k - 1)}$ in $U^{(k)}$, $k = 2, \ldots, K$, and
\begin{equation*}
b^{(k)} = p^{(k)}m(b)s^{(k)}, \ k = 1, \ldots, K,
\end{equation*}
where $p^{(k)}$ is a prefix of $b^{(k)}$, $s^{(k)}$ is a suffix of $b^{(k)}$, $p^{(k)}$ and $s^{(k)}$ are small pieces (any of them may be empty). In particular,
\begin{equation*}
\SPM(b) - 2 \leqslant \SPM(b^{(k)}) \leqslant \SPM(b) + 2, \ k = 1, \ldots, K.
\end{equation*}
\end{proposition}
\begin{proof}
The occurrence $a_h^{(1)}$ begins either from the left or from the right of the beginning of $b = b^{(1)}$ in $U = U^{(1)}$. Assume that $a_h^{(1)}$ begins from the left of the beginning of $b^{(1)}$. The replacements $a_h^{(1)} \mapsto a_j^{(1)}, \ldots, a_h^{(K - 1)} \mapsto a_j^{(K - 1)}$ can be separated into groups. Namely, there exist $k_1, \ldots, k_t$ such that $a_h^{(k)}$ begins from the left of the beginning of $b^{(k)}$ in $U^{(k)}$ when $k = 1, \ldots, k_1 - 1$, $a_h^{(k)}$ begins from the right of the beginning of $b^{(k)}$ in $U^{(k)}$ when $k = k_1, \ldots, k_2 - 1$, etc., ($k_t = K$). We call $t - 1$ \emph{the number of switchings of sides in the sequence}
\begin{equation*}
(U = U^{(1)} \mapsto \ldots \mapsto U^{(K)}, ((a_h^{(1)}, a_j^{(1)}), \ldots, (a_h^{(K - 1)}, a_j^{(K - 1)})), (b = b^{(1)}, \ldots, b^{(K)})).
\end{equation*}
We prove Proposition~\ref{admissible_replacements_left_right} by induction on $t$.

Clearly, the sequence
\begin{equation*}
(U^{(1)} \mapsto \ldots \mapsto U^{(k_1)}, ((a_h^{(1)}, a_j^{(1)}), \ldots, (a_h^{(k_1 - 1)}, a_j^{(k_1 - 1)})), (b^{(1)}, \ldots, b^{(k_1)}))
\end{equation*}
is a $(b, U)$-admissible sequence. By definition, we have $t(b) = m(b)s$, where $s$ is a small piece (possibly empty). Hence, $\SPM(t(b)) \geqslant \SPM(m(b)) \geqslant 3$. Recall that $a_h^{(k)}$ begins from the left of the beginning of $b^{(k)}$ when $k = 1, \ldots, k_1 - 1$. Therefore, by Lemma~\ref{admissible_replacements_left}, we obtain that $b^{(k)}$ is a single image of $b^{(k - 1)}$ in $U^{(k)}$, $k = 2, \ldots, k_1$, and
\begin{align*}
&b^{(k)} = p^{(k)}t(b), \textit{ where } p^{(k)} \textit{ is a small piece (possibly empty)},\\
&U^{(k)} = A^{(k)}b^{(k)}B, \textit{ where } A^{(k)} \textit{ is a prefix of } U^{(k)},\ B \textit{ is a suffix of } U^{(k)},\\
&k = 1, \ldots, k_1.
\end{align*}
So,
\begin{equation*}
b^{(k)} = p^{(k)}t(b) = p^{(k)}m(b)s, \ k = 1, \ldots, k_1.
\end{equation*}
If $t = 1$ (the number of switchings of sides in the initial sequence is equal to~$0$), then $k_1 = K$. Thus, if $t = 1$ and $a_h^{(1)}$ starts from the left of the beginning of $b^{(1)}$, we are done. The case when $t = 1$ and $a_h^{(1)}$ starts from the right of the beginning of $b^{(1)}$ is studied in the same way.

Assume $t > 1$. For now, we have
\begin{equation*}
U^{(k)} = A^{(k)}b^{(k)}B = A^{(k)}p^{(k)}t(b)B = A^{(k)}p^{(k)}m(b)sB,\ k = 1, \ldots, k_1.
\end{equation*}
By Lemma~\ref{admissible_replacements_left}, every $a_h^{(k)}$ for $k = 1, \ldots, k_1 - 1$ is contained in the subword $A^{(k)}p^{(k)}$ of $U^{(k)}$. Since $\SPM(m(b)) \geqslant 3$, every $a_h^{(k)}$ is separated from the right neighbours of $b^{(k)}$ in $U^{(k)}$ for $k = 1, \ldots, k_1 - 1$. Therefore, configuration of the type ~\ref{long_right_neighbour}---\ref{short_right_neighbour} for $b^{(1)}$ and its right neighbours in $U^{(1)}$ is the same as for $b^{(k_1)}$ and its right neighbours in $U^{(k_1)}$. So, since $i(b) = p^{(1)}m(b)$, we obtain $i(b^{(k_1)}) = p^{(k_1)}m(b)$. Hence, $\SPM(i(b^{(k_1)})) \geqslant \SPM(m(b)) \geqslant 3$. Clearly,
\begin{equation*}
(U^{(k_1)} \mapsto \ldots \mapsto U^{(k_2)}, ((a_h^{(k_1)}, a_j^{(k_1)}), \ldots, (a_h^{(k_2 - 1)}, a_j^{(k_2 - 1)})), (b^{(k_1)}, \ldots, b^{(k_2)}))
\end{equation*}
is a $(b^{(k_1)}, U^{(k_1)})$-admissible sequence. Recall that $a_h^{(k)}$ begins from the right of the beginning of $b^{(k)}$ when $k = k_1, \ldots, k_2 - 1$, and $U^{(k_1)} = A^{(k_1)}b^{(k_1)}B$. Therefore, it follows from Lemma~\ref{admissible_replacements_right} that $b^{(k)}$ is a single image of $b^{(k - 1)}$ in $U^{(k)}$, $k = k_1 + 1, \ldots, k_2$, and
\begin{align*}
&b^{(k)} = i(b^{(k_1)})s^{(k)}, \textit{ where } s^{(k)} \textit{ is a small piece (possibly empty)},\\
&U^{(k)} =  A^{(k_1)}b^{(k)}B^{(k)}, \textit{ where } A^{(k_1)} \textit{ is a prefix of } U^{(k)},\ B^{(k)} \textit{ is a suffix of } U^{(k)},\\
&k = k_1, \ldots, k_2.
\end{align*}
So,
\begin{equation*}
b^{(k)} =  i(b^{(k_1)})s^{(k)} = p^{(k_1)}m(b)s^{(k)}, \ k_1, \ldots, k_2.
\end{equation*}
If $t = 2$ (the number of switchings of sides in the initial sequence is equal to~$1$), then $k_2 = K$. Thus, if $t = 2$ and $a_h^{(1)}$ starts from the left of the beginning of $b^{(1)}$, we are done. The case when $t = 2$ and $a_h^{(1)}$ starts from the right of the beginning of $b^{(1)}$ is studied in the same way. Further we will argue by induction on $t$, the cases~$t = 1$ and~$t = 2$ are the basis of induction.

Assume $t > 2$. Let us do the step of induction. Notice that we already studied the form of $b^{(k)}$ for $k \leqslant k_2$. So, it remains to consider only $k > k_2$. Since $t > 2$, there exists $k_3$ such that $a_h^{(k)}$ begins from the left of the beginning of $b^{(k)}$ if $k = k_2, \ldots, k_3 - 1$. We proved above that
\begin{equation*}
U^{(k_2)} =  A^{(k_1)}b^{(k_2)}B^{(k_2)} = A^{(k_1)}p^{(k_1)}m(b)s^{(k_2)}B^{(k_2)}.
\end{equation*}
Notice, that both $p^{(k_1)}$ and $s^{(k_2)}$ can be empty. Then $b^{(k_2)} = m(b)$. At the same time it may happen that we have configuration~\ref{long_left_neighbour} for $b^{(k_2)}$ in $U^{(k_2)}$ and $b^{(k_2)}$ has an overlap $e$ with its neighbour of $\SPM$-measure $\geqslant \tau - 3$. Then $b^{(k_2)} = et(b^{(k_2)}) = m(b)$ and $\SPM(t(b^{(k_2)})) \leqslant \SPM(m(b))$.
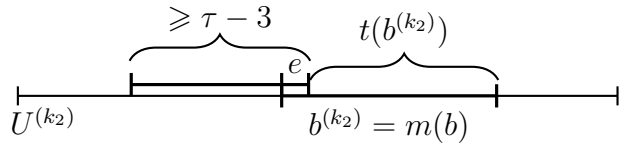
\begin{center}
\begin{tikzpicture}
\node at (0, 0) [below, xshift=10] {$U^{(k_2)}$};
\draw[|-|, black, thick] (0,0)--(8,0);
\draw[|-, black, very thick] (1.5, 0.15)--(3.5, 0.15);
\draw[|-|, black, very thick] (3.5, 0.15)--(3.9, 0.15) node[midway, above] {$e$};
\draw[|-|, black, very thick] (3.5, 0)--(6.4, 0) node[midway, below] {$b^{(k_2)} = m(b)$};
\draw [thick, decorate, decoration={brace, amplitude=10pt, raise=6pt}] (3.9, 0) to node[midway, above, yshift=14pt] {$t(b^{(k_2)})$} (6.4, 0);
\draw [thick, decorate, decoration={brace, amplitude=10pt, raise=9pt}] (1.5, 0.15) to node[midway, above, yshift=16pt] {$\geqslant \tau - 3$} (3.9, 0.15);
\end{tikzpicture}
\end{center}
So, in this case we may have $\SPM(t(b^{(k_2)})) < 3$. Then we can not apply Lemma~\ref{admissible_replacements_left} to the sequence
\begin{equation*}
(U^{(k_2)} \mapsto \ldots \mapsto U^{(k_3)}, ((a_h^{(k_2)}, a_j^{(k_2)}), \ldots, (a_h^{(k_3 - 1)}, a_j^{(k_3 - 1)})), (b^{(k_2)}, \ldots, b^{(k_3)})),
\end{equation*}
as we do above. Notice that even if $\SPM(t(b^{(k_2)})) \geqslant 3$, we do not obtain the desired result if we use Lemma~\ref{admissible_replacements_left} in this case. Informally speaking, we need to prove that the unchangeable part of $b^{(k)}$ is equal to $m(b)$. However, if we directly apply Lemma~\ref{admissible_replacements_left} in this case, we obtain  $b^{(k)} = q^{(k)}t(b^{(2)})$, where $q^{(k)}$ is a small piece (possibly empty). That is, the unchangeable part of $b^{(k)}$ becomes smaller than $m(b)$, because $t(b^{(k_2)})$ is a proper subword of $m(b)$. So, we will argue in a different way.

Consider the range of indices $k = k_1, \ldots, k_2$. By Lemma~\ref{admissible_replacements_right}, we have
\begin{equation}
\label{right_replacements}
U^{(k)} = A^{(k_1)}b^{(k)}B^{(k)} = A^{(k_1)}p^{(k_1)}m(b)s^{(k)}B^{(k)}, \ k = k_1, \ldots, k_2,
\end{equation}
and every $a_h^{(k)}$ is contained in $s^{(k)}B^{(k)}$, $k = k_1, \ldots, k_2 - 1$. Let us replace the prefix of $U^{(k)}$ that is equal to $A^{(k_1)}p^{(k_1)}$ by $A^{(1)}p^{(1)}$. Then we obtain monomials
\begin{equation*}
V^{(k)} = A^{(1)}p^{(1)}m(b)s^{(k)}B^{(k)},\ k = k_1, \ldots, k_2.
\end{equation*}
Since $a_h^{(k)}$ is contained in $s^{(k)}B^{(k)}$ and $a_h^{(k)}$ is a maximal occurrence in $U^{(k)}$, one can easily see that $a_h^{(k)}$ is a maximal occurrence in $V^{(k)}$. So, $a_h^{(k)} \mapsto a_j^{(k)}$ can be considered as a replacement in $V^{(k)}$. Since $U^{(k + 1)} = A^{(k_1)}p^{(k_1)}m(b)s^{(k + 1)}B^{(k + 1)}$ is a resulting monomial of the replacement $a_h^{(k)} \mapsto a_j^{(k)}$ in $U^{(k)}$, we see that $V^{(k + 1)} = A^{(1)}p^{(1)}m(b)s^{(k + 1)}B^{(k + 1)}$ is a resulting monomial of the replacement $a_h^{(k)} \mapsto a_j^{(k)}$ in $V^{(k)}$. Since $\SPM(m(b)) \geqslant 3$ and $a_h^{(k)} \mapsto a_j^{(k)}$ is an admissible replacement in $U^{(k)}$, it follows from Lemma~\ref{admissible_replacement_in_subword} that $a_h^{(k)} \mapsto a_j^{(k)}$ is an admissible replacement in $m(b)s^{(k)}B^{(k)}$, $k = k_1, \ldots, k_2 - 1$. Then, again by Lemma~\ref{admissible_replacement_in_subword}, we obtain that $a_h^{(k)} \mapsto a_j^{(k)}$ is an admissible replacement in $V^{(k)}$, $k = k_1, \ldots, k_2 - 1$.

Consider the range of indices $k =1, \ldots, k_1$. By Lemma~\ref{admissible_replacements_left}, we have
\begin{equation*}
U^{(k)} = A^{(k)}b^{(k)}B = A^{(k)}p^{(k)}m(b)sB,\ k = 1, \ldots, k_1,
\end{equation*}
and every $a_h^{(k)}$ is contained in $ A^{(k)}p^{(k)}$, $k = 1, \ldots, k_1 - 1$. According to the above denotations (see~\eqref{right_replacements}), we have $B = B^{(k_1)}$ and $s = s^{(k_1)}$. So, we obtain
\begin{equation*}
U^{(k)} = A^{(k)}p^{(k)}m(b)sB = A^{(k)}p^{(k)}m(b)s^{(k_1)}B^{(k_1)},\ k = 1, \ldots, k_1.
\end{equation*}
Let us replace the suffix of $U^{(k)}$ that is equal to $s^{(k_1)}B^{(k_1)}$ by $s^{(k_2)}B^{(k_2)}$. Then we obtain monomials
\begin{equation*}
W^{(k)} = A^{(k)}p^{(k)}m(b)s^{(k_2)}B^{(k_2)},\ k = 1, \ldots, k_1.
\end{equation*}
By the same argument as above, we obtain that $a_h^{(k)}$ is a maximal occurrence in $W^{(k)}$, $a_h^{(k)} \mapsto a_j^{(k)}$ is an admissible replacement in $W^{(k)}$, and $W^{(k + 1)}$ is the resulting monomial, $k = 1, \ldots, k_1 - 1$.

Let
\begin{align*}
&c^{(k)} = p^{(1)}m(b)s^{(k)},\ k = k_1, \ldots, k_2,\\
&d^{(k)} = p^{(k)}m(b)s^{(k_2)}, \ k = 1, \ldots, k_1.
\end{align*}
Since $m(b)$ is not a small piece and $b^{(k)} = p^{(k)}m(b)s^{(k)}$ is a maximal occurrence in $U^{(k)}$, we obtain
\begin{align*}
&c^{(k)} \in \mo{V^{(k)}},\\
&d^{(k)} \in \mo{W^{(k)}}.
\end{align*}
Recall that $a_h^{(k)} \in \longmo{U^{(k)}}$ is contained in $s^{(k)}B^{(k)}$, and $b^{(k + 1)} = p^{(k_1)}m(b)s^{(k + 1)}$ is an image of $b^{(k)} = p^{(k_1)}m(b)s^{(k)}$ in $U^{(k + 1)}$ for $k = k_1, \ldots, k_2 - 1$. Since $m(b)$ is not a small piece, this implies that $c^{(k + 1)} = p^{(1)}m(b)s^{(k + 1)}$ is an image of $c^{(k)}= p^{(1)}m(b)s^{(k)}$ in $V^{(k + 1)}$ for $k = k_1, \ldots, k_2 - 1$. We also know that $a_h^{(k)} \in \longmo{U^{(k)}}$ is contained in $A^{(k)}p^{(k)}$, and $b^{(k + 1)} = p^{(k + 1)}m(b)s^{(k_1)}$ is an image of $b^{(k)} = p^{(k)}m(b)s^{(k_1)}$ in $U^{(k + 1)}$ for $k = 1, \ldots, k_1 - 1$. Since $m(b)$ is not a small piece, this implies that $d^{(k + 1)} = p^{(k + 1)}m(b)s^{(k_2)}$ is an image of $d^{(k)}= p^{(k)}m(b)s^{(k_2)}$ in $W^{(k + 1)}$ for $k = 1, \ldots, k_1 - 1$.

It follows directly from the corresponding definitions that $U^{(1)} = V^{(k_1)}$, $b^{(1)} = c^{(k_1)}$, and $U^{(k_2)} = W^{(k_1)}$, $ b^{(k_2)} = d^{(k_1)}$. So, we obtain three sequences. Namely,
\begin{align*}
(&(U^{(1)} = V^{(k_1)} \mapsto \ldots \mapsto V^{(k_2)}),\\
&((a_h^{(k_1)}, a_j^{(k_1)}), \ldots, (a_h^{(k_2 - 1)}, a_j^{(k_2 - 1)})),\\
&(b^{(1)} = c^{(1)}, \ldots, c^{(k_2)}))
\end{align*}
is a $(b^{(1)}, U^{(1)})$-admissible sequence,
\begin{align*}
(&(W^{(1)} \mapsto \ldots \mapsto W^{(k_1)} = U^{(k_2)}),\\
&((a_h^{(1)}, a_j^{(1)}), \ldots, (a_h^{(k_1 - 1)}, a_j^{(k_1 - 1)})),\\
&(d^{(1)}, \ldots, d^{(k_1)} = b^{(k_2)})
\end{align*}
is a $(d^{(1)}, W^{(1)})$-admissible sequence,
\begin{align*}
(&(U^{(k_2)} \mapsto \ldots \mapsto U^{(K)}),\\
&((a_h^{(k_2)}, a_j^{(k_2)}), \ldots, (a_h^{(K - 1)}, a_j^{(K - 1)})),\\
&(b^{(k_2)}, \ldots, b^{(K)}))
\end{align*}
is a $(b^{(k_2)}, U^{(k_2)})$-admissible sequence. Notice also that
\begin{align*}
&V^{(k_2)} = A^{(1)}p^{(1)}m(b)s^{(k_2)}B^{(k_2)} = W^{(1)},\\
&c^{(k_2)} =p^{(1)}m(b)s^{(k_2)} = d^{(1)}.
\end{align*}
Hence, we can glue these sequences and obtain a $(b, U)$-admissible sequence. Namely, the sequence
\begin{align}
\begin{split}
\label{swiched_sequence}
(&X^{(1)} \ldots \mapsto X^{(K)},\\
&((x_h^{(1)}, x_j^{(1)}), \ldots, (x_h^{(K - 1)}, x_j^{(K - 1)})),\\
&(e^{(1)}, \ldots, e^{(K)})),
\end{split}
\end{align}
where
\begin{align*}
&X^{(l)} = V^{(l + k_1 - 1)} \textit{ when } l = 1, \ldots, k_2 - k_1,\\
&X^{(l)} = W^{(l - k_2 + k_1)} \textit{ when } l = k_2 - k_1 + 1, \ldots, k_2 - 1,\\
&X^{(l)} = U^{(l)} \textit{ when } l = k_2, \ldots, K,\\
&x_h^{(l)} = a_h^{(l + k_1 - 1)},\ x_j^{(l)} = a_j^{(l + k_1 - 1)} \textit{ when } l = 1, \ldots, k_2 - k_1,\\
&x_h^{(l)} = a_h^{(l - k_ 2 + k_1)},\ x_j^{(l)} = a_h^{(l - k_ 2 + k_1)} \textit{ when } l = k_2 - k_1 + 1, \ldots, k_2 - 1,\\
&x_h^{(l)} = a_h^{(l)},\ x_j^{(l)} = a_h^{(l)} \textit{ when } l = k_2, \ldots, K - 1,\\
&e^{(l)} = c^{(l + k_1 - 1)} \textit{ when } l = 1, \ldots, k_2 - k_1,\\
&e^{(l)} = d^{(l - k_2 + k_1)} \textit{ when } l = k_2 - k_1 + 1, \ldots, k_2 - 1,\\
&e^{(l)} = b^{(l)} \textit{ when } l = k_2, \ldots, K,
\end{align*}
is a $(b, U)$-admissible sequence. We see that $x_h^{(l)}$ begins from the left of the beginning of $e^{(l)}$ when $l = k_1, \ldots, k_2 - 1$, and also when $l = k_2, \ldots, k_3 - 1$. Therefore, since the initial sequence~\eqref{initial_sequence} has $t - 1$ changings of sides, by the construction, the sequence~\eqref{swiched_sequence} has $t - 2$ changings of sides. Recall that $X^{(1)} = U^{(1)}$, and $e^{(1)} = b^{(1)}$. Hence, $m(e^{(1)}) = m(b^{(1)}) = m(b)$ and $\SPM(m(e^{(1)})) \geqslant 3$. So, we can apply the induction hypothesis to the sequence~\eqref{swiched_sequence} and obtain
\begin{multline*}
e^{(k)}= p^{(k)}m(b)s^{(k)} \textit{ when } k = k_2 + 1, \ldots, K,\\
\textit{where }p^{(k)} \textit{ and } s^{(k)} \textit{are small pieces (possibly empty)}.
\end{multline*}
Since $e^{(k)} = b^{(k)}$ when $k = k_2 + 1, \ldots, K$, we are done.
\end{proof}

\subsection{Virtual members of the chart and their properties}
\label{virtual_members_subsection}
\begin{definition}
\label{virtual_members_def}
Let $U$ be a monomial, $b \in \mo{U}$. We call $b$ \emph{a virtual member of the chart of $U$} if
\begin{enumerate}
\item
$\SPM(b) \geqslant \tau - 2$.
\item
There exists a $(b, U)$-admissible sequence
\begin{equation*}
(U = U^{(1)} \mapsto \ldots \mapsto U^{(K)}, ((a_h^{(1)}, a_j^{(1)}), \ldots, (a_h^{(K - 1)}, a_j^{(K - 1)})), (b = b^{(1)}, \ldots, b^{(K)}))
\end{equation*}
such that $\SPM(b^{(K)}) \geqslant \tau$.
\end{enumerate}

We denote the set of all virtual members of the chart of $U$ by $\virt{U}$\label{virt_def}. We denote the number of virtual members of the chart of $U$ by $\nvirt{U}$, that is, $\vert \virt{U}\vert = \nvirt{U}$\label{nvirt_def}.
\end{definition}

\begin{lemma}
\label{virtual_members_corresponding}
Assume $U_h$ is a monomial, $a_h \in \mo{U_h}$. Let $a_h \mapsto a_j$ be an admissible replacement in $U_h$, $U_j$ be the resulting monomial. Assume $b^{\prime}$ is a virtual member of the chart of $U_j$ and $b^{\prime} \neq a_j$. Then there exists a unique $b \in \nfc{U_h}$ such that $b^{\prime}$ is an image of $b$ in $U_j$. Moreover, $b$ is a virtual member of the chart of $U_h$, $b \neq a_h$, and $b^{\prime}$ is a single image of $b$ in $U_j$.

Under the same conditions assume $a_j \in \mo{U_j}$. Consider the replacement $a_j \mapsto a_h$ in $U_j$. Clearly, then $U_h$ is the resulting monomial. Then the element $b \in \nfc{U_h}$ obtained above is a single image of $b^{\prime}$ in $U_h$.
\end{lemma}
\begin{proof}
We assume that $b^{\prime}$ starts from the left of the beginning of $a_j$. The case when $b^{\prime}$ starts from the right of the beginning of $a_j$ is considered in the similar (but simpler) way.

Assume $U_h = La_hR$, $U_j = La_jR$. Let $\widehat{b}_1$ be the intersection of $b^{\prime}$ and $L$, $\widehat{b}_2$ be the intersection of $b$ and $R$. First let us show that $a_j$ is not contained inside $b^{\prime}$. Assume the contrary. Then $a_j$ is a small piece (see Section~\ref{mt_configurations}, ~\ref{a_j_keep_structure}).
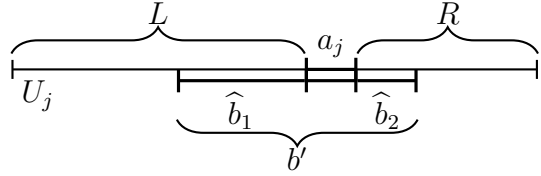
\begin{center}
\begin{tikzpicture}
\node at (0, 0) [below, xshift=10] {$U_j$};
\draw[|-|, black, thick] (0,0)--(7,0);
\draw[|-, black, very thick] (2.2, -0.15)--(3.9, -0.15) node[midway, below] {$\widehat{b}_1$};
\draw[|-|, black, very thick] (3.9, -0.15)--(4.6, -0.15);
\draw[|-|, black, very thick] (3.9, 0)--(4.6, 0) node[midway, above] {$a_j$};
\draw[-|, black, very thick] (4.6, -0.15)--(5.4, -0.15) node[midway, below] {$\widehat{b}_2$};
\draw [thick, decorate, decoration={brace, amplitude=10pt, raise=19pt, mirror}] (2.2, 0) to node[midway, below, yshift=-25pt] {$b^{\prime}$} (5.4, 0);
\draw [thick, decorate, decoration={brace, amplitude=10pt, raise=6pt}] (0, 0) to node[midway, above, yshift=13pt] {$L$} (3.9, 0);
\draw [thick, decorate, decoration={brace, amplitude=10pt, raise=6pt}] (4.6, 0) to node[midway, above, yshift=13pt] {$R$} (7, 0);
\end{tikzpicture}
\end{center}
We have
\begin{equation*}
\SPM(b^{\prime}) \leqslant \SPM(\widehat{b}_1) + \SPM(\widehat{b}_2) + \SPM(a_j) = \SPM(\widehat{b}_1) + \SPM(\widehat{b}_2) + 1.
\end{equation*}
Assume $\SPM(\widehat{b}_1) \leqslant 2$ and $\SPM(\widehat{b}_2) \leqslant 2$. Then we obtain $\SPM(b^{\prime}) \leqslant 2 + 2 + 1 = 5$. Since $\tau \geqslant 10$, we see that $\SPM(b^{\prime}) \leqslant 5 < 8 \leqslant \tau - 2$. However, since $b^{\prime}$ is a virtual member of the chart of $U_j$, we have $\SPM(b^{\prime}) \geqslant \tau - 2$, a contradiction. Therefore, $\SPM(\widehat{b}_1) \geqslant 3$ or $\SPM(\widehat{b}_2) \geqslant 3$. To be definite, assume $\SPM(\widehat{b}_1) \geqslant 3$. Then, by Lemma~\ref{non_covered_inverse_image}, there exists $c \in \longmo{U_h} \setminus \lbrace a_h \rbrace$ such that $b^{\prime}$ is an image of $c$. Since $a_j$ is covered by $b^{\prime}$, we obtain that $a_j$ is covered by an image of $c \in \longmo{U_h} \setminus \lbrace a_h\rbrace$. A contradiction with the condition that $a_h \mapsto a_j$ is an admissible replacement in $U_h$. Thus, $a_j$ is not contained inside $b^{\prime}$.

Since $a_j$ is not contained inside $b^{\prime}$, $\widehat{b}_2$ is empty, and $b^{\prime}$ and $a_j$ can have the following positions:
\begin{itemize}
\item
$b^{\prime}$ and $a_j$ are separated.
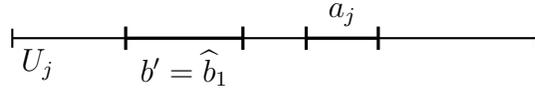
\begin{center}
\begin{tikzpicture}
\node at (0, 0) [below, xshift=10] {$U_j$};
\draw[|-|, black, thick] (0,0)--(7,0);
\draw[|-|, black, very thick] (1.5, 0)--(3.1, 0) node[midway, below] {$b^{\prime} = \widehat{b}_1$};
\draw[|-|, black, very thick] (3.9, 0)--(4.9, 0) node[midway, above] {$a_j$};
\end{tikzpicture}
\end{center}
Then $\widehat{b}_1 = b^{\prime}$. So, $\SPM(\widehat{b}_1) = \SPM(b^{\prime}) \geqslant \tau - 2$.
\item
$b^{\prime}$ and $a_j$ touch at a point.
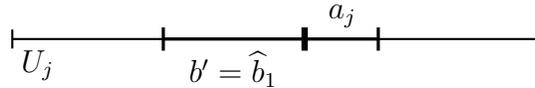
\begin{center}
\begin{tikzpicture}
\node at (0, 0) [below, xshift=10] {$U_j$};
\draw[|-|, black, thick] (0,0)--(7,0);
\draw[|-|, black, very thick] (2, 0)--(3.9, 0) node[midway, below] {$b^{\prime} = \widehat{b}_1$};
\draw[|-|, black, very thick] (3.9, 0)--(4.9, 0) node[midway, above] {$a_j$};
\end{tikzpicture}
\end{center}
Then $\widehat{b}_1 = b^{\prime}$. So, $\SPM(\widehat{b}_1) = \SPM(b^{\prime}) \geqslant \tau - 2$.
\item
$b^{\prime}$ and $a_j$ have an overlap $d$.
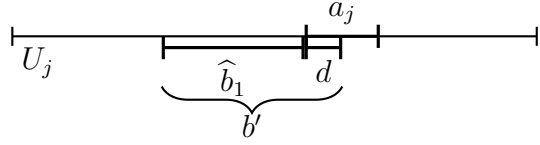
\begin{center}
\begin{tikzpicture}
\node at (0, 0) [below, xshift=10] {$U_j$};
\draw[|-|, black, thick] (0,0)--(7,0);
\draw[|-|, black, very thick] (2, -0.15)--(3.9, -0.15) node[midway, below] {$\widehat{b}_1$};
\draw[|-|, black, very thick] (3.9, -0.15)--(4.4, -0.15) node[midway, below] {$d$};
\draw[|-|, black, very thick] (3.9, 0)--(4.9, 0) node[midway, above] {$a_j$};
\draw [thick, decorate, decoration={brace, amplitude=10pt, raise=19pt, mirror}] (2, 0) to node[midway, below, yshift=-25pt] {$b^{\prime}$} (4.4, 0);
\end{tikzpicture}
\end{center}
Then $d$ is a small piece. We have $b^{\prime} = \widehat{b}_1d$, where $d$ is a suffix of $b^{\prime}$. Hence, $\SPM(\widehat{b}_1) \geqslant \SPM(b^{\prime}) - 1 \geqslant \tau - 2 - 1 = \tau - 3$.
\end{itemize}
We see that in all cases
\begin{equation*}
\SPM(\widehat{b}_1) \geqslant \tau - 3 \geqslant 7 > 3.
\end{equation*}
Recall that we can consider $\widehat{b}_1$ as an occurrence in $U_h$. Then, by Lemma~\ref{non_covered_inverse_image}, there exists a unique $b \in \nfc{U_h}$ such that $b^{\prime}$ is an image of $b$. Moreover, $b \neq a_h$,  $b^{\prime}$ is a single image of $b$ in $U_j$, and $b = \widehat{b}_1Y$, where $Y$ is either empty, or $Y$ is a small piece. It remains to show that the obtained $b$ is a virtual member of the chart of~$U_h$.

First notice the following. Since $b^{\prime}$ is a virtual member of the chart of $U_j$, there exists a $(b^{\prime}, U_j)$-admissible sequence
\begin{equation*}
(U_j = U^{(1)} \mapsto \ldots \mapsto U^{(K)}, ((a_h^{(1)}, a_j^{(1)}), \ldots, (a_h^{(K - 1)}, a_j^{(K - 1)})), (b^{\prime} = b^{(1)}, \ldots, b^{(K)}))
\end{equation*}
such that $\SPM(b^{(K)}) \geqslant \tau$. Consider the sequence
\begin{equation*}
U_h \mapsto U_j = U^{(1)} \mapsto \ldots \mapsto U^{(K)}.
\end{equation*}
Since $a_h \mapsto a_j$ is an admissible replacement in $U_h$, and $b \neq a_h$, and $b^{\prime}$ is an image of $b$, we obtain
\begin{align}
\begin{split}
\label{admissible_sequence}
(&U_h \mapsto U_j = U^{(1)} \mapsto \ldots \mapsto U^{(K)},\\
&((a_h, a_j), (a_h^{(1)}, a_j^{(1)}), \ldots, (a_h^{(K - 1)}, a_j^{(K - 1)})),\\
&(b, b^{\prime} = b^{(1)}, \ldots, b^{(K)}))
\end{split}
\end{align}
is a $(b, U)$-admissible sequence.

Assume $\SPM(b) \geqslant \tau - 2$. Then, since~\eqref{admissible_sequence} is a $(b, U)$-admissible sequence and $\SPM(b^{(K)}) \geqslant \tau$, by definition, we see that $b$ is a virtual member of the chart of~$U_h$.

Assume $\SPM(b) < \tau - 2$. Recall that $b = \widehat{b}_1Y$, where $Y$ is either empty, or $Y$ is a small piece. Hence, $\SPM(b) \geqslant \SPM(\widehat{b}_1)$. We proved above that $\SPM(\widehat{b}_1) \geqslant \tau - 3$. Therefore, $\SPM(b) \geqslant \tau - 3$. This means that the only possibility is $\SPM(b) = \tau - 3$. Recall that $\SPM(b) - 2 \leqslant \SPM(m(b)) \leqslant \SPM(b)$ (see page~\pageref{m_b_def_beginning}). So, in our case we have
\begin{equation*}
\SPM(m(b)) \geqslant \SPM(b) - 2 = \tau - 3 - 2 = \tau - 5  \geqslant 5 > 3.
\end{equation*}
Since~\eqref{admissible_sequence} is a $(b, U)$-admissible sequence, it follows from Corollary~\ref{admissible_replacements_left_right} that
\begin{equation*}
\SPM(b^{(K)}) \leqslant \SPM(m(b)) + 2 \leqslant \SPM(b) + 2 = \tau - 3 + 2 = \tau - 1.
\end{equation*}
This is a contradiction, since $\SPM(b^{(K)}) \geqslant \tau$. Therefore, the case $\SPM(b) < \tau - 2$ is not possible.

Consider the opposite replacement, namely, the replacement $a_j \mapsto a_h$ in $U_j$. Obviously, $U_h$ is the resulting monomial. Then it follows directly from the construction of $b$ that $b$ is a single image of $b^{\prime}$ in $U_h$. This completes the proof.
\end{proof}

\begin{corollary}
\label{virtual_members_number}
Assume $U_h$ is a monomial, $a_h$ is a virtual member of the chart of $U_h$. Let $a_h \mapsto a_j$ be an admissible replacement in $U_h$, $U_j$ be the resulting monomial. Then $\nvirt{U_j} \leqslant \nvirt{U_h}$. If, moreover, $a_j$ is not a virtual member of the chart of $U_j$, then $\nvirt{U_j} < \nvirt{U_h}$.
\end{corollary}
\begin{proof}
Let $b^{\prime}$ be a virtual member of the chart of $U_j$ such that $b^{\prime} \neq a_j$. Then, by Lemma~\ref{virtual_members_corresponding}, there exists $b$, a virtual member of the chart of $U_h$, such that $b\neq a_h$ and $b^{\prime}$ is a single image of $b$. This means that there exists $\mathcal{V}^{\prime} \subseteq \mathcal{V}(U_h) \setminus \lbrace a_h\rbrace$ such that taking images is a surjective mapping from $\mathcal{V}^{\prime}$ to $\mathcal{V}(U_j) \setminus \lbrace a_j\rbrace$. Hence,
\begin{equation*}
\vert \mathcal{V}^{\prime}\vert \geqslant \vert \mathcal{V}(U_j) \setminus \lbrace a_j\rbrace \vert.
\end{equation*}
Obviously, we have
\begin{equation*}
\vert \mathcal{V}(U_j) \setminus \lbrace a_j\rbrace \vert =
\begin{cases}
\nvirt{U_j} &\textit{ if } a_j \textit{ is not a virtual member}\\
&\textit{ of the chart of } U_j,\\
\nvirt{U_j} - 1 &\textit{ if } a_j \textit{ is a virtual member}\\
&\textit{ of the chart of } U_j.
\end{cases}
\end{equation*}
Since $a_h$ is a virtual member of the chart of $U_h$, we have
\begin{equation*}
\vert \mathcal{V}(U_h) \setminus \lbrace a_h\rbrace \vert = \nvirt{U_h} - 1.
\end{equation*}
Therefore,
\begin{equation*}
\nvirt{U_h} - 1 = \vert \mathcal{V}(U_h) \setminus \lbrace a_h\rbrace \vert \geqslant \vert \mathcal{V}^{\prime}\vert \geqslant \vert \mathcal{V}(U_j) \setminus \lbrace a_j\rbrace \vert \geqslant \nvirt{U_j} - 1.
\end{equation*}
So, we see that $\nvirt{U_h} \geqslant \nvirt{U_j}$.

Now assume that $a_j$ is not a virtual member of the chart of $U_j$. Then $\mathcal{V}(U_j) \setminus \lbrace a_j\rbrace = \mathcal{V}(U_j)$ and $\vert \mathcal{V}(U_j) \setminus \lbrace a_j\rbrace \vert = \nvirt{U_j}$. Therefore, we obtain
\begin{equation*}
\nvirt{U_h} - 1 = \vert \mathcal{V}(U_h) \setminus \lbrace a_h\rbrace \vert \geqslant \vert \mathcal{V}^{\prime}\vert \geqslant \vert \mathcal{V}(U_j) \setminus \lbrace a_j\rbrace \vert = \nvirt{U_j}.
\end{equation*}
Thus, $\nvirt{U_h} > \nvirt{U_j}$.
\end{proof}

Assume $U_h$ is a monomial. We can enumerate its virtual members of the chart according to their starting positions in $U_h$ from left to right. Then the $i$-th virtual member of the chart of $U_h$ (or the virtual member of the chart of $U_h$ that is located at $i$-th place) is defined.
\begin{corollary}
\label{full_virtual_members_corresponding}
Assume $U_h$ is a monomial, $a_h$ is a virtual member of the chart of $U_h$. Let $a_h$ and $a_j$ be incident monomials. Consider the replacement $a_h \mapsto a_j$ in $U_h$. Let $U_j$ be the resulting monomial. Suppose $a_j$ is a virtual member of the chart of $U_j$. Then taking  images gives a bijective correspondence between all the virtual members of the chart of $U_h$ and all the virtual members of the chart of $U_j$. Moreover, the $i$-th virtual member of the chart of $U_h$ goes to the $i$-th virtual member of the chart of $U_j$. In particular, $\nvirt{U_h} = \nvirt{U_j}$.
\end{corollary}
\begin{proof}
In order to show that taking images is a bijective mapping between $\mathcal{V}(U_h)$ and $\mathcal{V}(U_j)$ we have to prove two properties:
\begin{enumerate}
\item
If $b$ is a virtual member of the chart of $U_h$, then it has a single image in $U_j$ and this image is a virtual member of the chart. Then we obtain that taking images is a mapping from $\mathcal{V}(U_h)$ to $\mathcal{V}(U_j)$.
\item
If $b^{\prime}$ is a virtual member of the chart of $U_j$, then there exists $b$, a unique virtual member of the chart of $U_h$, such that $b^{\prime}$ is an image of $b$. Then we obtain that taking images is a bijective mapping from $\mathcal{V}(U_h)$ to $\mathcal{V}(U_j)$.
\end{enumerate}

Let us prove statement~$1$. Let $b$ be a virtual member of the chart of $U_h$. First assume $b = a_h$. Since $a_j$ is a virtual member of the chart, $a_j \in \mo{U_j}$. Then it follows from the definition of images that $a_j$ is a single image of $a_h$ in $U_j$.

Recall that since $a_h$ is a virtual member of the chart of $U_h$, we have $\SPM(a_h) \geqslant \tau - 2 \geqslant 8 > 2$. Therefore, $a_h \in \nfc{U_h}$. By the same argument, $\SPM(a_j) \geqslant\tau - 2$ and $a_j \in \nfc{U_j}$.

Assume $b \neq a_h$. Consider the replacement $a_j \mapsto a_h$ in $U_j$. Then, clearly, $U_h$ is the resulting monomial. Since $\SPM(a_j) \geqslant  \tau - 2$ and $a_h \in \nfc{U_h}$, we obtain that $a_j \mapsto a_h$ is an admissible replacement in $U_j$. We have $b$ a virtual member of the chart of $U_h$. Let us apply Lemma~\ref{virtual_members_corresponding} to the replacement $a_j \mapsto a_h$ in $U_j$. Then we obtain that there exists $b^{\prime}$,  a virtual member of the chart of $U_j$, such that $b$ is an image of $b^{\prime}$. Moreover, we can consider the opposite replacement, namely, $a_h \mapsto a_j$ in $U_h$. Then, by the second part of Lemma~\ref{virtual_members_corresponding}, $b^{\prime}$ is a single image of $b$ in $U_j$. So, statement~$1$ is proved.

Let us prove statement~$2$. Let $b^{\prime}$ be a virtual member of the chart of $U_j$. First assume $b^{\prime} = a_j$. Since $a_j \in \mo{U_j}$, we obtain that $a_j$ is an image of $a_h$. Assume $c \in \mo{U_h}$ and $c \neq a_h$. Then it follows directly from the definition of images, that an arbitrary image of $c$ is not equal to $a_j$. Therefore, $a_h$ is a unique virtual member of the chart of $U_h$ such that $a_j$ is its image is $U_j$.

Now assume $b^{\prime} \neq a_j$. Since $\SPM(a_h) \geqslant  \tau - 2$ and $a_j \in \nfc{U_j}$, we obtain that $a_h \mapsto a_j$ is an admissible replacement in $U_h$. Then it follows from Lemma~\ref{virtual_members_corresponding} that there exists a unique $b \in \nfc{U_h}$ such that $b^{\prime}$ is an image of $b$ in $U_j$. Moreover, $b$ is a virtual member of the chart of $U_h$. Obviously, every virtual member of the chart of $U_h$ belongs to $\nfc{U_h}$. Hence, $b$ is a unique virtual member of the chart of $U_h$ such that $b^{\prime}$ is an image of $b$. So, statement~$2$ is proved.


Finally we obtain that taking  images is a bijective mapping between $\mathcal{V}(U_h)$ and $\mathcal{V}(U_j)$. It follows directly from the definition of images that the $i$-th virtual member of the chart of $U_h$ goes to the $i$-th virtual member of the chart of $U_j$ under this mapping. This completes the proof of Corollary~\ref{full_virtual_members_corresponding}.
\end{proof}

In the same way as in Section~\ref{ideal_linear_descr_section}, we define a linear subspace of $\Galg$ of linear dependencies induced by multi-turns of virtual members of the chart of monomials. For every monomial of $\Fr$ we do all multi-turns of all virtual members of the chart and we consider all layouts of these multi-turns (see Definition~\ref{multiturn_def}). As a result, we obtain the set of expressions
\begin{multline}
\label{linear_dep_members2}
\GDp^{\prime} = \Bigg\lbrace \sum\limits_{j = 1}^n \alpha_jU_j \mid U_j \in \Fr, \textit{there exists an index }1\leqslant h \leqslant n \textit{ such that}\\
U_h \mapsto \sum_{\substack{j = 1 \\ j\neq h}}^{n} (-\alpha_h^{-1} \alpha_j U_j) \textit{ is a multi-turn of a virtual member of the chart of } U_h\Bigg\rbrace.
\end{multline}
We denote the linear span of $\GDp^{\prime}$ by $\langle\GDp^{\prime}\rangle$.

\begin{proposition}
\label{the_ideal_characterisation2}
The linear subspace $\langle\GDp^{\prime} \rangle \subseteq \Qalg$ is equal to the ideal $\Ideal$.
\end{proposition}
\begin{proof}
Let us show that $\GDp^{\prime} = \GDp$. Assume $U_h$ is a monomial and $a_h$ is a member of the chart of $U_h$. Then, by definition, $\SPM(a_h) \geqslant \tau$. Therefore, $a_h$ is a virtual member of the chart of $U_h$. So, every multi-turn of a member or the chart of $U_h$ is a multi-turn of a virtual member of the chart of $U_h$ as well. Thus, $\GDp \subseteq \GDp^{\prime}$.

Now let $a_h$ be a virtual member of the chart of $U_h$, $U_h = La_hR$. Assume $U_h \mapsto \sum_{\substack{j = 1 \\ j\neq h}}^{n} (-\alpha_h^{-1} \alpha_j U_j)$ is a multi-turn of $a_h$ in $U_h$. Then $\sum_{j = 1}^n \alpha_j U_j$ is a layout of this multi-turn. By definition, $\sum_{j = 1}^n \alpha_j U_j \in \GDp^{\prime}$ and every element of $\GDp^{\prime}$ can be obtained in such a way (with its own $U_h$ and $a_h$). We have to show that $\sum_{j = 1}^n \alpha_j U_j \in \GDp$.

Assume the multi-turn $U_h \mapsto \sum_{\substack{j = 1 \\ j\neq h}}^{n} (-\alpha_h^{-1} \alpha_j U_j)$ comes from an elementary multi-turn $a_h \mapsto \sum_{\substack{j = 1 \\ j\neq h}}^{n} (-\alpha_h^{-1} \alpha_j a_j)$. Then $U_j = La_jR$, $1 \leqslant j \leqslant n$. By definition, $\sum_{j = 1}^n \alpha_j a_j \in \Rel$. From Small Cancellation Axiom it follows that there exists a monomial $a_{j_0}$, $1 \leqslant j_{0} \leqslant n$, such that $\SPM(a_{j_0}) \geqslant \tau + 1$. Then $a_{j_0}$ is a member of the chart of $U_{j_0} = La_{j_0}R$, and
\begin{equation*}
U_{j_0} = La_{j_0}R \mapsto \sum_{\substack{j = 1 \\ j\neq j_0}}^{n} (-\alpha_h^{-1} \alpha_j La_jR) = \sum_{\substack{j = 1 \\ j\neq j_0}}^{n} (-\alpha_h^{-1} \alpha_j U_j)
\end{equation*}
is a multi-turn of a member of the chart of $U_{j_0}$. Clearly, $\sum_{j = 1}^n \alpha_j La_jR = \sum_{j = 1}^n \alpha_j U_j$ is a layout of this multi-turn. So, $\sum_{j = 1}^n \alpha_j U_j \in \GDp$. Thus, $\GDp^{\prime} \subseteq \GDp$.

So, finally we obtain $\GDp^{\prime} = \GDp$. Hence, by Proposition~\ref{the_ideal_characterisation}, $\langle\GDp^{\prime}\rangle = \Ideal$.
\end{proof}

The next proposition aggregates results of Subsection~\ref{coverings_subsection} and Subsection~\ref{virtual_members_subsection}.
\begin{proposition}
\label{non_increasing_parameter}
Assume $U_h$ is a monomial, $a_h$ is a virtual member of the chart of $U_h$. Let $a_h$ and $a_j$ be incident monomials. Consider the replacement $a_h \mapsto a_j$ in $U_h$. Let $U_j$ be the resulting monomial. If $a_j$ is a virtual member of the chart of $U_j$, then $\mincov{U_h} = \mincov{U_j}$ and $\nvirt{U_h} = \nvirt{U_j}$. If $a_j$ is not a virtual member of the chart of $U_j$, then either $\mincov{U_j} < \mincov{U_h}$, or $\mincov{U_j} = \mincov{U_h}$ and $\nvirt{U_j} < \nvirt{U_h}$.
\end{proposition}
\begin{proof}
First assume $a_j$ is a virtual member of the chart of $U_j$. Since $a_h$ is a virtual member of the chart of $U_h$, by definition, $\SPM(a_h) \geqslant \tau - 2$. Therefore, $a_h \in \nfc{U_h}$. Then if we consider the replacement $a_h \mapsto a_j$ in $U_h$, it follows from Lemma~\ref{minimal_coverings_property} that $\mincov{U_j} \leqslant \mincov{U_h}$. Since $a_j$ is a virtual member of the chart of $U_j$, by definition, $\SPM(a_j) \geqslant \tau - 2$. So, $a_j \in \nfc{U_j}$. Hence, we can also consider the replacement $a_j \mapsto a_h$ in $U_j$. Then, by Lemma~\ref{minimal_coverings_property}, $\mincov{U_h} \leqslant \mincov{U_j}$. So, $\mincov{U_h} = \mincov{U_j}$.

If $a_j$ is a virtual member of the chart of $U_j$, it follows directly from Corollary~\ref{full_virtual_members_corresponding} that $\nvirt{U_h} = \nvirt{U_j}$.

Now assume $a_j$ is not virtual member of the chart of $U_j$. By Lemma~\ref{minimal_coverings_property}, we have $\mincov{U_j} \leqslant \mincov{U_h}$. So, it remains to prove that if $\mincov{U_h} = \mincov{U_j}$, then $\nvirt{U_j} < \nvirt{U_h}$. Indeed, let $\mincov{U_h} = \mincov{U_j}$. If $a_h \mapsto a_j$ is not an admissible replacement in $U_h$, then it follows from Lemma~\ref{minimal_coverings_property} that $\mincov{U_j} < \mincov{U_h}$. Hence, $a_h \mapsto a_j$ is an admissible replacement in $U_h$ in the case under consideration. Then, since $a_j$ is not a virtual member of the chart of $U_j$, it follows from Corollary~\ref{virtual_members_number} that $\nvirt{U_j} < \nvirt{U_h}$. This completes the proof.
\end{proof}

\section{Transformation of a given monomial. Derived monomials}
\label{derived_monom_section}
\subsection{Derived monomials and $f$-characteristics of monomials}
\label{derived_monom_def_section}
\begin{definition}[\textbf{derived monomials}]
\label{derived_monomials}
Consider the following transformations of monomials:
\begin{enumerate}[label={(r\arabic*)}]
\item
\label{repl1}
Replacements of a virtual member of the chart by an incident monomial non-equal to~$1$ (see Definition~\ref{incident_momom}). Recall that in this case the result is always a reduced monomial (see Section~\ref{mt_configurations},~\ref{a_j_keep_structure}).
\item
\label{repl2}
Replacements of a virtual member of the chart by an incident monomial equal to~$1$ and the further cancellations (in order to obtain the reduced monomial).
\end{enumerate}
Starting with a certain monomial $U$ we consecutively apply transformations \ref{repl1}, \ref{repl2}. All the monomials that we obtain after some sequence of transformations \ref{repl1}, \ref{repl2} (including the monomial $U$ itself) are called \emph{derived monomials of $U$}.
\end{definition}

\begin{definition}
\label{f_char_def}
Let $U$ be a monomial. We introduce the following numerical characteristic of $U$ ($f$-characteristic of a monomial):
\begin{equation}
\label{f_char}
f(U) = (\mincov{U}, \nvirt{U}),
\end{equation}
where $\mincov{U}$ is the number of elements in a minimal covering of $U$, $\nvirt{U}$ is the number of virtual members of the chart of $U$. If $U_1$ and $U_2$ are monomials, we say that $f(U_1) < f(U_2)$ if and only if either $\mincov{U_1} < \mincov{U_2}$, or $\mincov{U_1} = \mincov{U_2}$ and $\nvirt{U_1} < \nvirt{U_2}$.
\end{definition}

The characteristic $f$ satisfies the following important property.
\begin{lemma}
\label{estimation_value_property}
Assume $U$ and $Z$ are monomials, $Z$ is a derived monomial of $U$. Then $f(Z) \leqslant f(U)$. Moreover, $f(Z) < f(U)$ if and only if in the corresponding sequence of replacements there exists at least one replacement of the form $La_hR \mapsto La_jR$ such that $a_h$ is a virtual member of the chart of $La_hR$ and $a_j$ is not a virtual member of the chart of $La_jR$.
\end{lemma}
\begin{proof}
The statement follows directly from Proposition~\ref{non_increasing_parameter}.
\end{proof}

In order to study the structure of $\Qalg$ as a vector space, we need to study interactions between linear generators $\GDp^{\prime}$ of the ideal $\Ideal$, namely, interactions between layouts of multi-turns of virtual members of the chart. The simplest example is as follows. Assume $T_1$ and $T_2$ are two layouts of multi-turns of virtual members of the chart of a monomial $U$. Then we have an interaction between them in the following sense: the monomial $U$ cancels out in the linear combination $T_1 - T_2$. Notice that, by Definition~\ref{derived_monomials}, all monomials of $T_1$ and $T_2$ are derived monomials of $U$. So, in order to study interactions of layouts, we need to study derived monomials.

In the simple example above the monomials of $T_1$ and $T_2$ are obtained from $U$ by a sequence of replacements of the form~\ref{repl1} and~\ref{repl2} of length one. However, we can consider more complicated linear combinations of layouts. Assume $U^{\prime}$ is some monomial in $T_1$ different from $U$. Assume $T_3$ is a layout of a multi-turn of a virtual member of the chart of $U^{\prime}$. Consider the linear combination $T_1 - T_2 - T_3$. Then $U$ and $U^{\prime}$ cancel out in this linear combination. However, the monomials of $T_3$ are obtained from $U$, using a sequence of replacements of length two. We can keep going this way and consider more and more complicated linear combinations of layouts of multi-turns. Then we deal with layouts that contain derived monomials of $U$  obtained by longer and longer sequences of replacements~\ref{repl1} and~\ref{repl2}, starting from $U$. So, we have to study sequences of replacements~\ref{repl1} and~\ref{repl2} of arbitrary length.

In Subsection~\ref{replacements_by_incident_section} we deal with sequences of replacements~\ref{repl1} and~\ref{repl2} of length two. This is our basic case. Namely, we consider two replacements in two different positions of the chart of $U$ and define their consecutive performing. For instance, assume that in the above example $f(U) = f(U^{\prime})$ and $T_1$ and $T_3$ come from two different positions of the chart of $U$ and of the chart of $U^{\prime}$. Then the monomials of $T_3$ are obtained from $U$ in this way. In Subsection~\ref{replacements_by_U_incident_section} we move on to the general case and study sequences of replacements~\ref{repl1} and~\ref{repl2} of arbitrary length.


\subsection{Replacements of virtual members of the chart by incident monomials}
\label{replacements_by_incident_section}
To be clear, in this section we use more explicit symbols for resulting monomials of replacements of maximal occurrences by incident monomials. Namely, let $Z$ be an arbitrary monomial, $c_h \in \mo{Z}$. Let $c_h$ and $c_j$ be incident monomials. Consider the replacement $c_h \mapsto c_j$ in $Z$. In this section we denote the resulting monomial of the replacement $c_h \mapsto c_j$ in $Z$ and the further cancellations (if there are any) by $Z[c_h \mapsto c_j]$.

Let $U$ be a monomial. Let $c_h$ and $d_h$ be virtual members of the chart of $U$. Let $c_h$ and $c_j$, $d_h$ and $d_j$ be incident monomials (see Remark~\ref{incident_same_indices_remark} about same indices for different pairs of incident monomials). We consider the replacements $c_h \mapsto c_j$ and $d_h \mapsto d_j$ in $U$.


Assume that $c_j$ is not a small piece, $d_j$ can be of arbitrary $\SPM$-measure. We want to perform the replacements $c_h \mapsto c_j$ and $d_h \mapsto d_j$ consecutively in the following sense. Let $\widetilde{d}_h$ be the image of $d_h$ in $U[c_h \mapsto c_j]$. Since $c_j$ is not a small piece, from the results of Section~\ref{mt_configurations} (see~\ref{keep_structure_case_l1},~\ref{keep_structure_case_l2},~\ref{keep_structure_case_r1},~\ref{keep_structure_case_r2}) it follows that
\begin{equation}
\label{image_of_right}
\widetilde{d}_h = \begin{cases}
d_h &\textit{if } c_h \textit{ and } d_h \textit{ are separated},\\
e^{\prime}(e^{-1}\cdot d_h) &\textit{if } c_h \textit{ and } d_h \textit{ are not separated},\\
&e \textit{ is the overlap of } c_h \textit{ and } d_h\\
&\textit{(empty if } c_h \textit{ and } d_h \textit{ touch at a point)},\\
&e^{\prime} \textit{ is the overlap of } c_j \textit{ and } \widetilde{d}_h\\
&\textit{(empty if } c_j \textit{ and } \widetilde{d}_h \textit{ touch at a point)}.
\end{cases}
\end{equation}
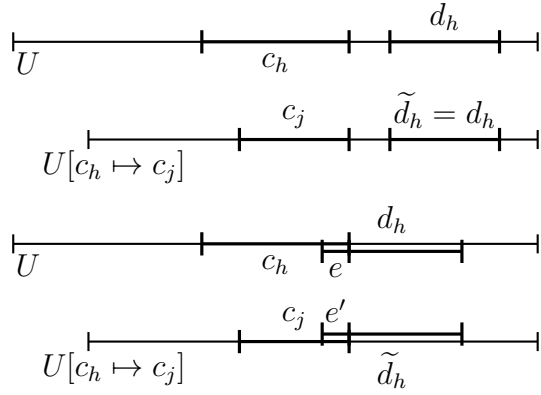
\begin{center}
\begin{tikzpicture}
\draw[|-|, black, thick] (-1.5,0)--(5.5,0);
\node[below, xshift=6] at (-1.5, 0) {$U$};
\draw[|-|, black, very thick] (1,0)--(3,0) node[midway, below] {$c_h$};
\draw[|-|, black, very thick] (3.5,0)--(5,0) node[midway, above] {$d_h$};

\draw[|-|, black, thick] (-0.5, -1.3)--(5.5, -1.3);
\node[below, xshift=10] at (-0.5, -1.3) {$U[c_h \mapsto c_j]$};
\draw[|-|, black, very thick] (1.5, -1.3)--(3, -1.3) node[midway, above] {$c_j$};
\draw[|-|, black, very thick] (3.5, -1.3)--(5, -1.3) node[midway, above] {$\widetilde{d}_h = d_h$};
\end{tikzpicture}

\begin{tikzpicture}
\draw[|-|, black, thick] (-1.5, 0)--(5.5, 0);
\node[below, xshift=6] at (-1.5, 0) {$U$};
\draw[|-|, black, very thick] (1,0)--(3,0) node[midway, below] {$c_h$};
\draw[-|, black, very thick] (3,-0.1)--(4.5,-0.1);
\draw[|-|, black, very thick] (2.6,-0.1)--(3,-0.1) node[midway, below] {$e$};
\path (2.6, 0)--(4.5, 0) node[midway, above] {$d_h$};

\draw[|-|, black, thick] (-0.5, -1.3)--(5.5, -1.3);
\node[below, xshift=10] at (-0.5, -1.3) {$U[c_h \mapsto c_j]$};
\draw[|-|, black, very thick] (1.5, -1.3)--(3,-1.3) node[midway, above] {$c_j$};
\draw[-|, black, very thick] (3, -1.2)--(4.5, -1.2);
\draw[|-|, black, very thick] (2.6, -1.2)--(3, -1.2) node[midway, above] {$e^{\prime}$};
\path (2.6, -1.3)--(4.5, -1.3) node[midway, below] {$\widetilde{d}_h$};
\end{tikzpicture}
\end{center}
Let us perform the same transformation of $d_j$ as was done on $d_h$ in order to obtain $\widetilde{d}_h$, and put
\begin{equation}
\label{corresponding_element_right}
\dbtilde{d}_j = \begin{cases}
d_j &\textit{if } c_h \textit{ and } d_h \textit{ are separated},\\
e^{\prime}\cdot e^{-1}\cdot d_j, &\textit{if } c_h \textit{ and } d_h \textit{ are not separated},\\
&e \textit{ is the overlap of } c_h \textit{ and } d_h\\
&\textit{(empty if } c_h \textit{ and } d_h \textit{ touch at a point)},\\
&e^{\prime} \textit{ is the overlap of } c_j \textit{ and } \widetilde{d}_h\\
&\textit{(empty if } c_j \textit{ and } \widetilde{d}_h \textit{ touch at a point)}.
\end{cases}
\end{equation}
Then, by Compatibility Axiom, we obtain that $\widetilde{d}_h$ and $\dbtilde{d}_j$ are incident monomials (see Lemma~\ref{second_replacement_incident} for a detailed proof). So, we have a replacement of incident monomials $\widetilde{d}_h \mapsto \dbtilde{d}_j$ in $U[c_h \mapsto c_j]$. We perform this replacement and obtain the monomial $U[c_h \mapsto c_j][\widetilde{d}_h \mapsto \dbtilde{d}_j]$. Notice that we do not claim that $\widetilde{d}_h$ is necessarily a virtual member of the chart of $U[c_h \mapsto c_j]$.

Assume that $d_j$ is not a small piece, $c_j$ can be of arbitrary $\SPM$-measure. We can do the same procedure, starting from the replacement $d_h \mapsto d_j$. Let $\widetilde{c}_h$ be the image of $c_h$ in $U[d_h \mapsto d_j]$. Since $d_j$ is not a small piece, we obtain as above
\begin{equation}
\label{image_of_left}
\widetilde{c}_h = \begin{cases}
c_h &\textit{if } c_h \textit{ and } d_h \textit{ are separated},\\
(c_h\cdot e^{-1})e^{\prime\prime} &\textit{if } c_h \textit{ and } d_h \textit{ are not separated},\\
&e \textit{ is the overlap of } c_h \textit{ and } d_h\\
&\textit{(empty if } c_h \textit{ and } d_h \textit{ touch at a point)},\\
&e^{\prime\prime} \textit{ is the overlap of } \widetilde{c}_h \textit{ and } d_j\\
&\textit{(empty if } \widetilde{c}_h \textit{ and } d_j \textit{ touch at a point)}.
\end{cases}
\end{equation}
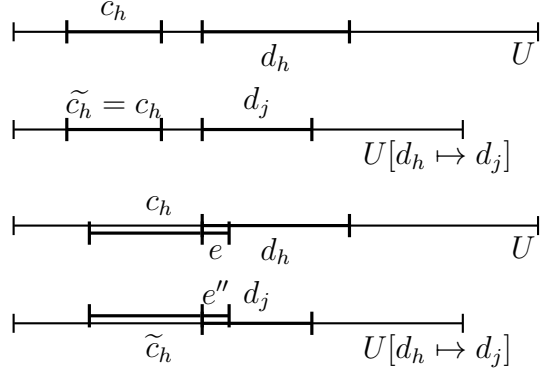
\begin{center}
\begin{tikzpicture}
\draw[|-|, black, thick] (0,0)--(7,0);
\node[below, xshift=-6] at (7, 0) {$U$};
\draw[|-|, black, very thick] (2.5,0)--(4.5,0) node[midway, below] {$d_h$};
\draw[|-|, black, very thick] (0.7,0)--(2,0) node[midway, above] {$c_h$};

\draw[|-|, black, thick] (0,-1.3)--(6,-1.3);
\node[below, xshift=-10] at (6, -1.3) {$U[d_h\mapsto d_j]$};
\draw[|-|, black, very thick] (2.5,-1.3)--(4,-1.3) node[midway, above] {$d_j$};
\draw[|-|, black, very thick] (0.7,-1.3)--(2,-1.3) node[midway, above] {$\widetilde{c_h} = c_h$};
\end{tikzpicture}

\begin{tikzpicture}
\draw[|-|, black, thick] (0, 0)--(7, 0);
\node[below, xshift=-6] at (7, 0) {$U$};
\draw[|-|, black, very thick] (2.5, 0)--(4.5, 0) node[midway, below] {$d_h$};
\draw[|-, black, very thick] (1, -0.1)--(2.5, -0.1);
\draw[|-|, black, very thick] (2.5, -0.1)--(2.9, -0.1) node[midway, below] {$e$};
\path (1, 0)--(2.9, 0) node[midway, above] {$c_h$};

\draw[|-|, black, thick] (0,-1.3)--(6,-1.3);
\node[below, xshift=-10] at (6, -1.3) {$U[d_h\mapsto d_j]$};
\draw[|-|, black, very thick] (2.5,-1.3)--(4,-1.3) node[midway, above] {$d_j$};
\draw[|-, black, very thick] (1, -1.2)--(2.5, -1.2) ;
\draw[|-|, black, very thick] (2.5, -1.2)--(2.9, -1.2) node[midway, above] {$e^{\prime\prime}$};
\path (1, -1.3)--(2.9, -1.3) node[midway, below] {$\widetilde{c}_h$};
\end{tikzpicture}
\end{center}
Let us perform the same transformation of $c_j$ as was done on $c_h$ in order to obtain $\widetilde{c}_h$, and put
\begin{equation}
\label{corresponding_element_left}
\dbtilde{c}_j  = \begin{cases}
c_j &\textit{if } c_h \textit{ and } d_h \textit{ are separated},\\
c_j \cdot e^{-1}\cdot e^{\prime\prime} &\textit{if } c_h \textit{ and } d_h \textit{ are not separated},\\
&e \textit{ is the overlap of } c_h \textit{ and } d_h\\
&\textit{(empty if } c_h \textit{ and } d_h \textit{ touch at a point)},\\
&e^{\prime\prime} \textit{ is the overlap of } \widetilde{c}_h \textit{ and } d_j\\
&\textit{(empty if } \widetilde{c}_h \textit{ and } d_j \textit{ touch at a point)}.
\end{cases}
\end{equation}
Then, by Compatibility Axiom, we obtain that $\widetilde{c}_h$ and $\dbtilde{c}_j$ are incident monomials (see Lemma~\ref{second_replacement_incident} for a detailed proof). So, we have a replacement of incident monomials $\widetilde{c}_h \mapsto \dbtilde{c}_j$ in $U[d_h \mapsto d_j]$. We perform this replacement and obtain the monomial $U[d_h \mapsto d_j][\widetilde{c}_h \mapsto \dbtilde{c}_j]$. Notice that we do not claim that $\widetilde{c}_h$ is necessarily a virtual member of the chart of $U[d_h \mapsto d_j]$.

\begin{lemma}
\label{second_replacement_incident}
Let $U$ be a monomial. Let $c_h$ and $d_h$ be virtual members of the chart of $U$. Assume $c_h$ starts from the left of the beginning of $d_h$. Let $c_h$, $c_j$ and $d_h$, $d_j$ be incident monomials.

Assume $c_j$ is not a small piece. Let $\widetilde{d}_h$ be the image of $d_h$ in $U[c_h\mapsto c_j]$, $\widetilde{d}_h$ is defined by formula~\eqref{image_of_right}. Assume $\dbtilde{d}_j$ is defined by formula~\eqref{corresponding_element_right}. Then $\widetilde{d}_h$, $\dbtilde{d}_j$ are incident monomials. In particular $\dbtilde{d}_j \in \Mon$.

Assume $d_j$ is not a small piece. Let $\widetilde{c}_h$ be the image of $c_h$ in $U[d_h\mapsto d_j]$, $\widetilde{c}_h$ be defined by formula~\eqref{image_of_left}. Assume $\dbtilde{c}_j$ is defined by formula~\eqref{corresponding_element_left}. Then $\widetilde{c}_h$, $\dbtilde{c}_j$ are incident monomials. In particular, $\dbtilde{c}_j \in \Mon$.
\end{lemma}
\begin{proof}
Let us prove the first statement of Lemma~\ref{second_replacement_incident}. So, in what follows we assume that $c_j$ is not a small piece. The second statement is proved similarly.

Assume that $c_h$ and $d_h$ are separated. Then $\widetilde{c}_h = c_h$, $\dbtilde{c}_j = c_j$, $\widetilde{c}_h = c_h$, $\dbtilde{c}_j = c_j$. So, the statement is trivial.

Assume that $c_h$ and $d_h$ are not separated. Let us show that $\widetilde{d}_h$, $\dbtilde{d}_j$ are incident monomials. We see that $e$ is a prefix of $d_h$ (possibly empty). So, by Lemma~\ref{compatibility_for_incedent}, we obtain that $e^{-1}\cdot d_h$ and $e^{-1}\cdot d_j$ are incident monomials.

Since $e$ is a small piece and $\SPM(d_h) \geqslant \tau - 2$, we see that $\SPM(e^{-1}\cdot d_h) \geqslant \tau - 3 \geqslant 7$. That is, $e^{-1}\cdot d_h$ is not a small piece. Notice that, by definition, $e^{\prime}$ and $e^{-1}\cdot d_h$ have no cancellations from the left and $\widetilde{d}_h = e^{\prime}(e^{-1}\cdot d_h) \in \Mon$. Hence, it follows from  Lemma~\ref{compatibility_for_incedent} that $e^{\prime}(e^{-1}\cdot d_h)$ and $e^{\prime}\cdot e^{-1}\cdot d_j$ are incident monomials. That is, $\widetilde{d}_h$ and $\dbtilde{d}_j$ are incident monomials. In particular, $\dbtilde{d}_j$ is a monomial of a polynomial from $\Rel$. Hence, $\dbtilde{d}_j \in \Mon$.
\end{proof}

\begin{lemma}
\label{two_replacements_incident}
Let $U$ be a monomial. Let $c_h$ and $d_h$ be virtual members of the chart of $U$. Assume $c_h$ starts from the left of the beginning of $d_h$. Let $c_h$, $c_j$ and $d_h$, $d_j$ be incident monomials. Assume $\SPM(c_j) \geqslant 4$, $\SPM(d_j) \geqslant 4$. Let $\widetilde{d}_h$ be the image of $d_h$ in $U[c_h\mapsto c_j]$, $\widetilde{d}_h$ is defined by formula~\eqref{image_of_right}. Assume $\dbtilde{d}_j$ is defined by formula~\eqref{corresponding_element_right}.  Let $\widetilde{c}_h$ be the image of $c_h$ in $U[d_h\mapsto d_j]$, $\widetilde{c}_h$ be defined by formula~\eqref{image_of_left}. Assume $\dbtilde{c}_j$ is defined by formula~\eqref{corresponding_element_left}. Then the following properties hold:
\begin{enumerate}[label={(\arabic*)}]
\item
\label{two_replacements_incident2}
$U[c_h \mapsto c_j][\widetilde{d}_h \mapsto \dbtilde{d}_j] = U[d_h \mapsto d_j][\widetilde{c}_h \mapsto \dbtilde{c}_j]$.

In what follows in Lemma~\ref{two_replacements_incident} we denote the monomial $U[c_h \mapsto c_j][\widetilde{d}_h \mapsto \dbtilde{d}_j] = U[d_h \mapsto d_j][\widetilde{c}_h \mapsto \dbtilde{c}_j]$ by $Z$.

\item
\label{two_replacements_incident3}
The occurrence $\dbtilde{d}_j$ is the single image of $d_j$ in $Z$ under the replacement $\widetilde{c}_h \mapsto \dbtilde{c}_j$ in $U[d_h \mapsto d_j]$. Similarly, $\dbtilde{c}_j$ is the single image of $c_j$ in $Z$ under the replacement $\widetilde{d}_h \mapsto \dbtilde{d}_j$ in $U[c_h \mapsto c_j]$.
\item
\label{two_replacements_incident4}
$\dbtilde{c}_j \in \mo{Z}$, $\dbtilde{d}_j \in \mo{Z}$, $\dbtilde{c}_j$ and $\dbtilde{d}_j$ are different maximal occurrences in $Z$.
\item
\label{two_replacements_incident5}
$\SPM(c_j) - 1 \leqslant \SPM(\dbtilde{c}_j) \leqslant \SPM(c_j) + 1$, $\SPM(d_j) - 1 \leqslant \SPM(\dbtilde{d}_j) \leqslant \SPM(d_j) + 1$.
\item
\label{two_replacements_incident6}
Since $\dbtilde{c}_j, \dbtilde{d}_j \in \mo{Z}$ and $\SPM(\dbtilde{c}_j) \geqslant 3$, $\SPM(\dbtilde{d}_j) \geqslant 3$, we consider the replacements in $Z$ in the other direction. Namely, we consider the replacements $\dbtilde{c}_j \mapsto \widetilde{c}_h$ and $\dbtilde{d}_j \mapsto \widetilde{d}_h$ in $Z$. The resulting monomial of the replacement $\dbtilde{c}_j \mapsto \widetilde{c}_h$ in $Z$ is, clearly, $U[d_h \mapsto d_j]$. The resulting monomial of the replacement $\dbtilde{d}_j \mapsto \widetilde{d}_h$ in $Z$ is, clearly, $U[c_h \mapsto c_j]$. Then $d_j$ is the single image of $\dbtilde{d}_j$ in $U[d_h \mapsto d_j]$ under the replacement $\dbtilde{c}_j \mapsto \widetilde{c}_h$ in $Z$. Similarly, $c_j$ is the single image of $\dbtilde{c}_j$ in $U[c_h \mapsto c_j]$ under the replacement $\dbtilde{d}_j \mapsto \widetilde{d}_h$ in $Z$.
\end{enumerate}
\end{lemma}
\begin{proof}
~\paragraph*{(1)} Assume $c_h$ and $d_h$ are separated in $U$. Then $U$ can be written as $U = Lc_hMd_hR$. We obtain  $U[c_h \mapsto c_j] = Lc_jMd_hR$. Since $c_j \neq 1$, the monomial $Lc_jMd_hR$ is reduced. Similarly, we obtain  $U[d_h \mapsto d_j] = Lc_hMd_jR$. Since $d_j \neq 1$, the monomial $Lc_hMd_jR$ is reduced. Hence,
\begin{align*}
&U[c_h \mapsto c_j][\widetilde{d}_h \mapsto \dbtilde{d}_j] = Lc_jMd_jR,\\
&U[d_h \mapsto d_j][\widetilde{c}_h \mapsto \dbtilde{c}_j] = Lc_jMd_jR.
\end{align*}
Therefore,
\begin{equation*}
U[c_h \mapsto c_j][\widetilde{d}_h \mapsto \dbtilde{d}_j] = U[d_h \mapsto d_j][\widetilde{c}_h \mapsto \dbtilde{c}_j].
\end{equation*}

Assume that $c_h$ and $d_h$ are not separated. Recall that $e$ is the overlap of $c_h$ and $d_h$. We assume that $e$ is empty if $c_h$ and $d_h$ touch at a point. Let $\widehat{c}_h = c_h\cdot e^{-1}$ and $\widehat{d}_h = e^{-1}\cdot d_h$. Then $U$ can be written in the form $U = L\widehat{c}_he\widehat{d}_hR$.
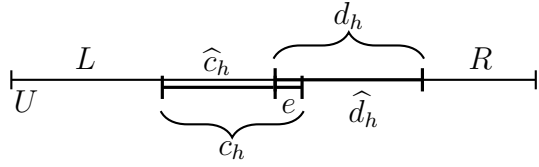
\begin{center}
\begin{tikzpicture}
\draw[|-|, black, thick] (0, 0)--(7, 0);
\node[below, xshift=6] at (0, 0) {$U$};
\draw[|-|, black, very thick] (3.5, 0)--(5.5, 0);
\path (3.9, 0)--(5.5, 0) node[midway, below] {$\widehat{d}_h$};
\draw [thick, decorate, decoration={brace, amplitude=10pt, raise=9pt, mirror}] (2, -0.1) to node[midway, below, yshift=-16pt] {$c_h$} (3.9, -0.1);
\draw[|-, black, very thick] (2, -0.1)--(3.5, -0.1) node[midway, above] {$\widehat{c}_h$};
\draw[|-|, black, very thick] (3.5, -0.1)--(3.9, -0.1) node[midway, below] {$e$};
\draw [thick, decorate, decoration={brace, amplitude=10pt, raise=8pt}] (3.5, 0) to node[midway, above, yshift=15pt] {$d_h$} (5.5, 0);
\path (0, 0)--(2, 0) node[midway, above] {$L$};
\path (5.5, 0)--(7, 0) node[midway, above] {$R$};
\end{tikzpicture}
\end{center}

Clearly, we have $U[c_h \mapsto c_j] = Lc_j\widehat{d}_hR$. Recall that $e^{\prime}$ is the overlap of $c_j$ and $\widetilde{d}_h$. So, we can represent $U[c_h \mapsto c_j]$ in the following way:
\begin{equation*}
U[c_h \mapsto c_j] = Lc_j\widehat{d}_hR =  L(c_j\cdot {e^{\prime}}^{-1})e^{\prime}\widehat{d}_hR = L(c_j\cdot {e^{\prime}}^{-1})\widetilde{d}_hR.
\end{equation*}
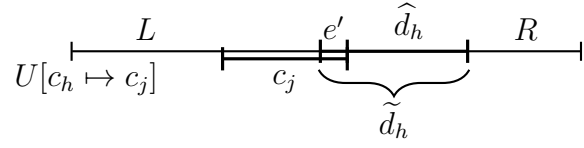
\begin{center}
\begin{tikzpicture}
\draw[|-|, black, thick] (0, 0)--(6.8, 0);
\node[below, xshift=6] at (0, 0) {$U[c_h \mapsto c_j]$};
\draw[-|, black, very thick] (3.7, 0)--(5.3, 0) node[midway, above] {$\widehat{d}_h$};
\draw[|-|, black, very thick] (2, -0.1)--(3.7, -0.1) node[midway, below] {$c_j$};
\draw[|-|, black, very thick] (3.3, 0)--(3.7, 0) node[midway, above] {$e^{\prime}$};
\draw [thick, decorate, decoration={brace, amplitude=10pt, raise=7pt, mirror}] (3.3, 0) to node[midway, below, yshift=-15pt] {$\widetilde{d}_h$} (5.3, 0);
\path (0, 0)--(2, 0) node[midway, above] {$L$};
\path (5.3, 0)--(6.8, 0) node[midway, above] {$R$};
\end{tikzpicture}
\end{center}
Therefore, we obtain
\begin{multline*}
U[c_h \mapsto c_j][\widetilde{d}_h \mapsto \dbtilde{d}_j] = L(c_j\cdot {e^{\prime}}^{-1})\dbtilde{d}_jR =\\
= L(c_j\cdot {e^{\prime}}^{-1})(e^{\prime}\cdot e^{-1}\cdot d_j)R = Lc_j\cdot e^{-1}\cdot d_jR.
\end{multline*}

Similarly, we have $U[d_h \mapsto d_j] = L\widehat{c}_hd_jR$. Recall that $e^{\prime\prime}$ is the overlap of $\widetilde{c}_h$ and $d_j$. Hence, we can represent $U[d_h \mapsto d_j]$ in the following way:
\begin{equation*}
U[d_h \mapsto d_j] = L\widehat{c}_hd_jR= L\widehat{c}_he^{\prime\prime}({e^{\prime\prime}}^{-1}\cdot d_j)R = L\widetilde{c}_h({e^{\prime\prime}}^{-1}\cdot d_j)R.
\end{equation*}
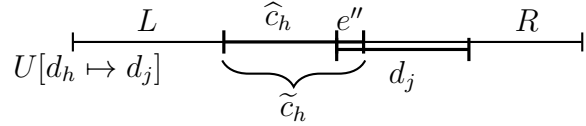
\begin{center}
\begin{tikzpicture}
\draw[|-|, black, thick] (0, 0)--(6.8, 0);
\node[below, xshift=6] at (0, 0) {$U[d_h \mapsto d_j]$};
\draw[|-|, black, very thick] (3.5, -0.1)--(5.3, -0.1) node[midway, below] {$d_j$};
\draw [thick, decorate, decoration={brace, amplitude=10pt, raise=7pt, mirror}] (2, 0) to node[midway, below, yshift=-15pt] {$\widetilde{c}_h$} (3.9, 0);
\draw[|-, black, very thick] (2, 0)--(3.5, 0) node[midway, above] {$\widehat{c}_h$};
\draw[|-|, black, very thick] (3.5, 0)--(3.9, 0) node[midway, above] {$e^{\prime\prime}$};
\path (0, 0)--(2, 0) node[midway, above] {$L$};
\path (5.3, 0)--(6.8, 0) node[midway, above] {$R$};
\end{tikzpicture}
\end{center}
Hence, we obtain
\begin{multline*}
U[d_h \mapsto d_j][\widetilde{c}_h \mapsto \dbtilde{c}_j] = L\dbtilde{c}_j({e^{\prime\prime}}^{-1}\cdot d_j)R =\\
= L(c_j\cdot e^{-1}\cdot e^{\prime\prime})({e^{\prime\prime}}^{-1}\cdot d_j)R = Lc_j\cdot e^{-1\cdot }d_jR.
\end{multline*}
Combining the results, we see that
\begin{equation*}
U[c_h \mapsto c_j][\widetilde{d}_h \mapsto \dbtilde{d}_j] = U[d_h \mapsto d_j][\widetilde{c}_h \mapsto \dbtilde{c}_j].
\end{equation*}

~\paragraph*{(2)} Assume $c_h$ and $d_h$ are separated. In this case $\dbtilde{c}_j = c_j$, $\dbtilde{d}_j = d_j$. Assume $U = Lc_hMd_hR$. Since $c_j$ is not a small piece, we obtain that $c_j$ is a maximal occurrence in $Lc_jMd_hR$ and $Lc_jMd_jR$. Similarly, since $d_j$ is not a small piece, we obtain that $d_j$ is a maximal occurrence in $Lc_hMd_jR$ and $Lc_jMd_jR$. So, the statement is obvious.

Assume $c_h$ and $d_h$ are not separated. Recall that $e$ is the overlap of $c_h$ and $d_h$. We assume that $e$ is empty if $c_h$ and $d_h$ touch at a point.

Consider the monomial $Z = U[c_h \mapsto c_j][\widetilde{d}_h \mapsto \dbtilde{d}_j] = U[d_h \mapsto d_j][\widetilde{c}_h \mapsto \dbtilde{c}_j]$. Recall that $e^{\prime}$ is the overlap of $c_j$ and $\widetilde{d}_h$ in $U[c_h \mapsto c_j]$. Let $\widehat{c}_j = c_j\cdot {e^{\prime}}^{-1}$. As above, we represent $Z$ in the form
\begin{equation*}
Z = U[c_h \mapsto c_j][\widetilde{d}_h \mapsto \dbtilde{d}_j] = L(c_j\cdot {e^{\prime}}^{-1})\dbtilde{d}_jR = L\widehat{c}_j\dbtilde{d}_jR.
\end{equation*}
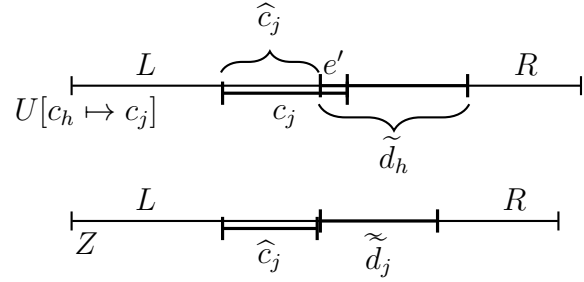
\begin{center}
\begin{tikzpicture}
\draw[|-|, black, thick] (0, 0)--(6.8, 0);
\node[below, xshift=6] at (0, 0) {$U[c_h \mapsto c_j]$};
\draw[-|, black, very thick] (3.7, 0)--(5.3, 0);
\draw[|-|, black, very thick] (2, -0.1)--(3.7, -0.1) node[midway, below] {$c_j$};
\draw[|-|, black, very thick] (3.3, 0)--(3.7, 0) node[midway, above] {$e^{\prime}$};
\draw [thick, decorate, decoration={brace, amplitude=8pt, raise=5pt}] (2, 0) to node[midway, above, yshift=15pt] {$\widehat{c}_j$} (3.3, 0);
\draw [thick, decorate, decoration={brace, amplitude=10pt, raise=7pt, mirror}] (3.3, 0) to node[midway, below, yshift=-15pt] {$\widetilde{d}_h$} (5.3, 0);
\path (0, 0)--(2, 0) node[midway, above] {$L$};
\path (5.3, 0)--(6.8, 0) node[midway, above] {$R$};

\draw[|-|, black, thick] (0, -1.8)--(6.5, -1.8);
\node[below, xshift=6] at (0, -1.8) {$Z$};
\draw[|-|, black, very thick] (3.3, -1.8)--(4.9, -1.8) node[midway, below] {$\dbtilde{d}_j$};
\draw[|-|, black, very thick] (2, -0.1 -1.8)--(3.3, -0.1-1.8) node[midway, below] {$\widehat{c}_j$};
\path (0, -1.8)--(2, -1.8) node[midway, above] {$L$};
\path (5.3, -1.8)--(6.5, -1.8) node[midway, above] {$R$};
\end{tikzpicture}
\end{center}
Recall that $e^{\prime\prime}$ is the overlap of $\widetilde{c}_h$ and $d_j$ in $U[d_h \mapsto d_j]$. Let $\widehat{d}_j = {e^{\prime\prime}}^{-1} \cdot d_j$. Similarly, we have
\begin{equation*}
Z = U[d_h \mapsto d_j][\widetilde{c}_h \mapsto \dbtilde{c}_j] = L\dbtilde{c}_j({e^{\prime\prime}}^{-1}\cdot d_j)R = L\dbtilde{c}_j\widehat{d}_jR.
\end{equation*}
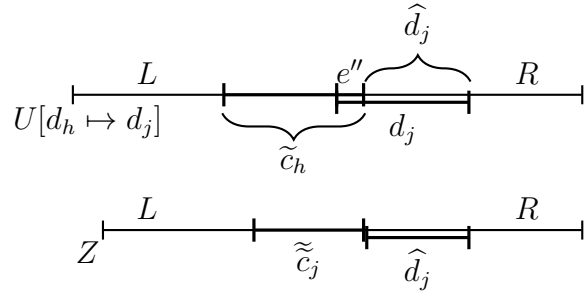
\begin{center}
\begin{tikzpicture}
\draw[|-|, black, thick] (0, 0)--(6.8, 0);
\node[below, xshift=6] at (0, 0) {$U[d_h \mapsto d_j]$};
\draw[|-|, black, very thick] (3.5, -0.1)--(5.3, -0.1) node[midway, below] {$d_j$};
\draw [thick, decorate, decoration={brace, amplitude=10pt, raise=7pt, mirror}] (2, 0) to node[midway, below, yshift=-15pt] {$\widetilde{c}_h$} (3.9, 0);
\draw [thick, decorate, decoration={brace, amplitude=10pt, raise=5pt}] (3.9, 0) to node[midway, above, yshift=15pt] {$\widehat{d}_j$} (5.3, 0);
\draw[|-, black, very thick] (2, 0)--(3.5, 0);
\draw[|-|, black, very thick] (3.5, 0)--(3.9, 0) node[midway, above] {$e^{\prime\prime}$};
\path (0, 0)--(2, 0) node[midway, above] {$L$};
\path (5.3, 0)--(6.8, 0) node[midway, above] {$R$};

\draw[|-|, black, thick] (0 + 0.4, -1.8)--(6.8, -1.8);
\node[below, xshift=6] at (0, -1.8) {$Z$};
\draw[|-|, black, very thick] (3.9, -0.1-1.8)--(5.3, -0.1-1.8) node[midway, below] {$\widehat{d}_j$};
\draw[|-|, black, very thick] (2 + 0.4, -1.8)--(3.9, -1.8) node[midway, below] {$\dbtilde{c}_j$};
\path (0, -1.8)--(2, -1.8) node[midway, above] {$L$};
\path (5.3, -1.8)--(6.8, -1.8) node[midway, above] {$R$};
\end{tikzpicture}
\end{center}
On the one hand, $Z = L\widehat{c}_j\dbtilde{d}_jR$. On the other hand, $Z = L\dbtilde{c}_j\widehat{d}_jR$. Therefore, $\widehat{c}_j\dbtilde{d}_j = \dbtilde{c}_j\widehat{d}_j$.

Notice that we do not make any additional assumptions on $\SPM$-measure of $e^{-1}$, ${e^{\prime}}^{-1}$ and ${e^{\prime\prime}}^{-1}$. In particular, we do not specially assume that they are small pieces. Therefore, initially, looking at formulas~\eqref{corresponding_element_right} and~\eqref{corresponding_element_left}, we do not have any special restrictions on $\SPM$-measure of $\dbtilde{c}_j$ and $\dbtilde{d}_j$.

Assume that $\widehat{c}_j\dbtilde{d}_j \in \Mon$. This implies that $\dbtilde{d}_j$ is a small piece (see Section~\ref{mt_configurations},~\ref{a_j_keep_structure}). Since $e^{\prime}$ is a suffix of $c_j$, by Lemma~\ref{compatibility_for_incedent}, we obtain that $c_h\cdot {e^{\prime}}^{-1}$ and $c_j \cdot {e^{\prime}}^{-1} = \widehat{c}_j$ are incident monomials. Since $\SPM(c_j) \geqslant 4$, we see that $\SPM(\widehat{c}_j) \geqslant 4 - 1 = 3$. That is, $\widehat{c}_j$ is not a small piece. Hence, since $\widehat{c}_j\dbtilde{d}_j \in \Mon$, it follows from Lemma~\ref{compatibility_for_incedent} that $c_h\cdot {e^{\prime}}^{-1}\cdot \dbtilde{d}_j$ and $\widehat{c}_j\dbtilde{d}_j$ are incident monomials. In particular, $c_h\cdot {e^{\prime}}^{-1}\cdot \dbtilde{d}_j\in \Mon$. Let $\widehat{c}_h = c_h\cdot e^{-1}$. So, we have
\begin{equation*}
c_h\cdot {e^{\prime}}^{-1}\cdot \dbtilde{d}_j = c_h\cdot {e^{\prime}}^{-1}\cdot e^{\prime}\cdot e^{-1}\cdot d_j = \widehat{c}_hd_j \in \Mon.
\end{equation*}
Therefore, $\widetilde{c}_h = \widehat{c}_he^{\prime\prime}$ and $d_j$ merge in $U[d_h \mapsto d_j]$. But this is not possible, since $d_j$ is not a small piece. A contradiction. Therefore, $\dbtilde{c}_j\widehat{d}_j = \widehat{c}_j\dbtilde{d}_j \notin \Mon$.

Since $\dbtilde{c}_j\widehat{d}_j = \widehat{c}_j\dbtilde{d}_j \notin \Mon$, we see that neither $\dbtilde{c}_j = 1$, nor $\dbtilde{d}_j = 1$. Hence, there are three possibilities for the monomial $\dbtilde{c}_j\widehat{d}_j = \widehat{c}_j\dbtilde{d}_j$:
\begin{itemize}
\item
$\dbtilde{c}_j = \widehat{c}_j$ and $\widehat{d}_j = \dbtilde{d}_j$;
\begin{center}
\begin{tikzpicture}
\draw[|-|, black, very thick] (0, 0)--(2, 0) node[midway, above] {$\widehat{c}_j$};
\draw[|-|, black, very thick] (2, 0)--(4, 0) node[midway, above] {$\dbtilde{d}_j$};
\draw[|-|, black, very thick] (0, -1)--(2, -1) node[midway, above] {$\dbtilde{c}_j$};
\draw[|-|, black, very thick] (2, -1)--(4, -1) node[midway, above] {$\widehat{d}_j$};
\end{tikzpicture}
\end{center}
\item
$\widehat{c}_j$ is a proper prefix of $\dbtilde{c}_j$, $\widehat{d}_j$ is a proper suffix of $\dbtilde{d}_j$;
\begin{center}
\begin{tikzpicture}
\draw[|-|, black, very thick] (0, 0)--(1.7, 0) node[midway, above] {$\widehat{c}_j$};
\draw[|-|, black, very thick] (1.7, 0)--(4, 0) node[midway, above] {$\dbtilde{d}_j$};
\draw[|-|, black, very thick] (0, -1)--(2.3, -1) node[midway, above] {$\dbtilde{c}_j$};
\draw[|-|, black, very thick] (2.3, -1)--(4, -1) node[midway, above] {$\widehat{d}_j$};
\end{tikzpicture}
\end{center}
\item
$\dbtilde{c}_j$ is a proper prefix of $\widehat{c}_j$, $\dbtilde{d}_j$ is a proper suffix of $\widehat{d}_j$.
\begin{center}
\begin{tikzpicture}
\draw[|-|, black, very thick] (0, 0)--(2.3, 0) node[midway, above] {$\widehat{c}_j$};
\draw[|-|, black, very thick] (2.3, 0)--(4, 0) node[midway, above] {$\dbtilde{d}_j$};
\draw[|-|, black, very thick] (0, -1)--(1.7, -1) node[midway, above] {$\dbtilde{c}_j$};
\draw[|-|, black, very thick] (1.7, -1)--(4, -1) node[midway, above] {$\widehat{d}_j$};
\end{tikzpicture}
\end{center}
\end{itemize}

Consider the first and the second cases. Since $\SPM(c_j) \geqslant 4$, we see that $\SPM(\widehat{c}_j) \geqslant 4 - 1 = 3$. Since $\SPM(d_j) \geqslant 4$, we see that $\SPM(\widehat{d}_j) \geqslant 4 - 1 = 3$. Since $\widehat{c}_j$ is a subword of $\dbtilde{c}_j$, and $\widehat{d}_j$ is a subword of $\dbtilde{d}_j$, we obtain  $\SPM(\dbtilde{c}_j) \geqslant \SPM(\widehat{c}_j) \geqslant 3$ and $\SPM(\dbtilde{d}_j) \geqslant \SPM(\widehat{d}_j) \geqslant 3$. So, $\dbtilde{c}_j$ and $\dbtilde{d}_j$ are not small pieces. Therefore, it follows from the results of Section~\ref{mt_configurations} that $\dbtilde{c}_j$ and $\dbtilde{d}_j$ are maximal occurrences in $Z$ (see~\ref{a_j_keep_structure}). Clearly, $\dbtilde{c}_j$ and $\dbtilde{d}_j$ are two different maximal occurrences in~$Z$.

Consider the third case. In this case $\widehat{c}_j$ and $\widehat{d}_j$ have a non-empty overlap. If this overlap is not a small piece, then $\widehat{c}_j$ and $\widehat{d}_j$ merge to one monomial of $\Mon$. So, $\dbtilde{c}_j\widehat{d}_j = \widehat{c}_j\dbtilde{d}_j \in \Mon$, a contradiction. Therefore, the overlap of $\widehat{c}_j$ and $\widehat{d}_j$ is a small piece. Since $\SPM(\widehat{c}_j) \geqslant 3$ and $\SPM(\widehat{c}_j) \geqslant 3$, this implies that $\SPM(\dbtilde{c}_j) \geqslant 2$ and $\SPM(\dbtilde{d}_j) \geqslant 2$. That is, $\dbtilde{c}_j$ and $\dbtilde{d}_j$ are not small pieces. Therefore, it follows from the results of Section~\ref{mt_configurations} that $\dbtilde{c}_j$ and $\dbtilde{d}_j$ are maximal occurrences in $Z$ (see~\ref{a_j_keep_structure}). This contradicts with the assumption that $\dbtilde{c}_j$ is a proper subword of $\widehat{c}_j$, and $\dbtilde{d}_j$ is a proper subword of $\widehat{d}_j$. So, the third case is not possible.

By definition, an image of $c_j$ in $Z$ is a maximal occurrence in $Z$ that contains $\widehat{c}_j$. We proved above that $\dbtilde{c}_j$ is a maximal occurrence in $Z$ that contains $\widehat{c}_j$. Therefore, $\dbtilde{c}_j$ is an image of $c_j$ in $Z$. Since $\widehat{c}_j$ is not a small piece, $\dbtilde{c}_j$ is the single maximal occurrence in $Z$ that contains $\widehat{c}_j$. Therefore, $\dbtilde{c}_j$ is the single image of $c_j$ in $Z$. In the same way one can show that $\dbtilde{d}_j$ is the single image of $d_j$ in $Z$.

~\paragraph*{(3)} In fact, we proved this property in the previous statement.

~\paragraph*{(4)} The proof of statement~(3) implies that $\dbtilde{c}_j$ differs from $c_j$ by at most one small piece at the end. So, $\SPM(c_j) - 1 \leqslant \SPM(\dbtilde{c}_j) \leqslant \SPM(c_j) + 1$. Similarly, $\dbtilde{d}_j$ differs from $d_j$ by at most one small piece at the beginning. So, $\SPM(d_j) - 1 \leqslant \SPM(\dbtilde{d}_j) \leqslant \SPM(d_j) + 1$.

~\paragraph*{(5)} If $c_h$ and $d_h$ are separated, then the statement is obvious.

Assume $c_h$ and $d_h$ are not separated. We use the same notations as in the proof of statement~\ref{two_replacements_incident3}. So, we have $Z = L\dbtilde{c}_j\widehat{d}_jR$, $\widehat{d}_j$ is a terminal subword of $\dbtilde{d}_j$. By statement~\ref{two_replacements_incident4}, we have $\dbtilde{c}_j, \dbtilde{d}_j \in \mo{Z}$. Since $d_j$ is longer than $\widehat{d}_j$ by at most one small piece, we have $\SPM(\widehat{d}_j) \geqslant \SPM(d_j) - 1 \geqslant 4 - 1 = 3$.
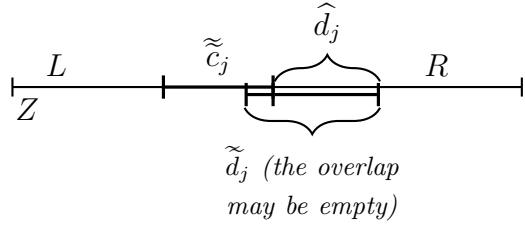
\begin{center}
\begin{tikzpicture}
\draw[|-|, black, thick] (0, 0)--(6.8, 0);
\node[below, xshift=6] at (0, 0) {$Z$};
\draw[-|, black, very thick] (3.9 -0.4, -0.1)--(5.3 -0.4, -0.1);
\draw[|-|, black, very thick] (3.5 -0.4, -0.1)--(3.9 -0.4, -0.1);
\draw[|-|, black, very thick] (2, 0)--(3.9 - 0.4, 0) node[midway, above] {$\dbtilde{c}_j$};
\path (0 - 0.4, 0)--(2 -0.4, 0) node[midway, above] {$L$};
\path (5.3 -0.4, 0)--(6.8 -0.4, 0) node[midway, above] {$R$};
\draw [thick, decorate, decoration={brace, amplitude=10pt, raise=5pt, mirror}] (3.5 - 0.4, -0.1) to node[text width=4cm, align=center, midway, below, yshift=-15pt] {\textit{\footnotesize{$\dbtilde{d}_j$ (the overlap may be empty)}}} (5.3 - 0.4, -0.1);
\draw [thick, decorate, decoration={brace, amplitude=10pt, raise=5pt}] (3.5, -0.1) to node[midway, above, yshift=15pt] {$\widehat{d}_j$} (4.9, -0.1);
\end{tikzpicture}
\end{center}
We have $U[d_h \mapsto d_j] = L\widetilde{c}_h\widehat{d}_jR$, $\widehat{d}_j$ is a terminal subword of $d_j$.
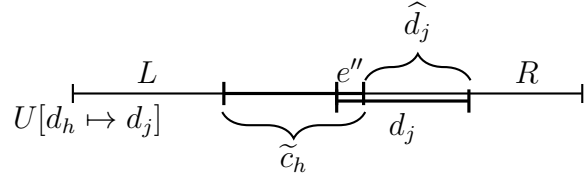
\begin{center}
\begin{tikzpicture}
\draw[|-|, black, thick] (0, 0)--(6.8, 0);
\node[below, xshift=6] at (0, 0) {$U[d_h \mapsto d_j]$};
\draw[|-|, black, very thick] (3.5, -0.1)--(5.3, -0.1) node[midway, below] {$d_j$};
\draw [thick, decorate, decoration={brace, amplitude=10pt, raise=7pt, mirror}] (2, 0) to node[midway, below, yshift=-15pt] {$\widetilde{c}_h$} (3.9, 0);
\draw [thick, decorate, decoration={brace, amplitude=10pt, raise=5pt}] (3.9, 0) to node[midway, above, yshift=15pt] {$\widehat{d}_j$} (5.3, 0);
\draw[|-, black, very thick] (2, 0)--(3.5, 0);
\draw[|-|, black, very thick] (3.5, 0)--(3.9, 0) node[midway, above] {$e^{\prime\prime}$};
\path (0, 0)--(2, 0) node[midway, above] {$L$};
\path (5.3, 0)--(6.8, 0) node[midway, above] {$R$};
\end{tikzpicture}
\end{center}
By definition, an image of $\dbtilde{d}_j$ in $U[d_h \mapsto d_j]$ under the replacement $\dbtilde{c}_j \mapsto \widetilde{c}_h$ in $Z$ is a maximal occurrence in $Z$ that contains $\widehat{d}_j$. Therefore, $d_j$ is an image of $\dbtilde{d}_j$ in $U[d_h \mapsto d_j]$ under the replacement $\dbtilde{c}_j \mapsto \widetilde{c}_h$ in $Z$. Since $\widehat{d}_j$ is not a small piece, $d_j$ is the single image of $\dbtilde{d}_j$ in $U[d_h \mapsto d_j]$ under the replacement $\dbtilde{c}_j \mapsto \widetilde{c}_h$ in $Z$.

Similarly, we obtain that $c_j$ is the single image of $\dbtilde{c}_j$ in $U[c_h \mapsto c_j]$ under the replacement $\dbtilde{d}_j \mapsto \widetilde{d}_h$ in $Z$.
\end{proof}



\begin{lemma}
\label{virtual_members_replacements_stability}
Let $U$ be a monomial. Let $c_h$ and $d_h$ be virtual members of the chart of $U$. Assume $c_h$ starts from the left of the beginning of $d_h$. Let $c_h$, $c_j$ and $d_h$, $d_j$ be incident monomials.

Assume $c_j$ is a virtual member of the chart of $U[c_h \mapsto c_j]$. Let $\widetilde{d}_h$ be the image of $d_h$ in $U[c_h\mapsto c_j]$, $\widetilde{d}_h$ is defined by formula~\eqref{image_of_right}. Assume $\dbtilde{d}_j$ is defined by formula~\eqref{corresponding_element_right}. Then $\dbtilde{d}_j$ is a virtual member of the chart of $U[c_h \mapsto c_j][\widetilde{d}_h \mapsto \dbtilde{d}_j]$ if and only if $d_j$ is a virtual member of the chart of $U[d_h \mapsto d_j]$.

Assume $d_j$ is a virtual member of the chart of $U[d_h \mapsto d_j]$. Let $\widetilde{c}_h$ be the image of $c_h$ in $U[d_h\mapsto d_j]$, $\widetilde{c}_h$ is defined by formula~\eqref{image_of_left}. Assume $\dbtilde{c}_j$ is defined by formula~\eqref{corresponding_element_left}. Then $\dbtilde{c}_j$ is a virtual member of the chart of $U[d_h \mapsto d_j][\widetilde{c}_h \mapsto \dbtilde{c}_j]$ if and only if $c_j$ is a virtual member of the chart of $U[c_h \mapsto c_j]$.
\end{lemma}
\begin{proof}
Let us prove the first part of Lemma~\ref{virtual_members_replacements_stability}. The second part is proved similarly. So, we suppose that $c_j$ is a virtual member of the chart of $U[c_h \mapsto c_j]$.

Assume that $d_j$ is a virtual member of the chart of $U[d_h \mapsto d_j]$. Let us denote the monomial $U[c_h \mapsto c_j][\widetilde{d}_h \mapsto \dbtilde{d}_j] = U[d_h \mapsto d_j][\widetilde{c}_h \mapsto \dbtilde{c}_j]$ by $Z$. Let us show that $\dbtilde{d}_j$ is a virtual member of the chart of $Z$.

First we prove that $\SPM(\dbtilde{d}_j) \geqslant \tau - 2$. Assume $c_h$ and $d_h$ are separated, then $\dbtilde{d}_j = d_j$. Since $d_j$ is a virtual member of the chart of $U[d_h \mapsto d_j]$, we have $\SPM(d_j) \geqslant \tau - 2$. Therefore, $\SPM(\dbtilde{d}_j) = \SPM(d_j) \geqslant \tau - 2$. Assume $c_h$ and $d_h$ are not separated. Consider the monomial $U[d_h \mapsto d_j]$. Recall that $e^{\prime\prime}$ is the overlap of $\widetilde{c}_h$ and $d_j$ ($e^{\prime\prime}$ is empty if $\widetilde{c}_h$ and $d_j$ touch at a point). Let $\widehat{d}_j = {e^{\prime\prime}}^{-1}\cdot d_j$.
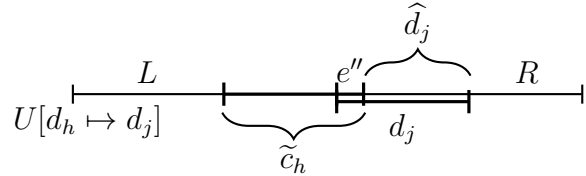
\begin{center}
\begin{tikzpicture}
\draw[|-|, black, thick] (0, 0)--(6.8, 0);
\node[below, xshift=6] at (0, 0) {$U[d_h \mapsto d_j]$};
\draw[|-|, black, very thick] (3.5, -0.1)--(5.3, -0.1) node[midway, below] {$d_j$};
\draw [thick, decorate, decoration={brace, amplitude=10pt, raise=7pt, mirror}] (2, 0) to node[midway, below, yshift=-15pt] {$\widetilde{c}_h$} (3.9, 0);
\draw [thick, decorate, decoration={brace, amplitude=10pt, raise=7pt}] (3.9, 0) to node[midway, above, yshift=15pt] {$\widehat{d}_j$} (5.3, 0);
\draw[|-, black, very thick] (2, 0)--(3.5, 0);
\draw[|-|, black, very thick] (3.5, 0)--(3.9, 0) node[midway, above] {$e^{\prime\prime}$};
\path (0, 0)--(2, 0) node[midway, above] {$L$};
\path (5.3, 0)--(6.8, 0) node[midway, above] {$R$};
\end{tikzpicture}
\end{center}
Since $d_j$ is a virtual member of the chart of $U[d_h \mapsto d_j]$, we have $\SPM(d_j) \geqslant \tau - 2$. Hence, we see that $\SPM(\widehat{d}_j) \geqslant \tau - 3$.

By statement~\ref{two_replacements_incident3} of Lemma~\ref{two_replacements_incident}, $\dbtilde{d}_j$ is the image of $d_j$ in $Z$ under the replacement $\widetilde{c}_h \mapsto \dbtilde{c}_j$ in $U[d_h \mapsto d_j]$. So, clearly, $\dbtilde{d}_j$ contains $\widehat{d}_j$. Therefore, $\SPM(\dbtilde{d}_j) \geqslant \SPM(\widehat{d}_j)$.

First assume $\SPM(\widehat{d}_j) \geqslant \tau - 2$. Then we immediately obtain $\SPM(\dbtilde{d}_j) \geqslant \SPM(\widehat{d}_j) \geqslant \tau - 2$.

Now assume that $\SPM(\widehat{d}_j) = \tau - 3$. Since $\SPM(c_h) \geqslant \tau - 2$ and $\widetilde{c}_h$ differs from $c_h$ by at most one small piece at the end, we see that $\SPM(\widetilde{c}_h) \geqslant \tau - 3$. Therefore $m(d_j)$ is an initial subword of $\widehat{d}_j$ (see the definition of a subword $m(\cdot)$ on page~\pageref{middle_subword_def}). So, $\SPM(m(d_j)) \leqslant \SPM(\widehat{d}_j) = \tau - 3$. Since $d_j$ is a virtual member of the chart of $U[d_h \mapsto d_j]$, there exists a $(U[d_h \mapsto d_j], d_j)$-admissible sequence such that the final image of $d_j$ in this sequence is of $\SPM$-measure $\geqslant \tau$. Let us denote this image by $d_j^{(K)}$. So, on the one hand, $\SPM(d_j^{(K)}) \geqslant \tau$. On the other hand, it follows from Propositions~\ref{admissible_replacements_left_right} that $d_j^{(K)} = p^{(K)}m(d_j)s^{(K)}$, where $p^{(K)}$ and $s^{(K)}$ are small pieces (possibly empty). So, we have
\begin{align*}
\SPM(d_j^{(K)}) &= \SPM(p^{(K)}m(d_j)s^{(K)}) \leqslant \SPM(p^{(K)}) + \SPM(m(d_j)) + \SPM(s^{(K)}) \leqslant\\
&\leqslant 1 + \tau - 3 + 1 = \tau - 1.
\end{align*}
A contradiction. Hence, $\SPM$-measure of $\widehat{d}_j$ is always $\geqslant \tau - 2$. Thus, $\SPM(\dbtilde{d}_j) \geqslant \SPM(\widehat{d}_j) \geqslant \tau - 2$.

One can prove similarly that $\SPM(\dbtilde{c}_j) \geqslant \tau - 2$.

Now let us prove that $\dbtilde{d}_j$ is a virtual member of the chart of $Z$. Consider the replacement $\dbtilde{c}_j \mapsto \widetilde{c}_h$ in $Z$. Clearly, $U[d_h \mapsto d_j]$ is the resulting monomial. We proved above that $\SPM(\dbtilde{c}_j) \geqslant \tau - 2$ and $\SPM(\widetilde{c}_h) \geqslant \tau - 3$. Therefore, $\dbtilde{c}_j \mapsto \widetilde{c}_h$ is an admissible replacement in $Z$. By statement~\ref{two_replacements_incident6} of Lemma~\ref{two_replacements_incident}, $d_j$ is the single image of $\dbtilde{d}_j$ in $U[d_h \mapsto d_j]$. Since $d_j$ is a virtual member of the chart of $U[d_h \mapsto d_j]$, there exists a $(U[d_h \mapsto d_j], d_j)$-admissible sequence such that the final image of $d_j$ in this sequence is of $\SPM$-measure $\geqslant \tau$. So, if we add the replacement $\dbtilde{c}_j \mapsto \widetilde{c}_h$ in $Z$ to the beginning of this sequence, we obtain a $(Z, \dbtilde{d}_j)$-admissible sequence such that the final image of $\dbtilde{d}_j$ in this sequence is of $\SPM$-measure $\geqslant \tau$. Combining this with $\SPM(\dbtilde{d}_j) \geqslant \tau - 2$, we see that $\dbtilde{d}_j$ is a virtual member of the chart of $Z$.

Let us prove the first part of Lemma~\ref{virtual_members_replacements_stability} in the other direction. We suppose that $\dbtilde{d}_j$ is a virtual member of the chart of $Z = U[c_h \mapsto c_j][\widetilde{d}_h \mapsto \dbtilde{d}_j]$. Let us show that $d_j$ is a virtual member of the chart of $U[d_h \mapsto d_j]$.

Consider the monomial $U[c_h \mapsto c_j]$. Since $c_h$ and $d_h$ are virtual members of the chart of $U$, and $c_j$ is a virtual member of the chart of $U[c_h \mapsto c_j]$, by Corollary~\ref{full_virtual_members_corresponding}, we obtain that $\widetilde{d}_h$ is a virtual member of the chart of $U[c_h \mapsto c_j]$.

Denote the monomial $U[c_h \mapsto c_j]$ by $W$. Denote $\widetilde{d}_h$ by $x_h$ and $\dbtilde{d}_j$ by $x_j$. Then $Z = W[x_h \mapsto x_j]$, $x_h$ is a virtual member of the chart of $W$, and $x_j$ is a virtual member of the chart of $W[x_h \mapsto x_j]$.
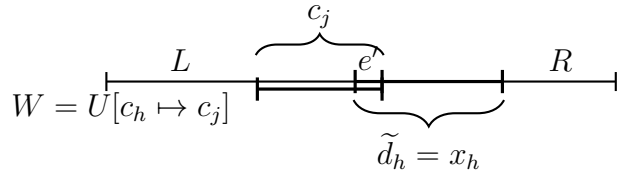
\begin{center}
\begin{tikzpicture}
\draw[|-|, black, thick] (0, 0)--(6.8, 0);
\node[below, xshift=6] at (0, 0) {$W = U[c_h \mapsto c_j]$};
\draw[-|, black, very thick] (3.7, 0)--(5.3, 0);
\draw[|-|, black, very thick] (2, -0.1)--(3.7, -0.1);
\draw[|-|, black, very thick] (3.3, 0)--(3.7, 0) node[midway, above] {$e^{\prime}$};
\draw [thick, decorate, decoration={brace, amplitude=8pt, raise=10pt}] (2, 0) to node[midway, above, yshift=16pt] {$c_j$} (3.7, 0);
\draw [thick, decorate, decoration={brace, amplitude=10pt, raise=7pt, mirror}] (3.3, 0) to node[midway, below, yshift=-15pt] {$\widetilde{d}_h = x_h$} (5.3, 0);
\path (0, 0)--(2, 0) node[midway, above] {$L$};
\path (5.3, 0)--(6.8, 0) node[midway, above] {$R$};
\end{tikzpicture}
\end{center}

By our initial assumption, $c_j$ is a virtual member of the chart of $W = U[c_h \mapsto c_j]$. Consider the replacement $c_j \mapsto c_h$ in $W$. Clearly, the resulting monomial $W[c_j \mapsto c_h]$ is equal to $U$. Hence, $c_h$ is a virtual member of the chart of $W[c_j \mapsto c_h]$. So, we obtain that $x_h$ is a virtual member of the chart of $W$, $x_j$ is a virtual member of the chart of $W[x_h \mapsto x_j]$, $c_j$ is a virtual member of the chart of $W$, $c_h$ is a virtual member of the chart of $W[c_j \mapsto c_h]$. Therefore, we can apply the part of Lemma~\ref{virtual_members_replacements_stability} that is proved above to the monomial $W$ and the replacements $c_j \mapsto c_h$ and $x_h \mapsto x_j$ in $W$.

Denote the image of $x_h$ in $W[c_j \mapsto c_h]$ by $\widetilde{x}_h$. Clearly, $\widetilde{x}_h$ is equal to $d_h$.
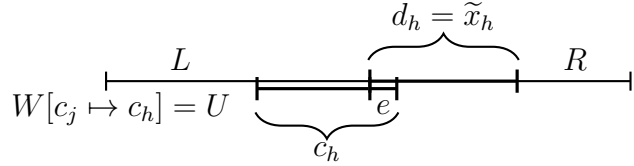
\begin{center}
\begin{tikzpicture}
\draw[|-|, black, thick] (0, 0)--(7, 0);
\node[below, xshift=6] at (0, 0) {$W[c_j \mapsto c_h] = U$};
\draw[|-|, black, very thick] (3.5, 0)--(5.5, 0);
\path (3.9, 0)--(5.5, 0);
\draw [thick, decorate, decoration={brace, amplitude=10pt, raise=9pt, mirror}] (2, -0.1) to node[midway, below, yshift=-16pt] {$c_h$} (3.9, -0.1);
\draw[|-, black, very thick] (2, -0.1)--(3.5, -0.1);
\draw[|-|, black, very thick] (3.5, -0.1)--(3.9, -0.1) node[midway, below] {$e$};
\draw [thick, decorate, decoration={brace, amplitude=10pt, raise=8pt}] (3.5, 0) to node[midway, above, yshift=15pt] {$d_h = \widetilde{x}_h$} (5.5, 0);
\path (0, 0)--(2, 0) node[midway, above] {$L$};
\path (5.5, 0)--(7, 0) node[midway, above] {$R$};
\end{tikzpicture}
\end{center}
Therefore, $\widetilde{x}_h = e\cdot {e^{\prime}}^{-1}\cdot x_h$. We put $\dbtilde{x}_j = e\cdot {e^{\prime}}^{-1}\cdot x_j$. Then, by the part of Lemma~\ref{virtual_members_replacements_stability} that is proved above, we obtain that $\dbtilde{x}_j$ is a virtual member of the chart of $W[c_j\mapsto c_h][\widetilde{x}_h \mapsto \dbtilde{x}_j]$. But we have
\begin{equation*}
\dbtilde{x}_j = e\cdot {e^{\prime}}^{-1}\cdot x_j = e\cdot {e^{\prime}}^{-1}\cdot \dbtilde{d}_j = e\cdot {e^{\prime}}^{-1}\cdot e^{\prime}\cdot e^{-1} \cdot d_j = d_j,
\end{equation*}
and $W[c_j\mapsto c_h][\widetilde{x}_h \mapsto \dbtilde{x}_j] = U[\widetilde{x}_h \mapsto \dbtilde{x}_j] = U[d_h \mapsto d_j]$. That is, $d_j$ is a virtual member of the chart of $U[d_h \mapsto d_j]$. This completes the proof.
\end{proof}

\subsection{Replacements of virtual members of the chart of a monomial $U$ by $U$-incident monomials}
\label{replacements_by_U_incident_section}
\begin{definition}
\label{U_incident_monomials}
Let $U$ be a monomial. Let $u \in \mo{U}$ and $U = LuR$. Let $a \in \Mon$. Assume there exists a sequence of monomials $m_1, \ldots, m_{n + 1} \in \Mon$ such that
\begin{enumerate}
\item
\label{U_incident_condition1}
$u = m_1$, $a = m_{n + 1}$;
\item
\label{U_incident_condition2}
$m_{i}$ and $m_{i + 1}$ are incident monomials for all $i = 1, \ldots, n$;
\item
\label{U_incident_condition3}
$m_i$ is a virtual member of the chart of $Lm_iR$ for $i = 2, \ldots, n$.
\end{enumerate}
Then the monomials $u$ and $a$ are called \emph{$U$-incident}.

Notice that if $u$ and $a$ are incident monomials, they are $U$-incident monomials as well. In this case we have the sequence of monomials of length $2$: $u = m_1$, $a = m_2$. In particular, $u$ is $U$-incident to itself. In this case we have the sequence of monomials of length $2$ as well: $u = m_1$, $u = m_2$.

If $u$ is a virtual member of the chart of $U$, then it follows directly from Definition~\ref{derived_monomials} that the monomial $LaR$ is a derived monomial of~$U$.
\medskip

Notice that if $u$ and $a$ are $U$-incident monomials and $a \in \mo{LaR}$, then $a$ and $u$ are $LaR$-incident-monomials.
\medskip

In what follows we consider replacements in $U$ such that $U = LuR$ goes to $LaR$, where $u$ and $a$ are $U$-incident monomials. We denote these replacements by $LuR \leadsto LaR$ or $u \leadsto a$\label{u_incident_replacements_def} in order to distinguish them from replacements of incident monomials. In other words, $LuR \leadsto LaR$ can be considered as a short notion for the sequence of replacements
\begin{equation*}
LuR= Lm_1R \mapsto Lm_2R\mapsto \ldots \mapsto Lm_nR = LaR.
\end{equation*}

Similarly to the above, we denote the resulting monomial of the replacement $u\leadsto a$ in $U$ and the further cancellations (if there are any) by $U[u\leadsto a]$.
\end{definition}


Let us notice that the complete analogue of property~\ref{a_j_keep_structure} from Section~\ref{mt_configurations} holds for replacements by $U$-incident monomials. Namely, we have the following statement.
\begin{lemma}
\label{U_incident_keep_structure}
Let $U$ be a monomial, $u \in \mo{U}$, $U = LuR$. Let $u$ and $a$ be $U$-incident monomials. Assume $a$ is not a small piece. Then the monomial $LaR$ is reduced and $a$ is a maximal occurrence in $LaR$.
\end{lemma}
\begin{proof}
Since $u$ and $a$ are $U$-incident monomials, there exists a sequence of monomials $m_1, \ldots, m_{n + 1}$ that satisfies the conditions of Definition~\ref{U_incident_monomials}. In particular, $u = m_1$, $a = m_{n + 1}$, and $\SPM(m_i) \geqslant \tau - 2$ for all $i = 2, \ldots, n$. Let us prove Lemma~\ref{U_incident_keep_structure} by induction on $n$. If $n = 1$, that is, $u$ and $a$ are incident monomials, the statement was proved in Section~\ref{mt_configurations}.

Assume $n > 1$. Let us split the replacement $u\leadsto a$ into two replacements $u \mapsto m_2$ and $m_2 \leadsto a$. Since $\SPM(m_2) \geqslant \tau - 2$, by the induction hypothesis, we obtain that $Lm_2R$ is a reduced monomial and $m_2$ is a maximal occurrence in $Lm_2R$. Now we consider the replacement $m_2 \leadsto a$ in $Lm_2R$. Since $a$ is not a small piece, by the induction hypothesis, we see that $LaR$ is a reduced monomial and $a$ is a maximal occurrence in $LaR$. This completes the proof.
\end{proof}

Notice that the whole list of cases in Section~\ref{mt_configurations} was based on two properties. The first is that we obtain a reduced monomial after we replace a maximal occurrence in $U$ by an incident monomial different from $1$. The second is property~\ref{a_j_keep_structure}. Therefore, using Lemma~\ref{U_incident_keep_structure}, we have the following statement.
\begin{corollary}
\label{replacements_by_U_incident_list}
We have literally the same list of possibilities for replacements of maximal occurrences in $U$ by $U$-incident monomials as we have in Section~\ref{mt_configurations} for replacements of maximal occurrences in $U$ by incident monomials.
\end{corollary}

Let $U$ be a monomial, $u$ be a maximal occurrence in $U$, $\SPM(u) \geqslant 3$, $U = LuR$. Assume $b$ is a maximal occurrence in $U$. Let $u$ and $a$ be $U$-incident monomials. We perform the replacement $u\leadsto a$ in $U$, $LaR$ is the resulting monomial. In the same way as in Definition~\ref{images_def}, we define a set of images of $b$ in the resulting monomial $LbR$.

Let $u^{(i_1)}$ and $u^{(i_2)}$ be virtual members of the chart of $U$. Assume $a^{(i_1)}$ starts from the left of the beginning of $a^{(i_2)}$. Assume $u^{(i_1)}$, $a^{(i_1)}$, and $u^{(i_2)}$, $a^{(i_2)}$ are $U$-incident monomials. We consider the replacements $u^{(i_1)} \leadsto a^{(i_1)}$ and $u^{(i_2)} \leadsto a^{(i_2)}$ in $U$. Assume $a^{(i_1)}$ is not a small piece. Let $\widetilde{u}^{(i_2)}$ be the image of $u^{(i_2)}$ in $U[u^{(i_1)} \leadsto a^{(i_1)}]$. Then, combining the results of Section~\ref{mt_configurations} and Corollary~\ref{replacements_by_U_incident_list}, we see that
\begin{equation}
\label{image_of_right_U_incident}
\widetilde{u}^{(i_2)} = \begin{cases}
u^{(i_2)} &\textit{if } u^{(i_1)} \textit{ and } u^{(i_2)} \textit{ are separated},\\
e^{\prime}(e^{-1}\cdot u^{(i_2)}) &\textit{if } u^{(i_1)} \textit{ and } u^{(i_2)} \textit{ are not separated},\\
&e \textit{ is the overlap of } u^{(i_1)} \textit{ and } u^{(i_2)}\\
&\textit{(empty if }u^{(i_1)} \textit{ and } u^{(i_2)} \textit{ touch at a point)},\\
&e^{\prime} \textit{ is the overlap of } a^{(i_1)} \textit{ and } \widetilde{u}^{(i_2)}\\
&\textit{(empty if } a^{(i_1)} \textit{ and } \widetilde{u}^{(i_2)} \textit{ touch at a point)}.
\end{cases}
\end{equation}
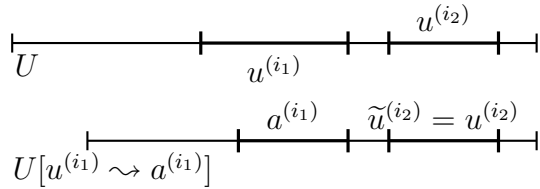
\begin{center}
\begin{tikzpicture}
\draw[|-|, black, thick] (-1.5,0)--(5.5,0);
\node[below, xshift=6] at (-1.5, 0) {$U$};
\draw[|-|, black, very thick] (1,0)--(3,0) node[midway, below] {$u^{(i_1)}$};
\draw[|-|, black, very thick] (3.5,0)--(5,0) node[midway, above] {$u^{(i_2)}$};

\draw[|-|, black, thick] (-0.5, -1.3)--(5.5, -1.3);
\node[below, xshift=10] at (-0.5, -1.3) {$U[u^{(i_1)}\leadsto a^{(i_1)}]$};
\draw[|-|, black, very thick] (1.5, -1.3)--(3, -1.3) node[midway, above] {$a^{(i_1)}$};
\draw[|-|, black, very thick] (3.5, -1.3)--(5, -1.3) node[midway, above] {$\widetilde{u}^{(i_2)} = u^{(i_2)}$};
\end{tikzpicture}

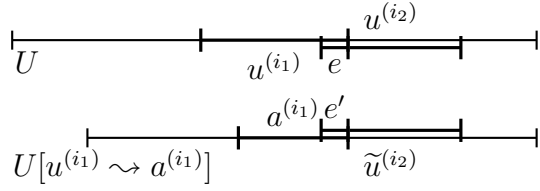
\begin{tikzpicture}
\draw[|-|, black, thick] (-1.5, 0)--(5.5, 0);
\node[below, xshift=6] at (-1.5, 0) {$U$};
\draw[|-|, black, very thick] (1,0)--(3,0) node[midway, below] {$u^{(i_1)}$};
\draw[-|, black, very thick] (3,-0.1)--(4.5,-0.1);
\draw[|-|, black, very thick] (2.6,-0.1)--(3,-0.1) node[midway, below] {$e$};
\path (2.6, 0)--(4.5, 0) node[midway, above] {$u^{(i_2)}$};

\draw[|-|, black, thick] (-0.5, -1.3)--(5.5, -1.3);
\node[below, xshift=10] at (-0.5, -1.3) {$U[u^{(i_1)}\leadsto a^{(i_1)}]$};
\draw[|-|, black, very thick] (1.5, -1.3)--(3,-1.3) node[midway, above] {$a^{(i_1)}$};
\draw[-|, black, very thick] (3, -1.2)--(4.5, -1.2);
\draw[|-|, black, very thick] (2.6, -1.2)--(3, -1.2) node[midway, above] {$e^{\prime}$};
\path (2.6, -1.3)--(4.5, -1.3) node[midway, below] {$\widetilde{u}^{(i_2)}$};
\end{tikzpicture}
\end{center}
Assume that $a^{(i_2)}$ is not a small piece. In the same way as~\eqref{image_of_left}, we have
\begin{equation}
\label{image_of_left_U_incident}
\widetilde{u}^{(i_1)} = \begin{cases}
u^{(i_1)} &\textit{if } u^{(i_1)} \textit{ and } u^{(i_2)} \textit{ are separated},\\
(u^{(i_1)}\cdot e^{-1})e^{\prime\prime} &\textit{if } u^{(i_1)} \textit{ and } u^{(i_2)} \textit{ are not separated},\\
&e \textit{ is the overlap of } u^{(i_1)} \textit{ and } u^{(i_2)}\\
&\textit{(empty if } u^{(i_1)} \textit{ and } u^{(i_2)} \textit{ touch at a point)},\\
&e^{\prime\prime} \textit{ is the overlap of } \widetilde{u}^{(i_1)} \textit{ and } a^{(i_2)}\\
&\textit{(empty if } \widetilde{u}^{(i_1)} \textit{ and } a^{(i_2)} \textit{ touch at a point)}.
\end{cases}
\end{equation}
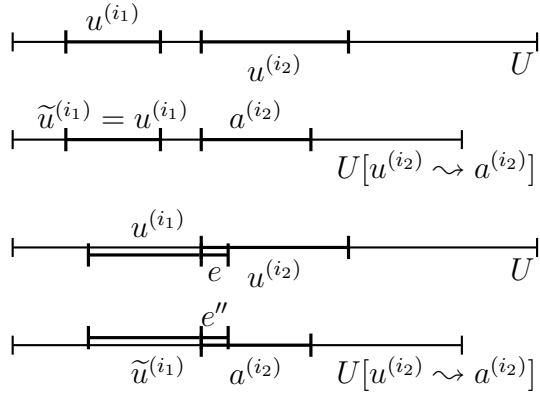
\begin{center}
\begin{tikzpicture}
\draw[|-|, black, thick] (0,0)--(7,0);
\node[below, xshift=-6] at (7, 0) {$U$};
\draw[|-|, black, very thick] (2.5,0)--(4.5,0) node[midway, below] {$u^{(i_2)}$};
\draw[|-|, black, very thick] (0.7,0)--(2,0) node[midway, above] {$u^{(i_1)}$};

\draw[|-|, black, thick] (0,-1.3)--(6,-1.3);
\node[below, xshift=-10] at (6, -1.3) {$U[u^{(i_2)}\leadsto a^{(i_2)}]$};
\draw[|-|, black, very thick] (2.5,-1.3)--(4,-1.3) node[midway, above] {$a^{(i_2)}$};
\draw[|-|, black, very thick] (0.7,-1.3)--(2,-1.3) node[midway, above] {$\widetilde{u}^{(i_1)} = u^{(i_1)}$};
\end{tikzpicture}

\begin{tikzpicture}
\draw[|-|, black, thick] (0, 0)--(7, 0);
\node[below, xshift=-6] at (7, 0) {$U$};
\draw[|-|, black, very thick] (2.5, 0)--(4.5, 0) node[midway, below] {$u^{(i_2)}$};
\draw[|-, black, very thick] (1, -0.1)--(2.5, -0.1);
\draw[|-|, black, very thick] (2.5, -0.1)--(2.9, -0.1) node[midway, below] {$e$};
\path (1, 0)--(2.9, 0) node[midway, above] {$u^{(i_1)}$};

\draw[|-|, black, thick] (0,-1.3)--(6,-1.3);
\node[below, xshift=-10] at (6, -1.3) {$U[u^{(i_2)}\leadsto a^{(i_2)}]$};
\draw[|-|, black, very thick] (2.5,-1.3)--(4,-1.3) node[midway, below] {$a^{(i_2)}$};
\draw[|-, black, very thick] (1, -1.2)--(2.5, -1.2) ;
\draw[|-|, black, very thick] (2.5, -1.2)--(2.9, -1.2) node[midway, above] {$e^{\prime\prime}$};
\path (1, -1.3)--(2.9, -1.3) node[midway, below] {$\widetilde{u}^{(i_1)}$};
\end{tikzpicture}
\end{center}
Similarly to~\eqref{corresponding_element_right}, we put
\begin{equation}
\label{corresponding_element_U_incident_right}
\dbtilde{a}^{(i_2)}= \begin{cases}
a^{(i_2)} &\textit{if } u^{(i_1)}\textit{ and } u^{(i_2)} \textit{ are separated},\\
e^{\prime}\cdot e^{-1}\cdot a^{(i_2)} &\textit{if } u^{(i_1)} \textit{ and } u^{(i_2)} \textit{ are not separated},\\
&e \textit{ is the overlap of } u^{(i_1)} \textit{ and } u^{(i_2)}\\
&\textit{(empty if } u^{(i_1)} \textit{ and } u^{(i_2)} \textit{ touch at a point)},\\
&e^{\prime} \textit{ is the overlap of } a^{(i_1)} \textit{ and } \widetilde{u}^{(i_2)}\\
&\textit{(empty if } a^{(i_1)} \textit{ and } \widetilde{u}^{(i_2)} \textit{ touch at a point)}.
\end{cases}\\
\end{equation}
Similarly to~\eqref{corresponding_element_left}, we put
\begin{equation}
\label{corresponding_element_U_incident_left}
\dbtilde{a}^{(i_1)} = \begin{cases}
a^{(i_1)} &\textit{if } u^{(i_1)} \textit{ and } u^{(i_2)} \textit{ are separated},\\
a^{(i_1)} \cdot e^{-1}\cdot e^{\prime\prime} &\textit{if } u^{(i_1)} \textit{ and } u^{(i_2)} \textit{ are not separated},\\
&e \textit{ is the overlap of } u^{(i_1)} \textit{ and } u^{(i_2)}\\
&\textit{(empty if } u^{(i_1)} \textit{ and } u^{(i_2)} \textit{ touch at a point)},\\
&e^{\prime\prime} \textit{ is the overlap of } \widetilde{u}^{(i_1)} \textit{ and } a^{(i_2)}\\
&\textit{(empty if } \widetilde{u}^{(i_1)} \textit{ and } a^{(i_2)} \textit{ touch at a point)}.
\end{cases}
\end{equation}

We want to perform the replacements $u^{(i_1)} \leadsto a^{(i_1)}$ and $u^{(i_2)} \leadsto a^{(i_2)}$ consecutively, similarly to replacements of incident monomials in Subsection~\ref{replacements_by_incident_section}. However, the monomials $\widetilde{u}^{(i_1)}$ and $\dbtilde{a}^{(i_1)}$ are not necessarily $U[u^{(i_2)} \leadsto a^{(i_2)}]$-incident monomials. Similarly, $\widetilde{u}^{(i_2)}$ and $\dbtilde{a}^{(i_2)}$ are not necessarily $U[u^{(i_1)} \leadsto a^{(i_1)}]$-incident monomials. So, we can not perform the replacements $u^{(i_1)} \leadsto a^{(i_1)}$ and $u^{(i_2)} \leadsto a^{(i_2)}$ consecutively under the same conditions as we have in Subsection~\ref{replacements_by_incident_section}. But it turns out that if $a^{(i_1)}$ is a virtual member of the chart of $U[u^{(i_1)} \leadsto a^{(i_1)}]$, or $a^{(i_2)}$ is a virtual member of the chart of $U[u^{(i_1)} \leadsto a^{(i_1)}]$, then we obtain properties of replacements by $U$-incident monomials similar to the properties that we have for replacements by incident monomials. So, first we consider this particular case.

\begin{corollary}
\label{virtual_members_U_incident_replacements}
Let $U$ be a monomial. Let $u^{(i_1)}$ and $u^{(i_2)}$ be virtual members of the chart of $U$. Assume $u^{(i_1)}$ starts from the left of the beginning of $u^{(i_2)}$. Let $u^{(i_1)}$, $a^{(i_1)}$, and $u^{(i_2)}$, $a^{(i_2)}$ be $U$-incident monomials.

Assume $a^{(i_1)}$ is a virtual member of the chart of $U[u^{(i_1)}\leadsto a^{(i_1)}]$. Assume $\widetilde{u}^{(i_2)}$ is the image of $u^{(i_2)}$ in $U[u^{(i_1)}\leadsto a^{(i_1)}]$, $\widetilde{u}^{(i_2)}$ is defined by formula~\eqref{image_of_right_U_incident}. Let $\dbtilde{a}^{(i_2)}$ be defined by formula~\eqref{corresponding_element_U_incident_right}. Then we have
\begin{enumerate}[label={(\arabic*)}]
\item
\label{virtual_members_U_incident_replacements1}
$\widetilde{u}^{(i_2)}$ is a virtual member of the chart of $U[u^{(i_1)}\leadsto a^{(i_1)}]$.
\item
\label{virtual_members_U_incident_replacements2}
$\widetilde{u}^{(i_2)}$ and $\dbtilde{a}^{(i_2)}$ are $U[u^{(i_1)}\leadsto a^{(i_1)}]$-incident monomials. In particular, $\dbtilde{a}^{(i_2)} \in \Mon$.
\item
\label{virtual_members_U_incident_replacements3}
$\dbtilde{a}^{(i_2)}$ is a virtual member of the chart of $U[u^{(i_1)}\leadsto a^{(i_1)}][\widetilde{u}^{(i_2)}\leadsto \dbtilde{a}^{(i_2)}]$ if and only if $a^{(i_2)}$ is a virtual member of the chart of $U[u^{(i_2)}\leadsto a^{(i_2)}]$.
\end{enumerate}

Assume $a^{(i_2)}$ is a virtual member of the chart of $U[u^{(i_2)}\leadsto a^{(i_2)}]$. Assume $\widetilde{u}^{(i_1)}$ is the image of $u^{(i_1)}$ in $U[u^{(i_2)}\leadsto a^{(i_2)}]$, $\widetilde{u}^{(i_1)}$ is defined by formula~\eqref{image_of_left_U_incident}. Let $\dbtilde{a}^{(i_2)}$ be defined by formula~\eqref{corresponding_element_U_incident_right}. Then we have
\begin{enumerate}[label={(\arabic*)}]
\setcounter{enumi}{3}
\item
\label{virtual_members_U_incident_replacements4}
$\widetilde{u}^{(i_1)}$ is a virtual member of the chart of $U[u^{(i_2)}\leadsto a^{(i_2)}]$.
\item
\label{virtual_members_U_incident_replacements5}
$\widetilde{u}^{(i_1)}$ and $\dbtilde{a}^{(i_1)}$ are $U[u^{(i_2)}\leadsto a^{(i_2)}]$-incident monomials. In particular, $\dbtilde{a}^{(i_1)} \in \Mon$.
\item
\label{virtual_members_U_incident_replacements6}
$\dbtilde{a}^{(i_1)}$ is a virtual member of the chart of $U[u^{(i_2)}\leadsto a^{(i_2)}][\widetilde{u}^{(i_1)}\leadsto \dbtilde{a}^{(i_1)}]$ if and only if $a^{(i_1)}$ is a virtual member of the chart of $U[u^{(i_1)}\leadsto a^{(i_1)}]$.
\end{enumerate}

Assume $a^{(i_1)}$ is a virtual member of the chart of $U[u^{(i_1)}\leadsto a^{(i_1)}]$, and $a^{(i_2)}$ is a virtual member of the chart of $U[u^{(i_2)}\leadsto a^{(i_2)}]$. Then
\begin{enumerate}[label={(\arabic*)}]
\setcounter{enumi}{6}
\item
\label{virtual_members_U_incident_replacements7}
$U[u^{(i_1)}\leadsto a^{(i_1)}][\widetilde{u}^{(i_2)} \leadsto \dbtilde{a}^{(i_2)}] = U[u^{(i_2)}\leadsto a^{(i_2)}][\widetilde{u}^{(i_1)} \leadsto \dbtilde{a}^{(i_1)}]$.
\end{enumerate}
\end{corollary}
\begin{proof}
~\paragraph*{(7)} At this moment, we can not yet claim that $\widetilde{u}^{(i_2)}$ and $\dbtilde{a}^{(i_2)}$ are $U[u^{(i_1)}\leadsto a^{(i_1)}]$-incident monomials. In the same way, we can not yet claim that $\widetilde{u}^{(i_1)}$ and $\dbtilde{a}^{(i_1)}$ are $U[u^{(i_2)}\leadsto a^{(i_2)}]$-incident monomials. However, we still can do formally the replacement $\widetilde{u}^{(i_2)} \leadsto \dbtilde{a}^{(i_2)}$ in $U[u^{(i_1)}\leadsto a^{(i_1)}]$ and the replacement $\widetilde{u}^{(i_1)} \leadsto \dbtilde{a}^{(i_1)}$ in $U[u^{(i_2)}\leadsto a^{(i_2)}]$. Then literally in the same way as in statement~\ref{two_replacements_incident2} of Lemma~\ref{two_replacements_incident}, one can show that
\begin{equation*}
U[u^{(i_1)}\leadsto a^{(i_1)}][\widetilde{u}^{(i_2)} \leadsto \dbtilde{a}^{(i_2)}] = U[u^{(i_2)}\leadsto a^{(i_2)}][\widetilde{u}^{(i_1)} \leadsto \dbtilde{a}^{(i_1)}].
\end{equation*}

~\paragraph*{(1)} Since $u^{(i_1)}$ and $a^{(i_1)}$ are $U$-incident monomials, by definition, there exists a sequence of monomials of $m_1^{(i_1)}, \ldots, m^{(i_1)}_{n_1 + 1}$ that satisfies conditions of Definition~\ref{U_incident_monomials}. In particular, $m_1^{(i_1)} = u^{(i_1)}$ and $m^{(i_1)}_{n_1 + 1} = a^{(i_1)}$. We prove statement~\ref{virtual_members_U_incident_replacements1} of Corollary~\ref{virtual_members_U_incident_replacements} by induction on $n_1$.

If $n_1 = 1$, then the statement follows from Corollary~\ref{full_virtual_members_corresponding}.

Assume $n_1 > 1$. Then first we perform the replacement $u^{(i_1)} \leadsto m^{(i_1)}_{n_1}$. Let $v^{(i_2)}$ be the image of $u^{(i_2)}$ in $U[u^{(i_1)} \leadsto m^{(i_1)}_{n_1}]$. Denote the monomial $U[u^{(i_1)} \leadsto m^{(i_1)}_{n_1}]$ by $Z$. Since $m^{(i_1)}_{n_1}$ is a virtual member of the chart of $Z$ and $u^{(i_2)}$ is a virtual member of the chart of $U$, by the induction hypothesis, we obtain that $v^{(i_2)}$ is a virtual member of the chart of $Z$.
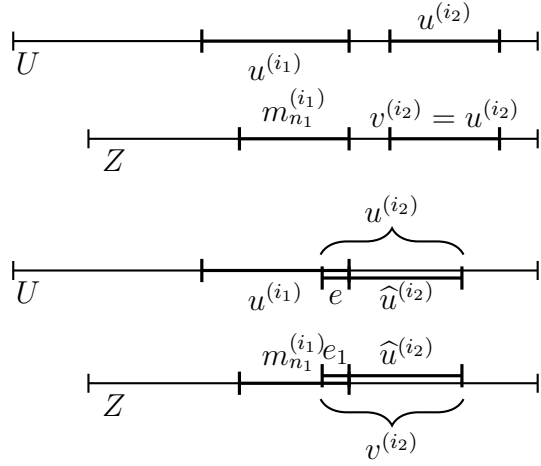
\begin{center}
\begin{tikzpicture}
\draw[|-|, black, thick] (-1.5,0)--(5.5,0);
\node[below, xshift=6] at (-1.5, 0) {$U$};
\draw[|-|, black, very thick] (1,0)--(3,0) node[midway, below] {$u^{(i_1)}$};
\draw[|-|, black, very thick] (3.5,0)--(5,0) node[midway, above] {$u^{(i_2)}$};

\draw[|-|, black, thick] (-0.5, -1.3)--(5.5, -1.3);
\node[below, xshift=10] at (-0.5, -1.3) {$Z$};
\draw[|-|, black, very thick] (1.5, -1.3)--(3, -1.3) node[midway, above] {$m_{n_1}^{(i_1)}$};
\draw[|-|, black, very thick] (3.5, -1.3)--(5, -1.3) node[midway, above] {$v^{(i_2)} = u^{(i_2)}$};
\end{tikzpicture}

\begin{tikzpicture}
\draw[|-|, black, thick] (-1.5, 0)--(5.5, 0);
\node[below, xshift=6] at (-1.5, 0) {$U$};
\draw[|-|, black, very thick] (1,0)--(3,0) node[midway, below] {$u^{(i_1)}$};
\draw[-|, black, very thick] (3,-0.1)--(4.5,-0.1);
\draw[|-|, black, very thick] (2.6,-0.1)--(3,-0.1) node[midway, below] {$e$};
\path (3, 0)--(4.5, 0) node[midway, below] {$\widehat{u}^{(i_2)}$};
\draw [thick, decorate, decoration={brace, amplitude=10pt, raise=6pt}] (2.6, 0) to node[midway, above, yshift=14pt] {$u^{(i_2)}$} (4.5, 0);

\draw[|-|, black, thick] (-0.5, -1.5)--(5.5, -1.5);
\node[below, xshift=10] at (-0.5, -1.5) {$Z$};
\draw[|-|, black, very thick] (1.5, -1.5)--(3,-1.5) node[midway, above] {$m_{n_1}^{(i_1)}$};
\draw[-|, black, very thick] (3, -1.4)--(4.5, -1.4);
\draw[|-|, black, very thick] (2.6, -1.4)--(3, -1.4) node[midway, above] {$e_1$};
\path (3, -1.5)--(4.5, -1.5) node[midway, above, yshift=3] {$\widehat{u}^{(i_2)}$};
\draw [thick, decorate, decoration={brace, amplitude=10pt, raise=6pt, mirror}] (2.6, -1.5) to node[midway, below, yshift=-14pt] {$v^{(i_2)}$} (4.5, -1.5);
\end{tikzpicture}
\end{center}
After that we perform the replacement $m^{(i_1)}_{n_1} \mapsto a^{(i_1)}$ in $Z$. Clearly, the resulting monomial $Z[m^{(i_1)}_{n_1} \mapsto a^{(i_1)}] = U[u^{(i_1)} \leadsto a^{(i_1)}]$, and the image of $v^{(i_2)}$ in $Z[m^{(i_1)}_{n_1} \mapsto a^{(i_1)}]$ is equal to $\widetilde{u}^{(i_2)}$.
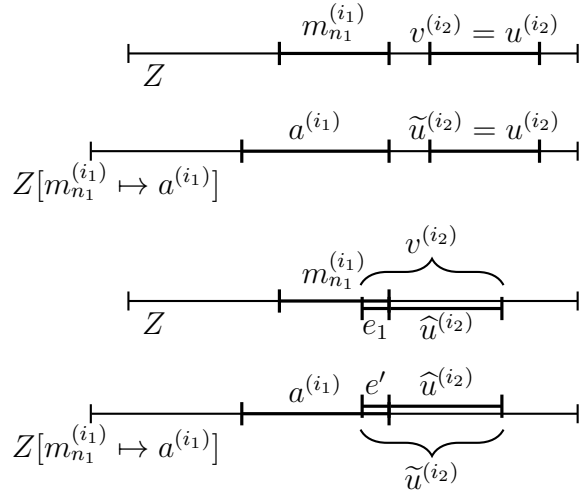
\begin{center}
\begin{tikzpicture}
\draw[|-|, black, thick] (-0.5, 0)--(5.5, 0);
\node[below, xshift=10] at (-0.5, 0) {$Z$};
\draw[|-|, black, very thick] (1.5, 0)--(3, 0) node[midway, above] {$m_{n_1}^{(i_1)}$};
\draw[|-|, black, very thick] (3.5, 0)--(5, 0) node[midway, above] {$v^{(i_2)} = u^{(i_2)}$};

\draw[|-|, black, thick] (-1, -1.3)--(5.5, -1.3);
\node[below, xshift=10] at (-1, -1.3) {$Z[m^{(i_1)}_{n_1} \mapsto a^{(i_1)}]$};
\draw[|-|, black, very thick] (1, -1.3)--(3, -1.3) node[midway, above] {$a^{(i_1)}$};
\draw[|-|, black, very thick] (3.5, -1.3)--(5, -1.3) node[midway, above] {$\widetilde{u}^{(i_2)} = u^{(i_2)}$};
\end{tikzpicture}

\begin{tikzpicture}
\draw[|-|, black, thick] (-0.5, 0)--(5.5, 0);
\node[below, xshift=10] at (-0.5, 0) {$Z$};
\draw[|-|, black, very thick] (1.5, 0)--(3, 0) node[midway, above] {$m_{n_1}^{(i_1)}$};
\draw[-|, black, very thick] (3, -0.1)--(4.5, -0.1);
\draw[|-|, black, very thick] (2.6, -0.1)--(3, -0.1) node[midway, below] {$e_1$};
\path (3, 0)--(4.5, 0) node[midway, below] {$\widehat{u}^{(i_2)}$};
\draw [thick, decorate, decoration={brace, amplitude=10pt, raise=6pt}] (2.6, 0) to node[midway, above, yshift=14pt] {$v^{(i_2)}$} (4.5, 0);

\draw[|-|, black, thick] (-1, -1.5)--(5.5, -1.5);
\node[below, xshift=10] at (-1, -1.5) {$Z[m^{(i_1)}_{n_1} \mapsto a^{(i_1)}]$};
\draw[|-|, black, very thick] (1, -1.5)--(3,-1.5) node[midway, above] {$a^{(i_1)}$};
\draw[-|, black, very thick] (3, -1.4)--(4.5, -1.4);
\draw[|-|, black, very thick] (2.6, -1.4)--(3, -1.4) node[midway, above] {$e^{\prime}$};
\path (3, -1.5)--(4.5, -1.5) node[midway, above, yshift=3] {$\widehat{u}^{(i_2)}$};
\draw [thick, decorate, decoration={brace, amplitude=10pt, raise=6pt, mirror}] (2.6, -1.5) to node[midway, below, yshift=-14pt] {$\widetilde{u}^{(i_2)}$} (4.5, -1.5);
\end{tikzpicture}
\end{center}
Since $a^{(i_1)}$ is a virtual member of the chart of $Z[m^{(i_1)}_{n_1} \mapsto a^{(i_1)}]$ and $v^{(i_2)}$ is a virtual member of the chart of $Z$, by the induction hypothesis, we have that $\widetilde{u}^{(i_2)}$ is a virtual member of the chart $U[u^{(i_1)} \leadsto a^{(i_1)}]$.

~\paragraph*{(4)} Statement~\ref{virtual_members_U_incident_replacements4} of Corollary~\ref{virtual_members_U_incident_replacements} is proved in the same way as statement~\ref{virtual_members_U_incident_replacements1}.

~\paragraph*{(2), (3), (5), (6)} Statements~\ref{virtual_members_U_incident_replacements2},~\ref{virtual_members_U_incident_replacements3},~\ref{virtual_members_U_incident_replacements5},~\ref{virtual_members_U_incident_replacements6} are proved together.

Since $u^{(i_1)}$ and $a^{(i_1)}$ are $U$-incident monomials, by definition, there exists a sequence of monomials of $m_1^{(i_1)}, \ldots, m^{(i_1)}_{n_1 + 1}$ that satisfies the conditions of Definition~\ref{U_incident_monomials}. In particular, $m_1^{(i_1)} = u^{(i_1)}$ and $m^{(i_1)}_{n_1 + 1} = a^{(i_1)}$. Since $u^{(i_2)}$ and $a^{(i_2)}$ are $U$-incident monomials, by definition, there exists a sequence of monomials of $m_1^{(i_2)}, \ldots, m^{(i_2)}_{n_2 + 1}$ that satisfies the conditions of Definition~\ref{U_incident_monomials}. In particular, $m_1^{(i_2)} = u^{(i_2)}$ and $m^{(i_2)}_{n_2 + 1} = a^{(i_2)}$. We prove statements~\ref{virtual_members_U_incident_replacements2},~\ref{virtual_members_U_incident_replacements3},~\ref{virtual_members_U_incident_replacements5},~\ref{virtual_members_U_incident_replacements6} of Corollary~\ref{virtual_members_U_incident_replacements} by induction on $n_1 + n_2$.

We consider only the case when $u^{(i_1)}$ and $u^{(i_2)}$ are not separated. The case when $u^{(i_1)}$ and $u^{(i_2)}$ are separated is considered in a similar (but simpler) way.

If $n_1 = 1$ and $n_2 = 1$, then these statements follow from Lemma~\ref{second_replacement_incident} and Lemma~\ref{virtual_members_replacements_stability}.

Assume $a^{(i_1)}$ is a virtual member of the chart of $U[u^{(i_1)} \leadsto a^{(i_1)}]$. First suppose that $n_2 = 1$, that is, $u^{(i_2)}$ and $a^{(i_2)}$ are incident monomials. Let us show that $\widetilde{u}^{(i_2)}$ and $\dbtilde{a}^{(i_2)}$ are incident monomials and that $a^{(i_2)}$ is a virtual member of the chart of $U[u^{(i_2)} \leadsto a^{(i_2)}]$ if and only if $\dbtilde{a}^{(i_2)}$ is a virtual member of the chart of  $U[u^{(i_1)} \leadsto a^{(i_1)}][\widetilde{u}^{(i_2)} \leadsto \dbtilde{a}^{(i_2)}]$. We prove it by induction on $n_1$.

We consider the replacement $u^{(i_1)} \leadsto a^{(i_1)}$ and split it into two replacements: $u^{(i_1)} \leadsto m_{n_1}^{(i_1)}$ and $m_{n_1}^{(i_1)} \mapsto a^{(i_1)}$. Denote the monomial $U[u^{(i_1)} \leadsto m_{n_1}^{(i_1)}]$ by $W$. Let $v^{(i_2)}$ be the image of $u^{(i_2)}$ in $W$. Then
\begin{equation*}
v^{(i_2)} = e_1\left(e^{-1}\cdot u^{(i_2)}\right), \textit{ where } e_1 \textit{ is the overlap of } m_{n_1}^{(i_1)} \textit{ and } v^{(i_2)}.
\end{equation*}
We put
\begin{equation*}
b^{(i_2)} = e_1 \cdot e^{-1}\cdot a^{(i_2)}.
\end{equation*}
Since $u^{(i_1)}$ and $m_{n_1}^{(i_1)}$ are connected by the sequence of monomials of length $n_1 - 1$ and $m_{n_1}^{(i_1)}$ is a virtual member of the chart of $W$, by the induction hypothesis, we obtain that $v^{(i_2)}$ and $b^{(i_2)}$ are incident monomials and $b^{(i_2)}$ is a virtual member of the chart of $W[v^{(i_2)} \mapsto b^{(i_2)}] = U[u^{(i_1)} \leadsto m_{n_1}^{(i_1)}][v^{(i_2)} \mapsto b^{(i_2)}]$ if and only if $a^{(i_2)}$ is a virtual member of the chart of $U[u^{(i_2)} \leadsto a^{(i_2)}]$.

Now consider the replacement $m_{n_1}^{(i_1)} \mapsto a^{(i_1)}$ in $W$. The resulting monomial is equal to $W[m_{n_1}^{(i_1)} \mapsto a^{(i_1)}] = U[u^{(i_1)} \leadsto a^{(i_1)}]$. So, $a^{(i_1)}$ is a virtual member of the chart of $W[m_{n_1}^{(i_1)} \mapsto a^{(i_1)}]$. Clearly, the image of $v^{(i_2)}$ in $W[m_{n_1}^{(i_1)} \mapsto a^{(i_1)}]$ is equal to
\begin{equation*}
\widetilde{u}^{(i_2)} = e^{\prime}\left(e_1^{-1} \cdot v^{(i_2)}\right).
\end{equation*}
If we do the same transformation of $b^{(i_2)}$, then we obtain
\begin{equation*}
e^{\prime} \cdot e_1^{-1} \cdot b^{(i_2)} = e^{\prime} \cdot e_1^{-1} \cdot e_1 \cdot e^{-1}\cdot a^{(i_2)} = e^{\prime} \cdot e^{-1}\cdot a^{(i_2)} = \dbtilde{a}^{(i_2)}.
\end{equation*}
Since $m_{n_1}^{(i_1)}$ and $a^{(i_1)}$ are connected by the sequence of monomials of length $1$ and $a^{(i_1)}$ is a virtual member of the chart of $W[m_{n_1}^{(i_1)} \mapsto a^{(i_1)}]$, by the induction hypothesis, we obtain that $\widetilde{u}^{(i_2)}$ and $\dbtilde{a}^{(i_2)}$ are incident monomials and $\dbtilde{a}^{(i_2)}$ is a virtual member of the chart of
\begin{equation*}
W[m_{n_1}^{(i_1)} \mapsto a^{(i_1)}][\widetilde{u}^{(i_2)} \mapsto \dbtilde{a}^{(i_2)}] = U[u^{(i_1)} \leadsto a^{(i_1)}][\widetilde{u}^{(i_2)} \mapsto \dbtilde{a}^{(i_2)}]
\end{equation*}
if and only if $b^{(i_2)}$ is a virtual member of the chart of $W[v^{(i_2)} \mapsto b^{(i_2)}]$.

Combining the above results, we obtain that $a^{(i_2)}$ is a virtual member of the chart of $U[u^{(i_2)} \mapsto a^{(i_2)}]$ if and only if $\dbtilde{a}^{(i_2)}$ is a virtual member of the chart of  $U[u^{(i_1)} \leadsto a^{(i_1)}][\widetilde{u}^{(i_2)} \mapsto \dbtilde{a}^{(i_2)}]$. So far, we are done with the case when $a^{(i_1)}$ is a virtual member of the chart of $U[u^{(i_1)} \leadsto a^{(i_1)}]$ and $n_2 = 1$.

Statements~\ref{virtual_members_U_incident_replacements5} and~\ref{virtual_members_U_incident_replacements6} for the case when $n_1 = 1$ are considered in the same way.

Now again let $a^{(i_1)}$ be a virtual member of the chart of $U[u^{(i_1)} \leadsto a^{(i_1)}]$. We assume that $n_2 > 1$. Let us show that $\widetilde{u}^{(i_2)}$ and $\dbtilde{a}^{(i_2)}$ are $U[u^{(i_1)} \leadsto a^{(i_1)}]$-incident monomials and that $a^{(i_2)}$ is a virtual member of the chart of $U[u^{(i_2)} \leadsto a^{(i_2)}]$ if and only if $\dbtilde{a}^{(i_2)}$ is a virtual member of the chart of  $U[u^{(i_1)} \leadsto a^{(i_1)}][\widetilde{u}^{(i_2)} \leadsto \dbtilde{a}^{(i_2)}]$. We prove it by induction on $n_1 + n_2$. That is, we suppose that statements~\ref{virtual_members_U_incident_replacements2},~\ref{virtual_members_U_incident_replacements3},~\ref{virtual_members_U_incident_replacements5},~\ref{virtual_members_U_incident_replacements6} are proved for total length of the corresponding sequences of monomials smaller than $n_1 + n_2$.

We consider the replacement $u^{(i_2)} \leadsto a^{(i_2)}$ in $U$ and split it into two replacements: $u^{(i_2)} \leadsto m_{n_2}^{(i_2)}$ and $m_{n_2}^{(i_2)} \mapsto a^{(i_2)}$. First we perform the replacement $u^{(i_2)} \leadsto m_{n_2}^{(i_2)}$ in $U$. Let $v^{(i_1)}$ be the image of $u^{(i_1)}$ in $U[u^{(i_2)} \leadsto m_{n_2}^{(i_2)}]$. Then
\begin{equation*}
v^{(i_1)} = \left(u^{(i_1)}\cdot e^{-1}\right)e_2, \textit{ where } e_2 \textit{ is the overlap of } v^{(i_1)} \textit{ and } m_{n_2}^{(i_2)}.
\end{equation*}
We put
\begin{equation*}
b^{(i_1)} = a^{(i_1)} \cdot e^{-1}\cdot e_2.
\end{equation*}
Since $u^{(i_2)}$ and $m_{n_2}^{(i_2)}$ are connected by the sequence of monomials of length $n_2 - 1$ and $m_{n_2}^{(i_2)}$ is a virtual member of the chart of $U[u^{(i_2)} \leadsto m_{n_2}^{(i_2)}]$, by the induction hypothesis, we obtain that $v^{(i_1)}$ and $b^{(i_1)}$ are $U[u^{(i_2)} \leadsto m_{n_2}^{(i_2)}]$-incident monomials and $b^{(i_1)}$ is a virtual member of the chart of $U[u^{(i_2)} \leadsto m_{n_2}^{(i_2)}][v^{(i_1)} \leadsto b^{(i_1)}]$.

Denote the monomial $U[u^{(i_2)} \leadsto m_{n_2}^{(i_2)}]$ by $Z$. We have the replacements $m_{n_2}^{(i_2)} \leadsto a^{(i_2)}$ and $v^{(i_1)} \leadsto b^{(i_1)}$ in $Z$. Let $\widetilde{m}_{n_2}^{(i_2)}$ be the image of $m_{n_2}^{(i_2)}$ in $Z[v^{(i_1)} \leadsto b^{(i_1)}] = U[u^{(i_2)} \leadsto m_{n_2}^{(i_2)}][v^{(i_1)} \leadsto b^{(i_1)}]$. Since $e_2$ is the overlap of $v^{(i_1)}$ and $m_{n_2}^{(i_2)}$, we see that
\begin{equation}
\label{m_transformation}
\widetilde{m}_{n_2}^{(i_2)} = e_2^{\prime}\left(e_2^{-1} \cdot m_{n_2}^{(i_2)}\right), \textit{ where } e_2^{\prime} \textit{ is the overlap of } b^{(i_1)} \textit{ and } \widetilde{m}_{n_2}^{(i_2)}.
\end{equation}

Let us show that $e_2^{\prime} \cdot e_2^{-1} = e^{\prime}\cdot e$. Let
\begin{equation*}
\dbtilde{m}_{n_2}^{(i_2)} = e^{\prime}\left(e^{-1}\cdot m_{n_2}^{(i_2)}\right).
\end{equation*}
Since $u^{(i_2)}$ and $m_{n_2}^{(i_2)}$ are $U$-incident monomials and $a^{(i_1)}$ is a virtual member of the chart of $U[u^{(i_1)} \leadsto a^{(i_1)}]$, by the induction hypothesis, we obtain that $\widetilde{u}^{(i_2)}$ and $\dbtilde{m}_{n_2}^{(i_2)}$ are $U[u^{(i_1)} \leadsto a^{(i_1)}]$-incident monomials. It follows from the induction hypothesis and statement~\ref{virtual_members_U_incident_replacements7} that
\begin{equation*}
U[u^{(i_2)} \leadsto m_{n_2}^{(i_2)}][v^{(i_1)} \leadsto b^{(i_1)}] = U[u^{(i_1)} \leadsto a^{(i_1)}][\widetilde{u}^{(i_2)} \leadsto \dbtilde{m}_{n_2}^{(i_2)}].
\end{equation*}
Denote these monomials by $Q$. Assume $u^{(i_2)}$ is the $j$-th virtual member of the chart of $U$. Then, since we replace only virtual members of the chart by virtual members of the chart in order to obtain the right-hand and the left-hand monomial in the above equality, by the induction hypothesis, we see that $\dbtilde{m}_{n_2}^{(i_2)}$ is the $j$-th virtual member of the chart of $Q$ and $\widetilde{m}_{n_2}^{(i_2)}$ is the $j$-th virtual member of the chart of $Q$. Therefore, $\widetilde{m}_{n_2}^{(i_2)} = \dbtilde{m}_{n_2}^{(i_2)}$. Hence, $e^{\prime}\cdot e = e_2^{\prime} \cdot e_2^{-1}$.

So, if we do the same transformation of $a^{(i_2)}$ as we do in~\eqref{m_transformation} over $m_{n_2}^{(i_2)}$, then we obtain $\dbtilde{a}^{(i_2)}$ as a result. Since $m_{n_2}^{(i_2)}$ and $a^{(i_2)}$ are incident monomials, by the basis of induction (case $n_2 = 1$) we obtain that $\widetilde{m}_{n_2}^{(i_2)}$ and $\dbtilde{a}^{(i_2)}$ are incident monomials. So, we see that $\widetilde{u}^{(i_2)} \leadsto  \widetilde{m}_{n_2}^{(i_2)}$ and $\widetilde{m}_{n_2}^{(i_2)} \mapsto \dbtilde{a}^{(i_2)}$ is a sequence of replacements in one position of the chart of $U[u^{(i_1)} \leadsto a^{(i_1)}]$. Since $\widetilde{m}_{n_2}^{(i_2)}$ is a virtual member of the chart of $U[u^{(i_1)} \leadsto a^{(i_1)}][\widetilde{u}^{(i_2)} \leadsto \widetilde{m}_{n_2}^{(i_2)}]$, we obtain that $\widetilde{u}^{(i_2)}$ and $\dbtilde{a}^{(i_2)}$ are $U[u^{(i_1)} \leadsto a^{(i_1)}]$-incident monomials.

By the basis of induction (case $n_2 = 1$), we obtain that $a^{(i_2)}$ is a virtual member of the chart of $Z[m_{n_2}^{(i_2)} \mapsto a^{(i_2)}]$ if and only if $\dbtilde{a}^{(i_2)}$ is a virtual member of the chart of $Z[v^{(i_1)} \leadsto b^{(i_1)}][\widetilde{m}_{n_2}^{(i_2)} \mapsto \dbtilde{a}^{(i_2)}]$. We proved above that
\begin{equation*}
Z[v^{(i_1)} \leadsto b^{(i_1)}] = U[u^{(i_2)} \leadsto m_{n_2}^{(i_2)}][v^{(i_1)} \leadsto b^{(i_1)}] = U[u^{(i_1)} \leadsto a^{(i_1)}][\widetilde{u}^{(i_2)} \leadsto \widetilde{m}_{n_2}^{(i_2)}].
\end{equation*}
Therefore,
\begin{equation*}
Z[v^{(i_1)} \leadsto b^{(i_1)}][\widetilde{m}_{n_2}^{(i_2)} \mapsto \dbtilde{a}^{(i_2)}] = U[u^{(i_1)} \leadsto a^{(i_1)}][\widetilde{u}^{(i_2)} \leadsto \dbtilde{a}^{(i_2)}].
\end{equation*}
Clearly, $Z[m_{n_2}^{(i_2)} \mapsto a^{(i_2)}] = U[u^{(i_2)} \leadsto a^{(i_2)}]$. This completes the proof of statements~\ref{virtual_members_U_incident_replacements2} and~\ref{virtual_members_U_incident_replacements3}.

Statements~\ref{virtual_members_U_incident_replacements5} and~\ref{virtual_members_U_incident_replacements6} for the case when $n_1 > 1$ are considered in the same way.
\end{proof}

Now we go to the principle case which serves as a basis for the argument in Subsection~\ref{tensor_product_section}.
\begin{definition}
\label{consecutive_replacements_def}
Let $U$ be a monomial. Let $u^{(1)}, \ldots, u^{(m)}$ be all the different virtual members of the chart of $U$ enumerated from left to right. Assume $\lbrace i_1, \ldots, i_k\rbrace \subseteq \lbrace 1, \ldots, m\rbrace$. Let $u^{(i_s)}$ and $a^{(i_s)}$ be $U$-incident monomials, $s = 1, \ldots, k$ (we allow that some $u^{(i_s)} = a^{(i_s)}$). Consider the replacements $u^{(i_s)} \leadsto a^{(i_s)}$ in $U$, $s = 1, \ldots, k$. Assume $a^{(i_s)}$ is a virtual member of the chart of $U[u^{(i_s)} \leadsto a^{(i_s)}]$ for all $s = 1, \ldots, k$. Notice that since $u^{(i_s)}$ is the $i_s$-th virtual member of the chart of $U$, we obtain that $a^{(i_s)}$ is the $i_s$-th virtual member of the chart of $U[u^{(i_s)} \leadsto a^{(i_s)}]$.

Let $j_1, \ldots, j_k$ be a permutation of the set $\lbrace i_1, \ldots, i_k\rbrace$. We define a sequence of replacements that we call \emph{consecutive performing  the replacements $u^{(i_s)} \leadsto a^{(i_s)}$, $s = 1, \ldots, k$, starting from $U$, in order $j_1,\ldots, j_k$}. We define it by induction on number $k$ of replacements.

Assume $k = 1$. Then $U \leadsto U[u^{(j_1)} \leadsto a^{(j_1)}]$ is the required sequence.

Assume $k > 1$. Then first we perform the replacement $u^{(j_1)} \leadsto a^{(j_1)}$ in $U$. Let $\widetilde{u}^{(j)}$ be the image of $u^{(j)}$ in $U[u^{(j_1)} \leadsto a^{(j_1)}]$, $j = j_2, \ldots, j_k$, defined by formulas~\eqref{image_of_right_U_incident} and~\eqref{image_of_left_U_incident}. Let $\widetilde{u}^{(j)} \leadsto \dbtilde{a}^{(j)}$ be the replacement in $U[u^{(j_1)} \leadsto a^{(j_1)}]$ that corresponds to the replacement $u^{(j)} \leadsto a^{(j)}$ in $U$, $j = j_2, \ldots, j_k$. Every element $\dbtilde{a}^{(j)}$, $j = j_2, \ldots, j_k$, is defined by formula~\eqref{corresponding_element_U_incident_right} or~\eqref{corresponding_element_U_incident_left}. Actually, at most two replacements may essentially change.

The following properties needed hold because $a^{(j_1)}$ is a virtual member of the chart of $U[u^{(j_1)} \leadsto a^{(j_1)}]$. By Corollary~\ref{virtual_members_U_incident_replacements}, we see that $\widetilde{u}^{(j)}$ is a virtual member of the chart of $U[u^{(j_1)} \leadsto a^{(j_1)}]$, $\widetilde{u}^{(j)}$ and $\dbtilde{a}^{(j)}$ are $U[u^{(j_1)} \leadsto a^{(j_1)}]$-incident monomials, and $\dbtilde{a}^{(j)}$ is a virtual member of the chart of $U[u^{(j_1)} \leadsto a^{(j_1)}][\widetilde{u}^{(j)} \leadsto \dbtilde{a}^{(j)}]$, $j = j_2, \ldots, j_k$. So, consecutive performing  the replacements $\widetilde{u}^{(j)} \leadsto \dbtilde{a}^{(j)}$, $j = j_2, \ldots, j_k$, starting from $U[u^{(j_1)} \leadsto a^{(j_1)}]$ in order $j_2, \ldots, j_k$ is defined by the induction hypothesis. We add the replacement $u^{(j_1)} \leadsto a^{(j_1)}$ in $U$ to the beginning of this sequence and obtain the desired sequence of replacements.
\end{definition}

\begin{definition}
\label{middle_subwords_chart}
Let $U$ be a monomial and $a \in \longmo{U}$. Similarly to subwords $t(a)$, $i(a)$, and $m(a)$, which are defined in Subsection~\ref{admissible_replacements_section} with the use of overlaps with neighbours from $\longmo{U}$, we define subwords $t_v(a)$, $i_v(a)$, and $m_v(a)$ with the use of overlaps with neighbouring virtual members of the chart. Namely,
\begin{equation*}
a = et_v(a),
\end{equation*}
where $e$ is the overlap of $a$ and the left neighbouring virtual member of the chart of $U$, $e$ is empty, if there are no virtual members of the chart of $U$ that start from the left of the beginning of $a$, or $a$ and the left neighbouring virtual member of the chart of $U$ are separated or touch at a point;
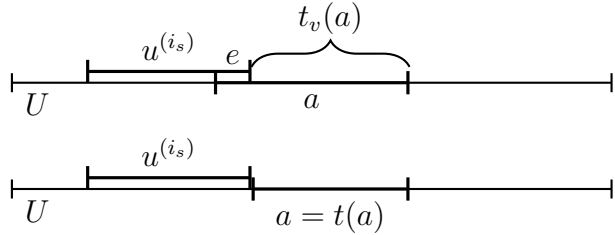
\begin{center}
\begin{tikzpicture}
\node at (0, 0) [below, xshift=10] {$U$};
\draw[|-|, black, thick] (0,0)--(8,0);
\draw[|-|, black, very thick] (1, 0.15)--(3.2, 0.15) node[midway, above] {$u^{(i_s)}$};
\draw[|-|, black, very thick] (2.7, 0)--(5.3, 0) node[midway, below] {$a$};
\draw [thick, decorate, decoration={brace, amplitude=10pt, raise=6pt}] (3.2, 0) to node[midway, above, yshift=14pt] {$t_v(a)$} (5.3, 0);
\path (2.7, 0)--(3.2, 0) node[midway, above, yshift=3] {$e$};
\end{tikzpicture}

\begin{tikzpicture}
\node at (0, 0) [below, xshift=10] {$U$};
\draw[|-|, black, thick] (0,0)--(8,0);
\draw[|-|, black, very thick] (1, 0.15)--(3.2, 0.15) node[midway, above] {$u^{(i_s)}$};
\draw[|-|, black, very thick] (3.2, 0)--(5.3, 0) node[midway, below] {$a = t(a)$};
\end{tikzpicture}
\end{center}
\begin{equation*}
a = i_v(a)f,
\end{equation*}
where $f$ is the overlap of $a$ and the right neighbouring virtual member of the chart of $U$, $f$ is empty, if there are no virtual members of the chart of $U$ that start from the right of the beginning of $a$, or $a$ and the right neighbouring virtual member of the chart of $U$ are separated or touch at a point;
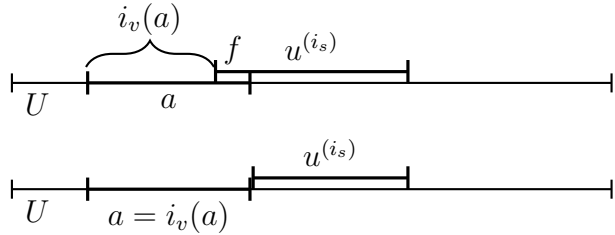
\begin{center}
\begin{tikzpicture}
\node at (0, 0) [below, xshift=10] {$U$};
\draw[|-|, black, thick] (0,0)--(8,0);
\draw[|-|, black, very thick] (1, 0)--(3.2, 0) node[midway, below] {$a$};
\draw[|-|, black, very thick] (2.7, 0.15)--(5.3, 0.15) node[midway, above] {$u^{(i_s)}$};
\draw [thick, decorate, decoration={brace, amplitude=10pt, raise=6pt}] (1, 0) to node[midway, above, yshift=14pt] {$i_v(a)$} (2.7, 0);
\path (2.7, 0)--(3.2, 0) node[midway, above, yshift=3] {$f$};
\end{tikzpicture}

\begin{tikzpicture}
\node at (0, 0) [below, xshift=10] {$U$};
\draw[|-|, black, thick] (0,0)--(8,0);
\draw[|-|, black, very thick] (1, 0)--(3.2, 0) node[midway, below] {$a = i_v(a)$};
\draw[|-|, black, very thick] (3.2, 0.15)--(5.3, 0.15) node[midway, above] {$u^{(i_s)}$};
\end{tikzpicture}
\end{center}
and finally
\begin{equation*}
a = em_v(a)f,
\end{equation*}
\begin{center}
\begin{tikzpicture}
\node[text width=11cm, align=left] at (1.6, -1) {\small\baselineskip=10pt either $i_s = i_r + 1$ if $a$ is not a virtual member of the chart of $U$, or $i_s = i_r + 2$ if $a$ is a virtual member of the chart of $U$.};
\end{tikzpicture}

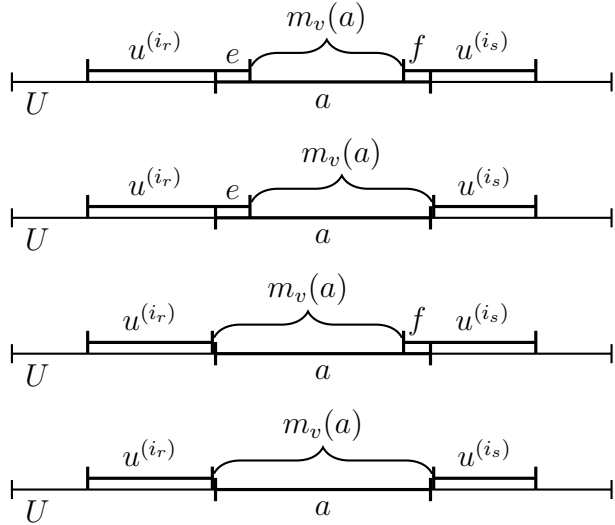
\begin{tikzpicture}
\node at (0, 0) [below, xshift=10] {$U$};
\draw[|-|, black, thick] (0,0)--(8,0);
\draw[|-|, black, very thick] (1, 0.15)--(3.2, 0.15) node[midway, above, xshift=-5] {$u^{(i_r)}$};
\draw[|-|, black, very thick] (2.7, 0)--(5.6, 0) node[midway, below] {$a$};
\draw[|-|, black, very thick] (5.2, 0.15)--(7, 0.15) node[midway, above, xshift=5] {$u^{(i_s)}$};
\draw [thick, decorate, decoration={brace, amplitude=10pt, raise=6pt}] (3.2, 0) to node[midway, above, yshift=14pt] {$m_v(a)$} (5.2, 0);
\path (2.7, 0)--(3.2, 0) node[midway, above, yshift=3] {$e$};
\path (5.2, 0)--(5.6, 0) node[midway, above, yshift=3] {$f$};
\end{tikzpicture}

\begin{tikzpicture}
\node at (0, 0) [below, xshift=10] {$U$};
\draw[|-|, black, thick] (0,0)--(8,0);
\draw[|-|, black, very thick] (1, 0.15)--(3.2, 0.15) node[midway, above, xshift=-5] {$u^{(i_r)}$};
\draw[|-|, black, very thick] (2.7, 0)--(5.6, 0) node[midway, below] {$a$};
\draw[|-|, black, very thick] (5.6, 0.15)--(7, 0.15) node[midway, above] {$u^{(i_s)}$};
\draw [thick, decorate, decoration={brace, amplitude=10pt, raise=6pt}] (3.2, 0) to node[midway, above, yshift=14pt] {$m_v(a)$} (5.6, 0);
\path (2.7, 0)--(3.2, 0) node[midway, above, yshift=3] {$e$};
\end{tikzpicture}

\begin{tikzpicture}
\node at (0, 0) [below, xshift=10] {$U$};
\draw[|-|, black, thick] (0,0)--(8,0);
\draw[|-|, black, very thick] (1, 0.15)--(2.7, 0.15) node[midway, above] {$u^{(i_r)}$};
\draw[|-|, black, very thick] (2.7, 0)--(5.6, 0) node[midway, below] {$a$};
\draw[|-|, black, very thick] (5.2, 0.15)--(7, 0.15) node[midway, above, xshift=5] {$u^{(i_s)}$};
\draw [thick, decorate, decoration={brace, amplitude=10pt, raise=6pt}] (2.7, 0) to node[midway, above, yshift=14pt] {$m_v(a)$} (5.2, 0);
\path (5.2, 0)--(5.6, 0) node[midway, above, yshift=3] {$f$};
\end{tikzpicture}

\begin{tikzpicture}
\node at (0, 0) [below, xshift=10] {$U$};
\draw[|-|, black, thick] (0,0)--(8,0);
\draw[|-|, black, very thick] (1, 0.15)--(2.7, 0.15) node[midway, above] {$u^{(i_r)}$};
\draw[|-|, black, very thick] (2.7, 0)--(5.6, 0) node[midway, below] {$a$};
\draw[|-|, black, very thick] (5.6, 0.15)--(7, 0.15) node[midway, above] {$u^{(i_s)}$};
\draw [thick, decorate, decoration={brace, amplitude=10pt, raise=6pt}] (2.7, 0) to node[midway, above, yshift=14pt] {$m_v(a)$} (5.6, 0);
\end{tikzpicture}
\end{center}
\end{definition}

First let us notice obvious properties of consecutive replacements.
\begin{remark}
\label{consecurive_repl_basic_properties}
Assume that we are under the conditions of Definition~\ref{consecutive_replacements_def}. When we consecutively perform replacements, at every step we have replacements that are already performed and the rest of the replacements that we need to perform. Assume $Z$ is the resulting monomial after some step. Let us emphasise that for every initial replacement in $U$ from the second set there exists a corresponding replacement in $Z$. Consider the corresponding replacements in more detail.

Assume that we consecutively performed the first $t - 1$ replacements $u^{(j_1)} \leadsto a^{(j_1)}, \ldots, u^{(j_{t - 1})} \leadsto a^{(j_{t - 1})}$, starting from $U$. Let $Z$ be the resulting monomial. Assume $\widetilde{u}^{(j_t)}$ is the $j_t$-th virtual member of the chart of $Z$. Let $\widetilde{u}^{(j_t)} \leadsto \dbtilde{a}^{(j_t)}$ be the next replacement in $Z$ in the sequence of consecutive replacements (the element $\dbtilde{a}^{(j_t)}$ corresponds to the element $a^{(j_t)}$). Then the following properties hold.
\begin{enumerate}
\item
It follows from Corollary~\ref{full_virtual_members_corresponding} that the structure of the chart is preserved completely after every replacement of a virtual member of the chart by a virtual member of the chart. That is, since $u^{(j_t)}$ is the $j_t$-th virtual member of the chart of $U$, we obtain that $\widetilde{u}^{(j_t)}$ is the $j_t$-th virtual member of the chart of $Z$. Since $a^{(j_t)}$ is a $j_t$-th virtual member of the chart of $U[u^{(j_t)} \leadsto a^{(j_t)}]$, we see that $\dbtilde{a}^{(j_t)}$ is the $j_t$-th virtual member of the chart of $Z[\widetilde{u}^{(j_t)} \leadsto \dbtilde{a}^{(j_t)}]$. So, index numbers of virtual members of the chart do not shift in a resulting monomial after every step.
\item
By definition, every virtual member of the chart is of $\SPM$-measure $\geqslant \tau - 2$. So, every  replacement above is, in particular, a replacement of a maximal occurrence of $\SPM$-measure $\geqslant \tau - 2$ by an element of $\Mon$ of $\SPM$-measure $\geqslant \tau - 2$. Therefore, it follows from the results of Section~\ref{mt_configurations} that $u^{(j_t)}$ in $U$ may differ from $\widetilde{u}^{(j_t)}$ in $Z$ only by  small pieces at the beginning and at the end.

In more detail, we have the following. Since we replace only virtual members of the chart by virtual members of the chart, we obtain
\begin{align*}
&u^{(j_t)} = em_v(u^{(j_t)})f,\ \widetilde{u}^{(j_t)} = \widetilde{e}m_v(u^{(j_t)})\widetilde{f},\\
&\widetilde{u}^{(j_t)} = \widetilde{e}\left(e^{-1}\cdot u^{(j_t)}\cdot f^{-1}\right)\widetilde{f},\ \dbtilde{a}^{(j_t)} = \widetilde{e}\cdot e^{-1}\cdot a^{(j_t)}\cdot f^{-1}\cdot\widetilde{f},
\end{align*}
where $e$ is an overlap of $u^{(j_t)}$ and its left neighbouring virtual member of the chart of $U$, and $\widetilde{e}$ is an overlap of $\widetilde{u}^{(j_t)}$ and its left neighbouring virtual member of the chart of $Z$, $f$ is an overlap of $u^{(j_t)}$ and its right neighbouring virtual member of the chart of $U$, and $\widetilde{f}$ is an overlap of $\widetilde{u}^{(j_t)}$ and its right neighbouring virtual member of the chart of $Z$.

Moreover, if all the replacements $u^{(j_1)} \leadsto a^{(j_1)}, \ldots, u^{(j_{t - 1})} \leadsto a^{(j_{t - 1})}$ are from the left of $u^{(j_t)}$ and its images (that is, $j_1, \ldots, j_{t - 1} < j_t$), then
\begin{align*}
&u^{(j_t)} = et_v(u^{(j_t)}),\ \widetilde{u}^{(j_t)} = \widetilde{e}t_v(u^{(j_t)}),\\
&\widetilde{u}^{(j_t)} = \widetilde{e}\left(e^{-1}\cdot u^{(j_t)}\right),\ \dbtilde{a}^{(j_t)} = \widetilde{e}\cdot e^{-1}\cdot a^{(j_t)}.
\end{align*}
If all the replacements $u^{(j_1)} \leadsto a^{(j_1)}, \ldots, u^{(j_{t - 1})} \leadsto a^{(j_{t - 1})}$ are from the right of $u^{(j_t)}$ and its images (that is, $j_1, \ldots, j_{t - 1} > j_t$), then
\begin{align*}
&u^{(j_t)} = i_v(u^{(j_t)})f,\ \widetilde{u}^{(j_t)} = i_v(u^{(j_t)})\widetilde{f},\\
&\widetilde{u}^{(j_t)} = \left(u^{(j_t)}\cdot f^{-1}\right)\widetilde{f},\ \dbtilde{a}^{(j_t)} = a^{(j_t)}\cdot f^{-1}\cdot \widetilde{f}.
\end{align*}

In what follows we say that \emph{the replacement $\widetilde{u}^{(j_t)} \leadsto \dbtilde{a}^{(j_t)}$ in $Z$ corresponds to the replacement $u^{(j_t)} \leadsto a^{(j_t)}$ in $U$.}
\end{enumerate}
\end{remark}

\begin{lemma}
\label{replace_consecutively}
Let $U$ be a monomial. Let $u^{(1)}, \ldots, u^{(m)}$ be all the different virtual members of the chart of $U$ enumerated from left to right. Assume $\lbrace i_1, \ldots, i_k\rbrace \subseteq \lbrace 1, \ldots, m\rbrace$. Let $u^{(i_s)}$ and $a^{(i_s)}$ be $U$-incident monomials, $s = 1, \ldots, k$ (we allow that some $u^{(i_s)} = a^{(i_s)}$). We consider the replacements $u^{(i_s)} \leadsto a^{(i_s)}$ in $U$, $s = 1, \ldots, k$. Assume $a^{(i_s)}$ is a virtual member of the chart of $U[u^{(i_s)} \leadsto a^{(i_s)}]$ for all $s = 1, \ldots, k$. We consecutively perform the replacements $u^{(i_s)} \leadsto a^{(i_s)}$, $s = 1, \ldots, k$, starting from $U$, in some order. Then the resulting monomial does not depend on the order of performing  the replacements. Moreover, the resulting monomial is a derived monomial of $U$ with the same $f$-characteristic as $U$.
\end{lemma}
\begin{proof}
We prove Lemma~\ref{replace_consecutively} by induction on $k$. Assume $k = 2$. Let $\widetilde{u}^{(j_2)}$ be defined by~\eqref{image_of_right_U_incident}, $\dbtilde{a}^{(j_2)}$ be defined by~\eqref{corresponding_element_U_incident_right}. Let $\widetilde{u}^{(j_1)}$ be defined by~\eqref{image_of_left_U_incident}, $\dbtilde{a}^{(j_1)}$ be defined by~\eqref{corresponding_element_U_incident_left}. It follows from statement~\ref{virtual_members_U_incident_replacements7} of Corollary~\ref{virtual_members_U_incident_replacements} that $U[u^{(j_1)} \leadsto a^{(j_1)}][\widetilde{u}^{(j_2)} \leadsto \dbtilde{a}^{(j_2)}] = U[u^{(j_2)} \leadsto a^{(j_2)}][\widetilde{u}^{(j_1)} \leadsto \dbtilde{a}^{(j_1)}]$. So, the basis of induction is proved.

Assume $k > 2$. Let $j_1, \ldots, j_k$ and $h_1, \ldots, h_k$ be two orders of performing  the replacements. First assume that $j_1 = h_1$. Then $u^{(j_1)} = u^{(h_1)}$, $a^{(j_1)} = a^{(h_1)}$, $U[u^{(j_1)} \leadsto a^{(j_1)}] = U[u^{(h_1)} \leadsto a^{(h_1)}]$ and the set of the rest of the replacements is the same for the both monomials. By the induction hypothesis, we can perform the rest  of the replacements in $U[u^{(j_1)} \leadsto a^{(j_1)}] = U[u^{(h_1)} \leadsto a^{(h_1)}]$ in any order and obtain the same resulting monomial. Hence, we obtain the same resulting monomials when we perform the replacements in order $j_1, \ldots, j_k$ and in order $h_1, \ldots, h_k$.

Assume $j_1 \neq h_1$. Let us perform the first replacement $u^{(j_1)} \leadsto a^{(j_1)}$ in $U$, and the first replacement $u^{(h_1)} \leadsto a^{(h_1)}$ in $U$. By the induction hypothesis, the rest of the replacements in $U[u^{(j_1)} \leadsto a^{(j_1)}]$ and the rest of the replacements in $U[u^{(h_1)} \leadsto a^{(h_1)}]$ can be performed in any order. So, we can continue with the replacement in any position of $j_2, \ldots, j_k$ in $U[u^{(j_1)} \leadsto a^{(j_1)}]$, and with the replacement in any position of $h_2, \ldots, h_k$ in $U[u^{(h_1)} \leadsto a^{(h_1)}]$. Notice that $h_1 \in \lbrace j_2, \ldots, j_k \rbrace$ and $j_1 \in \lbrace h_2, \ldots, h_k \rbrace$. Assume $h_1 = j_{r_0}$ and $j_1 = h_{s_0}$. Let us start from the replacement in the position $h_1$ in $U[u^{(j_1)} \leadsto a^{(j_1)}]$, and from the replacement in the position $j_1$ in $U[u^{(h_1)} \leadsto a^{(h_1)}]$. By the induction hypothesis, we obtain that the result of performing  the replacements, starting from $U[u^{(j_1)} \leadsto a^{(j_1)}]$, in order $j_2, \ldots, j_k$ is the same as in order $j_2, h_1 = j_{r_0}, j_2, \ldots, j_{r_0 - 1}, j_{r_0 + 1}, \ldots, j_k$. Hence, the result of performing  the replacements, starting from $U$, in order $j_1, \ldots, j_k$ is the same as in order $j_1, h_1 = j_{r_0}, j_2, \ldots, j_{r_0 - 1}, j_{r_0 + 1}, \ldots, j_k$. In the same way, we obtain that the result of performing  the replacements, starting from $U[u^{(h_1)} \leadsto a^{(h_1)}]$, in order $h_2, \ldots, h_k$ is the same as in order $h_2, j_1 = h_{s_0}, h_2, \ldots, h_{s_0 - 1}, h_{s_0 + 1}, \ldots, h_k$. Therefore, the result of performing  the replacements, starting from $U$, in order $h_1, \ldots, h_k$ is the same as in order $h_1, j_1 = h_{s_0}, h_2, \ldots, h_{s_0 - 1}, h_{s_0 + 1}, \ldots, h_k$.

Let $\widetilde{u}^{(h_1)} \leadsto \dbtilde{a}^{(h_1)}$ be the replacement in $U[u^{(j_1)} \leadsto a^{(j_1)}]$ that corresponds to the replacement $u^{(h_1)} \leadsto a^{(h_1)}$ in $U$, and let $\widetilde{u}^{(j_1)} \leadsto \dbtilde{a}^{(j_1)}$ be the replacement in $U[u^{(h_1)} \leadsto a^{(h_1)}]$ that corresponds to the replacement $u^{(j_1)} \leadsto a^{(j_1)}$ in $U$ ($\widetilde{u}^{(j_2)}$ is defined by~\eqref{image_of_right_U_incident}, $\dbtilde{a}^{(j_2)}$ is defined by~\eqref{corresponding_element_U_incident_right}, $\widetilde{u}^{(j_1)}$ is defined by~\eqref{image_of_left_U_incident}, $\dbtilde{a}^{(j_1)}$ is defined by~\eqref{corresponding_element_U_incident_left}). It follows from statement~\ref{virtual_members_U_incident_replacements7} of Corollary~\ref{virtual_members_U_incident_replacements} that
\begin{equation*}
U[u^{(j_1)} \leadsto a^{(j_1)}][\widetilde{u}^{(h_1)} \leadsto \dbtilde{a}^{(h_1)}] = U[u^{(h_1)} \leadsto a^{(h_1)}][\widetilde{u}^{(j_1)} \leadsto \dbtilde{a}^{(j_1)}].
\end{equation*}
Denote this monomial by $Z$. Obviously,
\begin{equation*}
\lbrace j_2, \ldots, j_{r_0 - 1}, j_{r_0 + 1}, \ldots, j_k \rbrace = \lbrace h_2, \ldots, h_{s_0 - 1}, h_{s_0 + 1}, \ldots, h_k\rbrace.
\end{equation*}
One can show that the set of the rest of the replacements in $Z$ is the same after the replacements $u^{(j_1)} \leadsto a^{(j_1)}$ in $U$ and $\widetilde{u}^{(h_1)} \leadsto \dbtilde{a}^{(h_1)}$ in $U[u^{(j_1)} \leadsto a^{(j_1)}]$, and after the replacements $u^{(h_1)} \leadsto a^{(h_1)}$ in $U$ and $\widetilde{u}^{(j_1)} \leadsto \dbtilde{a}^{(j_1)}$ in $U[u^{(h_1)} \leadsto a^{(h_1)}]$. Therefore, by the induction hypothesis, the resulting monomial of performing  these replacements, starting from $Z$, is the same in order $j_2, \ldots, j_{r_0 - 1}, j_{r_0 + 1}, \ldots, j_k$ and in order $h_2, \ldots, h_{s_0 - 1}, h_{s_0 + 1}, \ldots, h_k$. Since $Z$ is the resulting monomial of the consecutive replacements $u^{(j_1)} \leadsto a^{(j_1)}$ and $u^{(h_1)} \leadsto a^{(h_1)}$ in $U$ in order $j_1, h_1 =j_{r_0}$ and in order $h_1, j_1 = h_{s_0}$, we see that the result of performing  the initial replacements, starting from $U$, in order $j_1, h_1 = j_{r_0}, j_2, \ldots, j_{r_0 - 1}, j_{r_0 + 1}, \ldots, j_k$ is the same as in order $h_1, j_1 = h_{s_0}, h_2, \ldots, h_{s_0 - 1}, h_{s_0 + 1}, \ldots, h_k$.

We proved above that the result of performing  the initial replacements, starting from $U$, in order $j_1, h_1 = j_{r_0}, j_2, \ldots, j_{r_0 - 1}, j_{r_0 + 1}, \ldots, j_k$ is the same as in order $j_1, \ldots, j_k$, and that the result of performing the initial replacements, starting from $U$, in order $h_1, j_1 = h_{s_0}, h_2, \ldots, h_{s_0 - 1}, h_{s_0 + 1}, \ldots, h_k$ is the same as in order $h_1, \ldots, h_k$. Combining this with the above result, we obtain that the results of the initial replacements, starting from $U$, in orders $j_1, \ldots, j_k$ and $h_1, \ldots, h_k$ are equal.

By Definition~\ref{consecutive_replacements_def}, we replace only virtual members of the chart of monomials. Hence, the resulting monomial is a derived monomial of $U$. Moreover, by Definition~\ref{consecutive_replacements_def}, we perform only replacements of virtual members of the chart by virtual members of the chart. Such replacements preserve $f$-characteristic of monomials. Thus, the resulting monomial is of the same $f$-characteristic as $U$.
\end{proof}

So, different order of consecutive replacements gives the same result. In what follows when we speak about consecutive replacements, we omit in what order we perform them.

The next lemma follows from the fact that consecutive replacements can be performed in any order.
\begin{lemma}
\label{same_small_changes}
Let $U$ be a monomial. Let $u^{(1)}, \ldots, u^{(m)}$ be all the different virtual members of the chart of $U$ enumerated from left to right. Assume $\lbrace i_1, \ldots, i_k\rbrace \subseteq \lbrace 1, \ldots, m\rbrace$. Let $u^{(i_s)}$ and $a^{(i_s)}$ be $U$-incident monomials, $s = 1, \ldots, k$ (we allow that some $u^{(i_s)} = a^{(i_s)}$). We consider the replacements $u^{(i_s)} \leadsto a^{(i_s)}$ in $U$, $s = 1, \ldots, k$. Assume $a^{(i_s)}$ is a virtual member of the chart of $U[u^{(i_s)} \leadsto a^{(i_s)}]$ for all $s = 1, \ldots, k$. Let $i_0$ be one of indices $i_1, \ldots, i_s$. We consecutively perform the replacements $u^{(i_s)} \leadsto a^{(i_s)}$, $i_s \neq i_0$, starting from $U$. Let $Z$ be the resulting monomial. Let $\widetilde{u}^{(i_0)}$ be the $i_0$-th virtual member of the chart of $Z$. According to Remark~\ref{consecurive_repl_basic_properties}, we have
\begin{equation*}
u^{(i_0)} = em_v(u^{(i_0)})f,\ \widetilde{u}^{(i_0)} = \widetilde{e}m_v(u^{(i_0)})\widetilde{f},
\end{equation*}
where $e,\widetilde{e}, f, \widetilde{f}$ are small pieces. Let $\widetilde{u}^{(i_0)} \leadsto \dbtilde{a}^{(i_0)}$ be the replacement in $Z$ that corresponds to the replacement $u^{(i_0)} \leadsto a^{(i_0)}$ in $U$. Then
\begin{equation*}
a^{(i_0)} = e^{\prime}m_v(a^{(i_0)})f^{\prime},\ \dbtilde{a}^{(i_0)} = \widetilde{e}^{\prime}m_v(a^{(i_0)})\widetilde{f}^{\prime},
\end{equation*}
where $e^{\prime}, f^{\prime}, \widetilde{e}^{\prime}, \widetilde{f}^{\prime}$ are small pieces, and
\begin{equation*}
\widetilde{e}\cdot e^{-1} = \widetilde{e}^{\prime} \cdot {e^{\prime}}^{-1},\ f^{-1}\cdot \widetilde{f} = {f^{\prime}}^{-1}\cdot \widetilde{f}^{\prime}.
\end{equation*}
\end{lemma}
\begin{proof}
Let $W$ be the resulting monomial of consecutive performing all the replacements $u^{(i_s)} \leadsto a^{(i_s)}$, $s = 1, \ldots, k$. Then $\dbtilde{a}^{(i_0)}$ is the $i_0$-th virtual member of the chart of $W$. Let us consecutively perform the replacements $u^{(i_s)} \leadsto a^{(i_s)}$, $s = 1, \ldots, k$, starting from the replacement $u^{(i_0)} \leadsto a^{(i_0)}$. Then we see that $W$ is obtained from $U[u^{(i_0)} \leadsto a^{(i_0)}]$ as a result of replacements of virtual members of the chart by virtual members of the chart in positions different from $i_0$. By the initial assumptions, $a^{(i_0)}$ is the $i_0$-th virtual member of the chart of $U[u^{(i_0)} \leadsto a^{(i_0)}]$. Therefore, $\dbtilde{a}^{(i_0)}$ may differ from $a^{(i_0)}$ only by a small piece at the beginning and by a small piece at the end. Namely,
\begin{equation*}
a^{(i_0)} = e^{\prime}m_v(a^{(i_0)})f^{\prime},\ \dbtilde{a}^{(i_0)} = \widetilde{e}^{\prime}m_v(a^{(i_0)})\widetilde{f}^{\prime},
\end{equation*}
where $e^{\prime}, f^{\prime}, \widetilde{e}^{\prime}, \widetilde{f}^{\prime}$ are small pieces.

Notice that changes of $u^{(i_0)}$ after the consecutive replacements $u^{(i_s)} \leadsto a^{(i_s)}$, $s = 1, \ldots, k$, $i_s \neq i_0$, from the left side and from the right side are independent. This happens because the replacements from the left side of $i_0$ do not influence  the form of the replacements from the right side of $i_0$, and $m_v(u^{(i_0)})$, which stays unchanged, is not a small piece. Namely, they are independent in the following sense. Let us perform all the replacements $u^{(i_s)} \leadsto a^{(i_s)}$, $s = 1, \ldots, k$, such that $i_s < i_0$, starting from $U$. Let $V$ be the resulting monomial. Then the $i_0$-th virtual member of the chart of $V$ is of the form
\begin{equation*}
\widetilde{e}t_v(u^{(i_0)}) = \widetilde{e}m_v(u^{(i_0)})f.
\end{equation*}
In the same way, if we consecutively perform all the replacements $u^{(i_s)} \leadsto a^{(i_s)}$, $s = 1, \ldots, k$, such that $i_s > i_0$, starting from $U$, then the $i_0$-th virtual member of the chart of the resulting monomial is equal to
\begin{equation*}
i_v(u^{(i_0)})\widetilde{f} = em_v(u^{(i_0)})\widetilde{f}.
\end{equation*}
However, for now we are interested only in the replacements $u^{(i_s)} \leadsto a^{(i_s)}$, $s = 1, \ldots, k$, in $U$ such that $i_s < i_0$.

Similarly, changes of $a^{(i_0)}$ after the corresponding replacements in $U[u^{(i_0)} \leadsto a^{(i_0)}]$ are independent from the left side and from the right side. Namely, consider all the replacements in $U[u^{(i_0)} \leadsto a^{(i_0)}]$ that correspond to the replacements $u^{(i_s)} \leadsto a^{(i_s)}$, $s = 1, \ldots, k$, such that $i_s < i_0$. Let us consecutively perform them, starting from $U[u^{(i_0)} \leadsto a^{(i_0)}]$, and let $\widetilde{V}$ be the resulting monomial. Similarly, since $a^{(i_0)}$ is the $i_0$-th virtual member of the chart of $U[u^{(i_0)} \leadsto a^{(i_0)}]$, the $i_0$-th virtual member of the chart of $\widetilde{V}$ is of the form
\begin{equation*}
\widetilde{e}^{\prime}t_v(a^{(i_0)}) = \widetilde{e}^{\prime}m_v(a^{(i_0)})f^{\prime}.
\end{equation*}

On the other hand, by definition, $\widetilde{V}$ is the resulting monomial of the consecutive replacements $u^{(i_s)} \leadsto a^{(i_s)}$, $s = 1, \ldots, k$, such that $i_s \leqslant i_0$, starting from $U$. Assume that first we perform the replacements $u^{(i_s)} \leadsto a^{(i_s)}$, $s = 1, \ldots, k$, such that $i_s < i_0$. Recall that $V$ is the resulting monomial. Since the $i_0$-th virtual member of the chart of $V$ is of the form
\begin{equation*}
\widetilde{e}t_v(u^{(i_0)}) = \widetilde{e}\left(e^{-1}\cdot u^{(i_0)}\right),
\end{equation*}
we see that the replacement in $V$ that corresponds to the replacement $u^{(i_0)} \leadsto a^{(i_0)}$ in $U$ is of the form
\begin{equation*}
\widetilde{e}t_v(u^{(i_0)}) \leadsto \widetilde{e}\cdot e^{-1}\cdot a^{(i_0)},
\end{equation*}
wherein $V[\widetilde{e}t_v(u^{(i_0)}) \leadsto \widetilde{e}\cdot e^{-1}\cdot a^{(i_0)}] = \widetilde{V}$ and $\widetilde{e}\cdot e^{-1}\cdot a^{(i_0)}$ is the $i_0$-th virtual member of the chart of $\widetilde{V}$. However, we proved above that the $i_0$-th virtual member of the chart of $\widetilde{V}$ is of the form $\widetilde{e}^{\prime}t_v(a^{(i_0)})$. Therefore,
\begin{equation*}
\widetilde{e}\cdot e^{-1}\cdot a^{(i_0)} = \widetilde{e}^{\prime}t_v(a^{(i_0)}).
\end{equation*}
Since $a^{(i_0)} = e^{\prime}m_v(a^{(i_0)})f^{\prime} = e^{\prime}t_v(a^{(i_0)})$, we obtain
\begin{equation*}
\widetilde{e}^{\prime}t_v(a^{(i_0)}) = \widetilde{e}\cdot e^{-1}\cdot a^{(i_0)} = \widetilde{e}\cdot e^{-1}\cdot e^{\prime} t_v(a^{(i_0)}).
\end{equation*}
Hence, $\widetilde{e}^{\prime} = \widetilde{e}\cdot e^{-1}\cdot e^{\prime}$ and $\widetilde{e}^{\prime} \cdot {e^{\prime}}^{-1}= \widetilde{e}\cdot e^{-1}$.

Similarly, considering the replacements $u^{(i_s)} \leadsto a^{(i_s)}$, $s = 1, \ldots, k$, such that $i_s \geqslant i_0$, one can show that $f^{-1}\cdot \widetilde{f} = {f^{\prime}}^{-1}\cdot \widetilde{f}^{\prime}$.
\end{proof}

\begin{remark}
There are two comments about Lemma~\ref{same_small_changes}. So, assume we are under the conditions of Lemma~\ref{same_small_changes} and use the same notations. Let us consecutively perform all the replacements $u^{(i_s)} \leadsto a^{(i_s)}$, $s = 1, \ldots, k$, starting from $U$, and let $W$ be the resulting monomial.
\begin{enumerate}
\item
Since $a^{(i_0)}$ is a virtual member of the chart of $U[u^{(i_0)} \leadsto a^{(i_0)}]$ and we replace only virtual members of the chart by virtual members of the chart, we obtain that $\dbtilde{a}^{(i_0)}$ is the $i_0$-th virtual member of the chart of $W$. Since $i_0$ is an arbitrary index and the resulting monomial $W$ does not depend on the order of the replacements, we see that $i_s$-th virtual member of the chart of $W$ differs from $a^{(i_s)}$ by at most a small piece at the beginning and a small piece at the end for all $s = 1, \ldots, k$.

\item
By definition, $e^{\prime}$ is an overlap of $a^{(i_0)}$ and the left neighbouring virtual member of the chart of $U[u^{(i_0)} \leadsto a^{(i_0)}]$, $f^{\prime}$ is an overlap of $a^{(i_0)}$ and the right neighbouring virtual member of the chart of $U[u^{(i_0)} \leadsto a^{(i_0)}]$.

Since we replace only virtual members of the chart by virtual members of the chart, one can easily see that $\widetilde{e}^{\prime}$ is an overlap of $\dbtilde{a}^{(i_0)}$ and the left neighbouring virtual member of the chart of $W$, $\widetilde{f}^{\prime}$ is an overlap of $\dbtilde{a}^{(i_0)}$ and the right neighbouring virtual member of the chart of $W$.
\end{enumerate}
\end{remark}

\begin{lemma}
\label{create_initial_replacenets}
Let $U$ be a monomial. Assume $W$ is a derived monomial of $U$ and $f(U) = f(W)$. Then there exists a set of replacements $u^{(i_s)}\leadsto a^{(i_s)}$ in $U$, $s = 1, \ldots, k$, such that $u^{(i_s)}$ and $a^{(i_s)}$ are $U$-incident monomials, $u^{(i_s)}$ is the $i_s$-th virtual member of the chart of $U$, $a^{(i_s)}$ is the $i_s$-th virtual member of the chart of $U[u^{(i_s)}\leadsto a^{(i_s)}]$, and if we consecutively perform these replacements, then $W$ is the resulting monomial.
\end{lemma}
\begin{proof}
Since $W$ is a derived monomial of $U$ and $f(U) = f(W)$, there exists a sequence of replacements
\begin{equation}
\label{derived_sequence_for_initial_replacements}
U = Z_1 \overset{\phi_1}{\longmapsto} \ldots \overset{\phi_{t}}{\longmapsto} Z_{t + 1} = W
\end{equation}
such that every $\phi_l$ is a replacement of a virtual member of the chart of $W_l$ by an incident monomial such that it is a virtual member of the chart of $W_{l + 1}$. Let us prove Lemma~\ref{create_initial_replacenets} by induction on $t$. If $t = 1$, then the statement is trivial.

Assume $t > 1$. Consider the initial part of the sequence~\eqref{derived_sequence_for_initial_replacements}. Namely, we consider the sequence
\begin{equation*}
U = Z_1 \overset{\phi_1}{\longmapsto} \ldots \overset{\phi_{t - 1}}{\longmapsto}Z_{t}.
\end{equation*}
By the induction hypothesis, there exists a set of replacements $u^{(i_s)}\leadsto a^{(i_s)}$ in $U$, $s = 1, \ldots, n$, that satisfies the conditions of Definition~\ref{consecutive_replacements_def}, and such that $Z_{t}$ is the resulting monomial of their consecutive performing. Let
\begin{equation}
\label{direct_sequcence}
U = X_1 \leadsto X_2 \leadsto \ldots \leadsto X_{n +1} = Z_{t}
\end{equation}
be a consecutive performing the replacements $u^{(i_s)}\leadsto a^{(i_s)}$, starting from $U$, $s = 1, \ldots, n$.

Consider the last transformation $Z_{t} \overset{\phi_{t}}{\longmapsto} Z_{t + 1}$. Assume $\phi_{t}$ is a replacement of $j$-th virtual member of the chart of $Z_{t}$.

First assume that $j \neq i_s$ for all $s = 1, \ldots, n$. Let the replacement $\phi_{t}$ be of the form $\widetilde{u}^{(j)} \mapsto \dbtilde{c}^{(j)}$, where $\widetilde{u}^{(j)}$ is a $j$-th virtual member of the chart of $Z_{t}$, and $\dbtilde{c}^{(j)}$ is the $j$-th virtual member of the chart of $Z_{t + 1} = Z_t[\widetilde{u}^{(j)} \mapsto \dbtilde{c}^{(j)}]$. Let $u^{(j)}$ be the $j$-th virtual member of the chart of $U$. Since in sequence~\eqref{direct_sequcence} we replace only virtual members of the chart by virtual members of the chart, we see that
\begin{equation*}
u^{(j)} = em_v(u^{(j)})f,\ \widetilde{u}^{(j)} = \widetilde{e}m_v(u^{(j)})\widetilde{f},
\end{equation*}
where $e, f, \widetilde{e}, \widetilde{f}$ are small pieces. Let us put
\begin{equation*}
c^{(j)} = e\cdot \widetilde{e}^{-1}\cdot \dbtilde{c}^{(j)}\cdot \widetilde{f}^{-1}\cdot f.
\end{equation*}
Since in the sequence~\eqref{direct_sequcence} we replace only virtual members of the chart by virtual members of the chart, one is $X_l$-incident to the next one at every step, $l = 1, \ldots, n$, we can perform the corresponding reverse replacements in reverse order. Namely, $Z_{t}$ and $U$ are connected by the following sequence of replacements of virtual members of the chart by virtual members of the chart, one is $X_l$-incident to the next one at every step, $l = n + 1, \ldots, 2$:
\begin{equation}
\label{inverse_sequcence}
Z_{t} = X_{n + 1} \leadsto \ldots \leadsto X_1 = U.
\end{equation}
Therefore, since $\widetilde{u}^{(j)}$ and $\dbtilde{c}^{(j)}$ are $Z_{t}$-incident monomials, it follows from statements~\ref{virtual_members_U_incident_replacements2} and~\ref{virtual_members_U_incident_replacements5} of Corollary~\ref{virtual_members_U_incident_replacements} that $u^{(j)}$ and $c^{(j)}$ are $U$-incident monomials. Since $\dbtilde{c}^{(j)}$ is a virtual member of the chart of $Z_{t + 1}$, it follows from statements~\ref{virtual_members_U_incident_replacements3} and~\ref{virtual_members_U_incident_replacements6} of Corollary~\ref{virtual_members_U_incident_replacements} that $c^{(j)}$ is a virtual member of the chart of $U[u^{(j)} \leadsto c^{(j)}]$. By construction, we obtain that $\widetilde{u}^{(j)} \leadsto \dbtilde{c}^{(j)}$ in $Z_{t}$ corresponds to the replacement $u^{(j)} \leadsto c^{(j)}$ in $U$ after consecutive performing the replacements $u^{(i_s)}\leadsto a^{(i_s)}$, $s = 1, \ldots, n$. Therefore, $u^{(i_s)}\leadsto a^{(i_s)}$, $s = 1, \ldots, n$, together with $u^{(j)} \leadsto c^{(j)}$, is the required set of replacements in $U$.

Now assume that $j = i_{s_0}$ for some $1 \leqslant s_0 \leqslant n$. Then assume that sequence~\eqref{direct_sequcence} is performing  the replacements $u^{(i_s)}\leadsto a^{(i_s)}$, $s = 1, \ldots, n$, in order $j_1, \ldots, j_{n}$ such that $j_{n} = i_{s_0} = j$. Let $\widetilde{u}^{(j)}$ be the $j$-th virtual member of the chart of $X_{n}$. Let the last replacement $X_{n} \leadsto X_{n + 1}$ be of the form $\widetilde{u}^{(j)} \leadsto \dbtilde{a}^{(j)}$ in $X_{n}$. Then $\dbtilde{a}^{(j)}$ is $j$-th virtual member of the chart of $X_{n + 1} = Z_{t}$. We assumed that the last replacement $X_{n + 1} = Z_{t} \overset{\phi_{t}}{\longmapsto} Z_{t + 1}$ is a replacement of the $j$-th virtual member of the chart of $Z_t$. So, it is a replacement of $\dbtilde{a}^{(j)}$ by an incident monomial. Assume this replacement is of the form $\dbtilde{a}^{(j)} \mapsto \dbtilde{c}^{(j)}$, where $\dbtilde{c}^{(j)}$ is a virtual member of the chart of $Z_{t + 1} = Z_t[\dbtilde{a}^{(j)} \mapsto \dbtilde{c}^{(j)}]$.

We glue the replacement $\widetilde{u}^{(j)} \leadsto \dbtilde{a}^{(j)}$ in $X_{n}$ and the replacement $\dbtilde{a}^{(j)} \mapsto \dbtilde{c}^{(j)}$ in $Z_{t}$. Then we obtain the replacement $\widetilde{u}^{(j)} \leadsto \dbtilde{c}^{(j)}$ in $X_n$, and $Z_{t + 1} = W$ is the resulting monomial. Since $\dbtilde{a}^{(j)}$ is a virtual member of the chart of $X_{n + 1} = Z_{t}$, we see that $\widetilde{u}^{(j)}$ and  $\dbtilde{c}^{(j)}$ are $X_{n}$-incident monomials. As above, we can construct the replacement $u^{(j)} \leadsto c^{(j)}$ in $U$ such that $u^{(j)}$ and $c^{(j)}$ are $U$-incident monomials, and $c^{(j)}$ is a virtual member of the chart of $U[u^{(j)} \leadsto c^{(j)}]$, and the replacement $\widetilde{u}^{(j)} \leadsto \dbtilde{c}^{(j)}$ in $X_n$ corresponds to the replacement $u^{(j)} \leadsto c^{(j)}$ in $U$ after consecutive performing  the replacements $u^{(i_s)}\leadsto a^{(i_s)}$, $s = 1, \ldots, n$, $s \neq s_0$. Therefore, $u^{(i_s)}\leadsto a^{(i_s)}$, $s = 1, \ldots, n$, $s \neq s_0$, together with $u^{(i_{s_0})} \leadsto c^{(i_{s_0})}$ is the required set of replacements in $U$.
\end{proof}

\medskip

As above, let $U$ be a monomial and let $u^{(i_1)}$ and $u^{(i_2)}$ be virtual members of the chart of $U$. Assume $u^{(i_1)}$ starts from the left of the beginning of $u^{(i_2)}$. Let $u^{(i_1)}$, $a^{(i_1)}$ and $u^{(i_2)}$, $a^{(i_2)}$ be $U$-incident monomials. Now we study a general case with two replacements. Namely, we consider the replacements $u^{(i_1)} \leadsto a^{(i_1)}$ and $u^{(i_2)} \leadsto a^{(i_2)}$ in $U$ such that $a^{(i_1)}$ is not necessary a virtual member of the chart of $U[u^{(i_1)} \leadsto a^{(i_1)}]$, and $a^{(i_2)}$ is not necessary a virtual member of the chart of $U[u^{(i_2)} \leadsto a^{(i_2)}]$. Let us do the first replacement $u^{(i_1)} \leadsto a^{(i_1)}$ in $U$ and consider what happens with the second replacement.

Let $U = L^{(i_1)}u^{(i_1)}R^{(i_1)}$ and $\widetilde{u}^{(i_2)}$ be the image of $u^{(i_2)}$ in $U[u^{(i_1)} \leadsto a^{(i_1)}] = L^{(i_1)}a^{(i_1)}R^{(i_1)}$. The first problem that may appear is as follows\label{go_lower_difficulties}. Assume $a^{(i_1)} = 1$. Then the monomial $U[u^{(i_1)} \leadsto a^{(i_1)}] = L^{(i_1)}a^{(i_1)}R^{(i_1)} = L^{(i_1)}\cdot R^{(i_1)}$ does not have to be reduced. Hence, $u^{(i_2)}$ may have an empty image in $U[u^{(i_1)} \leadsto a^{(i_1)}]$ because of cancellations. Then the replacement of $\widetilde{u}^{(i_2)}$ is not defined.

We deal with this issue as follows. Let $\widehat{u}^{(i_2)}$ be the intersection of $u^{(i_2)}$ and $R^{(i_1)}$. Let us put
\begin{equation}
\label{corresponding_replacement_in_intersection_right}
\dbhat{a}^{(i_2)} = \begin{cases}
a^{(i_2)} &\textit{if } u^{(i_1)} \textit{ and } u^{(i_2)} \textit{ are separated or touch at a point}\\
&\textit{(that is, }  u^{(i_2)} = \widehat{u}^{(i_2)} \textit{)};\\
e^{-1}\cdot a^{(i_2)} &\textit{if } u^{(i_1)} \textit{ and } u^{(i_2)} \textit{ have an overlap } e\\
&\textit{(that is, }  u^{(i_2)} = e\widehat{u}^{(i_2)} \textit{)}.\\
\end{cases}
\end{equation}
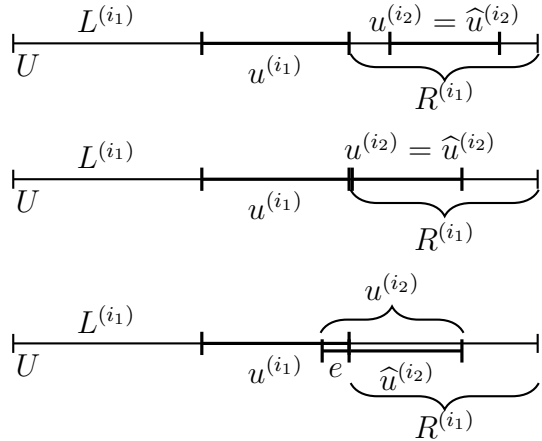
\begin{center}
\begin{tikzpicture}
\draw[|-|, black, thick] (-1.5,0)--(5.5,0);
\node[below, xshift=6] at (-1.5, 0) {$U$};
\draw[|-|, black, very thick] (1,0)--(3,0) node[midway, below] {$u^{(i_1)}$};
\draw[|-|, black, very thick] (3.5,0)--(5,0) node[midway, above] {$u^{(i_2)} = \widehat{u}^{(i_2)}$};
\draw [thick, decorate, decoration={brace, amplitude=10pt, raise=4pt, mirror}] (3, 0) to node[midway, below, yshift=-10pt] {$R^{(i_1)}$} (5.5, 0);
\path (-1.5, 0)--(1, 0) node[midway, above] {$L^{(i_1)}$};
\end{tikzpicture}

\begin{tikzpicture}
\draw[|-|, black, thick] (-1.5,0)--(5.5,0);
\node[below, xshift=6] at (-1.5, 0) {$U$};
\draw[|-|, black, very thick] (1,0)--(3,0) node[midway, below] {$u^{(i_1)}$};
\draw[|-|, black, very thick] (3,0)--(4.5,0) node[midway, above, xshift=5, yshift=3] {$u^{(i_2)} = \widehat{u}^{(i_2)}$};
\draw [thick, decorate, decoration={brace, amplitude=10pt, raise=4pt, mirror}] (3, 0) to node[midway, below, yshift=-12pt] {$R^{(i_1)}$} (5.5, 0);
\path (-1.5, 0)--(1, 0) node[midway, above] {$L^{(i_1)}$};
\end{tikzpicture}

\begin{tikzpicture}
\draw[|-|, black, thick] (-1.5, 0)--(5.5, 0);
\node[below, xshift=6] at (-1.5, 0) {$U$};
\draw[|-|, black, very thick] (1,0)--(3,0) node[midway, below] {$u^{(i_1)}$};
\draw[-|, black, very thick] (3,-0.1)--(4.5,-0.1);
\draw[|-|, black, very thick] (2.6,-0.1)--(3,-0.1) node[midway, below] {$e$};
\path (3, 0)--(4.5, 0) node[midway, below, yshift=-3] {$\widehat{u}^{(i_2)}$};
\draw [thick, decorate, decoration={brace, amplitude=10pt, raise=4pt}] (2.6, 0) to node[midway, above, yshift=12pt] {$u^{(i_2)}$} (4.5, 0);
\draw [thick, decorate, decoration={brace, amplitude=10pt, raise=15pt, mirror}] (3, 0) to node[midway, below, yshift=-21pt] {$R^{(i_1)}$} (5.5, 0);
\path (-1.5, 0)--(1, 0) node[midway, above] {$L^{(i_1)}$};
\end{tikzpicture}
\end{center}
Instead of replacements of maximal occurrences, we perform the following sequence of replacements, starting from $U$.
\begin{enumerate}
\label{replacement_in_intersection_process}
\item
We start with the replacement $u^{(i_1)} \leadsto a^{(i_1)}$ in $U$. The resulting monomial is $U[u^{(i_1)} \leadsto a^{(i_1)}] = L^{(i_1)}a^{(i_1)}R^{(i_1)}$. If $a^{(i_1)} = 1$ and $L^{(i_1)}\cdot R^{(i_1)}$ is not reduced, then we do not perform the cancellations.

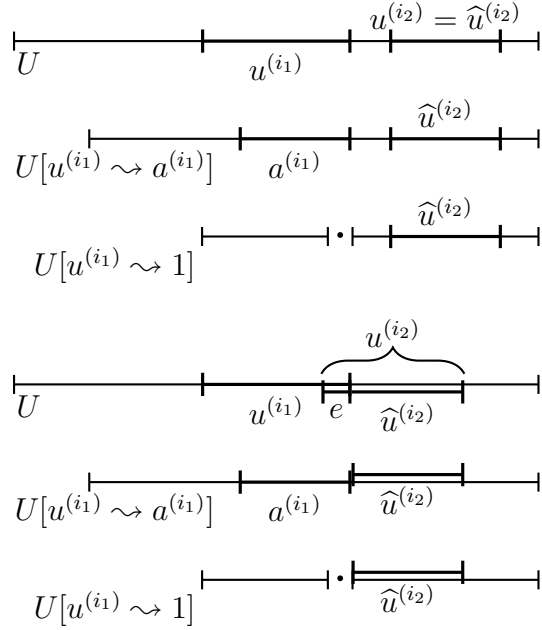
\begin{center}
\begin{tikzpicture}
\draw[|-|, black, thick] (-1.5,0)--(5.5,0);
\node[below, xshift=6] at (-1.5, 0) {$U$};
\draw[|-|, black, very thick] (1,0)--(3,0) node[midway, below] {$u^{(i_1)}$};
\draw[|-|, black, very thick] (3.5,0)--(5,0) node[midway, above] {$u^{(i_2)} = \widehat{u}^{(i_2)}$};

\draw[|-|, black, thick] (-0.5, -1.3)--(5.5, -1.3);
\node[below, xshift=10] at (-0.5, -1.3) {$U[u^{(i_1)}\leadsto a^{(i_1)}]$};
\draw[|-|, black, very thick] (1.5, -1.3)--(3, -1.3) node[midway, below] {$a^{(i_1)}$};
\draw[|-|, black, very thick] (3.5, -1.3)--(5, -1.3) node[midway, above] {$\widehat{u}^{(i_2)}$};

\draw[|-|, black, thick] (1, -2.6)--(2.7, -2.6);
\draw[|-|, black, thick] (3, -2.6)--(5.5, -2.6);
\node at (2.85, -2.6) {\LARGE $\cdot$};
\node[below, xshift=10] at (-0.5, -2.6) {$U[u^{(i_1)}\leadsto 1]$};
\draw[|-|, black, very thick] (3.5, -2.6)--(5, -2.6) node[midway, above] {$\widehat{u}^{(i_2)}$};
\end{tikzpicture}

\begin{tikzpicture}
\draw[|-|, black, thick] (-1.5, 0)--(5.5, 0);
\node[below, xshift=6] at (-1.5, 0) {$U$};
\draw[|-|, black, very thick] (1,0)--(3,0) node[midway, below] {$u^{(i_1)}$};
\draw[-|, black, very thick] (3,-0.1)--(4.5,-0.1);
\draw[|-|, black, very thick] (2.6,-0.1)--(3,-0.1) node[midway, below] {$e$};
\path (3, 0)--(4.5, 0) node[midway, below, yshift=-3] {$\widehat{u}^{(i_2)}$};
\draw [thick, decorate, decoration={brace, amplitude=10pt, raise=4pt}] (2.6, 0) to node[midway, above, yshift=10pt] {$u^{(i_2)}$} (4.5, 0);

\draw[|-|, black, thick] (-0.5, -1.3)--(5.5, -1.3);
\node[below, xshift=10] at (-0.5, -1.3) {$U[u^{(i_1)}\leadsto a^{(i_1)}]$};
\draw[|-|, black, very thick] (1.5, -1.3)--(3,-1.3) node[midway, below] {$a^{(i_1)}$};
\draw[|-|, black, very thick] (3, -1.2)--(4.5, -1.2) node[midway, below] {$\widehat{u}^{(i_2)}$};

\draw[|-|, black, thick] (1, -2.6)--(2.7, -2.6);
\draw[|-|, black, thick] (3, -2.6)--(5.5, -2.6);
\node at (2.85, -2.6) {\LARGE $\cdot$};
\node[below, xshift=10] at (-0.5, -2.6) {$U[u^{(i_1)}\leadsto 1]$};
\draw[|-|, black, very thick] (3, -2.5)--(4.5, -2.5) node[midway, below] {$\widehat{u}^{(i_2)}$};
\end{tikzpicture}

\begin{tikzpicture}
\node at (0, 0) {{\LARGE $\cdot$} \textit {on the above pictures means possibility of cancellations}};
\end{tikzpicture}
\end{center}
\item
We do the replacement $\widehat{u}^{(i_2)} \leadsto \dbhat{a}^{(i_2)}$ in $R^{(i_1)}$. If $\dbhat{a}^{(i_2)} = 1$ and the monomial that we obtain after the replacement $\widehat{u}^{(i_2)} \leadsto \dbhat{a}^{(i_2)}$ in $R^{(i_1)}$ is not reduced, then we perform the cancellations in this monomial. Let $R^{(i_1)}[\widehat{u}^{(i_2)} \leadsto \dbhat{a}^{(i_2)}]$ be the resulting monomial after the cancellations.

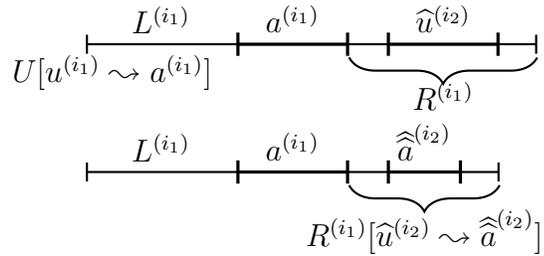
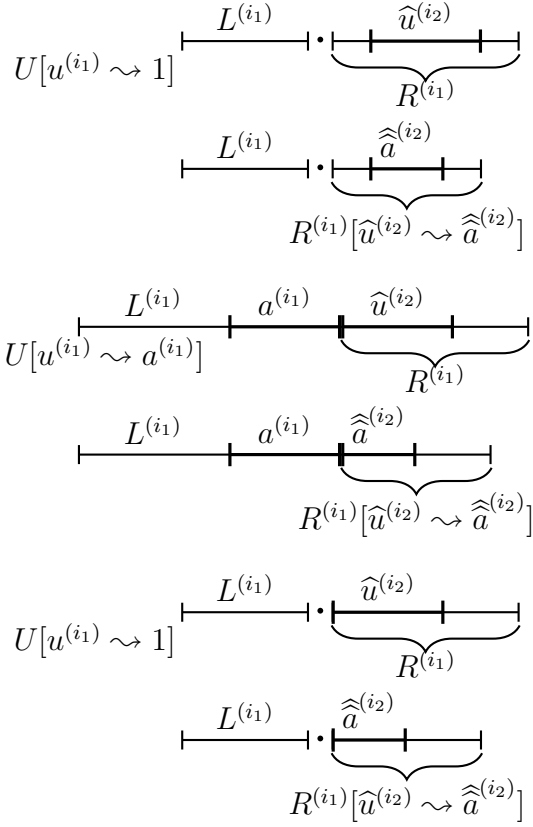
\begin{center}
\begin{tikzpicture}
\draw[|-|, black, thick] (-0.5, 0)--(5.5, 0);
\node[below, xshift=10] at (-0.5, 0) {$U[u^{(i_1)}\leadsto a^{(i_1)}]$};
\draw[|-|, black, very thick] (1.5, 0)--(3, 0) node[midway, above] {$a^{(i_1)}$};
\draw[|-|, black, very thick] (3.5, 0)--(5, 0) node[midway, above] {$\widehat{u}^{(i_2)}$};
\draw [thick, decorate, decoration={brace, amplitude=10pt, raise=5pt, mirror}] (3, 0) to node[midway, below, yshift=-10pt] {$R^{(i_1)}$} (5.5, 0);
\path (-0.5, 0)--(1.5, 0) node[midway, above] {$L^{(i_1)}$};

\draw[|-|, black, thick] (-0.5, -1.7)--(5, -1.7);
\draw[|-|, black, very thick] (1.5, -1.7)--(3, -1.7) node[midway, above] {$a^{(i_1)}$};
\draw[|-|, black, very thick] (3.5, -1.7)--(4.5, -1.7) node[midway, above] {$\dbhat{a}^{(i_2)}$};
\draw [thick, decorate, decoration={brace, amplitude=10pt, raise=5pt, mirror}] (3, -1.7) to node[midway, below, yshift=-10pt] {$R^{(i_1)}[\widehat{u}^{(i_2)} \leadsto \dbhat{a}^{(i_2)}]$} (5, -1.7);
\path (-0.5, -1.7)--(1.5, -1.7) node[midway, above] {$L^{(i_1)}$};
\end{tikzpicture}

\begin{tikzpicture}
\draw[|-|, black, thick] (1, 0)--(2.7, 0);
\draw[|-|, black, thick] (3, 0)--(5.5, 0);
\node at (2.85, 0) {\LARGE $\cdot$};
\node[below, xshift=10] at (-0.5, 0) {$U[u^{(i_1)}\leadsto 1]$};
\draw[|-|, black, very thick] (3.5, 0)--(5, 0) node[midway, above] {$\widehat{u}^{(i_2)}$};
\draw [thick, decorate, decoration={brace, amplitude=10pt, raise=5pt, mirror}] (3, 0) to node[midway, below, yshift=-10pt] {$R^{(i_1)}$} (5.5, 0);
\path (1, 0)--(2.7, 0) node[midway, above] {$L^{(i_1)}$};

\draw[|-|, black, thick] (1, -1.7)--(2.7, -1.7);
\draw[|-|, black, thick] (3, -1.7)--(5, -1.7);
\node at (2.85, -1.7) {\LARGE $\cdot$};
\draw[|-|, black, very thick] (3.5, -1.7)--(4.5, -1.7) node[midway, above] {$\dbhat{a}^{(i_2)}$};
\draw [thick, decorate, decoration={brace, amplitude=10pt, raise=5pt, mirror}] (3, -1.7) to node[midway, below, yshift=-10pt] {$R^{(i_1)}[\widehat{u}^{(i_2)} \leadsto \dbhat{a}^{(i_2)}]$} (5, -1.7);
\path (1, -1.7)--(2.7, -1.7) node[midway, above] {$L^{(i_1)}$};
\end{tikzpicture}

\begin{tikzpicture}
\draw[|-|, black, thick] (-0.5, 0)--(5.5, 0);
\node[below, xshift=10] at (-0.5, 0) {$U[u^{(i_1)}\leadsto a^{(i_1)}]$};
\draw[|-|, black, very thick] (1.5, 0)--(3, 0) node[midway, above] {$a^{(i_1)}$};
\draw[|-|, black, very thick] (3, 0)--(4.5, 0) node[midway, above] {$\widehat{u}^{(i_2)}$};
\draw [thick, decorate, decoration={brace, amplitude=10pt, raise=5pt, mirror}] (3, 0) to node[midway, below, yshift=-10pt] {$R^{(i_1)}$} (5.5, 0);
\path (-0.5, 0)--(1.5, 0) node[midway, above] {$L^{(i_1)}$};

\draw[|-|, black, thick] (-0.5, -1.7)--(5, -1.7);
\draw[|-|, black, very thick] (1.5, -1.7)--(3, -1.7) node[midway, above] {$a^{(i_1)}$};
\draw[|-|, black, very thick] (3, -1.7)--(4, -1.7) node[midway, above] {$\dbhat{a}^{(i_2)}$};
\draw [thick, decorate, decoration={brace, amplitude=10pt, raise=5pt, mirror}] (3, -1.7) to node[midway, below, yshift=-10pt] {$R^{(i_1)}[\widehat{u}^{(i_2)} \leadsto \dbhat{a}^{(i_2)}]$} (5, -1.7);
\path (-0.5, -1.7)--(1.5, -1.7) node[midway, above] {$L^{(i_1)}$};
\end{tikzpicture}

\begin{tikzpicture}
\draw[|-|, black, thick] (1, 0)--(2.7, 0);
\draw[|-|, black, thick] (3, 0)--(5.5, 0);
\node at (2.85, 0) {\LARGE $\cdot$};
\node[below, xshift=10] at (-0.5, 0) {$U[u^{(i_1)}\leadsto 1]$};
\draw[|-|, black, very thick] (3, 0)--(4.5, 0) node[midway, above] {$\widehat{u}^{(i_2)}$};
\draw [thick, decorate, decoration={brace, amplitude=10pt, raise=5pt, mirror}] (3, 0) to node[midway, below, yshift=-10pt] {$R^{(i_1)}$} (5.5, 0);
\path (1, 0)--(2.7, 0) node[midway, above] {$L^{(i_1)}$};

\draw[|-|, black, thick] (1, -1.7)--(2.7, -1.7);
\draw[|-|, black, thick] (3, -1.7)--(5, -1.7);
\node at (2.85, -1.7) {\LARGE $\cdot$};
\draw[|-|, black, very thick] (3, -1.7)--(4, -1.7) node[midway, above] {$\dbhat{a}^{(i_2)}$};
\draw [thick, decorate, decoration={brace, amplitude=10pt, raise=5pt, mirror}] (3, -1.7) to node[midway, below, yshift=-10pt] {$R^{(i_1)}[\widehat{u}^{(i_2)} \leadsto \dbhat{a}^{(i_2)}]$} (5, -1.7);
\path (1, -1.7)--(2.7, -1.7) node[midway, above] {$L^{(i_1)}$};
\end{tikzpicture}

\begin{tikzpicture}
\node at (0, 0) {{\LARGE $\cdot$} \textit {on the above pictures means possibility of cancellations}};
\end{tikzpicture}
\end{center}

Notice that the replacement $\widehat{u}^{(i_2)} \leadsto \dbhat{a}^{(i_2)}$ in $R^{(i_1)}$ is a formal procedure. We do not claim neither that $\widehat{u}^{(i_2)}$ is a virtual member of the chart of $R^{(i_1)}$, nor that $\widehat{u}^{(i_2)}$ and $\dbhat{a}^{(i_2)}$ are $R^{(i_1)}$-incident monomials.
\item
After the second replacement we obtain the monomial $\left(L^{(i_1)}a^{(i_1)}\right)\cdot R^{(i_1)}[\widehat{u}^{(i_2)} \leadsto \dbhat{a}^{(i_2)}]$. The last step is performing  the cancellations in this monomial if there are any.
\end{enumerate}

Let $U = L^{(i_2)}u^{(i_2)}R^{(i_2)}$. Let $\widehat{u}^{(i_1)}$ be the intersection of $u^{(i_1)}$ and $L^{(i_2)}$. Similarly to~\eqref{corresponding_replacement_in_intersection_right}, we define
\begin{equation}
\label{corresponding_replacement_in_intersection_left}
\dbhat{a}^{(i_1)} = \begin{cases}
a^{(i_1)} &\textit{if } u^{(i_1)} \textit{ and } u^{(i_2)} \textit{ are separated or touch at a point}\\
&\textit{(that is, }  u^{(i_2)} = \widehat{u}^{(i_2)} \textit{)};\\
a^{(i_1)}\cdot e^{-1} &\textit{if } u^{(i_1)} \textit{ and } u^{(i_2)} \textit{ have an overlap } e\\
&\textit{(that is, }  u^{(i_1)} = \widehat{u}^{(i_1)}e \textit{)}.\\
\end{cases}
\end{equation}
Clearly, we can also do the above process, starting from the replacement $a^{(i_2)} \leadsto u^{(i_2)}$ in $U$. But then at the second step we deal with the intersection of $u^{(i_1)}$ and $L^{(i_2)}$ and perform the replacement $\widehat{u}^{(i_1)} \leadsto \dbhat{a}^{(i_1)}$ in $L^{(i_2)}$.

Let us explain why the process of replacements introduced above agrees with our idea of replacements of maximal occurrences. Assume that the monomial $L^{(i_1)}a^{(i_1)}R^{(i_1)}$ does not have cancellations. Then $u^{(i_2)}$ has a non-empty image in $U[u^{(i_1)} \leadsto a^{(i_1)}] = L^{(i_1)}a^{(i_1)}R^{(i_1)}$. As above, let $\widetilde{u}^{(i_2)}$ be the image of $u^{(i_2)}$ in $U[u^{(i_1)} \leadsto a^{(i_1)}]$ and let $\widehat{u}^{(i_2)}$ be the intersection of $u^{(i_2)}$ and $R^{(i_1)}$. It follows from the results of Section~\ref{mt_configurations} that
\begin{equation}
\label{image_of_right_general}
\widetilde{u}^{(i_2)} = \begin{cases}
u^{(i_2)} = \widehat{u}^{(i_2)} &\textit{if } u^{(i_1)} \textit{ and } u^{(i_2)} \textit{ are separated},\\
X(e^{-1}\cdot u^{(i_2)}) = X\widehat{u}^{(i_2)} &\textit{if } u^{(i_1)} \textit{ and } u^{(i_2)} \textit{ are not separated},\\
&e \textit{ is the overlap of } u^{(i_1)} \textit{ and } u^{(i_2)}\\
&\textit{(empty if }u^{(i_1)} \textit{ and } u^{(i_2)} \textit{ touch at a point)},\\
&X \textit{ is a maximal prolongation of } \widehat{u}^{(i_2)}\\
&\textit{to the left in } L^{(i_1)}a^{(i_1)}R^{(i_1)}.
\end{cases}
\end{equation}
Similar to~\eqref{corresponding_element_U_incident_right}, let us put
\begin{equation}
\label{corresponding_element_right_general}
\dbtilde{a}^{(i_2)}= \begin{cases}
a^{(i_2)} &\textit{if } u^{(i_1)}\textit{ and } u^{(i_2)} \textit{ are separated},\\
X\cdot e^{-1}\cdot a^{(i_2)} &\textit{if } u^{(i_1)} \textit{ and } u^{(i_2)} \textit{ are not separated},\\
&e \textit{ is the overlap of } u^{(i_1)} \textit{ and } u^{(i_2)}\\
&\textit{(empty if } u^{(i_1)} \textit{ and } u^{(i_2)} \textit{ touch at a point)},\\
&X \textit{ is a maximal prolongation of } \widehat{u}^{(i_2)}\\
&\textit{to the left in } L^{(i_1)}a^{(i_1)}R^{(i_1)}.
\end{cases}
\end{equation}
We can not claim that $\widetilde{u}^{(i_2)}$ and $\dbtilde{a}^{(i_2)}$ are necessarily $U[u^{(i_1)} \leadsto a^{(i_1)}]$-incident monomials. So, let us specially assume that $\widetilde{u}^{(i_2)}$ and $\dbtilde{a}^{(i_2)}$ are $U[u^{(i_1)} \leadsto a^{(i_1)}]$-incident monomials. After the replacement $u^{(i_1)} \leadsto a^{(i_1)}$ in $U$, we can perform either the replacement $\widetilde{u}^{(i_2)} \leadsto \dbtilde{a}^{(i_2)}$ in $U[u^{(i_1)} \leadsto a^{(i_1)}]$ and obtain the monomial $U[u^{(i_1)} \leadsto a^{(i_1)}][\widetilde{u}^{(i_2)} \leadsto \dbtilde{a}^{(i_2)}]$, or the replacement $\widehat{u}^{(i_2)} \leadsto \dbhat{a}^{(i_2)}$ in $R^{(i_1)}$ and the further cancellations in the monomial $\left(L^{(i_1)}a^{(i_1)}\right)\cdot R^{(i_1)}[\widehat{u}^{(i_2)} \leadsto \dbhat{a}^{(i_2)}]$. Let us show that in both cases we obtain the same resulting monomial.

\begin{lemma}
\label{intersection_and_maximal_replacement}
Let $U$ be a monomial, and let $u^{(i_1)}$ and $u^{(i_2)}$ be virtual members of the chart of $U$. Assume $u^{(i_1)}$ starts from the left of the beginning of $u^{(i_2)}$. Let $u^{(i_1)}$, $a^{(i_1)}$ and $u^{(i_2)}$, $a^{(i_2)}$ be $U$-incident monomials.

Assume $U = L^{(i_1)}u^{(i_1)}R^{(i_1)}$. Let $\widehat{u}^{(i_2)}$ be the intersection of $u^{(i_2)}$ and $R^{(i_1)}$, let $\dbhat{a}^{(i_2)}$ be defined by formula~\eqref{corresponding_replacement_in_intersection_right}. Assume $L^{(i_1)}a^{(i_1)}R^{(i_1)}$ does not have cancellations. Let $\widetilde{u}^{(i_2)}$ be the image of $u^{(i_2)}$ in $U[u^{(i_1)} \leadsto a^{(i_1)}]$ ($\widetilde{u}^{(i_2)}$ be defined by formula~\eqref{image_of_right_general}), $\dbtilde{a}^{(i_2)}$ be defined by formula~\eqref{corresponding_element_right_general}. Assume $\widetilde{u}^{(i_2)}$ and $\dbtilde{a}^{(i_2)}$ are $U[u^{(i_1)} \leadsto a^{(i_1)}]$-incident monomials. Then we obtain
\begin{equation*}
U[u^{(i_1)} \leadsto a^{(i_1)}][\widetilde{u}^{(i_2)} \leadsto \dbtilde{a}^{(i_2)}] = \left(L^{(i_1)}a^{(i_1)}\right)\cdot R^{(i_1)}[\widehat{u}^{(i_2)} \leadsto \dbhat{a}^{(i_2)}].
\end{equation*}

Similarly, let $\widehat{u}^{(i_1)}$ be the intersection of $u^{(i_1)}$ and $L^{(i_2)}$, let $\dbhat{a}^{(i_1)}$ be defined by formula~\eqref{corresponding_replacement_in_intersection_left}. Assume $L^{(i_2)}a^{(i_2)}R^{(i_2)}$ does not have cancellations. Let $\widetilde{u}^{(i_1)}$ be the image of $u^{(i_1)}$ in $U[u^{(i_2)} \leadsto a^{(i_2)}]$ ($\widetilde{u}^{(i_2)}$ be defined similarly to~\eqref{image_of_right_general}), $\dbtilde{a}^{(i_2)}$ be defined similarly to~\eqref{corresponding_element_right_general}. Assume $\widetilde{u}^{(i_1)}$ and $\dbtilde{a}^{(i_1)}$ are $U[u^{(i_2)} \leadsto a^{(i_2)}]$-incident monomials. Then we obtain
\begin{equation*}
U[u^{(i_2)} \leadsto a^{(i_2)}][\widetilde{u}^{(i_1)} \leadsto \dbtilde{a}^{(i_1)}] = L^{(i_2)}[\widehat{u}^{(i_1)} \leadsto \dbhat{a}^{(i_1)}]\cdot \left(a^{(i_2)}R^{(i_2)}\right).
\end{equation*}
\end{lemma}
\begin{proof}
Let us prove the first part of Lemma~\ref{intersection_and_maximal_replacement}. The second part is proved analogously. Denote the monomial $\left(L^{(i_1)}a^{(i_1)}\right)\cdot R^{(i_1)}[\widehat{u}^{(i_2)} \leadsto \dbhat{a}^{(i_2)}]$ by $W$.

If $u^{(i_1)}$ and $u^{(i_2)}$ are separated, then $\widetilde{u}^{(i_2)} = \widehat{u}^{(i_2)} = u^{(i_2)}$ are the same occurrence in $U[u^{(i_1)} \leadsto a^{(i_1)}]$. Then we have $\dbtilde{a}^{(i_2)} = \dbhat{a}^{(i_2)} = a^{(i_2)}$. So, $U[u^{(i_1)} \leadsto a^{(i_1)}][\widetilde{u}^{(i_2)} \leadsto \dbtilde{a}^{(i_2)}] = W$.

Assume $u^{(i_1)}$ and $u^{(i_2)}$ are not separated. Let $\widetilde{L}$ be a prefix of $U[u^{(i_1)} \leadsto a^{(i_1)}]$ that ends at the beginning point of $X$. Let $U = L^{(i_2)}u^{(i_2)}R^{(i_2)}$. Then, clearly, $U[u^{(i_1)} \leadsto a^{(i_1)}]$ can be written as follows
\begin{equation*}
U[u^{(i_1)} \leadsto a^{(i_1)}] = L^{(i_1)}a^{(i_1)}R^{(i_1)} = L^{(i_1)}a^{(i_1)}\widehat{u}^{(i_2)}R^{(i_2)} = \widetilde{L}X\widehat{u}^{(i_2)}R^{(i_2)}.
\end{equation*}
Hence, on the one hand, we see that
\begin{equation*}
U[u^{(i_1)} \leadsto a^{(i_1)}][\widetilde{u}^{(i_2)} \leadsto \dbtilde{a}^{(i_2)}] = \widetilde{L}\dbtilde{a}^{(i_2)}R^{(i_2)} = \widetilde{L}\left(X\cdot e^{-1}\cdot a^{(i_2)}\right)R^{(i_2)}.
\end{equation*}
On the other hand, we have
\begin{align*}
W &=\left(L^{(i_1)}a^{(i_1)}\right)\cdot R^{(i_1)}[\widehat{u}^{(i_2)} \leadsto \dbhat{a}^{(i_2)}] = \left(L^{(i_1)}a^{(i_1)}\right)\cdot \dbhat{a}^{(i_2)}R^{(i_2)} =\\
&= \left(L^{(i_1)}a^{(i_1)}\right)\cdot \left(e^{-1}\cdot a^{(i_2)} \right)R^{(i_2)} = \widetilde{L}X\cdot \left(e^{-1}\cdot a^{(i_2)} \right)R^{(i_2)}.
\end{align*}
Therefore, $U[u^{(i_1)} \leadsto a^{(i_1)}][\widetilde{u}^{(i_2)} \leadsto \dbtilde{a}^{(i_2)}] = W$.
\end{proof}
So, if $L^{(i_1)}a^{(i_1)}R^{(i_1)}$ is a reduced monomial and $\widetilde{u}^{(i_2)}$, $\dbtilde{a}^{(i_2)}$ are $U[u^{(i_1)} \leadsto a^{(i_1)}]$-incident monomials, then the replacement $\widehat{u}^{(i_2)} \leadsto \dbhat{a}^{(i_2)}$ in $R^{(i_1)}$ and the further cancellations in $\left(L^{(i_1)}a^{(i_1)}\right)\cdot R^{(i_1)}[\widehat{u}^{(i_2)} \leadsto \dbhat{a}^{(i_2)}]$ agree with our general idea of replacements of maximal occurrences.

Now assume that $a^{(i_1)}$ is not a virtual member of the chart of $U[u^{(i_1)} \leadsto a^{(i_1)}]$, or $a^{(i_2)}$ is not a virtual member of the chart of $U[u^{(i_2)} \leadsto a^{(i_2)}]$. Our goal is to prove that in this case the monomial $\left(L^{(i_1)}a^{(i_1)}\right)\cdot R^{(i_1)}[\widehat{u}^{(i_2)} \leadsto \dbhat{a}^{(i_2)}]$ has smaller $f$-characteristic than $U$.

\begin{remark}
\label{first_monomial_reduced_difficulties}
To be definite, assume that $a^{(i_1)}$ is not a virtual member of the chart of $U[u^{(i_1)} \leadsto a^{(i_1)}]$. Notice that even if the monomial $L^{(i_1)}a^{(i_1)}R^{(i_1)}$, which we obtain after the first replacement $u^{(i_1)} \leadsto a^{(i_1)}$ in $U$, is reduced, we still have a number of issues with the second replacement, which corresponds to the replacement $u^{(i_2)} \leadsto a^{(i_2)}$ in $U$. So, the statement that the monomial $\left(L^{(i_1)}a^{(i_1)}\right)\cdot R^{(i_1)}[\widehat{u}^{(i_2)} \leadsto \dbhat{a}^{(i_2)}]$ has smaller $f$-characteristic than $U$ is non-trivial.

Namely, the situation is as follows. Suppose that $L^{(i_1)}a^{(i_1)}R^{(i_1)}$ has no cancellations. Then $u^{(i_2)}$ always has a non-empty image in $L^{(i_1)}a^{(i_1)}R^{(i_1)}$. So, let $\widetilde{u}^{(i_2)}$ be the image of $u^{(i_2)}$ in $L^{(i_1)}a^{(i_1)}R^{(i_1)}$ ($\widetilde{u}^{(i_2)}$ be defined by formula~\eqref{image_of_right_general}). Let $\dbtilde{a}^{(i_2)}$ be defined by formula~\eqref{corresponding_element_right_general}. Then the following problems may occur.
\begin{itemize}
\item
$\widetilde{u}^{(i_2)}$ is not necessarily a virtual member of the chart of $U[u^{(i_1)} \leadsto a^{(i_1)}]$.
\item
$\widetilde{u}^{(i_2)}$ and $\dbtilde{a}^{(i_2)}$ may not be $U[u^{(i_1)} \leadsto a^{(i_1)}]$-incident monomials.
\end{itemize}
\end{remark}

The following lemma shows that if $a^{(i_1)}$ is not a virtual member of the chart of $U[u^{(i_1)} \leadsto a^{(i_1)}]$, or $a^{(i_2)}$ is not a virtual member of the chart of $U[u^{(i_2)} \leadsto a^{(i_2)}]$, then we do not claim that we always obtain a derived monomial of $U$. However, we always obtain a monomial with smaller $f$-characteristic than $U$. In Subsection~\ref{tensor_product_section} we will see that this property is sufficient for the further argument.

\begin{lemma}
\label{replacenets_in_intersection_U_incident}
Let $U$ be a monomial. Let $u^{(i_1)}$ and $u^{(i_2)}$ be virtual members of the chart of $U$. Assume $u^{(i_1)}$ starts from the left of the beginning of $u^{(i_2)}$. Let $u^{(i_1)}, a^{(i_1)}$, and $u^{(i_2)}, a^{(i_2)}$ be $U$-incident monomials. Assume $U = L^{(i_1)}u^{(i_1)}R^{(i_1)} = L^{(i_2)}u^{(i_2)}R^{(i_2)}$. Let $\widehat{u}^{(i_2)}$ be the intersection of $u^{(i_2)}$ and $R^{(i_1)}$, let $\dbhat{a}^{(i_2)}$ be defined by formula~\eqref{corresponding_replacement_in_intersection_right}. Let $\widehat{u}^{(i_1)}$ be the intersection of $u^{(i_1)}$ and $L^{(i_2)}$, let $\dbhat{a}^{(i_1)}$ be defined by formula~\eqref{corresponding_replacement_in_intersection_left}. Then the following properties hold.
\begin{enumerate}[label={(\arabic*)}]
\item
\label{replacenets_in_intersection_U_incident1}
$\left(L^{(i_1)}a^{(i_1)}\right)\cdot R^{(i_1)}[\widehat{u}^{(i_2)} \leadsto \dbhat{a}^{(i_2)}] = L^{(i_2)}[\widehat{u}^{(i_1)} \leadsto \dbhat{a}^{(i_1)}]\cdot \left(a^{(i_2)}R^{(i_2)}\right)$.
\item
\label{replacenets_in_intersection_U_incident2}
Let us denote the resulting monomial from statement~\ref{replacenets_in_intersection_U_incident1} by $W$.

If $a^{(i_1)}$ is a virtual member of the chart of $U[u^{(i_1)}\leadsto a^{(i_1)}]$ and $a^{(i_2)}$ is a virtual member of the chart of $U[u^{(i_2)}\leadsto a^{(i_2)}]$, then $W$ is a derived monomial of $U$ and $f(W) = f(U)$.

If $a^{(i_1)}$ is not a virtual member of the chart of $U[u^{(i_1)}\leadsto a^{(i_1)}]$ or $a^{(i_2)}$ is not a virtual member of the chart of $U[u^{(i_2)}\leadsto a^{(i_2)}]$, then $f(W) < f(U)$.
\end{enumerate}
\begin{proof}
Assume that $u^{(i_1)}$ and $u^{(i_2)}$ are separated, and $U = L^{(i_1)}u^{(i_1)}Mu^{(i_2)}R^{(i_2)}$. Then $\left(L^{(i_1)}a^{(i_1)}\right)\cdot R^{(i_1)}[\widehat{u}^{(i_2)} \leadsto \dbhat{a}^{(i_2)}]$ is the result of the cancellations in the monomial $L^{(i_1)}a^{(i_1)}\cdot \left(Ma^{(i_2)}R^{(i_2)}\right)$, where cancellations are performed, starting from the right parenthesis. Similarly, $L^{(i_2)}[\widehat{u}^{(i_1)} \leadsto \dbhat{a}^{(i_1)}]\cdot \left(a^{(i_2)}R^{(i_2)}\right)$ is the result of the cancellations in the monomial $\left(L^{(i_1)}a^{(i_1)}M\right)\cdot a^{(i_2)}R^{(i_2)}$, where cancellations are performed, starting from the left parenthesis. Obviously, the results after all the cancellations are the same.

Assume $u^{(i_1)}$ and $u^{(i_2)}$ are not separated. Let $e$ be the overlap of $u^{(i_1)}$ and $u^{(i_2)}$ ($e$ is empty if $u^{(i_1)}$ and $u^{(i_2)}$ touch at a point). Then we have
\begin{equation*}
\dbhat{a}^{(i_2)} = e^{-1}\cdot a^{(i_2)},\ \dbhat{a}^{(i_1)} = a^{(i_2)} \cdot e^{-1}.
\end{equation*}
Hence, we obtain
\begin{align*}
&\left(L^{(i_1)}a^{(i_1)}\right)\cdot R^{(i_1)}[\widehat{u}^{(i_2)} \leadsto \dbhat{a}^{(i_2)}] = \left(L^{(i_1)}a^{(i_1)}\right)\cdot e^{-1} \cdot a^{(i_2)}R^{(i_2)},\\
&L^{(i_2)}[\widehat{u}^{(i_1)} \leadsto \dbhat{a}^{(i_1)}]\cdot \left(a^{(i_2)}R^{(i_2)}\right) = L^{(i_1)}a^{(i_1)}\cdot e^{-1} \cdot\left( a^{(i_2)}R^{(i_2)}\right).
\end{align*}
Thus, we see that
\begin{equation*}
\left(L^{(i_1)}a^{(i_1)}\right)\cdot R^{(i_1)}[\widehat{u}^{(i_2)} \leadsto \dbhat{a}^{(i_2)}] = L^{(i_2)}[\widehat{u}^{(i_1)} \leadsto \dbhat{a}^{(i_1)}]\cdot \left(a^{(i_2)}R^{(i_2)}\right).
\end{equation*}
The first part of Lemma~\ref{replacenets_in_intersection_U_incident} is proved.

Assume $a^{(i_1)}$ is a virtual member of the chart of $U[u^{(i_1)}\leadsto a^{(i_1)}]$ and $a^{(i_2)}$ is a virtual member of the chart of $U[u^{(i_2)}\leadsto a^{(i_2)}]$. Let $\widetilde{u}^{(i_2)}$ be the image of $u^{(i_2)}$ in $U[u^{(i_1)}\leadsto a^{(i_1)}]$ ($\widetilde{u}^{(i_2)}$ be defined by formula~\eqref{image_of_right_U_incident}), let $\dbtilde{a}^{(i_2)}$ be defined by formula~\eqref{corresponding_element_U_incident_right}. Since $a^{(i_1)} \neq 1$, the monomial $L^{(i_1)}a^{(i_1)}R^{(i_1)}$ is reduced. Since $a^{(i_1)}$ is a virtual member of the chart of $U[u^{(i_1)} \leadsto a^{(i_1)}]$, it follows from statement~\ref{virtual_members_U_incident_replacements2} of Corollary~\ref{virtual_members_U_incident_replacements} that $\widetilde{u}^{(i_2)}$ and $\dbtilde{a}^{(i_2)}$ are $U[u^{(i_1)} \leadsto a^{(i_1)}]$-incident monomials. Therefore, it follows from Lemma~\ref{intersection_and_maximal_replacement} that $W = U[u^{(i_1)}\leadsto a^{(i_1)}][\widetilde{u}^{(i_2)}\leadsto \dbtilde{a}^{(i_2)}]$. Since $a^{(i_1)}$ is a virtual member of the chart of $U[u^{(i_1)}\leadsto a^{(i_1)}]$ and $a^{(i_2)}$ is a virtual member of the chart of $U[u^{(i_2)}\leadsto a^{(i_2)}]$, it follows from Lemma~\ref{replace_consecutively} that $U[u^{(i_1)}\leadsto a^{(i_1)}][\widetilde{u}^{(i_2)}\leadsto \dbtilde{a}^{(i_2)}]$ is a derived monomial of $U$ with the same $f$-characteristic as $U$.

Now assume that $a^{(i_1)}$ is a virtual member of the chart of $U[u^{(i_1)}\leadsto a^{(i_1)}]$, and $a^{(i_2)}$ is not a virtual member of the chart of $U[u^{(i_2)}\leadsto a^{(i_2)}]$. As above, using statement~\ref{virtual_members_U_incident_replacements2} of Corollary~\ref{virtual_members_U_incident_replacements} and  Lemma~\ref{intersection_and_maximal_replacement}, we obtain $W = U[u^{(i_1)}\leadsto a^{(i_1)}][\widetilde{u}^{(i_2)}\leadsto \dbtilde{a}^{(i_2)}]$, where $\widetilde{u}^{(i_2)}$ is the image of $u^{(i_2)}$ in $U[u^{(i_1)}\leadsto a^{(i_1)}]$ and $\dbtilde{a}^{(i_2)}$ is defined by formula~\eqref{corresponding_element_U_incident_right}. It follows from statement~\ref{virtual_members_U_incident_replacements3} of Corollary~\ref{virtual_members_U_incident_replacements} that $\dbtilde{a}^{(i_2)}$ is not a virtual member of the chart of $U[u^{(i_1)}\leadsto a^{(i_1)}][\widetilde{u}^{(i_2)}\leadsto \dbtilde{a}^{(i_2)}]$. Therefore,
\begin{equation*}
f(U) > f(U[u^{(i_1)}\leadsto a^{(i_1)}][\widetilde{u}^{(i_2)}\leadsto \dbtilde{a}^{(i_2)}]) = f(W).
\end{equation*}

The case when $a^{(i_1)}$ is a not a virtual member of the chart of $U[u^{(i_1)}\leadsto a^{(i_1)}]$ and $a^{(i_2)}$ is a virtual member of the chart of $U[u^{(i_2)}\leadsto a^{(i_2)}]$ is studied similarly.

Now assume that $a^{(i_1)}$ is  not a virtual member of the chart of $U[u^{(i_1)}\leadsto a^{(i_1)}]$ and $a^{(i_2)}$ is not a virtual member of the chart of $U[u^{(i_2)}\leadsto a^{(i_2)}]$. In particular, this means that $\SPM(a^{(i_1)}) < \tau$ and  $\SPM(a^{(i_2)}) < \tau$.

First let us consider the following particular case. Recall that $\widehat{u}^{(i_2)}$ is the intersection of $u^{(i_2)}$ and $R^{(i_1)}$. We assume that $u^{(i_2)}$ and $a^{(i_2)}$ are incident monomials, and $\SPM(\widehat{u}^{(i_2)}) \geqslant \tau$. At the end we will show how to reduce the general case to this particular situation.

First assume that $L^{(i_1)}a^{(i_1)}R^{(i_1)}$ is reduced. Let $\widetilde{u}^{(i_2)}$ be the image of $u^{(i_2)}$ in $L^{(i_1)}a^{(i_1)}R^{(i_1)}$ ($\widetilde{u}^{(i_2)}$ be defined by formula~\eqref{image_of_right_general}). Therefore, since $\SPM(\widehat{u}^{(i_2)}) \geqslant \tau$, we see that $\SPM(\widetilde{u}^{(i_2)}) \geqslant \tau$. That is, $\widetilde{u}^{(i_2)}$ is a virtual member of the chart of $L^{(i_1)}a^{(i_1)}R^{(i_1)}$.

Let $\dbtilde{a}^{(i_2)}$ be defined by formula~\eqref{corresponding_element_right_general}. Since $u^{(i_2)}$ and $a^{(i_2)}$ are incident monomials and $\widehat{u}^{(i_2)}$ is not a small piece, it follows from Lemma~\ref{compatibility_for_incedent} that $\widetilde{u}^{(i_2)}$ and $\dbtilde{a}^{(i_2)}$ are incident monomials. Therefore, since $\widetilde{u}^{(i_2)}$ is a virtual member of the chart of $L^{(i_1)}a^{(i_1)}R^{(i_1)}$, it follows from Lemma~\ref{estimation_value_property} that the replacement $\widetilde{u}^{(i_2)} \mapsto \dbtilde{a}^{(i_2)}$ in $L^{(i_1)}a^{(i_1)}R^{(i_1)}$ does not increase $f$-characteristic. That is,
\begin{equation*}
f(U[u^{(i_1)} \leadsto a^{(i_1)}][\widetilde{u}^{(i_2)} \mapsto \dbtilde{a}^{(i_2)}]) \leqslant f(L^{(i_1)}a^{(i_1)}R^{(i_1)}).
\end{equation*}
However, since $u^{(i_1)}$ is a virtual member of the chart of $U$ and $a^{(i_1)}$ is not a virtual member of the chart of $L^{(i_1)}a^{(i_1)}R^{(i_1)}$, we see that $f(L^{(i_1)}a^{(i_1)}R^{(i_1)}) < f(U)$. Since the monomial $L^{(i_1)}a^{(i_1)}R^{(i_1)}$ is reduced and $\widetilde{u}^{(i_2)}$ and $\dbtilde{a}^{(i_2)}$ are incident monomials, by Lemma~\ref{intersection_and_maximal_replacement}, we have $W = U[u^{(i_1)} \leadsto a^{(i_1)}][\widetilde{u}^{(i_2)} \mapsto \dbtilde{a}^{(i_2)}]$. Hence,
\begin{equation*}
f(W) = f(U[u^{(i_1)} \leadsto a^{(i_1)}][\widetilde{u}^{(i_2)} \mapsto \dbtilde{a}^{(i_2)}]) \leqslant f(L^{(i_1)}a^{(i_1)}R^{(i_1)}) < f(U).
\end{equation*}

Now assume that $a^{(i_1)} = 1$ and $L^{(i_1)}\cdot R^{(i_1)}$ is not necessarily reduced. Notice that the replacement $\widehat{u}^{(i_2)} \mapsto \dbhat{a}^{(i_2)}$ in $R^{(i_1)}$ is a replacement of an element of $\longmo{R^{(i_1)}} \subseteq \nfc{R^{(i_1)}}$. Hence, similarly to the last part of Lemma~\ref{minimal_coverings_property}, one can prove that $\mincov{W} < \mincov{U}$. Therefore, $f(W) < f(U)$.

Now let us show the reduction of the general case to the particular case considered above. Initially we have the replacements $u^{(i_1)} \leadsto a^{(i_1)}$ and $u^{(i_2)} \leadsto a^{(i_2)}$ in $U$. Since $u^{(i_2)}$ and $a^{(i_2)}$ are $U$-incident monomials, there exists a sequence of elements $m_1^{(i_2)}, \ldots, m_{s + 1}^{(i_2)}$ that satisfies the conditions of Definition~\ref{U_incident_monomials}. Assume $\SPM(m_s^{(i_2)}) < \tau + 1$. Since $m_s^{(i_2)}$ and $a^{(i_2)}$ are incident monomials and $\SPM(a^{(i_2)}) \leqslant \tau$, it follows from Small Cancellation Axiom that there exists a monomial $m^{\prime} \in \Mon$ such that $m^{\prime}$ is incident to $m_s^{(i_2)}$ and $a^{(i_2)}$, and wherein $\SPM(m^{\prime}) \geqslant \tau + 1$. We can insert this $m^{\prime}$ to the sequence $m_1^{(i_2)}, \ldots, m_{s + 1}^{(i_2)}$ before the last monomial. Therefore, without loss of generality in what follows we assume that $\SPM(m_{s}^{(i_2)}) \geqslant \tau + 1$.

Roughly speaking, instead of two initial replacements in $U$ we consider three replacements: $u^{(i_1)} \leadsto a^{(i_1)}$, $u^{(i_2)} \leadsto m_s^{(i_2)}$, $m_s^{(i_2)} \mapsto a^{(i_2)}$. We perform the replacement $u^{(i_2)} \leadsto m_s^{(i_2)}$ in $U$. Let us denote the resulting monomial $ U[u^{(i_2)} \leadsto m_s^{(i_2)}]$ by $Z$. Since $m_s^{(i_2)}$ is a virtual member of the chart of $Z$, we have $f(Z) = f(U)$. We define the replacements in $Z$ that correspond to $u^{(i_1)} \leadsto a^{(i_1)}$ and $m_s^{(i_2)} \mapsto a^{(i_2)}$. On the one hand, if we perform them, using the process defined above, then $W$ is the resulting monomial. On the other hand, they satisfy the conditions of the particular case considered above. Hence, $f(W) < f(Z)$. Therefore, finally we obtain $f(W) < f(Z) = f(U)$.

Let us explain the situation in more detail. Let $v^{(i_1)}$ be the image of $u^{(i_1)}$ in $U[u^{(i_2)} \leadsto m_s^{(i_2)}] = Z$. We have the replacement $m_s^{(i_2)} \mapsto a^{(i_2)}$ in $Z$ (without any changes of $a^{(i_2)}$). Let $v^{(i_1)} \leadsto b^{(i_1)}$ be the replacement in $Z$ that corresponds to the replacement $u^{(i_1)} \leadsto a^{(i_1)}$ in $U$. Let $Z = \widetilde{L}^{(i_1)}v^{(i_1)}\widetilde{R}^{(i_1)}$. Clearly, $\widetilde{L}^{(i_1)} = L^{(i_1)}$, so, $Z = L^{(i_1)}v^{(i_1)}\widetilde{R}^{(i_1)}$. Let $\widehat{m}_s^{(i_2)}$ be the intersection of $m_s^{(i_2)}$ and $\widetilde{R}^{(i_1)}$. Let $\widehat{m}_s^{(i_2)} \mapsto \dbhat{b}^{(i_2)}$ be the replacement in $\widetilde{R}^{(i_1)}$ that corresponds to the replacement $m_s^{(i_2)} \mapsto a^{(i_2)}$ in $Z$. First let us show that
\begin{equation*}
\left(L^{(i_1)}b^{(i_1)}\right) \cdot \widetilde{R}^{(i_1)}[\widehat{m}_s^{(i_2)} \mapsto \dbhat{b}^{(i_2)}] = W.
\end{equation*}

Assume $u^{(i_1)}$ and $u^{(i_2)}$ are separated. Then $U = L^{(i_1)}u^{(i_1)}Mu^{(i_2)}R^{(i_2)}$, and we have
\begin{equation*}
Z = U[u^{(i_2)} \leadsto m_s^{(i_2)}] = L^{(i_1)}u^{(i_1)}Mm_s^{(i_2)}R^{(i_2)}.
\end{equation*}
We have $v^{(i_1)} = u^{(i_1)}$, $b^{(i_1)} = a^{(i_1)}$, $\widehat{m}_s^{(i_2)} = m_s^{(i_2)}$, $\dbhat{b}^{(i_2)} = a^{(i_2)}$. Therefore,
\begin{equation*}
\left(L^{(i_1)}b^{(i_1)}\right) \cdot \widetilde{R}^{(i_1)}[\widehat{m}_s^{(i_2)} \mapsto \dbhat{b}^{(i_2)}] = L^{(i_1)}a^{(i_1)}\cdot \left(Ma^{(i_2)}R^{(i_2)}\right) = W.
\end{equation*}

Assume $u^{(i_1)}$ and $u^{(i_2)}$ are not separated. Then $v^{(i_1)}$ and $m_s^{(i_2)}$ are not separated as well. Assume $u^{(i_1)}$ and $u^{(i_2)}$ have an overlap $e$ ($e$ is empty if $u^{(i_1)}$ and $u^{(i_2)}$ touch at a point), $v^{(i_1)}$ and $m_s^{(i_2)}$ have an overlap $d$ ($d$ is empty if $v^{(i_1)}$ and $m_s^{(i_2)}$ touch at a point). Then we see that
\begin{align*}
&v^{(i_1)} = \left(u^{(i_1)}\cdot e^{-1}\right)d,\ b^{(i_1)} = a^{(i_1)}\cdot e^{-1}\cdot d,\\
&\widehat{m}_s^{(i_2)} = d^{-1}\cdot m_s^{(i_2)},\ \dbhat{b}^{(i_2)} = d^{-1}\cdot a^{(i_2)}.
\end{align*}
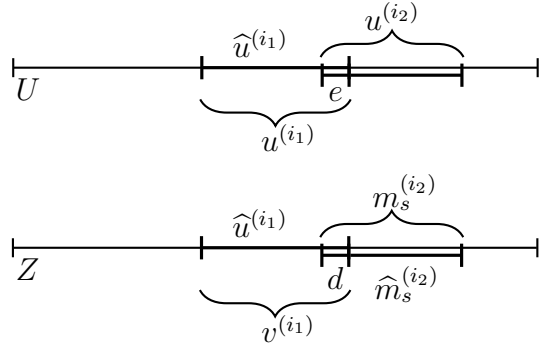
\begin{center}
\begin{tikzpicture}
\draw[|-|, black, thick] (-1.5, 0)--(5.5, 0);
\node[below, xshift=6] at (-1.5, 0) {$U$};
\draw[|-|, black, very thick] (1,0)--(3,0);
\path (1, 0)--(2.6, 0) node[midway, above] {$\widehat{u}^{(i_1)}$};
\draw[-|, black, very thick] (3,-0.1)--(4.5,-0.1);
\draw[|-|, black, very thick] (2.6,-0.1)--(3,-0.1) node[midway, below] {$e$};
\draw [thick, decorate, decoration={brace, amplitude=10pt, raise=4pt}] (2.6, 0) to node[midway, above, yshift=10pt] {$u^{(i_2)}$} (4.5, 0);
\draw [thick, decorate, decoration={brace, amplitude=10pt, raise=12pt, mirror}] (1, 0) to node[midway, above, yshift=-36pt, xshift=5] {$u^{(i_1)}$} (3, 0);
\end{tikzpicture}

\begin{tikzpicture}
\draw[|-|, black, thick] (-1.5, 0)--(5.5, 0);
\node[below, xshift=6] at (-1.5, 0) {$Z$};
\draw[|-|, black, very thick] (1,0)--(3,0);
\path (1, 0)--(2.6, 0) node[midway, above] {$\widehat{u}^{(i_1)}$};
\draw[-|, black, very thick] (3,-0.1)--(4.5,-0.1);
\draw[|-|, black, very thick] (2.6,-0.1)--(3,-0.1) node[midway, below] {$d$};
\path (3, 0)--(4.5, 0) node[midway, below, yshift=-3] {$\widehat{m}_s^{(i_2)}$};
\draw [thick, decorate, decoration={brace, amplitude=10pt, raise=4pt}] (2.6, 0) to node[midway, above, yshift=10pt, xshift=5] {$m_s^{(i_2)}$} (4.5, 0);
\draw [thick, decorate, decoration={brace, amplitude=10pt, raise=16pt, mirror}] (1, 0) to node[midway, above, yshift=-40pt, xshift=5] {$v^{(i_1)}$} (3, 0);
\end{tikzpicture}
\end{center}
Clearly, $Z$ can be written as $Z = L^{(i_1)}v^{(i_1)}\widehat{m}_s^{(i_2)}R^{(i_2)}$. Hence, we obtain
\begin{align*}
&\left(L^{(i_1)}b^{(i_1)}\right) \cdot \widetilde{R}^{(i_1)}[\widehat{m}_s^{(i_2)} \mapsto \dbhat{b}^{(i_2)}] = \left(L^{(i_1)}a^{(i_1)}\cdot e^{-1}\cdot d\right)\cdot d^{-1}\cdot a^{(i_2)} R^{(i_2)} =\\
&= L^{(i_1)}a^{(i_1)}\cdot e^{-1}\cdot a^{(i_2)} R^{(i_2)} = \left(L^{(i_1)}a^{(i_1)}\right) \cdot \dbhat{a}^{(i_2)}R^{(i_2)} =\\
&= \left(L^{(i_1)}a^{(i_1)}\right) \cdot R^{(i_1)}[\widehat{u}^{(i_1)} \leadsto \dbhat{a}^{(i_2)}] = W.
\end{align*}

Let us show that the replacements $m_{s}^{(i_2)} \mapsto a^{(i_2)}$ and $v^{(i_1)} \leadsto b^{(i_1)}$ in $Z$ satisfy the conditions of the particular case that we studied above.

Since $a^{(i_1)}$ is not a virtual member of the chart of $U[u^{(i_1)} \leadsto a^{(i_1)}]$ and $m_s^{(i_2)}$ is a virtual member of the chart of $U[u^{(i_2)} \leadsto m_s^{(i_2)}] = Z$, it follows from statement~\ref{virtual_members_U_incident_replacements6} of Corollary~\ref{virtual_members_U_incident_replacements} that $b^{(i_1)}$ is not a virtual member of the chart of $U[u^{(i_2)} \leadsto m_s^{(i_2)}][v^{(i_1)} \leadsto b^{(i_1)}] = Z[v^{(i_1)} \leadsto b^{(i_1)}]$.

Since $Z = U[u^{(i_2)} \leadsto m_s^{(i_2)}]$, we obviously have $U[u^{(i_2)} \leadsto a^{(i_2)}] = Z[m_{s}^{(i_2)} \mapsto a^{(i_2)}]$. So, since $a^{(i_2)}$ is not a virtual member of the chart of $U[u^{(i_2)} \leadsto a^{(i_2)}]$, we obviously obtain that $a^{(i_2)}$ is not a virtual member of the chart of $Z[m_{s}^{(i_2)} \mapsto a^{(i_2)}]$.

Since $\SPM(m_s^{(i_2)}) \geqslant \tau + 1$, we see that $\SPM(\widehat{m}_s^{(i_2)}) \geqslant \tau$.

Therefore, all the above conditions are satisfied. This completes the proof.
\end{proof}
\end{lemma}

\section{The structure of $\Qalg$ as a vector space: filtration, grading and tensor products}
\label{structure_calc}

\subsection{The filtration on the space $\Galg$}
\label{filtration_def_subsection}
We define an increasing filtration on $\Galg$, using $f$-characteristics of monomials. Consider the following correspondence $t$ between $\lbrace 0\rbrace \cup \mathbb{N}$ and values of $f$-characteristics of monomials. We put $t(0) = (0, 0)$. Assume $t(n) = (r, s)$, then we put
\begin{equation*}
t(n + 1) = \begin{cases}
(r, s + 1) &\textit{ if } r > s,\\
(r + 1, 0) &\textit{ if } r = s.
\end{cases}
\end{equation*}

We define an increasing filtration on $\Galg$ by the following rule:
\begin{equation}
\label{filtration_def}
\Ft_n(\Galg) = \left\langle\lbrace Z \mid Z \in \Fr, f(Z) \leqslant t(n)\rbrace\right\rangle,
\end{equation}
where $\langle \cdot \rangle$ is a linear span. That is, $\Ft_n(\Galg)$ is generated by all monomials with $f$-characteristics not greater than $t(n)$. Since, by definition, $t(n - 1) < t(n)$, we see that $\Ft_{n - 1}(\Galg) \subseteq \Ft_n(\Galg)$.

Since $\mincov{Z} \geqslant \nvirt{Z}$ for every monomial $Z$, we see that $t(n)$ covers all possible values of $f$-characteristics of monomials while $n$ varies over all values from $\lbrace 0\rbrace \cup \mathbb{N}$. Therefore, we really have
\begin{equation*}
\Galg = \bigcup_{n = 0}^{\infty} \Ft_n(\Galg).
\end{equation*}

Assume $Z$ is a monomial that belongs to $\Ft_n(\Galg)$. By Lemma~\ref{estimation_value_property}, $f$-characteristics of derived monomials of $Z$ is not greater that $f(Z)$. Hence, all derived monomials of $Z$ belong to $\Ft_n(\Galg)$. That is, $\Ft_n(\Galg)$ is closed under taking derived monomials.

\begin{definition}
\label{derived_monomials_space_def}
Let $U$ be a monomial. By $\DSpace{U}$, we denote a linear subspace of $\Galg$ generated by all the derived monomials of $U$.
\end{definition}

We consider derived monomials of $U$ and distinguish two sets of monomials:
\begin{align}
\label{upper_monomials}
\nonumber
\DMUp{U} = \lbrace Z &\mid Z \textit{ is a derived monomial of } U \\
&\textit{ such that } f(Z) = f(U) \rbrace;\\
\label{lower_monomials}
\nonumber
\DMLow{U} = \lbrace Z &\mid Z \textit{ is a derived monomial of } U \\
&\textit{ such that } f(Z) < f(U) \rbrace.
\end{align}
Notice that the set $\DMUp{U}$ is always non-empty, since $U \in \DMUp{U}$. On the other hand, the set $\DMLow{U}$ may be empty. Clearly, the set of all derived monomials of $U$ is equal to $\DMUp{U} \cup \DMLow{U}$. Therefore, we have
\begin{equation*}
\DSpace{U} = \left\langle \DMUp{U} \sqcup \DMLow{U} \right\rangle.
\end{equation*}
\medskip

\begin{definition}
\label{low_space_def}
Let $Y$ be a linear subspace of $\Galg$ generated by monomials, closed under taking derived monomials and such that $f$-characteristics of its monomials is bounded. Then, evidently, there exists a unique number $n$ such that $Y \subseteq \Ft_n(\Galg)$ and $Y \nsubseteq \Ft_{n - 1}(\Galg)$. We put
\begin{equation*}
\Low(Y) = Y \cap \Ft_{n - 1}(\Galg).
\end{equation*}
If $Y \subseteq \Ft_0(\Galg)$, then we put $\Low(Y) = 0$.
\medskip

Since both spaces $Y$ and $\Ft_{n - 1}(\Galg)$ are generated by monomials and closed under taking derived monomials, we obtain that $\Low(Y)$ is generated by monomials and closed under taking derived monomials as well. Moreover, if $Y \neq 0$, then $Y \neq \Low(Y)$.
\end{definition}

\begin{lemma}
\label{intersection_with_low_equality}
Let $U$ be a monomial. Assume $U \in \Ft_{n}(\Galg) \setminus \Ft_{n - 1}(\Galg)$. If $\DMLow{U}$ is non empty, then
\begin{equation*}
\Low\DSpace{U} = \langle \DMLow{U}\rangle.
\end{equation*}
If $\DMLow{U}$ is empty, then $\Low\DSpace{U} = 0$.
\end{lemma}
\begin{proof}
Since $U \in \Ft_{n}(\Galg)$, every monomial with $f$-characteristic smaller than $f(U)$ is contained in $\Ft_{n - 1}(\Galg)$. Therefore, $\langle \DMLow{U}\rangle \subseteq \Ft_{n - 1}(U)$. So, we obviously have
\begin{equation*}
\langle \DMLow{U}\rangle \subseteq \Low\DSpace{U} = \Ft_{n - 1}(\Galg) \cap \DSpace{U}.
\end{equation*}

Consider some arbitrary element $\sum_{i = 1}^l\alpha_i Z_i \in \Ft_{n - 1}(\Galg) \cap \DSpace{U}$, where $Z_i$ are monomials. Here we mean that the sum is additively reduced. Since $\Ft_{n - 1}(\Galg)$ is generated by monomials, every $Z_i$ belongs to $\Ft_{n - 1}(\Galg)$. Similarly, since $\DSpace{U}$ is generated by monomials, every $Z_i$ belongs to $\DSpace{U}$. So, on the one hand, every $Z_i$ is a derived monomial of $U$. On the other hand, since $Z_i \in \Ft_{n - 1}(\Galg)$, we have $f(Z_i) \leqslant t(n - 1) < t(n)$. Since $U \in \Ft_{n}(\Galg) \setminus \Ft_{n - 1}(\Galg)$, we see that $f(U) = t(n)$. Therefore, $f(Z_i) < f(U)$. By definition, this means that $Z_i \in \langle \DMLow{U}\rangle$, so, $\sum_{i = 1}^l\alpha_i Z_i \in  \langle \DMLow{U}\rangle$. Thus, $\Ft_{n - 1}(\Galg) \cap \DSpace{U} \subseteq  \langle \DMLow{U}\rangle$. Lemma~\ref{intersection_with_low_equality} is proved.
\end{proof}

In the sequel, we will widely use the following simple lemma.
\begin{lemma}
\label{derived_spaces_equality}
Let $U$ be a monomial, $Z$ be a derived monomial of $U$, $Y \subseteq \Galg$ be a linear space generated by monomials and closed under taking derived monomials. Then the following statements hold:
\begin{enumerate}
\item[$(1)$]
$Z \in \DMUp{U}$ if and only if $\DSpace{Z} = \DSpace{U}$.
\item[$(2)$]
If $Z \in Y$ and $Z \in \DMUp{U}$, then $\DSpace{U} \subseteq Y$.
\end{enumerate}
\end{lemma}
\begin{proof}
Let us prove statement $(1)$. Assume $Z \in \DMUp{U}$ and let us show that $\DSpace{Z} = \DSpace{U}$.

Since $Z$ is a derived monomial of $U$, clearly, every derived monomial of $Z$ is a derived monomial of $U$ as well. Hence, $\DSpace{Z} \subseteq \DSpace{U}$.

Since $Z$ is a derived monomial of $U$, there exists a sequence of transformations $\phi_1, \ldots, \phi_{s}$ of type~\ref{repl1}, \ref{repl2} (see Definition~\ref{derived_monomials}) such that
\begin{equation*}
U \overset{\phi_1}{\longmapsto} \ldots \overset{\phi_{s}}{\longmapsto} Z.
\end{equation*}
Moreover, Lemma~\ref{estimation_value_property} implies that every transformation $\phi_l$ preserves $f$-characteristic of the corresponding monomial. That is, every transformation $\phi_l$ is of the form $L^{(l)}a^{(l)}_hR^{(l)} \overset{\phi_{l}}{\longmapsto} L^{(l)}a^{(l)}_jR^{(l)}$, where $a_h^{(l)}$ is a virtual member of the chart of $L^{(l)}a^{(l)}_hR^{(l)}$ and $a_j^{(l)}$ is a virtual member of the chart of $L^{(l)}a^{(l)}_jR^{(l)}$, $a_h^{(l)}$ and $a_j^{(l)}$ are incident monomials. In particular, every transformation $\phi_l$ is of type~\ref{repl1}. Then the opposite replacement $L^{(l)}a^{(l)}_jR^{(l)} \overset{\phi_{l}}{\longmapsto} L^{(l)}a^{(l)}_hR^{(l)}$ is a transformation of type~\ref{repl1} as well. Denote it by $\phi^{-1}_l$. Then we obtain the following sequence of transformations of type~\ref{repl1}:
\begin{equation*}
Z \overset{\phi_s^{-1}}{\longmapsto} \ldots \overset{\phi_1^{-1}}{\longmapsto} U.
\end{equation*}
Therefore, by definition, $U$ is a derived monomial of $Z$. Hence, $\DSpace{U} \subseteq \DSpace{Z}$. Thus, $\DSpace{U} = \DSpace{Z}$.

Assume $\DSpace{Z} = \DSpace{U}$, let us show that $Z \in \DMUp{U}$. We see that $Z$ is a derived monomial of $U$, hence, $f(Z) \leqslant f(U)$. Similarly, $U$ is a derived monomial of $Z$, hence, $f(U) \leqslant f(Z)$. So, $f(U) = f(Z)$. By definition, this means that $Z \in \DMUp{U}$.

Let us prove statement $(2)$. Since $Z \in Y$ and $Y$ is closed under taking derived monomials, we obtain $\DSpace{Z} \subseteq Y$. However, since $Z \in \DMUp{U}$, it follows from statement~$(1)$ that $\DSpace{Z} = \DSpace{U}$. Therefore, $\DSpace{U} \subseteq Y$. This completes the proof.
\end{proof}

\begin{definition}
\label{subspace_of_dependencies}
Suppose $Y$ is a subspace of $\Galg$ linearly generated by a set of monomials and closed under taking derived monomials.
\begin{multline*}
\textit{We denote the set of all the layouts of multi-turns }\\
\textit{of virtual members of the chart of monomials of } Y \textit{ by } \GDp^{\prime}(Y).
\end{multline*}
\label{gdp_y_def}

Let $Z$ be a monomial. A layout a multi-turn of a virtual member of the chart of $Z$ is, in fact, a linear combination of derived monomials of $Z$. Therefore, if $Z \in Y$, a layout of a multi-turn of a virtual member of the chart of $Z$ belongs to $Y$, because $Y$ is closed under taking derived monomials. Hence, $\GDp^{\prime}(Y) \subseteq Y$. So, we can consider a subspace of $Y$ linearly generated by $\GDp^{\prime}(Y)$. We call this space \emph{the subspace of dependencies on $Y$} and denote it by $\Dp(Y)$. That is,
\begin{equation*}
\Dp(Y) = \langle \GDp^{\prime}(Y)\rangle.
\end{equation*}

Note that if a monomial $Z \in Y$ has no virtual members of the chart, then there are no multi-turns of virtual members of the chart of~$Z$. So, if $Y$ consists only of monomials with no virtual members of the chart, then $\GDp^{\prime}(Y) = \varnothing$. In this case, by definition, we put $\Dp(Y) = 0$.
\end{definition}

Using the above notions, we see that $\GDp^{\prime} = \GDp^{\prime}(\Galg)$ (see~\eqref{linear_dep_members2} on page~\pageref{linear_dep_members2}) and
\begin{equation*}
\Ideal = \langle \GDp^{\prime}\rangle = \langle \GDp^{\prime}(\Galg) \rangle = \Dp(\Galg)
\end{equation*}
(see~Proposition~\ref{the_ideal_characterisation2}).

\begin{remark}
Let $Y_1$ and $Y_2$ be subspaces of $\Galg$ linearly generated by a set of monomials and closed under taking derived monomials. Assume $Y_1 \subseteq Y_2$. Then it follows directly from Definition~\ref{subspace_of_dependencies} that $\GDp^{\prime}(Y_1) \subseteq \GDp^{\prime}(Y_2)$ and $\Dp(Y_1) \subseteq \Dp(Y_2)$.
\end{remark}

\begin{lemma}
\label{layouts_subspace}
Suppose $Y_1$ and $Y_2$ are subspaces of $\Galg$ linearly generated by monomials and closed under taking derived monomials, $Y_1 \subseteq Y_2$. Then $Y_1 \cap \GDp^{\prime}(Y_2) = \GDp^{\prime}(Y_1)$. In particular, $Y \cap \GDp^{\prime} = Y \cap \GDp^{\prime} (\Galg)= \GDp^{\prime}(Y)$.
\end{lemma}
\begin{proof}
Let $T \in \GDp^{\prime}(Y_1)$. Then $T \in Y_1$ because $Y_1$ is closed under taking derived monomials. On the other hand, since $Y_1 \subseteq Y_2$, $\GDp^{\prime}(Y_1) \subseteq \GDp^{\prime}(Y_2)$. Hence, $T \in \GDp^{\prime}(Y_2)$. So, $T \in Y_1 \cap \GDp^{\prime}(Y_2)$. That is, $\GDp^{\prime}(Y_1) \subseteq Y_1 \cap \GDp^{\prime}(Y_2)$.

Let $T \in Y_1 \cap \GDp^{\prime}(Y_2)$. Let $Z \in Y_2$ be a monomial such that $T$ is a layout of a multi-turn of a virtual member of the chart of $Z$. Since $T \in Y_1$ and $Y_1$ is generated by monomials, all the monomials of $T$ belong to $Y_1$. In particular, $Z \in Y_1$. Thus, $T \in \GDp^{\prime}(Y_1)$. That is, $Y_1 \cap \GDp^{\prime}(Y_2) \subseteq \GDp^{\prime}(Y_1)$.
\end{proof}

\begin{lemma}[\textbf{Main Lemma}]
\label{fall_through_linear_dep_one}
Let $U$ be an arbitrary monomial, $U \in \Ft_n(\Galg) \setminus \Ft_{n - 1}(\Galg)$. Then
\begin{equation}
\label{main_cancel_linear_dep}
\Dp\DSpace{U} \cap \Low\DSpace{U} \subseteq \Dp(\Ft_{n - 1}(\Galg)).
\end{equation}
\end{lemma}
We will prove this lemma in Subsection~\ref{tensor_product_section}. First let us apply it in order to prove the following key statement.

\begin{proposition}
\label{fall_through_linear_dep}
We have
\begin{equation}
\label{cancel_linear_dep}
\Dp(\Ft_n(\Galg)) \cap \Ft_{n - 1}(\Galg) = \Dp(\Ft_{n - 1}(\Galg)).
\end{equation}
\end{proposition}
\begin{proof}
We assume that Lemma~\ref{fall_through_linear_dep_one} is proved. Let us prove Proposition~\ref{fall_through_linear_dep}.

Since $\Ft_n(\Galg) \supseteq \Ft_{n - 1}(\Galg)$, we evidently have $\Dp(\Ft_n(\Galg)) \cap \Ft_{n - 1}(\Galg) \supseteq \Dp(\Ft_{n - 1}(\Galg))$. Hence, we need to show that
\begin{equation*}
\Dp(\Ft_n(\Galg)) \cap \Ft_{n - 1}(\Galg) \subseteq \Dp(\Ft_{n - 1}(\Galg)).
\end{equation*}

Clearly, we have
\begin{align*}
&\Dp(\Ft_n(\Galg)) = \Dp(\Ft_{n - 1}(\Galg)) + Y^{(n)},\\
&\textit{where } Y^{(n)} = \langle\lbrace T \mid T \in \GDp^{\prime}(\Ft_n(\Galg)), T \notin \Ft_{n - 1}(\Galg)\rbrace\rangle.
\end{align*}
That is, the space $Y^{(n)}$ is linearly generated by all the layouts of multi-turns of virtual members of the chart of all the monomials of $\Ft_n(\Galg) \setminus \Ft_{n - 1}(\Galg)$. Since $\Dp(\Ft_{n - 1}(\Galg)) \subseteq \Ft_{n - 1}(\Galg)$, we have
\begin{align*}
\Dp(\Ft_n(\Galg))\cap \Ft_{n - 1}(\Galg) &= (\Dp(\Ft_{n - 1}(\Galg)) + Y^{(n)}) \cap \Ft_{n - 1}(\Galg) =\\
&= \Dp(\Ft_{n - 1}(\Galg)) + Y^{(n)} \cap \Ft_{n - 1}(\Galg).
\end{align*}
So, we only need to prove that
\begin{equation*}
Y^{(n)} \cap \Ft_{n - 1}(\Galg) \subseteq \Dp(\Ft_{n - 1}(\Galg)).
\end{equation*}

Let $W \in Y^{(n)} \cap \Ft_{n - 1}(\Galg)$ be a non-zero element. Then
\begin{equation*}
W = \sum\limits_{i = 1}^l \gamma_i T_i,\ \gamma_i \neq 0,
\end{equation*}
where every $T_i \in \GDp^{\prime}(\Ft_n(\Galg))$ and $T_i \notin \Ft_{n - 1}(\Galg)$, and
\begin{equation}
\label{dep_sum_initial}
\sum\limits_{i = 1}^l \gamma_i T_i \in \Ft_{n - 1}(\Galg).
\end{equation}
We will prove that $W \in \Dp(\Ft_{n - 1}(\Galg))$.

Since $\Ft_{n - 1}(\Galg)$ is generated by monomials and the the set of all monomials is a basis of $\Galg$, every additively reduced linear combination of monomials that belongs to $\Ft_{n - 1}(\Galg)$ contains only monomials from $\Ft_{n - 1}(\Galg)$. Hence, all monomials of $\Ft_n(\Galg) \setminus \Ft_{n - 1}(\Galg)$ have to cancel out in~\eqref{dep_sum_initial}.

Let $Z_i \in \Ft_n(\Galg) \setminus \Ft_{n - 1}(\Galg)$ be a monomial such that $T_i$ is a layout of a multi-turn of a virtual member of the chart of $Z_i$, $i = 1, \ldots, l$. Then, obviously, $T_i \in \Dp\DSpace{Z_i}$, $i = 1, \ldots, l$. Consider the space $\sum_{i = 1}^l \DSpace{Z_i}$. If $\DSpace{Z_{i^{\prime}}} \subseteq \DSpace{Z_{i^{\prime\prime}}}$ for $i^{\prime} \neq i^{\prime\prime}$, then, clearly,
\begin{equation*}
\sum_{i = 1}^l \DSpace{Z_i} = \sum_{\substack{i = 1\\i \neq i^{\prime}}}^l \DSpace{Z_i}.
\end{equation*}
Hence, we can choose a subset $\lbrace Z_{i_1}, \ldots, Z_{i_m}\rbrace \subseteq \lbrace Z_1, \ldots, Z_l\rbrace$ such that
\begin{equation*}
\sum_{i = 1}^l \DSpace{Z_i} = \sum_{j = 1}^m \DSpace{Z_{i_j}} \textit{ and } \DSpace{Z_{i_j}} \nsubseteq \DSpace{Z_{i_{j^{\prime}}}} \textit{ whenever } j \neq j^{\prime}.
\end{equation*}
Then every $Z_i$, $i = 1, \ldots, l$, is a derived monomial of some $Z_{i_j}$, $j = 1, \ldots, m$, and $\DSpace{Z_i} \subseteq \DSpace{Z_{i_j}}$.

Since $T_i \in \Dp\DSpace{Z_i}$, $T_i$ belongs to some $\Dp\DSpace{Z_{i_j}}$, $j = 1, \ldots, m$. So, the layouts $\lbrace T_1, \ldots, T_l\rbrace$ can be separated into $m$ groups $\lbrace T_1^j, \ldots, T_{l_j}^j\rbrace$ such that $\lbrace T_1^j, \ldots, T_{l_j}^j\rbrace \subseteq \Dp\DSpace{Z_{i_j}}$. So, we have
\begin{equation*}
\lbrace T_1, \ldots, T_l\rbrace = \bigsqcup_{j = 1}^m\lbrace T_1^j, \ldots, T_{l_j}^j\rbrace.
\end{equation*}
Hence,
\begin{equation}
\label{dep_sum}
\sum\limits_{i = 1}^l \gamma_i T_i = \sum\limits_{j = 1}^m \sum\limits_{t = 1}^{l_j} \gamma_t^j T_t^j \in \Ft_{n - 1}(\Galg).
\end{equation}

Clearly, every monomial in the sum~\eqref{dep_sum} is a derived monomial of $Z_{i_j}$ for some $j = 1, \dots, m$. Recall that $Z_{i_j} \in \Ft_n(\Galg) \setminus \Ft_{n - 1}(\Galg)$. Hence, by definition, $f$-characteristics of monomials of $\DMLow{Z_{i_j}}$ are smaller than $t(n)$. Therefore, we see that
\begin{equation*}
\DMLow{Z_{i_j}} \subseteq \Ft_{n - 1}(\Galg).
\end{equation*}
So, if a monomial in the sum~\eqref{dep_sum} belongs to $\Ft_n(\Galg) \setminus \Ft_{n - 1}(\Galg)$, then it belongs to $\DMUp{Z_{i_j}}$ for some $j = 1, \dots, m$.

Let $Z$ be a monomial and $Z \in \DMUp{Z_{i_j}}$. Assume $Z$ belongs to $\Ft_{n - 1}(\Galg)$. Then it follows from Lemma~\ref{derived_spaces_equality} that $\DSpace{Z_{i_j}} \subseteq \Ft_{n - 1}(\Galg)$. This contradicts the assumption that $Z_{i_j} \in \Ft_n(\Galg) \setminus \Ft_{n - 1}(\Galg)$. Hence, $Z \in \Ft_n(\Galg) \setminus \Ft_{n - 1}(\Galg)$.

Combining the above results, we see that any monomial in the sum~\eqref{dep_sum} belongs to $\Ft_n(\Galg) \setminus \Ft_{n - 1}(\Galg)$ if and only if it belongs to $\DMUp{Z_{i_j}}$ for some $j = 1, \ldots, m$. Hence, the monomials of $\DMUp{Z_{i_j}}$, $j = 1, \ldots, m$, have to cancel out in the sum~\eqref{dep_sum}.

Again assume that a monomial $Z \in \DMUp{Z_{i_j}}$. If $Z$ belongs to $\DSpace{Z_{i_{j^{\prime}}}}$ for some $j^{\prime} \neq j$, then it follows from Lemma~\ref{derived_spaces_equality} that $\DSpace{Z_{i_j}} \subseteq \DSpace{Z_{i_{j^{\prime}}}}$. However, we have chosen the monomials $Z_{i_1}, \ldots, Z_{i_m}$ with the property $\DSpace{Z_{i_j}} \nsubseteq \DSpace{Z_{i_{j^{\prime}}}}$ whenever $j \neq j^{\prime}$. Hence, $Z \notin \DSpace{Z_{i_{j^{\prime}}}}$ if $j \neq j^{\prime}$. Thus,
\begin{equation*}
\DMUp{Z_{i_j}} \cap \DMUp{Z_{i_{j^{\prime}}}} = \varnothing \textit{ if } j \neq j^{\prime}.
\end{equation*}

Assume $Z$ is a monomial of some $T_t^j$ and $Z \in \DMUp{Z_{i_j}}$. Then it follows from the above argument that $Z$ is not contained in any $T_{t^{\prime}}^{j^{\prime}}$ whenever $j \neq j^{\prime}$. We proved that the monomials of $\DMUp{Z_{i_j}}$, $j = 1, \ldots, m$, have to cancel out in the sum~\eqref{dep_sum}. Therefore, the monomials of $\DMUp{Z_{i_j}}$ have to cancel out completely in the corresponding sum $\sum_{t = 1}^{l_j} \gamma_t^j T_t^j$. As a result we obtain
\begin{equation*}
\sum_{t = 1}^{l_j} \gamma_t^j T_t^j \in \Low\DSpace{Z_{i_j}} \textit{ for every } j = 1, \ldots, m.
\end{equation*}
Since $T_t^j \in \Dp\DSpace{Z_{i_j}}$, we have $\sum_{t = 1}^{l_j} \gamma_t^j T_t^j \in \Dp\DSpace{Z_{i_j}}$. Hence, we see that
\begin{equation*}
\sum_{t = 1}^{l_j} \gamma_t^j T_t^j \in \Dp\DSpace{Z_{i_j}} \cap \Low\DSpace{Z_{i_j}}.
\end{equation*}
So, by Lemma~\ref{fall_through_linear_dep_one}, we have
\begin{equation*}
\sum\limits_{t = 1}^{l_j} \gamma_t^j T_t^j \in \Dp(\Ft_{n - 1}(\Galg)) \textit{ for all } j = 1, \ldots, m.
\end{equation*}
Thus, we obtain
\begin{equation*}
\sum\limits_{i = 1}^l \gamma_i T_i = \sum\limits_{j = 1}^m\sum\limits_{t = 1}^{l_j} \gamma_t^j T_t^j \in \Dp(\Ft_{n - 1}(\Galg)).
\end{equation*}
Proposition~\ref{fall_through_linear_dep} is proved.
\end{proof}

Using Proposition~\ref{fall_through_linear_dep}, we obtain the following important statement.
\begin{proposition}
\label{fall_to_smaller_subspace}
Suppose $X, Y$ are subspaces of $\Galg$ generated by monomials and closed under taking derived monomials, $Y \subseteq X$. Then $\Dp(X)\cap Y = \Dp(Y)$.
\end{proposition}
\begin{proof}
Since $\Dp(Y) \subseteq \Dp(X)$ and $\Dp(Y) \subseteq Y$, we have $\Dp(Y) \subseteq \Dp(X)\cap Y$. Let us show that $\Dp(X)\cap Y \subseteq \Dp(Y)$.

Since $\Dp(X) \subseteq \Dp(\Galg)$, we evidently have
\begin{equation*}
\Dp(X) \cap Y \subseteq \Dp(\Galg) \cap Y.
\end{equation*}
We will show that $\Dp(\Galg) \cap Y \subseteq \Dp(Y)$. This implies $\Dp(X) \cap Y \subseteq \Dp(Y)$.

Let $T \in \Dp(\Galg)\cap Y$. Then
\begin{equation*}
T = \gamma_1T_1 +\ldots + \gamma_lT_l, \textit{ where } T_1, \ldots, T_l  \in \GDp^{\prime},\ \gamma_1, \ldots, \gamma_l \in \fld.
\end{equation*}
Let us show that $T \in \Dp(Y)$. Assume that some $T_{i_0} \in Y$, $1 \leqslant i_0 \leqslant l$. Since $Y$ is generated by monomials, we obtain that every monomial of $T_{i_0}$ belongs to $Y$. Therefore, $T_{i_0} \in \GDp^{\prime}(Y) \subseteq \Dp(Y)$. Hence, it is sufficient to prove that $T - \gamma_{i_0}T_{i_0} \in \Dp(Y)$. So, further we assume that every $T_i\notin Y$, $i = 1, \ldots, l$.

Denote by $X^{\prime}$ the linear space generated by all the monomials of $T_i$, $i = 1, \ldots, l$, and their derived monomials. Since $X^{\prime}$ is generated by derived monomials of a finite number of monomials, $f$-characteristics of monomials from $X^{\prime}$ is bounded. Therefore, there exists a unique $N$ such that
\begin{equation*}
X^{\prime} \subseteq \Ft_N(\Galg) \textit{ and } X^{\prime} \nsubseteq \Ft_{N - 1}(\Galg),
\end{equation*}
and the subspace $\Low(X^{\prime})$ is defined. Let us prove that $\gamma_1T_1 +\ldots + \gamma_lT_l \in \Dp(Y)$ by induction on $N$.

First we do the step of induction. Let $Z_i$ be a monomial such that $T_i$ is a layout of a multi-turn of a virtual member of the chart of $Z_i$. If $Z_i \in \Low(X^{\prime})$, then all the monomials of $T_i$ belong to $\Low(X^{\prime})$, because $\Low(X^{\prime})$ is closed under taking derived monomials. Assume $Z_i \in X^{\prime} \setminus \Low(X^{\prime})$. Let $Z$ be an arbitrary monomial of $T_i$ such that $Z \in X^{\prime} \setminus \Low(X^{\prime})$. Since $Z \in X^{\prime} \setminus \Low(X^{\prime})$, we have $Z \in \DMUp{Z_i}$. Therefore, Lemma~\ref{derived_spaces_equality} yields that $\DSpace{Z_i} = \DSpace{Z}$. That is, $Z_i$ is a derived monomial of $Z$. Therefore, all the monomials of $T_i$ are derived monomials of $Z$. Assume $Z \in Y$. Then all its derived monomials belong to $Y$ because $Y$ is closed under taking derived monomials. Hence, all the monomials of $T_i$ are contained in $Y$. This contradicts to our assumption that $T_i \notin Y$. So, $Z \notin Y$. Therefore, the monomials of $T_i$ that are contained in $X^{\prime} \setminus \Low(X^{\prime})$ are not contained in $Y$.

Since $\gamma_1T_1 +\ldots + \gamma_lT_l \in Y$ and $Y$ is generated by monomials, the monomials of $\gamma_1T_1 +\ldots + \gamma_lT_l$ that remain after the additive cancellations belong to~$Y$. Since monomials of $X^{\prime} \setminus \Low(X^{\prime})$ are not contained in $Y$, the monomials of $X^{\prime} \setminus \Low(X^{\prime})$ cancel in the sum $\gamma_1T_1 +\ldots + \gamma_lT_l$. Therefore, $\gamma_1T_1 +\ldots + \gamma_lT_l \in \Low(X^{\prime})$. Since $X^{\prime} \subseteq \Ft_N(\Galg)$ and $ X^{\prime} \nsubseteq \Ft_{N - 1}(\Galg)$, by definition, we see that $\Low(X^{\prime}) = X^{\prime}\cap \Ft_{N - 1}(\Galg) \subseteq \Ft_{N - 1}(\Galg)$. Hence,
\begin{equation*}
\gamma_1T_1 +\ldots + \gamma_lT_l \in \Ft_{N - 1}(\Galg).
\end{equation*}
Since $X^{\prime} \subseteq \Ft_N(\Galg)$, we obviously have
\begin{equation*}
\gamma_1T_1 +\ldots + \gamma_lT_l \in \Dp(\Ft_{N}(\Galg)).
\end{equation*}
Therefore, we see that
\begin{equation*}
\gamma_1T_1 +\ldots + \gamma_lT_l \in \Ft_{N - 1}(\Galg) \cap \Dp(\Ft_{N}(\Galg)).
\end{equation*}
So, Proposition~\ref{fall_through_linear_dep} implies
\begin{equation*}
\gamma_1T_1 +\ldots + \gamma_lT_l\in \Dp(\Ft_{N - 1}(\Galg)).
\end{equation*}

First assume that $\Dp(\Ft_{N - 1}(\Galg)) \neq 0$. Then we obtain
\begin{equation*}
\gamma_1T_1 +\ldots + \gamma_lT_l = \delta_1T_1^{\prime} +\ldots + \delta_{l^{\prime}}T_{l^{\prime}}^{\prime},
\end{equation*}
where $T_i^{\prime} \in \GDp^{\prime}(\Ft_{N - 1}(\Galg))$, $i = 1, \ldots, l^{\prime}$. Therefore, we have
\begin{equation*}
\gamma_1T_1 +\ldots + \gamma_lT_l = \delta_1T_1^{\prime} +\ldots + \delta_{l^{\prime}}T_{l^{\prime}}^{\prime} = \sum\limits_{T^{\prime}_i \in Y} \delta_i T^{\prime}_i + \sum\limits_{T^{\prime}_i \notin Y} \delta_i T^{\prime}_i.
\end{equation*}
As above, every element of the first sum $\sum_{T^{\prime}_i \in Y} \delta_i T^{\prime}_i $ belongs to $\Dp(Y)$. Since $T_i^{\prime} \in \GDp^{\prime}(\Ft_{N - 1}(\Galg))$, we obtain that all the monomials of $T_i^{\prime}$, $i = 1, \ldots, l^{\prime}$, and their derived monomials belong to $\Ft_{N - 1}(\Galg)$. Hence, the second sum $\sum_{T^{\prime}_i \notin Y} T^{\prime}_i$ belongs to $\Dp(Y)$ by the induction hypothesis. Thus, $\gamma_1T_1 +\ldots + \gamma_lT_l \in \Dp(Y)$.

Now assume that $\Dp(\Ft_{N - 1}(\Galg)) = 0$. Then, evidently, $\gamma_1T_1 +\ldots + \gamma_lT_l = 0 \in \Dp(Y)$. Clearly, in this case for every $n \leqslant N - 1$, we have $\Dp(\Ft_n(\Galg)) = 0$ and $\GDp^{\prime}(\Ft_n(\Galg)) = \varnothing$. So, this case at the same time is our basis of induction.

Thus, we obtain $\Dp(\Galg) \cap Y \subseteq \Dp(Y)$. This concludes the proof.
\end{proof}

\subsection{Tensor products. The proof of the Main Lemma}
\label{tensor_product_section}
The main goal of this section is to prove Lemma~\ref{fall_through_linear_dep_one}.

Assume $U$ is a monomial. Derived monomials of $U$ are defined with the use of certain sequences of replacements of virtual members of the chart (see Definition~\ref{derived_monomials}). When we perform replacements that preserve $f$-characteristics of monomials, they preserve, roughly speaking, the structure of the chart. Moreover, there is no interaction between the replaced occurrence and the separated virtual members of the chart and there is a very small interaction between the replaced occurrence and its neighbours. This kind of behaviour provides the idea to consider a tensor product of linear spaces that correspond to each place of the chart of $U$.

Assume a monomial $U$ has $m$ virtual members of the chart, that is, $\nvirt{U} = m$. We enumerate all the virtual members of the chart of $U$ from left to right. Let $u^{(i)}$ be the $i$-th virtual member of the chart of $U$. We define a linear space $A_i[U]$ by the following formula
\begin{equation}
\label{tens_prod_components}
A_i[U] = \left\langle\left\lbrace a^{(i)} \mid u^{(i)} \textit{ and } a^{(i)} \textit{ are } U\textit{-incident monomials}\right\rbrace\right\rangle.
\end{equation}

Suppose $U = L^{(i)}u^{(i)}R^{(i)}$. We define two sets of monomials $ME_i[U] \subseteq A_i[U]$ and $ML_i[U] \subseteq A_i[U]$ by the following rule:
\begin{align}
\begin{split}
\label{tens_prod_components_up_monomials}
ME_i[U] = &\left\lbrace a^{(i)} \mid u^{(i)} \textit{ and } a^{(i)} \textit{ are } U\textit{-incident monomials},\right.\\
&\left.L^{(i)}a^{(i)}R^{(i)}\in \DMUp{U}\right\rbrace;
\end{split}
\\
\begin{split}
\label{tens_prod_components_low_monomials}
ML_i[U] = &\left\lbrace a^{(i)} \mid u^{(i)} \textit{ and } a^{(i)}\textit{ are } U\textit{-incident monomials},\right.\\
&\left.L^{(i)}a^{(i)}R^{(i)} \in \DMLow{U}\right\rbrace.
\end{split}
\end{align}
We define the subspace $L_i[U] \subseteq A_i[U]$ by the formula
\begin{equation}
\label{tens_prod_components_low}
L_i[U] = \left\langle ML_i[U]\right\rangle.
\end{equation}
Let us put
\begin{equation*}
\Low(A_1[U]\otimes\ldots\otimes A_m[U]) = \sum\limits_{i = 1}^{m} A_1[U]\otimes \ldots \otimes L_i[U] \otimes \ldots \otimes A_m[U].
\end{equation*}

We put
\begin{equation}
\label{tens_prod_dep}
D_i[U] = \left\langle\Rel \cap A_i[U]\right\rangle.
\end{equation}

\begin{remark}
Let us show that $A_i[U]$ is closed under taking derived monomials.

First of all notice that if $a \in \Mon$, then $a$ is a virtual member of the chart of itself if and only if $\SPM(a) \geqslant \tau$. Indeed, since $a$ is a single maximal occurrence inside itself, there are no $(a, a)$-admissible sequences.

Let $a^{(i)}$ be a monomial from $A_i[U]$. It follows from the above that if $\SPM(a^{(i)}) < \tau$, then it does not have any derived monomials. Assume $\SPM(a^{(i)}) \geqslant \tau$, and assume $b^{(i)}$ is a derived monomial of $a^{(i)}$. By definition, this means that there exists a sequence of transformations of type~\ref{repl1} and~\ref{repl2} such that $b^{(i)}$ is the resulting monomial. In order to obtain a derived monomial, we replace only virtual members of the chart by incident monomials. Since $a^{(i)}\in \Mon$, after every such replacement, except the last one, the resulting monomial is an element of $\Mon$ such that it is a virtual member of the chart of itself. However, we proved above that this holds if and only if its $\SPM$-measure is $\geqslant \tau$. Therefore, there exists a sequence of monomials $m^{(i)}_1, \ldots, m_k^{(i)}$ such that $a^{(i)} = m^{(i)}_1$, $b^{(i)} = m^{(i)}_k$, $m^{(i)}_j$ and $m^{(i)}_{j + 1}$ are incident monomials for $j = 1, \ldots, k - 1$, and $\SPM(m^{(i)}_j) \geqslant \tau$ for $j = 1, \ldots, k - 1$.

Let us return to the monomial $U = L^{(i)}u^{(i)}R^{(i)}$. Since $m^{(i)}_j$ and $m^{(i)}_{j + 1}$ are incident monomials and they are not small pieces for $j = 1, \ldots , k - 2 $, it follows from the results of Section~\ref{mt_configurations} (see~\ref{a_j_keep_structure}) that $m^{(i)}_j$ is a maximal occurrence in $L^{(i)}m^{(i)}_jR^{(i)}$. Since $\SPM(m^{(i)}_j) \geqslant \tau$ for $j = 1, \ldots, k - 1$, we see that $m^{(i)}_j$ is a virtual member of the chart of $L^{(i)}m^{(i)}_jR^{(i)}$ for $j = 1, \ldots, k - 1$. Hence, by definition, $a^{(i)}$ and $b^{(i)}$ are $U$-incident monomials. Thus, $b^{(i)} \in A_i[U]$, and, therefore, $A_i[U]$ is closed under taking derived monomials.
\end{remark}

\begin{remark}
Since the space $A_i[U]$ is closed under taking derived monomials, the space $\Dp(A_i[U])$ is defined. Let $\sum_{j = 1}^l \alpha_ja^{(i)}_j \in \Rel \cap A_i[U]$. It follows from Small Cancellation Axiom that in the polynomial $\sum_{j = 1}^l \alpha_ja^{(i)}_j$ there exists a monomial of $\SPM$-measure $\geqslant \tau + 1$. By definition, this monomial is a virtual member of the chart of itself. Obviously, $\sum_{j = 1}^l \alpha_ja^{(i)}_j$ is a layout of a multi-turn of this monomial. Therefore, $\sum_{j = 1}^l \alpha_ja^{(i)}_j \in \Dp(A_i[U])$. Hence, $D_i[U] \subseteq \Dp(A_i[U])$.

Let $\sum_{j = 1}^l \alpha_ja^{(i)}_j$ be a layout of a multi-turn of a virtual member of the chart of a monomial from $A_i[U]$. Since $A_i[U]$ is closed under taking derived monomials, every $a^{(i)}_j \in A_i[U]$. Since every $a^{(i)}_j \in \Mon$, by the definition of layouts, we obtain  $\sum_{j = 1}^l \alpha_ja^{(i)}_j \in \Rel$. Combining these facts, we obtain $\sum_{j = 1}^l \alpha_ja^{(i)}_j \in D_i[U]$. Since $\Dp(A_i[U])$ is linearly generated by the layouts of multi-turns of virtual members of the chart of the monomials from $A_i[U]$, we have $\Dp(A_i[U]) \subseteq D_i[U]$.

Thus, we finally see that $D_i[U] = \Dp(A_i[U])$.
\end{remark}

\begin{lemma}
\label{glue_U_incident}
Let $U$ be a monomial, $u$ be a virtual member of the chart of $U$, $U = LuR$. Let $u$ and $b$ be $U$-incident monomials, and $u$ and $c$ be $U$-incident monomials. Assume $b \in \mo{LbR}$. Then $b$ and $c$ are $LbR$-incident monomials.
\end{lemma}
\begin{proof}
Since $u$ and $b$ are $U$-incident monomials, there exists a sequence of monomials $m_1^{(1)}, \ldots, m_{s_1}^{(1)}$, $m_1^{(1)} = u$, $m_{s_1}^{(1)} = b$ that satisfies the conditions of Definition~\ref{U_incident_monomials}. Since $u$ and $c$ are $U$-incident monomials, there exists a sequence of monomials $m_1^{(2)}, \ldots, m_{s_2}^{(2)}$, $m_1^{(2)} = u$, $m_{s_2}^{(2)} = c$ that satisfies the conditions of Definition~\ref{U_incident_monomials}. Since  $m_1^{(1)} = m_2^{(2)} = u$ is a virtual member of the chart of $U = LuR$, we obtain that the sequence
\begin{equation*}
m_{s_1}^{(1)} = b, m_{s_1 - 1}, \ldots, m_1^{(1)} = m_1^{(2)}, m_2^{(2)}, \ldots, m_{s_2}^{(2)} = c
\end{equation*}
satisfies the conditions of Definition~\ref{U_incident_monomials}. Thus, since $b \in \mo{LbR}$, we obtain that $b$ and $c$ are $LbR$-incident monomials.
\end{proof}

\begin{definition}
\label{mu_def}
Let $U$ be a monomial, $\nvirt{U} = m$, $U \in \Ft_{n}(\Galg) \setminus \Ft_{n - 1}(\Galg)$. Assume $u^{(i)}$ is the $i$-th virtual member of the chart of $U$, $i = 1, \ldots, m$. We construct a linear mapping
\begin{equation*}
\mu[U] : A_1[U] \otimes \ldots \otimes A_m[U] \to \DSpace{U} +  \Ft_{n - 1}(\Galg).
\end{equation*}
By definition, $A_i[U]$ is generated by monomials. Obviously, all monomials of $A_i[U]$ is a basis of $A_i[U]$. Hence, all the elements $a^{(1)} \otimes\ldots\otimes a^{(m)}$, where $a^{(i)}$ is a monomial of $A_i[U]$, is a basis of $A_1[U]\otimes \ldots \otimes A_m[U]$. We define $\mu[U]$ on these basis elements of $A_1[U]\otimes \ldots \otimes A_m[U]$ and then extend $\mu[U]$ by linearity on the whole space $A_1[U]\otimes \ldots \otimes A_m[U]$.

Informally speaking, the general idea is as follows. The element $a^{(1)} \otimes\ldots\otimes a^{(m)}$ encodes the replacements $u^{(1)} \leadsto a^{(1)}, \ldots, u^{(m)} \leadsto a^{(m)}$ in $U$. We want to perform them consecutively and to call the resulting monomial $\mu[U]\left(a^{(1)} \otimes\ldots\otimes a^{(m)}\right)$. However, in Subsection~\ref{replacements_by_U_incident_section} we defined consecutive performing  the replacements only for the case when $u^{(i)} \leadsto a^{(i)}$ preserves $f$-characteristic. That is, we considered only the case when every $a^{(i)}$ is a virtual member of the chart of the corresponding resulting monomial (see Definition~\ref{consecutive_replacements_def}). Also in Subsection~\ref{replacements_by_U_incident_section} we described possible difficulties that take place if some $u^{(i)} \leadsto a^{(i)}$ does not preserve $f$-characteristic (see page~\pageref{go_lower_difficulties} and Remark~\ref{first_monomial_reduced_difficulties}) and how we work with them (see procedure at page~\pageref{replacement_in_intersection_process}). So, in the definition of $\mu[U]$ we consider separately a number of cases.

Now let us precisely define $\mu[U]$ on the basis elements of $A_1[U]\otimes \ldots \otimes A_m[U]$. We distinguish between the following four possibilities.
\begin{enumerate}[label=\textbf{Case~\arabic*}, ref=Case~\arabic*]
\item
\label{mu_case_all_up}
Consider an element $a^{(1)} \otimes\ldots\otimes a^{(m)}$ such that every $a^{(i)} \in ME_i[U]$. Recall that $u^{(i)}$ and $a^{(i)}$ are $U$-incident monomials. Consider the replacements $u^{(1)} \leadsto a^{(1)}, \ldots, u^{(m)} \leadsto a^{(m)}$ in $U$. By definition, we put $\mu[U]\left(a^{(1)} \otimes\ldots\otimes a^{(m)}\right)$ to be the result of consecutive performing the replacements $u^{(1)} \leadsto a^{(1)}, \ldots, u^{(m)} \leadsto a^{(m)}$, starting from $U$ (see Definition~\ref{consecutive_replacements_def}). Recall that, by Lemma~\ref{replace_consecutively}, we can perform them in any order and obtain the same result.

It follows from Lemma~\ref{replace_consecutively} that the monomial $\mu[U]\left(a^{(1)} \otimes\ldots\otimes a^{(m)}\right)$ is a derived monomial of $U$ with the same $f$-characteristic as $U$. That is, $\mu[U]\left(a^{(1)} \otimes\ldots\otimes a^{(m)}\right) \in \DMUp{U}$.
\item
\label{mu_case_one_low}
Assume $1 \leqslant i_0 \leqslant m$ is a position in the chart of $U$. Consider an element $a^{(1)} \otimes\ldots\otimes a^{(m)}$ such that $a^{(i)} \in ME_i[U]$ if $i \neq i_0$, and $a^{(i_0)} \in ML_{i_0}[U]$. Then first we consecutively perform the replacements $u^{(i)} \leadsto a^{(i)}$, $i = 1, \ldots, m$, $i \neq i_0$, starting from $U$. Recall that, by Lemma~\ref{replace_consecutively}, we can perform them in any order and obtain the same result. Let $Z$ be the resulting monomial.

Assume $\widetilde{u}^{(i_0)}$ is the $i_0$-th virtual member of the chart of $Z$. Let $\widetilde{u}^{(i_0)} \leadsto \dbtilde{a}^{(i_0)}$ be the replacement in $Z$ that corresponds to the replacement $u^{(i_0)} \leadsto a^{(i_0)}$ in $U$ (see Remark~\ref{consecurive_repl_basic_properties}). Since $a^{(i)} \in ME_i[U]$ for $i \neq i_0$, it follows from statements~\ref{virtual_members_U_incident_replacements2} and~\ref{virtual_members_U_incident_replacements5} of Corollary~\ref{virtual_members_U_incident_replacements} that $\widetilde{u}^{(i_0)}$ and $\dbtilde{a}^{(i_0)}$ are $Z$-incident monomials. Finally we perform the replacement $\widetilde{u}^{(i_0)} \leadsto \dbtilde{a}^{(i_0)}$ in $Z$ and do the cancellations if there are any (cancellations may occur only if $\dbtilde{a}^{(i_0)} = 1$). The resulting monomial of the replacement $\widetilde{u}^{(i_0)} \leadsto \dbtilde{a}^{(i_0)}$ in $Z$ is, by definition, $\mu[U]\left(a^{(1)} \otimes\ldots\otimes a^{(m)}\right)$.

By definition of consecutive performing  the replacements (see Definition~\ref{consecutive_replacements_def}), the monomial $Z$ is a derived monomial of $U$. Since $\widetilde{u}^{(i_0)}$ and $\dbtilde{a}^{(i_0)}$ are $Z$-incident monomials, we obtain, by definition, that $\mu[U]\left(a^{(1)} \otimes\ldots\otimes a^{(m)}\right)$ is a derived monomial of $U$. Furthermore, since $a^{(i_0)} \in ML_{i_0}[U]$, statements~\ref{virtual_members_U_incident_replacements3} and~\ref{virtual_members_U_incident_replacements6} of Corollary~\ref{virtual_members_U_incident_replacements} imply that $\dbtilde{a}^{(i_0)}$ is not a virtual member of the chart of $Z[\widetilde{u}^{(i_0)} \leadsto \dbtilde{a}^{(i_0)}]$. That is, the replacement $\widetilde{u}^{(i_0)} \leadsto \dbtilde{a}^{(i_0)}$ in $Z$ decreases $f$-characteristic of the resulting monomial. Therefore, $\mu[U]\left(a^{(1)} \otimes\ldots\otimes a^{(m)}\right)$ belongs to $\DMLow{U} \subseteq \Ft_{n - 1}(\Galg)$.
\item
\label{mu_case_two_low}
Assume $1 \leqslant i_1 < i_2 \leqslant m$ are two different positions in the chart of~$U$. Consider an element $a^{(1)} \otimes\ldots\otimes a^{(m)}$ such that $a^{(i)} \in ME_i[U]$ if $i \neq i_1, i_2$, and $a^{(i_1)} \in ML_{i_1}[U]$, $a^{(i_2)} \in ML_{i_2}[U]$.

First we consecutively perform the replacements $u^{(i)} \leadsto a^{(i)}$, $i = 1, \ldots, m$, $i \neq i_1, i_2$, starting from $U$. Let $Z$ be the resulting monomial.

Assume $\widetilde{u}^{(i_1)}$ is the $i_1$-th virtual member of the chart of $Z$, and $\widetilde{u}^{(i_2)}$ is the $i_2$-th virtual member of the chart of $Z$. Let $\widetilde{u}^{(i_1)} \leadsto \dbtilde{a}^{(i_1)}$ be the replacement in $Z$ that corresponds to the replacement $u^{(i_1)} \leadsto a^{(i_1)}$ in $U$, and $\widetilde{u}^{(i_2)} \leadsto \dbtilde{a}^{(i_2)}$ be the replacement in $Z$ that corresponds to the replacement $u^{(i_2)} \leadsto a^{(i_2)}$ in $U$ (see Remark~\ref{consecurive_repl_basic_properties}). We will do the procedure that we introduced in the second part of Subsection~\ref{replacements_by_U_incident_section} (see page~\pageref{replacement_in_intersection_process}) with these two replacements.

Namely, we do the following process, starting from $Z$. Since $a^{(i)} \in ME_i[U]$ for $i \neq i_1, i_2$, it follows from statements~\ref{virtual_members_U_incident_replacements2} and~\ref{virtual_members_U_incident_replacements5} of Corollary~\ref{virtual_members_U_incident_replacements} that $\widetilde{u}^{(i_1)}$ and $\dbtilde{a}^{(i_1)}$ are $Z$-incident monomials. We do the replacement $\widetilde{u}^{(i_1)} \leadsto \dbtilde{a}^{(i_1)}$ in $Z$. If $\widetilde{a}^{(i_1)} = 1$, there may be cancellations in the resulting monomial. In this case we do not do the cancellations after the replacement.

Let $Z = \widetilde{L}\widetilde{u}^{(i_1)}\widetilde{R}$, where $\widetilde{u}^{(i_1)}$ is the $i_1$-th virtual member of the chart of $Z$. Then the resulting monomial of the replacement $\widetilde{u}^{(i_1)} \leadsto \dbtilde{a}^{(i_1)}$ in $Z$ is equal to $\widetilde{L}\widetilde{a}^{(i_1)}\widetilde{R}$. Let $\widehat{u}^{(i_2)}$ be the intersection of $\widetilde{R}$ and $\widetilde{u}^{(i_2)}$. Let $\widehat{u}^{(i_2)} \leadsto \dbhat{a}^{(i_2)}$ be the replacement in $\widetilde{R}$ that corresponds to the replacement $\widetilde{u}^{(i_2)} \leadsto \dbtilde{a}^{(i_2)}$ in $Z$ (see~\eqref{corresponding_replacement_in_intersection_right}). We do the replacement $\widehat{u}^{(i_2)} \leadsto \dbhat{a}^{(i_2)}$ in $\widetilde{R}$, let $\widetilde{R}[\widehat{u}^{(i_2)} \leadsto \dbhat{a}^{(i_2)}]$ be the resulting monomial. After that we do the cancellations in the monomial $\left(\widetilde{L}\dbtilde{a}^{(i_1)}\right)\cdot \widetilde{R}[\widehat{u}^{(i_2)} \leadsto \dbhat{a}^{(i_2)}]$ if there are any. The resulting monomial is, by definition, $\mu[U]\left(a^{(1)} \otimes\ldots\otimes a^{(m)}\right)$.

By Lemma~\ref{replace_consecutively}, we obtain that $Z$ is a derived monomial of $U$ with the same $f$-characteristic as $U$. That is, $Z \in \Ft_{n}(\Galg)$. It follows from statement~\ref{replacenets_in_intersection_U_incident2} of Lemma~\ref{replacenets_in_intersection_U_incident} that the monomial $\mu[U]\left(a^{(1)} \otimes\ldots\otimes a^{(m)}\right)$ obtained by the above process in $Z$ belongs to $\Ft_{n - 1}(\Galg)$.

Clearly, we can do the above process in $Z$, starting from the position $i_2$ (then in the second step we deal with the intersection of $\widetilde{u}^{(i_1)}$ with~$\widetilde{L}$). It follows from statement~\ref{replacenets_in_intersection_U_incident1} of Lemma~\ref{replacenets_in_intersection_U_incident} that we obtain the same resulting monomial as if we start from the position $i_1$.
\item
\label{mu_case_many_low}
Consider an element $a^{(1)} \otimes\ldots\otimes a^{(m)}$ such that there are more than two $a^{(i)} \in ML_i[U]$. We do not need to preserve full information about these elements. So, by definition, we put $\mu[U]\left(a^{(1)} \otimes\ldots\otimes a^{(m)}\right) = 0$ in this case.
\end{enumerate}
\end{definition}

The following lemma states properties of $\mu[U]$. We will use them in order to prove Lemma~\ref{fall_through_linear_dep_one}. In order to prove statement~\ref{mu_properties_up_correspondence} we use Isolation Axiom. Here we prove statement~\ref{mu_properties_up_correspondence} under assumption that right-sided Isolation Axiom holds. The case when left-sided Isolation Axiom holds is considered similarly.

\begin{lemma}
\label{mu_properties}
Let $U$ be a monomial. Assume $U \in \Ft_{n}(\Galg) \setminus \Ft_{n - 1}(\Galg)$ and $U$ has $m$ virtual members of the chart. For the matter of convenience, we denote all the layouts of all the multi-turns of $i$-th virtual members of the chart of the monomials of $\DMUp{U}$ by $\GDp^{(i)}[U]$. Let $\mu[U]$ be a mapping defined above by Definition~\ref{mu_def}. Then $\mu[U]$ possesses the following properties.
\begin{enumerate}[label={(\arabic*)}]
\item
\label{mu_properties_low_correspondence}
$\mu[U](\Low(A_1[U]\otimes \ldots \otimes A_m[U])) \subseteq\Ft_{n - 1}(\Galg)$.
\item
\label{mu_properties_up_correspondence}
Let $a^{(1)}\otimes \ldots \otimes a^{(m)} \in A_1[U]\otimes \ldots \otimes A_m[U]$ such that every $a^{(i)} \in ME_i[U]$. Then $\mu[U]\left(a^{(1)}\otimes \ldots \otimes a^{(m)}\right) \in \DMUp{U}$. Moreover, $\mu[U]$ gives a bijective correspondence between all the elements $a^{(1)}\otimes \ldots \otimes a^{(m)}$ such that every $a^{(i)} \in ME_i[U]$ and all the monomials of $\DMUp{U}$.
\item
\label{mu_properties_up_dep_image}
Let $a^{(1)}\otimes \ldots \otimes t^{(i_0)}\otimes \ldots \otimes a^{(m)} \in A_1[U]\otimes \ldots \otimes D_{i_0}[U]\otimes \ldots \otimes A_m[U]$ such that $a^{(i)} \in ME_i[U]$ for $i \neq i_0$ and $t^{(i_0)} \in \Rel \cap D_{i_0}[U]$. Then $\mu[U]\left(a^{(1)}\otimes \ldots \otimes t^{(i_0)}\otimes \ldots \otimes a^{(m)}\right) \in \GDp^{(i_0)}[U]$.
\item
\label{mu_properties_dep_inverse_image}
Let $T \in \GDp^{\prime}(\DSpace{U})$ such that $T \notin \Low\DSpace{U}$. Then there exists an element $a^{(1)}\otimes \ldots \otimes t^{(i_0)}\otimes \ldots \otimes a^{(m)} \in A_1[U]\otimes \ldots \otimes D_{i_0}[U]\otimes \ldots \otimes A_m[U]$ such that $a^{(i)} \in ME_i[U]$  for $i \neq i_0$, $t^{(i_0)} \in \Rel \cap D_{i_0}[U]$, and $\mu[U]\left(a^{(1)}\otimes \ldots \otimes t^{(i_0)}\otimes \ldots \otimes a^{(m)}\right) = T$.
\item
\label{mu_properties_dep_low_image}
Let $a^{(1)}\otimes \ldots \otimes t^{(i_1)}\otimes \ldots \otimes a^{(m)} \in A_1[U]\otimes \ldots \otimes D_{i_1}[U]\otimes\ldots \otimes A_m[U]$ such that $a^{(i)} \in ME_i[U]$  for $i \neq i_1, i_2$, and $a^{(i_2)} \in ML_{i_2}[U]$. Then $\mu[U]\left(a^{(1)}\otimes \ldots \otimes t^{(i_1)}\otimes \ldots \otimes a^{(m)}\right) \in \Dp(\Ft_{n - 1}(\Galg))$. Moreover, if $t^{(i_1)} \in \Rel \cap D_{i_1}[U]$, then $\mu[U]\left(a^{(1)}\otimes \ldots \otimes t^{(i_1)}\otimes \ldots \otimes a^{(m)}\right) \in \GDp^{\prime}(\Ft_{n - 1}(\Galg))$.
\end{enumerate}
\end{lemma}
\begin{proof}
Let us fix the notations. Namely, let $u^{(1)}, \ldots, u^{(m)}$ be all the virtual members of the chart of $U$ enumerated from left to right. We use the same notations for resulting monomials of replacements as we used in Subsection~\ref{replacements_by_U_incident_section}. Namely, let $Z$ be a monomial, $a \in \mo{Z}$, and $a, b$ be $Z$-incident monomials. Then the resulting monomial of the replacement $a\leadsto b$ in $Z$ is denoted by~$Z[a\leadsto b]$.

~\paragraph*{(1)} Let $a^{(1)}\otimes \ldots \otimes a^{(m)} \in \Low(A_1[U]\otimes \ldots \otimes A_m[U])$, where all $a^{(i)}\in A_i[U]$ are monomials, $i = 1, \ldots, m$. This means that there exist positions $i_1, \ldots, i_k$ such that $a^{i_j} \in ML_{i_j}[U]$, $j = 1, \ldots, k$. So, $\mu[U]\left(a^{(1)}\otimes \ldots \otimes a^{(m)}\right)$ is defined by~\ref{mu_case_one_low}---\ref{mu_case_many_low} of Definition~\ref{mu_def}. We already noticed in the definition of $\mu[U]$ that in~\ref{mu_case_one_low} and~\ref{mu_case_two_low} we have $\mu[U]\left(a^{(1)}\otimes \ldots \otimes a^{(m)}\right) \in \Ft_{n - 1}(\Galg)$. In~\ref{mu_case_many_low}, by definition, $\mu[U]\left(a^{(1)}\otimes \ldots \otimes a^{(m)}\right)$ is equal to $0$, so, it is also contained in $\Ft_{n - 1}(\Galg)$. Thus, $\mu[U](\Low(A_1[U]\otimes \ldots \otimes A_m[U])) \subseteq \Ft_{n - 1}(\Galg)$.

~\paragraph*{(2)}
Assume $a^{(1)}\otimes \ldots \otimes a^{(m)} \in A_1[U]\otimes \ldots \otimes A_m[U]$, where all $a^{(i)}\in ME_i[U]$. We noticed in~\ref{mu_case_all_up} of the definition of $\mu[U]$ that then $\mu[U]\left(a^{(1)}\otimes \ldots \otimes a^{(m)}\right) \in \DMUp{U}$. Hence, the restriction of $\mu[U]$ on the set of all elements $a^{(1)}\otimes \ldots \otimes a^{(m)}$ such that every $a^{(i)}\in ME_i[U]$ is a mapping to $\DMUp{U}$. Let us show that this mapping is bijective.

First we show that the restriction of $\mu[U]$ is a surjective mapping. Let $Z \in \DMUp{U}$. Then it follows from Lemma~\ref{create_initial_replacenets} that there exists a set of replacements $u^{(i)} \leadsto b^{(i)}$, $i = 1, \ldots, m$, such that $u^{(i)}$ and $b^{(i)}$ are $U$-incident monomials, $b^{(i)}$ is a virtual member of the chart of $U[u^{(i)} \leadsto b^{(i)}]$, and if we perform these replacements consecutively, then $Z$ is the resulting monomial. Since every $b^{(i)}$ is a virtual member of the chart of $U[u^{(i)} \leadsto b^{(i)}]$, we see that $b^{(i)} \in ME_i[U]$. Therefore, $\mu[U]\left(b^{(1)}\otimes\ldots\otimes b^{(m)}\right)$ is defined using~\ref{mu_case_all_up} of Defintion~\ref{mu_def}. Hence, we see that $Z = \mu[U]\left(b^{(1)}\otimes\ldots\otimes b^{(m)}\right)$. Thus, $\mu[U]$ is a surjective mapping from the set of all elements $a^{(1)}\otimes \ldots \otimes a^{(m)}$ such that every $a^{(i)}\in ME_i[U]$ to $\DMUp{U}$.

Assume $b^{(1)}\otimes \ldots \otimes b^{(m)}, c^{(1)}\otimes \ldots \otimes c^{(m)}$ are two elements of $A_1[U]\otimes \ldots \otimes A_m[U]$ such that $b^{(i)}, c^{(i)} \in ME_i[U]$, and $b^{(1)}\otimes \ldots \otimes b^{(m)} \neq c^{(1)}\otimes \ldots \otimes c^{(m)}$. Let us show that
\begin{equation*}
\mu[U]\left(b^{(1)}\otimes \ldots \otimes b^{(m)}\right) \neq \mu[U]\left(c^{(1)}\otimes \ldots \otimes c^{(m)}\right).
\end{equation*}

Since $b^{(1)}\otimes \ldots\otimes b^{(m)} \neq c^{(1)}\otimes \ldots \otimes c^{(m)}$, not all $b^{(i)} = c^{(i)}$. Assume $m = 1$. Then $b^{(1)} \neq c^{(1)}$. We have
\begin{align*}
&\mu[U]\left(b^{(1)}\right) = U[u^{(1)} \leadsto b^{(1)}],\\
&\mu[U]\left(c^{(1)}\right) = U[u^{(1)} \leadsto c^{(1)}].
\end{align*}
Since $b^{(1)} \neq c^{(1)}$, we obviously see that $U[u^{(1)} \leadsto b^{(1)}] \neq U[u^{(1)} \leadsto c^{(1)}]$. Hence, $\mu[U]\left(b^{(1)}\right) \neq \mu[U]\left(c^{(1)}\right)$. So far, we are done with the case $m = 1$.

Assume $m > 1$. Let $i_0$ be the first position such that $b^{(i)}$ and $c^{(i)}$ differ from each other. That is, $b^{(i_0)} \neq c^{(i_0)}$ and $b^{(i)} = c^{(i)}$ if $i = 1, \ldots, i_{0} - 1$. Since $b^{(i)}, c^{(i)} \in ME_i[U]$ for all $i = 1, \ldots, m$, the monomials $\mu[U]\left(b^{(1)}\otimes \ldots \otimes b^{(m)}\right)$ and $\mu[U]\left(c^{(1)}\otimes \ldots \otimes c^{(m)}\right)$ are defined, using \ref{mu_case_all_up} of Definition~\ref{mu_def}. Recall that, according to Lemma~\ref{replace_consecutively}, in~\ref{mu_case_all_up} of the definition of $\mu[U]$ the replacements can be done in any order. First let us consecutively perform the replacements $u^{(i)} \leadsto b^{(i)}$ and $u^{(i)} \leadsto c^{(i)}$, starting from $U$, for $i = 1, \ldots, i_0 - 1$. Since $b^{(i)} = c^{(i)}$ for $i = 1, \ldots, i_0 - 1$, we obviously see that the resulting monomial is the same. Denote this resulting monomial by $Z$.

Let $\widetilde{u}^{(i_0)}$ be $i_0$-th virtual member of the chart of $Z$. Let $\widetilde{u}^{(i_0)} \leadsto \dbtilde{b}^{(i_0)}$ be the replacement in $Z$ that corresponds to the replacement $u^{(i_0)} \leadsto b^{(i_0)}$ in $U$. Let $\widetilde{u}^{(i_0)} \leadsto \dbtilde{c}^{(i_0)}$ be the replacement in $Z$ that corresponds to the replacement $u^{(i_0)} \leadsto c^{(i_0)}$ in $U$. Since all replacements $u^{(i)} \leadsto b^{(i)} = c^{(i)}$ for $i = 1, \ldots, i_0 - 1$ are from the left of the beginning of $u^{(i_0)}$, we obtain $\widetilde{u}^{(i_0)} = e^{\prime}\left(e^{-1} \cdot u^{(i_0)}\right)$, where $e$ and $e^{\prime}$ are small pieces. Then, according to formula~\eqref{corresponding_element_U_incident_right}, we have
\begin{equation*}
\dbtilde{b}^{(i_0)} = e^{\prime}\cdot e^{-1}\cdot b^{(i_0)}, \ \dbtilde{c}^{(i_0)} = e^{\prime}\cdot e^{-1}\cdot c^{(i_0)}.
\end{equation*}
Therefore, since $b^{(i_0)} \neq c^{(i_0)}$, we obtain
 $\dbtilde{b}^{(i_0)} \neq \dbtilde{c}^{(i_0)}$.

Assume $i_0 = m$. In this case
\begin{align*}
\mu[U]\left(b^{(1)}\otimes \ldots \otimes b^{(m)}\right) = Z[\widetilde{u}^{(i_0)} \leadsto \dbtilde{b}^{(i_0)}],\\
\mu[U]\left(c^{(1)}\otimes \ldots \otimes c^{(m)}\right) = Z[\widetilde{u}^{(i_0)} \leadsto \dbtilde{c}^{(i_0)}].
\end{align*}
Therefore, since $\dbtilde{b}^{(i_0)} \neq \dbtilde{c}^{(i_0)}$, we obtain $Z[\widetilde{u}^{(i_0)} \leadsto \dbtilde{b}^{(i_0)}] \neq Z[\widetilde{u}^{(i_0)} \leadsto \dbtilde{c}^{(i_0)}]$. Thus, $\mu[U]\left(b^{(1)}\otimes \ldots \otimes b^{(m)}\right) \neq \mu[U]\left(c^{(1)}\otimes \ldots \otimes c^{(m)}\right)$. So far, we are done with the case $i_0 = m$.

Assume $i_0 < m$. Clearly, $u^{(i_0 + 1)}$ can be considered as an occurrence in $Z$ and $u^{(i_0 + 1)}$ is the $(i_0 + 1)$-th virtual member of the chart of $Z$. Assume $\widetilde{u}^{(i_0)}$ and $u^{(i_0 + 1)}$ are separated in $Z$. Then there are no further changes of $\dbtilde{b}^{(i_0)}$ after the rest of the replacements in $Z[\widetilde{u}^{(i_0)} \leadsto \dbtilde{b}^{(i_0)}]$. Similarly, there are no further changes of $\dbtilde{c}^{(i_0)}$ after the rest of the replacements in $Z[\widetilde{u}^{(i_0)} \leadsto \dbtilde{c}^{(i_0)}]$. Hence, $\dbtilde{b}^{(i_0)}$ is the $i_0$-th virtual member of the chart of $\mu[U]\left(b^{(1)}\otimes \ldots \otimes b^{(m)}\right)$ and $\dbtilde{c}^{(i_0)}$ is the $i_0$-th virtual member of the chart of $\mu[U]\left(c^{(1)}\otimes \ldots \otimes c^{(m)}\right)$. Since $\dbtilde{b}^{(i_0)} \neq \dbtilde{c}^{(i_0)}$, we see that $\mu[U]\left(b^{(1)}\otimes \ldots \otimes b^{(m)}\right) \neq \mu[U]\left(c^{(1)}\otimes \ldots \otimes c^{(m)}\right)$ in this case.

Assume $\widetilde{u}^{(i_0)}$ and $u^{(i_0 + 1)}$ are not separated in $Z$. Let $f$ be the overlap of $\widetilde{u}^{(i_0)}$ and $u^{(i_0 + 1)}$ in $Z$ ($f$ is empty if $\widetilde{u}^{(i_0)}$ and $u^{(i_0 + 1)}$ touch at a point). Let $\widehat{u}^{(i_0 + 1)} = f^{-1}\cdot u^{(i_0 + 1)}$.
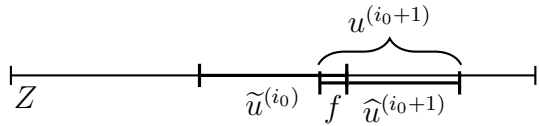
\begin{center}
\begin{tikzpicture}
\draw[|-|, black, thick] (-1.5, 0)--(5.5, 0);
\node[below, xshift=6] at (-1.5, 0) {$Z$};
\draw[|-|, black, very thick] (1,0)--(3,0) node[midway, below] {$\widetilde{u}^{(i_0)}$};
\draw[-|, black, very thick] (3,-0.1)--(4.5,-0.1);
\draw[|-|, black, very thick] (2.6,-0.1)--(3,-0.1) node[midway, below] {$f$};
\path (3, 0)--(4.5, 0) node[midway, below, yshift=-3] {$\widehat{u}^{(i_0 + 1)}$};
\draw [thick, decorate, decoration={brace, amplitude=10pt, raise=4pt}] (2.6, 0) to node[midway, above, yshift=12pt] {$u^{(i_0 + 1)}$} (4.5, 0);
\end{tikzpicture}
\end{center}
Let us perform the replacement $\widetilde{u}^{(i_0)} \leadsto \dbtilde{b}^{(i_0)}$ in $Z$ and consider the resulting monomial $Z[\widetilde{u}^{(i_0)} \leadsto \dbtilde{b}^{(i_0)}]$. Let $f_1$ be the overlap of $\dbtilde{b}^{(i_0)}$ and the image of $u^{(i_0 + 1)}$ in $Z[\widetilde{u}^{(i_0)} \leadsto \dbtilde{b}^{(i_0)}]$ ($f_1$ is empty if $\dbtilde{b}^{(i_0)}$ and the image of $u^{(i_0 + 1)}$ in $Z[\widetilde{u}^{(i_0)} \leadsto \dbtilde{b}^{(i_0)}]$ touch at a point). We can write $Z[\widetilde{u}^{(i_0)} \leadsto \dbtilde{b}^{(i_0)}]$ in the form
\begin{equation}
\label{first_different_replacement1}
Z[\widetilde{u}^{(i_0)} \leadsto \dbtilde{b}^{(i_0)}] = \widetilde{L}\dbhat{b}^{(i_0)}f_1\widehat{u}^{(i_0 + 1)}R.
\end{equation}
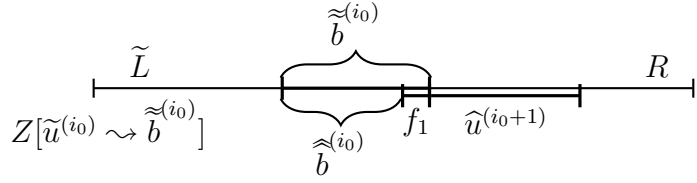
\begin{center}
\begin{tikzpicture}
\draw[|-|, black, thick] (-1.5, 0)--(6.5, 0);
\node[below, xshift=6] at (-1.5, 0) {$Z[\widetilde{u}^{(i_0)} \leadsto \dbtilde{b}^{(i_0)}]$};
\draw[|-|, black, very thick] (1,0)--(3,0);
\draw[-|, black, very thick] (3,-0.1)--(5,-0.1);
\draw[|-|, black, very thick] (2.6,-0.1)--(3,-0.1) node[midway, below] {$f_1$};
\path (3, 0)--(5, 0) node[midway, below, yshift=-3] {$\widehat{u}^{(i_0 + 1)}$};
\draw [thick, decorate, decoration={brace, amplitude=10pt, raise=4pt}] (1, 0) to node[midway, above, yshift=12pt] {$\dbtilde{b}^{(i_0)}$} (3, 0);
\draw [thick, decorate, decoration={brace, amplitude=10pt, raise=4pt, mirror}] (1, 0) to node[midway, below, yshift=-12pt] {$\dbhat{b}^{(i_0)}$} (2.6, 0);
\path (-1.5, 0)--(1, 0) node[near start, above] {$\widetilde{L}$};
\path (4.5, 0)--(6.5, 0) node[near end, above] {$R$};
\end{tikzpicture}
\end{center}
Similarly, we can write $Z[\widetilde{u}^{(i_0)} \leadsto \dbtilde{c}^{(i_0)}]$ in the form
\begin{equation}
\label{first_different_replacement2}
Z[\widetilde{u}^{(i_0)} \leadsto \dbtilde{c}^{(i_0)}] = \widetilde{L}\dbhat{c}^{(i_0)}f_2\widehat{u}^{(i_0 + 1)}R,
\end{equation}
where $f_2$ is the overlap of $\dbtilde{c}^{(i_0)}$ and the image of $u^{(i_0 + 1)}$ in $Z[\widetilde{u}^{(i_0)} \leadsto \dbtilde{c}^{(i_0)}]$ ($f_2$ is empty if $\dbtilde{c}^{(i_0)}$ and the image of $u^{(i_0 + 1)}$ in $Z[\widetilde{u}^{(i_0)} \leadsto \dbtilde{c}^{(i_0)}]$ touch at a point).
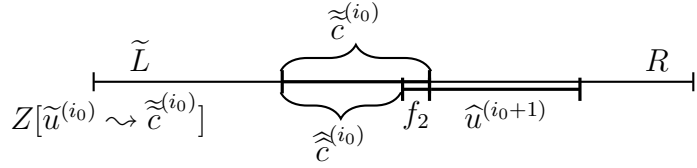
\begin{center}
\begin{tikzpicture}
\draw[|-|, black, thick] (-1.5, 0)--(6.5, 0);
\node[below, xshift=6] at (-1.5, 0) {$Z[\widetilde{u}^{(i_0)} \leadsto \dbtilde{c}^{(i_0)}]$};
\draw[|-|, black, very thick] (1,0)--(3,0);
\draw[-|, black, very thick] (3,-0.1)--(5,-0.1);
\draw[|-|, black, very thick] (2.6,-0.1)--(3,-0.1) node[midway, below] {$f_2$};
\path (3, 0)--(5, 0) node[midway, below, yshift=-3] {$\widehat{u}^{(i_0 + 1)}$};
\draw [thick, decorate, decoration={brace, amplitude=10pt, raise=4pt}] (1, 0) to node[midway, above, yshift=12pt] {$\dbtilde{c}^{(i_0)}$} (3, 0);
\draw [thick, decorate, decoration={brace, amplitude=10pt, raise=4pt, mirror}] (1, 0) to node[midway, below, yshift=-12pt] {$\dbhat{c}^{(i_0)}$} (2.6, 0);
\path (-1.5, 0)--(1, 0) node[near start, above] {$\widetilde{L}$};
\path (4.5, 0)--(6.5, 0) node[near end, above] {$R$};
\end{tikzpicture}
\end{center}

If $m > i_0 + 1$, then, clearly, $u^{(i)}$, $i = i_0 + 2, \ldots, m$, can be considered as occurrences in $Z[\widetilde{u}^{(i_0)} \leadsto \dbtilde{b}^{(i_0)}]$ and in $Z[\widetilde{u}^{(i_0)} \leadsto \dbtilde{c}^{(i_0)}]$. Clearly, $u^{(i)}$, $i = i_0 + 2, \ldots, m$, are virtual members of the chart of $Z[\widetilde{u}^{(i_0)} \leadsto \dbtilde{b}^{(i_0)}]$ and of $Z[\widetilde{u}^{(i_0)} \leadsto \dbtilde{c}^{(i_0)}]$. Let us consecutively perform the replacements $u^{(i)} \leadsto b^{(i)}$ for $i = i_0 + 2, \ldots, m$ (if $m > i_0 + 1$), starting from $Z[\widetilde{u}^{(i_0)} \leadsto \dbtilde{b}^{(i_0)}]$, and let $V$ be the resulting monomial. Let us consecutively perform the replacements $u^{(i)} \leadsto c^{(i)}$ for $i = i_0 + 2, \ldots, m$ (if $m > i_0 + 1$), starting from $Z[\widetilde{u}^{(i_0)} \leadsto \dbtilde{c}^{(i_0)}]$, and let $W$ be the resulting monomial. Then it follows from~\eqref{first_different_replacement1} and~\eqref{first_different_replacement2} that
\begin{align}
\label{before_last_replacement1}
&V = \widetilde{L}\dbhat{b}^{(i_0)}f_1m_v(u^{(i_0 + 1)})d_1R_1,\\
\label{before_last_replacement2}
&W = \widetilde{L}\dbhat{c}^{(i_0)}f_2m_v(u^{(i_0 + 1)})d_2R_2,
\end{align}
where $f_1m_v(u^{(i_0 + 1)})d_1$ is the $(i_0 + 1)$-th virtual member of the chart of $V$, $f_1m_v(u^{(i_0 + 1)})d_2$ is the $(i_0 + 1)$-th virtual member of the chart of $W$, $d_1$ is the overlap of $f_1m_v(u^{(i_0 + 1)})d_1$ and the $(i_0 + 2)$-th virtual member of the chart of $V$,  $d_2$ is the overlap of $f_2m_v(u^{(i_0 + 1)})d_2$ and the $(i_0 + 2)$-th virtual member of the chart of $W$ ($d_1$, $d_2$ may be empty).
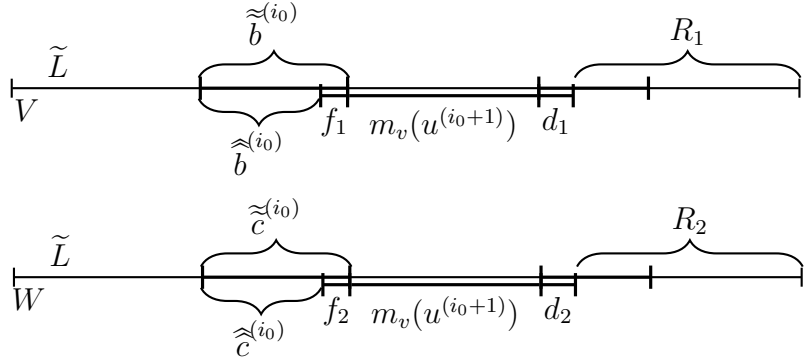
\begin{center}
\begin{tikzpicture}
\draw[|-|, black, thick] (-1.5, 0)--(9, 0);
\node[below, xshift=6] at (-1.5, 0) {$V$};
\draw[|-|, black, very thick] (1,0)--(3,0);
\draw[-, black, very thick] (3,-0.1)--(5.5,-0.1);
\draw[|-|, black, very thick] (2.6,-0.1)--(3,-0.1) node[midway, below] {$f_1$};
\path (3, 0)--(5.5, 0) node[midway, below, yshift=-3] {$m_v(u^{(i_0 + 1)})$};
\draw [thick, decorate, decoration={brace, amplitude=10pt, raise=4pt}] (1, 0) to node[midway, above, yshift=12pt] {$\dbtilde{b}^{(i_0)}$} (3, 0);
\draw [thick, decorate, decoration={brace, amplitude=10pt, raise=4pt, mirror}] (1, 0) to node[midway, below, yshift=-12pt] {$\dbhat{b}^{(i_0)}$} (2.6, 0);
\path (-1.5, 0)--(1, 0) node[near start, above] {$\widetilde{L}$};
\draw[|-|, black, very thick] (5.5,-0.1)--(6,-0.1) node[midway, below] {$d_1$};
\draw[|-|, black, very thick] (5.5,0)--(7,0);
\draw [thick, decorate, decoration={brace, amplitude=10pt, raise=4pt}] (6, 0) to node[midway, above, yshift=12pt] {$R_1$} (9, 0);
\end{tikzpicture}

\begin{tikzpicture}
\draw[|-|, black, thick] (-1.5, 0)--(9, 0);
\node[below, xshift=6] at (-1.5, 0) {$W$};
\draw[|-|, black, very thick] (1,0)--(3,0);
\draw[-, black, very thick] (3,-0.1)--(5.5,-0.1);
\draw[|-|, black, very thick] (2.6,-0.1)--(3,-0.1) node[midway, below] {$f_2$};
\path (3, 0)--(5.5, 0) node[midway, below, yshift=-3] {$m_v(u^{(i_0 + 1)})$};
\draw [thick, decorate, decoration={brace, amplitude=10pt, raise=4pt}] (1, 0) to node[midway, above, yshift=12pt] {$\dbtilde{c}^{(i_0)}$} (3, 0);
\draw [thick, decorate, decoration={brace, amplitude=10pt, raise=4pt, mirror}] (1, 0) to node[midway, below, yshift=-12pt] {$\dbhat{c}^{(i_0)}$} (2.6, 0);
\path (-1.5, 0)--(1, 0) node[near start, above] {$\widetilde{L}$};
\draw[|-|, black, very thick] (5.5,-0.1)--(6,-0.1) node[midway, below] {$d_2$};
\draw[|-|, black, very thick] (5.5,0)--(7,0);
\draw [thick, decorate, decoration={brace, amplitude=10pt, raise=4pt}] (6, 0) to node[midway, above, yshift=12pt] {$R_2$} (9, 0);
\end{tikzpicture}
\end{center}
If $m = i_0 + 1$ or the $(i_0 + 1)$-th and the $(i_0 + 2)$-th virtual members of the chart of $U$ are separated, then, obviously, $m_v(u^{(i_0 + 1)})d_1 = m_v(u^{(i_0 + 1)})d_2 = \widehat{u}^{(i_0 + 1)}$. Let us put
\begin{equation*}
v^{(i_0 + 1)} = f_1m_v(u^{(i_0 + 1)})d_1, \ w^{(i_0 + 1)} = f_2m_v(u^{(i_0 + 1)})d_2.
\end{equation*}
Let $v^{(i_0 + 1)} \leadsto x^{(i_0 + 1)}$ be the replacement in $V$ that corresponds to the replacement $u^{(i_0 + 1)} \leadsto b^{(i_0 + 1)}$ in $U$, let $w^{(i_0 + 1)} \leadsto y^{(i_0 + 1)}$ be the replacement in $W$ that corresponds to the replacement $u^{(i_0 + 1)} \leadsto c^{(i_0 + 1)}$ in $U$. Then we have
\begin{align*}
&\mu[U]\left(b^{(1)}\otimes \ldots \otimes b^{(m)}\right) = V[v^{(i_0 + 1)} \leadsto x^{(i_0 + 1)}],\\
&\mu[U]\left(c^{(1)}\otimes \ldots \otimes c^{(m)}\right) = W[w^{(i_0 + 1)} \leadsto y^{(i_0 + 1)}].
\end{align*}

Assume that
\begin{equation*}
\mu[U]\left(b^{(1)}\otimes \ldots \otimes b^{(m)}\right) = \mu[U]\left(c^{(1)}\otimes \ldots \otimes c^{(m)}\right).
\end{equation*}
Then
\begin{equation}
\label{equal_after_last}
V[v^{(i_0 + 1)} \leadsto x^{(i_0 + 1)}] = W[w^{(i_0 + 1)} \leadsto y^{(i_0 + 1)}].
\end{equation}
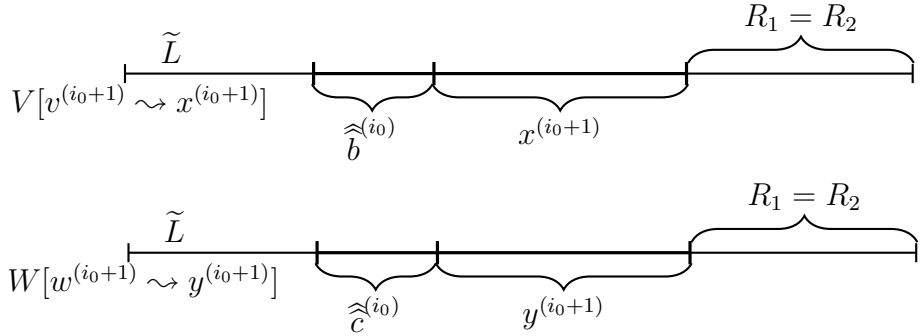
\begin{center}
\begin{tikzpicture}
\draw[|-|, black, thick] (-1.5, 0)--(9, 0);
\node[below, xshift=6] at (-1.5, 0) {$V[v^{(i_0 + 1)} \leadsto x^{(i_0 + 1)}]$};
\draw[|-, black, very thick] (1,0)--(2.6,0);
\draw[|-|, black, very thick] (2.6,0)--(6,0);
\draw [thick, decorate, decoration={brace, amplitude=10pt, raise=4pt, mirror}] (1, 0) to node[midway, below, yshift=-12pt] {$\dbhat{b}^{(i_0)}$} (2.6, 0);
\path (-1.5, 0)--(1, 0) node[near start, above] {$\widetilde{L}$};
\draw [thick, decorate, decoration={brace, amplitude=10pt, raise=4pt}] (6, 0) to node[midway, above, yshift=12pt] {$R_1 = R_2$} (9, 0);
\draw [thick, decorate, decoration={brace, amplitude=10pt, raise=4pt, mirror}] (2.6, 0) to node[midway, below, yshift=-12pt] {$x^{(i_0 + 1)}$} (6, 0);
\end{tikzpicture}

\begin{tikzpicture}
\draw[|-|, black, thick] (-1.5, 0)--(9, 0);
\node[below, xshift=6] at (-1.5, 0) {$W[w^{(i_0 + 1)} \leadsto y^{(i_0 + 1)}]$};
\draw[|-, black, very thick] (1,0)--(2.6,0);
\draw[|-|, black, very thick] (2.6,0)--(6,0);
\draw [thick, decorate, decoration={brace, amplitude=10pt, raise=4pt, mirror}] (1, 0) to node[midway, below, yshift=-12pt] {$\dbhat{c}^{(i_0)}$} (2.6, 0);
\path (-1.5, 0)--(1, 0) node[near start, above] {$\widetilde{L}$};
\draw [thick, decorate, decoration={brace, amplitude=10pt, raise=4pt}] (6, 0) to node[midway, above, yshift=12pt] {$R_1 = R_2$} (9, 0);
\draw [thick, decorate, decoration={brace, amplitude=10pt, raise=4pt, mirror}] (2.6, 0) to node[midway, below, yshift=-12pt] {$y^{(i_0 + 1)}$} (6, 0);
\end{tikzpicture}
\end{center}
Since $x^{(i_0 + 1)}$ is the $(i_0 + 1)$-th virtual member of the chart of $V[v^{(i_0 + 1)} \leadsto x^{(i_0 + 1)}]$ and $y^{(i_0 + 1)}$ is the $(i_0 + 1)$-th virtual member of the chart of $W[w^{(i_0 + 1)} \leadsto y^{(i_0 + 1)}]$, we see that $x^{(i_0 + 1)} = y^{(i_0 + 1)}$ as occurrences. That is, they are equal as words and they start at the same position in $V[v^{(i_0 + 1)} \leadsto x^{(i_0 + 1)}] = W[w^{(i_0 + 1)} \leadsto y^{(i_0 + 1)}]$. Therefore, \eqref{equal_after_last} implies that the monomials $V[v^{(i_0 + 1)} \leadsto x^{(i_0 + 1)}]$ and $W[w^{(i_0 + 1)} \leadsto y^{(i_0 + 1)}]$ have equal prefix that end at the beginning point of $x^{(i_0 + 1)} = y^{(i_0 + 1)}$. On the one hand, it follows from~\eqref{before_last_replacement1} that this prefix is equal to $\widetilde{L}\dbhat{b}^{(i_0)}$. On the other hand, it follows from~\eqref{before_last_replacement2} that this prefix is equal to $\widetilde{L}\dbhat{c}^{(i_0)}$. Therefore, $\dbhat{b}^{(i_0)} = \dbhat{c}^{(i_0)}$.

We have that $v^{(i_0 + 1)}$ and $x^{(i_0 + 1)}$ are $V$-incident monomials and $w^{(i_0 + 1)}$ and $y^{(i_0 + 1)}$ are $W$-incident monomials. Hence, $x^{(i_0 + 1)}$ and $v^{(i_0 + 1)}$ are $V[v^{(i_0 + 1)} \leadsto x^{(i_0 + 1)}]$-incident monomials and $y^{(i_0 + 1)}$ and $w^{(i_0 + 1)}$ are $W[w^{(i_0 + 1)} \leadsto y^{(i_0 + 1)}]$-incident monomials. Therefore, combining~\eqref{equal_after_last}, Lemma~\ref{glue_U_incident} and the fact that $x^{(i_0 + 1)} = y^{(i_0 + 1)}$ are virtual members of the chart of the corresponding monomials, we get that that $v^{(i_0 + 1)}$ and $w^{(i_0 + 1)}$ are $V$-incident monomials.

We have that $\widetilde{u}^{(i_0)}$ is a virtual member of the chart of $Z$, $\widetilde{u}^{(i_0)}$ and $\dbtilde{b}^{(i_0)}$ are $Z$-incident monomials, and $\widetilde{u}^{(i_0)}$ and $\dbtilde{c}^{(i_0)}$ are $Z$-incident monomials. Therefore, by Lemma~\ref{glue_U_incident}, we obtain that $\dbtilde{b}^{(i_0)}$ and $\dbtilde{c}^{(i_0)}$ are $Z[\widetilde{u}^{(i_0)} \leadsto \dbtilde{b}^{(i_0)}]$-incident monomials.

Recall that $\dbtilde{b}^{(i_0)} = \dbhat{b}^{(i_0)}f_1$ and $\dbtilde{c}^{(i_0)} = \dbhat{c}^{(i_0)}f_2$. Consider the monomials
\begin{align*}
\dbhat{b}^{(i_0)}f_1m_v(u^{(i_0 + 1)})d_1 = \dbtilde{b}^{(i_0)}m_v(u^{(i_0 + 1)})d_1,\\
\dbhat{c}^{(i_0)}f_2m_v(u^{(i_0 + 1)})d_2= \dbtilde{c}^{(i_0)}m_v(u^{(i_0 + 1)})d_2.
\end{align*}
The following properties hold.
\begin{enumerate}
\item
We proved above that $\dbtilde{b}^{(i_0)}$ and $\dbtilde{c}^{(i_0)}$ are $Z[\widetilde{u}^{(i_0)} \leadsto \dbtilde{b}^{(i_0)}]$-incident monomials. Hence, there exists a sequence of monomials $m_1^{(i_0)}, \ldots, m_{k + 1}^{(i_0)}$ that satisfies the conditions of Definition~\ref{U_incident_monomials}. In particular, $m_1^{(i_0)} = \dbtilde{b}^{(i_0)}$, $m_{k + 1}^{(i_0)} = \dbtilde{c}^{(i_0)}$, $m_t^{(i_0)}$ and $m_{t + 1}^{(i_0)}$ are incident monomials for $t = 1, \ldots, k$, and $\SPM(m_t^{(i_0)}) \geqslant \tau - 2$ for $t = 2, \ldots, k$.
\item
Since $\dbtilde{b}^{(i_0)} = m_1^{(i_0)}$ and $\dbtilde{c}^{(i_0)} = m_{k + 1}^{(i_0)}$ are virtual members of the chart of the corresponding monomials, we obtain  $\SPM(\dbtilde{b}^{(i_0)}) \geqslant \tau - 2$ and $\SPM(\dbtilde{c}^{(i_0)}) \geqslant \tau - 2$.
\item
Since $u^{(i_0 + 1)}$ is a virtual member of the chart of $U$, we obtain  $\SPM(m_v(u^{(i_0 + 1)})) \geqslant \tau - 2$.
\item
$f_1m_v(u^{(i_0 + 1)})$ is a maximal occurrence in $\dbhat{b}^{(i_0)}f_1m_v(u^{(i_0 + 1)})$;

$f_2m_v(u^{(i_0 + 1)})$ is a maximal occurrence in $\dbhat{c}^{(i_0)}f_2m_v(u^{(i_0 + 1)})$.
\item
$d_1$ and $d_2$ are small pieces.
\item
We proved above that $v^{(i_0 + 1)} = f_1m_v(u^{(i_0 + 1)})d_1$ and $w^{(i_0 + 1)} = f_2m_v(u^{(i_0 + 1)})d_2$ are $V$-incident monomials. Therefore, there exists a sequence of monomials $m_1^{(i_0 + 1)}, \ldots, m_{l + 1}^{(i_0 + 1)}$ that satisfies the conditions of Definition~\ref{U_incident_monomials}. In particular, $m_1^{(i_0 + 1)} = v^{(i_0 + 1)}$, $m_{l + 1}^{(i_0 + 1)} = w^{(i_0 + 1)}$, $m_t^{(i_0 + 1)}$ and $m_{t + 1}^{(i_0 + 1)}$ are incident monomials for $t = 1, \ldots, l$, and $\SPM(m_t^{(i_0 + 1)}) \geqslant \tau - 2$ for $t = 2, \ldots, l$.
\item
Since $v^{(i_0 + 1)} = f_1m_v(u^{(i_0 + 1)})d_1$ and $w^{(i_0 + 1)} = f_2m_v(u^{(i_0 + 1)})d_2$ are virtual members of the chart of the corresponding monomials, we obtain  $\SPM(v^{(i_0 + 1)}) \geqslant \tau - 2$ and $\SPM(w^{(i_0 + 1)}) \geqslant \tau - 2$.
\end{enumerate}
So, the monomials $\dbtilde{b}^{(i_0)}$, $\dbtilde{c}^{(i_0)}$, together with $m_v(u^{(i_0 + 1)})$ satisfy the initial conditions of right-sided Isolation Axiom. However, we proved above that
\begin{equation*}
\dbtilde{b}^{(i_0)}\cdot f_1^{-1} = \dbhat{b}^{(i_0)} = \dbhat{c}^{(i_0)} = \dbtilde{c}^{(i_0)}\cdot f_2^{-1}.
\end{equation*}
Therefore, $\dbtilde{b}^{(i_0)}$, $\dbtilde{c}^{(i_0)}$, together with $m_v(u^{(i_0 + 1)})$ violate right-sided Isolation Axiom. A contradiction. Thus,
\begin{equation*}
\mu[U]\left(b^{(1)}\otimes \ldots \otimes b^{(m)}\right) \neq \mu[U]\left(c^{(1)}\otimes \ldots \otimes c^{(m)}\right).
\end{equation*}
Statement~\ref{mu_properties_up_correspondence} is proved.

~\paragraph*{(3)} Let $1 \leqslant i_0 \leqslant m$ be a position in the chart of $U$. We consider an element $a^{(1)}\otimes\ldots\otimes t^{(i_0)}\otimes\ldots \otimes a^{(m)} \in A_1[U]\otimes\ldots\otimes D_{i_0}[U]\otimes\ldots\otimes A_m[U]$ such that $a^{(i)} \in ME_i[U]$ for $i \neq i_0$, and $t^{(i_0)} \in D_{i_0}[U] \cap \Rel$.

Assume $t^{(i_0)} = \sum_{j= 1}^k \beta_j b_j^{(i_0)}$. Let
\begin{equation*}
W_j = \mu[U]\left(a^{(1)}\otimes\ldots\otimes b^{(i_0)}_j\otimes\ldots\otimes a^{(m)}\right),\ j = 1, \ldots, k.
\end{equation*}
Clearly, if $b^{(i_0)}_j \in ME_{i_0}[U]$, then $W_j$ is defined using \ref{mu_case_all_up} of Defintion~\ref{mu_def}; if $b^{(i_0)}_j \in ML_{i_0}[U]$, then $W_j$ is defined using \ref{mu_case_one_low} of Defintion~\ref{mu_def}. Notice that in both cases the replacement in the position $i_0$ can be performed last.

Let us consecutively perform the replacements $u^{(i)} \leadsto a^{(i)}$ for $i \neq i_0$, starting from $U$, and let $Z$ be the resulting monomial. Let $\widetilde{u}^{(i_0)}$ be the $i_0$-th virtual member of the chart of $Z$. Then, according to Remark~\ref{consecurive_repl_basic_properties}, we have
\begin{align*}
&u^{(i_0)} = em_v(u^{(i_0)})f,\ \widetilde{u}^{(i_0)} = \widetilde{e}m_v(u^{(i_0)})\widetilde{f},\textit{ and}\\
&\widetilde{u}^{(i_0)} = \widetilde{e}\left(e^{-1}\cdot u^{(i_0)} \cdot f^{-1}\right)\widetilde{f},
\end{align*}
where $e, f, \widetilde{e}, \widetilde{f}$ are small pieces. Let $\widetilde{u}^{(i_0)} \leadsto \dbtilde{b}^{(i_0)}_j$ be the replacement in $Z$ that corresponds to the replacement $u^{(i_0)} \leadsto b^{(i_0)}_j$ in $U$, $j = 1, \ldots, k$. Then, by definition,
\begin{equation*}
\dbtilde{b}^{(i_0)}_j = \widetilde{e} \cdot e^{-1}\cdot b^{(i_0)}_j \cdot f^{-1}\cdot \widetilde{f},\ j = 1, \ldots, k
\end{equation*}
(see also Remark~\ref{consecurive_repl_basic_properties}). It follows from Small Cancellation Axiom that $t^{(i_0)}= \sum_{j= 1}^k \beta_j b_j^{(i_0)}$ contains a monomial of $\SPM$-measure $\geqslant \tau + 1$. Without loss of generality, we can assume that $\SPM(b_1^{(i_0)}) \geqslant \tau + 1$. Then $b_1^{(i_0)}$ is a virtual member of the chart of $U[u^{(i_0)}\leadsto b^{(i_0)}]$. Hence, Lemma~\ref{same_small_changes} implies
\begin{align*}
&b^{(i_0)}_1 = e^{\prime}m_v(b^{(i_0)}_1)f^{\prime},\ \dbtilde{b}^{(i_0)}_1 = \widetilde{e}^{\prime}m_v(b^{(i_0)}_1)\widetilde{f}^{\prime},\\
&\dbtilde{b}^{(i_0)}_1 = \widetilde{e}^{\prime} \left({e^{\prime}}^{-1}\cdot b^{(i_0)}_1 \cdot {f^{\prime}}^{-1}\right)\cdot \widetilde{f}^{\prime},\\
&\widetilde{e} \cdot e^{-1} = \widetilde{e}^{\prime} \cdot {e^{\prime}}^{-1},\  f^{-1}\cdot \widetilde{f} = {f^{\prime}}^{-1}\cdot \widetilde{f}^{\prime},
\end{align*}
where $e^{\prime}, f^{\prime}, \widetilde{e}^{\prime}, \widetilde{f}^{\prime}$ are small pieces.
Hence, we have
\begin{equation*}
\dbtilde{b}^{(i_0)}_j = \widetilde{e}^{\prime} \cdot {e^{\prime}}^{-1}\cdot b^{(i_0)}_j \cdot {f^{\prime}}^{-1}\cdot \widetilde{f}^{\prime},
\end{equation*}
and we obtain
\begin{align}
\begin{split}
\label{t_t_tilde_correspondnce}
\widetilde{t}^{(i_0)} &= \sum_{j = 1}^k \beta_j \dbtilde{b}^{(i_0)}_j = \sum_{j = 1}^k\widetilde{e}^{\prime} \cdot{e^{\prime}}^{-1}\cdot b^{(i_0)}_j \cdot {f^{\prime}}^{-1}\cdot \widetilde{f}^{\prime} =\\
&= \widetilde{e}^{\prime} \cdot{e^{\prime}}^{-1} \cdot t^{(i_0)}\cdot {f^{\prime}}^{-1} \cdot \widetilde{f}^{\prime}.
\end{split}
\end{align}

Since $b^{(i_0)}_1 = e^{\prime}m_v(b^{(i_0)}_1)f^{\prime}$ is of $\SPM$-measure $\geqslant \tau + 1$ and $e^{\prime}$ and $f^{\prime}$ are small pieces, we see that $m_v(b^{(i_0)}_1) = {e^{\prime}}^{-1}\cdot b^{(i_0)}_1 \cdot {f^{\prime}}^{-1}$ is not a small piece. Lemma~\ref{same_small_changes} implies  $\dbtilde{b}^{(i_0)}_1 = \widetilde{e}^{\prime}m_v(b^{(i_0)}_1)\widetilde{f}^{\prime} \in \Mon$. Therefore, it follows from Corollary~\ref{not_sp_cancellation_prolongation} that
\begin{equation*}
\widetilde{t}^{(i_0)} = \sum_{j = 1}^k \beta_j \dbtilde{b}^{(i_0)}_j = \widetilde{e}^{\prime} \cdot{e^{\prime}}^{-1} \cdot t^{(i_0)}\cdot {f^{\prime}}^{-1} \cdot \widetilde{f}^{\prime} \in \Rel.
\end{equation*}

Let $Z = L\widetilde{u}^{(i_0)}R$. Since $Z$ is the resulting monomial of the consecutive replacements $u^{(i)} \leadsto a^{(i)}$ for $i \neq i_0$, starting from $U$, and $\widetilde{u}^{(i_0)} \leadsto \dbtilde{b}^{(i_0)}_j$ is the replacement in $Z$ that corresponds to the replacement $u^{(i_0)} \leadsto b^{(i_0)}_j$ in $U$, we obtain  $W_j = Z[\widetilde{u}^{(i_0)} \leadsto \dbtilde{b}^{(i_0)}_j] = L\dbtilde{b}^{(i_0)}_jR$. Therefore,
\begin{equation}
\label{image_of_dependence_formula}
\mu[U]\left(a^{(1)}\otimes\ldots\otimes t^{(i_0)}\otimes\ldots \otimes a^{(m)}\right) = \sum_{j = 1}^k \beta_jW_j = \sum_{j = 1}^k \beta_j L\dbtilde{b}^{(i_0)}_jR.
\end{equation}

Recall that we assumed that $\SPM(b_1^{(i_0)}) \geqslant \tau + 1$. Then $b_1^{(i_0)}$ is a virtual member of the chart of $U[u^{(i_0)} \leadsto b_1^{(i_0)}]$. Therefore, $b_1^{(i_0)} \in ME_{i_0}[U]$. Since $a^{(i)}\in ME_i[U]$ for $i\neq i_0$, statement~\ref{mu_properties_up_correspondence} of Lemma~\ref{mu_properties} implies
\begin{equation*}
W_1 = \mu[U]\left(a^{(1)}\otimes\ldots\otimes b^{(i_0)}_1\otimes\ldots\otimes a^{(m)}\right) \in \DMUp{U}.
\end{equation*}
It follows from statements~\ref{virtual_members_U_incident_replacements3} and~\ref{virtual_members_U_incident_replacements6} Corollary~\ref{virtual_members_U_incident_replacements} that $\dbtilde{b}^{(i_0)}_1$ is the $i_0$-th virtual member of the chart of $W_1 = L\dbtilde{b}^{(i_0)}_1R$ (see also Remark~\ref{consecurive_repl_basic_properties}). Combining this with~\eqref{image_of_dependence_formula}, we see that $\sum_{j = 1}^k \beta_jW_j$ is a layout of a multi-turn of the $i_0$-th virtual member of the chart of the monomial $W_1 \in \DMUp{U}$. Thus, $\mu[U]\left(a^{(1)}\otimes\ldots\otimes t^{(i_0)}\otimes\ldots \otimes a^{(m)}\right) \in \GDp^{(i_0)}[U]$.

~\paragraph*{(4)} Let $1 \leqslant i_0 \leqslant m$ be a position in the chart of $U$. Let $T \in \GDp^{(i_0)}[U]$. Then
\begin{equation*}
T = \sum_{j = 1}^k \beta_j W_j = \sum_{j = 1}^k\beta_j L\dbtilde{b}^{(i_0)}_jR,
\end{equation*}
where $\sum_{j = 1}^k\beta_j \dbtilde{b}^{(i_0)}_j \in \Rel$, $W_h \in \DMUp{U}$ for some $1 \leqslant h \leqslant k$, and $\dbtilde{b}_{h}^{(i_0)}$ is the $i_0$-th virtual member of the chart of $W_h$. Without loss of generality we can assume that $h = 1$. So, $W_1 \in \DMUp{U}$.

Since $W_1 \in \DMUp{U}$, it follows from the statement~\ref{mu_properties_up_correspondence} of Lemma~\ref{mu_properties} that there exists an element $a^{(1)}\otimes \ldots \otimes b_1^{(i_0)}\otimes\ldots\otimes a^{(m)} \in A_1[U]\otimes\ldots\otimes A_m[U]$ such that every $a^{(i)} \in ME_i[U]$ for $i \neq i_0$, $b_1^{(i_0)} \in ME_{i_0}[U]$, and
\begin{equation*}
W_1 = L\dbtilde{b}_1^{(i_0)}R = \mu[U]\left(a^{(1)}\otimes \ldots \otimes b_1^{(i_0)}\otimes\ldots\otimes a^{(m)}\right).
\end{equation*}

Let us consecutively perform the replacements $u^{(i)}\leadsto a^{(i)}$, $i \neq i_0$. Let $Z$ be the resulting monomial. Assume $\widetilde{u}^{(i_0)}$ is the $i_0$-th virtual member of the chart of $Z$. Let $\widetilde{u}^{(i_0)} \leadsto c^{(i_0)}$ be the replacement in $Z$ that corresponds to the replacement $u^{(i_0)} \leadsto b_1^{(i_0)}$ in $U$. By the \ref{mu_case_all_up} of Definition~\ref{mu_def}, we obtain
\begin{equation*}
Z[\widetilde{u}^{(i_0)} \leadsto c^{(i_0)}] = \mu[U]\left(a^{(1)}\otimes \ldots \otimes b_1^{(i_0)}\otimes\ldots\otimes a^{(m)}\right) = W_1 = L\dbtilde{b}^{(i_0)}_1R.
\end{equation*}
Since $c^{(i_0)}$ is the $i_0$-th virtual member of the chart of $Z[\widetilde{u}^{(i_0)} \leadsto c^{(i_0)}]$ and $\dbtilde{b}^{(i_0)}_1$ is the $i_0$-th virtual member of the chart of $L\dbtilde{b}^{(i_0)}_1R$, the last equality implies  $c^{(i_0)} = \dbtilde{b}^{(i_0)}_1$ and $Z = L\widetilde{u}^{(i_0)}R$.

It follows from Remark~\ref{consecurive_repl_basic_properties} and Lemma~\ref{same_small_changes} that
\begin{align}
&u^{(i_0)} = em_v(u^{(i_0)})f,\ \widetilde{u}^{(i_0)} = \widetilde{e}m_v(u^{(i_0)})\widetilde{f},\nonumber\\
&\widetilde{u}^{(i_0)} = \widetilde{e}\left(e^{-1}\cdot u^{(i_0)} \cdot f^{-1}\right)\widetilde{f},\nonumber\\
&b^{(i_0)}_1 = e^{\prime}m_v(b^{(i_0)}_1)f^{\prime},\ \dbtilde{b}^{(i_0)}_1 = \widetilde{e}^{\prime}m_v(b^{(i_0)}_1)\widetilde{f}^{\prime},\nonumber\\
&\dbtilde{b}^{(i_0)}_1 = \widetilde{e}^{\prime} \left({e^{\prime}}^{-1}\cdot b^{(i_0)}_1 \cdot {f^{\prime}}^{-1}\right)\cdot \widetilde{f}^{\prime},\nonumber\\
\label{equal_trnsformation}
&\widetilde{e} \cdot e^{-1} = \widetilde{e}^{\prime} \cdot {e^{\prime}}^{-1},\  f^{-1}\cdot \widetilde{f} = {f^{\prime}}^{-1}\cdot \widetilde{f}^{\prime},
\end{align}
where $e, f, \widetilde{e}, \widetilde{f}, e^{\prime}, f^{\prime}, \widetilde{e}^{\prime}, \widetilde{f}^{\prime}$ are small pieces. Let us put
\begin{equation*}
b^{(i_0)}_j = e^{\prime}\cdot \widetilde{e^{\prime}}^{-1}\cdot \dbtilde{b}^{(i_0)}_j \cdot \widetilde{f^{\prime}}^{-1}\cdot f^{\prime},\ j = 2, \ldots, k,
\end{equation*}
and consider
\begin{equation*}
t^{(i_0)} = \sum_{j = 1}^k b^{(i_0)}_j = \sum_{j = 1}^k e^{\prime}\cdot \widetilde{e^{\prime}}^{-1}\cdot \dbtilde{b}^{(i_0)}_j \cdot \widetilde{f^{\prime}}^{-1}\cdot f^{\prime}.
\end{equation*}

Since $\dbtilde{b}^{(i_0)}_1$ is a virtual member of the chart of $W_1$, we have $\SPM(\dbtilde{b}^{(i_0)}_1) \geqslant \tau - 2$. Since $\widetilde{e}^{\prime}$ and $\widetilde{f}^{\prime}$ are small pieces, we see that $m_v(b^{(i_0)}_1) = \widetilde{e^{\prime}}^{-1}\cdot \dbtilde{b}^{(i_0)}_1 \cdot \widetilde{f^{\prime}}^{-1}$ is not a small piece. So, since $\sum_{j = 1}^k\beta_j \dbtilde{b}^{(i_0)}_j \in \Rel$ and $b^{(i_0)}_1 = e^{\prime}m_v(b^{(i_0)}_1)f^{\prime} \in \Mon$, it follows from Corollary~\ref{not_sp_cancellation_prolongation} that
\begin{equation*}
t^{(i_0)} = \sum_{j = 1}^k e^{\prime}\cdot \widetilde{e^{\prime}}^{-1}\cdot \dbtilde{b^{\prime}}^{(i_0)}_j \cdot \widetilde{f^{\prime}}^{-1}\cdot f^{\prime} \in \Rel.
\end{equation*}

Since $t^{(i_0)} \in \Rel$, we see that $b_j^{(i_0)} \in \Mon$ for all $j = 1, \ldots, k$, and $b^{(i_0)}_1$ and $b^{(i_0)}_j$ are incident monomials for $j = 2,\ldots, k$. Since $b^{(i_0)}_1 \in A_{i_0}[U]$, by definition, $u^{(i_0)}$ and $b^{(i_0)}_1$ are $U$-incident monomials. Combining these two observations with the fact that $b^{(i_0)}_1$ is a virtual member of the chart of $U[u^{(i_0)} \leadsto b^{(i_0)}_1]$, we see that $u^{(i_0)}$ and $b^{(i_0)}_j$ are $U$-incident monomials for $j = 2,\ldots, k$. That is, $b^{(i_0)}_j \in A_{i_0}[U]$.

Let us show that
\begin{equation*}
W_j = L\dbtilde{b}^{(i_0)}_jR = \mu[U]\left(a^{(1)}\otimes \ldots \otimes b_j^{(i_0)}\otimes\ldots\otimes a^{(m)}\right), j = 2, \ldots, k.
\end{equation*}
Recall that $Z$ is the resulting monomial of the consecutive replacements $u^{(i)}\leadsto a^{(i)}$, $i \neq i_0$, starting from $U$, and $\widetilde{u}^{(i_0)}$ is the $i_0$-th virtual member of the chart of $Z$, $Z = L\widetilde{u}^{(i_0)}R$. Consider the replacement in $Z$ that corresponds to the replacement $u^{(i_0)} \leadsto b^{(i_0)}_j$ in $U$. Since $\widetilde{u}^{(i_0)} = \widetilde{e}\left(e^{-1}\cdot u^{(i_0)} \cdot f^{-1}\right)\widetilde{f}$, the corresponding replacement in $Z$ is of the form
\begin{equation*}
\widetilde{u}^{(i_0)} \leadsto \widetilde{e} \cdot e^{-1}\cdot b^{(i_0)}_j \cdot  f^{-1} \cdot \widetilde{f}.
\end{equation*}
By the equality $b^{(i_0)}_j = e^{\prime}\cdot \widetilde{e^{\prime}}^{-1}\cdot \dbtilde{b}^{(i_0)}_j \cdot \widetilde{f^{\prime}}^{-1}\cdot f^{\prime}$ and formula~\eqref{equal_trnsformation}, we have
\begin{equation*}
\widetilde{e} \cdot e^{-1}\cdot b^{(i_0)}_j \cdot  f^{-1} \cdot \widetilde{f} = \widetilde{e} \cdot e^{-1}\cdot \left(e^{\prime}\cdot \widetilde{e^{\prime}}^{-1}\cdot \dbtilde{b}^{(i_0)}_j \cdot \widetilde{f^{\prime}}^{-1}\cdot f^{\prime}\right)\cdot  f^{-1} \cdot \widetilde{f} = \dbtilde{b}^{(i_0)}_j.
\end{equation*}
Hence, the replacement $\widetilde{u}^{(i_0)} \leadsto \dbtilde{b}^{(i_0)}_j$ in $Z$ corresponds to the replacement $u^{(i_0)} \leadsto b^{(i_0)}_j$ in $U$. So, it follows from \ref{mu_case_all_up} and \ref{mu_case_one_low} of Definition~\ref{mu_def} that
\begin{equation*}
\mu[U]\left(a^{(1)}\otimes \ldots \otimes b_j^{(i_0)}\otimes\ldots\otimes a^{(m)}\right) = Z[\widetilde{u}^{(i_0)} \leadsto \dbtilde{b}^{(i_0)}_j].
\end{equation*}
Since $Z = L\widetilde{u}^{(i_0)}R$, we obtain
\begin{equation*}
\mu[U]\left(a^{(1)}\otimes \ldots \otimes b_j^{(i_0)}\otimes\ldots\otimes a^{(m)}\right) = Z[\widetilde{u}^{(i_0)} \leadsto \dbtilde{b}^{(i_0)}_j] = L\dbtilde{b}^{(i_0)}_jR = W_j.
\end{equation*}

Since $t^{(i_0)} = \sum_{j = 1}^k b^{(i_0)}_j \in \Rel$ and $b^{(i_0)}_j \in A_{i_0}[U]$ for all $j = 1,\ldots, k$, we have $t_{i_0} \in D_{i_0}[U] \cap \Rel$.  It follows from the above that
\begin{equation*}
\mu[U]\left(a^{(1)}\otimes\ldots\otimes t^{(i_0)} \otimes\ldots\otimes a^{(m)}\right) = \sum_{j = 1}^k \beta_j W_j = T.
\end{equation*}
Thus, $a^{(1)}\otimes\ldots\otimes t^{(i_0)} \otimes\ldots\otimes a^{(m)}$ is the desired element.

~\paragraph*{(5)} Let $a^{(1)}\otimes\ldots\otimes t^{(i_1)}\otimes\ldots\otimes a^{(m)} \in A_1[U]\otimes\ldots\otimes D_{i_1}[U]\otimes\ldots\otimes A_m[U]$ such that $a^{(i)} \in ME_i[U]$ for $i \neq i_1, i_2$, and $a^{(i_2)} \in ML_{i_2}[U]$. First we consider the case $i_1 < i_2$.

Since $\Rel \cap D_{i_1}[U]$ is a set of generators of $D_{i_1}[U]$ and $\GDp^{\prime}(\Ft_{n - 1}(\Galg))$ is a set of generators of $\Dp(\Ft_{n - 1}(\Galg))$, it is enough to prove that
 \begin{equation*}
\mu[U]\left(a^{(1)}\otimes\ldots\otimes t^{(i_1)}\otimes\ldots\otimes a^{(m)}\right) \in \GDp^{\prime}(\Ft_{n - 1}(\Galg)) \textit{ for } t^{(i_1)} \in \Rel\cap D_{i_1}[U].
 \end{equation*}

So, we assume that $t^{(i_1)} \in \Rel\cap D_{i_1}[U]$. Let $t^{(i_1)}= \sum_{j = 1}^k\beta_jb_j^{(i_1)}$, and let
\begin{equation*}
W_j = \mu[U]\left(a^{(1)}\otimes\ldots\otimes b^{(i_1)}_j\otimes\ldots\otimes a^{(m)}\right).
\end{equation*}
If $b^{(i_1)}_j \in ME_{i_1}[U]$, then $W_j$ is defined using \ref{mu_case_one_low} of Definition~\ref{mu_def}. If $b^{(i_1)}_j \in ML_{i_1}[U]$, then $W_j$ is defined using \ref{mu_case_two_low} of Definition~\ref{mu_def}.

Let us consecutively do the replacements $u^{(i)} \leadsto a^{(i)}$ for $i \neq i_1, i_2$. Let $Z$ be the resulting monomial. Let $\widetilde{u}^{(i_1)}$ be the $i_1$-th virtual member of the chart of $Z$. Then, according to Remark~\ref{consecurive_repl_basic_properties}, we see that
\begin{align*}
&u^{(i_1)} = em_v\left(u^{(i_1)}\right)f,\ \widetilde{u}^{(i_1)} = \widetilde{e}m_v\left(u^{(i_1)}\right)\widetilde{f},\textit{ and}\\
&\widetilde{u}^{(i_1)} = \widetilde{e}\left(e^{-1}\cdot u^{(i_1)} \cdot f^{-1}\right)\widetilde{f},
\end{align*}
where $e, f, \widetilde{e}, \widetilde{f}$ are small pieces. Then the replacement in $Z$ that corresponds to the replacement $u^{(i_1)} \leadsto b^{(i_1)}_j$ in $U$ is of the form $\widetilde{u}^{(i_1)} \leadsto \dbtilde{b}^{(i_1)}_j$, where
\begin{equation*}
\dbtilde{b}^{(i_1)}_j = \widetilde{e} \cdot e^{-1}\cdot b^{(i_1)}_j \cdot f^{-1}\cdot \widetilde{f},\ j = 1, \ldots, k
\end{equation*}
(see Remark~\ref{consecurive_repl_basic_properties}).

Let $\widetilde{u}^{(i_2)}$ be the $i_2$-th virtual member of the chart of $Z$. Let $\widetilde{u}^{(i_2)} \leadsto \dbtilde{a}^{(i_2)}$ be the replacement in $Z$ that corresponds to the replacement $u^{(i_2)} \leadsto a^{(i_2)}$ in $U$. Let $Z = L\widetilde{u}^{(i_1)}R$. Let $\widehat{u}^{(i_2)}$ be the intersection of $\widetilde{u}^{(i_2)}$ with $R$. As before, $\widehat{u}^{(i_2)}$ can be considered as an occurrence in the monomial $Z[\widetilde{u}^{(i_1)} \leadsto \dbtilde{b}^{(i_1)}_j] = L\dbtilde{b}^{(i_1)}_jR$. Let $\widehat{u}^{(i_2)} \leadsto \dbhat{a}^{(i_2)}$ be the replacement in $R$ that corresponds to the replacement $\widetilde{u}^{(i_2)} \leadsto \dbtilde{a}^{(i_2)}$ in $Z$ (see~\eqref{corresponding_replacement_in_intersection_right}).

Assume $b_j^{(i_1)} \in ML_{i_1}[U]$, that is, $W_j$ is defined using \ref{mu_case_two_low} of Definition~\ref{mu_def}. Assume that we start with the replacement in the position $i_1$. Then, by definition, in order to obtain $W_j$, first we replace $\widehat{u}^{(i_2)}$ in $R$ by $\dbhat{a}^{(i_2)}$. Denote the corresponding resulting monomial $R[\widehat{u}^{(i_2)} \leadsto \dbhat{a}^{(i_2)}]$ by $\widetilde{R}$. After that we do the cancellations in the monomial $L\dbtilde{b}^{(i_1)}_j\cdot \widetilde{R}$.

Assume $b_j^{(i_1)} \in ME_{i_1}[U]$, that is, $W_j$ is defined using \ref{mu_case_one_low} on Definition~\ref{mu_def}. Then, clearly, the last replacement in~\ref{mu_case_one_low} of Defintion~\ref{mu_def} (in the monomial $L\dbtilde{b}^{(i_1)}_jR$) can be represented as the replacement $\widehat{u}^{(i_2)} \leadsto \dbhat{a}^{(i_2)}$ in $R$ and the further cancellations in the monomial
\begin{equation*}
\left(L\dbtilde{b}^{(i_1)}_j\right)\cdot R[\widehat{u}^{(i_2)} \leadsto \dbhat{a}^{(i_2)}] = \left(L\dbtilde{b}^{(i_1)}_j\right)\cdot \widetilde{R}
\end{equation*}
(see also Remark~\ref{first_monomial_reduced_difficulties}).

So, in both cases (when $b_j^{(i_1)} \in ML_{i_1}[U]$ and when $b_j^{(i_1)} \in ME_{i_1}[U]$) we have $W_j = \left(L\dbtilde{b}^{(i_1)}_j\right)\cdot \widetilde{R}$.

Using the same argument as in statement~\ref{mu_properties_up_dep_image} of Lemma~\ref{mu_properties}, we obtain  $\sum_{j = 1}^k \beta_j\dbtilde{b}_j^{(i_1)} \in \Rel$. By the same argument as in Proposition~\ref{the_ideal_characterisation}, we see that
\begin{equation*}
\sum_{j = 1}^k\beta_j W_j =  \sum_{j = 1}^k\beta_jL\dbtilde{b}_j^{(i_1)}\cdot\widetilde{R} \in \GDp^{\prime}(\Galg).
\end{equation*}
We noticed in Definition~\ref{mu_def} that in \ref{mu_case_one_low}---\ref{mu_case_many_low} resulting monomials belong to $\Ft_{n - 1}(\Galg)$. Therefore, every $W_j \in \Ft_{n - 1}(\Galg)$. As a result we see that
\begin{equation*}
\mu[U]\left(a^{(1)}\otimes\ldots\otimes t^{(i_1)}\otimes\ldots\otimes a^{(m)}\right) = \sum_{j = 1}^k\beta_j W_j \in \Dp(\Ft_{n - 1}(\Galg)).
\end{equation*}

Recall that in~\ref{mu_case_two_low} of Definition~\ref{mu_def} it does not matter whether we start two last replacements from the position with a smaller number or with a bigger number. So, the case $i_1 > i_2$ is considered in the same way as above (but in this case we need to deal with the intersection of $\widetilde{u}^{(i_2)}$ with $L$ in $Z = L \widetilde{u}^{(i_1)}R$).

Lemma~\ref{mu_properties} is proved.
\end{proof}

\begin{remark}
\label{isolation_axiom_remark}
Throughout this paper we use Isolation Axiom only in the argument of statement~\ref{mu_properties_up_correspondence} of Lemma~\ref{mu_properties}, in order to prove that the corresponding mapping is injective. Let us notice that, of course, Isolation Axiom is only a sufficient condition of injectivity of this mapping, not a criterion. One can state another sufficient condition if it would be more suitable in a particular case. We state the current version of Isolation Axiom because it sounds transparent, not too difficult for verification in a particular case, and holds for particular cases that we are interested in (see Section~\ref{examples_section}).
\end{remark}

Let us go back to consideration of the space $\Dp\DSpace{U}$. Assume $T \in \GDp^{\prime}(\DSpace{U})$. Then, by definition, $T$ is a layout of a multi-turn of a virtual member of the chart of some monomial $Z \in\DSpace{U}$. If $Z \in \DMUp{U}$, then it follows from Lemma~\ref{intersection_with_low_equality} that $T \notin \Low\DSpace{U}$. Assume $Z \in \Low\DSpace{U}$. Since $\Low\DSpace{U}$ is closed under taking derived monomials, we see that $T \in \Dp(\Low\DSpace{U})$. Therefore, since $\Dp\DSpace{U} = \langle \GDp^{\prime}(\DSpace{U})\rangle$, we obtain
\begin{align*}
&\Dp\DSpace{U} = \Dp(\Low\DSpace{U}) + \EDp\DSpace{U},\\
&\textit{where }\EDp\DSpace{U} = \left\langle \lbrace T \mid T \in \GDp^{\prime}(\DSpace{U}), T \notin \Low\DSpace{U} \rbrace\right\rangle.
\end{align*}
That is, by definition,  the linear space $\EDp\DSpace{U}$\label{edp_def} is generated by all the layouts of all the multi-turns of monomials of $\DMUp{U}$.

However, in order to prove Lemma~\ref{fall_through_linear_dep_one}, we consider a bigger set of generators of $\EDp\DSpace{U}$. Namely, for every $1 \leqslant i \leqslant m = \nvirt{U}$ let us consider the set of elements as follows:\label{tilde_gdp_def}
\begin{align*}
\widetilde{\GDp}^{(i)}[U] = \bigg\lbrace \mu[U]\left(a^{(1)}\otimes\right.&\left. \ldots \otimes \widetilde{t}^{(i)}\otimes \ldots \otimes a^{(m)}\right) \mid\\
&a^{(i^{\prime})} \in ME_{i^{\prime}}[U] \textit{ if } i^{\prime} \neq i,\  \widetilde{t}^{(i)} \in D_i[U] \setminus \lbrace 0\rbrace\bigg\rbrace.
\end{align*}
Let us put
\begin{equation*}
\widetilde{\GDp}[U] = \bigcup_{i = 1}^{\nvirt{U}} \widetilde{\GDp}^{(i)}[U].
\end{equation*}

\begin{lemma}
\label{EDp_bigger_set_of_generators}
The following properties hold.
\begin{enumerate}[label={(\arabic*)}]
\item
\label{EDp_bigger_set_of_generators1}
For every $T \in \widetilde{\GDp}[U]$ we have $T \notin \Low\DSpace{U}$.
\item
\label{EDp_bigger_set_of_generators2}
Let $\EDp\DSpace{U} = \left\langle \left\lbrace T \mid T \in \GDp^{\prime}(\DSpace{U}), T \notin \Low\DSpace{U} \right\rbrace\right\rangle$.\\Then $\EDp\DSpace{U} = \left\langle \widetilde{\GDp}[U]\right\rangle$.
\end{enumerate}
\end{lemma}
\begin{proof}
Let us prove statement~\ref{EDp_bigger_set_of_generators1}. Let $T \in \widetilde{\GDp}[U]$. Then, by definition,
\begin{equation*}
T = \mu[U]\left(a^{(1)}\otimes \ldots \otimes \widetilde{t}^{(i)}\otimes \ldots \otimes a^{(m)}\right),
\end{equation*}
where $a^{(i^{\prime})} \in ME_{i^{\prime}}[U]$ for $i^{\prime} \neq i$, and $\widetilde{t}^{(i)} \in D_i[U]$, $\widetilde{t}^{(i)} \neq 0$. Let
\begin{equation}
\label{dependence_one_position}
\widetilde{t}^{(i)} = \sum_{k = 1}^l \beta_k b^{(i)}_k,
\end{equation}
where $b^{(i)}_k$ are monomials, $b^{(i)}_k \in A_i[U]$ and the sum in the right-hand side is additively reduced. Then we have
\begin{align}
\begin{split}
\label{dep_inverse_image}
T = \mu[U]\left(a^{(1)}\otimes \ldots \otimes \widetilde{t}^{(i)}\otimes \ldots \otimes a^{(m)}\right)& =\\
\sum_{k = 1}^l\beta_k\mu[U]\left(a^{(1)}\right.&\left.\otimes \ldots \otimes b^{(i)}_k\otimes \ldots \otimes a^{(m)}\right).
\end{split}
\end{align}
It follows from Lemma~\ref{mu_properties} that $\mu[U]\left(a^{(1)}\otimes \ldots \otimes b^{(i)}_k\otimes \ldots \otimes a^{(m)}\right)$ belongs to $\DMUp{U}$ if and only if $b^{(i)}_k \in ME_i[U]$. Moreover, if $b^{(i)}_{k_1}, b^{(i)}_{k_2} \in ME_i[U]$ and $b^{(i)}_{k_1}\neq b^{(i)}_{k_2}$, then
\begin{equation*}
\mu[U]\left(a^{(1)}\otimes \ldots \otimes b^{(i)}_{k_1}\otimes \ldots \otimes a^{(m)}\right) \neq \mu[U]\left(a^{(1)}\otimes \ldots \otimes b^{(i)}_{k_2}\otimes \ldots \otimes a^{(m)}\right).
\end{equation*}
So, if there exists $b^{(i)}_k$ such that $b^{(i)}_k \in ME_{i}[U]$, then monomials of $\DMUp{U}$ can not cancel in the right-hand side of~\eqref{dep_inverse_image}.

Since $\widetilde{t}^{(i)} \in D_i[U]$, by definition, $\widetilde{t}^{(i)}$ is a linear combination of elements of $\Rel$. Hence, it follows from Small Cancellation Axiom that there exists at least one monomial $b^{(i)}_k$ in~\eqref{dependence_one_position} of $\SPM$-measure $\geqslant \tau + 1$. Therefore, using the above result, we obtain that monomials of $\DMUp{U}$ do not cancel in the right-hand side of~\eqref{dep_inverse_image}. Thus, $T \notin \Low\DSpace{U}$. Statement~\ref{EDp_bigger_set_of_generators1} is proved.

Let us prove statement~\ref{EDp_bigger_set_of_generators2}. Let $T$ be a generator of the space $\EDp\DSpace{U} $, that is, $T \in \GDp^{\prime}(\DSpace{U})$, $T \notin \Low\DSpace{U}$. Then it follows from statement~\ref{mu_properties_dep_inverse_image} of Lemma~\ref{mu_properties} that there exists an element $a^{(1)}\otimes \ldots \otimes t^{(i)}\otimes\ldots \otimes a^{(m)}$ such that $a^{(i^{\prime})} \in ME_{i^{\prime}}[U]$ for $i^{\prime} \neq i$, $t^{(i)} \in \Rel \cap D_i[U]$, and
\begin{equation*}
T = \mu[U]\left(a^{(1)}\otimes \ldots \otimes t^{(i)}\otimes\ldots \otimes a^{(m)}\right).
\end{equation*}
This means that $T \in \widetilde{\GDp}[U]$. So, $\EDp\DSpace{U} \subseteq \left\langle \widetilde{\GDp}[U]\right\rangle$.

Now let $T \in \widetilde{\GDp}[U]$. Then, by definition,
\begin{equation*}
T = \mu[U]\left(a^{(1)}\otimes \ldots \otimes \widetilde{t}^{(i)}\otimes \ldots \otimes a^{(m)}\right),
\end{equation*}
where $a^{(i^{\prime})} \in ME_{i^{\prime}}[U]$ for $i^{\prime} \neq i$, and $\widetilde{t}^{(i)} \in D_i[U]$, $\widetilde{t}^{(i)} \neq 0$. Since $\widetilde{t}^{(i)} \in D_i[U]$, by the definition of $D_i[U]$, we have
\begin{equation*}
\widetilde{t}^{(i)} = \sum_{k = 1}^s \gamma_k t_k,
\end{equation*}
where every $t_k \in \Rel \cap D_i[U]$. So,
\begin{equation*}
T = \sum_{k = 1}^s \gamma_k \mu[U]\left(a^{(1)}\otimes \ldots \otimes t_k\otimes \ldots \otimes a^{(m)}\right).
\end{equation*}
It follows from statement~\ref{mu_properties_up_dep_image} of Lemma~\ref{mu_properties} that $\mu[U]\left(a^{(1)}\otimes \ldots \otimes t_k\otimes \ldots \otimes a^{(m)}\right) \in \GDp^{\prime}(\DSpace{U})$. Statement~\ref{EDp_bigger_set_of_generators1} of Lemma~\ref{EDp_bigger_set_of_generators} implies  $\mu[U]\left(a^{(1)}\otimes \ldots \otimes t_k\otimes \ldots \otimes a^{(m)}\right) \notin \Low\DSpace{U}$. Therefore, $\mu[U]\left(a^{(1)}\otimes \ldots \otimes t_k\otimes \ldots \otimes a^{(m)}\right) \in \EDp\DSpace{U}$. Hence, $T \in \EDp\DSpace{U}$ and, so, $\left\langle \widetilde{\GDp}[U]\right\rangle \subseteq \EDp\DSpace{U}$. This completes the proof.
\end{proof}

\begin{proof}[\textbf{Proof of Lemma~\ref{fall_through_linear_dep_one}}]
Let $U$ be a monomial with $m$ virtual members of the chart, $U \in \Ft_{n}(\Galg) \setminus \Ft_{n - 1}(\Galg)$. Recall that we need to prove that
\begin{equation*}
\Dp\DSpace{U} \cap \Low\DSpace{U} \subseteq \Dp(\Ft_{n - 1}(\Galg)).
\end{equation*}
We noticed above that
\begin{align*}
&\Dp\DSpace{U} = \Dp(\Low\DSpace{U}) + \EDp\DSpace{U},\\
&\textit{where }\EDp = \left\langle \left\lbrace T \mid T \in \GDp^{\prime}(\DSpace{U}), T \notin \Low\DSpace{U} \right\rbrace\right\rangle.
\end{align*}
We have $\Dp(\Low\DSpace{U}) \subseteq \Low\DSpace{U}$. Hence, we see that
\begin{align*}
\Dp\DSpace{U} \cap \Low\DSpace{U} &= (\Dp(\Low\DSpace{U}) + \EDp\DSpace{U})\cap \Low\DSpace{U} =\\
&= \Dp(\Low\DSpace{U}) + \EDp\DSpace{U}\cap \Low\DSpace{U}.
\end{align*}
Obviously, $\Dp(\Low\DSpace{U}) \subseteq \Dp(\Ft_{n - 1}(\Galg))$. So, we need to prove that
\begin{equation*}
\EDp\DSpace{U} \cap \Low\DSpace{U} \subseteq \Dp(\Ft_{n - 1}(\Galg)).
\end{equation*}

Let us define a linear order on the monomials of $\DMUp{U}$. Since $ME_i[U]$ is a finite or countable set, we can introduce a linear order on the monomials of $ME_i[U]$ without infinite decreasing chains. For instance, the lexicographical order on $ME_i[U]$ satisfies this condition. So, in what follows we use the lexicographical order on $ME_i[U]$. We extend this ordering lexicographically on all the elements $a^{(1)}\otimes\ldots\otimes a^{(m)}$ such that $a^{(i)} \in ME_i[U]$ for all $i = 1, \ldots, m$. By the statement~\ref{mu_properties_up_correspondence} of Lemma~\ref{mu_properties}, $\mu[U]$ gives a bijective correspondence between all the elements $a^{(1)}\otimes\ldots\otimes a^{(m)}$ such that $a^{(i)} \in ME_i[U]$ for all $i = 1, \ldots, m$, and all the monomials of $\DMUp{U}$. Therefore, a linear order on the elements $a^{(1)}\otimes\ldots\otimes a^{(m)}$ such that $a^{(i)} \in ME_i[U]$ for all $i = 1, \ldots, m$, induces a linear order on $\DMUp{U}$. Clearly, this order on $\DMUp{U}$ does not have infinite decreasing chains. So, we can use this order for an induction.

Let $W \in \EDp\DSpace{U}\cap \Low\DSpace{U}$. By statement~\ref{EDp_bigger_set_of_generators2} of Lemma~\ref{EDp_bigger_set_of_generators}, we obtain  $\EDp\DSpace{U} = \left\langle \widetilde{\GDp}[U]\right\rangle$. Hence, we have
\begin{align*}
W = \sum_{q = 1}^l \gamma_q \widetilde{T}_q, \ &\textit{where } \widetilde{T}_q \textit{ are elements of } \widetilde{\GDp}[U],\\
&\textit{and } \sum\limits_{q = 1}^l \gamma_q \widetilde{T}_q \in \Low\DSpace{U}.
\end{align*}
It follows from statement~\ref{EDp_bigger_set_of_generators1} of Lemma~\ref{EDp_bigger_set_of_generators} that every $\widetilde{T}_q$ does not belong to $\Low\DSpace{U}$. Therefore, every $\widetilde{T}_q$ contains monomials of $\DMUp{U}$ in its additively reduced representation. We call the biggest monomial of $\DMUp{U}$ in additively reduced representation of $\widetilde{T}_q$ \emph{the highest monomial of $\widetilde{T}_q$}. We need to show that $\sum_{q = 1}^l \gamma_q \widetilde{T}_q \in \Dp(\Ft_{n - 1}(\Galg))$. We will prove this by induction on the biggest monomial among all highest monomials of $\widetilde{T}_q$ in the sum.

Let us make the step of induction. Let $X$ be the biggest monomial among all highest monomials of $\widetilde{T}_q$, $q = 1, \ldots, l$. Without loss of generality we can assume that $\widetilde{T}_1, \ldots, \widetilde{T}_{\widetilde{l}}$ have the highest monomial $X$, and $\widetilde{T}_{\widetilde{l} + 1}, \ldots, \widetilde{T}_{l}$ have smaller highest monomials (clearly, it is possible that $\widetilde{l} = l$). We will show that
\begin{equation}
\label{sum_with_smaller_leading_monomial}
\sum_{q = 1}^{\widetilde{l}} \gamma_q \widetilde{T}_q = \sum_{p} \gamma^{\prime}_p \widetilde{T}_p^{\prime} + Q,
\end{equation}
where every $\widetilde{T}_p^{\prime} \in \widetilde{\GDp}[U]$, the highest monomial of every $\widetilde{T}_p^{\prime}$ is smaller than $X$, and $Q \in \Dp(\Ft_{n - 1}(\Galg))$ (the sum $\sum_{p} \gamma_p^{\prime} \widetilde{T}_p^{\prime}$ may be empty).

Assume~\eqref{sum_with_smaller_leading_monomial} is proved, then we obtain
\begin{equation*}
W = \sum_{q = 1}^l \gamma_q\widetilde{T}_q = \sum_{q = 1}^{\widetilde{l}} \gamma_q\widetilde{T}_q + \sum_{q = \widetilde{l} + 1}^l \gamma_q\widetilde{T}_q = \sum_{p} \gamma_p^{\prime} \widetilde{T}_p^{\prime} + Q + \sum_{q = \widetilde{l} + 1}^l \gamma_q\widetilde{T}_q.
\end{equation*}
Since all monomials of $\DMUp{U}$ cancel out in $\sum_{q = 1}^l \gamma_q\widetilde{T}_q$ and $Q \in \Ft_{n - 1}(\Galg)$, we obtain that all monomials of $\DMUp{U}$ cancel out in the sum $\sum_{p} \gamma_p^{\prime} \widetilde{T}_p^{\prime} + \sum_{q = \widetilde{l} + 1}^l \gamma_q\widetilde{T}_q$. However, the highest monomial of every $\widetilde{T}_p^{\prime}$ and $\widetilde{T}_q$ for $q = \widetilde{l} + 1, \ldots, l$ is smaller than $X$. Therefore, by the induction hypothesis, we have
\begin{equation*}
\sum_{p} \gamma_p^{\prime} \widetilde{T}_p^{\prime} + \sum_{q = \widetilde{l} + 1}^l \gamma_q\widetilde{T}_q\in \Dp(\Ft_{n - 1}(\Galg)).
\end{equation*}
Thus, $W \in \Dp(\Ft_{n - 1}(\Galg))$.

Let us prove~\eqref{sum_with_smaller_leading_monomial}. Since $W = \sum_{q = 1}^l \gamma_q\widetilde{T}_q\in \Low\DSpace{U}$ and $\Low\DSpace{U}$ is generated by the monomials of $\DMLow{U}$, all monomials of $\DMUp{U}$ cancel out in the sum $\sum_{q = 1}^l
\gamma_q\widetilde{T}_q$. In particular, the monomial $X$ cancels out. Since $\widetilde{T}_{\widetilde{l} + 1}, \ldots, \widetilde{T}_{l}$ have the highest monomials smaller than $X$, the monomial $X$ is contained in none of $\widetilde{T}_{\widetilde{l} + 1}, \ldots, \widetilde{T}_{l}$. Therefore, $X$ cancels out in the sum~$\sum_{q = 1}^{\widetilde{l}} \gamma_q \widetilde{T}_q$.

Let $\delta_q$ be the coefficient of $X$ in $\widetilde{T}_q$, $q = 1, \ldots, \widetilde{l}$. Then we obtain
\begin{align*}
\sum_{q = 1}^{\widetilde{l}} \gamma_q \widetilde{T}_q =\ &\gamma_1\delta_1(\delta_1^{-1}\widetilde{T}_1 - \delta_2^{-1}\widetilde{T}_2) + (\gamma_1\delta_1 + \gamma_2\delta_2)(\delta_2^{-1}\widetilde{T}_2 - \delta_3^{-1}\widetilde{T}_3) + \ldots\\
&+ (\gamma_1\delta_1 + \ldots + \gamma_{\widetilde{l}}\delta_{\widetilde{l}})\delta_{\widetilde{l}}^{-1}\widetilde{T}_{\widetilde{l}}.
\end{align*}
Notice that $\gamma_1\delta_1 + \ldots + \gamma_{\widetilde{l}}\delta_{\widetilde{l}}$ is a coefficient of $X$ in $\sum_{q = 1}^{\widetilde{l}} \gamma_q \widetilde{T}_q$. Since $X$ cancels out in this sum, we see that $\gamma_1\delta_1 + \ldots + \gamma_{\widetilde{l}}\delta_{\widetilde{l}} = 0$. Therefore,
\begin{equation}
\label{pair_of_polynomials_representation}
\sum_{q = 1}^{\widetilde{l}} \gamma_q \widetilde{T}_q = \sum_{\substack{r, s = 1\\r < s}}^{\widetilde{l}}\nu_{r, s}(\delta_r^{-1}\widetilde{T}_r - \delta_s^{-1}\widetilde{T}_s).
\end{equation}
The monomial $X$ has the coefficient $1$ in every $\delta_q^{-1}\widetilde{T}_q$. Therefore, $X$ cancels out in every $\delta_r^{-1}\widetilde{T}_r - \delta_s^{-1}\widetilde{T}_s$.

Clearly, $\delta_r^{-1}\widetilde{T}_r, \delta_s^{-1}\widetilde{T}_s \in \widetilde{\GDp}[U]$. Hence, by definition, there exists $v^{(1)}\otimes \ldots \otimes t^{(i_r)}_r\otimes \ldots \otimes v^{(m)} \in A_1[U]\otimes\ldots\otimes D_{i_r}[U]\otimes\ldots\otimes A_m[U]$ such that $v^{(i)} \in ME_i[U]$ for $i \neq i_r$, $t^{(i_r)}_r \in D_{i_r}[U]$, and
\begin{equation*}
\mu[U]\left(v^{(1)}\otimes \ldots \otimes t^{(i_r)}_r\otimes \ldots \otimes v^{(m)}\right) = \delta_r^{-1}\widetilde{T}_r.
\end{equation*}
Similarly, there exists $w^{(1)}\otimes \ldots \otimes t^{(i_s)}_s\otimes \ldots \otimes w^{(m)} \in A_1[U]\otimes\ldots\otimes D_{i_s}[U]\otimes\ldots\otimes A_m[U]$ such that $w^{(i)} \in ME_i[U]$ for $i \neq i_s$, $t^{(i_s)}_r \in D_{i_s}[U]$, and
\begin{equation*}
\mu[U]\left(w^{(1)}\otimes \ldots \otimes t^{(i_s)}_s\otimes \ldots \otimes w^{(m)}\right) = \delta_s^{-1}\widetilde{T}_s.
\end{equation*}
It follows from the statement~\ref{mu_properties_up_correspondence} of Lemma~\ref{mu_properties} that there exists an element $x^{(1)}\otimes\ldots\otimes x^{(m)} \in A_1[U]\otimes\ldots\otimes A_m[U]$ such that every $x^{(i)} \in ME_i[U]$ and
\begin{equation*}
\mu[U]\left(x^{(1)}\otimes\ldots\otimes x^{(m)}\right) = X.
\end{equation*}
Since $X$ is a monomial of $\delta_r^{-1}\widetilde{T}_r$ and $\delta_s^{-1}\widetilde{T}_s$ and $\mu[U]$ is a bijective correspondence between the elements $a^{(1)}\otimes\ldots\otimes a^{(m)}$ such that $a^{(i)} \in ME_i[U]$ and the monomials of $\DMUp{U}$, we see that $v^{(i)} = w^{(i)} = x^{(i)}$ for $i \neq i_r, i_s$, and $x^{(i_r)}$ is a monomial of $t^{(i_r)}_r$, $x^{(i_s)}$ is a monomial of $t^{(i_s)}_s$. So, we have
\begin{align*}
&\delta_r^{-1}\widetilde{T}_r = \mu[U]\left(x^{(1)}\otimes \ldots \otimes t^{(i_r)}_r\otimes \ldots \otimes x^{(m)}\right),\\
&\delta_s^{-1}\widetilde{T}_s = \mu[U]\left(x^{(1)}\otimes \ldots \otimes t^{(i_s)}_s\otimes \ldots \otimes x^{(m)}\right).
\end{align*}
Since $X$ is the highest monomial of $\delta_r^{-1}\widetilde{T}_r$ and $\delta_s^{-1}\widetilde{T}_s$, we see that $x^{(i_r)}$ is the lexicographically biggest monomial of $ME_{i_r}[U]$ in $t^{(i_r)}_r$, and $x^{(i_s)}$ is the lexicographically biggest monomial of $ME_{i_s}[U]$ in $t^{(i_s)}_s$.

First consider the case $i_r = i_s = i_0$. Then we have
\begin{align*}
&\delta_r^{-1}\widetilde{T}_r - \delta_r^{-1}\widetilde{T}_s =\\
&= \mu[U]\left(x^{(1)}\otimes\ldots\otimes t^{(i_0)}_r\otimes\ldots\otimes x^{(m)}\right) - \mu[U]\left(x^{(1)}\otimes\ldots\otimes t^{(i_0)}_s\otimes\ldots\otimes x^{(m)}\right) =\\
&= \mu[U]\left(x^{(1)}\otimes\ldots\otimes \left(t^{(i_0)}_r - t^{(i_0)}_s\right) \otimes\ldots\otimes x^{(m)}\right).
\end{align*}
If $t^{(i_0)}_r = t^{(i_0)}_s$, then $\delta_r^{-1}\widetilde{T}_r - \delta_r^{-1}\widetilde{T}_s = 0$, hence, we can delete it from the sum~\eqref{pair_of_polynomials_representation}. Assume $t^{(i_0)}_r \neq t^{(i_0)}_s$. Then, by definition, $\delta_r^{-1}\widetilde{T}_r - \delta_r^{-1}\widetilde{T}_s \in \widetilde{\GDp}[U]$. Recall that  $X$ cancels out in $\delta_r^{-1}\widetilde{T}_r - \delta_s^{-1}\widetilde{T}_s$. Therefore,
\begin{align}
\begin{split}
\label{the_same_position_result}
&\delta_r^{-1}\widetilde{T}_r - \delta_r^{-1}\widetilde{T}_s \in \widetilde{\GDp}[U] \textit{ and }\\
&\delta_r^{-1}\widetilde{T}_r - \delta_s^{-1}\widetilde{T}_s \textit{ has the smaller highest monomial than } X.
\end{split}
\end{align}

Now consider the case $i_r \neq i_s$. To be definite, assume $i_r < i_s$. Let
\begin{equation*}
t^{(i_r)}_r = \sum_{j = 1}^{n_r}\beta_j b_j^{(i_r)},\ t^{(i_s)}_s = \sum_{k = 1}^{n_s}\eta_k c_k^{(i_s)},
\end{equation*}
where $b_j^{(i_r)} \in A_{i_r}[U]$, $j = 1, \ldots, n_r$, and $c_k^{(i_s)} \in A_{i_s}[U]$, $k = 1, \ldots, n_s$, are monomials, and the right-hand sides are additively reduced. Without loss of generality we can assume that $b_1^{(i_r)}$ is the lexicographically biggest monomial of $ME_{i_r}[U]$ in $t^{(i_r)}_r$, and $c_1^{(i_s)}$ is the lexicographically biggest monomial of $ME_{i_s}[U]$ in $t^{(i_s)}_s$. According to the above, this means that $b_1^{(i_r)} = x^{(i_r)}$ and $c_1^{(i_s)} = x^{(i_s)}$. Since the coefficient of $X$ in $\delta_r^{-1}\widetilde{T}_r$ and in $\delta_s^{-1}\widetilde{T}_s$ is equal to $1$, we see that $\beta_1 = 1$ and $\eta_1 = 1$.

We consider elements $x^{(1)}\otimes\ldots \otimes t^{(i_r)}_r\otimes\ldots\otimes x^{(m)}$ and $x^{(1)}\otimes \ldots \otimes t^{(i_s)}_s\otimes \ldots \otimes x^{(m)}$. Since we have changes only in positions $i_r$ and $i_s$, let us focus only on these positions in the further auxiliary calculations. Since the coefficients $\beta_1 = 1$ and $\eta_1 = 1$, one can easily see that
\begin{align*}
&t^{(i_r)}_r\otimes x^{(i_s)} - x^{(i_r)}\otimes t^{(i_s)}_s = t^{(i_r)}_r\otimes c_1^{(i_s)} - b_1^{(i_r)}\otimes t^{(i_s)}_s =\\
&= \sum_{j = 1}^{n_r}\beta_j b_j^{(i_r)}\otimes c_1^{(i_s)} - \sum_{k = 1}^{n_s}\eta_k b^{(i_r)}_1\otimes c^{(i_s)}_k =\\
&= b_1^{(i_r)} \otimes c_1^{(i_r)} + \sum_{j = 2}^{n_r}\beta_j b_j^{(i_r)}\otimes c_1^{(i_s)} -  b_1^{(i_r)} \otimes c_1^{(i_r)} - \sum_{k = 2}^{n_s}\eta_k b^{(i_r)}_1\otimes c^{(i_s)}_k =\\
&= \sum_{j = 2}^{n_r}\beta_j b_j^{(i_r)}\otimes c_1^{(i_s)} - \sum_{k = 2}^{n_s}\eta_k b^{(i_r)}_1\otimes c^{(i_s)}_k.
\end{align*}
On the other hand, we have
\begin{align*}
\sum_{j = 2}^{n_r}\beta_j b_j^{(i_r)}&\otimes t^{(i_s)}_s - \sum_{k = 2}^{n_s} \eta_k t^{(i_r)}_r\otimes c^{(i_s)}_k =\\
&= \sum_{j = 2}^{n_r}\beta_j b_j^{(i_r)}\otimes \left(\sum_{k = 1}^{n_s}\eta_k c_k^{(i_s)}\right) - \sum_{k = 2}^{n_s} \eta_k \left(\sum_{j = 1}^{n_r}\beta_j b_j^{(i_r)}\right)\otimes c^{(i_s)} =\\
&= \sum_{j = 2}^{n_r}\sum_{k = 1}^{n_s}\beta_j\eta_k b_j^{(i_r)}\otimes c^{(i_s)}_k - \sum_{k = 2}^{n_s}\sum_{j = 1}^{n_r} \eta_k\beta_j  b^{(i_r)}_j\otimes c^{(i_s)}_k =\\
=\sum_{j = 2}^{n_r}\beta_j\eta_1& b_j^{(i_r)}\otimes c^{(i_s)}_1 - \sum_{k = 2}^{n_s} \eta_k\beta_1 b^{(i_r)}_1\otimes c^{(i_s)}_k = \sum_{j = 2}^{n_r}\beta_j b_j^{(i_r)}\otimes c^{(i_s)}_1 - \sum_{k = 2}^{n_s} \eta_k b^{(i_r)}_1\otimes c^{(i_s)}_k.
\end{align*}
Combining the above results, we see that
\begin{align*}
&t^{(i_r)}_r\otimes x^{(i_s)} - x^{(i_r)}\otimes t^{(i_s)}_s = \sum_{j = 2}^{n_r}\beta_j b_j^{(i_r)}\otimes t^{(i_s)}_s - \sum_{k = 2}^{n_s} \eta_k t^{(i_r)}_r\otimes c^{(i_s)}_k.
\end{align*}
So, we have
\begin{align*}
\delta_r^{-1}\widetilde{T}_r - \delta_s^{-1}\widetilde{T}_s =\ \ \ \ \ \ \ \ \ \ \ \ \ \ \ \ \ \ \ \ \ \ \ \ &\\
\mu[U]\left(x^{(1)}\otimes\ldots \otimes t^{(i_r)}_r\otimes \ldots \otimes x^{(m)}\right) &- \mu[U]\left(x^{(1)}\otimes\ldots \otimes t^{(i_s)}_s\otimes \ldots \otimes x^{(m)}\right) =\\
\mu[U]\left(x^{(1)}\otimes\ldots \otimes t^{(i_r)}_r\otimes \ldots \otimes x^{(i_s)}\otimes \right.&\ldots \otimes x^{(m)} -\\
- x^{(1)}\otimes&\left.\ldots \otimes x^{(i_r)}\otimes \ldots\otimes t^{(i_s)}_s\otimes \ldots \otimes x^{(m)}\right) =\\
= \sum_{j = 2}^{n_r}\beta_j \mu[U]\left(x^{(1)}\otimes\ldots \otimes b_j^{(i_r)}\otimes\ldots\right.&\left. \otimes t^{(i_s)}_s\otimes \ldots \otimes x^{(m)}\right) -\\
- \sum_{k = 2}^{n_s} \eta_k \mu[U]\left(x^{(1)}\right.&\left.\otimes \ldots \otimes t^{(i_r)}_r\otimes \ldots \otimes c^{(i_s)}_k\otimes \ldots \otimes x^{(m)}\right).
\end{align*}

Consider the element $\mu[U]\left(x^{(1)}\otimes \ldots \otimes b_j^{(i_r)}\otimes\ldots\otimes t^{(i_s)}_s\otimes \ldots \otimes x^{(m)}\right)$, $j \neq 1$. We distinguish two possibilities: $b_j^{(i_r)} \in ML_{i_r}[U]$ and $b_j^{(i_r)} \in ME_{i_r}[U]$. First assume that $b_j^{(i_r)} \in ML_{i_r}[U]$. Then, by the statement~\ref{mu_properties_dep_low_image} of Lemma~\ref{mu_properties}, we obtain
\begin{equation*}
\mu[U]\left(x^{(1)}\otimes \ldots \otimes b_j^{(i_r)}\otimes t^{(i_s)}_s\otimes \ldots \otimes x^{(m)}\right) \in \Dp(\Ft_{n - 1}(\Galg)).
\end{equation*}
Now assume that $b_j^{(i_r)} \in ME_{i_r}[U]$. Then, by definition, we have
\begin{equation*}
\mu[U]\left(x^{(1)}\otimes \ldots \otimes b_j^{(i_r)}\otimes\ldots\otimes t^{(i_s)}_s\otimes \ldots \otimes x^{(m)}\right) \in \widetilde{T}[U].
\end{equation*}
We assumed above that $c_1^{(i_s)}$ is the lexicographically biggest monomial of $ME_{i_s}[U]$ in $t^{(i_s)}_s$. Hence, by the definition of our order on $\DMUp{U}$, we obtain that $\mu[U]\left(x^{(1)}\otimes \ldots \otimes b_j^{(i_r)}\otimes\ldots\otimes c_1^{(i_s)}\otimes \ldots \otimes x^{(m)}\right)$ is the highest monomial of $\mu[U]\left(x^{(1)}\otimes \ldots \otimes b_j^{(i_r)}\otimes\ldots\otimes t^{(i_s)}_s\otimes \ldots \otimes x^{(m)}\right)$. We also assumed that $b_1^{(i_r)}$ is the lexicographically biggest monomial of $ME_{i_r}[U]$ in $t^{(i_r)}_r$. Therefore, since $b_j^{(i_r)} \in ME_{i_r}[U]$, we see that $b_j^{(i_r)}$, $j \neq 1$, is lexicographically smaller than $b_1^{(i_r)}$. Recall that $x^{(i_r)} = b_1^{(i_r)}$ and $x^{(i_s)} = c_1^{(i_s)}$. That is, we have
\begin{equation*}
X = \mu[U]\left(x^{(1)}\otimes \ldots \otimes b_1^{(i_r)}\otimes\ldots\otimes c^{(i_s)}_1\otimes \ldots \otimes x^{(m)}\right).
\end{equation*}
Hence, the monomial $\mu[U]\left(x^{(1)}\otimes \ldots \otimes b_j^{(i_r)}\otimes\ldots\otimes c_1^{(i_s)}\otimes \ldots \otimes x^{(m)}\right)$, $j \neq 1$, is smaller than $X$ with respect to our order on $\DMUp{U}$. Thus,
\begin{multline*}
\mu[U]\left(x^{(1)}\otimes \ldots \otimes b_j^{(i_r)}\otimes\ldots\otimes t^{(i_s)}_s\otimes \ldots \otimes x^{(m)}\right) \in \widetilde{\GDp}[U]\\
\textit{ has smaller highest monomial than } X.
\end{multline*}

The element $\mu[U]\left(x^{(1)}\otimes \ldots \otimes t^{(i_r)}_r\otimes\ldots\otimes c^{(i_s)}_k\otimes \ldots \otimes x^{(m)}\right)$, $k \neq 1$, is studied similarly. Namely, if $c_k^{(i_s)} \in ML_{i_s}[U]$, then
\begin{equation*}
\mu[U]\left(x^{(1)}\otimes \ldots \otimes t^{(i_r)}_r\otimes c^{(i_s)}_k\otimes \ldots \otimes x^{(m)}\right) \in \Dp(\Ft_{n - 1}(\Galg)).
\end{equation*}
If $c_k^{(i_s)} \in ME_{i_s}[U]$, then
\begin{multline*}
\mu[U]\left(x^{(1)}\otimes \ldots \otimes t^{(i_r)}_r\otimes\ldots\otimes c^{(i_s)}_k\otimes \ldots \otimes x^{(m)}\right) \in \widetilde{\GDp}[U]\\
\textit{ has smaller highest monomial than } X.
\end{multline*}

So, for $i_r \neq i_s$ we eventually obtain
\begin{align}
\begin{split}
\label{different_positions_result}
\delta_r^{-1}\widetilde{T}_r - \delta_s^{-1}\widetilde{T}_s =\ \ \ \ \ \ \ \ \ \ \ \ \ \ \ \ \ &\\
\sum_{j = 2}^{n_r}\beta_j \mu[U]\left(x^{(1)}\otimes \ldots \otimes b_j^{(i_r)}\otimes\right.&\left.\ldots \otimes t^{(i_s)}_s\otimes \ldots \otimes x^{(m)}\right) -\\
- \sum_{k = 2}^{n_s} \eta_k \mu[U]\left(x^{(1)}\otimes \right.&\left.\ldots \otimes t^{(i_r)}_r\otimes \ldots \otimes c^{(i_s)}_k\otimes \ldots \otimes x^{(m)}\right) =\\
&\ \ \ \ \ \ \ \ \ \ \ \ \ \ \ \ \ \ \ \ \ = \sum_p\gamma_p^{r, s}\widetilde{T}_p^{r, s} + Q^{r, s},
\end{split}
\end{align}
where $\widetilde{T}_p^{r, s} \in \widetilde{\GDp}[U]$, every $\widetilde{T}_p^{r, s}$ has the highest monomial smaller than $X$, and $Q^{r, s} \in \Dp(\Ft_{n - 1}(\Galg))$.

Combining~\eqref{different_positions_result} with~\eqref{pair_of_polynomials_representation} and~\eqref{the_same_position_result}, we obtain the equality~\eqref{sum_with_smaller_leading_monomial}. So far, we are done with the step of induction.

Let us prove the basis of induction. Again let $X$ be the biggest monomial among the highest monomials of all $\widetilde{T}_q$, $q = 1, \ldots, l$. Assume that $X$ is the smallest monomial of $\DMUp{U}$ with respect to our order on $\DMUp{U}$. Then $X$ is the only monomial of $\DMUp{U}$ in every $\widetilde{T}_q$, $q = 1, \ldots, l$. We have to prove that
\begin{equation*}
\sum_{q = 1}^l \gamma_q\widetilde{T}_q \in \Dp(\Ft_{n - 1}(\Galg)).
\end{equation*}
Let us argue as in the step of induction and use the same notations. As above, we need to consider $\delta_r^{-1}\widetilde{T}_r - \delta_s^{-1}\widetilde{T}_s$. The equality~\eqref{pair_of_polynomials_representation} implies that it is enough to prove that $\delta_r^{-1}\widetilde{T}_r - \delta_s^{-1}\widetilde{T}_s \in \Dp(\Ft_{n - 1}(\Galg))$.

Notice that since $X$ is the only monomial of $\DMUp{U}$ in every $\widetilde{T}_q$, we see that
\begin{equation}
\label{ind_basis_low}
\delta_r^{-1}\widetilde{T}_r - \delta_s^{-1}\widetilde{T}_s \in \Low\DSpace{U} \subseteq \Ft_{n - 1}(\Galg).
\end{equation}
As above, let
\begin{align*}
&\delta_r^{-1}\widetilde{T}_r = \mu[U]\left(x^{(1)}\otimes \ldots \otimes t^{(i_r)}_r\otimes \ldots \otimes x^{(m)}\right),\\
&\delta_s^{-1}\widetilde{T}_s = \mu[U]\left(x^{(1)}\otimes \ldots \otimes t^{(i_s)}_s\otimes \ldots \otimes x^{(m)}\right),\\
&\textit{where } t^{(i_r)}_r = \sum_{j = 1}^{n_r}\beta_j b_j^{(i_r)},\ t^{(i_s)}_s = \sum_{k = 1}^{n_s}\eta_k c_k^{(i_s)},\\
&X = \mu[U]\left(x^{(1)}\otimes \ldots \otimes b^{(i_r)}_1\otimes\ldots\otimes c^{(i_s)}_1\otimes \ldots \otimes x^{(m)}\right).
\end{align*}
Since $X$ is the only monomial of $\DMUp{U}$ in every $\widetilde{T}_q$, the statement~\ref{mu_properties_up_correspondence} of Lemma~\ref{mu_properties} implies  $b_j^{(i_r)} \in ML_{i_r}[U]$ if $j \neq 1$, and $c_k^{(i_s)} \in ML_{i_s}[U]$ if $k \neq 1$.

First assume $i_r = i_s = i_0$. Then, as above,
\begin{equation*}
\delta_r^{-1}\widetilde{T}_r - \delta_s^{-1}\widetilde{T}_s = \mu[U]\left(x^{(1)}\otimes \ldots \otimes \left(t^{(i_0)}_r - t^{(i_0)}_s\right)\otimes \ldots \otimes x^{(m)}\right).
\end{equation*}
Assume $t^{(i_0)}_r - t^{(i_0)}_s \neq 0$. Then it follows from statement~\ref{EDp_bigger_set_of_generators1} of Lemma~\ref{EDp_bigger_set_of_generators} that $\delta_r^{-1}\widetilde{T}_r - \delta_s^{-1}\widetilde{T}_s \notin \Low\DSpace{U}$. A contradiction with~\eqref{ind_basis_low}. Therefore, if $i_r = i_s = i_0$, then $t^{(i_0)}_r - t^{(i_0)}_s = 0$ and $\delta_r^{-1}\widetilde{T}_r - \delta_s^{-1}\widetilde{T}_s = 0$.

Now assume $i_r \neq i_s$. To be definite, assume $i_r < i_s$.  As above, we obtain that $\delta_r^{-1}\widetilde{T}_r - \delta_s^{-1}\widetilde{T}_s$ is a linear combination of the elements
\begin{align*}
&\mu[U]\left(x^{(1)}\otimes \ldots \otimes t^{(i_r)}_r\otimes\ldots\otimes c^{(i_s)}_k\otimes \ldots \otimes x^{(m)}\right), k\neq 1,\\
&\mu[U]\left(x^{(1)}\otimes \ldots \otimes b_j^{(i_r)}\otimes\ldots\otimes t^{(i_s)}_s\otimes \ldots \otimes x^{(m)}\right), j \neq 1.
\end{align*}
Since $b_j^{(i_r)} \in ML_{i_r}[U]$ if $j \neq 1$ and $c_k^{(i_s)} \in ML_{i_s}[U]$ if $k \neq 1$, the statement~\ref{mu_properties_dep_low_image} of Lemma~\ref{mu_properties} implies
\begin{align*}
&\mu[U](a^{(1)}\otimes \ldots \otimes t^{(i_r)}_r\otimes\ldots\otimes c^{(i_s)}_k\otimes \ldots \otimes a^{(m)}) \in \Dp(\Ft_{n - 1}(\Galg)) \textit{ if }k \neq 1\\
&\mu[U](a^{(1)}\otimes \ldots \otimes b_j^{(i_r)}\otimes\ldots\otimes t^{(i_s)}_s\otimes \ldots \otimes a^{(m)}) \in \Dp(\Ft_{n - 1}(\Galg)) \textit{ if } j \neq 1.
\end{align*}
Thus, $\delta_r^{-1}\widetilde{T}_r - \delta_s^{-1}\widetilde{T}_s \in \Dp(\Ft_{n - 1}(\Galg))$. So, we are done with the basis of induction. This completes the proof of Lemma~\ref{fall_through_linear_dep_one}.
\end{proof}

\begin{proposition}
\label{correspondence_to_tensor_product}
Let $U$ be a monomial with $m$ virtual members of the chart. Suppose $A_i[U]$, $L_i[U] \subseteq A_i[U]$, $D_i[U] \subseteq A_i[U]$, $i = 1, \ldots, m$, are subspaces of $\Galg$ defined above by~\eqref{tens_prod_components},~\eqref{tens_prod_components_low}, and~\eqref{tens_prod_dep}. Then we have
\begin{multline*}
\DSpace{U} / (\Dp\DSpace{U} + \Low\DSpace{U}) \cong\\
\cong A_1[U]/(D_1[U] + L_1[U]) \otimes \ldots \otimes A_m[U]/(D_m[U] + L_m[U]).
\end{multline*}
\end{proposition}
\begin{proof}
Assume $U \in \Ft_n(\Galg) \setminus \Ft_{n - 1}(\Galg)$. Let
\begin{equation*}
\mu[U]: A_1[U]\otimes\ldots\otimes A_m[U] \to \DSpace{U} + \Ft_{n - 1}(\Galg)
\end{equation*}
be a linear mapping defined by Definition~\ref{mu_def}. This statement is, in fact, a corollary of Lemma~\ref{mu_properties}. It follows from statements~\ref{mu_properties_low_correspondence} and~\ref{mu_properties_up_dep_image} of Lemma~\ref{mu_properties} that
\begin{multline*}
\mu[U](A_1[U]\otimes \ldots \otimes D_i[U]\otimes \ldots\otimes A_m[U] + \Low(A_1[U]\otimes\ldots\otimes A_m[U])) \subseteq\\
\subseteq \Dp\DSpace{U} + \Ft_{n - 1}(\Galg) \textit{ for all } i = 1, \ldots, m.
\end{multline*}
Hence,
\begin{multline}
\label{mu_dependencies_plus_low}
\mu[U](\Dp(A_1[U]\otimes\ldots\otimes A_m[U]) + \Low(A_1[U]\otimes\ldots\otimes A_m[U])) \subseteq\\
\subseteq \Dp\DSpace{U} + \Ft_{n - 1}(\Galg).
\end{multline}
We define the mapping
\begin{multline*}
\overline{\mu}[U]: A_1[U]\otimes\ldots\otimes A_m[U] / (\Dp(A_1[U]\otimes\ldots\otimes A_m[U]) + \Low(A_1[U]\otimes\ldots\otimes A_m[U])) \to\\
\to (\DSpace{U} + \Ft_{n - 1}(\Galg)) / (\Dp\DSpace{U} + \Ft_{n - 1}(\Galg))
\end{multline*}
by the formula
\begin{multline*}
\overline{\mu}[U](W) + \Dp(A_1[U]\otimes\ldots\otimes A_m[U]) + \Low(A_1[U]\otimes\ldots\otimes A_m[U])) =\\
= \mu[U](W) + \Dp\DSpace{U} + \Ft_{n - 1}(\Galg),
\end{multline*}
where $W$ is an arbitrary element of $A_1[U]\otimes\ldots\otimes A_m[U]$. It follows from~\eqref{mu_dependencies_plus_low} that the mapping $\overline{\mu}[U]$ is a well-defined homomorphism of vector spaces. Let us show that $\overline{\mu}[U]$ is an isomorphism of vector spaces. That is, we have to show that $\overline{\mu}[U]$ is a bijective mapping.

Clearly, $(\DSpace{U} + \Ft_{n - 1}(\Galg)) / (\Dp\DSpace{U} + \Ft_{n - 1}(\Galg))$ is linearly generated by all the elements of the form $Z + \Dp\DSpace{U} + \Ft_{n - 1}(\Galg)$ such that $Z \in \DMUp{U}$. By the statement~\ref{mu_properties_up_correspondence} of Lemma~\ref{mu_properties}, for every $Z \in \DMUp{U}$ there exists an element $W \in A_1[U]\otimes\ldots\otimes A_m[U]$ such that $\mu[U](W) = Z$. Hence, $\overline{\mu}[U]$ is a surjective mapping.

Assume
\begin{equation*}
\overline{T}^{\prime} \in A_1[U]\otimes\ldots\otimes A_m[U] / (\Dp(A_1[U]\otimes\ldots\otimes A_m[U]) + \Low(A_1[U]\otimes\ldots\otimes A_m[U]))
\end{equation*}
and $\overline{\mu}[U](\overline{T}^{\prime}) = 0$. Let us show that $\overline{T}^{\prime} = 0$. We have
\begin{equation*}
\overline{T}^{\prime} = T^{\prime} + \Dp(A_1[U]\otimes\ldots\otimes A_m[U]) + \Low(A_1[U]\otimes\ldots\otimes A_m[U]),
\end{equation*}
where $T^{\prime} \in A_1[U]\otimes\ldots\otimes A_m[U]$. By the definition of $\overline{\mu}[U]$, we have
\begin{align*}
&0 = \overline{\mu}[U](\overline{T}^{\prime}) =\\
&= \overline{\mu}[U](T^{\prime} + \Dp(A_1[U]\otimes\ldots\otimes A_m[U]) + \Low(A_1[U]\otimes\ldots\otimes A_m[U])) =\\
&= \mu[U](T^{\prime}) + \Dp\DSpace{U} + \Ft_{n - 1}(\Galg).
\end{align*}
Therefore, $\mu[U](T^{\prime}) \in \Dp\DSpace{U} + \Ft_{n - 1}(\Galg)$. This means that
\begin{equation*}
\mu[U](T^{\prime}) = \sum_{s = 1}^l \gamma_sT_s + X,
\end{equation*}
where $T_s \in \GDp^{\prime}(\DSpace{U})$, $T_s \notin \Ft_{n - 1}(\Galg)$, $s = 1, \ldots, l$, and $X \in \Ft_{n - 1}(\Galg)$. Since $T_s \notin \Ft_{n - 1}(\Galg)$, Lemma~\ref{intersection_with_low_equality} implies  $T_s \notin \Low\DSpace{U}$. Hence, it follows from the statement~\ref{mu_properties_dep_inverse_image} of Lemma~\ref{mu_properties} that there exists $T^{\prime}_s \in \Dp(A_1[U]\otimes\ldots\otimes A_m[U])$ such that $\mu[U](T^{\prime}_s) = T_s$. Therefore,
\begin{equation*}
\mu[U](T^{\prime} - \sum_{s = 1}^l \gamma_sT_s^{\prime}) = X \in \Ft_{n - 1}(\Galg).
\end{equation*}

We have $T^{\prime} - \sum_{s = 1}^l \gamma_sT_s^{\prime} \in A_1[U]\otimes\ldots\otimes A_m[U]$. Let
\begin{equation}
\label{lower_sum_for_mu}
T^{\prime} - \sum_{s = 1}^l \gamma_sT_s^{\prime} = \sum_{j_1, \ldots, j_m}\alpha_{j_1, \ldots, j_m}a^{(1)}_{j_1} \otimes \ldots \otimes a^{(m)}_{j_m},
\end{equation}
where $a^{(i)}_{j_i}$ is a monomial of $A_i[U]$, $i = 1, \ldots, m$, and the right-hand sum is additively reduced. Assume $T^{\prime} - \sum_{s = 1}^l \gamma_sT_s^{\prime} \notin \Low(A_1[U]\otimes\ldots\otimes A_m[U])$. This means that not all $a^{(1)}_{j_1} \otimes \ldots \otimes a^{(m)}_{j_m}$ in the right-hand side of~\eqref{lower_sum_for_mu} belong to $\Low(A_1[U]\otimes\ldots\otimes A_m[U])$. By the statement~\ref{mu_properties_up_correspondence} of Lemma~\ref{mu_properties}, $\mu[U]$ gives a bijective correspondence between all the elements $b^{(1)}\otimes\ldots\otimes b^{(m)} \notin \Low(A_1[U]\otimes\ldots\otimes A_m[U])$ such that $b^{(i)} \in A_i[U]$, $i = 1, \ldots, m$, are monomials, and all the monomials of $\DMUp{U}$. By the statement~\ref{mu_properties_low_correspondence} of Lemma~\ref{mu_properties}, we have $\mu[U](\Low(A_1[U]\otimes\ldots\otimes A_m[U])) \subseteq \Ft_{n - 1}(\Galg)$. Notice that every monomial of $\DMUp{U}$ belongs to $\Ft_n(\Galg) \setminus  \Ft_{n - 1}(\Galg)$. Combining these statements, we see that if there exist elements $a^{(1)}_{j_1} \otimes \ldots \otimes a^{(m)}_{j_m} \notin \Low(A_1[U]\otimes\ldots\otimes A_m[U])$ in the right-hand side of~\eqref{lower_sum_for_mu}, then their images under the mapping $\mu[U]$ belong to $\DMUp{U}$ and can not additively cancel out in
\begin{equation*}
\mu[U]\left(\sum_{j_1, \ldots, j_m}\alpha_{j_1, \ldots, j_m}a^{(1)}_{j_1} \otimes \ldots \otimes a^{(m)}_{j_m}\right) = \sum_{j_1, \ldots, j_m}\alpha_{j_1, \ldots, j_m}\mu[U](a^{(1)}_{j_1} \otimes \ldots \otimes a^{(m)}_{j_m}).
\end{equation*}
So, since every monomial of $\DMUp{U}$ belongs to $\Ft_n(\Galg) \setminus  \Ft_{n - 1}(\Galg)$, we obtain
\begin{equation*}
\mu[U]\left(T^{\prime} - \sum_{s = 1}^l \gamma_sT_s^{\prime}\right) = \sum_{j_1, \ldots, j_m}\alpha_{j_1, \ldots, j_m}\mu[U](a^{(1)}_{j_1} \otimes \ldots \otimes a^{(m)}_{j_m}) \notin \Ft_{n - 1}(\Galg),
\end{equation*}
a contradiction. Therefore, $T^{\prime} - \sum_{s = 1}^l \gamma_sT_s^{\prime} \in \Low(A_1[U]\otimes\ldots\otimes A_m[U])$.

So, we proved that
\begin{equation*}
T^{\prime} = \sum_{s = 1}^l \gamma_sT_s^{\prime} + X^{\prime},\textit{ where } X^{\prime} \in \Low(A_1[U]\otimes\ldots\otimes A_m[U]).
\end{equation*}
Since $T_s^{\prime} \in \Dp(A_1[U]\otimes\ldots\otimes A_m[U])$, $s = 1, \ldots, l$, this yields
\begin{equation*}
T^{\prime} \in \Dp(A_1[U]\otimes\ldots\otimes A_m[U]) + \Low(A_1[U]\otimes\ldots\otimes A_m[U]).
\end{equation*}
Hence, finally we see that
\begin{equation*}
\overline{T}^{\prime} = T^{\prime} + \Dp(A_1[U]\otimes\ldots\otimes A_m[U]) + \Low(A_1[U]\otimes\ldots\otimes A_m[U]) = 0.
\end{equation*}
Thus, $\overline{\mu}[U]$ is an injective mapping.

Let us define the mapping
\begin{multline*}
\psi[U] : A_1[U]/(D_1[U] + L_1[U]) \otimes \ldots \otimes A_m[U]/(D_m[U] + L_m[U]) \to\\
\to A_1[U]\otimes\ldots\otimes A_m[U] / (\Dp(A_1[U]\otimes\ldots\otimes A_m[U]) + \Low(A_1[U]\otimes\ldots\otimes A_m[U])).
\end{multline*}
Let $\psi[U]$ take each element
\begin{equation*}
(a^{(1)} + D_1[U] + L_1[U])\otimes\ldots\otimes (a^{(m)} + D_m[U] + L_m[U])
\end{equation*}
such that $a^{(i)} \in ME_i[U]$, $i = 1, \ldots, m$, to
\begin{equation*}
a^{(1)}\otimes\ldots\otimes a^{(m)} + \Dp(A_1[U]\otimes\ldots\otimes A_m[U]) + \Low(A_1[U]\otimes\ldots\otimes A_m[U]),
\end{equation*}
and let it be extended linearly on the space $A_1[U]/(D_1[U] + L_1[U]) \otimes \ldots \otimes A_m[U]/(D_m[U] + L_m[U])$. One can easily show that $\psi[U]$ is well-defined and $\psi[U]$ is an isomorphism of vector spaces. Thus,
\begin{multline*}
A_1[U]/(D_1[U] + L_1[U]) \otimes \ldots \otimes A_m[U]/(D_m[U] + L_m[U]) \cong\\
\cong A_1[U]\otimes\ldots\otimes A_m[U] / (\Dp(A_1[U]\otimes\ldots\otimes A_m[U]) + \Low(A_1[U]\otimes\ldots\otimes A_m[U])).
\end{multline*}

Let us show that
\begin{equation*}
(\DSpace{U} + \Ft_{n - 1}(\Galg)) / (\Dp\DSpace{U} + \Ft_{n - 1}(\Galg)) \cong \DSpace{U} / (\Dp\DSpace{U} + \Low\DSpace{U}).
\end{equation*}
Indeed, by the isomorphism theorem, we have
\begin{multline*}
(\DSpace{U} + \Ft_{n - 1}(\Galg)) / (\Dp\DSpace{U} + \Ft_{n - 1}(\Galg)) =\\
= (\DSpace{U} + \Dp\DSpace{U} + \Ft_{n - 1}(\Galg)) / (\Dp\DSpace{U} + \Ft_{n - 1}(\Galg)) \cong\\
\cong \DSpace{U} / ((\Dp\DSpace{U} + \Ft_{n - 1}(\Galg))\cap \DSpace{U}).
\end{multline*}
Since $\Dp\DSpace{U} \subseteq \DSpace{U}$, we obtain
\begin{equation*}
(\Dp\DSpace{U} + \Ft_{n - 1}(\Galg))\cap \DSpace{U} = \Dp\DSpace{U} + \Ft_{n - 1}(\Galg) \cap \DSpace{U}.
\end{equation*}
By Lemma~\ref{intersection_with_low_equality}, we have $\Ft_{n - 1}(\Galg) \cap \DSpace{U} = \Low\DSpace{U}$. So, finally we see that
\begin{equation*}
(\Dp\DSpace{U} + \Ft_{n - 1}(\Galg))\cap \DSpace{U} = \Dp\DSpace{U} + \Low\DSpace{U}
\end{equation*}
and, therefore,
\begin{equation}
\label{U_quotient_spaces}
(\DSpace{U} + \Ft_{n - 1}(\Galg)) / (\Dp\DSpace{U} + \Ft_{n - 1}(\Galg)) \cong \DSpace{U} / (\Dp\DSpace{U} + \Low\DSpace{U}).
\end{equation}

Combining the above results we obtain
\begin{multline}
\label{tensor_product_final_isom}
A_1[U]/(D_1[U] + L_1[U]) \otimes \ldots \otimes A_m[U]/(D_m[U] + L_m[U]) \cong\\
\cong\DSpace{U} / (\Dp\DSpace{U} + \Low\DSpace{U}).
\end{multline}
Let $a^{(i)} \in ME_i[U]$, $i = 1, \ldots, m$. Then it follows from the definitions of $\overline{\mu}[U]$, $\psi[U]$, and the canonical isomorphism~\eqref{U_quotient_spaces} that the isomorphism of vector spaces~\eqref{tensor_product_final_isom} takes each
\begin{equation*}
(a^{(1)} + D_1[U] + L_1[U])\otimes\ldots\otimes(a^{(m)} + D_m[U] + L_m[U])
\end{equation*}
to
\begin{equation*}
\mu[U](a^{(1)}\otimes\ldots\otimes a^{(m)}) + \Dp\DSpace{U} + \Low\DSpace{U}.
\end{equation*}
This completes the proof.
\end{proof}

\subsection{The grading on the space $\Galg$}
\label{grading_subsection}
Recall that $\Dp(\Galg) = \Ideal$. The quotient space $\Qalg$ inherits the filtration from $\Galg$, namely,
\begin{equation*}
\Ft_n (\Qalg) = (\Ft_n(\Galg) + \Dp(\Galg)) / \Dp(\Galg) = (\Ft_n(\Galg) + \Ideal) / \Ideal.
\end{equation*}
We have the corresponding graded spaces\label{grading_def}
\begin{align*}
&\Gr(\Galg) = \bigoplus\limits_{n = 0}^{\infty}\Gr_n(\Galg), \textit{ where } \Gr_n(\Galg) = \Ft_n(\Galg) / \Ft_{n - 1}(\Galg),\\
&\Gr(\Qalg) = \bigoplus\limits_{n = 0}^{\infty}\Gr_n(\Qalg), \\
&\textit{where }\Gr_n(\Qalg) = \Ft_n (\Qalg) / \Ft_{n - 1} (\Qalg).
\end{align*}
It is well-known that $\Galg \cong \Gr(\Galg)$ and $\Qalg \cong \Gr(\Qalg)$ as vector spaces.
\medskip

The following theorem establishes the compatibility of the filtration and the corresponding grading on~$\Galg$ with the space of dependencies~$\Dp(\Galg)$.
\begin{theorem}
\label{structure_of_quotient_space}
\begin{equation*}
\Gr_n(\Qalg) \cong \Ft_n(\Galg) / (\Dp(\Ft_n(\Galg)) + \Ft_{n - 1}(\Galg)).
\end{equation*}
\end{theorem}
\begin{proof}
Using the isomorphism theorems, we obtain
\begin{align}
\label{canonical_isomorphisms_graded_component}
\Gr_n(\Qalg) = & \Ft_n (\Qalg) / \Ft_{n - 1} (\Qalg)\cong\\
& \cong (\Ft_n(\Galg) + \Dp(\Galg)) / (\Ft_{n - 1}(\Galg) + \Dp(\Galg)) =\nonumber\\
& = (\Ft_n(\Galg) + \Ft_{n - 1}(\Galg) + \Dp(\Galg)) / (\Ft_{n - 1}(\Galg) + \Dp(\Galg))\cong\nonumber\\
& \cong \Ft_n(\Galg) / (\Ft_n(\Galg) \cap (\Ft_{n - 1}(\Galg) + \Dp(\Galg))).\nonumber
\end{align}
Since $\Ft_{n - 1}(\Galg) \subseteq \Ft_n(\Galg)$, we have
\begin{equation}
\label{spaces_intersection}
\Ft_n(\Galg) \cap (\Ft_{n - 1}(\Galg) + \Dp(\Galg)) = \Ft_{n - 1}(\Galg) + \Ft_n(\Galg) \cap \Dp(\Galg).
\end{equation}
Therefore, we obtain
\begin{multline*}
\Ft_n(\Galg) / (\Ft_n(\Galg) \cap (\Ft_{n - 1}(\Galg) + \Dp(\Galg))) =\\
= \Ft_n(\Galg) / (\Ft_{n - 1}(\Galg) + \Ft_n(\Galg) \cap \Dp(\Galg)).
\end{multline*}

Since $\Galg$ and $\Ft_n(\Galg)$ are generated by monomials and are closed under taking derived monomials, from Proposition~\ref{fall_to_smaller_subspace} it follows that
\begin{equation*}
\Ft_n(\Galg) \cap \Dp(\Galg) = \Dp(\Ft_n(\Galg)).
\end{equation*}
Therefore,
\begin{equation*}
\Ft_{n - 1}(\Galg) + \Ft_n(\Galg) \cap \Dp(\Galg) = \Ft_{n - 1}(\Galg) + \Dp(\Ft_n(\Galg)).
\end{equation*}
Thus, combining~\eqref{canonical_isomorphisms_graded_component} and~\eqref{spaces_intersection}, we obtain
\begin{equation*}
\Gr_n(\Qalg) \cong \Ft_n(\Galg) / (\Dp(\Ft_n(\Galg)) + \Ft_{n - 1}(\Galg)).
\end{equation*}
\end{proof}

\begin{proposition}[Reducing to the cyclic case]
\label{component_subspaces_structure}
Assume $Y$ is a non-trivial subspace of $\Galg$ generated by monomials and is closed under taking derived monomials such that $f$-characteristics of monomials from $Y$ is bounded. Consider the set of spaces as follows
\begin{equation*}
\lbrace\DSpace{Z} \mid Z \in \Fr, Z \in Y \setminus \Low(Y)\rbrace.
\end{equation*}
Let $\lbrace V_i\rbrace_{i \in I}$ be all the different spaces from the above set (some different $Z$ may give the same $V_i$). Then
\begin{enumerate}[label={(\arabic*)}]
\item
\label{component_subspaces_structure_1}
\begin{equation}
\label{isomorphish_separate_spaces}
Y / (\Dp(Y) + \Low(Y)) \cong \bigoplus\limits_{i \in I}V_i/(\Dp(V_i) + \Low(V_i)).
\end{equation}
\item
\label{component_subspaces_structure_2}
Assume $Z \in Y$ is a monomial such that $Z \notin \Low(Y)$, and $Z \in V_i\cap V_j$. Then $i = j$. And then in~\eqref{isomorphish_separate_spaces}
\begin{equation*}
Z + \Dp(Y) + \Low(Y) \mapsto (0, \ldots, 0, \underbrace{Z + \Dp(V_i) + \Low(V_i)}_{\textit{i-th place}}, 0, \ldots).
\end{equation*}
\end{enumerate}
\end{proposition}
\begin{proof}
Let $\lbrace X_s \rbrace_{s \in S}$ be all the monomials of $Y$. Since $Y$ is generated by monomials, we have
\begin{equation*}
Y = \bigoplus_{s \in S} \langle X_s\rangle.
\end{equation*}
Since $\langle X_s\rangle \subseteq \DSpace{X_s}$, we also have
\begin{equation}
\label{monomial_spaces_sum}
Y = \sum_{s \in S} \DSpace{X_s}.
\end{equation}
Consider the quotient space $Y / \Low(Y)$. Using~\eqref{monomial_spaces_sum}, we obtain
\begin{equation}
\label{quotient_monomial_derived_spaces_sum}
Y  / \Low(Y) = \sum_{s \in S} (\DSpace{X_s} + \Low(Y)) / \Low(Y).
\end{equation}

Assume $X_s \in \Low(Y)$. Since $\Low(Y)$ is generated by monomials and closed under taking derived monomials, we have $\DSpace{X_s} \subseteq \Low(Y)$. So, $(\DSpace{X_s} + \Low(Y)) / \Low(Y)$ is trivial. Therefore, in the sum~\eqref{quotient_monomial_derived_spaces_sum} we can take only spaces $(\DSpace{X_s} + \Low(Y)) / \Low(Y)$ such that $X_s \notin \Low(Y)$ and obtain the same resulting space. Assume $\DSpace{X_{s_1}} = \DSpace{X_{s_2}}$, $s_1 \neq s_2$. Then, obviously, $(\DSpace{X_{s_1}} + \Low(Y)) / \Low(Y) = (\DSpace{X_{s_2}} + \Low(Y)) / \Low(Y)$. Clearly, we can take only one of them in the sum~\eqref{quotient_monomial_derived_spaces_sum}, and again obtain the same resulting space. Therefore, finally we have
\begin{equation}
\label{v_spaces_quotient_low}
Y  / \Low(Y) = \sum_{i \in I} (V_i + \Low(Y)) / \Low(Y).
\end{equation}

Let us show that~\eqref{v_spaces_quotient_low} is a direct sum. Let $\overline{T}_j \in (V_{i_j} + \Low(Y)) / \Low(Y)$, $j = 1, \ldots, l$, $i_j \neq i_{j^{\prime}}$ if $j \neq j^{\prime}$. Assume $\sum_{j = 1}^l\gamma_j\overline{T}_j = 0$. It is enough to show that $\overline{T}_j = 0$, $j = 1, \ldots, l$. Clearly,
\begin{equation*}
\overline{T}_j = T_j + \Low(Y),
\end{equation*}
where $T_j$ is a linear combination of monomials of $V_{i_j}$. Since $\sum_{j = 1}^l\gamma_j\overline{T}_j = 0$, we have
\begin{equation}
\label{up_mon_sum}
\sum_{j = 1}^l\gamma_j T_j \in \Low(Y).
\end{equation}
Our aim is to show that every $T_j \in \Low(Y)$, $j = 1, \ldots, l$.

Let $V_i = \DSpace{Z_i}$, where $Z_i \in \Fr$, $Z_i \in Y \setminus \Low(Y)$. Then we obtain $T_j \in V_{i_j} = \DSpace{Z_{i_j}}$. Since $\sum_{j = 1}^l\gamma_j T_j \in \Low(Y)$, all monomials of $T_j$, $j = 1, \ldots, l$, that belong to $Y \setminus \Low(Y)$ cancel out in the sum~\eqref{up_mon_sum}. Let  $X$ be a monomial and $X \in \DMUp{Z_{i_j}}$ (see~\eqref{upper_monomials}). Assume $X \in \Low(Y)$. Then it follows from Lemma~\ref{derived_spaces_equality} that $\DSpace{Z_{i_j}} \subseteq \Low(Y)$. Hence, $Z_{i_j} \in \Low(Y)$, a contradiction. Therefore, $\DMUp{Z_{i_j}} \subseteq Y \setminus \Low(Y)$. So, all monomials of $T_j$ that belong to $\DMUp{Z_{i_j}}$, $j = 1, \ldots, l$, cancel out in the sum~\eqref{up_mon_sum}.

Assume $Z$ is a monomial such that $Z \in Y \setminus \Low(Y)$. Assume $Z \in V_{i} = \DSpace{Z_{i}}$ and $Z \in V_{j} = \DSpace{Z_{j}}$, $i, j \in I$, $i \neq j$. Since $Z \notin \Low(Y)$, we have $Z \in \DMUp{Z_{i}}$ and $Z \in \DMUp{Z_{j}}$. However, then it follows from Lemma~\ref{derived_spaces_equality} that $\DSpace{Z_{i}} = \DSpace{Z}$ and $\DSpace{Z_{j}} = \DSpace{Z}$. Therefore, $\DSpace{Z_{i}} = \DSpace{Z_{j}}$, a contradiction. Thus, $Z$ belongs to precisely one space $V_i$, $i \in I$.

Assume $U$ is a monomial in $T_{j}$ and $U \in \DMUp{Z_{i_{j}}} \subseteq V_{i_{j}}$. Since $\DMUp{Z_{i_{j}}} \subseteq Y \setminus \Low(Y)$, we have $U \in Y \setminus \Low(Y)$. Hence, it follows from the above that $U$ is not contained in any $V_i$ for $i \neq i_{j}$. Let $1 \leqslant j^{\prime} \leqslant l$ and $j^{\prime} \neq j$. Then $U$ is not contained in $T_{j^{\prime}}$, since $T_{j^{\prime}}$ is a linear combination of monomials of $V_{i_{j^{\prime}}}$ and $i_j \neq i_{j^{\prime}}$. So, if $U$ is cancelled in the sum~\eqref{up_mon_sum}, then $U$ is cancelled in $T_{j}$. That is, $T_{j} \in \Low(V_{i_{j}}) \subseteq \Low(Y)$. Hence, every $T_j \in \Low(Y)$, $j = 1, \ldots, l$. Thus, the sum~\eqref{v_spaces_quotient_low} is a direct sum.

Consider the space $\Dp(Y)$. Recall that $\Dp(Y) = \langle \GDp^{\prime}(Y)\rangle$. Let $T \in \GDp^{\prime}(Y)$. That is, $T$ is a layout of a multi-turn of a virtual member of the chart of some monomial $Z \in Y$. By definition, all monomials of $T$ are derived monomials of $Z$. Therefore, since $\DSpace{Z}$ is closed under taking derived monomials, we have $T \in \GDp^{\prime}(\DSpace{Z}) \subseteq \Dp\DSpace{Z}$. Recall that $\lbrace X_s \rbrace_{s \in S}$ are all the monomials of $Y$. So, we obtain
\begin{equation}
\label{dep_monomials_sum}
\Dp(Y) = \sum_{s \in S} \Dp\DSpace{X_s}.
\end{equation}
Consider the space $(\Dp(Y) + \Low(Y)) / \Low(Y)$. By~\eqref{dep_monomials_sum}, we have
\begin{equation*}
(\Dp(Y) + \Low(Y)) / \Low(Y) = \sum_{s \in S} (\Dp\DSpace{X_s} + \Low(Y)) / \Low(Y).
\end{equation*}
Arguing as at the beginning of the proof, we obtain
\begin{equation}
\label{dep_sum_quotient_low}
(\Dp(Y) + \Low(Y)) / \Low(Y) = \sum\limits_{i \in I} (\Dp(V_i) + \Low(Y)) / \Low(Y).
\end{equation}
Since $(\Dp(V_i) + \Low(Y)) / \Low(Y) \subseteq (V_i + \Low(Y)) / \Low(Y)$ and~\eqref{v_spaces_quotient_low} is a direct sum, we obtain that~\eqref{dep_sum_quotient_low} is a direct sum as well.

Consider the space $Y / (\Dp(Y) + \Low(Y))$. It follows from the isomorphism theorem that
\begin{equation*}
Y / (\Dp(Y) + \Low(Y)) \cong (Y / \Low(Y)) / ((\Dp(Y) + \Low(Y)) / \Low(Y)).
\end{equation*}
We proved above that
\begin{align*}
&Y / \Low(Y) = \bigoplus\limits_{i \in I} (V_i + \Low(Y)) / \Low(Y),\\
&(\Dp(Y) + \Low(Y)) / \Low(Y) = \bigoplus\limits_{i \in I} (\Dp(V_i) + \Low(Y)) / \Low(Y).
\end{align*}
Hence, we have the following sequence of isomorphisms of vector spaces
\begin{align}
\label{direct_sum_components}
Y /& (\Dp(Y) + \Low(Y)) \cong (Y / \Low(Y)) / ((\Dp(Y) + \Low(Y)) / \Low(Y)) = \\
&= \left(\bigoplus\limits_{i \in I} (V_i + \Low(Y)) / \Low(Y) \right) / \left(\bigoplus\limits_{i \in I} (\Dp(V_i) + \Low(Y)) / \Low(Y) \right) \cong \nonumber \\
&\cong \bigoplus\limits_{i \in I} ((V_i + \Low(Y)) / \Low(Y)) / ((\Dp(V_i) + \Low(Y)) / \Low(Y)) \cong \nonumber \\
&\cong \bigoplus\limits_{i \in I} (V_i + \Low(Y)) / (\Dp(V_i) + \Low(Y)).\nonumber
\end{align}
By the isomorphism theorem, we obtain
\begin{align}
\label{single_component}
(V_i + \Low(Y)) / &(\Dp(V_i) + \Low(Y)) = (V_i + \Dp(V_i) + \Low(Y)) / (\Dp(V_i) + \Low(Y)) \cong \\
&\cong V_i / ((\Dp(V_i) + \Low(Y)) \cap V_i).\nonumber
\end{align}
Since $\Dp(V_i) \subseteq V_i $, we have
\begin{equation}
\label{cyclic_spaces_intersection}
(\Dp(V_i) + \Low(Y)) \cap V_i = \Dp(V_i) + \Low(Y) \cap V_i.
\end{equation}

Let us show that $\Low(Y) \cap V_i = \Low(V_i)$. Assume $Y \subseteq \Ft_n(\Galg)$ and $Y \nsubseteq \Ft_{n - 1}(\Galg)$.  By definition, $\Low(Y) = Y \cap \Ft_{n - 1}(\Galg)$, so,
\begin{equation*}
\Low(Y) \cap V_i = Y \cap \Ft_{n - 1}(\Galg) \cap V_i.
\end{equation*}
Recall that $V_i = \DSpace{Z_i}$, where $Z_i$ is a monomial, and $Z_i \in Y\setminus \Low(Y)$. Therefore, $V_i \subseteq Y \subseteq \Ft_{n}(\Galg)$ and $V_i \nsubseteq \Ft_{n - 1}(\Galg)$. Hence, $\Ft_{n - 1}(\Galg) \cap V_i = \Low(V_i)$. Since $\Low(V_i) \subseteq V_i \subseteq Y$, we finally obtain
\begin{equation*}
\Low(Y) \cap V_i = Y \cap \Ft_{n - 1}(\Galg) \cap V_i = Y \cap \Low(V_i) = \Low(V_i).
\end{equation*}

Applying the equality $\Low(Y) \cap V_i = \Low(V_i)$ to~\eqref{cyclic_spaces_intersection}, we see that
\begin{equation*}
(\Dp(V_i) + \Low(Y)) \cap V_i = \Dp(V_i) + \Low(V_i).
\end{equation*}
Thus, using~\eqref{direct_sum_components} and~\eqref{single_component}, we finally obtain
\begin{equation*}
Y / (\Dp(Y) + \Low(Y)) \cong \bigoplus\limits_{i \in I} V_i/(\Dp(V_i) + \Low(V_i)).
\end{equation*}
So, the first statement of Proposition~\ref{component_subspaces_structure} is proved.

Assume $Z$ is a monomial such that $Z \in Y \setminus \Low(Y)$. We proved above that $Z$ belongs to precisely one space $V_i$, $i \in I$. Then, obviously,
\begin{equation*}
Z + \Low(Y) \in (V_i + \Low(Y)) / \Low(Y).
\end{equation*}
Thus, by the isomorphism theorems, we obtain that the sequence of the canonical isomorphisms~\eqref{direct_sum_components} and~\eqref{single_component} acts on $Z + \Dp(Y) + \Low(Y)\in Y / (\Dp(Y) + \Low(Y))$ in the following way:
\begin{align*}
Z + \Dp(Y)& + \Low(Y) \mapsto Z + \Low(Y) + (\Dp(Y) + \Low(Y)) / \Low(Y) \mapsto \\
&\mapsto (0, \ldots, 0,  \underbrace{Z + L(Y)}_{\textit{i-th place}}, 0, \ldots) + \bigoplus\limits_{i \in I} (\Dp(V_i) + \Low(Y)) / \Low(Y) \mapsto\\
&\mapsto (0, \ldots, 0, \underbrace{Z + L(Y) + (\Dp(V_i) + \Low(Y)) / \Low(Y)}_{\textit{i-th place}}, 0, \ldots) \mapsto\\
&\mapsto (0, \ldots, 0, \underbrace{Z + \Dp(V_i) + \Low(Y)}_{\textit{i-th place}}, 0, \ldots) \mapsto \\
&\mapsto (0, \ldots, 0, \underbrace{Z + (\Dp(V_i) + \Low(Y)) \cap V_i}_{\textit{i-th place}}, 0, \ldots) = \\
&= (0, \ldots, 0,\underbrace{Z + \Dp(V_i) + \Low(V_i)}_{\textit{i-th place}}, 0, \ldots).
\end{align*}
So, the second statement of Proposition~\ref{component_subspaces_structure} is proved.
\end{proof}

\section{Construction of a basis of $\Qalg$ (ensuring the non-triviality of $\Qalg$)}
\label{basis_section}
\subsection{Non-triviality of $\Qalg$}
\label{quotient_space_non_triviality_subsection}
\begin{lemma}
\label{non_trivial_spaces_existence}
Let $\lbrace V_i\rbrace_{i \in I}$ be all the different spaces $\lbrace\DSpace{Z} \mid Z \in \Fr\rbrace$. Then not all spaces $V_i / (\Dp(V_i) + \Low(V_i))$, $i\in I$, are trivial. Namely, the space $\DSpace{X} / (\Dp\DSpace{X} + \Low\DSpace{X})$, where $X$ is a monomial with no virtual members of the chart, is always non-trivial, and is of dimension~$1$. In particular, $\DSpace{1} / (\Dp\DSpace{1} + \Low\DSpace{1}) \neq 0$, where $1$ is the empty word.
\end{lemma}
\begin{proof}
Let $X$ be a monomial with no virtual members of the chart. Then there are no derived monomials of $X$ except $X$ itself, and there are no multi-turns of virtual members of the chart of $X$. So, by definition, $\DSpace{X}$ is linearly generated by $X$ and, therefore, it is of dimension~$1$; $\Dp\DSpace{X} = 0$; $\Low\DSpace{X} = 0$. Therefore,
\begin{equation*}
\DSpace{X} / (\Dp\DSpace{X} + \Low\DSpace{X}) = \DSpace{X} = \langle X\rangle \neq 0,
\end{equation*}
and $\DSpace{X} / (\Dp\DSpace{X} + \Low\DSpace{X}) $ is of dimension~$1$.

By definition, the empty word~$1$ is a small piece. Therefore, $1$ has no virtual members of the chart. So, it follows from the above that $\DSpace{1} / (\Dp\DSpace{1} + \Low\DSpace{1}) \neq 0$.
\end{proof}

\begin{remark}
Notice that, by the definition of small pieces, there always exists at least one small piece. Namely, the monomial $1$ is always a small piece. Let us emphasise that this fact plays the crucial role in the argument of Lemma~\ref{non_trivial_spaces_existence}.
\end{remark}

Now we can prove that the quotient ring $\Qalg$ is non-trivial.
\begin{corollary}
\label{non_trivial_quotient}
The quotient ring $\Qalg$ is non-trivial.
\end{corollary}
\begin{proof}
Let $U$ be a monomial. Consider the space $\DSpace{U}$ and the corresponding subspace in $\Qalg$, namely, $(\DSpace{U} + \Ideal) / \Ideal$. From the isomorphism theorem it follows that
\begin{equation*}
(\DSpace{U} + \Ideal) / \Ideal \cong \DSpace{U} / (\DSpace{U}\cap \Ideal).
\end{equation*}
Recall that $\Ideal = \Dp(\Galg)$. From Proposition~\ref{fall_to_smaller_subspace} it follows that $\DSpace{U} \cap \Dp(\Galg) = \Dp\DSpace{U}$. Hence,
\begin{equation*}
(\DSpace{U} + \Ideal) / \Ideal \cong \DSpace{U} / \Dp\DSpace{U}.
\end{equation*}

By Lemma~\ref{non_trivial_spaces_existence}, there exists a space $\DSpace{U_0}$, $U_0 \in \Fr$, such that $\DSpace{U_0} / (\Dp\DSpace{U_0} + \Low\DSpace{U_0}) \neq 0$. Hence, we see that $\DSpace{U_0} / \Dp\DSpace{U_0} \neq 0$ and $(\DSpace{U_0} + \Ideal) / \Ideal \neq 0$. So, there exists a non-trivial subspace of $\Qalg$. Thus, $\Qalg$ itself is non-trivial.
\end{proof}

\subsection{Construction of a basis of $\Qalg$}
\label{construction_of_basis_subsection}
Now we show how to construct a basis of $\Qalg$. First we construct a basis for non-trivial graded components
\begin{equation*}
\Gr_n(\Qalg) = \Ft_n(\Qalg) / \Ft_{n - 1}(\Qalg),
\end{equation*}
where $\Ft_n(\Galg)$ is a filtration defined by~\eqref{filtration_def}.
\begin{proposition}
\label{graded_component_basis}
Let $n$ be a level of the filtration $\Ft_n(\Galg)$ defined by~\eqref{filtration_def}. We consider the set of spaces
\begin{equation*}
\lbrace\DSpace{Z} \mid Z \in \Fr, Z \in \Ft_n(\Galg)  \setminus \Ft_{n - 1}(\Galg) , \DSpace{Z} / (\Dp\DSpace{Z} + \Low\DSpace{Z}) \neq 0\rbrace.
\end{equation*}
Let $\lbrace V_{i}^{(n)}\rbrace_{i \in I^{(n)}}$ be all the different spaces from the above set. Then $\Gr_n(\Qalg)$ is non-trivial if and only if $\lbrace V_{i}^{(n)}\rbrace_{i \in I^{(n)}} \neq \varnothing$. If $\lbrace V_{i}^{(n)}\rbrace_{i \in I^{(n)}} \neq \varnothing$, then we have
\begin{equation*}
\Gr_n(\Qalg) \cong \bigoplus\limits_{i \in I^{(n)}} V_{i}^{(n)} /(\Dp(V_{i}^{(n)}) + \Low(V_{i}^{(n)})).
\end{equation*}

Assume $\lbrace \overline{W}^{(i, n)}_j\rbrace_j$ is a basis of $V_{i}^{(n)}/(\Dp(V_{i}^{(n)}) + \Low(V_{i}^{(n)}))$, $i \in I^{(n)}$. Let $W^{(i, n)}_j \in V_{i}^{(n)}$ be an arbitrary representative of the coset $\overline{W}^{(i, n)}_j$. Then
\begin{equation*}
\bigcup\limits_{i \in I^{(n)}}\left\lbrace W^{(i, n)}_j + \Ideal + \Ft_{n - 1}(\Qalg)\right\rbrace_j
\end{equation*}
is a basis of $\Gr_n(\Qalg)$.
\end{proposition}

\begin{proof}
While the statement is pretty obvious, we prefer to give a proof to recollect the previously stated facts.

Recall that we constructed the following sequence of canonical isomorphisms of vector spaces in Theorem~\ref{structure_of_quotient_space} (see~\eqref{canonical_isomorphisms_graded_component}):
\begin{align}
\label{graded_component_isomorphisms}
&\Gr_n(\Qalg) = \Ft_n(\Qalg) / \Ft_{n - 1}(\Qalg) = \\
&=((\Ft_n(\Galg) + \Ideal)/\Ideal) / ((\Ft_{n - 1}(\Galg) + \Ideal)/\Ideal) \to \nonumber\\
&\to (\Ft_n(\Galg) + \Ideal) / (\Ft_{n - 1}(\Galg) + \Ideal) = (\Ft_n(\Galg) + \Dp(\Galg)) / (\Ft_{n - 1}(\Galg) + \Dp(\Galg))\to\nonumber\\
&\to \Ft_n(\Galg) / (\Ft_{n - 1}(\Galg) + \Dp(\Galg) \cap \Ft_n(\Galg)) = \Ft_n(\Galg) / (\Ft_{n - 1}(\Galg) + \Dp(\Ft_n(\Galg))).\nonumber
\end{align}
It follows from Proposition~\ref{component_subspaces_structure} that if $\lbrace V_{i}^{(n)}\rbrace_{i \in I^{(n)}} = \varnothing$, then $\Ft_n(\Galg) / (\Ft_{n - 1}(\Galg) + \Dp(\Ft_n(\Galg)))$ is trivial; otherwise $\Ft_n(\Galg) / (\Ft_{n - 1}(\Galg) + \Dp(\Ft_n(\Galg)))$ is non-trivial, and we obtain
\begin{equation}
\label{graded_component_final_isomorphism}
\Ft_n(\Galg) / (\Ft_{n - 1}(\Galg) + \Dp(\Ft_n(\Galg))) \cong \bigoplus\limits_{i \in I^{(n)}} V_{i}^{(n)} /(\Dp(V_{i}^{(n)}) + \Low(V_{i}^{(n)})).
\end{equation}
Thus, $\Gr_n(\Qalg)$ is non-trivial if and only if $\lbrace V_{i}^{(n)}\rbrace_{i \in I^{(n)}} \neq \varnothing$. And if $\lbrace V_{i}^{(n)}\rbrace_{i \in I^{(n)}} \neq \varnothing$, then we have
\begin{equation}
\label{graded_component_isomorphism_combine_all}
\Gr_n(\Qalg) \cong \bigoplus\limits_{i \in I^{(n)}} V_{i}^{(n)} /(\Dp(V_{i}^{(n)}) + \Low(V_{i}^{(n)})).
\end{equation}

Let $W$ belong to some $V_{i}^{(n)}$, that is, $W$ is a linear combination of monomials of $V_{i}^{(n)}$. Then, clearly, the sequence of isomorphisms~\eqref{graded_component_isomorphisms} acts on
\begin{equation*}
W + \Ideal + \Ft_{n - 1}(\Qalg) \in \Gr_n(\Qalg) = \Ft_n(\Qalg) / \Ft_{n - 1}(\Qalg)
\end{equation*}
in the following way:
\begin{align*}
&W + \Ideal + \Ft_{n - 1}(\Qalg) \mapsto W + \Ideal + (\Ft_{n - 1}(\Galg) + \Ideal) \mapsto\\
&\mapsto  W + \Ft_{n - 1}(\Galg) + \Dp(\Ft_n(\Galg)).\nonumber
\end{align*}
By Proposition~\ref{component_subspaces_structure}, the isomorphism~\eqref{graded_component_final_isomorphism} acts on $W + \Ft_{n - 1}(\Galg) + \Dp(\Ft_n(\Galg))$ as follows:
\begin{equation*}
W + \Ft_{n - 1}(\Galg) + \Dp(\Ft_n(\Galg)) \mapsto (0, \ldots, 0, \underbrace{W + \Dp(V^{(n)}_{i}) + \Low(V^{(n)}_{i})}_{i-\textit{th place}}, 0, \ldots).
\end{equation*}
Combining two last mappings, we obtain that the isomorphism~\eqref{graded_component_isomorphism_combine_all} acts as follows:
\begin{equation}
\label{final_isom_fron_graded_component}
W + \Ideal + \Ft_{n - 1}(\Qalg) \mapsto (0, \ldots, 0, \underbrace{W + \Dp(V^{(n)}_{i}) + \Low(V^{(n)}_{i})}_{i-\textit{th place}}, 0, \ldots).
\end{equation}

Assume $\lbrace \overline{W}^{(i, n)}_j\rbrace_j$ is a basis of $V_{i}^{(n)} / (\Dp(V^{(n)}_{i}) + \Low(V^{(n)}_{i}))$. Let
\begin{equation*}
\overline{W}^{(i, n)}_j = W^{(i, n)}_j+ \Dp(V^{(n)}_{i}) + \Low(V^{(n)}_{i}),
\end{equation*}
where $W^{(i, n)}_j \in V^{(n)}_{i}$ is an arbitrary representative of the coset $\overline{W}^{(i, n)}_j$. Using~\eqref{final_isom_fron_graded_component}, we obtain
\begin{align*}
W^{(i, n)}_j + \Ideal + \Ft_{n - 1}(\Qalg) &\mapsto (0, \ldots, 0, \underbrace{W^{(i, n)}_j + \Dp(V^{(n)}_{i}) + \Low(V^{(n)}_{i})}_{i-\textit{th place}}, 0, \ldots) =\\
&=(0, \ldots, 0, \underbrace{\overline{W}^{(i, n)}_j}_{i-\textit{th place}}, 0, \ldots).
\end{align*}
Therefore, since
\begin{equation*}
\bigcup_{i \in I^{(n)}} \left\lbrace(0, \ldots, 0, \underbrace{\overline{W}^{(i, n)}_j}_{i-\textit{th place}}, 0, \ldots)\right\rbrace_j
\end{equation*}
is a basis of $\bigoplus\limits_{i \in I^{(n)}} V_{i}^{(n)} /(\Dp(V_{i}^{(n)}) + \Low(V_{i}^{(n)}))$, we see that
\begin{equation*}
\bigcup_{i \in I^{(n)}} \lbrace W^{(i,n)}_j + \Ideal + \Ft_{n - 1}(\Qalg) \rbrace_j
\end{equation*}
is a basis of $\Gr_n(\Qalg) = \Ft_n(\Qalg) / \Ft_{n - 1}(\Qalg)$. This completes the proof.
\end{proof}

The following theorem describes a structure of $\Qalg$ as a vector space and describes a basis of $\Qalg$.
\begin{theorem}
\label{whole_quotient_ring_structure}
Consider the set of spaces as follows
\begin{equation*}
\lbrace\DSpace{Z} \mid Z \in \Fr, \DSpace{Z} / (\Dp\DSpace{Z} + \Low\DSpace{Z}) \neq 0\rbrace.
\end{equation*}
Let $\lbrace V_i\rbrace_{i \in I}$ be all the different spaces from the above set ($V_{i_1} \neq V_{i_2}$ for $i_1 \neq i_2$). Then $\lbrace V_i\rbrace_{i \in I} \neq \varnothing$ and we have
\begin{equation*}
\Qalg \cong \bigoplus\limits_{i \in I} V_i/(\Dp(V_i) + \Low(V_i))
\end{equation*}
as vector spaces, and the right-hand side is explicitly described in Proposition~\ref{correspondence_to_tensor_product}.

Assume $\lbrace \overline{W}^{(i)}_j \rbrace_j$ is a basis of $V_i/(\Dp(V_i) + \Low(V_i))$, $i \in I$. Let $W^{(i)}_j \in V_i$ be an arbitrary representative of the coset $\overline{W}^{(i)}_j$. Then
\begin{equation*}
\bigcup\limits_{i \in I}\left\lbrace W^{(i)}_j  + \Ideal\right\rbrace_j
\end{equation*}
is a basis of $\Qalg$.
\end{theorem}
\begin{proof}
Recall that we have
\begin{multline}
\label{graded_space}
\Qalg \cong \Gr(\Qalg) = \bigoplus_{n = 0}^{\infty}\Gr_n(\Qalg),\\
\textit{ where } \Gr_n(\Qalg) = \Ft_n(\Qalg) / \Ft_{n - 1}(\Qalg).
\end{multline}
For every $n \in \lbrace 0 \rbrace \cup \mathbb{N}$ we consider the set of spaces
\begin{equation*}
\lbrace\DSpace{Z} \mid Z \in \Fr, Z \in \Ft_n(\Galg)  \setminus \Ft_{n - 1}(\Galg) , \DSpace{Z} / (\Dp\DSpace{Z} + \Low\DSpace{Z}) \neq 0\rbrace.
\end{equation*}
Let $\lbrace V_{i}^{(n)}\rbrace_{i \in I^{(n)}}$ be all the different spaces from the above set. Proposition~\ref{graded_component_basis} implies that $\Gr_n(\Qalg)$ is non-trivial if and only if $\lbrace V_{i}^{(n)}\rbrace_{i \in I^{(n)}} \neq \varnothing$; if $\lbrace V_{i}^{(n)}\rbrace_{i \in I^{(n)}} \neq \varnothing$, then
\begin{equation}
\label{graded_component_final_isomorphism_recall}
\Gr_n(\Qalg) \cong \bigoplus\limits_{i \in I^{(n)}} V_{i}^{(n)} /(\Dp(V_{i}^{(n)}) + \Low(V_{i}^{(n)})).
\end{equation}
Since every monomial belongs to $\Ft_n(\Galg) \setminus \Ft_{n - 1}(\Galg)$ for some $n$, we have
\begin{equation}
\label{all_levels_spaces_union}
\bigsqcup_{n = 0}^{\infty}\lbrace V_{i}^{(n)}\rbrace_{i \in I^{(n)}} = \lbrace V_{i}\rbrace_{i \in I}.
\end{equation}
Therefore, combining~\eqref{graded_space} and~\eqref{graded_component_final_isomorphism_recall}, we obtain
\begin{equation*}
\Gr(\Qalg) = \bigoplus\limits_{i \in I} V_i/(\Dp(V_i) + \Low(V_i)).
\end{equation*}

Let us construct a basis of $\Qalg$. Let $\lbrace \overline{W}^{(i, n)}_j \rbrace_j$ be a basis of $V_i^{(n)}$, $i \in I^{(n)}$. Then it follows from Proposition~\ref{graded_component_basis} that
\begin{equation*}
\bigcup_{i \in I^{(n)}} \lbrace W^{(i, n)}_j + \Ideal + \Ft_{n - 1}(\Galg) \rbrace_j
\end{equation*}
is a basis of $\Gr_n(\Qalg)$, where $W^{(i, n)}_j$ is an arbitrary representative of the coset $\overline{W}^{(i, n)}_j$. It is well-known that $\Qalg$ is isomorphic as a vector space to $\Gr(\Qalg) = \bigoplus_{n = 0}^{\infty}\Gr_n(\Qalg)$. Although there is no canonical isomorphism between them, we have the following correspondence. Assume $\lbrace \overline{e}_j^{(n)}\rbrace_j$ is a basis of a non-trivial graded component $\Gr_n(\Qalg)$. If $\Gr_n(\Qalg)$ is trivial, for the matter of convenience we assume that $\lbrace \overline{e}_j^{(n)}\rbrace_j = \varnothing$. Then $\cup_{n = 0}^{\infty} \lbrace e_j^{(n)}\rbrace_j$ is a basis of $\Qalg$, where $e_j^{(n)} \in \Ft_n(\Qalg)$ is an arbitrary representative of the coset $\overline{e}_j^{(n)}$. Therefore,
\begin{equation*}
\bigcup_{n = 0}^{\infty}\bigcup_{i \in I^{(n)}} \lbrace W^{(i,n)}_j + \Ideal\rbrace_j
\end{equation*}
is a basis of $\Qalg$.

Using~\eqref{all_levels_spaces_union}, we see that
\begin{equation*}
\bigcup_{n = 0}^{\infty}\bigcup_{i \in I^{(n)}} \lbrace W^{(i, n)}_j + \Ideal \rbrace_j = \bigcup\limits_{i \in I}\left\lbrace W^{(i)}_j  + \Ideal\right\rbrace_j.
\end{equation*}
Thus, $\bigcup\limits_{i \in I}\left\lbrace W^{(i)}_j  + \Ideal\right\rbrace_j$ is a basis of $\Qalg$. This completes the proof.
\end{proof}

\section{Greedy Algorithm based on $f$-characteristic and Ideal Membership Problem}
\label{algorithm_section}
We take \emph{an additive closure of the set of generators $\Rel$} in the following sense. Assume two polynomials $p, q \in \Rel$ have a common monomial of $\SPM$-measure $\geqslant \tau - 2$. Denote this monomial by $a$. Assume $\alpha$ is the coefficient of $a$ in $p$ and $\beta$ is the coefficient of $a$ in $q$. Then we add to $\Rel$ additively reduced linear combinations of the form $\gamma(\alpha^{-1}p - \beta^{-1}q)$, where $\gamma$ is an arbitrary element of the field $\fld$. Notice that $a$ cancels out in such linear combinations. As a result we obtain the set of generators of $\Ideal$ of the form
\begin{align*}
\Rel \cup \lbrace \gamma(\alpha^{-1}p - \beta^{-1}q) \mid& \gamma \in \fld, p, q \in \Rel,\  p, q \textit{ have a common monomial}\\
&\textit{of } \SPM\textit{-measure} \geqslant \tau - 2 \textit{ with a coefficient } \alpha \textit{ in } p\\
&\textit{and with a coefficient } \beta \textit{ in } q\rbrace.
\end{align*}
We repeat the same procedure for the obtained set, then for the set obtained after the second step, etc. Let $\Add(\Rel)$\label{add_def} be the union of sets of generators obtained after every step.

\begin{remark}
\label{strong_additive_closed}
It can be more convenient to consider a stronger additive closure of $\Rel$. Namely, assume two polynomials $p, q \in \Rel$ have a common monomial of $\SPM$-measure $\geqslant \tau - 2$. Then we add all their possible additively reduced linear combinations $\gamma p + \delta q$, $\gamma, \delta \in \fld$, and obtain a set of generators
\begin{align*}
\Rel \cup \lbrace \gamma p + \delta q \mid &\gamma, \delta \in \fld,\ p, q \in \Rel,\  p, q \textit{ have a common monomial}\\
&\textit{of } \SPM\textit{-measure} \geqslant \tau - 2 \textit{ with a coefficient }\rbrace.
\end{align*}
We repeat the same procedure for the obtained set, then for the set obtained after the second step, etc. The desired result is the union of sets of generators obtained after every step.
\end{remark}

Clearly, the set of all monomials of $\Add(\Rel)$ is equal to $\Mon$. We consider $\SPM$-measure on monomials of $\Add(\Rel)$ with respect to the initial set of small pieces $\SP$. Every element of $\Add(\Rel)$ is, in fact, a linear combination of elements of $\Rel$. Therefore, it follows from Small Cancellation Axiom that every linear combination of elements of $\Add(\Rel)$ after additive cancellations is either trivial, or contains a monomial of $\SPM$-measure $\geqslant \tau + 1$.

We do not claim that $\Add(\Rel)$ necessarily satisfies Compatibility Axiom.

The initial set of generators $\Rel$ itself can be additively closed in the above sense. That is, it is possible that $\Add(\Rel) = \Rel$. The initial set of generators even can be closed in the stronger sense formulated in Remark~\ref{strong_additive_closed}. Corresponding natural examples are given in Section~\ref{examples_section}.

In the same way as in Definition~\ref{multi_turn_def}, we can define elementary multi-turns and multi-turns with respect to the set $\Add(\Rel)$. Similarly to $\GDp^{\prime}$ (see~\eqref{linear_dep_members2}), we define the set of layouts $\GDp^{\prime\prime}$ with respect to the set $\Add(\Rel)$. Namely,
\begin{multline*}
\label{add_linear_dep_members}
\GDp^{\prime\prime} = \Bigg\lbrace \sum\limits_{j = 1}^s \alpha_j La_jR\mid L, R \in \Fr, \textit{there exists an index }1\leqslant h \leqslant s \textit{ such that}\\
a_h \textit{ is a virtual member of the chart of } La_hR, \textit{ and } \sum_{j = 1}^s \alpha_j a_j \in \Add(\Rel)\Bigg\rbrace.
\end{multline*}

The complete analogue of property~\ref{a_j_keep_structure} from Section~\ref{mt_configurations} holds for replacements by monomials that belong to one polynomial from $\Add(\Rel)$.
\begin{lemma}
\label{add_rel_incident_keep_structure}
Let $U$ be a monomial, $a \in \mo{U}$, $U = LaR$. Let $p \in \Add(\Rel)$ and $a, b$ be two monomials of $p$. Assume $b$ is not a small piece. Then the monomial $LbR$ is reduced and $b$ is a maximal occurrence in $LbR$.
\end{lemma}
\begin{proof}
Since $p \in \Add(\Rel)$, this means that $p$ belongs to $\Add(\Rel)$ after some number $N$ of steps. Let us prove Lemma~\ref{add_rel_incident_keep_structure} by induction on $N$.

Assume $N = 0$. Then $p \in \Rel$. Therefore, $a$ and $b$ are incident monomials. Hence, the statement of Lemma~\ref{add_rel_incident_keep_structure} follows from property~\ref{a_j_keep_structure} from Section~\ref{mt_configurations}.

Let us make the step of induction. We have $p = \gamma_1q_1 - \gamma_2q_2$, where $q_1$ and $q_2$ belong to $\Add(\Rel)$ after $N - 1$ steps, there exists a common monomial $c$ of $q_1$ and $q_2$ such that $\SPM(c) \geqslant \tau - 2$, and $c$ cancels out in $\gamma_1q_1 - \gamma_2q_2$. Since $a$ and $b$ are monomials of $p$, we get that $a$ belongs to at least one of $q_1$ and $q_2$, and $b$ belongs to at least one of $q_1$ and $q_2$. Assume $a$ and $b$ both belong  to $q_1$, or both belong  to $q_2$. Then the statement of Lemma~\ref{add_rel_incident_keep_structure} follows directly from the induction hypothesis.

Now assume that $a$ and $b$ belong to the different polynomials. Without loss of generality, we can assume that $a$ is contained in $q_1$ and $b$ is contained in $q_2$. First we consider the replacement $a \mapsto c$ in $LaR$. So, $a$ and $c$ are monomials of $q_1$, $q_1$ belongs to $\Add(\Rel)$ after $N - 1$ steps, and $c$ is not a small piece. Therefore, by the induction hypothesis, we obtain that the monomial $LcR$ is reduced and $c \in \mo{LcR}$. Now we consider the replacement $c \mapsto b$ in $LcR$. Similarly, $c$ and $b$ are monomials of $q_2$, $q_2$ belongs to $\Add(\Rel)$ after $N - 1$ steps, and $b$ is not a small piece. Therefore, since $c \in \mo{LcR}$, it follows from the induction hypothesis that the monomial $LbR$ is reduced and $b \in \mo{LbR}$. Lemma~\ref{add_rel_incident_keep_structure} is proved.
\end{proof}

\begin{definition}[\textbf{Order $<_f$}]
\label{f_char_order_def}
Let us define a linear order on monomials based on $f$-characteristic and denote it by $<_f$. Consider the set of spaces $\lbrace \DSpace{Z} \mid Z \in \Fr\rbrace$. We fix monomials $\lbrace Z_i\rbrace_{i \in I}$ such that $\lbrace \DSpace{Z_i}\rbrace_{i \in I}$ are all different spaces from the above set. Then it follows from Lemma~\ref{derived_spaces_equality} that
\begin{equation}
\label{split_by_upper_levels}
\Fr = \bigsqcup_{i \in I}\DMUp{Z_i}.
\end{equation}

First we order monomials $Z_i$ according to their $f$-characteristic. Then we linearly order monomials $Z_i$ with the same $f$-characteristic. We can order them in an arbitrary way. For example, we can take Deglex ordering.

Now we define a linear order on every $\DMUp{Z_i}$. Let $m_i = \nvirt{Z_i}$, $A^{(j)}[Z_i]$ be defined by formula~\eqref{tens_prod_components}, $ME_j[Z_i]$ be defined by formula~\eqref{tens_prod_components_up_monomials}, $j = 1, \ldots, m_i$. We order $ME_j[Z_i]$ first by $\SPM$-measure. Elements of  $ME_j[Z_i]$ which have equal $\SPM$-measure can be ordered in an arbitrary way. For instance, we can take Deglex ordering. We denote the obtained order on $ME_j[Z_i]$ by $<_{\SPM}$. After that we lexicographically order elements $a^{(1)}\otimes \ldots \otimes a^{(m_i)} \in A^{(1)}[Z_i]\otimes \ldots \otimes A^{(m_i)}[Z_i]$ with every $a^{(j)} \in ME_j[Z_i]$, using the order on $ME_j[Z_i]$ introduced above. That is,
\begin{align}
\label{order_inside_up_level}
\begin{split}
b^{(1)}\otimes \ldots \otimes b^{(m_i)} &< c^{(1)}\otimes \ldots \otimes c^{(m_i)} \Longleftrightarrow\\
&b^{(j)} = c^{(j)} \textit{ for } j < j_0 < m_i \textit{ and } b^{(j_0)} <_{\SPM} c^{(j_0)}.
\end{split}
\end{align}
By statement~\ref{mu_properties_up_correspondence} of Lemma~\ref{mu_properties}, we obtain that there exists a bijective correspondence between the monomials of $\DMUp{Z_i}$ and the elements $a^{(1)}\otimes \ldots \otimes a^{(m_i)} \in A^{(1)}[Z_i]\otimes \ldots \otimes A^{(m_i)}[Z_i]$ such that every $a^{(j)} \in ME_j[Z_i]$. Using this bijective correspondence, the order constructed above on the elements $a^{(1)}\otimes \ldots \otimes a^{(m_i)} \in A^{(1)}[Z_i]\otimes \ldots \otimes A^{(m_i)}[Z_i]$ with every $a^{(j)} \in ME_j[Z_i]$ induces  an order on $\DMUp{Z_i}$.

As a result, we define a linear order $<_f$ on all monomials as follows. Consider two monomials $U$ and $V$.
\begin{itemize}
\item
If $f(U) < f(V)$, then $U <_f V$.
\item
Assume $f(U) = f(V)$. It follows from~\eqref{split_by_upper_levels} that there exist monomials $Z_{i_1}$ and $Z_{i_2}$ such that $U \in \DMUp{Z_{i_1}}$ and $V \in \DMUp{Z_{i_2}}$.

Assume $i_1 \neq i_2$. Then either $Z_{i_1} < Z_{i_2}$, or $Z_{i_2} < Z_{i_1}$, where $<$ is Deglex ordering. If $Z_{i_1} < Z_{i_2}$, then $U <_f V$. If $Z_{i_2} < Z_{i_1}$, then $V <_f U$.
\item
Assume $i_1 = i_2$, that is $U, V \in \DMUp{Z_{i_1}}$. Then we order $U$ and $V$, using formula~\eqref{order_inside_up_level}.
\end{itemize}

Clearly, the linear order $<_f$ does not have infinite decreasing chains.
\end{definition}

Let $G = \langle \mathcal{X} \mid \Rel_G\rangle$ be a group given by generators and defining relations. Let $A = M_1M_2^{-1} \in \Rel_G$. Assume $LM_1R$ and $L\cdot M_2\cdot R$ are two words. Recall that the transition from $LM_1R$ to $L\cdot M_2\cdot R$
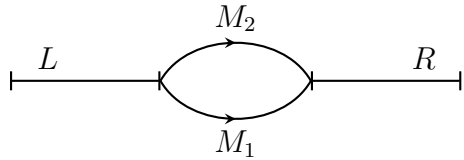
\begin{center}
\begin{tikzpicture}
\draw[|-|, black, thick] (0,0)--(2,0) node [near start, above] {$L$};
\draw[black, thick, arrow=0.5] (2,0) to [bend left=60] node [above] {$M_2$} (4,0);
\draw[black, thick, arrow=0.5] (2,0) to [bend right=60] node [below] {$M_1$} (4,0);
\draw[|-|, black, thick] (4,0)--(6,0) node [near end, above] {$R$};
\end{tikzpicture}
\end{center}
is called a \emph{turn} of an occurrence of the subrelation $M_1$ (to its complement~$M_2$), see~\cite{NA1}.

Let the group $G = \langle \mathcal{X} \mid \Rel_G\rangle$ be given. In what follows we assume that $\Rel_G$ is closed under taking cyclic shifts and inverses of relators, and every relator from $\Rel_G$ is a cyclically reduced word. Below is the procedure called \emph{Dehn's algorithm} (see~\cite{LS}). Generally speaking, Dehn's algorithm is a greedy algorithm on words based on the corresponding set of turns with respect to words length in generators. Let us explain it in more detail. Let $W$ be a word. Let $A = M_1M_2^{-1}$ be an element from $\Rel_G$ such that $M_1$ is an occurrence in $W$, $W = LM_1R$, and $\vert M_1\vert > \vert M_2\vert$, where $\vert \cdot\vert$ is the number of generators in a reduced word. The step of Dehn's algorithm is as follows:
\begin{equation*}
W = LM_1R \longmapsto L\cdot M_2\cdot R.
\end{equation*}
That is, we perform a turn of the occurrence $M_1$ to its complement~$M_2$. If there is no element of $\Rel_G$ with the required property, then the algorithm terminates.

The group $G = \langle \mathcal{X} \mid \Rel_G\rangle$ satisfies \emph{small cancellation condition $C^{\prime}(\frac{1}{6})$} if for every $R_{j_1} = cR_{j_1}^{\prime} \in \Rel_G$ and $R_{j_2} = cR_{j_2}^{\prime} \in \Rel_G$, we have $\vert c \vert < \frac{1}{6}\vert R_{j_1}\vert$ and $\vert c \vert < \frac{1}{6}\vert R_{j_2}\vert$ (see~\cite{LS}). The following theorem is one of main results of Small Cancellation Groups Theory (see~\cite{LS}).
\begin{theorem}
\label{dehns_algorithm_sc_groups}
Assume $G = \langle \mathcal{X} \mid \Rel_G\rangle$ satisfies condition $C^{\prime}(\frac{1}{6})$. Then a word $W$ is equal to $1$ in $G$ if and only if Dehn's algorithm, starting from $W$, terminates at $1$.
\end{theorem}

\begin{definition}
\label{greedy_alg_rings}
Based on order $<_f$ and on the set $\Add(\Rel)$, we define an algorithm with a black-box on $\Galg$. Every step of this algorithm is a reduction of the highest monomial of an element of $\Galg$ with respect to order $<_f$. So, in fact, we define  a greedy algorithm with respect to order $<_f$. We denote this procedure by $\GreedyAlg(<_f, \Add(\Rel))$. Namely, let $W_1, \ldots, W_k$ be different monomials, $\sum_{i = 1}^k \gamma_i W_i$, $\gamma_i \neq 0$, be an element of $\Galg$. Let $W_{i_0}$ be the highest monomial among $W_1, \ldots, W_k$ with respect to order $<_f$. Let $\sum_{j = 1}^s \alpha_j a_j$ be a polynomial from $\Add(\Rel)$ such that
\begin{enumerate}[label={(GA\arabic*)}]
\item
\label{reduct_polynomial1}
there exists an index $1 \leqslant h \leqslant s$ such that $a_h$ is a virtual member of the chart of $W_{i_0}$;
\item
\label{reduct_polynomial2}
$La_jR <_f La_hR$, where $W_{i_0} = La_hR$, for all $j \neq h$.
\end{enumerate}
Then the step of $\GreedyAlg( <_f, \Add(\Rel))$ is as follows:
\begin{align*}
&\sum_{i = 1}^k \gamma_i W_i = \sum_{\substack{i = 1 \\ i \neq i_0}}^s \gamma_i W_i + \gamma_{i_0} W_{i_0} \longmapsto \sum_{\substack{i = 1 \\ i \neq i_0}}^k \gamma_i W_i + \gamma_{i_0} \sum_{\substack{j = 1 \\ j \neq h}}^s (-\alpha_j\alpha_h^{-1} La_jR)\\
&\textit{ (and the further cancellations if there are any)}.
\end{align*}
That is, we perform a multi-turn of $W_{i_0}$, which comes from the elementary multi-turn $a_h \mapsto \sum_{\substack{j = 1 \\ j \neq h}}^s (-\alpha_j\alpha_h^{-1} a_j)$. If there is no polynomial in $\Add(\Rel)$ that satisfies conditions~\ref{reduct_polynomial1} and~\ref{reduct_polynomial2}, then the algorithm terminates.

A black-box answers whether there exits or does not exist a polynomial in $\Add(\Rel)$ that satisfies conditions~\ref{reduct_polynomial1} and~\ref{reduct_polynomial2}.

In other words, $\GreedyAlg( <_f, \Add(\Rel))$ for $\sum_{i = 1}^k \gamma_i W_i$ works in the following way. We take a polynomial $T$ from $\GDp^{\prime\prime}$ such that $W_{i_0}$ is its highest monomial. Let $\delta$ be the coefficient of $W_{i_0}$ in $T$. Then the step of $\GreedyAlg( <_f, \Add(\Rel))$ is as follows:
\begin{equation*}
\sum_{i = 1}^k \gamma_i W_i \longmapsto \sum_{i = 1}^k \gamma_i W_i - \gamma_{i_0}\delta^{-1}T.
\end{equation*}
If there is no element of $\GDp^{\prime\prime}$ such that $W_{i_0}$ is its highest monomial, then the algorithm terminates.

Notice that there may exist several elements of $\Add(\Rel)$ that satisfy conditions~\ref{reduct_polynomial1} and~\ref{reduct_polynomial2}. If this is the case, then we choose one of them arbitrarily. Hence, the algorithm $\GreedyAlg(<_f, \Add(\Rel))$ may have several possibilities at every step, so it may have different possible branches of execution. That is, $\GreedyAlg(<_f, \Add(\Rel))$ is a non-deterministic algorithm.

One can consider $\GreedyAlg(<_f, \Add(\Rel))$ as a generalization of Dehn's algorithm.
\end{definition}

In fact, inside Main Lemma (Lemma~\ref{fall_through_linear_dep_one}) we proved the following statement.
\begin{lemma}
\label{cancelling_of_highest_monomial}
Let $U$ be a monomial. Let $T_1, \ldots, T_n$ be elements of $\widetilde{\GDp}[U]$. Let $X$ be the biggest monomial of all highest monomials of $T_1, \ldots, T_n$ with respect to order~$<_f$. Assume $X$ additively cancels out in a linear combination $\sum_{j = 1}^n \delta_j T_j$, $\delta_j \neq 0$. Then
\begin{equation*}
\sum_{j = 1}^n \delta_j T_j = \sum_{r = 1}^{n^{\prime}} \delta_r^{\prime} T_r^{\prime},
\end{equation*}
where all $T_r^{\prime} \in \widetilde{\GDp}[U] \cup \GDp^{\prime}(\Low\DSpace{U})$, and all monomials of $T_1^{\prime}, \ldots, T_{n^{\prime}}^{\prime}$ are smaller than $X$ with respect to order~$<_f$.
\end{lemma}
\begin{proof}
See the proof of Lemma~\ref{fall_through_linear_dep_one}.
\end{proof}

\begin{lemma}
\label{sets_of_generators_conncetion}
Let $U$ be a monomial. Then for every element of $\widetilde{\GDp}[U]$ there exists an element of $\GDp^{\prime\prime}$ with the same highest monomial with respect to order~$<_f$.
\end{lemma}
\begin{proof}
Let $T \in \widetilde{\GDp}[U]$. Then we have $T \in \widetilde{\GDp}^{(i_0)}[U]$ for some $1 \leqslant i_0 \leqslant \nvirt{U}$. Hence,
\begin{equation*}
T = \mu[U]\left(a^{(1)}\otimes \ldots \otimes t\otimes\ldots \otimes a^{(m)}\right),
\end{equation*}
where $a^{(i)} \in ME_i[U]$ for $i \neq i_0$ and $t \in D_{i_0}[U]$. Since $t \in D_{i_0}[U]$, we have
\begin{equation*}
t = \sum_{j = 1}^k \delta_j t_j,
\end{equation*}
where $t_j \in \Rel \cap D_{i_0}[U]$. So,
\begin{equation}
\label{t_tilda_element_representation}
T = \sum_{j = 1}^k \delta_j T_j, \textit{ where } T_j = \mu[U]\left(a^{(1)}\otimes \ldots \otimes t_j\otimes\ldots\otimes a^{(m)}\right).
\end{equation}
Assume $\widetilde{t}_j$ corresponds to $t_j$ in $T_j$ (see the proof of statement~\ref{mu_properties_up_dep_image} of Lemma~\ref{mu_properties}, formula~\eqref{t_t_tilde_correspondnce}, for the details). Statement~\ref{mu_properties_up_dep_image} of Lemma~\ref{mu_properties} implies $\widetilde{t}_j \in \Rel$ and $T_j \in \GDp^{\prime}$. Since $\GDp^{\prime} \subseteq \GDp^{\prime\prime}$, we have $T_j \in \GDp^{\prime\prime}$.

Let us prove the following statement. Assume $T = \sum_{j = 1}^k \delta_j T_j$, where
\begin{itemize}
\item
$T_j = \mu[U]\left(a^{(1)}\otimes \ldots \otimes t_j\otimes\ldots\otimes a^{(m)}\right) = L\widetilde{t}_jR$, $\widetilde{t}_j$ corresponds to $t_j$ in $T_j$ (see the proof of statement~\ref{mu_properties_up_dep_image} of Lemma~\ref{mu_properties}, formula~\eqref{t_t_tilde_correspondnce}, for the details);
\item
$a^{(i)} \in ME_i[U]$ for $i \neq i_0$;
\item
$t_j \in \Add(\Rel) \cap D_{i_0}[U]$;
\item
$\widetilde{t}_j \in \Add(\Rel)$;
\item
$T_j \in \GDp^{\prime\prime}$.
\end{itemize}
Let $X$ be the highest monomial of $T$. Then there exists an element of $\GDp^{\prime\prime}$ with the same highest monomial $X$.

Looking at~\eqref{t_tilda_element_representation}, one can see that Lemma~\ref{sets_of_generators_conncetion} follows from the above statement. We prove it by induction on $k$.

Assume $k = 1$. Then $T = \delta_1T_1$. Therefore, $X$ is the highest monomial of $T_1 \in \GDp^{\prime\prime}$.

Assume $k > 1$. Let us make the step of induction. Let $x$ be the biggest monomial of the highest monomials of all $t_1, \ldots, t_k$. First assume that $x$ does not additively cancel out in $\sum_{j = 1}^k \delta_j t_j$. Obviously, $x$ is the biggest monomial of all monomials of $t_1, \ldots, t_k$. Therefore, by the definition of $<_f$,
\begin{equation*}
X = \mu[U]\left(a^{(1)}\otimes \ldots \otimes x\otimes\ldots\otimes a^{(m)}\right).
\end{equation*}
The monomial $x$ is the highest monomial of some $t_{j_0}$. Hence, $X$ is the highest monomial of $T_{j_0} = \mu[U]\left(a^{(1)}\otimes \ldots \otimes t_{j_0}\otimes\ldots\otimes a^{(m)}\right) \in \GDp^{\prime\prime}$.

Now assume that $x$ additively cancels out in $\sum_{j = 1}^k \delta_j t_j$. Without loss of generality, we can assume that $x$ is the highest monomial of $t_1, \ldots, t_{\widetilde{k}}$, and $t_{\widetilde{k} + 1},\ldots, t_k$ have smaller highest monomials. Then $x$ cancels out in $\sum_{j = 1}^{\widetilde{k}} \delta_j t_j$. Assume $\eta_j$ is the coefficient of $x$ in $t_j$. Then we have
\begin{multline*}
\sum_{j = 1}^{\widetilde{k}} \delta_j t_j = \delta_1\eta_1(\eta_1^{-1}t_1 - \eta_2^{-1}t_2) + (\delta_1\eta_1 + \delta_2\eta_2)(\eta_2^{-1}t_2 - \eta_3^{-1}t_3) + \ldots\\
\ldots + (\delta_1\eta_1 + \ldots + \delta_{\widetilde{k}}\eta_{\widetilde{k}})\eta_{\widetilde{k}}^{-1}t_{\widetilde{k}}.
\end{multline*}
Notice that $(\delta_1\eta_1 + \ldots + \delta_{\widetilde{k}}\eta_{\widetilde{k}})$ is the coefficient of $x$ in $\sum_{j = 1}^{\widetilde{k}} \delta_j t_j$. Therefore, $\delta_1\eta_1 + \ldots + \delta_{\widetilde{k}}\eta_{\widetilde{k}} = 0$. Let us put
\begin{align*}
s_{j, j + 1} = \eta_j^{-1}t_j - \eta_{j + 1}^{-1}t_{j + 1}.\\
\end{align*}
So, we have
\begin{equation*}
\sum_{j = 1}^{\widetilde{k}} \delta_j t_j = \delta_1\eta_1s_{1,2} + (\delta_1\eta_1 + \delta_2\eta_2)s_{2, 3} + \ldots + (\delta_1\eta_1 + \ldots + \delta_{\widetilde{k} - 1}\eta_{\widetilde{k} - 1})s_{\widetilde{k} - 1, \widetilde{k}}.
\end{equation*}

Consider $s_{j, j + 1}$ for $1 \leqslant j \leqslant \widetilde{k} - 1$ in detail. By Small Cancellation Axiom, every polynomial $t_j$ has a monomial of $\SPM$-measure $\geqslant \tau + 1$. Therefore, by the definition of $<_f$, we see that $\SPM(x) \geqslant \tau + 1$. So, since $x$ additively cancels out in $s_{j, j + 1}$, we get $s_{j, j + 1} \in \Add(\Rel)$.

We have
\begin{align*}
\eta_j^{-1}T_j - \eta_{j + 1}^{-1}T_{j + 1} &= \mu[U]\left(a^{(1)}\otimes\ldots \otimes(\eta_j^{-1}t_j - \eta_{j + 1}^{-1}t_{j + 1})\otimes \ldots \otimes a^{(m)} \right) =\\
&= \eta_j^{-1}L\widetilde{t}_jR - \eta_{j + 1}^{-1}L\widetilde{t}_{j + 1}R.
\end{align*}
Let us put
\begin{equation*}
\widetilde{s}_{j, j + 1} = \eta_j^{-1}\widetilde{t}_j - \eta_{j + 1}^{-1}\widetilde{t}_{j + 1}.
\end{equation*}
Then we have
\begin{equation*}
\eta_j^{-1}T_j - \eta_{j + 1}^{-1}T_{j + 1} = L\widetilde{s}_{j, j + 1}R.
\end{equation*}
Let $\widetilde{x}$ be a monomial that corresponds to $x$ in $\widetilde{t}_j$ and $\widetilde{t}_{j + 1}$. Then Lemma~\ref{same_small_changes} implies that $\widetilde{x}$ may differ from $x$ by at most a small piece at the beginning and by a small piece at the end. So, $\SPM(\widetilde{x}) \geqslant \SPM(x) - 2 \geqslant \tau - 1$. Therefore, since $\widetilde{t}_j, \widetilde{t}_{j + 1} \in \Add(\Rel)$, we obtain  $\widetilde{s}_{j, j + 1} \in \Add(\Rel)$. By Small Cancellation Axiom, $\widetilde{s}_{j, j + 1}$ contains a monomial $a$ of $\SPM$-measure $\geqslant \tau + 1$. Since $a$ is not a small piece and $T_j, T_{j + 1} \in \GDp^{\prime\prime}$, by Lemma~\ref{add_rel_incident_keep_structure}, $a$ is a maximal occurrence in $LaR$. Since $\SPM(a) \geqslant \tau$, $a$ is a virtual member of $LaR$. Therefore, $L\widetilde{s}_{j, j+1}R \in \GDp^{\prime\prime}$.

As a result, we obtain
\begin{equation}
\label{induction_step_sum}
T = \sum_{j = 1}^{\widetilde{k} - 1}\delta^{\prime}_j T_j^{\prime} + \sum_{j = \widetilde{k} + 1}^k \delta_jT_j,
\end{equation}
where $T_j^{\prime} = L\widetilde{s}_{j, j + 1}R$ and $\delta^{\prime}_j = \delta_1\eta_1 + \ldots + \delta_{j}\eta_{j}$. Moreover, we proved that
\begin{itemize}
\item
$T_j^{\prime} = \mu[U]\left(a^{(1)}\otimes \ldots \otimes s_{j, j + 1}\otimes\ldots\otimes a^{(m)}\right) = L\widetilde{s}_{j, j + 1}R$, $\widetilde{s}_{j, j + 1}$ corresponds to $s_{j, j + 1}$ in $T_j^{\prime}$;
\item
$s_{j, j + 1} \in \Add(\Rel) \cap D_{i_0}[U]$;
\item
$\widetilde{s}_{j, j + 1} \in \Add(\Rel)$;
\item
$T^{\prime}_j = L\widetilde{s}_{j, j+1}R \in \GDp^{\prime\prime}$.
\end{itemize}
Notice that sum~\eqref{induction_step_sum} contains $k - 1$ members. Thus, by the induction hypothesis, there exists an element of $\GDp^{\prime\prime}$ with the same highest monomial $X$ as $T$ has. This completes the proof.
\end{proof}

Using Lemma~\ref{cancelling_of_highest_monomial} and Lemma~\ref{sets_of_generators_conncetion}, we obtain the following statement.
\begin{lemma}
\label{possibility_of_reduction_step}
Assume $W_1, \ldots, W_k$ are different monomials, and an element $\sum_{i = 1}^k \gamma_i W_i$, $\gamma_i \neq 0$, belongs to $\Ideal$. Then it is possible to make a step of the algorithm $\GreedyAlg(<_f, \Add(\Rel))$ for $\sum_{i = 1}^k \gamma_i W_i$. Namely, let $W_{i_0}$ be the highest monomial among $W_1, \ldots, W_k$ with respect to order $<_f$. Then there exists an element of $\GDp^{\prime\prime}$ such that $W_{i_0}$ is its highest monomial.
\end{lemma}
\begin{proof}
Assume $\sum_{i = 1}^k \gamma_i W_i \in \Ft_N(\Galg) \setminus \Ft_{N - 1}(\Galg)$. Recall that $\Ideal = \Dp(\Galg)$. Then $\sum_{i = 1}^k \gamma_i W_i \in \Ft_N(\Galg) \cap \Dp(\Galg)$. It follows from Proposition~\ref{fall_to_smaller_subspace} that $\Ft_N(\Galg) \cap \Dp(\Galg) = \Dp(\Ft_{N}(\Galg))$. So, $\sum_{i = 1}^k \gamma_i W_i \in \Dp(\Ft_{N}(\Galg))$. Therefore,
\begin{equation*}
\sum_{i = 1}^k \gamma_i W_i = \sum_{j = 1}^n \delta_j T_j,
\end{equation*}
where $T_j \in \GDp^{\prime}(\Ft_N(\Galg))$.

Let $Z_j$ be a monomial such that $T_j$ is a layout of a multi-turn of a virtual member of the chart of $Z_j$, $j = 1, \ldots, n$. Then, obviously, $T_j \in \Dp\DSpace{Z_j}$. Consider the space $\sum_{j = 1}^n \DSpace{Z_j}$. If $\DSpace{Z_{j^{\prime}}} \subseteq \DSpace{Z_{j^{\prime\prime}}}$ for $j^{\prime} \neq j^{\prime\prime}$, then, clearly,
\begin{equation*}
\sum_{j = 1}^n \DSpace{Z_j} = \sum_{\substack{j = 1\\j \neq j^{\prime}}}^n \DSpace{Z_j}.
\end{equation*}
Hence, we can choose a subset $\lbrace Z_{j_1}, \ldots, Z_{j_t}\rbrace \subseteq \lbrace Z_1, \ldots, Z_n\rbrace$ such that
\begin{equation*}
\sum_{j = 1}^n \DSpace{Z_j} = \sum_{s = 1}^t \DSpace{Z_{j_s}} \textit{ and } \DSpace{Z_{j_s}} \nsubseteq \DSpace{Z_{j_{s^{\prime}}}} \textit{ whenever } s \neq s^{\prime}.
\end{equation*}
Then every $Z_j$, $j = 1, \ldots, n$, is a derived monomial of some $Z_{j_s}$, $s = 1, \ldots, t$.

Since $T_j \in \Dp\DSpace{Z_j}$, we have that $T_j$ belongs to some $\Dp\DSpace{Z_{j_s}}$, $s = 1, \ldots, t$. Therefore, we can separate $\lbrace T_1, \ldots, T_n\rbrace$ into $t$ groups as follows:
\begin{equation*}
\lbrace T_1, \ldots, T_n\rbrace = \bigsqcup_{s = 1}^t\lbrace T_1^{(s)}, \ldots, T_{n_s}^{(s)}\rbrace,
\end{equation*}
where $\lbrace T_1^{(s)}, \ldots, T_{n_s}^{(s)}\rbrace \subseteq \Dp\DSpace{Z_{j_s}}$. Hence,
\begin{equation*}
\sum\limits_{j = 1}^n \delta_j T_j = \sum\limits_{s = 1}^t \sum\limits_{p = 1}^{n_s} \delta_p^{(s)} T_p^{(s)}.
\end{equation*}

Applying consecutively Lemma~\ref{cancelling_of_highest_monomial}, we obtain
\begin{equation*}
\sum\limits_{p = 1}^{n_s} \delta_p^{(s)} T_p^{(s)} = \sum\limits_{q = 1}^{\widetilde{n}_s} \widetilde{\delta}_q^{(s)} \widetilde{T}_q^{(s)} + Q^{(s)},
\end{equation*}
where every $\widetilde{T}_q^{(s)} \in \widetilde{\GDp}[Z_{j_s}]$, $Q^{(s)} \in \DMLow{Z_{j_s}}$, and the biggest monomial of all highest monomials of $\widetilde{T}_1^{(s)}, \ldots, \widetilde{T}_{n_s}^{(s)}$ does not cancel. Since $\DMLow{Z_{j_s}} \subseteq \Ft_{N - 1}(\Galg)$, we have $Q^{(s)} \in \Ft_{N - 1}(\Galg)$. Since $\sum\limits_{j = 1}^n \delta_j T_j \notin \Ft_{N - 1}(\Galg)$, we obtain that at least one sum $\sum\limits_{q = 1}^{\widetilde{n}_s} \widetilde{\delta}_q^{(s)} \widetilde{T}_q^{(s)}$ is non-empty.

Lemma~\ref{EDp_bigger_set_of_generators} implies  $\widetilde{\GDp}[Z_{j_s}] \cap \Low\DSpace{Z_{j_s}} = \varnothing$. Therefore, every $\widetilde{T}_q^{(s)} \in \DSpace{Z_{j_s}}$ and $\widetilde{T}_q^{(s)}  \notin \Low\DSpace{Z_{j_s}}$. So, it follows from the definition of $<_f$ that the highest monomial of $\widetilde{T}_q^{(s)}$ is contained in $\DMUp{Z_{j_s}}$. Since $\DSpace{Z_{j_s}} \nsubseteq \DSpace{Z_{j_{s^{\prime}}}} $ whenever $s \neq s^{\prime}$,  it follows from Lemma~\ref{derived_spaces_equality} that
\begin{equation*}
\DMUp{Z_{j_s}} \cap \DSpace{Z_{j_{s^{\prime}}}} = \varnothing \textit{ if } s \neq s^{\prime}.
\end{equation*}
Therefore, the highest monomials of $\lbrace T_1^{(s)}, \ldots, T_{n_s}^{(s)}\rbrace$ for different $s$ can not cancel each other out. Hence, there exists $\widetilde{T}_{q_0}^{(s_0)}$ such that its highest monomial is equal to $W_{i_0}$. Then it follows from Lemma~\ref{sets_of_generators_conncetion} that there exists an element of $\GDp^{\prime\prime}$ such that its highest monomial is equal to $W_{i_0}$. This completes the proof.
\end{proof}

As a result, we have the following theorem.
\begin{theorem}
\label{possibility_of_reduction}
Assume $W_1, \ldots, W_k$ are different monomials. We take an element $\sum_{i = 1}^k \gamma_i W_i \in \Galg$, $\gamma_i \neq 0$. Then the following statements are equivalent:
\begin{enumerate}[label={(\arabic*)}]
\item
\label{possibility_of_reduction1}
some branch of the algorithm $\GreedyAlg(<_f, \Add(\Rel))$, starting from $\sum_{i = 1}^k \gamma_i W_i$, terminates at~$0$;
\item
\label{possibility_of_reduction2}
$\sum_{i = 1}^k \gamma_i W_i \in \Ideal$;
\item
\label{possibility_of_reduction3}
every branch of the algorithm $\GreedyAlg(<_f, \Add(\Rel))$, starting from $\sum_{i = 1}^k \gamma_i W_i$, terminates at~$0$;
\end{enumerate}
\end{theorem}
\begin{proof}
$\ref{possibility_of_reduction1} \Rightarrow \ref{possibility_of_reduction2}$. Assume some branch of the algorithm $\GreedyAlg(<_f, \Add(\Rel))$, starting from $\sum_{i = 1}^k \gamma_i W_i$, terminates at $0$. The $n$-th step of the branch of $\GreedyAlg(<_f, \Add(\Rel))$ is of the form
\begin{equation*}
\sum_{i = 1}^{k^{(n)}} \gamma_{i}^{(n)}W_i^{(n)} \longmapsto \sum_{i = 1}^{k^{(n)}} W_i^{(n)} - \delta^{(n)} T^{(n)},
\end{equation*}
where $\sum_{i = 1}^{k^{(1)}} \gamma_{i}^{(1)}W_i^{(1)} = \sum_{i = 1}^k \gamma_i W_i$, $\sum_{i = 1}^{k^{(n)}} \gamma_{i}^{(n)}W_i^{(n)}$ is the result of the step $n - 1$ for $n > 1$, and $T^{(n)} \in \GDp^{\prime\prime}$. Since $\GDp^{\prime\prime} \subseteq \Ideal$, we have $T^{(n)} \in \Ideal$. Hence, the result of the $n$-th step is of the form
\begin{equation*}
\sum_{i = 1}^k \gamma_i W_i - \sum_{r = 1}^n \delta^{(r)}T^{(r)},
\end{equation*}
where $T^{(1)}, \ldots, T^{(n)} \in \Ideal$. The branch of the algorithm $\GreedyAlg(<_f, \Add(\Rel))$, starting from $\sum_{i = 1}^k \gamma_i W_i$, terminates at $0$ after a step with some number $N$. Therefore, we obtain
\begin{equation*}
\sum_{i = 1}^k \gamma_i W_i  = \sum_{r = 1}^N \delta^{(r)}T^{(r)}.
\end{equation*}
Thus, $\sum_{i = 1}^k \gamma_i W_i \in \Ideal$.
\medskip

$\ref{possibility_of_reduction2} \Rightarrow \ref{possibility_of_reduction3}$.
Assume $\sum_{i = 1}^k \gamma_i W_i \in \Ideal$. We need to prove that every branch of the algorithm $\GreedyAlg(<_f, \Add(\Rel))$, starting from $\sum_{i = 1}^k \gamma_i W_i$, terminates at~$0$. Let us take an arbitrary branch of $\GreedyAlg(<_f, \Add(\Rel))$, starting from $\sum_{i = 1}^k \gamma_i W_i$, and show that it terminates at~$0$.

Assume the contrary, namely, assume that the branch of $\GreedyAlg(<_f, \Add(\Rel))$, starting from $\sum_{i = 1}^k \gamma_i W_i$, terminates at some non-zero element of~$\Galg$. Assume $\sum_{i = 1}^{k^{\prime}} \gamma_i^{\prime} W_i^{\prime} \in \Galg$, where $W_i^{\prime}$ are monomials, $\gamma_i^{\prime} \neq 0$, is this element. We proved in implication $\ref{possibility_of_reduction1} \Rightarrow \ref{possibility_of_reduction2}$ that a result of the $n$-th step of $\GreedyAlg(<_f, \Add(\Rel))$, starting from $\sum_{i = 1}^k \gamma_i W_i$, is of the form
\begin{equation*}
\sum_{i = 1}^k \gamma_i W_i - \sum_{r = 1}^n \delta^{(r)}T^{(r)},
\end{equation*}
where $T^{(1)}, \ldots, T^{(n)}$ are some elements of $\GDp^{\prime\prime} \subseteq \Ideal$. Therefore, since $\sum_{i = 1}^k \gamma_i W_i \in \Ideal$, we obtain
\begin{equation*}
\sum_{i = 1}^k \gamma_i W_i - \sum_{r = 1}^n \delta^{(r)}T^{(r)} \in \Ideal.
\end{equation*}
In particular, $\sum_{i = 1}^{k^{\prime}} \gamma_i^{\prime} W_i^{\prime} \in \Ideal$. Let $W_{i_0}^{\prime}$ be the highest monomial of $\sum_{i = 1}^{k^{\prime}} \gamma_i^{\prime} W_i^{\prime}$ with respect to order $<_f$. Then it follows from Lemma~\ref{possibility_of_reduction_step} that there exists $T \in \GDp^{\prime\prime}$ such that $W_{i_0}^{\prime}$ is the highest monomial of $T$. Let $\delta^{\prime}$ be the coefficient of $W_{i_0}^{\prime}$ in $T$. Hence, it is possible to do a step of $\GreedyAlg(<_f, \Add(\Rel))$ for $\sum_{i = 1}^{k^{\prime}} \gamma_i^{\prime} W_i^{\prime}$ of the form
\begin{equation*}
\sum_{i = 1}^{k^{\prime}} \gamma_i^{\prime} W_i^{\prime} \longmapsto \sum_{i = 1}^{k^{\prime}} \gamma_i^{\prime} W_i^{\prime} - \gamma_{i_0}^{\prime}{\delta^{\prime}}^{-1}T.
\end{equation*}
Therefore, $\GreedyAlg(<_f, \Add(\Rel))$  does not terminate at $\sum_{i = 1}^{k^{\prime}} \gamma_i^{\prime} W_i^{\prime}$. A contradiction. Thus, an arbitrary branch of $\GreedyAlg(<_f, \Add(\Rel))$, starting from $\sum_{i = 1}^k \gamma_i W_i$, terminates at~$0$.
\medskip

$\ref{possibility_of_reduction3} \Rightarrow \ref{possibility_of_reduction1}$. This implication is trivial.
\end{proof}

Theorem~\ref{possibility_of_reduction} implies the following statemement.
\begin{corollary}
\label{grobner_basis_ideal_membership}
The set $\Add(\Rel)$ is a Gr\"{o}bner basis of the ideal $\Ideal$ with respect to monomial ordering $<_f$, and $\GreedyAlg(<_f, \Add(\Rel))$ solves the Ideal Membership Problem for~$\Ideal$.
\end{corollary}

Let us show the following property of additively closed system of generators $\Rel$. We will use it in Section~\ref{examples_section}.
\begin{lemma}
\label{add_closed_incident_chain_of_incidend}
Assume $\Add(\Rel) = \Rel$. Let $m_1, \ldots, m_k$ be a sequence of monomials from $\Mon$ such that $m_i$, $m_{i + 1}$ for $i = 1, \ldots, k - 1$ are incident monomials and $\SPM(m_i) \geqslant \tau - 2$ for $i = 2, \ldots, k - 1$. Then $m_1$ and $m_k$ are incident monomials. In particular, the notions of incident monomials and $U$-incident monomials ($U$ is a monomial) coincide.
\end{lemma}
\begin{proof}
We prove Lemma~\ref{add_closed_incident_chain_of_incidend} by induction on $k$.

Assume $k = 2$, then $m_1$ and $m_2$ are incident monomials by the initial assumption.

Assume $k > 2$. Let us make a step of induction. Consider the first three monomials $m_1, m_2$ and $m_3$. First assume that $m_1 = m_3$. Then we obtain the sequence of monomials $m_1, m_4, \ldots, m_k$ such that consecutive monomials are incident, $\SPM(m_i) \geqslant \tau - 2$ for $i = 4, \ldots, k - 1$, and length of this sequence is equal to $k - 2$. So, $m_1$ and $m_k$ are incident monomials, by the induction hypothesis.

Now assume that $m_1 \neq m_3$. Since $m_1$ and $m_2$ are incident monomials, there exists a polynomial $p \in \Rel$ of the form
\begin{equation*}
p = \alpha_1m_1 + \alpha_2m_2 + \sum_{i = 3}^{n_1}\alpha_i a_i,
\end{equation*}
and exists a polynomial $q \in \Rel$ of the form
\begin{equation*}
q = \beta_1m_2 + \beta_2m_3 + \sum_{i = 3}^{n_2}\beta_i b_i.
\end{equation*}
Since $\Add(\Rel) = \Rel$ and $\SPM(m_2) \geqslant \tau - 2$, we get  $g = \alpha_2^{-1}p - \beta_2^{-1}q \in \Rel$. Obviously, $g$ contains monomials $m_1$ and $m_3$. So, $m_1$ and $m_3$ are incident monomials. Therefore, we obtain the sequence of monomials $m_1, m_3, \ldots, m_k$ such that consecutive monomials are incident, $\SPM(m_i) \geqslant \tau - 2$ for $i = 3, \ldots, k - 1$, and length of this sequence is equal to $k - 1$. Hence, $m_1$ and $m_k$ are incident monomials, by the induction hypothesis.
\end{proof}

\section{Examples}
\label{examples_section}
\subsection{Group algebras of small cancellation groups}
\label{group_algebras_of_sc_groups_section}
Let a group $G$ be given by generators and defining relations,
\begin{equation*}
G = \langle \mathcal X \mid \Rel_G\rangle,\ \Rel_G = \lbrace R_j\rbrace_{j \in J},
\end{equation*}
where $\Rel_G$ is closed under taking cyclic shifts of relators and taking inverses of relators, and every $R_j$ is a cyclically reduced word. The word $c$ is called \emph{a small piece with respect to $\Rel_G$ (in a group sense)} if there exist $R_{j_1}, R_{j_2} \in \Rel_G$ such that $R_{j_1} = cR_{j_1}^{\prime}$, $R_{j_2} = cR_{j_2}^{\prime}$ and $R_{j_1}^{\prime} \neq R_{j_2}^{\prime}$ as words in the corresponding free group. Obviously, every subword of a small piece is a small piece as well.

Since $\Rel_G$ is closed under taking inverses of relators, we get $R_{j_1}^{-1} = {R_{j_1}^{\prime}}^{-1}c^{-1} \in \Rel_G$ and $R_{j_2}^{-1} = {R_{j_2}^{\prime}}^{-1}c^{-1} \in \Rel_G$. Since $\Rel_G$ is closed under taking cyclic shifts of relators, we see that $c^{-1}{R_{j_1}^{\prime}}^{-1} \in \Rel_G$ and $c^{-1}{R_{j_2}^{\prime}}^{-1} \in \Rel_G$. Therefore, if $c$ is a small piece, then $c^{-1}$ is a small piece as well.

Let $\SP_{G}(\Rel_G)$ be a set of all small pieces with respect to $\Rel_G$ in the group sense.

\begin{remark}
Assume $c$ is a subword of some $R_j \in \Rel_G$ and $c$ is not a small piece. The following properties follow directly from the above definition and the fact that $\Rel_G$ is closed under taking cyclic shifts of relators.
\begin{itemize}
\item
If $c$ is a prefix of $R_{j_1} \in \Rel_G$ and $c$ is a prefix of $R_{j_2} \in \Rel_G$, then $R_{j_1} = R_{j_2}$.
\item
There are no occurrences of $c$ in other relators from $\Rel_G$ different from cyclic shifts of~$R_j$.
\end{itemize}
\end{remark}

\begin{remark}
Notice that even if $c$ is not a small piece and $c$ is a prefix of some $R_{j_1} \in \Rel_G$, $c$ can be a prefix of some $R_{j_2}^{n}$, where $R_{j_2} \in \Rel_G$, $R_{j_1} \neq R_{j_2}$, for $n \in \mathbb{Z}$, $\vert n\vert \geqslant 2$. If this happens, we get that ether $R_{j_2}$, or $R_{j_2}^{-1}$ is a proper prefix of $R_{j_1}$. Therefore, $R_{j_2}$ and $R_{_2}^{-1}$ are small pieces.
\end{remark}

\begin{lemma}
\label{occurrences_of_non_sp_in_powers}
Assume $c$ is a subword of some $R_j \in \Rel_G$ and $c$ is not a small piece. Assume there is more than one occurrence of $c$ in $R_j^m$, $m \in \mathbb{N}$. Then all the occurrences of $c$ are shifts of the first occurrence by a multiple of length of the smallest (with respect to length in the generators $\mathcal{X} \cup \mathcal{X}^{-1}$) period of~$R_j$.
\end{lemma}
\begin{proof}
First assume that there is a single occurrence of $c$ in $R_j$. Then all occurrences of $c$ in $R_j^m$ are shifts of $c$ by length of $R_j$. Otherwise, we can shift by a multiple of length of $R_j$ an occurrence of $c$ in $R_j^m$ that violates this property and obtain more than one occurrence of $c$ in $R_j$.

Now assume that there is more than one occurrence of $c$ in $R_j$. Let  $R_j = l_1cr_1 = l_2cr_2$, where $l_1$ is a proper prefix of $l_2$. Let $l_2 = l_1l^{\prime}$. Then $l^{\prime}$ is a shift between these two occurrences of $c$ in $R_j$. Let us show that $\vert l^{\prime}\vert$ is a multiple of length of the smallest period of $R_j$.

Consider a cyclic shift of $R_j$ by $l_1$ and denote it by $R_{j^{\prime}}$. That is, we have $R_{j^{\prime}} = cr_1l_1$. Since $\Rel_G$ is closed under cyclic shifts, we get  $cr_1l_1\in \Rel_G$ and $cr_2l_2 \in \Rel_G$. Since $c$ is not a small piece, we obtain that $cr_1l_1 = cr_2l_2$ as words in the corresponding free group. Therefore, $R_{j^{\prime}} = cr_2l_2$.

Since $R_j = l_2cr_2 = l_1r^{\prime}cR_2$, we have $R_{j^{\prime}} = l^{\prime}cr_2l_1$. We have $cr_2l_2 = cr_2l_1l^{\prime}$. Therefore, $cr_2l_2$ is a cyclic shift of $R_{j^{\prime}}$ by $l^{\prime}$. However, we proved above that $R_{j^{\prime}} = cr_2l_2$. Hence, $R_{j^{\prime}}$ is equal to its cyclic shift by $l^{\prime}$. Since we work with words in a free group, this implies that there exists $a \in \Fr$ such that $l^{\prime} = a^{n_1}$ and $cr_2l_1 = a^{n_2}$, $n_1, n_2 \in \mathbb{N}$.

Let $a$ be the shortest word with the property $l^{\prime} = a^{n_1}$ and $cr_2l_1 = a^{n_2}$, $n_1, n_2 \in \mathbb{N}$. Let us show that $a$ is the smallest period of $R_{j^{\prime}}$. Assume the contrary. Namely, we assume that there exists $b \in \Fr$ such that $b$ is a proper prefix of $a$ and $R_{j^{\prime}} = b^t$. Then, on the one hand, $R_{j^{\prime}} = a^{n}$, $n = n_1 + n_2 \geqslant 2$, on the other hand, $R_{j^{\prime}} = b^t$. One can show that if two power words have a common subword of length greater than the sum of their periods, then their periods are powers of the same word. Since $n \geqslant 2$ and $\vert a\vert > \vert b\vert$, we obtain that $a^{n}$ and $b^{t}$ have a common subword $a^n$ of length greater then $\vert a\vert + \vert b\vert$. Hence, $a$ and $b$ are powers of the same word. This contradicts the assumption that $a$ is the shortest word with the property $L^{\prime} = a^{n_1}$ and $cR_2L_1 = a^{n_2}$. So, $a$ is the smallest period of $R_{j^{\prime}}$.

Length of the smallest period of $R_{j}$ is equal to length of the smallest period of $R_{j^{\prime}}$, because $R_{j^{\prime}}$ is a cyclic shift of $R_j$. Therefore, length of $L^{\prime}$ in generators is a multiple of the length of the smallest period of $R_j$.

Clearly, the same property holds for positions of occurrences of $c$ in $R_j^m$ for $m \geqslant 2$. Lemma~\ref{occurrences_of_non_sp_in_powers} is proved.
\end{proof}

\begin{remark}
Assume $c$ is a subword of some $R_j \in \Rel_G$ and $c$ is not a small piece.  Lemma~\ref{occurrences_of_non_sp_in_powers} implies the following properties.
\begin{itemize}
\item
If $R_j$ is a primitive word, then there is a unique occurrence of $c$ in $R_j$.
\item
If $R_j$ is a proper power and $c$ is a subword of its smallest period, then there is a unique occurrence of $c$ in the smallest period of $R_j$.
\end{itemize}
\end{remark}

\begin{lemma}
\label{unique_relation_for_non_small_piece}
Assume a group $G = \left\langle \mathcal{X} \mid \Rel_G\right\rangle$. Let $R_{j_1}, R_{j_2} \in \Rel_G$. Assume $c$ is not a small piece and $c$ is a prefix of $R_{j_1}^{k_1}$ and $R_{j_2}^{k_2}$, where $k_1, k_2 \in \mathbb{N}$. Then either $R_{j_1}$ is a proper prefix of $R_{j_2}$, or $R_{j_2}$ is a proper prefix of $R_{j_1}$, or $R_{j_1} = R_{j_2}$.
\end{lemma}
\begin{proof}
Since $c$ is a prefix of $R_{j_1}^{k_1}$, we see that either $c$ is a prefix of $R_{j_1}$, or $R_{j_1}$ is a prefix of $c$. Similarly, since $c$ is a prefix of $R_{j_2}^{k_2}$, we see that either $c$ is a prefix of $R_{j_2}$, or $R_{j_2}$ is a prefix of $c$.

Assume that $c$ is a prefix of $R_{j_1}$ and $c$ is a prefix of $R_{j_1}$. So, $R_{j_1} = cR_{j_1}^{\prime}$ and $R_{j_2} = cR_{j_2}^{\prime}$. Then, by definition, $R_{j_1}^{\prime}= R_{j_2}^{\prime}$, because $c$ is not a small piece. Therefore, $R_{j_1} = R_{j_2}$.

Assume that $c$ is a prefix of $R_{j_1}$ and $R_{j_2}$ is a prefix of $c$. Then we obtain that $R_{j_2}$ is a prefix of $R_{j_1}$. That is, either $R_{j_2}= R_{j_1}$, or $R_{j_2}$ is a proper prefix of $R_{j_1}$. Assume $R_{j_1}$ is a prefix of $c$ and $c$ is a prefix of $R_{j_2}$. Then we obtain that $R_{j_1}$ is a prefix of $R_{j_2}$. That is, either $R_{j_1}= R_{j_2}$, or $R_{j_1}$ is a proper prefix of $R_{j_2}$.

Assume that $R_{j_1}$ is a prefix of $c$ and $R_{j_2}$ is a prefix of $c$. Then, clearly, either $R_{j_1} = R_{j_2}$, or $R_{j_1}$ is a proper prefix of $R_{j_2}$, or $R_{j_2}$ is a proper prefix of $R_{j_1}$. This completes the proof.
\end{proof}

\begin{definition}
\label{sc_group}
We say that a group $G = \langle \mathcal X \mid \Rel_G\rangle$ satisfies \emph{small cancellation condition~$C(m)$}\label{cond_cm_2} if every $R_j \in \Rel_G$ can not be written as a product of less than $m$ small pieces.
\end{definition}

We noticed above that if $c$ is a small piece, then $c^{-1}$ is a small piece as well. Therefore, a word $a$ can not be written as a product of less than $t$ small pieces if and only if $a^{-1}$ can not be written as a product of less than $t$ small pieces for every $t \in \mathbb{N}$.
\medskip

Assume $G = \langle \mathcal X \mid \Rel_G\rangle$, where $\Rel_G = \lbrace R_j\rbrace_{j \in J}$, $\Rel_G$ is closed under taking cyclic shifts of relators and taking inverses of relators, and every $R_j$ is a cyclically reduced word. For now, we do not make any additional assumptions about the group $G$. We consider its group algebra $\fld G$, where $\fld$ is a field. Let $\Fr$ be the free group with the set of free generators~$\mathcal{X}$. Consider an ideal $\Ideal$ of $\Galg$ that is generated as an ideal by the set $\lbrace  R_j - 1\rbrace_{j \in J}$. That is,
\begin{equation}
\label{ideal_group_algebra}
\Ideal = \left\langle \lbrace  R_j - 1\rbrace_{j \in J}\right\rangle_i.
\end{equation}

\begin{lemma}
\label{group_algebra_as_a_quotient}
We have $\fld G \cong \Qalg$.
\end{lemma}
\begin{proof}
By the universal property of $\Fr$ we obtain that there exists the canonical surjective group homomorphism $\phi_1: \Fr \to G$ such that $\ker \phi_1$ is a normal subgroup generated as a normal subgroup by $\Rel_G$. Clearly, we can linearly extend $\phi_1$ to $\Galg$ and obtain a linear mapping of vector spaces $\phi_2: \Galg \to \fld G$. Using the well-known definition of $\phi_1$, one can easily see that $\phi_2$ preserves the multiplication as well, so $\phi_2$ is a homomorphism of algebras.

By the isomorphism theorem, we obtain that $\fld G \cong \Galg / \ker\phi_2$. We need to prove that $\Ideal = \ker\phi_2$. Let $R_j \in \Rel_G$. We have
\begin{equation*}
\phi_2(R_j - 1) = \phi_1(R_j) - \phi_1(1) = 1 - 1 = 0.
\end{equation*}
Therefore, $\lbrace  R_j - 1\rbrace_{j \in J} \subseteq \ker\phi_2$. Since $\ker\phi_2 \triangleleft \Galg$, this implies  $\Ideal \subseteq \ker\phi_2$.

Let us show that $\ker\phi_2 \subseteq \Ideal$. Assume $\sum_{i = 1}^n\alpha_i A_i \in \ker\phi_2$, where $A_i \in \Fr$. Then
\begin{equation}
\label{kernel_element}
0 = \phi_2\left(\sum_{i = 1}^n\alpha_i A_i\right) = \sum_{i = 1}^n\alpha_i \phi_1(A_i).
\end{equation}
Since the elements of $G$ is a basis of $\fld G$, we see that all elements in the right-hand side of the above sum have to additively cancel. We split the monomials~$A_i$ into sets such that the elements of the same set have equal images under the mapping $\phi_1$ and the elements from different sets have different images. Let $A_{i_1}, \ldots, A_{i_t}$ be such a set. Since $\phi_1(A_{i_1}) = \ldots = \phi_1(A_{i_t})$ additively cancel with each other in~\eqref{kernel_element}, we see that $\sum_{s = 1}^{t}\alpha_{i_s} = 0$. Hence,
\begin{equation*}
\phi_{2}\left(\sum_{s = 1}^{t}\alpha_{i_s}A_{i_s}\right) = \sum_{s = 1}^{t}\alpha_{i_s}\phi_1(A_{i_s}) = \left(\sum_{s = 1}^{t}\alpha_{i_s} \right)\phi_1(A_{i_1}) = 0.
\end{equation*}
Therefore,
\begin{equation*}
\sum_{s = 1}^{t}\alpha_{i_s}A_{i_s} \in \ker\phi_2.
\end{equation*}
Thus, it is enough to consider elements $\sum_{i = 1}^n\alpha_i A_i\in \ker\phi_2$ such that all $A_i$ have the same image under the mapping $\phi_1$ and show that such elements belong to $\Ideal$.

So, we take $\sum_{i = 1}^n\alpha_i A_i \in \ker\phi_2$ such that all $A_i$ have the same image under the mapping $\phi_1$. Recall that, since
\begin{equation*}
0 = \phi_2\left(\sum_{i = 1}^n\alpha_i A_i \right) = \sum_{i = 1}^n\alpha_i \phi_1(A_i) = \left(\sum_{i = 1}^n\alpha_i\right)\phi_1(A_1),
\end{equation*}
we have $\sum_{i = 1}^n\alpha_i = 0$ in this case. Since all $A_i$ have the same image under the mapping $\phi_1$, we obtain that all $A_i$ belong to the same coset of $\ker\phi_1$. That is, every $A_i = A\cdot H_i$, where $A \in \Fr$, $H_i \in \ker\phi_1$. Hence, it is enough to prove that $\sum_{i = 1}^n\alpha_i H_i \in \Ideal$. Since $\ker\phi_1$ is generated as a normal subgroup by $\Rel_G = \lbrace R_j\rbrace_{j \in J}$, by definition, we have
\begin{equation*}
H_i = \prod_{s = 1}^{n_i}U_s^{(i)}\cdot R_{s}^{(i)}\cdot {U_s^{(i)}}^{-1}, \textit{ where } R_{s}^{(i)} \in \Rel_G,\ U_s^{(i)} \in \Fr.
\end{equation*}
Since $R_{s}^{(i)} - 1 \in \Ideal$, we have
\begin{multline*}
U_s^{(i)}\cdot R_{s}^{(i)}\cdot {U_s^{(i)}}^{-1} = U_s^{(i)}\cdot (R_{s}^{(i)} - 1 + 1)\cdot {U_s^{(i)}}^{-1} =\\
= U_s^{(i)}\cdot (R_{s}^{(i)} - 1)\cdot {U_s^{(i)}}^{-1} + U_s^{(i)}\cdot {U_s^{(i)}}^{-1} = 1 \mod \Ideal.
\end{multline*}
Therefore, every $H_i = 1 \mod \Ideal$. Hence,
\begin{equation*}
\sum_{i = 1}^n\alpha_i H_i = \sum_{i = 1}^n\alpha_i \cdot 1 \mod \Ideal  = 0.
\end{equation*}
That is, $\sum_{i = 1}^n\alpha_i H_i \in \Ideal$. Therefore,
\begin{equation*}
\sum_{i = 1}^n\alpha_i A_i = \sum_{i = 1}^n\alpha_i A\cdot H_i = A\cdot \left(\sum_{i = 1}^n\alpha_i H_i\right) \in \Ideal.
\end{equation*}
This completes the proof.
\end{proof}

In what follows we assume that the group $G = \left\langle \mathcal{X} \mid \Rel_G\right\rangle$ satisfies small cancellation condition $C(m)$. For now, we do not specify a value of the constant $m$. Let us study the quotient algebra $\Qalg$, where the ideal $\Ideal$ is defined by~\eqref{ideal_group_algebra}. Using the set of relators $\lbrace R_j - 1\rbrace_{j \in J}$, we will extend it and construct a set of relators $\Rel \supseteq \lbrace R_j - 1\rbrace_{j \in J}$ with the following properties:
\begin{enumerate}[label={(\roman*)}]
\item
\label{rel_prop1}
$\Rel$ generates as an ideal the same ideal $\Ideal$;
\item
\label{rel_prop2}
$\Rel$ satisfies Compatibility Axiom, Small Cancellation Axiom, and Isolation Axiom (both left-sided and right-sided);
\item
\label{rel_prop3}
$\SP_G(\Rel_G) = \SP_R(\Rel)$, where $\SP_R(\Rel)$ is a set of all small pieces with respect to $\Rel$ in the sense of Definition~\ref{sp}.
\end{enumerate}

~\paragraph*{$1^{\mathrm{o}}$} Let us look first at the set $\lbrace R_j - 1\rbrace_{j \in J}$. One can easily see that it does not satisfy Compatibility Axiom. We consider the Cayley graph of the group $G$ with respect to the set of generators $\mathcal{X}$. Then every $R_j \in \Rel_G$ corresponds to a closed path in the Cayley graph.
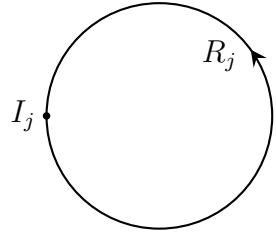
\begin{equation}
\label{relator_graph}
\begin{tikzpicture}
\draw[black, thick, bigarrow=0.1] (0, 0) circle (1.5);
\node at (0.8, 0.8) {$R_j$};
\node[circle, fill, inner sep=1] at (-1.5, 0) {};
\node[left] at (-1.5, 0) {$I_j$};
\end{tikzpicture}
\end{equation}
Here the point $I_j$ is the initial and the final point of the path that corresponds to $R_j$. Every cyclic shift of $R_j$ corresponds to a closed path in graph~\eqref{relator_graph} with some other initial and final point.

Let us do the first obvious step. Assume $R_j = a_j^{(s)}{b_j^{(s)}}^{-1}$. Then, clearly $a_j^{(s)} = b^{(s)}_j$ in $G$ and in $\Qalg$. We consider all binomial relations in $G$ and in $\Qalg$ of such a form. That is, we take all different points $F_j^{(s)}$ on graph~\eqref{relator_graph} and consider two different simple paths from $I_j$ to~$F_j^{(s)}$.
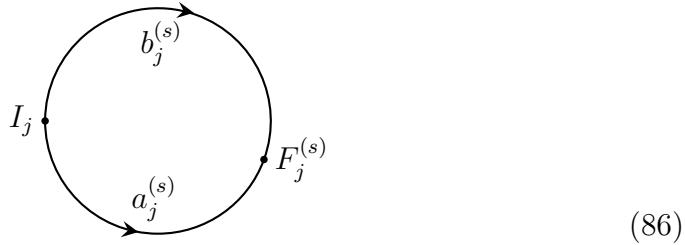
\begin{equation}
\label{relator_graph_two_points}
\begin{tikzpicture}
\node[circle, fill, inner sep=1] at (-1.5, 0) {};
\node[left] at (-1.5, 0) {$I_j$};
\draw[thick, bigarrow=0.5] (-1.5, 0) arc (-180:-20:1.5) node[midway, above, xshift=6] {$a_j^{(s)}$} coordinate (F);
\node[circle, fill, inner sep=1] at (F) {};
\node[right] at (F) {$F_j^{(s)}$};
\draw[thick, reversebigarrow=0.5] (F) arc (-20:180:1.5) node[midway, below, xshift=-6] {$b_j^{(s)}$};
\end{tikzpicture}
\end{equation}
We take the set of binomials
\begin{equation*}
\left\lbrace \gamma(a_j^{(s)} - b_j^{(s)}) \mid R_j \in \Rel_G,\  R_j = a_j^{(s)}{b_j^{(s)}}^{-1}, \gamma \in \fld\right\rbrace
\end{equation*}
and add them to the initial set $\lbrace R_j - 1\rbrace_{j \in J}$. Denote the obtained set by $\Rel_1$. Obviously, $\Rel_1 \subseteq \Ideal$, so, $\Rel_1$ generates as an ideal the same ideal $\Ideal$.

Since $\Rel_G$ is closed under taking cyclic shifts, one can see that Compatibility Axiom holds for the obtained set $\Rel_1$. However, if we consider $\SP_R(\Rel_1)$, a set of all small pieces with respect to $\Rel_1$ in the sense of Definition~\ref{sp}, we see that there are too many small pieces and $\Rel_1$ does not satisfy Small Cancellation Axiom. Namely, assume $R_j = a_jb_j \in \Rel_G$, both $a_j$ and $b_j$ are not empty. Since $\Rel_G$ is closed under cyclic shifts, we see that $b_ja_j \in \Rel_G$. Therefore, $a_jb_j - 1 \in \Rel \subseteq \Rel_1$ and $b_ja_j- 1 \in \Rel \subseteq \Rel_1$. Notice that monomials in every binomial from $\Rel_1$ do not have common prefixes and common suffixes. Hence, $b_ja_jb_j - b_j \notin \Rel_1$. Therefore, by Definition~\ref{sp}, $a_j$ is a small piece with respect to $\Rel_1$. Similarly, $b_j$ is a small piece with respect to $\Rel_1$. That is, every proper subword of every $R_j \in \Rel_G$ is a small piece with respect to $\Rel_1$. So, every $R_j \in \Rel_G$ is a product of not more than two small pieces with respect to $\Rel_1$ (regardless the constant $m$ in condition $C(m)$ for $G = \langle \mathcal{X} \mid \Rel_G\rangle$).

In order to deal with the above difficulties, we do a further extension of the set $\Rel_1$. Let $R_j \in \Rel_G$. As above, we consider the corresponding graph~\eqref{relator_graph} and take all different points $F_j^{(s)}$ in this graph. But now we consider all possible different paths with the initial point $I_j$ and the final point $F_j^{(s)}$, not only two different simple paths as we did above. Let $R_j = a_j^{(s)}{b_j^{(s)}}^{-1}$ (see~\eqref{relator_graph_two_points}). Then every such path corresponds to a monomial either of the form $R_j^na_j^{(s)}$, or of the form $R_j^{-n}b_j^{(s)}$, where $n \in \lbrace 0\rbrace\cup \mathbb{N}$. Clearly, all such monomials are equal in $G$ and in $\Qalg$. We take the set of binomials
\begin{align*}
\left\lbrace \gamma(c_j^{(s)} - d_j^{(s)}) \right.\mid &R_j \in \Rel_G, \ R_j = a_j^{(s)}{b_j^{(s)}}^{-1}, \gamma \in \fld,\\
&\left.c_j^{(s)}, d_j^{(s)} \in \lbrace R_j^na_j^{(s)}, R_j^{-n}b_j^{(s)} \mid n \in \lbrace 0\rbrace\cup \mathbb{N}\rbrace\right\rbrace
\end{align*}
and add them to the initial set $\lbrace R_j - 1\rbrace_{j \in J}$. We denote the obtained set by~$\Rel$. That is,
\begin{align}
\begin{split}
\label{group_algebra_relatios_set}
\Rel = \lbrace R_j - 1\rbrace_{j \in J} \bigcup \Big\lbrace &\gamma(c_j^{(s)} - d_j^{(s)})\mid R_j \in \Rel_G, \ R_j = a_j^{(s)}{b_j^{(s)}}^{-1}, \gamma \in \fld,\\
&c_j^{(s)}, d_j^{(s)} \in \lbrace R_j^na_j^{(s)}, R_j^{-n}b_j^{(s)} \mid n \in \lbrace 0\rbrace\cup \mathbb{N}\rbrace\Big\rbrace.
\end{split}
\end{align}
Let us show that $\Rel$ satisfies properties~\ref{rel_prop1}---\ref{rel_prop3}.

Obviously, $\Rel$ generates as an ideal the same ideal $\Ideal$. Since $\Rel_G$ is closed under taking cyclic shifts, we see that $\Rel$ satisfies Compatibility Axiom. Let us show that the rest of the properties are satisfied as well.

~\paragraph*{$2^{\mathrm{o}}$} Let us study $\SP_R(\Rel)$, a set of all small pieces with respect to $\Rel$ in the sense of Definition~\ref{sp}, and show that $\Rel$ satisfies Small Cancellation Axiom with the constant $\tau = \left[\frac{m}{2}\right] - 1$ ($\left[\frac{m}{2}\right]$ is the integer part of $\frac{m}{2}$).

In the following proposition we assume that $G = \left\langle \mathcal{X} \mid \Rel_G\right\rangle$ satisfies the following condition: for every two different $R_{j_1}, R_{j_2} \in \Rel_G$ neither $R_{j_1}$ is a subword of $R_{j_2}$, nor $R_{j_2}$ is a subword of~$R_{j_1}$. Notice that this condition is equivalent to small cancellation condition~$C(2)$.
\begin{proposition}
\label{sc_group_same_small_pieces}
Assume a group $G = \left\langle \mathcal{X} \mid \Rel_G\right\rangle$, where $\Rel_G$ is closed under taking cyclic shifts and  inverses of relators, and every $R_j \in \Rel_G$ is a cyclically reduced word. For every two different $R_{j_1}, R_{j_2} \in \Rel_G$ we assume that neither $R_{j_1}$ is a subword of $R_{j_2}$, nor $R_{j_2}$ is a subword of~$R_{j_1}$. Let $\Rel$ be defined by~\eqref{group_algebra_relatios_set}, $\SP_R(\Rel)$ be a set of all small pieces with respect to $\Rel$ in the sense of Definition~\ref{sp}. Then $\SP_G(\Rel_G) = \SP_R(\Rel)$.
\end{proposition}
\begin{proof}
Assume $c \in \Fr$, $c \notin \SP_G(\Rel_G)$. Let us show that $c \notin \SP_R(\Rel)$. Assume that $c$ is an occurrence in monomials in two polynomials $T_1, T_2 \in \Rel$. By the definition of $\Rel$, there exist $R_{j_1}, R_{j_2} \in \Rel_G$ such that
\begin{align*}
&T_1 = \gamma_1(c_1 - d_1),\ T_2 = \gamma_2(c_2 - d_2),\textit{ where}\\
&R_{j_1} = a_{j_1}b_{j_1}^{-1}, \ R_{j_2} = a_{j_2}b_{j_2}^{-1},\\
&c_1, d_1 \in \lbrace R_{j_1}^na_{j_1}, R_{j_1}^{-n}b_{j_1} \mid n \in \lbrace 0\rbrace\cup \mathbb{N}\rbrace,\\
&c_2, d_2 \in \lbrace R_{j_2}^na_{j_2}, R_{j_2}^{-n}b_{j_2} \mid n \in \lbrace 0\rbrace\cup \mathbb{N}\rbrace,\\
&\gamma_1, \gamma_2 \in \fld.
\end{align*}
Without loss of generality, we assume that $c$ is an occurrence in $c_1$ and in $c_2$. Let $c_1 = u_1cv_1$ and  $c_2 = u_2cv_2$. Then
\begin{align*}
T_1 = \gamma_1(c_1 - d_1) = \gamma_1(u_1cv_1 - d_1),\\
T_2 = \gamma_2(c_2 - d_2) = \gamma_2(u_2cv_2 - d_2).
\end{align*}
In order to show that $c \notin \SP_R(\Rel)$, we need to prove that
\begin{align*}
u_2\cdot u_1^{-1}\cdot T_1 = \gamma_1(u_2cv_1 - u_2\cdot u_1^{-1}\cdot d_1) \in \Rel,\\
T_1\cdot v_1^{-1}\cdot v_2 = \gamma_1(u_1cv_2 - d_1\cdot v_1^{-1}\cdot v_2) \in \Rel.
\end{align*}

Assume that $c$ is a proper subword of $R_{j_1}$ or of a cyclic shift of $R_{j_1}$. Then Lemma~\ref{occurrences_of_non_sp_in_powers} implies that positions of occurrences of $c$ in $R_{j_1}^m$, $m \in \mathbb{N}$, differ by a multiple of length of the smallest period of $R_{j_1}$. Assume $c$ contains $R_{j_1}$ or some cyclic shift of $R_{j_1}$. Then, similarly to Lemma~\ref{occurrences_of_non_sp_in_powers}, one can show that positions of occurrences of $c$ in $R_{j_1}^m$, $m \in \mathbb{N}$, also differ by a multiple of length of the smallest period of $R_{j_1}$. By the same reason, positions of occurrences of $c$ in $R_{j_2}^m$, $m \in \mathbb{N}$, differ by a multiple of length of the smallest period of~$R_{j_2}$.

Since $c$ is an occurrence in $c_1$, we get that $c$ is a prefix of $\widetilde{R}_{j_1}^{k_1}$, where $\widetilde{R}_{j_1}$ is a cyclic shift of $R_{j_1}$ or $R_{j_1}^{-1}$, $k_1 \in \mathbb{N}$. Similarly, since $c$ is an occurrence in $c_2$, we get that $c$ is a prefix of $\widetilde{R}_{j_2}^{k_2}$, where $\widetilde{R}_{j_2}$ is a cyclic shift of $R_{j_2}$ or $R_{j_2}^{-1}$, $k_2 \in \mathbb{N}$. Since $\Rel_G$ is closed under taking inverses and cyclic shifts of relators, we obtain  $\widetilde{R}_{j_1}, \widetilde{R}_{j_2} \in \Rel_G$. Since $c \notin \SP_G(\Rel_G)$, Lemma~\ref{unique_relation_for_non_small_piece} implies that either $\widetilde{R}_{j_1}$ is a proper prefix of $\widetilde{R}_{j_2}$, or $\widetilde{R}_{j_2}$ is a proper prefix of $\widetilde{R}_{j_1}$, or $\widetilde{R}_{j_2} = \widetilde{R}_{j_2}$. It follows from our initial assumptions that neither $\widetilde{R}_{j_1}$ is a proper subword of $\widetilde{R}_{j_2}$, nor $\widetilde{R}_{j_2}$ is a proper subword of $\widetilde{R}_{j_1}$. Hence, we get  $\widetilde{R}_{j_2} = \widetilde{R}_{j_2}$. Therefore, either $R_{j_2}$ is a cyclic shift of $R_{j_1}$, or $R_{j_2}^{-1}$ is a cyclic shift of $R_{j_1}$.

Consider a graph of the form~\eqref{relator_graph_two_points} for $T_1$.
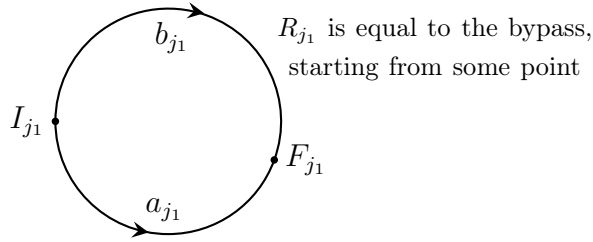
\begin{center}
\begin{tikzpicture}
\node[circle, fill, inner sep=1] at (-1.5, 0) {};
\node[left] at (-1.5, 0) {$I_{j_1}$};
\draw[thick, bigarrow=0.5] (-1.5, 0) arc (-180:-20:1.5) node[midway, above, xshift=6] {$a_{j_1}$} coordinate (F);
\node[circle, fill, inner sep=1] at (F) {};
\node[right] at (F) {$F_{j_1}$};
\draw[thick, reversebigarrow=0.5] (F) arc (-20:180:1.5) node[midway, below, xshift=-6] {$b_{j_1}$} coordinate[pos = 0.3] (R);
\node[text width=4.5cm, right, align=center] at (R) {\footnotesize $R_{j_1}$ is equal to the bypass, starting from some point};
\end{tikzpicture}
\end{center}
It follows from the latter observation that $c_2$, $d_2$ correspond to paths in the same graph with some other initial point $I_{j_1}$ and final point $F_{j_2}$. We proved above that the positions of all occurrences of $c$ in $R_{j_1}^m$, $m \in \mathbb{N}$, differ by a multiple of length of the smallest period of $R_{j_1}$. Hence, we can choose the points $I_{j_2}$ and $F_{j_2}$ such that the  paths in the above graph that correspond to $c_1 = u_1cv_1$ and $c_2 = u_2cv_2$ have a common subpath that corresponds to the occurrences of $c$ under consideration.
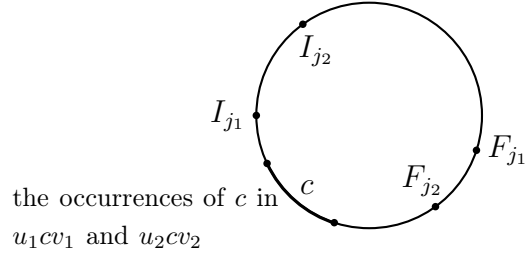
\begin{center}
\begin{tikzpicture}
\draw[thick] (1.5, 0) arc (0:360:1.5) coordinate[pos=0.5] (I1) coordinate[pos=0.95] (F1) coordinate[pos=0.35] (I2) coordinate[pos=0.85] (F2) coordinate[pos=0.57] (C1) coordinate[pos=0.7] (C2) coordinate[pos=0.625] (CM) ;
\node[circle, fill, inner sep=1] at (I1) {};
\node[left] at (I1) {$I_{j_1}$};
\node[circle, fill, inner sep=1] at (F1) {};
\node[right] at (F1) {$F_{j_1}$};
\node[circle, fill, inner sep=1] at (I2) {};
\node[below, xshift=5] at (I2) {$I_{j_2}$};
\node[circle, fill, inner sep=1] at (F2) {};
\node[above, xshift=-5] at (F2) {$F_{j_2}$};
\node[circle, fill, inner sep=1] at (C1) {};
\node[circle, fill, inner sep=1] at (C2) {};
\draw[very thick] (C1) arc (180 + 25: 180 + 73:1.5);
\node[above, right, yshift=3] at (CM) {$c$};
\node[below, xshift=-45, yshift=7, text width=4.2cm] at (CM) {\footnotesize the occurrences of $c$ in $u_1cv_1$ and $u_2cv_2$};
\end{tikzpicture}
\end{center}
Therefore, the monomials $u_2cv_1$ and $u_2\cdot u_1^{-1}\cdot d_1$ correspond to paths in the above graph with the initial point $I_{j_2}$ and the final point $F_{j_1}$. Similarly, the monomials $u_1cv_2$ and $d_2\cdot v_1^{-1}\cdot v_2$ correspond to paths with the initial point $I_{j_1}$ and the final point $F_{j_2}$. So, by definition~\eqref{group_algebra_relatios_set},
\begin{align*}
u_2\cdot u_1^{-1}\cdot T_1 = \gamma_1(u_2cv_1 - u_2\cdot u_1^{-1}\cdot d_1) \in \Rel,\\
T_1\cdot v_1^{-1}\cdot v_2 = \gamma_1(u_1cv_2 - d_1\cdot v_1^{-1}\cdot v_2) \in \Rel.
\end{align*}
Thus, $c \notin \SP_R(\Rel)$. So, we obtain  $\SP_G(\Rel_G) \supseteq \SP_R(\Rel)$.

Let us show that $\SP_G(\Rel_G) \subseteq \SP_R(\Rel)$. Assume $c \in \SP_G(\Rel_G)$. We need to prove that $c \in \SP_R(\Rel)$. If $c = 1$, then $c$ belongs to $\SP_R(\Rel)$, by definition. In what follows we assume that $c \neq 1$.

Since $c \in \SP_G(\Rel_G)$, there exist $R_{j_1}, R_{j_2} \in \Rel_G$ such that $R_{j_1} = cR_{j_1}^{\prime}$, $R_{j_2} = cR_{j_2}^{\prime}$, and $R_{j_1}^{\prime} \neq R_{j_2}^{\prime}$ as words in the corresponding free group. We can assume that $c$ is a maximal common prefix of $R_{j_1}$ and $R_{j_2}$, since $\SP_R(\Rel)$ is closed under taking subwords.

Assume $c \notin \SP_R(\Rel)$. By definition~\eqref{group_algebra_relatios_set} of $\Rel$, we have
\begin{align*}
R_{j_1} - 1 = cR_{j_1}^{\prime} - 1 \in \Rel,\\
R_{j_2} - 1 = cR_{j_2}^{\prime} - 1 \in \Rel.
\end{align*}
Since $c \notin \SP_R(\Rel)$, we obtain
\begin{equation*}
(cR_{j_1}^{\prime} - 1)\cdot {R_{j_1}^{\prime}}^{-1}\cdot R_{j_2}^{\prime} = cR_{j_2}^{\prime} - {R_{j_1}^{\prime}}^{-1}R_{j_2}^{\prime} = R_{j_2} - {R_{j_1}^{\prime}}^{-1}R_{j_2}^{\prime} \in \Rel.
\end{equation*}
By definition~\eqref{group_algebra_relatios_set} of $\Rel$, there exists $R_{j_3} \in \Rel_G$ such that
\begin{equation*}
R_{j_3} = a_{j_3}b_{j_3}^{-1} \textit{ and }R_{j_2}, {R_{j_1}^{\prime}}^{-1}R_{j_2}^{\prime} \in \lbrace R_{j_3}^na_{j_3}, R_{j_3}^{-n}b_{j_3} \mid n \in \lbrace 0\rbrace\cup \mathbb{N}\rbrace.
\end{equation*}
It follows from our initial assumptions that $R_{j_3}^{\pm 1}$ are not proper prefixes of $R_{j_2}$. Therefore, $R_{j_2}$ can not be of the form $R_{j_3}^na_{j_3}, R_{j_3}^{-n}b_{j_3}$ for $n > 0$. Similarly, $R_{j_2}$ is not a proper prefix of $R_{j_3}^{\pm 1}$. Hence, $R_{j_2}$ can not be of the form $a_{j_3}, b_{j_3}$ if at least one of $a_{j_3}$ and $b_{j_3}$ is not empty. Therefore, either $a_{j_3}$ is empty, or $b_{j_3}$ is empty.

Assume $a_{j_3}$ is empty. Then $b_{j_3} = R_{j_3}^{-1}$. Therefore,
\begin{equation*}
R_{j_2} \in \lbrace R_{j_3}^n, R_{j_3}^{-n} \mid n \in \lbrace 0\rbrace\cup \mathbb{N}\rbrace.
\end{equation*}
Hence, we get $R_{j_2} = R_{j_3}^{\pm m}$, $m \in \mathbb{N}$. Since $R_{j_3}^{\pm 1}$ are not proper prefixes of $R_{j_2}$, this implies  $R_{j_2} = R_{j_3}^{\pm 1}$. Similarly, if $b_{j_3}$ is empty, $R_{j_2} = R_{j_3}^{\pm 1}$.

Since $R_{j_2} = R_{j_3}^{\pm 1}$ and one of $a_{j_3}$ or $b_{j_3}$ is empty, we have
\begin{equation*}
{R_{j_1}^{\prime}}^{-1}R_{j_2}^{\prime} \in \lbrace R_{j_2}^n, R_{j_2}^{-n} \mid n \in \lbrace 0\rbrace\cup \mathbb{N}\rbrace.
\end{equation*}
Hence, ${R_{j_1}^{\prime}}^{-1}R_{j_2}^{\prime} = R_{j_2}^{\pm m}$, $m \in \mathbb{N}$. Notice that since neither $R_{j_1} = cR_{j_1}^{\prime}$ is a subword of $R_{j_2} = cR_{j_2}^{\prime}$, nor $R_{j_2}$ is a subword of $R_{j_1}$, we get that both $R_{j_1}^{\prime}$ and $R_{j_2}^{\prime}$ are not empty.

First assume that ${R_{j_1}^{\prime}}^{-1}R_{j_2}^{\prime} = R_{j_2}^m$. Since $c$ is a prefix of $R_{j_2}$, we get that $c$ and ${R_{j_1}^{\prime}}^{-1}$ have a common prefix. However, this is not possible, because $R_{j_1} = cR_{j_1}^{\prime}$ is a cyclically reduced word. Now assume that ${R_{j_1}^{\prime}}^{-1}R_{j_2}^{\prime} = R_{j_2}^{-m}$. Then $R_{j_2}^{\prime}$ and $R_{j_2}^{-1} = {R_{j_2}^{\prime}}^{-1}c^{-1}$ have a common suffix. This is not possible, because $R_{j_2} = cR_{j_2}^{\prime}$ is a cyclically reduced word. A contradiction. Therefore,
\begin{equation*}
(cR_{j_1}^{\prime} - 1)\cdot {R_{j_1}^{\prime}}^{-1}\cdot R_{j_2}^{\prime} = cR_{j_2}^{\prime} - {R_{j_1}^{\prime}}^{-1}R_{j_2}^{\prime} = R_{j_2} - {R_{j_1}^{\prime}}^{-1}R_{j_2}^{\prime} \notin \Rel.
\end{equation*}
Thus, $c \in \SP_R(\Rel)$. Proposition~\ref{sc_group_same_small_pieces} is proved.
\end{proof}
So, in what follows in the current section we do not need to distinguish small pieces in the group sense and in the ring sense.

\begin{proposition}
\label{small_cancellation_in_ring}
Assume a group $G = \left\langle \mathcal{X} \mid \Rel_G\right\rangle$, where $\Rel_G$ is closed under taking cyclic shifts and inverses of relators, and every $R_j \in \Rel_G$ is a cyclically reduced word. Assume the group $G$ satisfies small cancellation condition $C(m)$, $m \geqslant 2$. Then the set $\Rel$ defined by~\eqref{group_algebra_relatios_set} satisfies Small Cancellation Axiom with the constant $\tau = \left[\frac{m}{2}\right] - 1$.
\end{proposition}
\begin{proof}
Let $T \in \Rel$. Then there exists $R_j \in \Rel_G$ such that
\begin{align}
\begin{split}
\label{binomial_relation_def}
T = \gamma(c_j - d_j), &\textit{ where }\ R_j = a_j{b_j}^{-1}, \gamma \in \fld,\\
&c_j, d_j \in \left\lbrace R_j^na_j, R_j^{-n}b_j \mid n \in \lbrace 0\rbrace\cup \mathbb{N}\right\rbrace.
\end{split}
\end{align}
Recall that there exists a graph of the form~\eqref{relator_graph_two_points} such that the monomials $c_j$ and $d_j$ correspond to two different paths in this graph with the same initial points and the same final points.
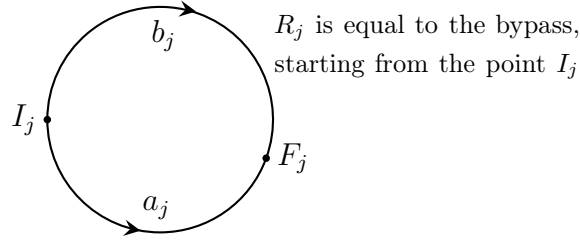
\begin{center}
\begin{tikzpicture}
\node[circle, fill, inner sep=1] at (-1.5, 0) {};
\node[left] at (-1.5, 0) {$I_j$};
\draw[thick, bigarrow=0.5] (-1.5, 0) arc (-180:-20:1.5) node[midway, above, xshift=6] {$a_j$} coordinate (F);
\node[circle, fill, inner sep=1] at (F) {};
\node[right] at (F) {$F_j$};
\draw[thick, reversebigarrow=0.5] (F) arc (-20:180:1.5) node[midway, below, xshift=-6] {$b_j$} coordinate[pos = 0.3] (R);
\node[text width=4.5cm, right, align=center] at (R) {\footnotesize $R_{j}$ is equal to the bypass, starting from the point $I_j$};
\end{tikzpicture}
\end{center}
Condition $C(m)$ implies that $R_j$ can not be written as a product of less than $m$ small pieces from  $\SP_G(\Rel_G)$. Therefore, at least one of $a_j$ and $b_j$ can not be written as a product of less than $\left[\frac{m}{2}\right]$ small pieces from  $\SP_G(\Rel_G)$. Hence, at least one of $c_j$ and $d_j$ can not be written as a product of less than $\left[\frac{m}{2}\right]$ small pieces from  $\SP_G(\Rel_G)$. It follows from Proposition~\ref{sc_group_same_small_pieces} that $\SP_G(\Rel_G) \supseteq \SP_R(\Rel)$. Hence, at least one of $c_j$ and $d_j$ can not be written as a product of less than $\left[\frac{m}{2}\right]$ small pieces from $\SP_R(\Rel)$. Therefore, Small Cancellation Axiom with the constant $\tau = \left[\frac{m}{2}\right] - 1$ holds for every single polynomial of $\Rel$.

Now consider an arbitrary linear combination $\sum_{r = 1}^n \gamma_rT_r$, where $T_r \in \Rel$, $\gamma_r \in \fld$. By definition~\eqref{group_algebra_relatios_set} of $\Rel$, we see that for every $T_r$ there exists $R_{j_r} \in \Rel_G$ and a fragmentation $R_{j_r} = a_{j_r}{b_{j_r}}^{-1}$ with property~\eqref{binomial_relation_def}. Let us collect all $T_r$ such that they correspond to the same $R_{j_r} \in \Rel_G$ up to taking inverse and the same fragmentation of $R_{j_r}$. Without loss of generality we can assume that the obtained sets are
\begin{equation}
\label{same_reletion_same_points}
\left\lbrace T_1, \ldots, T_{n_1} \right\rbrace, \ldots, \left\lbrace T_{n_t + 1}, \ldots, T_{n} \right\rbrace.
\end{equation}

Let us consider every set from~\eqref{same_reletion_same_points} and the corresponding linear combination separately. Without loss of generality we can consider the first set $\left\lbrace T_1, \ldots, T_{n_1}\right\rbrace$ and the corresponding linear combination $\sum_{r = 1}^{n_1} \gamma_rT_r$. So, all monomials of $T_1, \ldots, T_{n_1}$ belong to the set
\begin{equation*}
\left\lbrace R_{j_1}^na_{j_1}^{(s_1)}, R_{j_1}^{-n}b_{j_1}^{(s_1)} \mid n \in \lbrace 0\rbrace\cup \mathbb{N}\right\rbrace, \textit{ where } R_{j_1} \in \Rel_ G,\ R_{j_1} = a_{j_1}^{(s_1)}{b_{j_1}^{(s_1)}}^{-1}.
\end{equation*}
Since all $T_r$ are binomials, we obtain that either $\sum_{r = 1}^{n_1} \gamma_rT_r = 0$, or $\sum_{r = 1}^{n_1} \gamma_rT_r$ after additive cancellations contains at least two different monomials. Therefore, in the same way as above we obtain that at least one of these monomials can not be written as a product of less than $\left[\frac{m}{2}\right]$ small pieces from $\SP_R(\Rel)$. So, Small Cancellation Axiom with the constant $\tau = \left[\frac{m}{2}\right] - 1$ holds for every linear combination that corresponds to a single set from~\eqref{same_reletion_same_points}.

Let $T_{r_1}, T_{r_2}$ be two different polynomials from $\Rel$. Let $a$ be a monomial of $T_1$ and of $T_2$, and $a \notin \SP_R(\Rel)$. There exists $R_{t_1} \in \Rel_G$ that satisfies property~\eqref{binomial_relation_def} for $T_1$. There exists $R_{t_2} \in \Rel_G$ that satisfies property~\eqref{binomial_relation_def} for $T_2$. Therefore,
\begin{align}
\begin{split}
\label{common_not_sp}
a \in &\left\lbrace R_{t_1}^na_{t_1}^{(s_1)}, R_{t_1}^{-n}b_{t_1}^{(s_1)} \mid n \in \lbrace 0\rbrace\cup \mathbb{N}\right\rbrace \bigcap \left\lbrace R_{t_2}^na_{t_2}^{(s_2)}, R_{t_2}^{-n}b_{t_2}^{(s_2)} \mid n \in \lbrace 0\rbrace\cup \mathbb{N}\right\rbrace,\\
&\textit{where } R_{t_1} = a_{t_1}^{(s_1)}{b_{t_1}^{(s_1)}}^{-1},\ R_{t_2} = a_{t_2}^{(s_2)}{b_{t_2}^{(s_2)}}^{-1}.
\end{split}
\end{align}
In particular, either $a$ is a prefix of $R_{t_1}^{m_1}$, or $a$ is a prefix of $R_{t_1}^{-m_1}$ for some $m_1 \in \mathbb{N}$. Similarly, either $a$ is a prefix of $R_{t_2}^{m_2}$, or $a$ is a prefix of $R_{t_2}^{-m_2}$ for some $m_2 \in \mathbb{N}$. Recall that, by Proposition~\ref{sc_group_same_small_pieces}, we have $\SP_R(\Rel) = \SP_G(\Rel_G)$. Since $a \notin \SP_R(\Rel) = \SP_G(\Rel_G)$ and $G$ satisfies condition $C(2)$, it follows from Lemma~\ref{unique_relation_for_non_small_piece} that $R_{t_1} = R_{t_2}^{\pm 1}$.

If $R_{t_1} = R_{t_2}$, then it follows from~\eqref{common_not_sp} that $a_{t_1}^{(s_1)} = a_{t_2}^{(s_2)}$ and $b_{t_1}^{(s_1)} = b_{t_2}^{(s_2)}$. If $R_{t_1} = R_{t_2}^{-1}$, then it follows from~\eqref{common_not_sp} that $a_{t_1}^{(s_1)} = b_{t_2}^{(s_2)}$ and $b_{t_1}^{(s_1)} = a_{t_2}^{(s_2)}$. Therefore, $T_{r_1}, T_{r_2}$ belong to the same set from~\eqref{same_reletion_same_points}. So, non-small pieces that belong to different sets from~\eqref{same_reletion_same_points} are not equal to each other. Hence, if non-small pieces cancel in the linear combination $\sum_{r = 1}^n \gamma_rT_r$, then they cancel inside the linear combination $\sum_{r = n_i + 1}^{n_{i + 1}}\gamma_r T_r$ that corresponds to a single set from~\eqref{same_reletion_same_points}. Combining this with the above, we obtain that Small Cancellation Axiom with the constant $\tau = \left[\frac{m}{2}\right] - 1$ holds for~$\Rel$.
\end{proof}

From the very beginning for the argument in the presented paper we require the constant $\tau \geqslant 10$ in Small Cancellation Axiom. So, Proposition~\ref{small_cancellation_in_ring} implies that we obtain the required value of the constant $\tau$ for the group algebra $\fld G \cong \Qalg$ if $G = \left\langle \mathcal{X} \mid \Rel_G\right\rangle$ satisfies small cancellation condition $C(m)$ with $m \geqslant 22$.
\begin{corollary}
\label{small_cancellation_in_ring_min_const}
Assume a group $G = \left\langle \mathcal{X} \mid \Rel_G\right\rangle$, where $\Rel_G$ is closed under taking cyclic shifts and inverses of relators, and every $R_j \in \Rel_G$ is a cyclically reduced word. Assume the group $G$ satisfies small cancellation condition $C(m)$ with $m \geqslant 22$. Then the set $\Rel$ defined by~\eqref{group_algebra_relatios_set} satisfies Small Cancellation Axiom with the constant $\tau \geqslant 10$.
\end{corollary}

~\paragraph*{$3^{\mathrm{o}}$} As above, assume a group $G = \left\langle \mathcal{X} \mid \Rel_G\right\rangle$, where $\Rel_G$ is closed under taking cyclic shifts and inverses of relators, and every $R_j \in \Rel_G$ is a cyclically reduced word. Assume the group $G$ satisfies small cancellation condition $C(m)$, $m \geqslant 2$. We will prove that both Isolation Axioms (left-sided and right-sided) are satisfied for $\Rel$ defined by~\eqref{group_algebra_relatios_set}.

First notice that $\Add(\Rel) = \Rel$ (see the beginning of Section~\ref{algorithm_section} for the definition of $\Add(\Rel)$). Indeed, assume $T_1, T_2 \in \Rel$, $T_1 = \gamma_1(a - b_1)$, $T_1 = \gamma_2(a - b_2)$, and $\SPM(a) \geqslant \tau - 2$. Let us show that
\begin{equation*}
\gamma_1^{-1}T_1 - \gamma_2^{-1}T_2 = a - b_1 - (a - b_2) = b_2 - b_1 \in \Rel.
\end{equation*}
There exists $R_{j_1} \in \Rel_G$ and its fragmentation with property~\eqref{binomial_relation_def} for $T_1$. There exists $R_{j_2} \in \Rel_G$ and its fragmentation with property~\eqref{binomial_relation_def} for $T_2$. Since $a$ is not a small piece, we can argue in the same way as at the end of the proof of Proposition~\ref{small_cancellation_in_ring} and obtain that $R_{j_1} = R_{j_2}^{\pm 1}$ and the corresponding fragmentations are the same. Therefore, the monomials $a, b_1, b_2$ correspond to paths in the graph for $R_{j_1}$ of the form~\eqref{relator_graph} with the same initial point and the same final points. Thus, by definition, $b_2 - b_1 \in \Rel$.

Recall that, by definition, $\Mon$ is a set of monomials of polynomials from~$\Rel$. So, in the current case $\Mon$ consists of all subwords of powers of all relators from $\Rel_G$.

Let $m_1, m_2, \ldots, m_k$ be a sequence of monomials of $\Mon$ such that all of them are of $\SPM$-measure $\geqslant \tau - 2$ and the consecutive monomials are incident. Since $\Add(\Rel) = \Rel$, it follows from Lemma~\ref{add_closed_incident_chain_of_incidend} that $m_1$ and $m_k$ are incident monomials. So, we need to check Isolation Axiom only for incident monomials.

Let us show that right-sided  Isolation Axiom is satisfied for $\Rel$. Let $m_1, m_2 \in \Mon$ be two incident monomials, $m_1 \neq m_2$, $\SPM(m_1) \geqslant \tau - 2$, $\SPM(m_2) \geqslant \tau - 2$. We will prove even a slightly stronger condition. Namely, we take $a \in \Mon$ such that $\SPM(a) \geqslant \tau - 2$, $m_1a, m_2a \notin \Mon$, $m_1$ is a maximal occurrence in $m_1a$, $m_2$ is a maximal occurrence in $m_2a$ (that is, we omit the last condition on $a$ in right-sided  Isolation Axiom). Let $s_1a$ be a maximal occurrence in $m_1a$ that contains $a$, and $s_2a$ be a maximal occurrence in $m_2a$ that contains $a$. We will show that
\begin{equation*}
m_1\cdot s_1^{-1} \neq m_2 \cdot s_2^{-1}.
\end{equation*}

Assume the contrary, namely, assume that $m_1\cdot s_1^{-1} = m_2 \cdot s_2^{-1} = m$. Then we have $m_1 = ms_1 $, $m_2 = ms_2$. Since $s_1$ and $s_2$ are overlaps of maximal occurrences, $s_1$ and $s_2$ are small pieces. So, $m$ is not empty, because $\SPM(m_1) \geqslant \tau - 2 > 1$ and $\SPM(m_2) \geqslant \tau - 2 > 1$. Since $m_1$ and $m_2$ are two incident monomials, they are monomials of some polynomial $T \in \Rel$. So, by definition, there exists $R_j \in \Rel_G$ such that
\begin{equation*}
m_1, m_2 \in \left\lbrace R_{j}^na_{j}, R_{j}^{-n}b_{j} \mid n \in \lbrace 0\rbrace\cup \mathbb{N}\right\rbrace, \textit{ where } R_{j} = a_{j}b_{j}^{-1}.
\end{equation*}
So, the monomials $m_1$ and $m_2$ correspond to two paths in the graph of the form~\eqref{relator_graph} for $R_j$ with the same initial points and the same final points. Denote their initial point by $I_j$ and their final point by $F_j$.
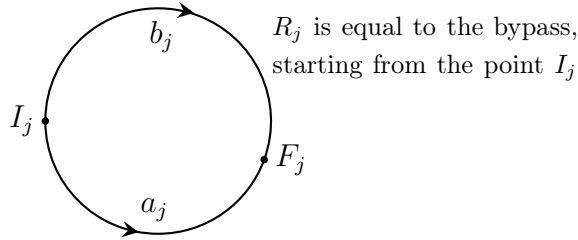
\begin{center}
\begin{tikzpicture}
\node[circle, fill, inner sep=1] at (-1.5, 0) {};
\node[left] at (-1.5, 0) {$I_j$};
\draw[thick, bigarrow=0.5] (-1.5, 0) arc (-180:-20:1.5) node[midway, above, xshift=6] {$a_j$} coordinate (F);
\node[circle, fill, inner sep=1] at (F) {};
\node[right] at (F) {$F_j$};
\draw[thick, reversebigarrow=0.5] (F) arc (-20:180:1.5) node[midway, below, xshift=-6] {$b_j$} coordinate[pos = 0.3] (R);
\node[text width=4.5cm, right, align=center] at (R) {\footnotesize $R_{j}$ is equal to the bypass, starting from the point $I_j$};
\end{tikzpicture}
\end{center}

Assume $m_1$ is of the form $R_j^{n}a_j$ and $m_2$ is of the form $R_j^{-n}b_j$. Then the first letters of $m_1$ and $m_2$ are different, because $R_j$ is a cyclically reduced word. Hence, it is not possible that $m_1 = ms_1$ and that $m_2 = ms_2$.

Now assume that $m_1 = R_j^{n_1}a_j$, $m_2 = R_j^{n_2}a_j$, where $n_1, n_2 \in \mathbb{N}$, $n_1 \neq n_2$. Then at least one of $s_1$ and $s_2$ contains a cyclic shift of $R_j$. Since $G$ satisfies condition $C(2)$, every cyclic shift of $R_j$ is not a small piece. So, this is not possible, because $s_1$ and $s_2$ are small pieces.

Finally assume that  $m_1 = R_j^{-n_1}b_j$ and $m_2 = R_j^{-n_2}b_j$ where $n_1, n_2 \in \mathbb{N}$, $n_1 \neq n_2$. Similarly, at least one of $s_1$ and $s_2$ contains a cyclic shift of $R_j^{-1}$. Since $G$ satisfies condition $C(2)$, every cyclic shift of $R_j^{-1}$ is not a small piece. So, again this is not possible, because $s_1$ and $s_2$ are small pieces. Thus, we obtain $m_1\cdot s_1^{-1} \neq m_2 \cdot s_2^{-1}$. Right-sided Isolation Axiom for $\Rel$ is proved. Left-sided Isolation Axiom for $\Rel$ is checked in the same way.

\subsection{Equating of a binomial to a single monomial in a group algebra of a free group}
\label{equating_to_monomial_subsection}
Here we refer to paper~\cite{AKPR}.

Let $\Fr$ be a free group with at least four free generators. We fix a set of free generators of $\Fr$. Let $w \in \Fr$ be an arbitrary cyclically reduced primitive word. Let $x$ and $y$ be letters from the set of free generators of $\Fr$ such that the initial and the final letter of $w$ and $w^{-1}$ differ from $x^{\pm 1}$ and $y^{\pm 1}$. Consider the word
\begin{equation*}
v = x^{n_1}y x^{n_1 + 1}y \cdots x^{n_2}y, \ n_1, n_2 \in \mathbb{N},
\end{equation*}
such that $\vert w \vert \ll n_1 \ll n_2$ (namely, $n_1 - \vert w \vert > 0$ and $n_2 - n_1 \geqslant 21$). This are assumptions that we fix for Subsection~\ref{equating_to_monomial_subsection}.

The word $v$ exhibits small cancellation properties, because a subword of $v^m$, $m \in \mathbb{Z}$, containing at least two letters $y^{\pm 1}$, appears in $v^{m}$ uniquely modulo a shift by multiple of $\vert v\vert$. Since the initial and the final letter of $w$ and $w^{-1}$ differ from $x^{\pm 1}$ and $y^{\pm 1}$, we obtain that $v^{\pm 1}w^{\pm 1}$ and $w^{\pm 1}v^{\pm 1}$ with any combination of signs have no cancellations.

We consider a group algebra $\Galg$, where $\fld$ is a field. Let $\Ideal$ be an ideal of $\Galg$ generated by the polynomial $v^{-1} - 1 - w$ as an ideal. That is,
\begin{equation*}
\Ideal = \langle v^{-1} - 1 - w \rangle_i.
\end{equation*}

Let us show how to construct a set of generators of $\Ideal$ as an ideal, starting from $v^{-1} - 1 - w$, that satisfies Compatibility Axiom, Small Cancellation Axiom and Isolation Axiom. Let us notice that in paper~\cite{AKPR} we take $k = \mathbb{Z}_2$ in order to simplify calculations. However, a set of generators of $\Ideal$ that we will construct satisfies all necessary conditions over an arbitrary field.

~\paragraph*{$1^{\mathrm{o}}$} Let us put $\Rel_1 = \lbrace v^{-1} - w - 1\rbrace$. One can see that $\Rel_1$ does not satisfy Compatibility Axiom. We can do the following natural procedure. We start adding to $\Rel_1$ relations of the form $z_1^{-1}\cdot (v^{-1} - 1 - w)$ such that $z_1$ is a prefix of $v^{-1}$ or $w$, and of the form $(v^{-1} - 1 - w)\cdot z_2^{-1}$ such that $z_2$ is a suffix of $v^{-1}$ or $w$. After that we repeat the same procedure with the obtained set of relations, then again with the set of relations obtained after the second step, etc. Then, by construction, the union of all sets obtained in this way satisfies Compatibility Axiom. However, we make another more explicit and transparent procedure. Let us explain this procedure in detail.

Similarly to a Cayley graph of a group, we consider oriented graphs with edges marked by generators of the group $\Fr$. Take such a graph of the form
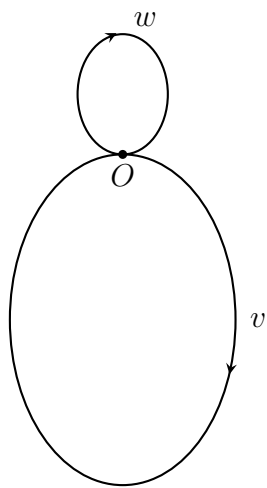
\begin{equation}
\begin{tikzpicture}
\begin{scope}
\coordinate (O) at ($(0,0)+(90:1.5 and 2.2)$);

\draw[black, thick, reversearrow=0.95] (0, 0) ellipse (1.5 and 2.2) node at (1.8, 0) {$v$};
\node[circle,fill,inner sep=1.2] at (O) {};
\path let \p1 = (O) in node at (\x1, \y1 - 8) {$O$};

\draw[black, thick, reversearrow=0.3] (0, 2.2+0.8) ellipse (0.6 and 0.8) node at (0.3, 2.2+1.8) {$w$};



\end{scope}
\end{tikzpicture}
\label{v_diagram_picture}
\end{equation}
(the word $w$ is written on the small arc, the word $v$ is written on the big arc). Assume that we have an oriented path in the graph~\eqref{v_diagram_picture}. We mean that it is possible to pass edges in the positive and negative direction. When we go along this path, we can write down the mark of an edge if we pass the edge in the positive direction, and we can write down the inverse to the mark of an edge if we pass the edge in the negative direction. So, speaking formally, every connected pair of vertices in~\eqref{v_diagram_picture} is connected by two oriented edges: one is marked by a free generator of $\Fr$ and another one is marked by the same generator in power $-1$. Hence, vertex $O$ is of degree $4$ and all other vertices are of degree $2$.

A path in the graph~\eqref{v_diagram_picture} with the initial and the final vertex $O$ naturally corresponds to a monomial over $v^{\pm 1}, w^{\pm 1}$ and vice versa. So, we can match every Laurent polynomial over $v, w$ with a collection of such paths. In particular, the polynomial $v^{-1} - 1 - w$ corresponds to the set of three paths in~\eqref{v_diagram_picture} that start and end at the point~$O$.

We can consider paths in the graph~\eqref{v_diagram_picture} with arbitrary initial and final points, which do not have to be equal to $O$. For instance, let $v = v_iv_mv_f$.
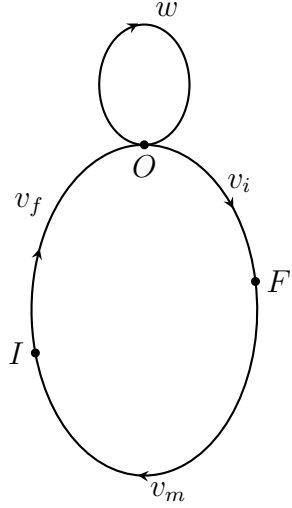
\begin{center}
\begin{tikzpicture}
\begin{scope}
\coordinate (O) at ($(90:1.5 and 2.2)$);
\coordinate (Si) at ($(10:1.5 and 2.2)$);
\coordinate (Se) at ($(195:1.5 and 2.2)$);

\coordinate (vi) at (0.96, 1.69);
\coordinate (vm) at (0.32, -2.15);
\coordinate (vf) at (-1.19, 1.33);

\draw[black, thick, arrow=0.6] (O) arc (90:10:1.5 and 2.2);
\node[right] at (vi) {$v_i$};
\draw[black, thick, arrow=0.6] (Si) arc (10:-360+195:1.5 and 2.2);
\node[below] at (vm) {$v_m$};
\draw[black, thick, arrow=0.4] (Se) arc (195:90:1.5 and 2.2);
\node[left, yshift=2] at (vf) {$v_f$};

\node[circle,fill,inner sep=1.2] at (O) {};
\path let \p1 = (O) in node at (\x1, \y1 - 8) {$O$};

\node[circle,fill,inner sep=1.2] at (Si) {};
\node[right] at (Si) {$F$};
\node[circle,fill,inner sep=1.2] at (Se) {};
\node[left] at (Se) {$I$};

\draw[black, thick, reversearrow=0.3] (0, 2.2+0.8) ellipse (0.6 and 0.8) node at (0.3, 2.2+1.8) {$w$};
\end{scope}
\end{tikzpicture}
\end{center}
Then monomials of the form $v_fM(v, w)v_i$, where $M(v, w)$ is a non-commutative reduced monomial over $v^{\pm 1}$ and $w^{\pm 1}$, correspond to paths in the above graph with the initial point~$I$ and the final point~$F$.

Assume that $z_1$ is a prefix of $v^{-1}$ or $w$. Then one can see that $z_1^{-1}\cdot (v^{-1} - 1 - w)$ corresponds to an agreed shifting of the initial point of the paths that correspond to $v^{-1} - 1 - w$ in graph~\eqref{v_diagram_picture}. Assume $z_2$ is a suffix of $v^{-1}$ or $w$. Similarly, one can see that $(v^{-1} - 1 - w)\cdot z_2^{-1}$ corresponds to an agreed shifting of the final point of the paths that correspond to $v^{-1} - 1 - w$ in graph~\eqref{v_diagram_picture}. Now we take all possible agreed shifts of the initial and the final points of the paths with the initial and the final point~$O$ that correspond to $v^{-1} - 1 - w$ in graph~\eqref{v_diagram_picture}, not only shifts by inverses of prefixes and suffixes of $v^{-1}$ and $w$, and add the corresponding polynomials to $\Rel_1$. Denote the obtained set of generators by $\Rel_2$. So,
\begin{align*}
\Rel_2 = \big\lbrace &\gamma M_L\cdot (v^{-1} - 1 - w)\cdot M_R \mid v = v_iv_mv_f, w = w_iw_mw_f, \gamma \in \fld,\\
&M_L \in \lbrace v_fM_1(v, w), v_i^{-1}M_1(v, w), w_fM_1(v, w), w_i^{-1}M_1(v, w)\rbrace,\\
&M_R \in \lbrace M_2(v, w)v_i, M_2(v, w)v_f^{-1}, M_2(v, w)w_i, M_2(v, w)w_f^{-1}\rbrace,\\
&M_1(v, w),  M_2(v, w) \textit{ are non-commutative reduced monomials over } v^{\pm 1}, w^{\pm 1}\big\rbrace.
\end{align*}
Clearly, $\Rel_2$ generates as an ideal the same ideal $\Ideal$.

Let us show that $\Rel_2$ satisfies Compatibility Axiom. Assume
\begin{equation*}
T = \gamma M_L\cdot (v^{-1} - 1 - w)\cdot M_R \in \Rel_2.
\end{equation*}
The monomial $M_L$ corresponds to a path in graph~\ref{v_diagram_picture} with some initial point $I$ and the final point $O$. The monomial $M_R$ corresponds to a path in graph~\ref{v_diagram_picture} with the initial point $O$ and some final point $F$. So, the polynomial $T$ corresponds to a collection of paths with the initial point $I$ and the final point $F$. Assume a letter $z^{-1}$ cancels from the left with some monomial in~$T$. Then $z^{-1}$ corresponds to a path with the final point $I$. Hence, $z^{-1}\cdot M_L$ corresponds to a path with the final point $O$. That is, $z^{-1}\cdot M_L$ belongs to the set $\lbrace v_fM(v, w), v_i^{-1}M(v, w), w_fM(v, w), w_i^{-1}M(v, w)\rbrace$. Therefore,
\begin{equation*}
z^{-1}\cdot T = \gamma (z^{-1}\cdot M_L) \cdot (v^{-1} - 1 - w)\cdot M_R \in \Rel_2.
\end{equation*}
Similarly, if a letter $z^{-1}$ cancels from the right with some monomial in $T$, we obtain  $T \cdot z^{-1} \in \Rel_2$.

Obviously, in the current case the set of monomials $\Mon$ for $\Rel_2$ consists of all subwords of non-commutative reduced monomials over $v^{\pm 1}, w^{\pm 1}$. We call such subwords \emph{$(v, w)$-generalized fractional powers}.

~\paragraph*{$2^{\mathrm{o}}$} Let $a$ be a subword of $M(v, w)$, where $M(v, w)$ is a non-commutative reduced monomial over $v^{\pm 1}, w^{\pm 1}$. Then, as we noticed above, $a$ corresponds to a path in graph~\eqref{v_diagram_picture}. There are two types of subwords of $M(v, w)$.
\begin{itemize}
\item
A subword that corresponds to a unique path in~\eqref{v_diagram_picture}. For example, $yx^ny$, $n_1 \leqslant n \leqslant n_2$.
\item
A subword that corresponds to more than one path in~\eqref{v_diagram_picture}. For example, $xy$.
\end{itemize}

\begin{proposition}
\label{sp_for_sc_binomial}
A word $c \in \Mon$ is not a small piece with respect to $\Rel_2$ (in the sence of Definition~\ref{sp}) if and only if $c$ corresponds to a unique path in graph~\eqref{v_diagram_picture}.
\end{proposition}
\begin{proof}
Assume $c$ corresponds to a unique path in graph~\eqref{v_diagram_picture} with the initial point $I$ and the final point $F$. Let us show that $c$ is not a small piece with respect to $\Rel_2$. Assume
\begin{align*}
&T_1 = \alpha\widehat{a}_1c\widehat{a}_2 + \sum_{j = 1}^{n_1}\alpha_j a_j \in \Rel_2,\\
&T_2 = \beta\widehat{b}_1c\widehat{b}_2 + \sum_{j = 1}^{n_2}\beta_j b_j \in \Rel_2.
\end{align*}
Then we need to show that
\begin{align*}
\widehat{b}_1\cdot \widehat{a}_1^{-1}\cdot T_1 = \alpha\widehat{b}_1c\widehat{a}_2 + \sum_{j = 1}^{n_1}\alpha_j \widehat{b}_1\cdot \widehat{a}_1^{-1}\cdot a_j \in \Rel_2,\\
T_1\cdot \widehat{a}_2^{-1} \cdot \widehat{b}_2 = \alpha\widehat{a}_1c\widehat{b}_2 + \sum_{j = 1}^{n_2}\alpha_j a_j \cdot \widehat{a}_2^{-1} \cdot \widehat{b}_2 \in \Rel_2.
\end{align*}

Since $\Rel_2$ satisfies Compatibility Axiom, we have $\widehat{a}_1^{-1}\cdot T_1\in \Rel_2$. So, $\widehat{a}_1^{-1}\cdot T_1 = \gamma M_L\cdot (v^{-1} - 1 - w)\cdot M_R$, where
\begin{align*}
&M_L \in \left\lbrace v_fM_1(v, w), v_i^{-1}M_1(v, w), w_fM_1(v, w), w_i^{-1}M_1(v, w)\right\rbrace,\\
&M_R \in \left\lbrace M_2(v, w)v_i, M_2(v, w)v_f^{-1}, M_2(v, w)w_i, M_2(v, w)w_f^{-1}\right\rbrace,\\
&\textit{where } M_1(v, w),  M_2(v, w) \textit{ are non-commutative}\\
&\textit{reduced monomials over } v^{\pm 1}, w^{\pm 1}.
\end{align*}
The key step is as follows. Assume $M_L$ corresponds to a path in graph~\eqref{v_diagram_picture} with some initial point $I_1$ and the final point $O$ and $M_R$ corresponds to a path in graph~\eqref{v_diagram_picture} with the initial point $O$ and some final point $F_1$. Then, on the one hand, all monomials of $\gamma M_L\cdot (v^{-1} - 1 - w)\cdot M_R$ correspond to a set of paths with the initial point $I_1$ and the final point $F_1$. On the other hand, the monomial $c$ corresponds to the unique path with the initial point~$I$. Therefore, every path that corresponds to the monomial $c\widehat{a}_2$ has the initial point $I$. So, since $c\widehat{a}_2$ is a monomial of $\gamma M_L\cdot (v^{-1} - 1 - w)\cdot M_R$, we obtain that $I_1$ has to be equal to $I$. That is, the monomial $M_L$ corresponds to a path in graph~\eqref{v_diagram_picture} with the initial point $I$ and the final point $O$.

The monomial $\widehat{b}_1c \in \Mon$ corresponds to some path in graph~\eqref{v_diagram_picture}. Since the monomial $c$ corresponds to the unique path with the initial point $I$, we obtain that $\widehat{b}_1$ corresponds to a path with some initial point $I_2$ and the final point $I$. Hence, the monomial $\widehat{b}_1\cdot M_L$ corresponds to a path with the initial point $I_2$ and the final point $O$. Therefore,
\begin{align*}
&\widehat{b}_1\cdot M_L \in \lbrace v_fM(v, w), v_i^{-1}M(v, w), w_fM(v, w), w_i^{-1}M(v, w)\rbrace,\\
&\textit{where } M(v, w) \textit{ is a non-commutative}\\
&\textit{reduced monomial over } v^{\pm 1}, w^{\pm 1}.
\end{align*}
So, $\gamma (\widehat{b}_1\cdot M_L)\cdot (v^{-1} - 1 - w)\cdot M_R \in \Rel_2$. Hence, we obtain
\begin{equation*}
\widehat{b}_1\cdot \widehat{a}_1^{-1}\cdot T_1 = \gamma (\widehat{b}_1\cdot M_L)\cdot (v^{-1} - 1 - w)\cdot M_R \in \Rel_2.
\end{equation*}
In the same way, one can show that $T_1\cdot \widehat{a}_2^{-1} \cdot \widehat{b}_2 \in \Rel_2$.

Let us show that if a word corresponds to more than one path in graph~\eqref{v_diagram_picture}, then this word is a small piece. The following words always correspond to more than one path:
\begin{enumerate}[label=(\arabic*)]
\item
\label{binomial_sp1}
$x^{n}$, $yx^{n}$, $x^{n}y$ for $n_1 \leqslant n < n_2$;
\item
\label{binomial_sp2}
$x^{-n}$, $y^{-1}x^{-n}$, $x^{-n}y^{-1}$ for $n_1 \leqslant n < n_2$;
\item
\label{binomial_sp3}
subwords of $w^{\pm m}$ that appear in $w^{\pm m}$ more than once modulo the period $w^{\pm 1}$;
\end{enumerate}
The following words may correspond to more than one path:
\begin{enumerate}[label=(\arabic*)]
\setcounter{enumi}{3}
\item
\label{binomial_sp4}
$v_fw_i,\ v_fw_f^{-1},\ v_i^{-1}w_i,\ v_i^{-1}w_f^{-1},$\\
$w_fv_i,\ w_i^{-1}v_i,\ w_fv_f^{-1},\ w_i^{-1}v_f^{-1},$\\
where $w_i$ is a prefix of $w$, $w_f$ is a suffix of $w$, $v_i$ is a prefix of $v$, $v_f$ is a suffix of $v$ such that $v_i$ and $v_f$ contain not more than one letter $y$.
\end{enumerate}
Let $c$ be a word of the form~\ref{binomial_sp1}---\ref{binomial_sp4}. Notice that for every type of $c$ it is enough to find a specific example of two polynomials from $\Rel_2$ that satisfy Definition~\ref{sp}. Let us take $c$ of the form~\ref{binomial_sp1} and construct a corresponding example of two polynomials. Examples for types~\ref{binomial_sp2}---\ref{binomial_sp4} can be produced similarly.

Let $v = v_iv_mv_f = v_i^{\prime}v_m^{\prime} v_f^{\prime}$ be two different fragmentations of $v$ such that $v_m = v_m^{\prime}$ are words of the form~\ref{binomial_sp1}. Since the fragmentations are different, we see that $v_i \neq v_i^{\prime}$ and $v_f \neq v_f^{\prime}$. Clearly, we can take such different fragmentations of $v$ for every word type~\ref{binomial_sp1}. We have
\begin{equation*}
v\cdot (v^{-1} - w - 1) = 1 - vw - v \in \Rel_2.
\end{equation*}
We take $T_1 = T_2 \in \Rel_2$, but consider different occurrences of $c = v_m = v_m^{\prime}$ in their monomials:
\begin{align*}
&T_1 = 1 - vw - v = 1 -  v_iv_mv_f - v_iv_mv_f w,\\
&T_2 = 1 - vw - v = 1 -  v_i^{\prime}v_m^{\prime}v_f^{\prime} - v_i^{\prime}v_m^{\prime}v_f^{\prime} w.
\end{align*}
Looking at the explicit form of $v$, we see that for every word of type~\ref{binomial_sp1} we can take fragmentations of $v$ such that both $v_i$ and $v_i^{\prime}$ contain more than two letters $y$. Then $v_i$ corresponds to a unique path in~\eqref{v_diagram_picture} and $v_i^{\prime}$ corresponds to a unique path in~\eqref{v_diagram_picture}. Clearly, these paths differ only by the final points. Do be definite, assume that $\vert v_i \vert > \vert v_i^{\prime}\vert$.
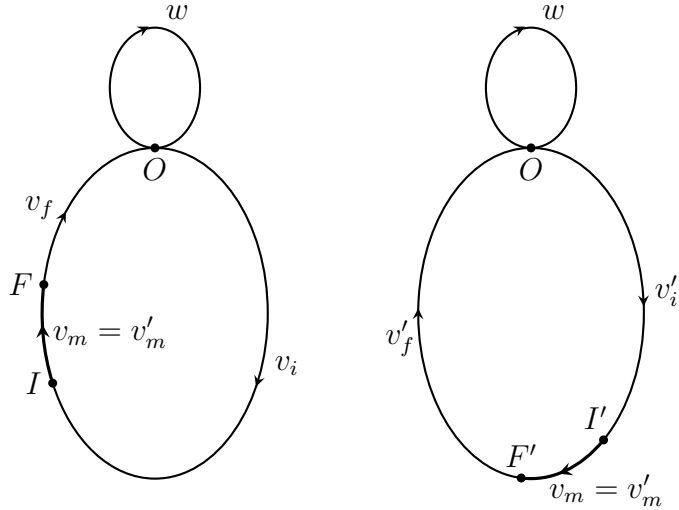
\begin{center}
\begin{tikzpicture}
\begin{scope}
\coordinate (O) at ($(90:1.5 and 2.2)$);
\coordinate (Si) at ($(-155:1.5 and 2.2)$);
\coordinate (Se) at ($(170:1.5 and 2.2)$);

\coordinate (vi) at (1.45, -0.7);
\coordinate (vm) at (-0.6, 0.1);
\coordinate (vf) at (-1.19, 1.33);

\draw[black, thick, arrow=0.5] (O) arc (90:-155:1.5 and 2.2);
\node[right] at (vi) {$v_i$};
\draw[black, very thick, arrow=0.6] (Si) arc (-155:-360+170:1.5 and 2.2);
\node[below] at (vm) {$v_m = v_m^{\prime}$};
\draw[black, thick, arrow=0.4] (Se) arc (170:90:1.5 and 2.2);
\node[left, yshift=2] at (vf) {$v_f$};

\node[circle,fill,inner sep=1.2] at (O) {};
\path let \p1 = (O) in node at (\x1, \y1 - 8) {$O$};

\node[circle,fill,inner sep=1.2] at (Si) {};
\node[left] at (Si) {$I$};
\node[circle,fill,inner sep=1.2] at (Se) {};
\node[left] at (Se) {$F$};

\draw[black, thick, reversearrow=0.3] (0, 2.2+0.8) ellipse (0.6 and 0.8) node at (0.3, 2.2+1.8) {$w$};
\end{scope}

\begin{scope}
\coordinate (O) at ($(5,0)+(90:1.5 and 2.2)$);
\coordinate (Si) at ($(5,0)+(-50:1.5 and 2.2)$);
\coordinate (Se) at ($(5,0)+(-95:1.5 and 2.2)$);

\coordinate (vi) at (5 + 1.5, 0.3);
\coordinate (vm) at (5 + 1, -2);
\coordinate (vf) at (5 + -1.4, -0.4);

\draw[black, thick, arrow=0.6] (O) arc (90:-50:1.5 and 2.2);
\node[right] at (vi) {$v_i^{\prime}$};
\draw[black, very thick, arrow=0.6] (Si) arc (-50:-95:1.5 and 2.2);
\node[below] at (vm) {$v_m = v_m^{\prime}$};
\draw[black, thick, arrow=0.5] (Se) arc (-95:90-360:1.5 and 2.2);
\node[left, yshift=2] at (vf) {$v_f^{\prime}$};

\node[circle,fill,inner sep=1.2] at (O) {};
\path let \p1 = (O) in node at (\x1, \y1 - 8) {$O$};

\node[circle,fill,inner sep=1.2] at (Si) {};
\node[above, xshift=-3] at (Si) {$I^{\prime}$};
\node[circle,fill,inner sep=1.2] at (Se) {};
\node[above] at (Se) {$F^{\prime}$};

\draw[black, thick, reversearrow=0.3] (5, 2.2+0.8) ellipse (0.6 and 0.8) node at (5.3, 2.2+1.8) {$w$};
\end{scope}
\end{tikzpicture}
\end{center}

We consider
\begin{align*}
v_i^{\prime} \cdot v_i^{-1}\cdot T_1 &= v_i^{\prime} \cdot v_i^{-1}\cdot (1 -  v - vw) = v_i^{\prime} \cdot v_i^{-1}\cdot (1 - v_iv_mv_f  - v_iv_mv_f w ) =\\
&= v_i^{\prime} \cdot v_i^{-1} - v_i^{\prime}v_mv_f  - v_i^{\prime}v_mv_f w.
\end{align*}
Let us show that $v_i^{\prime}v_mv_f w \notin \Mon$. Assume the contrary. The word $w$ corresponds to the unique path in~\eqref{v_diagram_picture} with the beginning and the end at the point $O$. Hence, $v_i^{\prime}v_mv_f$, which is a prefix of  $v_i^{\prime}v_mv_f w$, corresponds to a path with the end at the point $O$. Since $v_i^{\prime}$ corresponds to the unique path with the beginning point $O$, every path that corresponds to $v_i^{\prime}v_mv_f$ starts at the point $O$. Since the path that corresponds to $v_i^{\prime}$ is contained in $v$-arc and ends at point $I^{\prime} \neq O$, we obtain that every path that corresponds to $v_i^{\prime}v_mv_f$ is contained in $v$-arc. We have $\vert v\vert = \vert v_i\vert + \vert v_mv_f\vert = \vert v_i^{\prime}\vert + \vert v_m^{\prime}v_f^{\prime}\vert$. Since $\vert v_i \vert > \vert v_i^{\prime}\vert$, we get $\vert v_mv_f\vert < \vert v_m^{\prime}v_f^{\prime}\vert$. Combining these observations, we get that $v_i^{\prime}v_mv_f$ can not end at the point $O$ because of length. A contradiction. So, $v_i^{\prime}v_mv_f w \notin \Mon$.

Since $v_i^{\prime}v_mv_f w \notin \Mon$, we obtain that all the more $v_i^{\prime} \cdot v_i^{-1}\cdot T_1 \notin \Rel_2$. Thus, $v_m = v_m^{\prime}$ is a small piece with respect to~$\Rel_2$. Proposition~\ref{sp_for_sc_binomial} is proved.
\end{proof}

So, there are three types of words that are necessarily non-small pieces for an arbitrary $w$:
\begin{enumerate}[label=(NSP\arabic*)]
\item
\label{sp1}
subwords of $v^{\pm 1}$ that contain at least two letters $y^{\pm 1}$ and the subword~$x^{\pm n_2}$;
\item
\label{sp2}
$w^{\pm n}$ or its cyclic shift;
\item
\label{sp3}
subwords of $w^{\pm n}$ that appear in $w^{\pm n}$ uniquely modulo the period~$w^{\pm 1}$.
\end{enumerate}
Notice that words of the form
\begin{align*}
&v_fw_i,\ v_fw_f^{-1},\ v_i^{-1}w_i,\ v_i^{-1}w_f^{-1},\\
&w_fv_i,\ w_i^{-1}v_i,\ w_fv_f^{-1},\ w_i^{-1}v_f^{-1},
\end{align*}
where $v_i$ is a prefix of $v$, $v_f$ is a suffix of $v$, $w_i$ is a prefix of $w$, $w_f$ is a suffix of $w$, may be non-small pieces for some $w$ as well.

One can show that $\Rel_2$ satisfies Small Cancellation Axiom with the constant $\tau \geqslant 10$ and Isolation Axiom. So, we can work with $\Rel_2$. However, one can see that $\Add(\Rel_2) \neq \Rel_2$. Indeed, $1 - v - vw, v - v^2 - v^2w \in \Rel_2$, wherein
\begin{equation*}
1 - v - vw - (v - v^2 - v^2w) = 1 - vw + v^2 + v^2w \notin \Rel_2,
\end{equation*}
because all polynomials from $\Rel_2$ are trinomials. Let us further extend $\Rel_2$ in order to produce additively closed set of generators with the same set of small pieces.

\begin{lemma}
\label{polynomial_in_ideal}
Let $k(t)$ be the field of rational functions in one variable $t$ over the field~$\fld$. Let $P(x_1, x_2)$ be a non-commutative Laurent polynomial over the field $\fld$ such that $P((1 + t)^{-1}, t) = 0$ as an element of $k(t)$. Then $P(v, w) \in \Ideal$.
\end{lemma}
\begin{proof}
First assume that $P(x_1, x_2)$ is an arbitrary non-commutative Laurent polynomial over the field $\fld$ without any additional conditions. Let
\begin{equation*}
P(x_1, x_2) = \sum_{i,j}\eta_{ij}x_1^{n_i}x_2^{n_j} + \sum_{i,j}\zeta_{ij}x_2^{k_i}x_1^{k_j},\ n_i, n_j, k_i, k_j \in \mathbb{Z}.
\end{equation*}
Then
\begin{equation*}
P((1 + t)^{-1}, t) = \sum_{i,j}\eta_{ij}\frac{t^{n_j}}{(1 + t)^{n_i}} + \sum_{i,j}\zeta_{ij}\frac{t^{k_i}}{(1 + t)^{k_j}}.
\end{equation*}
Let us decompose every fraction of the form $\frac{t^n}{(1 + t)^m}$ in the above expression to elementary fractions and make all possible additive cancellations. Then we obtain
\begin{equation}
\label{elementary_fractions_decomposition}
P((1 + t)^{-1}, t) = \sum_r \gamma_{r}\frac{1}{(1 + t)^{n_r}} + \sum_s\delta_{s}t^{n_s}.
\end{equation}

Let us show that the last equality holds in $\Qalg$ if we replace $t$ by $w + \Ideal$ and $(1 + t)^{-1}$ by $v + \Ideal$. That is, we will show that
\begin{equation*}
P(v, w) + \Ideal = \sum_{i,j}\eta_{ij}v^{n_i}w^{n_j} + \sum_{i,j}\zeta_{ij}w^{k_i}v^{k_j} + \Ideal = \sum_r \gamma_{r}v^{n_r} + \sum_s\delta_{s}w^{n_s} + \Ideal.
\end{equation*}
Indeed, when we decompose $\frac{t^n}{(1 + t)^m}$ to elementary fractions, we use only identities that hold in an arbitrary ring (adding and subtracting  the same value and the binomial formula) and the equality
\begin{equation*}
(1 + t)^n\cdot \frac{1}{(1 + t)^n} = \frac{1}{(1 + t)^n}\cdot (1 + t)^n = 1, \ n \in \mathbb{N}.
\end{equation*}
Clearly, we have the same equality in $\Qalg$ if we replace $(1 + t)^{-1}$ by $v + \Ideal$ and $t$ by $w + \Ideal$. That is, $v^n\cdot (1 + w)^n + \Ideal = (1 + w)^n\cdot v^n + \Ideal = 1 + \Ideal$ in $\Qalg$. Therefore, if
\begin{equation*}
\frac{t^n}{(1 + t)^m} = \sum_r\alpha_{n_r}\frac{1}{(1 + t)^{n_r}} + \sum_s\beta_{s}t^{n_s}
\end{equation*}
is the decomposition to elementary fractions, the corresponding equality
\begin{equation*}
v^nw^m + \Ideal = \sum_r\alpha_{n_r}v^{n_r} + \sum\beta_{s}w^{n_s} + \Ideal
\end{equation*}
holds in $\Qalg$. So, we have in $\Qalg$
\begin{align}
\begin{split}
\label{elementary_fractions_decomposition_in_quotient}
P(v, w) + \Ideal &= \sum_{i,j}\eta_{ij}v^{n_i}w^{n_j} + \sum_{i,j}\zeta_{ij}w^{k_i}v^{k_j} + \Ideal =\\
&= \sum_r \gamma_{r}v^{n_r} + \sum_s\delta_{s}w^{n_s} + \Ideal.
\end{split}
\end{align}

Now assume that $P((1 + t)^{-1}, t) = 0$ as an element of $k(t)$. It is well known that $\left\lbrace \frac{1}{(1 + t)^n}, t^n \mid n \in \lbrace 0\rbrace\cup\mathbb{N} \right\rbrace$ is a set of linearly independent elements in~$k(t)$. So, since $P((1 + t)^{-1}, t) = 0$, we obtain that every $\gamma_r = 0$ and every $\delta_s = 0$ in~\eqref{elementary_fractions_decomposition}. Hence,~\eqref{elementary_fractions_decomposition_in_quotient} implies that $P(v, w) + \Ideal = \Ideal$ in $\Qalg$. Thus, $P(v, w) \in \Ideal$.
\end{proof}

We define
\begin{align}
\begin{split}
\Rel = \big\lbrace & M_L\cdot P(v, w)\cdot M_R \mid v = v_iv_mv_f, w = w_iw_mw_f,\\
&M_L \in \lbrace v_fM_1(v, w), v_i^{-1}M_1(v, w), w_fM_1(v, w), w_i^{-1}M_1(v, w)\rbrace,\\
&M_R \in \lbrace M_2(v, w)v_i, M_2(v, w)v_f^{-1}, M_2(v, w)w_i, M_2(v, w)w_f^{-1}\rbrace,\\
&M_1(v, w), M_2(v, w) \textit{ are non-commutative reduced monomials over } v^{\pm 1}, w^{\pm 1},\\
&P(x_1, x_2) \textit{ is a non-commutative Laurent polynomial over the field } k,\\
&P((1 + t)^{-1}, t) = 0 \textit{ as an elelent of } k(t)\big\rbrace.
\label{binom_additively_closed_relations}
\end{split}
\end{align}
Since $P(v, w) \in \Ideal$, we get that $\Rel$ generates as an ideal the same ideal $\Ideal$. Using the same argument that is used for $\Rel_2$, one can show that $\Rel$ satisfies Compatibility Axiom. Using the same argument as in Proposition~\ref{sp_for_sc_binomial}, one can show that the set of small pieces with respect to $\Rel$ is the same as the set of small pieces with respect to $\Rel_2$.

\begin{lemma}
\label{binom_additively_closed_system}
We have $\Add(\Rel) = \Rel$.
\end{lemma}
\begin{proof}
As usual, $\tau$ is a natural number $\geqslant 10$. Let
\begin{equation*}
T_1 = \alpha c + \sum_{j = 1}^{n_1} \alpha_j a_j,\ T_2 = \beta c + \sum_{j = }^{n_2} \beta_j b_j \in \Rel,
\end{equation*}
where $\SPM(c) \geqslant \tau - 2$. Let us show that $\gamma_1T_1 + \gamma_2T_2 \in \Rel$ for arbitrary $\gamma_1, \gamma_2 \in \fld$.

By the definition of $\Rel$, we have
\begin{equation*}
T_1 = M_L^{(1)}\cdot P_1(v, w)\cdot M_R^{(1)},\ T_2 = M_L^{(2)}\cdot P_2(v, w)\cdot M_R^{(2)}
\end{equation*}
(see~\eqref{binom_additively_closed_relations} for a definition of $M_L^{(1)}$, $P_1(v, w)$, $M_R^{(1)}$, $M_L^{(2)}$, $P_2(v, w)$, $M_R^{(2)})$). Since $\SPM(c) \geqslant \tau - 2$, $c$ is not a small piece. Therefore, it follows from Proposition~\ref{sp_for_sc_binomial} that $c$ corresponds to a unique path in graph~\eqref{v_diagram_picture} (particularly, with unique initial and final points). Denote the initial point of this path by $I$ and the final point by $F$. Notice that for incident monomials there exist paths in~\eqref{v_diagram_picture} with the same initial and the same final points. Therefore, $a_j$, $j = 1, \ldots, n_1$, and $b_j$, $j = 1, \ldots, n_2$, correspond to paths in graph~\eqref{v_diagram_picture} with the initial point $I$ and the final point $F$. Hence, there exits $\widetilde{M}_L \in \Mon$ that corresponds to a path with the initial point $O$ and the final point $I$, and there exists $\widetilde{M}_R \in \Mon$ that corresponds to a path with the initial point $F$ and the final point $O$ such that $\widetilde{M}_L\cdot T_1\cdot \widetilde{M}_R$ and $\widetilde{M}_L\cdot T_2\cdot \widetilde{M}_R$ correspond to collections of paths with the initial point $O$ and the final point $O$.

Consider monomials $\widetilde{M}_L\cdot M_L^{(1)}$, $M_R^{(1)}\cdot \widetilde{M}_R$, $\widetilde{M}_L\cdot M_L^{(2)}$, $M_R^{(2)}\cdot \widetilde{M}_R$. Since $c$ corresponds to a unique path in graph~\eqref{v_diagram_picture} with the initial point $I$ and the final point $F$, this implies that $M_L^{(1)}$ and $M_L^{(2)}$ correspond to paths with the initial point $I$ and the final point $O$. Similarly, $M_R^{(1)}$ and $M_R^{(2)}$ correspond to paths with the initial point $O$ and the final point $F$. Therefore, it follows from the above definition of  $\widetilde{M}_L$ and $\widetilde{M}_R$ that
\begin{align*}
&\widetilde{M}_L\cdot M_L^{(1)} = M_L^{(1)}(v, w),\ M_R^{(1)}\cdot \widetilde{M}_R = M_R^{(1)}(v, w),\\
&\widetilde{M}_L\cdot M_L^{(2)} = M_L^{(2)}(v, w),\ M_R^{(2)}\cdot \widetilde{M}_R = M_R^{(2)}(v, w)
\end{align*}
(possibly after the cancellations), where $M_L^{(1)}(x_1, x_2)$, $M_R^{(1)}(x_1, x_2)$, $M_L^{(1)}(x_1, x_2)$, $M_R^{(1)}(x_1, x_2)$ are non-commutative monomials in $x_1^{\pm 1}$, $x_2^{\pm 1}$.

Consider non-commutative Laurent polynomials
\begin{align*}
&Q_1(x_1, x_2) = M_L^{(1)}(x_1, x_2)\cdot P_1(x_1, x_2)\cdot M_R^{(1)}(x_1, x_2),\\
&Q_2(x_1, x_2) = M_L^{(2)}(x_1, x_2)\cdot P_2(x_1, x_2)\cdot M_R^{(2)}(x_1, x_2).
\end{align*}
Combining the above equalities, we obtain
\begin{align*}
Q_1(v, w) &=  M_L^{(1)}(v, w)\cdot P_1(v, w)\cdot M_R^{(1)}(v, w) =\\
&= \widetilde{M}_L\cdot M_L^{(1)}\cdot P_1(v, w) \cdot M_R^{(1)}\cdot \widetilde{M}_R = \widetilde{M}_L\cdot T_1\cdot \widetilde{M}_R,\\
Q_2(v, w) &=  M_L^{(2)}(v, w)\cdot P_2(v, w)\cdot M_R^{(2)}(v, w) =\\
&= \widetilde{M}_L\cdot M_L^{(2)}\cdot P_2(v, w) \cdot M_R^{(2)}\cdot \widetilde{M}_R = \widetilde{M}_L\cdot T_2\cdot \widetilde{M}_R.
\end{align*}

Since $P_1((1 + t)^{-1}, t) = 0$ and $P_2((1 + t)^{-1}, t) = 0$, we get $Q_1((1 + t)^{-1}, t) = 0$ and $Q_2((1 + t)^{-1}, t) = 0$. Let us take
\begin{equation*}
Q(x_1, x_2) = \gamma_1Q_1(x_1, x_2) + \gamma_2Q_2(x_1, x_2).
\end{equation*}
Then, evidently, $Q((1 + t)^{-1}, t) = 0$. We have
\begin{align*}
\gamma_1T_1 + \gamma_2T_2 &= \gamma_1 \widetilde{M}_L^{-1}\cdot Q_1(v, w)\cdot \widetilde{M}_R^{-1} + \gamma_2 \widetilde{M}_L^{-1}\cdot Q_2(v, w)\cdot \widetilde{M}_R^{-1} =\\
&= \widetilde{M}_L^{-1}\cdot Q(v, w)\cdot \widetilde{M}_R^{-1}.
\end{align*}
Thus, we obtain $\gamma_1T_1 + \gamma_2T_2 \in \Rel$.
\end{proof}

It is proved in paper~\cite{AKPR} that Small Cancellation Axiom with a constant $\tau \geqslant 10$ holds for $\Rel$ (see~\cite{AKPR}, Proposition~$3.1$, Transversality Condition). In fact, we proved in~\cite{AKPR} even more.
\begin{proposition}[Transversality Condition]
\label{transversality_condition}
Let $T_1, \ldots, T_n \in \Rel$, and $\sum_{j = 1}^n\gamma_j T_j$ be non-zero element of $\Galg$. Then $\sum_{j = 1}^n\gamma_j T_j$ after additive cancellations contains a monomial $A$ that contains separate subwords of $v^{\pm m}$ such that they contain in total $\geqslant \tau + 1$ letters from the set $\lbrace y, y^{-1}\rbrace$ ($\tau \geqslant 10$).
\end{proposition}
Although in  paper~\cite{AKPR} we worked with $k = \mathbb{Z}_2$, the argument in Transversality Condition works for an arbitrary field with very small changes.

\begin{remark}
\label{measure_number_of_y}
Let us also notice that in~\cite{AKPR} we use a measure $\SPM^{\prime}$ on monomials of $\Mon$ that slightly differs from $\SPM$-measure. In order to define $\SPM^{\prime}$ we count only letters $y^{\pm 1}$ in subwords of $v^{\pm m}$. Namely, let $u$ be a subword of a reduced monomial over $v^{\pm1}, w^{\pm 1}$, then $\SPM^{\prime}(u)$ is equal to the number of letters $y^{\pm 1}$ in total in all maximal occurrences of subwords of $v^{\pm m}$ in $u$. So, all subwords of $w^{\pm n}$ have $\SPM^{\prime}$ equal to $0$, and all subwords of $v^{\pm 1}$ of the form $x^{\pm n}$ have $\SPM^{\prime}$ equal to $0$.
\end{remark}
It is possible to use measure $\SPM^{\prime}$, because Small Cancellation Axiom holds for $\Rel$ in a stronger form, which is stated in Proposition~\ref{transversality_condition}.

~\paragraph*{$3^{\mathrm{o}}$} Let us check Isolation Axiom for the set $\Rel$. Let $m_1, m_2, \ldots, m_k$ be a sequence of monomials of $\Mon$ such that all of them are of $\SPM$-measure $\geqslant \tau - 2$ and the consecutive monomials are incident. Since $\Add(\Rel) = \Rel$, it follows from Lemma~\ref{add_closed_incident_chain_of_incidend} that $m_1$ and $m_k$ are incident monomials. So, we need to check Isolation Axiom only for incident monomials.

Notice that incident monomials correspond to paths in~\eqref{v_diagram_picture} with the same initial and the same final points.

Let us check right-sided Isolation Axiom. Let $m_1, m_2 \in \Mon$ be incident monomials, $m_1 \neq m_2$, $\SPM(m_1) \geqslant \tau - 2$, $\SPM(m_2) \geqslant \tau - 2$. We will prove even a slightly stronger condition. Namely, we take $a \in \Mon$ such that $\SPM(a) \geqslant \tau - 2$, $m_1a, m_2a \notin \Mon$, $m_1$ is a maximal occurrence in $m_1a$, $m_2$ is a maximal occurrence in $m_2a$ (that is, we omit the last condition on $a$ in right-sided  Isolation Axiom). Let $s_1a$ be a maximal occurrence in $m_1a$ that contains $a$, and $s_2a$ be a maximal occurrence in $m_2a$ that contains $a$. Then we will show that
\begin{equation*}
m_1\cdot s_1^{-1} \neq m_2 \cdot s_2^{-1}.
\end{equation*}

Assume the contrary, namely, assume that $m_1\cdot s_1^{-1} = m_2 \cdot s_2^{-1} = m$. Then we have $m_1 = ms_1$, $m_2 = ms_2$. By definition, $s_1$ and $s_2$ are overlaps of maximal occurrences. Hence, $s_1$ and $s_2$ are small pieces. Since $\SPM(m_1) \geqslant \tau - 2$ and $\SPM(m_2) \geqslant \tau - 2$, we obtain  $\SPM(m) \geqslant \tau - 2 - 1 = \tau - 3 \geqslant 7$. Therefore, $m$ is not a small piece. So, $m$ corresponds to a unique path in graph~\eqref{v_diagram_picture}. Hence, paths that correspond to $m_1$ and $m_2$ start with one path, which corresponds to $m$. Since $m_1$ and $m_2$ are incident monomials, they correspond to paths in~\eqref{v_diagram_picture} with the same initial and the same final points. Combining these facts, we obtain that $s_1$ and $s_2$ correspond to paths in graph~\eqref{v_diagram_picture} with the same initial and the same final points.

Since $ms_1 \neq ms_2$, we get  $s_1 \neq s_2$. We can write $s_1$ and $s_2$ in the form $s_1 = ss_1^{\prime}$ and $s_2 = ss_2^{\prime}$, where $s_1^{\prime}$ and $s_2^{\prime}$ do not have a common prefix. Then, in the same way as above, we obtain that $ms$ corresponds to a unique path in graph~\eqref{v_diagram_picture} and that $s_1^{\prime}$ and $s_2^{\prime}$ correspond to paths in graph~\eqref{v_diagram_picture} with the same initial and the same final points.

On the one hand, $s_1^{\prime}$ and $s_2^{\prime}$ do not have a common prefix, on the other hand, $s_1^{\prime}$ and $s_2^{\prime}$ prolong without any cancellations the path that corresponds to $ms$. Therefore, $s_1^{\prime}$ and $s_2^{\prime}$ can not start at a vertex of degree~$2$.
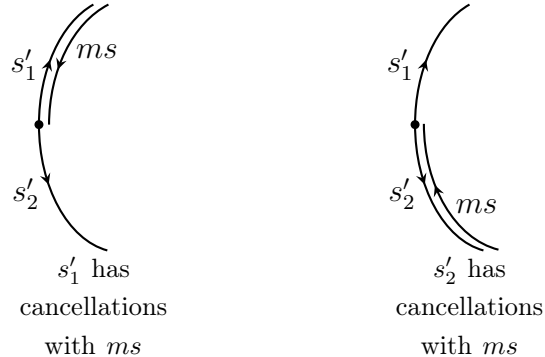
\begin{center}
\begin{tikzpicture}
\begin{scope}
\coordinate (O) at (0, 0);

\draw[black, thick, reversearrow=0.6] (O) arc (110:180:1.1 and 1.7) coordinate[pos = 1] (F) coordinate[pos = 0.5] (p1);
\draw[black, thick, arrow=0.4] (F) arc (180:260:1.1 and 1.7) coordinate[pos = 0.5] (p2);
\node[circle,fill,inner sep=1.2] at (F) {};
\node[left, yshift=-5, xshift=-2] at (p1) {$s_1^{\prime}$};
\node[left, yshift=5, xshift=-3] at (p2) {$s_2^{\prime}$};

\draw[black, thick, arrow=0.6] (0.2, 0) arc (110:180:1.2 and 1.7) coordinate[pos = 0.5] (mp);
\node[right] at (mp) {$ms$};

\node[text width=2cm, align=center] at (0, -4) {\footnotesize{$s_1^{\prime}$ has cancellations with $ms$}};
\end{scope}
\begin{scope}
\coordinate (O) at (5, 0);

\draw[black, thick, reversearrow=0.6] (O) arc (110:180:1.1 and 1.7) coordinate[pos = 1] (F) coordinate[pos = 0.5] (p1);
\draw[black, thick, arrow=0.4] (F) arc (180:260:1.1 and 1.7) coordinate[pos = 0.5] (p2) coordinate[pos = 1] (F1);
\node[circle,fill,inner sep=1.2] at (F) {};
\node[left, yshift=-5, xshift=-2] at (p1) {$s_1^{\prime}$};
\node[left, yshift=5, xshift=-3] at (p2) {$s_2^{\prime}$};

\draw[black, thick, arrow=0.6]  ($(F1) + (0.2, 0.01)$) arc (260:180:1.2 and 1.7) coordinate[pos = 0.5] (mp);
\node[right] at (mp) {$ms$};

\node[text width=2cm, align=center] at (5, -4) {\footnotesize{$s_2^{\prime}$ has cancellations with $ms$}};
\end{scope}

\end{tikzpicture}
\end{center}
So, $s_1^{\prime}$ and $s_2^{\prime}$ start at point~$O$.

Let us calculate possible forms of $s_1^{\prime}$ and $s_2^{\prime}$. Since $s_1$ and $s_2$ are small pieces, we get that $s_1^{\prime}$ and $s_2^{\prime}$ are small pieces as well. It follows from Definition~\ref{sp} that a small piece can not contain non-small pieces as subwords. Hence, $s_1^{\prime}$ and $s_2^{\prime}$ do not contain subwords that are non-small pieces. Recall that $s_1^{\prime}$ and $s_2^{\prime}$ end at the same point. Therefore, $s_1^{\prime}$ and $s_2^{\prime}$ can end only inside $w$-arc and be of the form $w_i$ and $w_f^{-1}$, where $w = w_iw_f$ and both $w_i$ and $w_f$ are non-empty. Otherwise, at least one of $s_1^{\prime}$ and $s_2^{\prime}$ contains $w^{\pm 1}$ or a subword of $v^{\pm 1}$ of $\SPM$-measure $\geqslant \frac{1}{2}\SPM(v)$. Denote their final point by $F$. To be definite, assume that $s_1^{\prime} = w_i$ and $s_2^{\prime} = w_f^{-1}$.
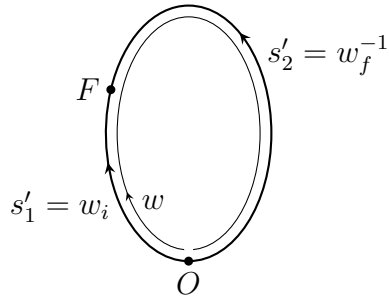
\begin{center}
\begin{tikzpicture}
\coordinate (O) at (0, 0);

\draw[black, thick, arrow=0.6] (O) arc (-90:160:1.1 and 1.7) coordinate[pos = 1] (F) coordinate[pos = 0.5] (wf);
\draw[black, thick, reversearrow=0.4] (F) arc (160:270:1.1 and 1.7) coordinate[pos = 0.5] (wi);
\node[circle,fill,inner sep=1.2] at (O) {};
\node[circle,fill,inner sep=1.2] at (F) {};
\node[right, yshift=2] at (wf) {$s_2^{\prime} = w_f^{-1}$};
\node[left] at (wi) {$s_1^{\prime} = w_i$};
\node[below] at (O) {$O$};
\node[left] at (F) {$F$};

\draw[black, arrow=0.15] (0, 0.15) arc (-90:-450:0.95 and 1.55) coordinate[pos = 0.15] (w);
\node[right] at (w) {$w$};
\node[white, circle,fill,inner sep=1.2] at (0, 0.15) {};
\end{tikzpicture}
\end{center}
In particular, we obtain that $s_1^{\prime}$ and $s_2^{\prime}$ have no common suffix.

Now let us look at $a$ and consider $s_1^{\prime}$ and $s_2^{\prime}$ as prolongations of $a$ to the left. Since $s_1^{\prime}$ and $s_2^{\prime}$ have no common suffix and they prolong $a$ from the left without any cancellations, we get that $a$ can not start at a vertex of degree~$2$. Hence, $a$ starts at point $O$. By the initial assumptions, $s_1^{\prime} = w_i$ starts with a letter different from $x^{\pm 1}$ and $y^{\pm 1}$. So, since $\vert s_1^{\prime}\vert < \vert w\vert < \vert v\vert$, $s_1^{\prime}$ can be contained only inside $w$-arc. Similarly, $s_2^{\prime} = w_f^{-1}$ starts with a letter different from $x^{\pm 1}$ and $y^{\pm 1}$. Hence, since $\vert s_2^{\prime}\vert < \vert w\vert < \vert v\vert$, $s_2^{\prime}$ can be contained only inside $w$-arc. Since $s_1^{\prime}$ and $s_2^{\prime}$ have no common suffix, they end by two different edges that come in vertex~$O$. Since
\begin{equation*}
\vert s_1^{\prime}\vert + \vert s_2^{\prime}\vert = \vert w_i\vert + \vert w_f^{-1}\vert = \vert w\vert,
\end{equation*}
we obtain that $s_1^{\prime}$ and $s_2^{\prime}$ start at one point at $w$-arc. Denote this point by $I$.
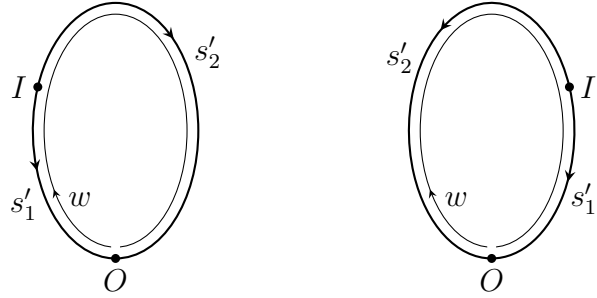
\begin{center}
\begin{tikzpicture}
\begin{scope}
\coordinate (O) at (0, 0);

\draw[black, thick, reversearrow=0.6] (O) arc (-90:160:1.1 and 1.7) coordinate[pos = 1] (I) coordinate[pos = 0.5] (wf);
\draw[black, thick, arrow=0.4] (I) arc (160:270:1.1 and 1.7) coordinate[pos = 0.5] (wi);
\node[circle,fill,inner sep=1.2] at (O) {};
\node[circle,fill,inner sep=1.2] at (I) {};
\node[right, yshift=2] at (wf) {$s_2^{\prime}$};
\node[left] at (wi) {$s_1^{\prime}$};
\node[below] at (O) {$O$};
\node[left] at (I) {$I$};

\draw[black, arrow=0.15] (0, 0.15) arc (-90:-450:0.95 and 1.55) coordinate[pos = 0.15] (w);
\node[right] at (w) {$w$};
\node[white, circle,fill,inner sep=1.2] at (0, 0.15) {};
\end{scope}
\begin{scope}
\coordinate (O) at (5, 0);

\draw[black, thick, reversearrow=0.6] (O) arc (-90:20:1.1 and 1.7) coordinate[pos = 1] (I) coordinate[pos = 0.5] (wf);
\draw[black, thick, arrow=0.4] (I) arc (20:270:1.1 and 1.7) coordinate[pos = 0.5] (wi);
\node[circle,fill,inner sep=1.2] at (O) {};
\node[circle,fill,inner sep=1.2] at (I) {};
\node[right, yshift=2] at (wf) {$s_1^{\prime}$};
\node[left] at (wi) {$s_2^{\prime}$};
\node[below] at (O) {$O$};
\node[right] at (I) {$I$};

\draw[black, arrow=0.15] (5, 0.15) arc (-90:-450:0.95 and 1.55) coordinate[pos = 0.15] (w);
\node[right] at (w) {$w$};
\node[white, circle,fill,inner sep=1.2] at (5, 0.15) {};
\end{scope}
\end{tikzpicture}
\end{center}
So, we obtain that either $w = {s_1^{\prime}}^{-1}s_2^{\prime} = w_i^{-1}w_f^{-1}$, or $w = {s_2^{\prime}}^{-1}s_1^{\prime} = w_fw_i$. Notice that $w_fw_i$ is a cyclic shift of $w = w_iw_f$. Therefore, the equality $w = w_fw_i$ is not possible in $\Fr$, because $w$ is a primitive word.

Now consider the case $w = {s_1^{\prime}}^{-1}s_2^{\prime} = w_i^{-1}w_f^{-1}$. We have $w = w_iw_f = w_i^{-1}w_f^{-1}$. Since $w_iw_f$ and $w_i^{-1}w_f^{-1}$ have no cancellations inside, this implies  $w_i = w_i^{-1}$ and $w_f = w_f^{-1}$ in the free group $\Fr$. Since at least one of $w_i$ and $w_f$ is not equal to $1$, this is not possible. A contradiction. Thus, $m_1\cdot s_1^{-1} \neq m_2 \cdot s_2^{-1}$. So, right-sided Isolation Axiom holds for $\Rel$. Left-sided Isolation Axiom holds for $\Rel$ as well and is checked similarly.

\section{Table of notations}
\renewcommand{\arraystretch}{1.4}
\begin{tabular}{ p{4cm} p{7cm} p{3.5cm} }
$\Fr$ & \textit{the free group} & page~\pageref{fr_def} \\
$S$ & \textit{the set of free generators of $\Fr$} & page~\pageref{alphabet_def}\\
$\Galg$ & \textit{the group algebra of $\Fr$ over a field $k$} & page~\pageref{galg_def}\\
$\Rel$& \textit{the set of relations in $\Galg$} &page~\pageref{rel_def}\\
$\Ideal$ & \textit{the ideal generated by $\Rel$ (as an ideal)}  & page~\pageref{ideal_def}\\
$\Mon$ & \textit{the set of summands of elements from~$\Rel$}  & page~\pageref{mon_def}\\
$AB$ \textit{and} $A\cdot B$ & & page~\pageref{product_def}\\
$\SP$ & \textit{the set of small pieces with respect to $\Rel$} & page~\pageref{set_sp_def}, see also page~\pageref{alphabet_def},

Definition~\ref{sp}\\
\emph{$\SPM$-measure} & & page~\pageref{measure_def}\\
$\tau$ & \textit{a fixed natural number $\geqslant 10$} & page~\pageref{tau_def}\\
\emph{small cancellation condition $C(m)$} &  & pages~\pageref{cond_cm_1},~\pageref{cond_cm_2} \\
$\GDp$ & \textit{the set of layouts of multi-turns of members of the chart of all monomials} & page~\pageref{linear_dep_members}, see also page~\pageref{multiturn_def},

Definition~\ref{multiturn_def}\\
$\mo{U}$ & \textit{the set of maximal occurrences in a monomial $U$} & page~\pageref{mo_def}, see also page~\pageref{max_occurrence_def},

Definition~\ref{max_occurrence_def} \\
$a_h \mapsto a_j$ & \textit{a replacement of an occurrence $a_h$ by $a_j$ such that $a_h$ and $a_j$ are incident-monomials} & page~\pageref{incidents_replacement_def}, see also page~\pageref{incident_momom},

Definition~\ref{incident_momom}\\
$\fc{U}$ & \textit{the set of elements of $\mo{U}$ fully covered by other elements of $\mo{U}$} & page~\pageref{fc_def} \\
$\nfc{U}$ & \textit{the set of elements of $\mo{U}$ not fully covered by other elements of $\mo{U}$}  & page~\pageref{nfc_def} \\
$\mathcal{C}_i(U)$ & \textit{a covering of a monomial $U$} & page~\pageref{cov_def}
\end{tabular}
\begin{tabular}{ p{4cm} p{7cm} p{3.5cm} }
$\mathcal{C}_i(Z, U)$ & \textit{the set of elements of $\mathcal{C}_i(U)$ that intersects with $Z$, where $Z$ is an occurrence in $U$} & page~\pageref{cov_intersection_def}\\
$\mincov{U}$ & \textit{the size of a minimal covering of $U$} & page~\pageref{mincov_def} \\
$\longmo{U}$ & \textit{the set of elements of $\mo{U}$ of $\SPM$-measure $\geqslant 3$} & page~\pageref{longmo_def}\\
$\shortmo{U}$ & \textit{the set of elements of $\mo{U}$ of $\SPM$-measure $\leqslant 2$} & page~\pageref{shortmo_def}\\
$\virt{U}$ & \textit{the set of virtual members of the chart of $U$} & page~\pageref{virt_def}, see also page~\pageref{virtual_members_def},

Definition~\ref{virtual_members_def}\\
$\nvirt{U}$ & \textit{the number of virtual members of the chart of $U$} & page~\pageref{nvirt_def}\\
$\GDp^{\prime}$ & \textit{the set of layouts of multi-turns of virtual members of the chart of all monomials} & page~\pageref{linear_dep_members2}\\
$f(U)$ & \textit{$f$-characteristic of a monomial $U$} & page~\pageref{f_char_def},

Definition~\ref{f_char_def}\\
$a \leadsto b$ & \textit{a replacement of an occurrence $a_h$ in a monomial $U$ by $a_j$ such that $a_h$ and $a_j$ are $U$-incident-monomials} & page~\pageref{u_incident_replacements_def}, see also page~\pageref{U_incident_monomials},

Definition~\ref{U_incident_monomials}\\
$\DSpace{U}$ & \textit{the subspace of $\Galg$ linearly generated by all derived monomials of $U$} & page~\pageref{derived_monomials_space_def},

Definition~\ref{derived_monomials_space_def}, see also page~\pageref{derived_monomials},

Definition~\ref{derived_monomials}\\
$\DMUp{U}$ & \textit{the set of derived monomials of $U$ with $f$-characteristic equal to $f(U)$} & page~\pageref{upper_monomials}\\
$\DMLow{U}$ & \textit{the set of derived monomials of $U$ with $f$-characteristic smaller than $f(U)$} & page~\pageref{lower_monomials}\\
$\Ft_n(\Galg)$ & \textit{the filtration on $\Galg$ based on $f$-characteristic} & page~\pageref{filtration_def}\\
$\Low(Y)$ & & page~\pageref{low_space_def}
\end{tabular}
\begin{tabular}{ p{4.5cm} p{7cm} p{3.5cm} }
$\Dp(Y)$ & \textit{the space of dependencies on $Y$, where $Y$ is a linear subspace of $\Galg$} & page~\pageref{subspace_of_dependencies},

Definition~\ref{subspace_of_dependencies}\\
$\GDp^{\prime}(Y)$ & \textit{the set of layouts of multi-turns of all monomials from $Y$, where $Y$ is a linear subspace of $\Galg$} & page~\pageref{gdp_y_def}\\
$A_i[U]$, $ME_i[U]$,\par $ML_i[U]$, $D_i[U]$ & \textit{subspaces of $\Galg$ linearly generated by special subsets of $\Mon$ that depend on a monomial $U$} & page~\pageref{tens_prod_components}\\
$\mu[U]$ &  & page~\pageref{mu_def},

Definition~\ref{mu_def}\\
$\GDp^{(i)}[U]$ & \textit{all layouts of all multi-turns of the $i$-th virtual members of the chart of the monomials of $\DMUp{U}$} & page~\pageref{mu_properties}, see

Lemma~\ref{mu_properties}\\
$\EDp\DSpace{U}$ & \textit{the set of dependencies that come from monomials of $\DMUp{U}$} & page~\pageref{edp_def}\\
$\widetilde{\GDp}^{(i)}[U]$, $\widetilde{\GDp}[U]$ & & page~\pageref{tilde_gdp_def}\\
$\Gr_n$ & \textit{a graded component that corresponds to the filtration $\Ft_n$} & page~\pageref{grading_def}\\
$\Add(\Rel)$ & \textit{the additive closure of $\Rel$} & page~\pageref{add_def}\\
$\GDp^{\prime\prime}$ & \textit{the set of layouts of multi-turns of virtual members of the chart of all monomials that come from $\Add(\Rel)$} & page~\pageref{add_def}\\
$<_f$ & \textit{the total ordering of monomials based on $f$-characteristic} & page~\pageref{f_char_order_def},

Definition~\ref{f_char_order_def}\\
$\GreedyAlg(<_f, \Add(\Rel))$ & & page~\pageref{greedy_alg_rings},

Definition~\ref{greedy_alg_rings}
\end{tabular}

{\bf Acknowledgements.} The research of the first, second and the third authors was supported by  ISF grant 1994/20  and the Emmy Noether Research Institute for Mathematics. The research of the first  and the forth authors was also supported by ISF fellowship. The research of the second author was  supported by the Russian science foundation, grant  17-11-01377.

We are very grateful to I.Kapovich, B.Kunyavskii and D.Osin   for invaluable cooperation.

\end{document}